\documentclass[reqno]{amsart}

\usepackage{graphicx, amsmath, amssymb, amsfonts, amsthm, stmaryrd, amscd}
\usepackage[dvipsnames]{xcolor}
\usepackage{tikz}

\usepackage{enumerate}
\usepackage[normalem]{ulem}
\usepackage{color}

\usepackage{ifpdf} \ifpdf \usepackage{hyperref} \usepackage{navigator} \fi

\usepackage{xparse}

\usepackage[all]{xy} \usepackage{enumerate} \usetikzlibrary{matrix,arrows,decorations.pathmorphing}

\usepackage{mathtools} \DeclarePairedDelimiter{\paren}{(}{)} \DeclarePairedDelimiter{\abs}{\lvert}{\rvert}

\makeatletter \newcommand*{\transpose}{%
  {\mathpalette\@transpose{}}%
} \newcommand*{\@transpose}[2]{%
  \raisebox{\depth}{$\m@th#1\intercal$}%
} \makeatother

\makeatletter \newcommand*{\da@rightarrow}{\mathchar"0\hexnumber@\symAMSa 4B } \newcommand*{\da@leftarrow}{\mathchar"0\hexnumber@\symAMSa 4C } \newcommand*{\xdashrightarrow}[2][]{%
  \mathrel{%
    \mathpalette{\da@xarrow{#1}{#2}{}\da@rightarrow{\,}{}}{}%
  }%
} \newcommand{\xdashleftarrow}[2][]{%
  \mathrel{%
    \mathpalette{\da@xarrow{#1}{#2}\da@leftarrow{}{}{\,}}{}%
  }%
} \newcommand*{\da@xarrow}[7]{%
  \sbox0{$\ifx#7\scriptstyle\scriptscriptstyle\else\scriptstyle\fi#5#1#6\m@th$}%
  \sbox2{$\ifx#7\scriptstyle\scriptscriptstyle\else\scriptstyle\fi#5#2#6\m@th$}%
  \sbox4{$#7\dabar@\m@th$}%
  \dimen@=\wd0 %
  \ifdim\wd2 >\dimen@ \dimen@=\wd2 %
  \fi \count@=2 %
  \def\da@bars{\dabar@\dabar@}%
  \@whiledim\count@\wd4<\dimen@\do{%
    \advance\count@\@ne \expandafter\def\expandafter\da@bars\expandafter{%
      \da@bars \dabar@ }%
  }%
  \mathrel{#3}%
  \mathrel{%
    \mathop{\da@bars}\limits \ifx\\#1\\%
    \else _{\copy0}%
    \fi \ifx\\#2\\%
    \else ^{\copy2}%
    \fi }%
  \mathrel{#4}%
} \makeatother
 \DeclareMathOperator{\sgn}{sgn} \DeclareMathOperator{\SL}{SL}  \DeclareMathOperator{\glLie}{\mathfrak{g}\mathfrak{l}}   \DeclareMathOperator{\GL}{GL}        \DeclareMathOperator{\htt}{ht}        \DeclareMathOperator{\ad}{ad}
\DeclareMathOperator{\Ad}{Ad}      \DeclareMathOperator{\Hom}{Hom}    \DeclareMathOperator{\h}{h}   \DeclareMathOperator{\End}{End}  \DeclareMathOperator{\image}{image}  \def\div{\operatorname{div}} \def\eps{\varepsilon}        \DeclareMathOperator{\PGL}{PGL}            
   \DeclareMathOperator{\dom}{dom}      \DeclareMathOperator{\U}{U}
\DeclareMathOperator{\trace}{trace}                                    \DeclareMathOperator{\dist}{dist}          \DeclareMathOperator{\diag}{diag} \DeclareMathOperator{\Sym}{Sym}    \DeclareMathOperator{\Lie}{Lie}    \DeclareMathOperator{\stab}{stab}   \newcommand{\git}{/\!\!/}   \def\O{\operatorname{O}}   \DeclareMathOperator{\Opp}{Op}   
\DeclareMathOperator{\sym}{sym}          
\DeclareMathOperator{\fin}{fin}   \DeclareMathOperator{\rank}{rank}     \DeclareMathOperator{\Eis}{Eis} \DeclareMathOperator{\reg}{reg}       \DeclareMathOperator{\Ind}{Ind}   \DeclareMathOperator{\ord}{ord}      \DeclareMathOperator{\res}{res}        \DeclareMathOperator{\ev}{eval}  \DeclareMathOperator{\eval}{eval}                \DeclareMathOperator{\vol}{vol}      \DeclareMathOperator{\supp}{supp}

\theoremstyle{plain} \newtheorem{theorem} {Theorem}

\newtheorem{corollary} [theorem] {Corollary}

\newtheorem{proposition} [theorem] {Proposition}

\theoremstyle{definition} \newtheorem{definition} [theorem] {Definition}  \newtheorem{assumptions} [theorem] {Assumptions} \newtheorem{example} [theorem] {Example}      \newtheorem{remark} [theorem] {Remark} \newtheorem{notation} [theorem] {Notation}

\newtheoremstyle{itplain} 
{6pt} 
{5pt\topsep} 
{\itshape} 
{} 
{\itshape} 
{.}  
{5pt plus 1pt minus 1pt} 
{} 

\theoremstyle{itplain} 
\newtheorem{lemma}[theorem]{Lemma}

\newtheorem*{lemma*}{Lemma} \newtheorem*{proposition*}{Proposition} \newtheorem*{definition*}{Definition} \newtheorem*{example*}{Example}

 \newtheorem*{remark*} {Remark}

\numberwithin{equation}{section} \numberwithin{theorem}{section}

\usepackage{multicol} \makeindex

  

\renewcommand{\geq}{\geqslant} \renewcommand{\leq}{\leqslant}

\DeclareMathOperator{\HACH}{HC}

\begin{document}
\author{Paul D. Nelson}
\email{nelson.paul.david@gmail.com}
\address{Aarhus University, Denmark}
\title{Bounds for standard $L$-functions}

\begin{abstract}
  Let $\pi$ be a cuspidal automorphic representation of a general linear group over the rational numbers.  We establish a subconvex bound for the standard $L$-function of $\pi$ in the $t$-aspect.  More generally, we address the spectral aspect in the case of uniform parameter growth.
\end{abstract}

\maketitle

\setcounter{tocdepth}{1} \tableofcontents
\section{Introduction}\label{sec:introduction}
Let $\pi$ be a unitary cuspidal automorphic representation of $\GL_n$ over $\mathbb{Q}$.  The standard $L$-functions $L(\pi,s)$ ($s \in \mathbb{C}$) were introduced by Godement--Jacquet \cite{MR0342495}, and later interpreted more broadly by Jacquet--Piatetski-Shapiro--Shalika as the $\GL_n \times \GL_1$ case of Rankin--Selberg convolutions for $\GL_n \times \GL_m$; we refer to \S\ref{sec:standard-l-factors}, \S\ref{sec:standard-l-functions} and references for further details, and also the summaries of \cite{MR546609}, \cite[\S2]{MR1395406} and \cite{MR2508768}.  Conjectures of Langlands \cite{MR0302614} predict that ``all $L$-functions'' arise as shifted products of such standard $L$-functions.

As a matter of convention, we exclude archimedean factors from the definition, so that $L(\pi,s)$ is defined for $s \in \mathbb{C}$ with $\Re(s) > 1$ \cite[I, Thm 5.3]{MR618323} by an Euler product
\begin{equation*}
  L(\pi,s) = \prod_p L(\pi_p,s),
\end{equation*}
and in general by meromorphic continuation.  In the case $n = 1$, these are the Dirichlet $L$-functions, which include the Riemann zeta function as a special case.  We focus here on the case $n \geq 2$, in which $L(\pi,s)$ is known to be entire.

\subsection{The subconvexity problem in the $t$-aspect}\label{sec:subc-probl-t}
For each natural number $n$, let $\delta_n$ denote the supremum of all real numbers $\delta$ with the property that for each unitary cuspidal automorphic representation $\pi$ of $\GL_n$ over $\mathbb{Q}$, there exists $C_0 \geq 0$ so that for all $t \in \mathbb{R}$,
\begin{equation}\label{eq:lpi-tfrac12-+}
  |L(\pi,\tfrac{1}{2} + i  t)| \leq C_0 (1 + |t|)^{(1 - \delta) n/4}.
\end{equation}
The convexity bound (see Lemma \ref{lem:convexity-bound}) implies that $\delta_n \geq 0$, while the generalized Lindel\"{o}f hypothesis predicts that $\delta_n = 1$.  The $t$-aspect subconvexity problem for $\GL_n$ (over $\mathbb{Q}$) is to show that $\delta_n > 0$.  For an overview of the history on this problem, we refer to \cite{MR1403937}, \cite[\S2]{MR1826269}, \cite{MichelVenkateshICM}, \cite[\S4, \S5]{MR2331346} and \cite{MR3966770}.  We survey below some of the existing literature, which has primarily addressed the case $n \leq 3$.  Our motivating result is the following.
\begin{theorem}\label{thm:t-aspect}
  For each natural number $n$, we have $\delta_n > 0$.  In fact, for $n \geq 2$,
  \begin{equation}\label{eq:delta_n-geq-frac1}
    \delta_n \geq \delta_n^\sharp := \frac{2}{ 3 n^5 - 2 n^4 - n^2 } > \frac{2}{3 n^5}.
  \end{equation}
  That is to say, for each unitary cuspidal automorphic representation $\pi$ of $\GL_n$ over $\mathbb{Q}$ and each $\delta < \delta_n^\sharp$, there exists $C_0 \geq 0$ so that $|L(\pi,\tfrac{1}{2} + i t)| \leq C_0 (1 + |t|)^{(1 - \delta) n/4}$.
\end{theorem}

We note that when $n=2$ or $3$, our estimate \eqref{eq:delta_n-geq-frac1} is not competitive with the existing numerical records.  On the other hand, we have prioritized the clarity and uniformity of our arguments over their numerical optimality.

\subsubsection*{Summary of existing literature}
It will be useful to introduce the modified quantity $\delta_n^{(*)}$ (resp.\ $\delta_{n}^{(**)}$) defined similarly, but restricting the quantification to $\pi$ that are unramified at every finite place (resp.\ that correspond to Hecke--Maass cusp forms on $\SL_n(\mathbb{Z})$).  Clearly $\delta_n^{(**)} \geq \delta_n^{(*)} \geq \delta_n$.  Many papers establish lower bounds for $\delta_n^{(*)}$ or $\delta_n^{(**)}$ using methods that are expected to extend to $\delta_n$, although such extensions have not consistently been carried out.

It has been known that $\delta_1 > 0$ since the classical work of Weyl, Hardy--Littlewood and Landau \cite{MR1511862, MR1544665}, who in fact established the stronger bound $\delta_1 \geq 1/3$.  The current numerical record, $\delta_1^{(*)} \geq 65/171 = 0.380\dotsb$, follows from work of Bombieri--Iwaniec \cite{MR881101} and Bourgain \cite{MR3556291}.

It was shown first by Good \cite{MR696884} and Meurman \cite{MR1058223} that $\delta_2^{(*)} > 0$, and in fact that $\delta_2^{(*)} \geq 1/3$.  Recently Aggarwal \cite{MR4032289} showed that $\delta_2 \geq 1/3$ (see also Booker--Milinovich--Ng \cite{MR3872849}).



It was shown first by Munshi \cite{MR3369905}, following earlier work of X. Li \cite{MR2753605}, that $\delta_3^{(**)} > 0$, and in fact that $\delta_3^{(**)} \geq 1/24$.  This has been improved to $\delta_3^{(**)} \geq 1/10$ by Aggarwal \cite{MR4270876} and to $\delta_3^{(**)} \geq 1/6$ by Aggarwal--Leung--Munshi \cite{2022arXiv220606517A}.

Sporadic cases of $t$-aspect subconvex bounds in higher rank have also been given, e.g., for Rankin--Selberg convolutions on $\GL_2 \times \GL_2$ \cite{MR2726097, 2021arXiv210112106B} or $\GL_3 \times \GL_2$ \cite{2018arXiv181000539M}.

Soundararajan \cite[Example 3]{MR2680497} established the following nontrivial bound, weaker than the assertion that $\delta_n > 0$: for each $n$ and $\pi$ as above, and each $\eps > 0$, there exists $C_0$ so that
\begin{equation}\label{eq:lpi-tfrac12-+-1}
  |L(\pi, \tfrac{1}{2} + i t)| \leq C_0 \frac{(1 + |t|)^{n/4} }{(\log (3 + |t|))^{1 - \eps }}.
\end{equation}

\subsection{Subconvex bounds assuming uniform parameter growth}\label{sec:subc-bounds-assum}

We will establish Theorem \ref{thm:t-aspect} in a more general form.

The archimedean factor attached to $L(\pi,s)$ may be written
\begin{equation}\label{eqn:L-infinity}
  L_\infty(\pi,s) =
  \prod_{j=1}^{n}
  \Gamma_\mathbb{R}(s + t_j),
  \quad
  \Gamma_{\mathbb{R}}(s) := \pi^{-s/2} \Gamma(s/2)
\end{equation}
for some complex numbers $t_j$, which we call the \emph{archimedean $L$-function parameters} of $\pi$.  The functional equation for $L(\pi,s)$ involves a natural number, called the \emph{finite conductor} of $\pi$, which we denote by $C(\pi_{\fin})$.  The \emph{analytic conductor} of $\pi$ \cite{MR1826269} is then defined to be the product $C(\pi) := C(\pi_{\fin}) C(\pi_\infty)$, where $C(\pi_\infty)$, the \emph{archimedean conductor}, is the product $\prod_{j=1}^n (3 + |t_j|)$ or some mild variant thereof.

Write $\mathcal{A}(n)$ for the set of unitary cuspidal automorphic representations of $\GL_n$ over $\mathbb{Q}$.  The convexity bound (Lemma \ref{lem:convexity-bound}) asserts that for each natural number $n$ and $\eps > 0$, there exists $C_0 \geq 0$ so that for each $\pi \in \mathcal{A}(n)$,
\begin{equation*}
  |L(\pi,\tfrac{1}{2})| \leq C_0 C(\pi)^{1/4+\eps}.
\end{equation*}

We say that a family $\mathcal{F} \in \mathcal{A}(n)$ satisfies a \emph{subconvex bound} if there exists $\delta > 0$ and $C_0 \geq 0$ so that for all $\pi \in \mathcal{F}$,
\begin{equation*}
  |L(\pi, \tfrac{1}{2} )| \leq C_0 C(\pi)^{(1 - \delta)/4}.
\end{equation*}
The generalized Lindel\"{o}f hypothesis predicts that every family (or equivalently, the ``universal family'' $\mathcal{F} = \mathcal{A}(n)$) satisfies a subconvex bound with $\delta = 1$.  It is known by the classical work of Weyl \cite{MR1511862}, Hardy--Littlewood, Landau \cite{MR1544665}, Burgess \cite{Bur57} and Heath-Brown \cite{MR485727} that for $n = 1$, every family satisfies a subconvex bound.  The more recent work of Michel--Venkatesh \cite{michel-2009} established the same conclusion when $n=2$ (indeed, for the analogous problem over any fixed number field).

Soundararajan--Thorner \cite{MR3953433} recently established the general logarithmic improvement
\begin{equation}\label{eq:lpi-tfrac12-leq-2}
  |L(\pi,\tfrac{1}{2})| \leq C_0 \frac{C(\pi)^{1/4}}{(\log C(\pi))^{1/10^7 n^3}}.
\end{equation}
Their work made unconditional an earlier result of Soundararajan \cite[Thm 1]{MR2680497}, which featured the stronger logarithmic improvement $(\log C(\pi))^{1-\eps}$, as in \eqref{eq:lpi-tfrac12-+-1}, but under the assumption of (a weakened form of) the generalized Ramanujan conjecture.

We focus in this paper on the \emph{spectral aspect}, in which one regards the archimedean conductor $C(\pi_\infty)$ as varying substantially and the finite conductor $C(\pi_{\fin})$ as either fixed or varying mildly.  We say that $\mathcal{F} \subseteq \mathcal{A}(n)$ satisfies a \emph{spectral aspect subconvex bound, with polynomial dependence upon the level}, if there exists $\delta > 0$, $C_0 \geq 0$ and $A \geq 0$ so that for all $\pi \in \mathcal{F}$,
\begin{equation}\label{eq:lpi-tfrac12-leq-1}
  |L(\pi,\tfrac{1}{2})| \leq C_0 C(\pi_\infty)^{(1-\delta)/4} C(\pi_{\fin})^{A}.
\end{equation}

The following definition has played an important (although often implicit) role in the literature on the subject, and also in this paper.  We say that $\mathcal{F} \subseteq \mathcal{A}(n)$ has \emph{uniform parameter growth} if there exists $c > 0$ so that for each $\pi \in \mathcal{F}$, the archimedean $L$-function parameters $t_1,\dotsc,t_n$ of $\pi$ satisfy
\begin{equation}\label{eq:c-leq-fract_j1}
  \frac{\min \{|t_1|, \dotsc, |t_n|\} }{1 + \max \{|t_1|,\dotsc,|t_n|\}} \geq c.
\end{equation}
Note that the LHS of \eqref{eq:c-leq-fract_j1} is obviously bounded by $1$, so this condition forces all of the $t_j$ to be of roughly the same size.  For instance, a result of Jutila--Motohashi \cite{MR2233686} may be formulated as follows: for each $\delta < 1/3$, there exists $C_0 \geq 0$ so that for each $\pi \in \mathcal{A}(2)$ of ``full level'' (i.e., $C(\pi_{\fin}) = 1$), we have
\begin{equation*}
  |L(\pi,\tfrac{1}{2})| \leq C_0 (1 + \max \{|t_1|, |t_2|\})^{(1 - \delta)/2},
\end{equation*}
where $\{t_1, t_2\}$ denote the archimedean $L$-function parameters of $\pi$.  This result implies that a family $\mathcal{F} \subseteq \mathcal{A}(2)$ with uniform parameter growth and full level satisfies a spectral aspect subconvex bound (with polynomial dependence upon the level, tautologically).  However, specialized to the ``conductor dropping'' case in which one of the two parameters $t_1$ or $t_2$ is much smaller than the other (e.g., $t_2 = \O(1)$ while $t_1 \rightarrow \infty$), this estimate does not improve upon the convexity bound.  That case was not addressed until the work of Michel--Venkatesh noted above, following an important breakthrough of Duke--Friedlander--Iwaniec \cite{MR1923476} (see also \cite{MR2382859, Mi04, MR3912807}).

For $n \geq 3$, the first general examples of families $\mathcal{F} \subseteq \mathcal{A}(n)$ proven to satisfy subconvex bounds, beyond the ``$t$-aspect'' setting discussed in \S\ref{sec:subc-probl-t} or special families arising from functorial lifts, were given by Blomer--Buttcane \cite{MR4203038, MR4039487}.  To describe their result, let us say that $\mathcal{F} \subseteq \mathcal{A}(n)$ \emph{avoids the walls} if there exists $c > 0$ so that for each $\pi \in \mathcal{F}$,
\begin{equation*}
  \frac{  \max_{j \neq k}|t_j - t_k|}{1 + \max \{|t_1|,\dotsc,|t_n|\}} \geq c.
\end{equation*}
Blomer--Buttcane's results may be formulated as follows: if $\mathcal{F} \subseteq \mathcal{A}(3)$ has uniform parameter growth, avoids the walls, and each $\pi \in \mathcal{F}$ has trivial central character, then $\mathcal{F}$ satisfies a spectral aspect subconvex bound, with polynomial dependence upon the level.  A preprint of Kumar--Mallesham--Singh \cite{2020arXiv200607819K} establishes a similar conclusion, but with a modified form of the wall-avoidance condition.

We will show the following:
\begin{theorem}\label{thm:main-result}
  Let $n$ be a natural number.
  \begin{enumerate}[(a)]
  \item \label{itm:standard:main-a} Every family $\mathcal{F} \subseteq \mathcal{A}(n)$ with uniform parameter growth satisfies a spectral aspect subconvex bound, with polynomial dependence upon the level.
  \item \label{itm:standard:main-b} More quantitatively, suppose that $n \geq 2$.  Let $\delta < \delta_n^\sharp$, with $\delta_n^\sharp$ as in \eqref{eq:delta_n-geq-frac1}.  Let $C_0 > c_0 > 0$.  There exists $C_1 \geq 0$ with the following property.  Let $T \geq 1$.  Let $\pi$ be a unitary cuspidal automorphic representation of $\GL_n$ over $\mathbb{Q}$ whose archimedean $L$-function parameters $t_1,\dotsc,t_n \in \mathbb{C}$ satisfy
    \begin{equation}\label{eq:c_1-t-leq}
      c_0 T \leq |t_j| \leq C_0 T \quad \text{ for each }  j \in \{1,\dotsc,n\}.
    \end{equation}
    Then
    \begin{equation}\label{eq:lpi-tfrac12-leq}
      |L(\pi, \tfrac{1}{2})| \leq C_1 C(\pi_{\fin})^{1/2} T^{(1-\delta)n/4}.
    \end{equation}
  \end{enumerate}
\end{theorem}

For the sake of illustration, we record the straightforward deduction from the latter result to the former:
\begin{proof}[Proof of Theorem \ref{thm:main-result}, part \eqref{itm:standard:main-a}, assuming Theorem \ref{thm:main-result}, part \eqref{itm:standard:main-b}]
  By the classical results mentioned above, we may assume that $n \geq 2$.  Let $\mathcal{F} \subseteq \mathcal{A}(n)$ be a family with uniform parameter growth.  We may then find $C_0 > c_0 > 0$ so that for each $\pi \in \mathcal{F}$, there exists $T \geq 1$ (e.g., $T = 1 + \max \{|t_1|, \dotsc, |t_n|\}$) so that the parameters $\{t_1,\dotsc,t_n\}$ of $\pi$ satisfy \eqref{eq:c_1-t-leq}.  The ratio $C(\pi_\infty) / T^n$ is then bounded, from above and below, by positive quantities depending only upon $\mathcal{F}$, so the conclusion of part \eqref{itm:standard:main-b} implies that for each $\delta < \delta_n^\sharp$, there exists $C_1 \geq 0$ so that for all $\pi \in \mathcal{F}$, the bound \eqref{eq:lpi-tfrac12-leq-1} holds with the given value of $\delta$ and with $A := 1/2$.
\end{proof}

It is similarly straightforward to see that Theorem \ref{thm:main-result} generalizes Theorem \ref{thm:t-aspect}:
\begin{proof}[Proof of Theorem \ref{thm:t-aspect}, assuming Theorem \ref{thm:main-result}]
  Let $\pi \in \mathcal{A}(n)$.  Let $t_1, \dotsc, t_n \in \mathbb{C}$ denote its parameters.  We may find $T_0 \geq 1$ so that whenever $t \in \mathbb{R}$ satisfies $|t| > T_0$, we have
  \begin{equation*}
    \tfrac{1}{2} \lvert t \rvert \leq |t_j + t | \leq 2 \lvert t \rvert,
  \end{equation*}
  say.  Then the unitary twist $\pi \otimes |\det|^{i t}$ satisfies the hypotheses of Theorem \ref{thm:main-result} with $c_0 = \tfrac{1}{2}$, $C_0 = 2$ and $T = \lvert t \rvert$.  Since $L(\pi \otimes |\det|^{it}, \tfrac{1}{2}) = L(\pi, \tfrac{1}{2} + i t)$, the conclusion of Theorem \ref{thm:main-result} implies that we may find $C_1 \geq 0$ so that whenever $|t| > T_0$, we have
  \begin{equation*}
    |L(\pi, \tfrac{1}{2} + i t)| \leq C_1 ( 1 + |t|)^{(1 - \delta) n/4}.
  \end{equation*}
  The bound \eqref{eq:lpi-tfrac12-+} then holds with the leading scalar (denoted there ``$C_0$'') given by
  \begin{equation*}
    C_2 := \max \left(\{C_1\} \cup \left\{ \frac{|L(\pi,\tfrac{1}{2} + i t)|}{(1 + |t|)^{(1 - \delta) n/4}} : |t| \leq T_0 \right\}\right).
  \end{equation*}
\end{proof}

\begin{remark}
  While the primary novelty of these last results is in addressing the case $n \geq 4$, the conclusion appears to be new even when $n=3$.  For instance, the spectral aspect subconvexity problem for $\GL_3$ in the degenerate case where $(t_1,t_2,t_3) = (T, T, - 2 T)$ is addressed by Theorem \ref{thm:main-result}, but not by the earlier works \cite{MR4203038, MR4039487, 2020arXiv200607819K}.  We note also that the proof of Theorem \ref{thm:main-result} treats the ``spectral aspect'' and ``$t$-aspect'' uniformly, whereas earlier papers covering those cases \cite{MR3369905, MR4270876, MR4203038, MR4039487, 2020arXiv200607819K} employ fundamentally different methodology.
\end{remark}

\begin{remark}
  Since our focus in this paper is on the spectral aspect, we have made no effort to optimize the exponent of $C(\pi_{\fin})^{1/2}$.  Indeed, we would have been satisfied with any estimate involving $C(\pi_{\fin})^A$ for some explicit $A \geq 0$.  The exponent $A = 1/2$ turns out to be achievable with essentially no effort beyond what would be required in a treatment of the full level case, where $C(\pi_{\fin}) = 1$.  It would be interesting and natural to strengthen this exponent to arbitrary $A > 1/4$, as in the convexity bound.
\end{remark}

\begin{remark}
  We suspect that the proof of Theorem \ref{thm:main-result} can be made polynomially effective with respect to the parameter $c_0$.  That is to say, we expect that one could provably take $C_1 = C_2 (1/c_0)^{C_3}$, where $C_2, C_3 > 0$ depend upon $C_0$ but not upon $c_0$.  To achieve this, it should suffice to give a strengthening of \cite[Theorem 15.2]{2020arXiv201202187N} that ``depends polynomially'' upon the quantity ``$\mathfrak{S}$'' in its statement.  We expect that an ineffective strengthening of the required form follows formally from the Nullstellensatz, an effective strengthening by quantifying each step of the proof.  Since the estimate \eqref{eq:lpi-tfrac12-leq} (with $C_1$ as just specified) is weaker than the convexity bound unless $c_0$ can be taken quite close to $1$, we expect that such a strengthening would yield the following slightly more natural statement, which would in turn represent the limit of the technique developed in this paper (as far as the spectral aspect is concerned, and ignoring the problem of optimizing $\delta$):

  \emph{ For each natural number $n$, there exists $C_0 > 0$ and $\delta > 0$ such that for each unitary cuspidal automorphic representation $\pi$ of $\GL_n$ over $\mathbb{Q}$ whose archimedean $L$-function parameters $t_1,\dotsc,t_n \in \mathbb{C}$ satisfy $|t_j| \leq T$ for all $j$, we have}
  \begin{equation}\label{eq:cuhjn2tftk}
    |L(\pi, \tfrac{1}{2})| \leq C_0 C(\pi_{\fin})^{1/2} T^{(1-\delta)n/4}.
  \end{equation}

  (In more detail, assume the following claimed polynomial effectivization of Theorem \ref{thm:main-result}: there are fixed quantities $\delta_0 > 0$ and $A \geq 0$ so that, writing $c_0 := T^{-1} \min_j |t_j| \leq 1$, we have
  \begin{equation}\label{eq:cuhjn2moma}
    L(\pi, \tfrac{1}{2})  \ll(1 / c_0)^{A} C(\pi_{\fin})^{1/2} T^{(1-\delta_0)n/4}.
  \end{equation}
  Recall also the convexity bound:
  \begin{equation}\label{eq:cuhjn2s7g2}
    L(\pi, \tfrac{1}{2}) \ll C(\pi_\infty)^{1/4 + o(1)} C(\pi_{\fin})^{1/4 + o(1)}.
  \end{equation}
  Fix $\delta > 0$ and $\eta > 0$.  If $c_0 \geq T^{- \eta}$, then $1 / c_0 \leq T^{\eta}$, so \eqref{eq:cuhjn2tftk} follows from \eqref{eq:cuhjn2moma} provided that
  \begin{equation}\label{eq:cuhjn5n62z}
    A \eta - \frac{\delta_0 n}{4} \leq - \frac{\delta n}{4}.
  \end{equation}
  On the other hand, if if $c_0 \leq T^{- \eta}$, then $C(\pi_\infty) \leq T^{n - \eta}$, so \eqref{eq:cuhjn2tftk} follows from \eqref{eq:cuhjn2s7g2} provided that
  \begin{equation}\label{eq:cuhjn5n3ni}
    -\frac{\eta}{4} < - \frac{\delta n}{4}.
  \end{equation}
  To cover both cases, we need to choose $\delta$ and $\eta$ so that both \eqref{eq:cuhjn5n62z} and \eqref{eq:cuhjn5n3ni} hold.  Take $\delta := \eta / 2 n$ (say), so that \eqref{eq:cuhjn5n3ni} holds.  Then \eqref{eq:cuhjn5n62z} becomes
  \begin{equation*}
    \biggl(A + \frac{1}{8}\biggr) \eta \leq \delta_0 \frac{n}{4},
  \end{equation*}
  which holds for $\eta$ small enough in terms of $A$ and $\delta_0$.)
\end{remark}

\subsection{Discussion of the proof}\label{sec:high-level-overview}

\subsubsection{Overview}

The closest result in the literature to Theorem \ref{thm:main-result} is the main result of \cite{2020arXiv201202187N}.  We summarize that result informally:
\begin{theorem}\label{thm:unitary-review}
  Let $G = \U_{n+1} \geq H = \U_n$ be unitary groups over a number field $F$, with $H$ anisotropic.  Let $S$ be a large enough finite set of places.  Let $\mathfrak{q} \in S$ be an archimedean place.  Let $\mathcal{F}$ be a family of pairs $(\pi,\sigma)$ of cuspidal automorphic representations of $G$ and $H$, respectively, with the following properties.
  \begin{itemize}
  \item The local components of $\pi$ and $\sigma$ are tempered at all places inside $S$, unramified outside $S$, and have bounded ramification at all places other than $\mathfrak{q}$.
  \item The archimedean parameters at $\mathfrak{q}$ of the base change Rankin--Selberg $L$-function $L(\pi \times \sigma, \tfrac{1}{2})$ satisfy the analogue of the ``uniform parameter growth'' condition \eqref{eq:c-leq-fract_j1}.
  \end{itemize}
  Then there is an explicit $\delta > 0$ so that, for some $C_0 \geq 0$, we have
  \begin{equation*}
    |L(\pi \times \sigma, \tfrac{1}{2})| \leq C_0 C(\pi \times \sigma)^{1/4-\delta}
  \end{equation*}
  for all $(\pi,\sigma) \in \mathcal{F}$.
\end{theorem}
Theorem \ref{thm:main-result} may be understood as a specialization of a hypothetical generalization of Theorem \ref{thm:unitary-review}: if we
\begin{itemize}
\item generalize Theorem \ref{thm:unitary-review} by dropping
  \begin{itemize}
  \item the anisotropy condition on $H$,
  \item the cuspidality condition on the automorphic representations, and
  \item the temperedness condition on their local components,
  \end{itemize}
\item specialize that generalization to the ``least anisotropic'' example $(G,H) = (\GL_{n+1}, \GL_n)$ and $F = \mathbb{Q}$, and then
\item specialize further by taking for $\sigma$ the ``least cuspidal'' automorphic representation, consisting of Eisenstein series with trivial parameters (i.e., induced from the trivial character of the Borel subgroup), so that the Rankin--Selberg convolution is given by
  \begin{equation*}
    L(\pi \times \sigma, \tfrac{1}{2}) = L(\pi, \tfrac{1}{2})^n,
  \end{equation*}
\end{itemize}
then we arrive at something like Theorem \ref{thm:main-result} (for $\mathrm{GL}_{n + 1}$).  In this paper, we address that specialization by combining the strategy and results of \cite{2020arXiv201202187N} with some new ideas that overcome the failure of temperedness, cuspidality and (most seriously) anisotropy.

We recall that the anisotropy of $H$ corresponds to the compactness of the quotients $H(\mathbb{Q}) \backslash H(\mathbb{A})$, $H(\mathbb{Z}) \backslash H(\mathbb{R})$.  There has been little work concerning the problem of bounding automorphic forms on higher-rank non-compact quotients, uniformly with respect to both
\begin{itemize}
\item the argument of the automorphic form as it tends into the cusp, and
\item the underlying parameters of the automorphic form.
\end{itemize}
Some exceptions include the recent work of Brumley--Templier \cite{MR4150475} and Blomer--Harcos--Maga \cite{MR3905609, MR4093914}.

We note that the present paper does not ``replace'' the work \cite{2020arXiv201202187N} treating the compact case.  Indeed, the main technical results of that work (concerning the volume estimates and refined microlocal calculus) are applied here as a ``black box'' and combined with some new ingredients.  On the other hand, the relative character asymptotics of \cite[\S19]{nelson-venkatesh-1}, which served as a crucial input to \cite{2020arXiv201202187N}, do not enter here; the analogous role in the present treatment is played by direct analysis of local zeta integrals.

We survey the local issues addressed by this paper in \S\ref{sec:intertw-oper-distr}, \S\ref{sec:asympt-local-zeta} and \S\ref{sec:asympt-kirill-model}.  The main global problem that we face and overcome is to give nontrivial bounds for the ``local $L^2$-norms'' of certain Eisenstein series, uniformly as we move away from the bulk of the fundamental domain.  This is needed to compensate for the growth of the geometric side of a pretrace formula for $\GL_n$, coming from a long sum over unipotent elements of $\GL_n(\mathbb{Z})$.  It is a familiar feature that trace-like formulas on non-compact quotients are complicated by the contributions of unipotent elements.  The Arthur--Selberg trace formula gives a well-known example where the taming of the unipotent contribution requires an inclusion-exclusion scheme involving Eisenstein spectra coming from Levi subgroups.  The present situation more closely resembles the Jacquet--Rallis relative trace formula for $(G, H \times H)$ (see \cite{MR2767518, MR4195660, MR3947975}, \cite[\S3]{MR4426741}), but differs in that, by working with wave packets (denoted $\Psi$ below) rather than automorphic forms, our integrals converge absolutely.  We are thus faced not with the qualitative-analytic problem of making sense of divergent integrals, but rather the quantitative-analytic problem of obtaining satisfactory estimates for convergent ones.  It was nevertheless unclear to us whether an inclusion-exclusion scheme would be necessary, or whether direct estimation would be possible.  We eventually succeeded via the latter route, in the manner indicated below in \S\ref{sec:theta-funct-posit}--\S\ref{sec:duality-1}.

In the remainder of this introductory section, we give a high-level overview, intended for experts, of the new features of this paper relative to the compact quotient case treated in \cite{2020arXiv201202187N}.  In particular, we assume some familiarity with the notion of ``coadjoint microlocalization'' described in \cite[\S1]{nelson-venkatesh-1} and applied in \cite{2020arXiv201202187N}.  Some readers may wish to first study \cite[\S2]{2020arXiv201202187N} for an informal overview of the compact case, and \cite[\S1]{nelson-venkatesh-1} for an overview of how microlocal analysis and the orbit method may be applied to the analysis of automorphic forms.  On the other hand, we have attempted to keep the body of this paper self-contained and to give complete details.  We emphasize that the following discussion is informal.  In particular, notation such as $\approx$ or $\lessapprox$ is not defined precisely.

\subsubsection{Setup}\label{sec:crb1dssqx4}
The main ideas may be explained in the full level case (i.e., $C(\pi_{\fin}) = 1$), where all the action takes place on the quotients
\begin{equation*} [G] := \Gamma \backslash G := \GL_{n+1}(\mathbb{Z}) \backslash \GL_{n+1}(\mathbb{R}), \quad [H] := \Gamma_H \backslash H := \GL_{n}(\mathbb{Z}) \backslash \GL_{n}(\mathbb{R}),
\end{equation*}
with $H$ embedded in $G$ via the upper-left corner map $h \mapsto \left(
  \begin{smallmatrix}
    h&0\\
    0     &1 \\
  \end{smallmatrix}
\right)$.  We denote by $\mathfrak{g}, \mathfrak{h}$ the Lie algebras and \index{Lie algebra!imaginary dual $\mathfrak{g}^\wedge = i \mathfrak{g}^*$} $\mathfrak{g}^\wedge, \mathfrak{h}^\wedge$ their imaginary duals.  We write $Z \leq G$ and $Z_H \leq H$ for the centers.

We regard $\pi$ as a space of cusp forms on $[G]$, transforming under $Z$ via some unitary character.  Similarly, $\sigma$ is the space of Eisenstein series on $[H]$ obtained by inducing the trivial character of the Borel subgroup.  The method of Jacquet--Piatetski-Shapiro--Shalika leads to an integral representation
\begin{equation}\label{eq:int-_h-varphi}
  \int_{[H]} \varphi \Psi \approx L(\pi,\tfrac{1}{2})^{n} Z(\varphi,\Psi)
  \qquad
  (\varphi \in \pi, \, \Psi \in \sigma),
\end{equation}
where $Z(\varphi,\Psi)$ denotes the local integral obtained by integrating the archimedean Whittaker functions of $\varphi$ and $\Psi$.  To bound $L(\pi,\tfrac{1}{2})$ using \eqref{eq:int-_h-varphi} (following \cite{MR2726097, venkatesh-2005, michel-2009}), we must choose vectors $\varphi$ and $\Psi$ for which we may prove that the local integral $Z(\varphi, \Psi)$ is not too small, while the global integral $\int_{[H]} \varphi \Psi$ is not too large.

Strictly speaking, we take $\Psi$ to be a ``wave packet'' or ``pseudo Eisenstein series'', obtained by integrating Eisenstein series induced from unramified characters of the Borel, that decays rapidly on $[H]$.  This choice has the effect of replacing $L(\pi,\tfrac{1}{2})^n$ by an average of $\prod_{j=1}^n L(\pi,\tfrac{1}{2} + s_j)$ over small $s = (s_1,\dotsc,s_n) \in \mathbb{C}^n$, but we will often ignore this point in this overview and write simply $L(\pi,\tfrac{1}{2})^n$ in place of the more complicated expression that actually arises.

We note the recent work of Tsuzuki \cite{2021arXiv210307835T}, which applies an integral representation similar to \eqref{eq:int-_h-varphi} to the non-vanishing problem for central $L$-values.

\subsubsection{Analytic test vectors}\label{sec:analyt-test-vect}
We choose vectors following \cite{nelson-venkatesh-1} (although with completely different implementation details, for reasons discussed below):

We retain the setting of Theorem \ref{thm:main-result} and \S\ref{sec:crb1dssqx4}; in particular, $T \geq 1$ is a parameter that captures the size of the archimedean $L$-function parameters for $\pi$.  We should think of $T$ as tending off to $\infty$, as the required estimates follow from the convexity bound when $T$ is bounded.  Let $\mathcal{O}_\pi \subseteq \mathfrak{g}^\wedge$ and $\mathcal{O}_\sigma \subseteq \mathfrak{h}^\wedge$ denote the associated coadjoint orbits.  (More precisely, we take these to be the loci of the infinitesimal characters -- see \S\ref{sec:infin-char}.)  We may choose $\tau \in \mathfrak{g}^\wedge$ and $\theta \in \mathfrak{h}^\wedge$ so that
\begin{equation*}
  T \tau \in \mathcal{O}_\pi, \quad
  T \theta \in \mathcal{O}_{\sigma},
  \quad
  \tau|_{\mathfrak{h}} = \theta.
\end{equation*}
Indeed, we may be quite explicit: since $\sigma$ consists of Eisenstein series with trivial parameters, $\mathcal{O}_\sigma$ is the nilcone, so we may take for $\theta$ any regular nilpotent Jordan block (i.e., of signature $(n)$).  We then solve for $\tau$ by an exercise in rational canonical form (\S\ref{sec:parameter-tau}).  For example, when $n = 3$, we might take
\begin{equation*}
  \theta = \begin{pmatrix}
    0 & 0 & 0 \\
    1 & 0 & 0 \\
    0 & 1 & 0
  \end{pmatrix},
  \quad
  \tau = \begin{pmatrix}
    0 & 0 & 0 & \ast \\
    1 & 0 & 0 & \ast \\
    0 & 1 & 0 & \ast \\
    0 & 0 & 1 & \ast
  \end{pmatrix}.
\end{equation*}
The asterisks here are essentially the $L$-function parameters of $\pi$ divided by $T$, so $\tau$ is normalized so as to remain bounded as $T \rightarrow \infty$.

We take for the cusp form $\varphi$ a unit vector microlocalized at $T \tau$.

We take for the pseudo Eisenstein series $\Psi$ a unit vector microlocalized at $T \theta$.  More precisely, $\Psi$ is described by a function
\begin{equation*}
  f : U_H \backslash H \rightarrow \mathbb{C},
\end{equation*}
where $U_H$ denotes a maximal unipotent subgroup.  For example, if we pretend that $\Gamma_H \backslash H = \SL_2(\mathbb{Z}) \backslash \SL_2(\mathbb{R})$, then we may identify $U_H \backslash H$ with the punctured plane $\mathbb{R}^2 - \{0\}$ and take $\Psi(g) = \sum_{v \in \mathbb{Z}^2 - \{0\}} f(v g)$.  The space $U_H \backslash H$ admits a right action by $H$ and a left action by the toral part $A_H$ of the normalizer of $U_H$.  In practice, we take $U_H$ lower-triangular and identify $A_H$ with the diagonal subgroup, thus, e.g., for $n=3$,
\begin{equation*}
  U_H = \begin{pmatrix}
    1 & 0 & 0 \\
    \ast & 1 & 0 \\
    \ast & \ast & 1
  \end{pmatrix},
  \quad
  A_H = \begin{pmatrix}
    \ast& 0 & 0 \\
    0 & \ast & 0 \\
    0 & 0 & \ast
  \end{pmatrix}.
\end{equation*}
We take $f$ to be:
\begin{itemize}
\item $L^2$-normalized,
\item microlocalized under $H$ at $T \theta$,
\item uniformly smooth under $A_H$, and
\item concentrated near $T^{-\rho^\vee}$, where $\rho^\vee$ denotes the half-sum of positive coroots for $U_H$.
\end{itemize}
The matrix $T^{- \rho^\vee}$ plays a central role throughout the analysis, so we record for the sake of illustration that when $n =2,3,4$, it is the diagonal matrix with entries
\begin{equation*}
  (T^{1/2}, T^{-1/2}),
  \quad
  (T, 1, T^{-1}),
  \quad
  (T^{3/2}, T^{1/2}, T^{-1/2}, T^{-3/2}),
\end{equation*}
respectively.

In the $\SL_2$ example, the resulting function on $\mathbb{R}^2 - \{0\}$ has the rough shape
\begin{equation*}
  f(x,y) \approx
  \begin{cases}
    T^{-1/4}
    e^{i T y/x}
    & \text{ if }
      x \asymp T^{1/2}, \, y \ll 1, \\
    0 & \text{ otherwise,}
  \end{cases}
\end{equation*}
or equivalently,
\begin{equation*}
  f(r \cos \alpha, r \sin \alpha) \approx
  \begin{cases}
    T^{-1/4}  e^{i T \alpha}
    & \text{ if } r \asymp T^{1/2}, \,  \alpha \ll T^{-1/2}, \\
    0 & \text{ otherwise.}
  \end{cases}
\end{equation*}
We assume moreover that $f$ is chosen so that $\Psi$ is ``nondegenerate'' in certain senses (e.g., orthogonal to the residual spectrum); in the $\SL_2$ example, this amounts to requiring that the radial Mellin transform of $f$ vanish at certain arguments.

Ignoring wave packet issues, we obtain with this choice
\begin{equation}\label{eq:zvarphi-psi-approx}
  Z(\varphi, \Psi) \approx T^{-n^2/4} \quad
  \left( \, = T^{- \dim(H) / 4} \right),
\end{equation}
leading to an integral representation roughly of the shape
\begin{equation*}
  L(\pi,\tfrac{1}{2})^{n}  \approx T^{n^2/4} \int_{[H]} \varphi \cdot \Psi.
\end{equation*}
\begin{remark}
  Strictly speaking, since we are working with a wave packet, what we really consider is a family of zeta integrals $Z(\varphi, \Psi, s)$ indexed by small elements $s \in \mathbb{C}^n$.  The concentration condition for $f$ may be understood as ensuring that the variation of this family with respect to the real part of $s$ matches up with the convexity bound for the $L$-functions (see \S\ref{sec:reduct-peri-bounds} for details).  This matching may be understood as the analogue in our setup of working with a ``balanced'' approximate functional equation.
\end{remark}

\subsubsection{Relative trace formula}\label{sec:relat-trace-form}
We construct, as in \cite[\S1.5.2]{2020arXiv201202187N}, a ``convolution kernel'' $\omega \in C_c^\infty(G)$ with the following properties.
\begin{enumerate}[(i)]
\item The associated integral operator $\pi(\omega)$ satisfies
  \begin{equation}\label{eq:piom-varphi-appr}
    \pi(\omega) \varphi \approx \varphi.
  \end{equation}
\item $\omega$ concentrates on
  \begin{equation*}
    \left\{ g \in G : g = 1 + O(T^{-\eps}), \quad \Ad^*(g) \tau = \tau + O(T^{-1/2-\eps}) \right\}.
  \end{equation*}
\item $\omega$ has $L^\infty$-norm $\lessapprox T^{n(n+1)/2} $ ($ \, = T^{\frac{1}{2} \dim(\mathcal{O}_\pi)} = T^{\frac{1}{2} \dim(G / G_\tau)} = T^{\dim(N)}$, for a maximal unipotent subgroup $N \leq G$).
\end{enumerate}
Informally, the spectral support of $\omega$ is ``$\pi + \O(1)$,'' and $\pi(\omega)$ is essentially the orthogonal projection onto $\varphi$.  The function $\omega$ has the same shape as its adjoint $\omega^*$ and self-convolution $\omega \ast \omega^{\ast}$, so we denote these simply by $\omega$ in this informal overview.

As in \cite[\S1.5.3]{2020arXiv201202187N}, the main point is to show that
\begin{equation}\label{eq:t-n2-int}
  T^{-n/2} \int_{h_1, h_2 \in [H]} \sum_{\gamma \in \Gamma } \bar{\Psi}(h_1) \Psi(h_2) \omega(h_1^{-1} \gamma h_2) \, d h_1 \, d h_2
  = (\text{main term}) + \O(T^{-\delta}).
\end{equation}
Indeed, the pretrace formula for $[G]$ (or Cauchy--Schwarz) applied to $\omega$ shows that the LHS of \eqref{eq:t-n2-int} majorizes an average of $|L(\pi,\tfrac{1}{2})|^{2 n}$ over a family of size $T^{n(n+1)/2}$.  In view of our uniform parameter growth assumption, the family size is roughly the fourth root of the conductor, so the amplification method will deliver a subconvex bound if we can give a sufficiently robust proof of \eqref{eq:t-n2-int}.  The sum over $\gamma \in H Z$ contributes an apparent main term, so we reduce to showing that the $\gamma \notin H Z$ contribute $\O(T^{-\delta})$.  (There would be additional contributions to the main term had we not passed to wave packets, cf.\ \cite{MR2200997}.)

\subsubsection{Bilinear forms estimates}
We record an informal statement of the ``bilinear forms'' estimate given by \cite[Thm 15.1]{2020arXiv201202187N}.  The proof of that estimate relies ultimately upon the volume bound established in \cite[\S16--17]{2020arXiv201202187N}, a special case of which was established independently by Marshall \cite[Prop 6.2]{2023arXiv2309.16667}.  Fix a compact neighborhood $\mathcal{H}$ of the identity element in $H$.  Let $\Psi_1$ and $\Psi_2$ be measurable functions, supported on $\mathcal{H}$, with $\Psi_2$ satisfying the approximate equivariance condition
\begin{equation*}
  |\Psi_2(h z)| \approx |\Psi_2(h)| \text{ for } (h,z) \in \mathcal{H} \times Z_H  \text{ with } z = \O(1).
\end{equation*}
Let $g \in G$ be ``not too close to $H Z$.''  Then
\begin{equation}\label{eq:t-n2-int-1}
  T^{-n/2} \int_{h_1, h_2 \in \mathcal{H} }
  \left\lvert
    \Psi_1(h_1) \Psi_2(h_2)
    \omega(h_1^{-1} g h_2) \right\rvert \, d h_1 \, d h_2
  \lessapprox
  T^{-1/2 + \eps}
  \|\Psi_1\|_{L^2(\mathcal{H})}
  \|\Psi_2\|_{L^2(\mathcal{H})}.
\end{equation}
The factor $T^{-1/2}$ improves upon what one might regard as the ``trivial estimate'' obtained by omitting that factor.  We refer to \cite[\S1.5.3--1.5.7, \S2.9--2.11]{2020arXiv201202187N} for an extended discussion of such estimates.

In the compact quotient case, the estimate \eqref{eq:t-n2-int-1} readily yields \eqref{eq:t-n2-int}: take $\mathcal{H}$ large enough to contain a fundamental domain for $[H]$ and observe that
\begin{itemize}
\item the number of $\gamma$ contributing to \eqref{eq:t-n2-int} is essentially $\O(1)$, and
\item by the discreteness of $\Gamma$, any element of $\Gamma - H Z$ is ``not too close to $H Z$'' in the required sense.
\end{itemize}

In the non-compact case, the same argument gives an adequate estimate for the contribution from $x, y$ lying in any fixed compact subset of $[H]$.  (Recall that $\Psi$ is a wave packet, so we have good control over its $L^2$-norm.)  Unfortunately, as $x$ and/or $y$ ``escape to $\infty$,'' the sum over $\gamma$ in \eqref{eq:t-n2-int} features a large number of unipotent elements, so unlike in the compact case, we cannot conclude immediately.

\subsubsection{Siegel domains}\label{sec:crb1eg63zb}

Every element of $[H]$ may be represented in the form $h = t x$ with $(t,x) \in \mathcal{T} \times \mathcal{H}$, where $\mathcal{T}$ denotes the set of diagonal matrices $t$ with entries $t_1 \geq \dotsb \geq t_n > 0$ and, as above, $\mathcal{H}$ is a large enough fixed compact set.  The Haar measure $d h$ on $[H]$ is comparable to
\begin{equation*}
  \frac{d t}{\delta_H(t)} \, d x,
  \quad \delta_H(t) := \delta_{N_H}(t) = t_1^{n-1} t_2^{n-3} \dotsb t_n^{1-n},
  \quad
  d t = \frac{d t_1}{t_1} \dotsb \frac{d t_n}{t_n}.
\end{equation*}
We majorize the sum over $\gamma \notin H Z$ in \eqref{eq:t-n2-int} by
\begin{equation*}
  T^{-n/2}   \int_{t, u \in \mathcal{T}  }
  \sum_{\gamma \in \Gamma  - H Z}
  \int_{x,y \in \mathcal{H} }
  \left\lvert
    \Psi(t x)
    \Psi(u y)
    \omega(x^{-1} t^{-1}  \gamma u y)
  \right\rvert
  \, d x \, d y \, \frac{d t \, d u}{\delta_H(t) \delta_H(u)}.
\end{equation*}
We can apply the bilinear forms estimate \eqref{eq:t-n2-int-1} with $g = t^{-1} \gamma u$.  (Many of these elements will be closer to $H Z$ than in the compact case, where discreteness forces them all to be far away, but the dominant contribution still comes from those that are not too close.)  Our task then reduces to showing that
\begin{equation*}
  T^{-1/2}
  \int_{t, u \in \mathcal{T} }
  \nu(t,u)
  \|\Psi \|_{L^2(t \mathcal{H})}
  \|\Psi \|_{L^2(u \mathcal{H})}
  \frac{d t \, d u}{\delta_H(t) \delta_H(u)} \ll T^{-\delta},
\end{equation*}
where $\nu(t,u)$ denotes the number of $\gamma \in \Gamma - H Z$ for which $t^{-1} \gamma u = \O(1)$.  By elementary counting arguments, we see that $\nu(t,u)$ vanishes unless $t \approx u$ (in this informal discussion, take ``$A \approx B$'' to mean that $A B^{-1}$ lies in some essentially fixed compact set), in which case $\nu(t,u) \approx \delta_H(t) t^\dagger$, where
\begin{equation}\label{eq:t-dagger-:=}
  t^\dagger := \prod_{j=1}^n \max(t_j, t_j^{-1}).
\end{equation}
After Cauchy--Schwarz, the main problem is then to show that
\begin{equation}\label{eq:t-14-int}
  T^{-1/2}
  \int_{t  \in \mathcal{T} }
  \|\Psi \|_{L^2(t \mathcal{H})}^2
  t^\dagger
  \frac{d t}{\delta_H(t)} \ll T^{-\delta}.
\end{equation}
By comparison, the $L^2$-normalization of $\Psi$ gives
\begin{equation*}
  \int_{t  \in \mathcal{T} }
  \|\Psi \|_{L^2(t \mathcal{H})}^2
  \frac{d t}{\delta_H(t)} \ll 1,
\end{equation*}
so the ``enemy'' is the additional factor $t^\dagger$.

The wave packet $\Psi$ concentrates near $\SL_n(\mathbb{R})$ (since $\det(T^{-\rho^\vee}) = 1$), so one can focus on the contribution to \eqref{eq:t-14-int} from when $\det(t) \approx 1$.  Since $\Psi$ is of rapid decay, the integrals \eqref{eq:t-14-int} decay rapidly in a qualitative sense, so ``very large'' values of $t$ are no problem; the ``intermediate'' values of size larger than $1$ must be addressed.  We refer to \cite[\S1.2]{MR4093914} for some hint as to the precise meaning of ``intermediate.''  In any event, the total volume of the set of $t$ that we must consider is logarithmic in $T$, so the basic problem is to show that
\begin{equation}\label{eq:t-in-mathcalt}
  t \in \mathcal{T},
  \det(t) = 1
  \implies
  T^{-1/2}
  \|\Psi \|_{L^2(t \mathcal{H})}^2
  t^\dagger  \ll
  \delta_H(t)
  T^{-\delta}.
\end{equation}

\begin{remark}
  We caution that $\Psi$ is not, e.g., bounded by $T^{o(1)}$ -- one can show that the Whittaker function of $\Psi$ contains large ``spikes,'' which, for the same reasons as in the work of Brumley--Templier \cite{MR4150475}, implies that it assumes large values near the cusp.  While this feature suggests that we should be careful, the failure of $\Psi$ to satisfy strong pointwise bounds does not directly contradict the required estimate \eqref{eq:t-in-mathcalt}, which concerns only the local average behavior.
\end{remark}

\subsubsection{Theta functions and positivity}\label{sec:theta-funct-posit}
The estimate \eqref{eq:t-in-mathcalt} is straightforward when $H = \GL_2(\mathbb{R})$, which is the case relevant for the subconvexity problem on $\GL_3$: every approach that we have tried works without difficulty.  Already when $H = \GL_3(\mathbb{R})$ (relevant for the subconvexity problem on $\GL_4$, hitherto unaddressed), an approach based on direct estimation of the Fourier expansion (e.g., as in \cite{MR3905609}) seems hopeless: optimal bounds for all terms in that expansion add up to a worse than required bound.  Several other natural avenues of attack (e.g., using the Rankin--Selberg integrals that represent $L(\sigma \times \tilde{\sigma}, 1+ \eps)$) also seem inadequate.  We will establish \eqref{eq:t-in-mathcalt} instead using ``pseudo mirabolic Eisenstein series'' or ``Siegel theta functions,'' which provide a more flexible source of nonnegative functions, taking large values in the cusp, than do their ``spectral'' counterparts relevant for the integral representation of $L(\sigma \times \tilde{\sigma}, 1 + \eps)$.

Let $\phi \in C_c^\infty(\mathbb{R}^n)$ be a nonnegative function supported on a dyadic annulus of radii $\approx t_n$.  Form the series
\begin{equation*}
  E(g) := \sideset{}{^*}\sum_{v \in \mathbb{Z}^n} \phi(v g),
\end{equation*}
where $\sum^*$ denotes a sum over primitive vectors.  By considering the contribution to that series from the $n$th standard basis vector, we see that $E(g) \gg 1$ for $g \in t \mathcal{H}$.  It follows that
\begin{equation*}
  \|\Psi \|^2_{L^2(t \mathcal{H})} \ll \delta_H(t) \int_{[H]} |\Psi|^2 E
\end{equation*}
The integral $\int_{[H]} |\Psi|^2 E$ is a ``pseudo Rankin--Selberg integral'' (i.e., a spectral expansion of $\Psi$ and $E$ yields $\GL_n \times \GL_n$ Rankin--Selberg integrals with Eisenstein factors).  It unfolds nicely enough.  With some work, we can estimate the associated ``pseudo local Rankin--Selberg integrals.''  We arrive in this way at the estimate
\begin{equation}\label{eq:psi-_l2t-mathcalh2}
  \|\Psi \|_{L^2(t \mathcal{H})}^2 \ll \delta_H(t) t_n^n T^{\eps}.
\end{equation}
We could have chosen $\phi$ differently, i.e., concentrated on a different dyadic annulus, but the estimate \eqref{eq:psi-_l2t-mathcalh2} turns out to be optimal among such choices.

The estimate \eqref{eq:psi-_l2t-mathcalh2} does not imply \eqref{eq:t-in-mathcalt}.  A typical ``problem case'' is when
\begin{equation}\label{eq:t-=-diagt_1}
  t = \diag(t_1,\dotsc,t_n),
  \quad
  t_1 = r^{n-1}, \quad
  t_2 = \dotsb = t_n = r^{-1}
\end{equation}
where $r$ is a small positive power of $T$.  In that case, $t_n^n = r^{-n}$ and $t^\dagger \approx r^{2 (n-1)}$, so to achieve \eqref{eq:t-in-mathcalt}, we would need $T^{-1/2} r^{-n + 2 (n - 1)} \ll T^{-\delta}$, or equivalently, $r^{n - 2} \ll T^{1/2 - \delta}$.  The proof of \cite[Thm 3]{MR4093914} suggests that the rapid decay of $\Psi$ kicks in only when $r$ is a bit larger than $T^{1/n}$.  Consider then, for instance, the contribution from when $r = T^{1/n}$.  The desired estimate is that $T^{1 - 2/n} \ll T^{1/2 - \delta}$, which holds only for $n \leq 3$.

\subsubsection{Duality}\label{sec:duality-1}
We can do better by exploiting an important symmetry of the problem, coming from the outer automorphism of $\GL_n$ defined by the inverse transpose map $g \mapsto g^{-\transpose}$.  (We came to appreciate the power of this symmetry while perusing \cite[\S3.3]{MR4093914}.)  We can argue either in terms of a dual mirabolic Eisenstein series
\begin{equation*}
  \tilde{E}(g) = \sideset{}{^*}\sum_{v \in \mathbb{Z}^n} \phi(v g^{-\transpose}),
\end{equation*}
or a dual Eisenstein series
\begin{equation*}
  \tilde{\Psi}(g) := \Psi(g^{-\transpose}).
\end{equation*}
Either way, we obtain the following bound, complementary to \eqref{eq:psi-_l2t-mathcalh2}:
\begin{equation}\label{eq:psi-_l2t-mathcalh2-2}
  \|\Psi \|_{L^2(t \mathcal{H})}^2 \ll \delta_H(t) t_1^{-n} T^{\eps}.
\end{equation}
It is easy to see that this estimate is highly effective in ``problem cases'' such as \eqref{eq:t-=-diagt_1}.  In fact, the two estimates \eqref{eq:psi-_l2t-mathcalh2} and \eqref{eq:psi-_l2t-mathcalh2-2} combine to cover all cases, due to elementary inequalities such as the following (see Lemma \ref{lem:standard:any-dominant-element}):
\begin{equation}\label{eq:cujyyxote2}
  \det(1) = 1 \implies t^\dagger \leq \max(t_1^n, t_n^{-n}).
\end{equation}
We thus obtain
\begin{equation}\label{eq:int-_t-in-Psi-t-dagger}
  \int_{t  \in \mathcal{T} }
  \|\Psi \|_{L^2(t \mathcal{H})}^2
  t^\dagger
  \frac{d t}{\delta_H(t)} \ll T^{\eps},
\end{equation}
and the proof of \eqref{eq:t-14-int} is, to first approximation, complete.

\begin{remark}
  While the argument outlined above seems natural in retrospect, it was unclear to us for some time whether the estimates \eqref{eq:t-14-int} and \eqref{eq:t-in-mathcalt} were even true, or whether more careful arguments following \eqref{eq:t-n2-int} (e.g., Poisson summation applied to the unipotent contributions, exploiting cancellation coming from the sum over $\gamma$, or inclusion-exclusion involving the non-cuspidal spectrum of $G$) would be needed to achieve a satisfactory estimate.  We elaborate briefly on one source of our concern.  We have referred to the geometric expansion \eqref{eq:t-n2-int} as describing a spectral average of $\left\lvert L(\pi,\tfrac{1}{2}) \right\rvert^{2 n}$ over $\pi$ in some family.  Let us pause to consider in more detail what that spectral average actually looks like.  (It does not arise explicitly in this work -- we argue via Cauchy--Schwarz, in the guise of the pretrace inequality stated in \cite[Lemma 5.5]{2020arXiv201202187N}, rather than an actual pretrace formula.)  The contribution from cuspidal $\pi$ is as we have described, modulo the wave packet issues noted previously which have the effect of replacing the $2n$th moment with an average of shifts.  One could likely show that the contribution of the residual spectrum is negligible, thanks to the specific shape of either $\omega$ or $\Psi$.  It remains to understand the contribution of the Eisenstein spectrum, which involves non-cuspidal generic automorphic representations $\pi$.  Let $\varphi \in \pi$ be a unit vector in the essential image of $\pi(\omega)$ (heuristically, $\varphi$ is essentially unique, since $\pi(\omega)$ is essentially of rank one).  Then the contribution of $\pi$ to the spectral average is essentially the squared magnitude $\left\lvert \int_{[H]} \varphi \Psi \right\rvert^2$ of the integral analogous to that in \eqref{eq:int-_h-varphi}.  There are no convergence issues here, since $\Psi$ is a wave packet.  However, since $\pi$ is non-cuspidal, unfolding as in \eqref{eq:int-_h-varphi} produces additional ``degenerate'' terms, coming from the non-generic parts of the Fourier expansion of $\varphi$.  These would be ``subtracted off'' in a proper discussion of regularized $\GL_{n+1} \times \GL_{n}$ zeta integrals \cite{MR3334892}, and it was not clear to us whether they might be large enough to spoil the estimate \eqref{eq:t-n2-int}, thereby precluding the possibility of \eqref{eq:t-14-int}.  (There is precedent for the spoiling role of contributions from the Eisenstein spectrum in the subconvexity literature for $\GL_2$ and $\GL_3$, see for instance \cite{MR1923476, MR3694659}.)  It would be possible to excise such terms from the spectral average (essentially via the Kuznetsov formula, roughly as in the work of Blomer--Buttcane \cite{MR4203038, MR4039487}), but at the cost of producing a more complicated geometric expansion \eqref{eq:t-n2-int}.
\end{remark}

\subsubsection{Intertwining operators and distributions related to Jacquet integrals}\label{sec:intertw-oper-distr}
For the cuspidal variant of the problem considered here, the proof of \eqref{eq:int-_t-in-Psi-t-dagger} is quite simple: see \S\ref{sec:simpler-variants}.  For the Eisenstein case of primary interest, there are some further difficulties, discussed at length in Part \ref{part:local-l2-growth}.  One of these difficulties is to understand the asymptotics of the various ``partial Whittaker transforms'' of $\Psi$, including the constant term.  Recalling that $\Psi$ is attached to a function $f \in C^\infty(U_H \backslash H)$, an important intermediate step is to understand the behavior of $f$ under the Fourier transforms, or normalized intertwining operators, attached to all elements of the Weyl group of $H$.  We manage to avoid this step altogether by carefully constructing $f$ so as to be exactly invariant under all normalized intertwining operators.  To do this, we take for $f$ the image, under a carefully-constructed convolution operator on $A_H \times H$, of a distribution on $U_H \backslash H$ closely related to the theory of Jacquet integrals.  Since (by definition) normalized intertwining operators preserve Jacquet integrals, they likewise preserve $f$ provided that we take our convolution operator on $A_H \times H$ to be invariant under the action of the Weyl group on the first factor.

\subsubsection{Asymptotics of local zeta integrals}\label{sec:asympt-local-zeta}
We have not yet commented on the key local estimates \eqref{eq:piom-varphi-appr} and \eqref{eq:zvarphi-psi-approx}.  The analogues of these already appeared in the compact case \cite{2020arXiv201202187N}.  In that case, we argued as follows.  We started with the test function $\omega$.  We showed that $\pi(\omega)$ behaves like a self-adjoint rank one idempotent, and so approximately corresponds to a scaling class of vectors $\varphi$.  The vectors $\varphi$ themselves did not appear explicitly.  In the setting of unitary groups, the local zeta integrals $Z(\varphi, \Psi)$ are not necessarily defined.  In their place, we have certain matrix coefficients integrals, defined at least in the tempered case, which may be understood as avatars of $|Z(\varphi,\Psi)|^2$.  We estimated the matrix coefficient integrals using the ``reduction to small elements'' technique of \cite[\S19]{nelson-venkatesh-1} and then ultimately via appeal to the Kirillov character formula for $\pi$ and $\sigma$.  More specifically, we proved such estimates on average for $\varphi$ in an implicit family defined by $\omega$.

The approach of the present paper is completely different.  The primary motivation for the difference is that, since we are working with wave packets, we require not just estimates for individual zeta integrals, but rather a family of integrals $Z(\varphi, \Psi, s)$ indexed by small elements $s \in \mathbb{C}^n$.  Off the ``unitary axis'' $(i \mathbb{R})^n$, the connection between such zeta integrals and matrix coefficient integrals is more tenuous, and does not seem to us to lend itself readily to the problem of estimation.  A secondary motivation for the departure from \cite{2020arXiv201202187N} is that we do not wish to impose temperedness assumptions on $\pi$.  The Kirillov character formulas for non-tempered representations (see, e.g., \cite{MR4131117}) are more unwieldy than in the tempered case, and in any event, the estimates of \cite[\S19]{nelson-venkatesh-1} are available only in the tempered case.

We will instead study the zeta integrals $Z(\varphi,\Psi,s)$ directly.  In doing so, it is very convenient to assume that the Whittaker function of $\varphi$ has compactly-supported image in the Kirillov model.  This assumption trivializes the problem of analytic continuation of such integrals, and also reduces the range of arguments for which we need to estimate the Whittaker functions associated to $\Psi$.  For this reason, we will view the vector $\varphi \in \pi$ as the primary object, rather than the test function $\omega$.

\subsubsection{Asymptotics of the Kirillov model}\label{sec:asympt-kirill-model}
While the adjustment of perspective noted in \S\ref{sec:asympt-local-zeta} reduces the complexity of analyzing the zeta integrals, it introduces a new difficulty: proving that $\omega$ has the required approximation property \eqref{eq:piom-varphi-appr}.  (Recall that this property was essentially automatic in the approach pursued in the compact case, by virtue of the implicit definition of $\varphi$ in terms of $\omega$.)  The mirabolic subgroup $P$ of $G$ acts in a very simple way on the Kirillov model.  We construct $\varphi$ so that it is microlocalized with respect to $P$ at the restriction of $T \tau$ to $\Lie(P)$.  We need to know that $\varphi$ is actually microlocalized with respect to $G$ at $T \tau$.  At a heuristic level, this is clear, since $T \tau$ is the only element of the coadjoint orbit $\mathcal{O}_\pi$ lying above its restriction to $\Lie(P)$.  It is a bit subtle to implement this heuristic rigorously.  The scenario that we need to rule out is that the ``frequency'' of $\varphi$ with respect to the $G$-action is significantly larger than $T$.  To do this, we first apply a quantitative form of Jacquet's arguments \cite[\S3]{MR2733072} concerning the archimedean Kirillov model, which allows us to show that $\varphi$ has ``frequency $T^{\O(1)}$.''  This serves as the first step in a bootstrapping argument, involving the microlocal calculi of \cite{nelson-venkatesh-1} and \cite{2020arXiv201202187N}, by which we eventually conclude that $\varphi$ has frequency $\O(T)$ and is microlocalized at $T \tau$, as required.  We record a more detailed overview of this part of the paper in \S\ref{sec:overview-asymp-kirillov}.  We note in particular that the Kirillov character formula makes no appearance in this paper, in contrast to its repeated usage in the previous works \cite{nelson-venkatesh-1, 2020arXiv201202187N}.

\subsection{Outlook}
We briefly comment on the scope of the method.

It would be natural to give a common generalization of Theorems \ref{thm:main-result} and \ref{thm:unitary-review}, valid for all unitary (or orthogonal) GGP pairs and all families of pairs of automorphic representations $(\pi,\sigma)$ with uniform parameter growth.  As an imperfect analogy, one might think of such a generalization as the ``convex hull'' of those theorems.  For example, the case of Rankin--Selberg $L$-functions attached to cusp forms on $\GL_{n+1} \times \GL_n$ can be approached either by dropping the anisotropy and temperedness assumptions from Theorem \ref{thm:unitary-review}, or by adapting Theorem \ref{thm:main-result} to the case that $\sigma$ is cuspidal, which should in most respects be simpler than the Eisenstein case treated here.  However, such a generalization would not be altogether formal, as the local estimates developed here for the Eisenstein case use that they are tempered and belong locally to the principal series.  We hope that our work might nevertheless serve as a useful template for such a generalization.

Over number fields with complex places, there are a couple of lemmas in the literature concerning $\GL_n(\mathbb{R})$ that would need to be generalized to $\GL_n(\mathbb{C})$ (e.g., \cite[Lemma 3]{MR2733072} and \cite[Lemma 5.2]{JN19a}) but we do not foresee significant difficulties there.  The depth aspect is often similar to the archimedean aspect, but with fewer technical complications due to the availability of the exact projectors in the Hecke algebra; one might hope that it could be treated similarly.

The ``horizontal'' level aspect (e.g., taking $\pi$ of prime level), the quadratic twist case of which may have applications to quadratic forms in view of the appearance of self-dual standard $L$-functions in formulas for Whittaker coefficients on higher rank metaplectic groups \cite{MR3267120, MR3619910, MR3649366}, would require new ideas concerning the construction of analytic test vectors, and, we imagine, new estimates for exponential sums over finite fields, substituting for the role played by Proposition \ref{prop:standard2:let-tau-belong}, \cite[Thm 4.2 (iii) and Thm 15.1]{2020arXiv201202187N} and \cite[\S16--17]{2020arXiv201202187N} for the ``vertical'' (archimedean or depth) aspects.  The works of Bykovskii \cite{MR1433344}, Blomer--Harcos \cite{MR2431250, MR3259045} and Fouvry--Kowalski--Michel \cite{MR3334236} (resp.\ Sharma \cite{2019arXiv190609493S} and Lin--Michel--Sawin \cite{2019arXiv191209473L}) concerning the case of $\GL_2 \times \GL_1$ (resp.\ $\GL_3 \times \GL_2$) give some hint as to what might be required.  As a first step towards generalizing the results of those works to higher rank, it might be interesting to revisit them from the perspective of integral representations (cf.\ \cite{2020arXiv200201993N, MR4206609, 2020arXiv200406791S}).

The outstanding open problem is to address the case in which $\pi$ is not of uniform parameter growth, often known as the conductor dropping case.  We have already noted the works \cite{MR1923476, MR2382859, Mi04, michel-2009, MR3912807} addressing that case (or its related level aspect analogue) when $n=2$.  An extension of Theorem \ref{thm:main-result} for $n=3$ or $n=6$ to the conductor dropping case in which one or two of the archimedean $L$-function parameters remains bounded would have well-known applications to quantum chaos, as described, for instance, in \cite{MR1826269, watson-2008, MR2680499, PDN-AP-AS-que, sarnak-progress-que, 2020arXiv200108704N}.  Such an extension should not be expected to follow via analysis of the specific moments considered in this paper.  Indeed, to obtain such an extension, we would need to further shrink the archimedean families that we consider.  Doing so is not analytically feasible, even when $n=2$: the families are already of ``size $\O(1)$'' with respect to the parameters.

\subsection{Organization}
The paper is intended to be read linearly, but other approaches may be preferred.  For this reason, we briefly indicate the organization.  In Part \ref{part:reduction-proof}, we begin by formulating two auxiliary results:
\begin{itemize}
\item a local result, Theorem \ref{thm:main-local-results}, asserting the existence of elements $\varphi$, $f$ and $\omega$ satisfying the desiderata indicated in the sketch of \S\ref{sec:high-level-overview}, with $\varphi$ constructed in terms of its archimedean Whittaker function $W$, and
\item a global result, Theorem \ref{thm:growth-eisenstein-nonstandard}, asserting the required $L^2$-local growth bounds \eqref{eq:int-_t-in-Psi-t-dagger} for Eisenstein series.
\end{itemize}
The proof of the former is given in Part \ref{part:constr-test-vect}, the latter in Part \ref{part:local-l2-growth}.  In the remainder of Part \ref{part:reduction-proof}, we deduce our main result, Theorem \ref{thm:main-result}, from these auxiliary results, via the expected global ingredients (amplification method, counting arguments, $\dotsc$).  Part \ref{part:asympt-analys-kirill} supplies an independent ingredient (discussed in \S\ref{sec:asympt-kirill-model}) required in the construction of test vectors.  The dependency graph of the parts of the paper is thus as follows:
\begin{itemize}
\item Part \ref{part:reduction-proof}: Division of the argument
  \begin{itemize}
  \item Part \ref{part:constr-test-vect}: Construction of test vectors
    \begin{itemize}
    \item Part \ref{part:asympt-analys-kirill}: Asymptotics of the Kirillov model
    \end{itemize}
  \item Part \ref{part:local-l2-growth}: Growth bounds for Eisenstein series
  \end{itemize}
\end{itemize}

\subsection*{Acknowledgments}
We would like to thank V. Blomer, R. Beuzart-Plessis, P. Garrett, G. Harcos, Yueke Hu, P. Sarnak, E. Lapid, Ph. Michel, Y. Sakellaridis, A. Venkatesh and Liyang Yang for helpful discussions and feedback.  We are especially grateful to the anonymous referees for their careful reading and for several useful comments and corrections.  Much of this paper was written while the author was an Assistant Professor at ETH Zurich, during the Hausdorff Trimester Program \emph{Harmonic Analysis and Analytic Number Theory} hosted by the Hausdorf Research Institute for Mathematics.  This paper was completed at the Institute for Advanced Study, where the author was supported by the National Science Foundation under Grant No. DMS-1926686.  The paper was revised at Aarhus University, where the author is supported by a research grant (VIL54509) from VILLUM FONDEN.  We are grateful to each of these institutions for their excellent working conditions.

\section{Preliminaries}\label{sec:preliminaries}

\subsection{Asymptotic notation and terminology}\label{sec:asympt-notat-term}
\subsubsection{Overview}\label{sec:asymptotic-notation-overview}
This paper, like its predecessor \cite{2020arXiv201202187N}, involves the iterated application of estimates involving functions whose support, size and oscillation depend upon many parameters.  To manage this iteration sensibly, we again adopt the formalism of nonstandard analysis.  Specifically, we adopt the axiomatization due to E. Nelson \cite{MR469763} (``internal set theory'', or IST for short), but with the terminological substitution of ``fixed'' for ``standard.''  This substitution has the primary effect of making this paper read (we hope) like a typical analytic number theory paper, so as to require little adjustment of habit for the reader who is familiar with such papers.  A secondary effect is to have a convenient verb form: we write ``fix $x$'' synonymously with ``let $x$ be fixed.''  A tertiary effect is to avoid conflict with ``standard $L$-functions.''  We briefly recall that axiomatization, referring to \emph{loc.\ cit.}\ for details and superior exposition.

In ZFC, everything is a set.  ZFC consists of the binary predicate ``$\in$'' (interpreted: ``is an element of'') and some axioms for working with it.  IST adjoins to ZFC an additional unary predicate, which we label ``fixed,'' and three axioms, recalled below, for working with that.  A statement of IST is \emph{internal} if it is really a statement of ZFC, i.e., if it does not use the predicate ``fixed,'' even implicitly; otherwise, it is \emph{external}.  The appendix of \cite{MR469763} establishes the \emph{conservation theorem}: any internal statement provable in IST is also provable in ZFC.  More precisely, any proof in IST may be mechanically converted to one in ZFC.  In some applications of nonstandard analysis, such conversion is not practical, but in this paper, it is: the proof of our main results could be rewritten in ZFC in a reasonable amount of time by a human.


We define customary asymptotic notation in terms of ``fixed,'' as in \cite[\S3.1]{2020arXiv201202187N}.  \index{asymptotic notation!$\ll, \lll, \gg, \ggg, \O(\cdot), o(\cdot)$} If $|A| \leq C |B|$ for some fixed $C \geq 0$, then we write $A \ll B$ or $B \gg A$ or $A = \O(B)$; otherwise, we write $B \lll A$ or $A \ggg B$ or $B = o(A)$.  We write $A \asymp B$ as shorthand for $A \ll B \ll A$.  We introduce some further notation below, and will introduce other definitions involving ``fixed'' throughout the paper.  The reader who is familiar with the customary use of such notation can likely skip or skim the remainder of \S\ref{sec:asympt-notat-term} on a first reading, referring back in case anything seems unclear or imprecise.

Informally, ``fixed'' quantities behave like ``absolute constants.''  In particular, ordinary mathematical objects ($\mathbb{Z}$, $\mathbb{R}$, $\sqrt{2}$, the set of cuspidal automorphic representations of $\GL_{100}(\mathbb{Q})$, $\dotsc$) are fixed.  On the other hand, any infinite set (and any non-fixed finite set) contains non-fixed elements.  To make this heuristic more precise, one could introduce an asymptotic parameter $T$ (tending off to $\infty$ in some directed set), take every object in one's mathematical universe to depend upon $T$ by default, and think of ``fixed'' as meaning ``independent of $T$.''  This is essentially what happens in the ultrafilter approach to nonstandard analysis, which in turn forms the basis for the proof of the conservation theorem of IST.  The definitions and results of this paper could be reformulated in such ``parameter-dependent'' language, much like in \cite{nelson-venkatesh-1}, but doing so seemed to us to ultimately reduce the clarity and precision.

As detailed in \cite[\S2]{MR469763}, any external statement may be mechanically converted to an equivalent internal statement (subject to a minor technical caveat that applies in ordinary mathematics, and in particular, this paper).  We illustrate this algorithm below in the discussion of \S\ref{sec:main-result:-reform} and in the proof of Lemma \ref{lem:four-transf-begin}.

The principal results upon which we base the proof of Theorem \ref{thm:main-result} are formulated in the language of IST.  We record their equivalent internal formulations in Appendix \ref{sec:elem-form}.  Those reformulations are clunkier, and we do not work with them in the body of the text; they are provided only for the sake of illustration.

\subsubsection{Axioms of IST}\label{sec:axioms-ist}
We now recall the axioms, referring again to \cite{MR469763} for details.  In what follows, take ``quantity'' to be a synonym for ``set,'' given that we are working in ZFC.
\begin{itemize}
\item \emph{Transfer}.  Let $A(x)$ be an internal statement involving a quantity $x$ and possibly other fixed quantities, but no non-fixed quantities.  Then $A(x)$ is true for all $x$ if and only if it is true for all fixed $x$.
\item \emph{Idealization}.  Let $A(x,y)$ be an internal statement involving quantities $x$, $y$ and possibly other quantities.  Then the following statements are equivalent.
  \begin{itemize}
  \item For each fixed finite set $x'$, there exists $y$ so that, for each $x \in x'$, the statement $A(x,y)$ is true.
  \item There exists $y$ so that for each fixed $x$, the statement $A(x,y)$ is true.
  \end{itemize}
\item \emph{Standardization}.  Let $A(z)$ be a statement involving the quantity $z$.  Let $x$ be a fixed set.  Then there is a fixed subset $y \subseteq x$ so that for each fixed $z \in x$, we have $z \in y$ if and only if $A(z)$ holds.
\end{itemize}

The idealization axiom simplifies when
\begin{itemize}
\item the variable $x$ of the statement $A(x,y)$ is restricted to some directed set $(X, \geq)$, and
\item the statement $A(x,y)$ is \emph{filtered} in $x$: if $x_1 \geq x_2$, then $A(x_1,y) \implies A(x_2,y)$.
\end{itemize}
In that case, it says that one can commute quantifiers: ``for each fixed $x$, there exists $y$ so that $A(x,y)$'' is equivalent to ``there exists $y$ so that for each fixed $x$, we have $A(x,y)$.''

A primary consequence of the standardization axiom is the following (see \cite[Thm 1.3]{MR469763}).  Suppose given fixed sets $X$ and $Y$ and a statement $A(x,y)$, concerning elements $x \in X$ and $y \in Y$, such that for each fixed $x \in X$, there is a fixed $y \in Y$ for which $A(x,y)$ holds.  Then there is a fixed function $\tilde{y} : X \rightarrow Y$ such that for all fixed $x \in X$, we have $A(x, \tilde{y}(x))$

Each of the above axioms has a dual form obtained by swapping ``for all'' with ``there exists.''

\subsubsection{Main result: reformulation}\label{sec:main-result:-reform}
Recall that the main result to be proved is part \eqref{itm:standard:main-b} of Theorem \ref{thm:main-result}.  We explain here the equivalence of that statement with the following one.
\begin{theorem}[Theorem \ref{thm:main-result}, reformulated]\label{thm:main-result-reformulated}
  Let $T \geq 1$.  Let $n \geq 2$ be fixed.  Let $\pi$ be a unitary cuspidal automorphic representation of $\GL_n$ over $\mathbb{Q}$ whose archimedean $L$-function parameters satisfy $t_j \asymp T$.  Then for each fixed $\delta < \delta_n^\sharp$, we have
  \begin{equation*}
    L(\pi, \tfrac{1}{2}) \ll C(\pi_{\fin})^{1/2} T^{(1-\delta)n/4}.
  \end{equation*}
\end{theorem}
Recall the universal quantifier $\forall$ (``for all'') and the existential quantifier $\exists$ (``there exists'').  Let us write $\forall^*$ and $\exists^*$ for the modified quantifiers obtained by restricting quantification to fixed elements (thus, e.g., $\forall^* x$ means ``for all $x$, if $x$ is fixed, then $\dotsb$.'')  The condition ``$t_j \asymp T$ for all $j$'' expands to
\begin{equation*}
  \forall j \, \exists^* C_0, c_0 \, : j \in \{1, \dotsc, n\} \implies C_0 > c_0 > 0 \text{ and }  c_0 T \leq |t_j| \leq C_0 T.
\end{equation*}
The conclusion is filtered with respect to $(C_0,c_0)$, so by idealization, the above is equivalent to
\begin{equation*}
  \exists^* C_0, c_0 \,   \forall j \, : j \in \{1, \dotsc, n\} \implies C_0 > c_0 > 0 \text{ and }  c_0 T \leq |t_j| \leq C_0 T.
\end{equation*}
Theorem \ref{thm:main-result-reformulated} is thus equivalent to
\begin{equation}\label{eq:forall-t-reducing-nonstandard-reformulation}
  \forall T \,
  \forall^* n \,
  \forall \pi \,
  \forall^* \delta , C_0, c_0 \,
  \exists^* C_1 : A,
\end{equation}
where $A$ denotes the internal statement ``if $T \geq 1$, $n$ is a natural number, $\pi$ is a unitary cuspidal automorphic representation of $\GL_n$ over $\mathbb{Q}$, $\delta < \delta_n^\sharp$, $C_0 > c_0 > 0$, $c_0 T \leq |t_j| \leq C_0 T$ for all $j$, and $\delta < \delta_n^\sharp$, then $|L(\pi,\tfrac{1}{2})| \leq C_1 C(\pi_{\fin})^{1/2} T^{(1-\delta) n/4}$.''  By Boolean logic, \eqref{eq:forall-t-reducing-nonstandard-reformulation} is equivalent to
\begin{equation*}
  \forall^* n, \delta , C_0, c_0  \,
  \forall T , \pi  \,
  \exists^* C_1 \, : A.
\end{equation*}
The internal statement $A$ is filtered with respect to $C_1$, so by idealization, we can commute its quantification with that for $(T,\pi)$.  This yields the equivalent statement
\begin{equation}\label{eq:forall-n-delta}
  \forall^* n, \delta , C_0, c_0  \,
  \exists^* C_1 \,
  \forall T , \pi  \,
  : A.
\end{equation}
By repeated application of transfer to the internal statement ``$\forall T, \pi \, : A$,'' we obtain the further reformulation
\begin{equation*}
  \forall n, \delta , C_0, c_0  \,
  \exists C_1 \,
  \forall T , \pi  \,
  : A.
\end{equation*}
The latter unwinds to the statement of part \eqref{itm:standard:main-b} of Theorem \ref{thm:main-result}.

In this way, we obtain the required equivalence with Theorem \ref{thm:main-result-reformulated}.  By the conservation theorem, our task reduces to verifying the latter.

\subsubsection{Classes}\label{sec:classes}
As in \cite[\S3.1.4]{2020arXiv201202187N}, by a \emph{class}, we mean a pair $(X,P)$, where $X$ is a set and $P$ is a boolean formula of IST that accepts an element $x \in X$ as its argument.  We note that if $P$ is internal, then an axiom of ZFC tells us that there exists a subset of $X$, denoted
\begin{equation}\label{eq:cumotvxqcc}
  \{x \in X : P(x)\},
\end{equation}
that consists of precisely those elements for which $P$ is true.  On the other hand, ZFC says nothing about external formulas.  For example, there is no subset of $\mathbb{R}$ consisting of precisely those elements $x$ for which $x = \O(1)$.  On the other hand, it makes perfect sense to form the pair $(\mathbb{R},P)$, where $P(x)$ denotes the predicate ``$x = \O(1)$.''  We refer to such a pair $(X,P)$ as ``the class of all $x \in X$ satisfying $P(x)$.''  For instance, we may speak of the class of all real numbers of the form $\O(1)$.  In the language of some axiomatizations (see, e.g., \cite[\S3]{MR469763}, \cite{MR723332} or \cite[\S3]{MR546178}), ``classes'' correspond to ``external subsets of internal sets,'' but since we will be working with classes in a purely syntactic way, it is not necessary to adopt such additional formalism.

Classes do not appear explicitly in Part \ref{part:reduction-proof}.  They play a significant organizing role in Parts \ref{part:constr-test-vect} and \ref{part:asympt-analys-kirill}, and a minor one in Part \ref{part:local-l2-growth}.

We may regard sets as classes (take for $P$ the predicate ``true'') and work with classes in much the same way we work with sets, except that we will never apply anything like the ``power set'' axiom to classes, e.g., by speaking of ``the class of all subclasses of $(\dotsb)$.''  An element $x$ belongs to the class $(X,P)$ if $x \in X$ and if $P(x)$ is true.  One class $(X,P)$ is contained in another class $(Y,Q)$ (and called, then, a ``subclass'') if for every $x \in X$ satisfying $P(x)$, we have $x \in Y$ and $Q(x)$.  For example, for $T \geq 1$, the class of all real numbers of the form $\O(T)$ is contained in the class of all real numbers of the form $\O(T^2)$, say.  This basic idea is routinely applied in analytic number theory papers, but since we will eventually be applying it in somewhat complicated ways, we have attempted to be precise.

Two classes are equal if they contain the same elements, or equivalently, if each is contained in the other.  We may define intersections or unions of classes, class maps $f : (X,P) \rightarrow (Y,Q)$ between classes (i.e., subclasses of the product class satisfying the map axioms), and so on.

We record a simple example of the sorts of classes we will work with.  This example will not play a direct role in the paper, and is perhaps better treated without such formalism, but illustrates the basic idea.
\begin{example}\label{exa:positive-parameter-t-scaling}
  For a positive parameter $T$, let $\mathcal{A}_T$ denote the subclass of the Schwartz space $\mathcal{S}(\mathbb{R})$ consisting of $f$ that satisfy, for all fixed $m,n \in \mathbb{Z}_{\geq 0}$, the $n$th derivative estimate
  \begin{equation*}
    f^{(n)}(x) \ll \frac{1}{T^{1+n}} \left(1 + \frac{|x|}{T}\right)^{-m}
  \end{equation*}
  for all $x \in \mathbb{R}$.

  Let $\mathcal{A}_T^\wedge$ denote the subclass of $\mathcal{S}(\mathbb{R})$ consisting of $f$ satisfying, for all fixed $m,n \in \mathbb{Z}_{\geq 0}$,
  \begin{equation*}
    f^{(n)}(x) \ll T^{n} \left( 1 + T |x| \right)^{-m}.
  \end{equation*}
  Informally, $\mathcal{A}_T$ consists of bump functions of mean $\ll 1$ that concentrate on $\O(T)$, while $\mathcal{A}_T^\wedge$ consists of bump functions of magnitude $\ll 1$ that concentrate on $\O(1/T)$.
\end{example}
\begin{lemma}\label{lem:four-transf-begin}
  The Fourier transform
  \begin{equation*}
    f \mapsto f^\wedge := [\xi \mapsto \int_{x \in \mathbb{R} } f(x) e^{i \xi x} \, d x]
  \end{equation*}
  induces a bijective class map $\mathcal{A}_T \rightarrow \mathcal{A}_T^\wedge$.
\end{lemma}
\begin{proof}
  This follows by repeated partial integration applied to the definition of the Fourier integral.  Alternatively, we may reduce to the case $T = 1$ by rescaling, and then apply transfer to relate the conclusion to the well-known fact that the Fourier transform defines a topological automorphism of the Schwartz space.  The ``two proofs'' are ultimately equivalent.  We indicate in more detail the latter argument, primarily as an excuse to illustrate further the algorithm for translating external statements into internal ones.  (We emphasize that this is a toy example.  The nonstandard language become more useful than in this example when the definitions of the classes involve, e.g., support conditions depending upon parameters.)

  In the case $T = 1$, we need to verify the implication
  \begin{align*}
    &\forall f \in \mathcal{S}(\mathbb{R}) :
      (  \forall^* k,  \ell \,
      \exists^* C_0 \,
      \forall x :
      (1 + |x|)^\ell |f^{(k)}(x)| \leq C_0
      )
    \\
    &\implies
      (  \forall^* m,  n \,
      \exists^* C_1 \,
      \forall y :
      (1 + |y|)^m |(f^\wedge)^{(n)}(y)| \leq C_1
      )
  \end{align*}
  as well as its converse.  By Fourier inversion and the symmetry of the problem, we can just focus on the forward implication.

  We may rewrite the implication to be proved more succinctly as
  \begin{equation*}
    \forall f \,
    :  (
    \forall^* \kappa \, \exists^* C_0 \, :  \kappa(f) \leq C_0
    ) \implies
    (
    \forall^* \mu \, \exists^* C_1 \,  : \mu(f^\wedge) \leq C_1
    ),
  \end{equation*}
  where $f \in \mathcal{S}(\mathbb{R})$, $\kappa$ and $\mu$ are (continuous) seminorms on $\mathcal{S}(\mathbb{R})$, and $C_0, C_1$ are nonnegative reals.  (Indeed,
  \begin{equation*}
    f \mapsto \max_{k \leq k'} \sup_{x \in \mathbb{R}} (1 + |x|)^{\ell} |f^{(k)}(x)|
  \end{equation*}
  defines a seminorm for each $(k,\ell) \in \mathbb{Z}_{\geq 0}^2$, and every seminorm is bounded by a constant multiple of such a seminorm.)  By standardization, we can permute quantification over $\kappa$ and $C_0$ at the cost of replacing $C_0$ by a fixed function that depends upon $\kappa$:
  \begin{equation*}
    \forall f \,
    :  (
    \exists^* \tilde{C}_0 \, \forall^* \kappa \, : \kappa(f) \leq \tilde{C}_0(\kappa)
    ) \implies
    (
    \forall^* \mu \, \exists^* C_1 \, :  \mu(f^\wedge) \leq C_1
    ).
  \end{equation*}
  By Boolean logic, the above is equivalent to
  \begin{equation*}
    \forall f \, \forall^* \tilde{C}_0, \mu \, \exists^* \kappa, C_1 \,
    : \kappa(f) \leq \tilde{C}_0(\kappa)
    \implies \mu(f^\wedge) \leq C_1.
  \end{equation*}
  The universal quantification over $f$ and $(\tilde{C}_0, \mu)$ may be commuted.  The innermost statement is then of the form $\forall f \exists^* \kappa, C_1 : (\dotsb)$, where $(\dotsb)$ is an internal statement that filters with respect to $(\kappa,C_1)$.  By idealization, we may commute quantification over $f$ and $(\kappa,C_1)$.  We arrive at the equivalent statement
  \begin{equation*}
    \forall^* \tilde{C}_0, \mu \, \exists^* \kappa, C_1 \,
    \forall f \,
    : \kappa(f) \leq \tilde{C}_0(\kappa)
    \implies \mu(f^\wedge) \leq C_1.
  \end{equation*}
  By transfer, this is equivalent to the internal statement obtained by replaced $\forall^*$ with $\forall$ and $\exists^*$ with $\exists$:
  \begin{equation*}
    \forall^* \tilde{C}_0, \mu \, \exists^* \kappa, C_1 \,
    \forall f \,
    : \kappa(f) \leq \tilde{C}_0(\kappa)
    \implies \mu(f^\wedge) \leq C_1.
  \end{equation*}
  The function $\tilde{C}_0$ describes a general bounded subset $\Omega$ of $\mathcal{S}(\mathbb{R})$, so we are left with the assertion that, for each bounded subset $\Omega$ of $\mathcal{S}(\mathbb{R})$ and seminorm $\mu$ on $\mathcal{S}(\mathbb{R})$, there exists $C_1 \geq 0$ so that $\mu(f^\wedge) \leq C_1$ for all $f \in \Omega$, which is the content of the continuity of $f \mapsto f^\wedge$.
\end{proof}

\subsubsection{Overspill}
For our purposes, the overspill principle is a theorem of IST \cite[Thm 1.1]{MR469763} which asserts that any subset of the natural numbers that contains every fixed natural number must contain some natural number $n$ that is not fixed, hence for which $n \ggg 1$.  One verifies in the same way that any subset of the positive reals that contains every fixed positive real must contain some positive real of the form $o(1)$.  These conclusions would fail if we relaxed ``subset'' to ``subclass'' in the sense described above (consider, e.g., the class of natural numbers of the form $\O(1)$), so it is important to pay attention to what is actually a set rather than merely a class.  We will use this principle primarily for ``$\eps$ clean-up.''  A sample application is given by the lemma below.

The notation $A \ll T^{o(1)}$ (together with variants, such as $A \ll B T^{o(1)}$) will occur frequently.  The direct meaning is that there is a fixed $C \geq 0$ and an $\eps \lll 1$ such that $|A| \leq C T^{\eps}$.  We will often make use of the following equivalent characterizations:
\begin{lemma}\label{lem:overspill-A-vs-T-eps}
  Let $T$ be a positive quantity with $T \ggg 1$.  Let $A$ be a scalar.  The following conditions are equivalent.
  \begin{enumerate}[(i)]
  \item $|A| \leq T^{o(1)}$.
  \item $A \ll T^{o(1)}$.
  \item For each fixed $\eps > 0$, we have $A \ll T^{\eps}$.
  \end{enumerate}
  Under the weaker assumption that $T \geq 1$, conditions (ii) and (iii) are equivalent.
\end{lemma}
\begin{proof}
  It is clear that (i) implies (ii) implies (iii).

  Assume (iii).  We claim first that $|A| \leq T^\eps$ for each fixed $\eps > 0$.  Indeed, our hypothesis (applied to $\eps/2$) gives $A \ll T^{\eps/2}$; since $T \ggg 1$, we have $T^{\eps/2} \ggg 1$, and so $|A| \leq T^{\eps}$, as required.  The claim implies that the set $\{\eps > 0 : |A| \leq T^{\eps}\}$ contains every fixed $\eps$.  By overspill, it must contain some $\eps = o(1)$.  Thus (i) holds.

  Assume now only that $T \geq 1$.  If $T \ggg 1$, then we see from the above that (ii) and (iii) are equivalent.  Otherwise $T \ll 1$, and it is easy to see that both (ii) and (iii) are equivalent to ``$A \ll 1$.''
\end{proof}

\subsubsection{Further notation and conventions}\label{sec:furth-notat-conv}
We conclude with some miscellaneity.

Let $F$ be a fixed local field, with normalized absolute value $|.| : F \rightarrow \mathbb{R}_{\geq 0}$.  Let $V$ be a fixed vector space (finite-dimensional) over $F$.  By a \emph{norm} on $V$, we mean a continuous function $|.| : V \rightarrow \mathbb{R}_{>0}$ such that $|v| = 0$ only if $v = 0$, and $|t v| = |t|\, |v|$ for all $(t,v) \in F \times V$.  It is well-known that for any two fixed norms $|.|_1$ and $|.|_2$ on $V$, we have $|x|_1 \asymp |x|_2$ for all $x \in V$.  Given $x,y \in F$ and a scalar $t$ (either in $F$ or complex), we write $x = y + \O(t)$ to denote that $|x-y| \ll |t|$ for some (equivalently, any) fixed norm $|.|$ on $V$.  We analogously define the notations $x = y + o(t)$, $x \asymp y$, $x \ll t$ and $x \lll t$.

We often write $T$ and $\h$ for positive parameters with $T \geq 1$ and $\h \leq 1$; typically, $T = 1/\h$ and $T \ggg 1$, or equivalently, $\h \lll 1$.  We use notation such as $A \ll B T^{-\infty}$ or $A = \O(B T^{-\infty})$ as shorthand for the assertion that for each fixed $m \geq 0$, we have $A \ll B T^{-m}$; more explicitly, for each fixed $m \geq 0$ there is a fixed $C \geq 0$ so that $|A| \leq C |B| T^{-m}$.  We analogously define the notation $A \ll B \h^\infty$.  Similarly, given a subclass $\mathcal{C}$ of a vector space, we write $T^{-\infty} \mathcal{C}$ for the intersection over all fixed $m \geq 0$ of the classes $T^{-m} \mathcal{C}$.  This notation will be used only in situations where that intersection is descending.

\subsection{General notation}
For a statement $S$, we define $1_S$ to be $1$ if $S$ is true and $0$ otherwise.  For a set $X$, we write $1_X$ for its characteristic function, thus $1_X(x)$ is $1$ or $0$ according as $x \in X$ or $x \notin X$; equivalently, $1_X(x) := 1_{x \in X}$.

We interchangeably write $\# X$ and $|X|$ for the cardinality of a set $X$.

\subsection{General preliminaries}\label{sec:gener-prel}

\subsubsection{Fields}\label{sec:crb1dyvdtv}
Let $F$ be a local or global number field (i.e., a local or global field of characteristic zero, the only case considered here).

In the global case, we write $\mathbb{A}$ for the adele ring of $F$ and $\mathbb{A}_f$ for the subring of finite adeles.  We use the letter $\mathfrak{p}$ to denote places of $F$, possibly archimedean, and denote by $F_\mathfrak{p}$ the associated completion.  We set $F_\infty := F \otimes_{\mathbb{Q}} \mathbb{R} = \prod_{\mathfrak{p} | \infty} F_\mathfrak{p}$, so $\mathbb{A} = F_\infty \times \mathbb{A}_f$.  When working over $F = \mathbb{Q}$, we write $p$ for a finite prime.

In the local case, we write $|.| : F \rightarrow \mathbb{R}_{\geq 0}$ for the normalized absolute value; in the global case, $|.| : \mathbb{A}^\times / F^\times \rightarrow \mathbb{R}^\times_+$ denotes the idelic absolute value.

In the local non-archimedean case, we denote by $\mathfrak{o}$ the ring of integers, by $\mathfrak{p}$ the maximal ideal, by $\varpi \in \mathfrak{p}$ a generator, and $q := \# \mathfrak{o}/\mathfrak{p}$.

\subsubsection{Reductive groups}\label{sec:cumos34zyo}
Let $G$ be a split connected reductive group over $F$.  The examples relevant for this paper are when $G$ is a finite product of general linear groups (e.g., arising as Levi factors of a given general linear group), but it will be convenient at times to work more generally -- for instance, some points can be explained most clearly by passing temporarily to the simpler group $G = \SL_2$.

In the local case, we identify $G$ with its set of $F$-points, and abbreviate $G := G(F)$.  In the global case, we denote by \index{adelic quotients!$[G]$} $[G] := G(F) \backslash G(\mathbb{A})$ the associated adelic quotient.

In the local case, we write \index{Lie algebra!Lie algebra $\Lie(G)$} $\Lie(G)$ for the Lie algebra.  In the local archimedean case, we write $\mathfrak{U}(G)$ for the universal enveloping algebra \index{Lie algebra!enveloping algebra $\mathfrak{U}$, center $\mathfrak{Z}$} of the complexification of $\Lie(G)$, and $\mathfrak{Z}(G)$ for the center of $\mathfrak{U}(G)$.

We occasionally take for $G$ a split connected reductive group over $\mathbb{Z}$ or over the ring of integers $\mathfrak{o}$ in a non-archimedean local field $F$; in such cases, by base restriction, we may understand $G$ as being defined also over $\mathbb{Q}$ or $F$, respectively.

We always assume that $G$ comes equipped with the following additional data:
\begin{itemize}
\item A closed linear embedding $G \hookrightarrow \SL_r$ for some $r$.
\item A split maximal torus \index{groups!maximal torus $A$} $A$.  We assume given a basis for its group of rational characters, hence an identification of $A$ with a product of copies of $\GL_1(F)$.
\item A Borel subgroup $B$ \index{groups!Borel subgroup $B = A N$} of $G$, containing $A$, with unipotent radical denoted $N$, so that $B = A N$.
\item Another Borel subgroup $Q = A U$ \index{groups!Borel subgroup $Q = A U$} containing $A$, which may or may not coincide with $B$.  (We will use $Q$ to define induced representations, and $N$ to define Whittaker functions -- it will be computationally convenient to keep the choices of $Q$ and $B$ decoupled.)
\item In the local (resp.\ global) case, a maximal compact subgroup \index{groups!maximal compact subgroup $K$} $K$ of $G$ (resp.\ $K = \prod_\mathfrak{p} K_\mathfrak{p}$ of $G(\mathbb{A})$) in good position relative to $B$ (resp.\ $B(\mathbb{A})$), so that $A \cap K$ (resp.\ $A(\mathbb{A}) \cap K$) is a maximal compact subgroup of $A$ (resp.\ $A(\mathbb{A})$) and the Iwasawa decomposition $G = B K$ (resp.\ $G(\mathbb{A}) = B(\mathbb{A}) K$) holds (see \cite[I.1.4]{MR1361168}, \cite[Thm A.1.1]{GetzHahn}).
\item In the local archimedean case, \index{Lie algebra!Lie algebra basis $\mathcal{B}(G)$} bases $\mathcal{B}(A)$ and $\mathcal{B}(G)$ for $\Lie(A)$ and $\Lie(G)$, respectively.
\end{itemize}
When $G$ is fixed, we may and shall assume that all of the above data is fixed, too.  We denote by $Z_G$ the center of $G$ and by $W_G$ the Weyl group for $(G,A)$.  When $G$ is clear from context, we abbreviate $Z := Z_G$ and $W := W_G$.

When working over $\mathbb{R}$, $\mathbb{Q}$ or $\mathbb{Z}$, we write $A(\mathbb{R})^0$ \index{groups!identity component $A(\mathbb{R})^0$} for the topologically connected component of the identity element of $A(\mathbb{R})$.

\subsubsection{General linear groups}\label{sec:gener-line-groups-1}
By a \emph{general linear group} $G$, we mean $G = \GL_n$ for some $n \in \mathbb{Z}_{\geq 0}$.  We take the closed embedding $G \hookrightarrow \SL_{2 n}$ given by $g \mapsto \diag(g, g^{-\transpose})$.  We assume that $A$ is the diagonal subgroup and $B$ is the upper-triangular Borel subgroup.  We will take $U$ to be either $N$ or its opposite, according to convenience.  It will often be convenient to use the letter ``$n$'' for other purposes (e.g., to denote elements of the group $N$) or to be flexible in the indexing of our general linear groups (e.g., writing $G = \GL_{n+1}$ or $G = \GL_{n-1}$, according to convenience).  To that end, we will make use of the formulas
\begin{equation*}
  \rank(G) = n, \quad \dim(G) = n^2, \quad
  \dim(N) = n(n-1)/2,
\end{equation*}
and so on (with ranks and dimensions computed over $F$, unless otherwise indicated).

By a \emph{general linear pair} $(G,H)$ over $F$, we mean a pair of $F$-groups of the form $(\GL_{n+1},\GL_n)$ for some natural number $n \in \mathbb{Z}_{\geq 0}$.  We regard $H$ as being included in $G$ as the upper left $n \times n$ block, and note that $H$ also embeds in $\bar{G} := G / Z$.  We denote by $A_H, B_H, N_H$ the subgroups of $H$ defined by intersecting the corresponding subgroups of $G$ with $H$.  We denote by $Z_H$ the center of $H$.  We will take for $Q_H$ a Borel subgroup of $H$ containing $A_H$, chosen according to convenience, and denote by $U_H$ its unipotent radical.  We define $\mathfrak{a}_H^*$ and $\mathfrak{a}_{H,\mathbb{C}}^*$ in terms of $A_H$ as in \eqref{eq:frak-a-star-X-A}.

\subsubsection{Norms}
In the local case, we define a norm $\|.\|_G$ on $G$, denoted simply $\|.\|$ when $G$ is clear from context, by pulling back the operator norm on the special linear group under the given embedding $G \hookrightarrow \SL_r$.  In the global case, we define the norm $\|.\|_{G(\mathbb{A})}$ on $G(\mathbb{A})$ to be the product $\|g\|_{G(\mathbb{A})} := \prod_{\mathfrak{p}} \|g_\mathfrak{p} \|_{G(F_\mathfrak{p})}$ of the corresponding local norms defined using the induced embeddings $G(F_\mathfrak{p}) \hookrightarrow \SL_r(F_\mathfrak{p})$.

\subsubsection{Haar measures}\label{sec:haar-measures}
Since this paper is ultimately concerned with estimates, the choice of Haar measure is unimportant, provided that we are consistent.  We record a standard choice sufficient for our purposes.

In local and global cases, we equip $K$ with the Haar probability measure.

In the local case, we assume given Haar measures on $U$, $N$, $A$, $Z$.  In the non-archimedean case, we assume that maximal compact subgroups of each of these groups have volume one.  In the archimedean case, we assume only that the given Haar measures are fixed when $G$ is fixed.  We then equip $G$ with the Haar measure defined using the Iwasawa decomposition:
\begin{equation}\label{eq:int-_g-in-3}
  \int_{g \in G} f(g) \, d g
  = \int_{u \in U}
  \int_{a \in A}
  \int_{k \in K}
  f(u a k)
  \, d k
  \, \frac{d a}{\delta_U(a)}
  \, d u.
\end{equation}
Equivalently, the multiplication map $A \times U \times K \rightarrow G$ is measure-preserving.  We obtain a quotient measure on $U \backslash G$, given by
\begin{equation}\label{eq:int-_g-in-U-bakslash-G}
  \int_{g \in U \backslash G} f(g) \, d g
  =
  \int_{a \in A}
  \int_{k \in K}
  f(a k)
  \, d k
  \, \frac{d a}{\delta_U(a)}.
\end{equation}
These last identities remain valid with $U$ replaced by $N$.

Sometimes, it will be convenient instead to suppose that the multiplication map $A \times U \times N \rightarrow G$ is measure-preserving.  We will indicate when this is the case.

In the global case, we choose factorizable Haar measures on $U(\mathbb{A})$, $N(\mathbb{A})$, $A(\mathbb{A})$, $Z(\mathbb{A})$ and $K$ so that $[U]$ and $K$ have volume one.  We equip $G(\mathbb{A})$ with the Haar measure defined by analogy to \eqref{eq:int-_g-in-3}.  We fix compatible factorizations of the Haar measures on each of these groups, so that for each $\mathfrak{p}$, the analogue of \eqref{eq:int-_g-in-3} for $G(F_\mathfrak{p})$ is valid and the Haar measures on the groups $U(F_\mathfrak{p})$, $N(F_\mathfrak{p})$, $A(F_\mathfrak{p})$, $Z(F_\mathfrak{p})$ and $K_\mathfrak{p}$ satisfy the desiderata of the previous paragraph.  We equip the discrete groups consisting of rational points with the counting measure, adelic quotients $[G]$, $[A]$, $[Z]$, etc., with the quotient measure.

\subsection{Characters}\label{sec:cnjeor2mwy}
Given a locally compact abelian group $V$, we write $V^\wedge$ for its Pontryagin dual.  This notation applies in particular when $F$ is local and $V$ is a finite-dimensional vector space over $F$, in which case $V^\wedge$ may be regarded as a vector space over $F$ of the same dimension.  It applies also in the local (resp.\ global) case to $A$ (resp.\ $[A]$).

We denote $X(A)$ the group of rational characters $\chi : A \rightarrow \GL_1$.  We set
\begin{equation}\label{eq:frak-a-star-X-A}
  \mathfrak{a}^* := X(A) \otimes_{\mathbb{Z}} \mathbb{R} \hookrightarrow \mathfrak{a}_{\mathbb{C}}^* := X(A) \otimes_{\mathbb{Z}} \mathbb{C}.
\end{equation}
We denote by $\mathfrak{a} \hookrightarrow \mathfrak{a}_{\mathbb{C}}$ the corresponding dual spaces.  For $s \in \mathfrak{a}_{\mathbb{C}}^*$ and $t \in \mathfrak{a}_{\mathbb{C}}$, we denote by $t(s) = \langle s, t \rangle = \langle t,s  \rangle \in \mathbb{C}$ the natural pairing.

In the local (resp.\ global) case, we denote \index{characters!character group $\mathfrak{X}(A)$, $\mathfrak{X}([A])$} by $\mathfrak{X}(A)$ (resp.\ $\mathfrak{X}([A])$) the group of continuous complex-valued characters $\chi$ of $A$ (resp.\ $[A]$).  It is naturally a complex manifold, locally isomorphic to $\mathfrak{a}^*_{\mathbb{C}}$.  The subgroup of unitary characters is $A^\wedge \leq \mathfrak{X}(A)$ (resp.\ $[A]^\wedge \leq \mathfrak{X}([A])$).

Each $s \in \mathfrak{a}_{\mathbb{C}}^*$ naturally defines an element $|.|^s : a \mapsto |a|^s$ of $\mathfrak{X}(A)$ (resp.\ $\mathfrak{X}([A])$).  When $F = \mathbb{R}$ and $a \in A(\mathbb{R})^0$, we abbreviate $a^s := |a|^s$.  We recall that $\chi \in \mathfrak{X}(A)$ (resp.\ $\mathfrak{X}([A])$) is \emph{unramified} if it is invariant under the maximal compact subgroup of $A$ (resp.\ $[A]$); in that case, $\chi$ is of the form $|.|^{s}$ for some $s \in \mathfrak{a}_{\mathbb{C}}^*$ (in the local non-archimedean case, non-uniquely).  We denote by \index{characters!unramified character group $\mathfrak{X}^e(A)$, $\mathfrak{X}^e([A])$} $\mathfrak{X}^e(A)$ (resp.\ $\mathfrak{X}^e([A])$) the subgroup consisting of unramified characters.  (The superscript ``$e$'' is motivated by the archimedean case, where it signifies an ``evenness'' condition.)

For each element $\chi$ of $\mathfrak{X}(A)$ (resp.\ $\mathfrak{X}([A])$), there is a unique element $\Re(\chi) \in \mathfrak{a}^*$, called the \emph{real part} of $\chi$, for which $|\chi| = |.|^{\Re(\chi)}$.  Thus $A^\wedge$ (resp.\ $[A]^\wedge$) consists of all $\chi$ with $\Re(\chi) = 0$.

We equip $A^\wedge$ (resp.\ $[A]^\wedge$) with the Haar measure $d \chi$ dual to the Haar measure $d a$ on $A$ (resp.\ $[A]$).  For each $\sigma \in \mathfrak{a}^*$, the ``vertical'' subset consisting of $\chi \in \mathfrak{X}(A)$ (resp.\ $\mathfrak{X}([A])$) with $\Re(\chi) = \sigma$ is a coset of $A^\wedge$ (resp.\ $[A]^\wedge$); we equip it with the Haar measure $d \chi$ obtained by transporting that on the latter group.  For a function $f$ on $\mathfrak{X}^e(A)$ (resp.\ $\mathfrak{X}^e([A])$), we denote by \index{characters!Haar measure $d \mu_A(s)$, $d \mu_{[A]}(s)$}
\begin{equation}\label{eq:int-_sigma-fs}
  \int_{(\sigma)} f(s) \, d \mu_A(s)
\end{equation}
the integral of $f(\chi) \, d \chi$ over $\chi$ in $\mathfrak{X}^e(A)$ (resp.\ $\mathfrak{X}^e([A])$) of real part $\sigma$.

We denote by $R$ the set of roots for $(G,A)$, by $R_U^+$ (resp.\ $R_N^+$) the subset of positive roots relative to $U$ (resp.\ $N$), and $\Delta_U$ (resp.\ $\Delta_N$) the associated set of simple roots.\index{characters!positive roots $R_U^+, R_N^+$} \index{characters!simple roots $\Delta_U, \Delta_N$}  For $\alpha \in R$, we denote by $\alpha^\vee \in \mathfrak{a}^*$ the corresponding coroot.

We say that $\sigma \in \mathfrak{a}^*$ is \emph{dominant} if it is dominant in the sense defined by $U$, i.e., if for each $\alpha \in R_U^+$, we have $\alpha^\vee(s) \geq 0$.  We say that $\sigma$ is \emph{strictly dominant} if the stronger condition $\alpha^\vee(s) > 0$ holds for all such $\alpha$.  Given $\sigma_1, \sigma_2 \in \mathfrak{a}^*$, we write
\begin{equation*}
  \sigma_1 \prec \sigma_2 \quad \text{ or } \quad \sigma_2 \succ \sigma_1
\end{equation*}
to denote that the difference $\sigma_2 - \sigma_1$ is strictly dominant. \index{characters!strict dominance relations $\succ, \prec$}

(When we wish to speak of dominance relative to $N$, we write more verbosely ``$N$-dominant,'' but this will rarely happen except when $N = U$.)

We denote by $\delta_U, \delta_N \in \mathfrak{X}(A)$ (resp.\ $\mathfrak{X}([A])$) the modulus characters attached to $Q = A U$ and $B = A N$, respectively.  We write $\rho_U, \rho_N \in \mathfrak{a}^*$ for the elements characterized by the identities $\delta_U = |.|^{2 \rho_U}$, $\delta_N = |.|^{2 \rho_N}$.  Then $\rho_U = (1/2) \sum_{\alpha \in R_U^+} \alpha$, and similarly for $\rho_N$.\index{characters!positive root half-sums $\rho_U, \rho_N$} \index{characters!positive coroot half-sums $\rho_U^\vee, \rho_N^\vee$}  \index{characters!modulus $\delta_U, \delta_N$}  We set $\rho_U^\vee := (1/2) \sum_{\alpha \in R_U^+} \alpha^\vee$, and define $\rho_N^\vee$ analogously.

\subsection{Quotients and actions}
Let $P$ be a parabolic subgroup of $G$ that contains $A$.  We may write $P = M_P U_P$, where $U_P$ is the unipotent radical of $P$ and $M_P$ is the unique Levi factor that contains $A$.  In the local case, we may form the quotient $U_P \backslash G$.  In the global case, we define \index{adelic quotients!$[G]_P, [G]_Q$}
\begin{equation*}
  [G]_P := P(F) U_P(\mathbb{A}) \backslash G(\mathbb{A}).
\end{equation*}
We will be almost exclusively concerned with the special cases $P = G$ and $P = Q$.  We note that, in the global case, $[G]_G = [G]$ and $[G]_Q = Q(F) U(\mathbb{A}) \backslash G(\mathbb{A}) = A(F) U(\mathbb{A}) \backslash G(\mathbb{A})$.

In the local (resp.\ global) case, given a complex-valued function $f$ on $U_P \backslash G$ (resp.\ $[G]_P$) and $g \in G$, we denote by $R(g) f$ the function on the same space obtained via \index{representations!translation action, right $R(g)$} right translation:
\begin{equation*}
  R(g) f(x) := f(x g).
\end{equation*}
More generally, given a measure $\gamma$ on the relevant group ($G$ in the local case, $G(\mathbb{A})$ in the global case), we write $R(\gamma) f$ for the corresponding integrated right-translation operator, whenever the integral
\begin{equation*}
  R(\gamma) f(x) := \int f(x g) \, d \gamma(g)
\end{equation*}
converges.  We use the same convention for representations equipped with a right action of the relevant group.
The spaces $U \backslash G$ (resp.\ $[G]_Q$) admit natural left action by $A$ (resp.\ $[A]$), and we denote by $L$ the associated normalized action on function spaces, namely, for $a \in A$ (resp.\ $[A]$),\footnote{ This formula defines an action because $A$ is abelian.  Such left actions are often defined using $f(a^{-1} x)$; we found the present normalization slightly more convenient.  } \index{representations!translation action, left, normalized $L(a)$}
\begin{equation}\label{eq:normalized-left-translation}
  L(a) f(x) :=
  \delta_U^{-1/2}(a)
  f(a x).
\end{equation}
In the local archimedean case, we denote also simply by $R$ (resp.\ $L$) the induced actions of $\Lie(G)$ and $\mathfrak{U}(G)$ (resp.\ $\Lie(A)$ and $\mathfrak{U}(A)$) on the corresponding space of smooth functions.  In the global setting, we use similar notation for the differential actions induced by $G(F_\infty)$ and $A(F_\infty)$.

\subsection{Schwartz spaces}\label{sec:schwartz-spaces}
We refer to Casselman \cite{MR1001613} and the thesis of Moore \cite{MooreOSU2018} for background, although we note that in non-archimedean aspects, we do not adopt the same definitions as the latter reference.

Continue to let $P = M_P U_P$ be a parabolic containing $A$.

Recall that in the local (resp.\ global) case, we have defined a norm $\|.\|$ on $G$ (resp.\ $G(\mathbb{A})$).  In the local case, we define a norm on $U_P \backslash G$ by
\begin{equation*}
  \|g\|_{U_P \backslash G} := \inf_{u \in U_P} \|u g\|.
\end{equation*}
In the global case, we define a norm on $[G]_P$, hence on $[G]_Q$ and $[G]$, by
\begin{equation}\label{eq:g_g-:=-inf_gamma}
  \|g\|_{[G]_P} := \inf_{h \in P(F) U_P(\mathbb{A})} \|h g\|,
\end{equation}

Suppose for the moment that $F$ is fixed (or that $F$ is non-archimedean local and ``everything is unramified'').  In the local case, for $g = u a k$ with $(u,a,k) \in U \times A \times K$, we then have $\|g\|_{U \backslash G} \asymp \|a\|$.  In the global case, we have under analogous hypotheses
\begin{equation*}
  \|g\| _{[G]_Q} \asymp \|a\|_{[A]} := \inf_{h \in A(F)} \|h a\|.
\end{equation*}
(The proofs of these estimates may be reduced, via the given linear embedding of $G$, to the case $G = \GL_n$, where they may be verified directly.)

Suppose for the moment that $F$ is local.  We define the Schwartz space $\mathcal{S}(U_P \backslash G)$ as follows.  In the non-archimedean case, it consists of locally constant compactly-supported functions $f : U_P \backslash G \rightarrow \mathbb{C}$.  In the archimedean case, it consists of smooth functions satisfying, for all $x \in \mathfrak{U}(G)$ and $m \in \mathbb{Z}_{\geq 0}$,
\begin{equation*}
  \sup_{g \in U_P \backslash G} \|g\|_{U_P \backslash G}^m | R(x) f(g) | < \infty.
\end{equation*}
In particular, by specializing to $P = G$ and $P = Q$, we obtain definitions of $\mathcal{S}(G)$ and $\mathcal{S}(U \backslash G)$. \index{function spaces!$\mathcal{S}(U_P \backslash G), \mathcal{S}(G), \mathcal{S}(U \backslash G)$}

In the global case, we define the Schwartz space $\mathcal{S}([G]_P)$ to be the union, over all open subgroups $J$ of $G(\mathbb{A}_f)$, of the space of smooth functions $f : [G]_P \rightarrow \mathbb{C}$ that are right-invariant under $J$ and satisfy, for each $m \geq 0$ and $x \in \mathfrak{U}(G(F_\infty))$,
\begin{equation*}
  \sup_{g \in [G]_P} \|g\|_{[G]_P}^m |R(x) f(g)| < \infty.
\end{equation*}
In particular, we obtain definitions of $\mathcal{S}([G])$ and $\mathcal{S}([G]_Q)$.  \index{function spaces!$\mathcal{S}([G]_P), \mathcal{S}([G]), \mathcal{S}([G]_Q)$}

In the local (resp.\ global) case, the space $\mathcal{S}(U \backslash G)$ (resp.\ $\mathcal{S}([G]_Q)$) admits a natural left action by $A$ (resp.\ $[A]$).  We denote by $\mathcal{S}^e(U \backslash G)$ (resp.\ $\mathcal{S}^e([G]_Q)$) the subspaces consisting of elements that are \emph{left-unramified}, i.e., left-invariant under the maximal compact subgroup of $A$ (resp.\ $[A]$).  \index{function spaces!$\mathcal{S}^e(U \backslash G), \mathcal{S}^e([G]_Q)$}

All of the above definitions apply in particular to $G = A$.  In particular, in the local case, we may define the Schwartz space $\mathcal{S}(A)$ and its subspace $\mathcal{S}^e(A)$ consisting of elements invariant under the maximal compact subgroup.  \index{function spaces!$\mathcal{S}(A), \mathcal{S}^e(A)$}

In the global case, we define the space $\mathcal{T}([G]_P)$ of functions of uniform moderate growth to be the union, over $d \in \mathbb{Z}_{\geq 0}$ and open subgroups $J$ of $G(\mathbb{A}_f)$, of the space of smooth functions $f : [G]_P \rightarrow \mathbb{C}$ that are right-invariant by $J$ and satisfy, for each $x \in \mathfrak{U}(G(F_\infty))$,
\begin{equation*}
  \sup_{g \in [G]_P} \|g\|^{-d}_{[G]_P} | R(x) \varphi(g)| < \infty.
\end{equation*}
\index{function spaces!$\mathcal{T}([G]_P)$}

Each of the function spaces just defined is naturally a Fr{\'e}chet space in the local archimedean case and an inductive limit of Fr{\'e}chet spaces otherwise.

\subsection{Principal series representations}\label{sec:induc-repr}
In the local (resp.\ global) case, for $\chi \in \mathfrak{X}(A)$ (resp.\ $\mathfrak{X}([A])$), we denote by \index{representations!principal series $\mathcal{I}(\chi)$}
\begin{equation*}
  \mathcal{I}(\chi) = \Ind_Q^G(\chi) \subseteq C^\infty(U \backslash G)
  \quad
  \text{
    (resp.\ $\mathcal{I}(\chi) = \Ind_{Q(\mathbb{A})}^{G(\mathbb{A})}(\chi) \subseteq C^\infty([G]_Q)$)
  }
\end{equation*}
the corresponding normalized principal series representation, consisting of smooth complex-valued functions $f$ on $G$ (resp.\ $G(\mathbb{A})$) satisfying $f(u a g) = \delta_U^{1/2}(a) \chi(a) f(g)$ for all $(u,a,g)$ in $U \times A \times G$ (resp.\ $U(\mathbb{A}) \times A(\mathbb{A}) \times G(\mathbb{A})$).  Here ``smooth'' carries the following convention:
\begin{itemize}
\item In the archimedean case, it means that all derivatives (left-invariant, say) exist.
\item In the non-archimedean (resp.\ global) case, it means that there is an open subgroup $J$ of $G$ (resp.\ of $G(\mathbb{A}_f)$) under which $f$ is right-invariant.
\end{itemize}

In the global case, for a prime $\mathfrak{p}$ of $F$, we denote by $\chi_\mathfrak{p}$ the local component of $\chi$ at $\mathfrak{p}$ and by $I_\mathfrak{p}(\chi_\mathfrak{p})$ (or simply $\mathcal{I}_\mathfrak{p}(\chi)$) the corresponding principal series representation of $G(F_\mathfrak{p})$; then $\mathcal{I}(\chi)$ identifies naturally with the restricted tensor product of the $\mathcal{I}_\mathfrak{p}(\chi_\mathfrak{p})$, restricted with respect to the spherical elements taking the value $1$ at the identity.

An invariant pairing between $\mathcal{I}(\chi)$ and $\mathcal{I}(\chi^{-1})$ is given by $(f_1,f_2) \mapsto \int_{K} f_1 f_2$.  In the local case, up to a constant factor depending upon measure normalizations, this coincides with the pairing given by $\int_N f_1 f_2$.

We denote by $\|.\|$ the norm on $\mathcal{I}(\chi)$ given by $\|f\|^2 := \int_{K} |f|^2$, which is invariant when $\Re(\chi) = 0$.

When $\chi$ is unramified, i.e., $\chi = |.|^s$ for some $s \in \mathfrak{a}_{\mathbb{C}}^*$, we often write simply \index{representations!principal series $\mathcal{I}(s)$}
\begin{equation*}
  \mathcal{I}(s) := \mathcal{I}(|.|^s)
\end{equation*}
in both local and global cases.  (This notational abbreviation has the potential to introduce confusion -- should $\mathcal{I}(1)$ denote the induction of the trivial character $1 = |.|^0$, or that of $|.|^1$? -- so we will be careful when using it.)

\subsection{Mellin analysis and Paley--Wiener theory}\label{sec:mell-analys-paley}
\subsubsection{Warm-up}
We briefly recall the case of the real line.  We say that a holomorphic function $\tilde{f} : \mathbb{C} \rightarrow \mathbb{C}$ has \emph{rapid vertical decay} if for each finite interval $[a,b] \subset \mathbb{R}$ and $m \geq 0$, the quantity
\begin{equation}\label{eq:sup_s-in-mathbbc}
  \sup_{s \in \mathbb{C} : \Re(s) \in [a,b]}
  (1 + |s|)^m |\tilde{f}|(s)
\end{equation}
is finite.  The Schwartz space $\mathcal{S}(\mathbb{R}^\times_+)$ is defined to consist of all smooth $f : \mathbb{R}^\times_+ \rightarrow \mathbb{C}$ with the property that for each $m, k \geq 0$, we have
\begin{equation}\label{eq:sup_y-in-mathbbrt}
  \sup_{y \in \mathbb{R}^\times_+} (y + y^{-1})^m \left\lvert \left( y \frac{d}{d y} \right)^k f(y) \right\rvert < \infty.
\end{equation}
For each $f \in \mathcal{S}(\mathbb{R}^\times_+)$, the Mellin transform $\tilde{f}(s) := \int_{y \in \mathbb{R}^\times_+} f(y) y^{-s} \, d^\times y$, where $d^\times y := \frac{d y}{y}$, converges absolutely, and defines a holomorphic function $\tilde{f}$ of rapid vertical decay; conversely, every such function $\tilde{f}$ arises in this way.  Moreover, the topology on $\mathcal{S}(\mathbb{R}^\times_+)$, as described by the seminorms \eqref{eq:sup_y-in-mathbbrt}, is equivalent to the topology described by the seminorms \eqref{eq:sup_s-in-mathbbc}.  Finally, we have
\begin{equation*}
  \int  _{y \in \mathbb{R}^\times_+} |f(y)|^2 \, d^\times y
  =
  \int_{\Re(s) = 0} \left\lvert \tilde{f}(s) \right\rvert^2 \, \frac{d s}{2 \pi i }.
\end{equation*}

The above assertions summarize parts of classical Paley--Wiener theory, expressed in terms of the Mellin transform.  The proofs consist of elementary complex analysis.  A similar theory relates decay properties of functions $f : \mathbb{Z} \rightarrow \mathbb{C}$ to the regularity of the associated Laurent series $\tilde{f}(z) = \sum_{n \in \mathbb{Z}} f(n) z^n$.

We temporarily drop the general notation of \S\ref{sec:preliminaries}.  In \S\ref{sec:tori} and \S\ref{sec:mellin-paley-affine-quotients}, we consider groups $A$ acting on spaces $X$ with the following properties:
\begin{itemize}
\item $A$ is isomorphic $\mathbb{R}^{n_1} \times \mathbb{Z}^{n_2} \times \Gamma$ for some $n_1, n_2 \geq 0$ and some profinite group $\Gamma$.
\item There is a compact group $K$ and a surjective $A$-equivariant homomorphism $A \times K \rightarrow X$ with compact kernel.
\end{itemize}
In each case, we describe a form of Paley--Wiener theory for the action of $A$ on certain function spaces on $X$.  The proofs remain elementary.  Analogous results appear throughout \cite[II.1]{MR1361168}.

\subsubsection{Tori}\label{sec:tori}
We now return to the general thread of \S\ref{sec:preliminaries}.  When $F$ is archimedean, we use the given linear embedding of $G$ to equip the real Lie group $\Lie(A)$, hence its real dual $\Lie(A)^*$, with a $W$-invariant norm $\|.\|$.  We denote by $d \chi \in \Lie(A)^*$ the differential at the identity of a character $\chi \in \mathfrak{X}(A)$, and define the analytic conductor $C(\chi) := 3 + \|d \chi\|$.  (It is easy to see that this definition is comparable to the ``standard'' one, obtained by identifying $A$ with a product of copies of $\GL_1(F)$, hence $\chi$ with a tuple of characters $\chi_j$ of $F^\times$, and multiplying together the analogous conductors of the $\chi_j$ as defined in \S\cite[\S3.1.8]{michel-2009}.)

We consider here the local case, which is all that we require.  We say that a complex-valued holomorphic function $\tilde{f}$ on $\mathfrak{X}(A)$ has \emph{rapid vertical decay} if
\begin{enumerate}[(i)]
\item in the non-archimedean case, there is an open subgroup $J$ of $A$ so that $\tilde{f}(\chi) = 0$ unless $\chi|_J$ is trivial, and
\item in the archimedean case, for each compact subset $\mathcal{D} \subseteq \mathfrak{a}^*$ and $m \geq 0$, the supremum
  \begin{equation}\label{eq:sup_chi-in-mathfr}
    \sup_{
      \substack{
        \chi \in \mathfrak{X}([A]) :
        \\
        \Re(\chi) \in \mathcal{D}
      }
    }
    C(\chi)^m
    |\tilde{f}(\chi)|
  \end{equation}
  is finite.
\end{enumerate}
For $f \in \mathcal{S}(A)$, the Mellin transform $\tilde{f}(\chi) := \int_{a \in A} f(\chi) \chi^{-1}(a) \, d a$ converges absolutely, and defines a holomorphic function $\tilde{f} : \mathfrak{X}(A) \rightarrow \mathbb{C}$ of rapid decay; conversely, every such function $\tilde{f}$ arises in this way.  In the archimedean case, the topology on $\mathcal{S}(A)$ is the same as that defined by the seminorms \eqref{eq:sup_chi-in-mathfr}.  \index{characters!Mellin transform $\tilde{f}(\chi), \tilde{f}(s)$}

We have $f \in \mathcal{S}^e(A)$ if and only if $\tilde{f}$ is supported on $\mathfrak{X}^e(A)$.  In that case, $\tilde{f}$ is described entirely by the numbers
\begin{equation*}
  \tilde{f}(s) := \tilde{f}(|.|^s).
\end{equation*}
(We will be careful in using such notation for the same reasons as indicated at the end of \S\ref{sec:induc-repr}.)

\subsubsection{Affine quotients}\label{sec:mellin-paley-affine-quotients}
In the local (resp.\ global) case, let $E$ be an open subset of $\mathfrak{X}(A)$ (resp.\ $\mathfrak{X}([A])$).  By a \emph{holomorphic family $\{f[\chi]\}_{\chi \in E}$ of vectors $f[\chi] \in \mathcal{I}(\chi)$}, we mean a family of such vectors that arises as the pointwise specialization $f[\chi](g) := f(\chi,g)$ of some complex-valued function $f$ on $E \times U \backslash G$ (resp.\ $E \times [G]_Q$) that is holomorphic in the first variable and smooth in the second variable, where ``smooth'' follows the same convention as in \S\ref{sec:induc-repr}, with the subgroup $J$ taken independent of the first argument.  This definition applies in particular when $E$ is the full character group $\mathfrak{X}(A)$ (resp.\ $\mathfrak{X}([A])$), in which case we write simply ``holomorphic family $\{f[\chi]\}$.''

We say that a holomorphic family $\{f[\chi]\}$ has \emph{rapid vertical decay} if
\begin{itemize}
\item in the non-archimedean local case, there is an open subgroup $J$ of $A$ so that $f[\chi] = 0$ unless $\chi|_J$ is trivial,
\item in the archimedean local case, for each compact subset $\mathcal{D} \subseteq \mathfrak{a}^*$, $m \geq 0$, and $x \in \mathfrak{U}(G)$, the supremum
  \begin{equation}\label{eq:sup_chi-in-mathfr-1}
    \sup_{
      \substack{
        \chi \in \mathfrak{X}(A) :  \\
        \Re(\chi) \in \mathcal{D}
      }
    } C(\chi_\infty)^m \sup_{k \in K} |R(x) f[\chi](k)|
  \end{equation}
  is finite, and
\item in the global case, there is an open subgroup $J$ of $A(\mathbb{A}_f)$ so that $f[\chi] = 0$ unless $\chi|_J$ is trivial, and for each compact subset $\mathcal{D} \subseteq \mathfrak{a}^*$, $m \geq 0$ and $x \in \mathfrak{U}(G(F_\infty))$, the supremum
  \begin{equation}\label{eq:sup_-substack-chi}
    \sup_{
      \substack{
        \chi \in \mathfrak{X}([A]) :  \\
        \Re(\chi) \in \mathcal{D}
      }
    }  C(\chi)^m \sup_{k \in K} |R(x) f[\chi](k)|
  \end{equation}
  is finite.
\end{itemize}
We note that in the archimedean local and global cases, replacing $\sup_{k \in K} |R(x) f[\chi](k)|$ with $\|R(x) f[\chi] \|$ yields the same definition, thanks to the Sobolev lemma for $K$.  We note also that the vanishing conditions imposed on $f[\chi]$ in the non-archimedean local and global cases can be seen to follow from our definition of ``holomorphic family''; we have stated it, redundantly, for emphasis.

In the local (resp.\ global) case, given $f \in \mathcal{S}(U \backslash G)$ (resp.\ $\mathcal{S}([G]_Q)$), we define the Mellin components $f[\chi] \in \mathcal{I}(\chi)$ by the absolutely convergent integrals \index{principal series!Mellin components $f[\chi], f[s]$}
\begin{equation*}
  f[\chi](g) := \int \delta_U^{1/2}(a)  \chi(a) f(a^{-1} g) \, d a
\end{equation*}
taken over $a \in A$ (resp.\ $[A]$).  These integrals define a holomorphic family of vectors, of rapid vertical decay; conversely, every such family arises in this way.  In the archimedean local (resp.\ global) case, the topology on $\mathcal{S}(U \backslash G)$ (resp.\ $\mathcal{S}([G]_Q)$) is described by the family of seminorms \eqref{eq:sup_chi-in-mathfr-1} (resp.\ \eqref{eq:sup_-substack-chi}).

We have for each $\sigma \in \mathfrak{a}^*$ the Mellin decomposition
\begin{equation}\label{eq:f-=-int}
  f = \int_{(\sigma)} f[\chi] \, d \chi.
\end{equation}
When $\sigma = 0$, this extends to a decomposition of $L^2(U \backslash G)$ (resp.\ $L^2([G]_Q)$) as the direct integral of the $\mathcal{I}(\chi)$ with $\Re(\chi) = 0$ (cf.\ \cite[\S3.4.9, \S3.5.2]{MR3468638}).  In particular, for $f_1, f_2 \in \mathcal{S}(U \backslash G)$ (resp.\ $\mathcal{S}([G]_Q)$), we have the Parseval formula
\begin{equation}\label{eq:leftlangle-f_1-f_2}
  \left\langle f_1, f_2 \right\rangle
  =
  \int_{(0)}
  \left\langle f_1[\chi], f_2[\chi] \right\rangle \, d \chi.
\end{equation}

If $f \in \mathcal{S}^e(U \backslash G)$ (resp.\ $\mathcal{S}^e([G]_Q)$), then $f[\chi] = 0$ unless $\chi$ is unramified, and so $f$ is again described by the vectors
\begin{equation*}
  f[s] := f[|.|^s] \in
  \mathcal{I}(s) := \mathcal{I} (|.|^s).
\end{equation*}

\subsection{Sobolev norms}\label{sec:local-sobolev-norms-prelims}
Here we restrict to the local archimedean case.  The topology on $\mathcal{S}^e(U \backslash G)$ is described by the following family of seminorms, indexed by compact subsets $\mathcal{D} \subseteq \mathfrak{a}^*$ and $\ell \in \mathbb{Z}_{\geq 0}$:  \index{Sobolev norms!$\nu_{\mathcal{D},\ell}$}
\begin{equation}\label{eq:nu_mathcald-ellf-:=}
  \nu_{\mathcal{D},\ell}(f) :=
  \sup_{
    \substack{
      s \in \mathfrak{a}_{\mathbb{C}}^* :  \\
      \Re(s) \in \mathcal{D}
    }
  }
  \sum_{m=0}^{\ell}
  \sum_{x_1,\dotsc,x_m \in \mathcal{B}(G)}
  (1 + |s|)^{\ell}
  \|R(x_1 \dotsb x_m) f[s] \|.
\end{equation}
Here $\|.\|$ denotes the norm on $\mathcal{I}(s)$, induced as in \S\ref{sec:induc-repr} from $L^2(K)$.  For each $T \geq 1$, we define a rescaled variant of such seminorms: \index{Sobolev norms!$\nu_{\mathcal{D},\ell,T}$}
\begin{equation*}
  \nu_{\mathcal{D},\ell,T}(f) :=
  \sup_{
    \substack{
      s \in \mathfrak{a}_{\mathbb{C}}^* :  \\
      \Re(s) \in \mathcal{D}
    }
  }
  \sum_{m=0}^{\ell}
  \sum_{x_1,\dotsc,x_m \in \mathcal{B}(G)}
  (1 + |s|)^{\ell}
  \frac{\|R(x_1 \dotsb x_m) f[s] \|}{T^{\left\langle \rho_U^\vee , \Re(s) \right\rangle + m}}.
\end{equation*}
The informal idea here is that $\nu_{\mathcal{D},\ell,T}(f) = \O(1)$ provided that $f$ concentrates near $T^{-\rho_U^\vee}$, is uniformly smooth under left translation by $A$, has $L^2$-norm $\O(1)$, and has ``frequency'' $\O(T)$ under right translation by $G$ (see, e.g., Lemmas \ref{lem:standard:frakE-equivalences} and \ref{lem:standard:let-f-in-2} for details).

Let $\pi$ be a unitary representation of $G$.  More precisely, we write $\pi$ for the space of smooth vectors in the underlying Hilbert space representation.  For each $d \in \mathbb{Z}_{\geq 0}$, we define a Sobolev norm $\mathcal{S}_d$ on $\pi$ by the formula: for $v \in \pi$,  \index{Sobolev norms!$\mathcal{S}_d$}
\begin{equation}\label{eq:mathc-:=-sum_x_1}
  \mathcal{S}_d(v)^2 := \sum_{x_1,\dotsb,x_d \in \mathcal{B}(G) \cup \{1\}} \| \pi(x_1 \dotsb x_d) v\|^2,
\end{equation}
where $\|.\|$ denotes the given norm on $\pi$ and $1$ denotes the unit element of the universal enveloping algebra.

For each $(d_1,d_2) \in \mathbb{Z}_{\geq 0} \times \mathbb{Z}_{\geq 0}$, we define a Sobolev norm $\mathcal{S}_{d_1,d_2}$ on $\mathcal{S}(U \backslash G)$ by the formula: for $f \in \mathcal{S}(U \backslash G)$,  \index{Sobolev norms!$\mathcal{S}_{d_1,d_2}$}
\begin{equation*}
  \mathcal{S}_{d_1,d_2}(f)^2 :=
  \sum_{
    \substack{
      x_1, \dotsb, x_{d_1} \in \mathcal{B}(A) \cup \{1\}  \\
      y_1, \dotsc, y_{d_2 } \in \mathcal{B}(G) \cup \{1\}
    }
  }
  \| L(x_1 \dotsb x_{d_1}) R(y_1 \dotsb y_{d_2}) f \|^2,
\end{equation*}
where $\|.\|$ denotes the norm on $L^2(U \backslash G)$.  For example, $\mathcal{S}_{0,d}(f) = \mathcal{S}_{d}(f)$ as defined via \eqref{eq:mathc-:=-sum_x_1}, taking for $\pi$ the right regular representation $R$ of $G$ on $\mathcal{S}(U \backslash G)$.  By the Parseval relation on $A$, we have for $f \in \mathcal{S}^e(U \backslash G)$ the relation
\begin{equation}\label{eq:sobolev-parseval-A-U-G}
  \mathcal{S}_{d_1,d_2}(f)^2 =
  \sum_{
    y_1, \dotsc, y_{d_2 } \in \mathcal{B}(G) \cup \{1\}
  }
  \int_{
    \substack{
      s \in \mathfrak{a}_{\mathbb{C}}^* :  \\
      \Re(s) = 0
    }
  }
  E_{d_1}(s)
  \| R(y_1 \dotsb y_{d_2}) f[s] \|^2 \, d \mu_A(s),
\end{equation}
where
\begin{equation}\label{eq:e_ds-:=-sum}
  E_d(s) :=
  \sum_{
    x_1, \dotsc, x_{d } \in \mathcal{B}(A) \cup \{1\}
  }
  \left\lvert \langle s, x_1 \rangle \dotsb \langle s, x_d \rangle \right\rvert^2.
\end{equation}

We note that for all $a \in A$ and $(d_1,d_2,f)$ as above,
\begin{equation}\label{eq:mathcals_d_1-d_2la-f}
  \mathcal{S}_{d_1,d_2}(L(a) f) = \mathcal{S}_{d_1,d_2}(f),
\end{equation}
as follows from the commutativity of $A$ and the unitarity of the operator $L(a)$ on $L^2(U \backslash G)$.

Given a parameter $T \geq 1$, we define the rescaled Sobolev norm $\mathcal{S}_{d_1,d_2,T}$ on $\mathcal{S}(U \backslash G)$ by \index{Sobolev norms!$\mathcal{S}_{d_1,d_2,T}$}
\begin{equation}\label{eq:mathcals_d_1-d_2-tf}
  \mathcal{S}_{d_1,d_2,T}(f) :=
  \max_{0 \leq k \leq d_2}
  T^{-k}
  \mathcal{S}_{d_1,k}(f).
\end{equation}

\begin{lemma}\label{lem:standard2:let-d_1-d_2}
  Let $d_1, d_2 \in \mathbb{Z}_{\geq 0}$.  There exists $C_0 \geq 0$ and $\ell \in \mathbb{Z}_{\geq 0}$ so that for all compact $\mathcal{D} \subseteq \mathfrak{a}^*$ containing the origin, and all $f \in \mathcal{S}^e(U \backslash G)$, we have
  \begin{equation*}
    \mathcal{S}_{d_1,d_2}(f) \leq C_0 \nu_{\mathcal{D},\ell}(f).
  \end{equation*}
\end{lemma}
\begin{proof}
  Immediate from \eqref{eq:nu_mathcald-ellf-:=} and \eqref{eq:sobolev-parseval-A-U-G}; we may take $\ell = d_1 + b$ for any $b \in \mathbb{Z}_{\geq 0}$ large enough that $\int_{s \in \mathfrak{a}_{\mathbb{C}}^* : \Re(s) = 0} (1 + |s|)^{-2 b} \, d \mu_A(s) < \infty$.
\end{proof}

\subsection{Zeta functions}\label{sec:local-zeta-functions}
In the local case, we denote by $\zeta_F$ the meromorphic function on $\mathbb{C}$ given by \index{zeta functions!$\zeta_F(s)$}
\begin{equation*}
  \zeta_F(s) :=
  \begin{cases}
    \pi^{-s/2} \Gamma(s/2) & \text{ if } F = \mathbb{R}, \\
    2 (2 \pi ) ^{- s} \Gamma (s) &  \text{ if } F = \mathbb{C}, \\
    (1 - q^{-s})^{-1} & \text{ if } F: \text{non-archimedean},
  \end{cases}
\end{equation*}
with $q$ the residue field cardinality (see \S\ref{sec:crb1dyvdtv}).  We write $\zeta_F(U,\cdot)$ for the meromorphic function on $\mathfrak{a}_{\mathbb{C}}^*$ given by
\begin{equation*}
  \zeta_F(U,s) := \prod_{\alpha \in R_U^+} \zeta_F(1 + \alpha^\vee(s)).
\end{equation*}
We note that $\zeta_F(U,s)$ depends only upon the character $\chi = |.|^s$ of $A$.  Thus the shorthand \index{zeta functions!$\zeta_F(U,\chi), \zeta_F(U,s)$}
\begin{equation*}
  \zeta_F(U,\chi) := \zeta_F(U,s)
\end{equation*}
defines a meromorphic function $\zeta_F(U,\cdot)$ on $\mathfrak{X}^e(A)$.

In the global case, we abbreviate $\zeta_\mathfrak{p} := \zeta_{F_\mathfrak{p}}$.  Given a finite set of places $S$ containing the archimedean ones, we write $\zeta_F^{(S)}(U,s)$ or $\zeta_F^{(S)}(U,\chi)$ for the product over $\mathfrak{p} \notin S$  of the corresponding local factors, with the product defined either via meromorphic continuation from the domain where $\Re(s)$ is strictly dominant or by the analogous product of partial Dedekind zeta functions.  When $F = \mathbb{Q}$ and $S = \{\infty\}$, we write simply $\zeta(U,s)$ or $\zeta(U,\chi)$.

\subsection{Additive characters}\label{sec:stand-nond-char}
We denote by $\psi_{\mathbb{Q}}$ the ``standard'' nontrivial unitary character of $\mathbb{A}/\mathbb{Q}$, characterized by $\psi_{\mathbb{Q}}|_{\mathbb{R}} = e^{2 \pi i (\cdot)}$.  For each completion $F$ of $\mathbb{Q}$, we write $\psi_F$ for the restriction of $\psi_{\mathbb{Q}}$ to $F$.  Each local field $F$ of characteristic zero is a finite extension of some completion $F_0$ of $\mathbb{Q}$; we define $\psi_F$ to be the pullback of $\psi_{F_0}$ via the trace map $F \rightarrow F_0$.

We denote by $\U(1) := \{z \in \mathbb{C}^\times : |z| = 1\}$ the unitary group in one variable.

In the local (resp.\ global) case, let $\psi : N \rightarrow \U(1)$ (resp.\ $\psi : N(\mathbb{A}) \rightarrow \U(1)$) be a nondegenerate unitary character, trivial on $N(F)$ in the global case.  For example, if $G = \GL_r$, then
\begin{equation}\label{eq:psi-n-vs-psi-F-whittaker}
  \psi(n) = \psi_F\left(\sum_{i=1}^{r-1} c_i n_{i,i+1}\right)
\end{equation}
for some $c_1,\dotsc,c_{r-1} \in F^\times$.

Suppose $G = \GL_r$.  In the local (resp.\ global) case, the \emph{standard} nondegenerate unitary character $\psi$ of $N$ (resp.\ of $N(\mathbb{A})$, trivial on $N(F)$) is then defined to be
\begin{equation}\label{eq:psin-=-psi_f}
  \psi(n) = \psi_{F} \left( \sum_{i=1}^{r-1} n_{i,i+1} \right).
\end{equation}

Suppose now that $F$ is local and that $G$ arises from some split reductive group over $\mathfrak{o}$, which we continue to denote by $G$.  We may speak then of $\psi$ being \emph{unramified}.  For example, if $G = \GL_r$ and $\psi_F$ itself is unramified (i.e., $F = \mathbb{Q}_p$ or an unramified extension), then this says that each $|c_i| = 1$ in the expression \eqref{eq:psi-n-vs-psi-F-whittaker}.

For a reductive group $G$ over $\mathbb{Z}$ and a nondegenerate unitary character $\psi$ of $N(\mathbb{A})$, trivial on $N(\mathbb{Q})$, we say that $\psi$ is \emph{unramified} if its restriction to each $N(\mathbb{Q}_p)$ is unramified.  For example, when $G = \GL_n$, the standard character and its inverse are unramified, as is any $\psi$ for which each $c_i = \pm 1$ in \eqref{eq:psi-n-vs-psi-F-whittaker}.

\subsection{Whittaker functions}\label{sec:whittaker-functions-1}
Let $\psi$ be as in \S\ref{sec:stand-nond-char}.

In the local case, we denote by $C^\infty(N \backslash G, \psi)$ the space of smooth functions $W : G \rightarrow \mathbb{C}$ satisfying $W(u g) = \psi(u) W(g)$ for all $(u,g) \in N \times G$.  (We are using the letter $W$ to denote both the Weyl group and Whittaker functions, but it should always be clear by context what we mean.)  In the global case, we analogously define $C^\infty(N(\mathbb{A}) \backslash G(\mathbb{A}), \psi)$.

We denote by $w_{\bar{U},N} \in W$ the unique Weyl group element with $w_{\bar{U},N} N = \bar{U} w_{\bar{U},N}$, where $\bar{U}$ denotes the maximal unipotent subgroup opposite to $U$.  We choose a representative for $w_{\overline{U}, N}$ in $G(F)$.  The definitions that follow depend upon this choice.  When $G$ is a general linear group (the case ultimately relevant for this paper), we always take the standard choice, given by a permutation matrix; more generally, if $U$ is opposite to $N$, then we take for $w_{\overline{U}, N}$ the identity element of $G(F)$.  In the local case, $w_{\overline{U},N}$ then lies in the standard maximal compact subgroup.

In the local case, for $\chi \in \mathfrak{X}(A)$ and $f \in \mathcal{I}(\chi)$, we denote by
\begin{equation*}
  W[f,\psi] \in C^\infty(N \backslash G, \psi)
\end{equation*}
the corresponding Jacquet integral, defined for $\Re(\chi) \succ 0$ by the absolutely convergent integral \index{Jacquet integrals and distributions!$W[f,\psi]$}
\begin{equation*}
  W[f, \psi](g) := \int_{n \in N} \psi^{-1}(n)
  f(w_{\bar{U},N} n g) \, d n
\end{equation*}
and in general by analytic continuation along a holomorphic family of vectors (see \cite[\S3]{MR716292}, \cite[\S7.2]{MR2533003}, \cite[\S15]{MR1170566}, \cite{MR563369}).

In the global case, for a continuous function $\varphi$ on $[G]$, we define
\begin{equation*}
  W[\varphi,\psi] \in C^\infty(N(\mathbb{A}) \backslash G(\mathbb{A}), \psi)
\end{equation*}
by the analogous global integral:
\begin{equation*}
  W[\varphi, \psi](g) := \int_{n \in [N]} \psi^{-1}(n)
  \varphi(w_{\bar{U},N} n g) \, d n.
\end{equation*}

Assume now that $F$ is non-archimedean local, that $G$ arises from some split reductive group over $\mathfrak{o}$, which we continue to denote by $G$, that $K = G(\mathfrak{o})$, and that $w_{\bar{U},N} \in K$.  Assume that $\psi$ is unramified, $\chi$ is unramified, and $f^0[\chi]$ is the unramified vector with $f^0[\chi](1) = 1$.  Temporarily write $W_\chi := W[f^0[\chi], \psi]$ for the corresponding unramified Whittaker function.  Recall our convention that the Haar measure on $N$ assigns volume one to a maximal compact subgroup.  The Casselman--Shalika formula \cite[Thm 5.4]{MR581582} then reads
\begin{equation}\label{eq:wf0chi-psi1-casselman-shalika}
  W_\chi(1)
  =
  \frac{1}{\zeta_F(U,\chi)},
\end{equation}
where the RHS is as defined in \S\ref{sec:local-zeta-functions}.


We recall the general formula for $W_\chi$, which will be used only in the proof of Lemma \ref{lem:sub-gln:we-have-begin}.  To that end, recall that $G$ admits a complex dual group $G^\vee(\mathbb{C})$ which contains a maximal torus $A^\vee(\mathbb{C})$ dual to $A$: cocharacters $\lambda : F^\times \rightarrow A$ (i.e., algebraic homomorphisms) identify with characters $\lambda : A^\vee(\mathbb{C}) \rightarrow \mathbb{C}^\times$.  Define the (representative for the) Satake parameter $\alpha_\chi \in A^\vee(\mathbb{C})$ by requiring that
\begin{equation*}
  \alpha_\chi^\lambda = \chi(\varpi^\lambda),
\end{equation*}
where a superscripted $\lambda$ denotes the image under the (co)character $\lambda$.  In view of the Iwasawa decomposition, $W_\chi$ is determined by its values on elements of the form $\varpi^{\lambda}$.  By \emph{loc.\ cit.},\footnote{ The present formulation of Casselman's results is well-known, see for instance \cite[p.\ 504-505]{MR1431508} for the case $\mathrm{GL}_2$.} we have $W_\chi(\varpi^\lambda) = 0$ unless $\lambda$ is $N$-dominant, in which case
\begin{equation*}
  W_\chi(\varpi^\lambda) = \frac{\delta_N^{1/2}(\varpi^\lambda) \chi_\lambda(\alpha_\chi)}{\zeta_F(U, \chi)},
\end{equation*}
where $\chi_\lambda$ denotes the character of the irreducible algebraic representation of the complex dual group $G^\vee(\mathbb{C})$ of highest weight $\lambda$.

\begin{example}\label{example:cuhi8y36c4}
  Suppose $G = \mathrm{GL}_n(F)$, so that $G^\vee(\mathbb{C}) = \mathrm{GL}_n(\mathbb{C})$.  Then $\lambda$ identifies with an element of $\mathbb{Z}^n$ and $\varpi^\lambda = \diag(\varpi^{\lambda_1}, \dotsc, \varpi^{\lambda_n})$.  Take $N$ upper-triangular.  Then $\lambda$ is $N$-dominant precisely when $\lambda_1 \geq \dotsb \geq \lambda_n$.  Writing $\chi = \chi_s : \diag(a_1,\dotsc,a_n) \mapsto \prod \lvert a_i \rvert^{s_i}$, we have $\alpha_{\chi_s} = \diag(q^{- s_1}, \dotsc, q^{- s_n})$.  Finally, for $N$-dominant $\lambda$, $\chi_\lambda(\alpha_{\chi_s})$ is the Schur polynomial for $\lambda$ evaluated at $(q^{- s_1}, \dotsc, q^{- s_n})$.  (Compare with Shintani \cite{MR407208}.)
\end{example}

\subsection{Intertwining operators}\label{sec:local-intertw-oper-four}

\subsubsection{Basic definition}\label{sec:basic-definition-intertwining-ops}
Let $w \in W$.

In the local case, for $\chi \in \mathfrak{X}(A)$, we denote by
\begin{equation*}
  M_w[\chi] : \mathcal{I}(\chi) \rightarrow \mathcal{I}(w \chi),
\end{equation*}
or simply $M_w$ when $\chi$ is clear from context, the standard local intertwining operator, defined initially for $\Re(\chi) \succ 0$ by \index{intertwining operators!standard, $M_w$}
\begin{equation*}
  M_w f[\chi](g) = \int_{u \in
    w U w^{-1} \cap U \backslash  U} f[\chi](w^{-1} u g) \, d u
\end{equation*}
and in general via meromorphic continuation (see \cite[10.1.14]{MR1170566}, \cite[\S2.2]{MR610479}).  We note that the definition depends upon a lift of $w$ to $G$.

In the global case, for a function $f$ on $[G]_Q$, we denote by $M_w f : [G]_Q \rightarrow \mathbb{C}$ the image of $f$ under the standard intertwining operator $M_w$, defined by the usual integral, if it converges absolutely:
\begin{equation*}
  M_w f(g) = \int_{(w U w^{-1} \cap U \backslash U)(\mathbb{A})}
  f(w^{-1} u g) \, d u.
\end{equation*}
Here we lift $w$ to an element of $G(F)$.  We recall the basic properties of such integrals (see \cite{MR546601, 2019arXiv191102342B, MR2402686}).  They converge if $f \in \mathcal{S}([G]_Q)$, in which case they induce continuous maps $\mathcal{S}([G]_Q) \rightarrow \mathcal{T}([G]_Q)$.  They also converge if $f \in \mathcal{I}(\chi)$ for some $\chi$ with $\Re(\chi) \succ \rho_U$.  For a holomorphic family $\{f[\chi]\}_{\chi \in \mathfrak{X}([A])}$, the assignment $\chi \mapsto M_w f[\chi]$ admits a meromorphic continuation to all $\chi$, inducing a meromorphic family of maps
\begin{equation*}
  M_w := M_w[\chi] : \mathcal{I}(\chi) \rightarrow \mathcal{I}(w \chi)
\end{equation*}
This family of operators is holomorphic in a neighborhood of the subgroup $[A]^\wedge$ of unitary characters, and is unitary on that axis: for $\chi \in [A]^\wedge$ and $f_1, f_2 \in \mathcal{I}(\chi)$, we have
\begin{equation}\label{eq:langle-f_1-f_2}
  \langle f_1, f_2 \rangle = \langle M_w f_1, M_w f_2 \rangle.
\end{equation}
For $w_1, w_2 \in W$ and $\chi \in \mathfrak{X}(A)$, we have the composition formula
\begin{equation}\label{eq:m_w_1w_2-chi-circ}
  M_{w_1}[w_2 \chi] \circ M_{w_2}[\chi] = M_{w_1 w_2}[\chi] : \mathcal{I}(\chi) \rightarrow \mathcal{I}(w_1 w_2 \chi).
\end{equation}

\subsubsection{Fourier transforms}\label{sec:fourier-transforms}
In the local case, given a nondegenerate unitary character $\psi$ of $N$, we denote by \index{intertwining operators!normalized, $\mathcal{F}_{w,\psi}$}
\begin{equation*}
  \mathcal{F}_{w,\psi}[\chi] : \mathcal{I}(\chi) \rightarrow \mathcal{I}(w \chi),
\end{equation*}
or simply $\mathcal{F}_{w}[\chi]$, $\mathcal{F}_{w,\psi}$ or $\mathcal{F}_w$ when $\psi$ and/or $\chi$ are understood by context, the normalized intertwining operator as in \cite[p336]{MR610479}, i.e., the scalar multiple of $M_w[\chi]$ for which
\begin{equation}\label{eq:wmathcalf_w-fchi-psi}
  W[\mathcal{F}_w f, \psi] = W[f,\psi] \quad \text{ for all } f \in \mathcal{I}(\chi).
\end{equation}
We note that, while $M_w$ depends upon a choice of representative for $w$, the normalized operators $\mathcal{F}_{w}$ are independent of that choice, although they do depend upon the representatives $w_{\bar{U},N}$ used to define our Whittaker functionals.

By \cite[Prop 3.1.4]{MR610479}, we have
\begin{equation}\label{eq:mathc-1-mathc}
  \mathcal{F}_{w^{-1}} \mathcal{F}_w  f = f
  \quad \text{ for all } f \in \mathcal{I}(\chi)
\end{equation}
and
\begin{equation}\label{eq:f_1chi-f_2chi-1}
  (f_1, f_2) =
  (\mathcal{F}_w f_1, \mathcal{F}_w f_2)
  \quad
  \text{ for all } (f_1,f_2,w) \in \mathcal{I}(\chi) \times \mathcal{I}(\chi^{-1}) \times W.
\end{equation}
In particular, $\mathcal{F}_w[\chi]$ is unitary when $\chi$ is unitary, and so the family $\{\mathcal{F}_w[\chi] : \chi \in A^\wedge \}$ induces a unitary operator $\mathcal{F}_w$ on $L^2(U \backslash G)$ (cf.\ \cite[\S3.5]{MR3468638}, \cite{MR1694894,MR1988971,MR3969881}).  It is also known that
\begin{equation*}
  \mathcal{F}_{w_1} \circ \mathcal{F}_{w_2} = \mathcal{F}_{w_1 w_2},
\end{equation*}
see, e.g., \cite[Thm 2.1, $(R_2)$]{MR999488}.

\begin{lemma}\label{lem:standard2:let-w-in}
  Let $w \in W$ and $f \in \mathcal{S}(U \backslash G)$.  The following are equivalent:
  \begin{enumerate}[(i)]
  \item \label{itm:standard2:F-w-unitary} The unitary operator $\mathcal{F}_w$ on $L^2(U \backslash G)$ satisfies $\mathcal{F}_w f = f$.
  \item \label{itm:standard2:F-w-mellin} For almost all $\chi \in \mathfrak{X}(A)$, the operators $\mathcal{F}_w : \mathcal{I}(\chi) \rightarrow \mathcal{I}(w \chi)$ satisfy $\mathcal{F}_w f[\chi] = f [ w \chi]$.
  \end{enumerate}
\end{lemma}
\begin{proof}
  By definition, $\mathcal{F}_w f$ is the element of $L^2(U \backslash G)$ that, with respect to the direct integral decomposition noted in \eqref{eq:f-=-int}, has components $(\mathcal{F}_w f)[w \chi]$ given, for almost all $\chi$ with $\Re(\chi) = 0$, by $\mathcal{F}_w f[\chi]$.  The identity $\mathcal{F}_w f = f$ is thus equivalent to:
  \begin{itemize}
  \item $\mathcal{F}_w f[\chi] = f[w \chi]$ for almost all $\chi$ with $\Re(\chi) = 0$.
  \end{itemize}
  This last identity should be understood initially as one of elements of $L^2(K)$, but both sides are smooth under right translation by $G$, so we can interpret it pointwise.  Both sides of this last identity vary meromorphically with respect to $\chi$, so it is equivalent to ask that the same identity hold for almost all $\chi \in \mathfrak{X}(A)$.
\end{proof}

For $w \in W$, we denote by \index{function spaces!$\mathcal{S}(U \backslash G)^w, \mathcal{S}^e(U \backslash G)^w$}
\begin{equation}\label{eq:mathc-backsl-gw}
  \mathcal{S}(U \backslash G)^w
\end{equation}
the space of all $f \in \mathcal{S}(U \backslash G)$ satisfying the equivalent conditions of Lemma \ref{lem:standard2:let-w-in}.  We define $\mathcal{S}^e(U \backslash G)^w$ analogously, and then define $\mathcal{S}(U \backslash G)^W$ and $\mathcal{S}^e(U \backslash G)^W$ by taking the intersection over all $w \in W$.  \index{function spaces!$\mathcal{S}(U \backslash G)^W, \mathcal{S}^e(U \backslash G)^W$}

\begin{remark}\label{rmk:we-note-that-dependence-psi-w}
  We note that, despite what the notation $\mathcal{S}(U \backslash G)^w$ suggests, these spaces depend upon $\psi$, since $\mathcal{F}_{w,\psi}$ does.  We might write more verbosely $\mathcal{S}(U \backslash G)^{w, \psi}$, etc., to reflect this dependence.  In what follows, the choice of $\psi$ will always be specified or clear from context.

  When $G$ is a general linear group, we adopt the convention that $S(U \backslash G)^W := \mathcal{S}(U \backslash G)^{W,\psi}$, with $\psi$ the standard character $\psi$ as in \eqref{eq:psin-=-psi_f}.  We adopt the analogous convention concerning $\mathcal{S}^e(U \backslash G)^W$.  Using that the element
  \begin{equation*}
    \eta = \diag(1,-1,1,-1,\dotsc) \in A \cap K
  \end{equation*}
  conjugates $\psi$ to $\psi^{-1}$, we may check that
  \begin{equation}\label{eq:mathc-backsl-gw-1}
    \mathcal{S}^e(U \backslash G)^W := \mathcal{S}^e(U \backslash G)^{W,\psi} = \mathcal{S}^e(U \backslash G)^{W,\psi^{-1}}.
  \end{equation}
\end{remark}

\begin{remark}
  The Mellin-theoretic condition \eqref{itm:standard2:F-w-mellin} is the relevant one for our purposes.  We do not make direct use of the $L^2$-based condition \eqref{itm:standard2:F-w-unitary}.  We have included Lemma \ref{lem:standard2:let-w-in} primarily to simplify our notation: condition \eqref{itm:standard2:F-w-unitary} is somewhat more concise than \eqref{itm:standard2:F-w-mellin} and better justifies the use of the notation \eqref{eq:mathc-backsl-gw}.  The above discussion could likely be formulated more naturally in terms of extended Schwartz spaces (see Remark \ref{rmk:one-could-likely}).
\end{remark}

\begin{example}
  Suppose $G = \SL_2$, $U$ is lower-triangular, $N$ is upper-triangular, and $\psi : N \rightarrow \U(1)$ is given by $\begin{pmatrix}
    1 & x \\
    0 & 1
  \end{pmatrix} \mapsto \psi_F(x)$. Then $U \backslash G$ identifies, via the ``top row'' map, with the punctured plane $F^2 - \{0\}$.  The Haar measure on $U \backslash G$ identifies with the restriction of the Haar measure on the plane $F^2$.  For the nontrivial element $w \in W$, one can check that the unitary operator $\mathcal{F}_{\psi} := \mathcal{F}_{w,\psi}$ identifies with the symplectic Fourier transform on $L^2(F^2)$, given by
  \begin{equation*}
    \mathcal{F}_{\psi} f(x,y) =
    \int_{u, v \in F}
    f(u, v) \psi_F(u y - x v)
    \, d u \, d v
  \end{equation*}
  with $d u$ and $d v$ the $\psi_F$-self-dual Haar measures.  Related calculations appear in \cite[\S3.4--\S3.6]{MR3468638}.
\end{example}

\subsubsection{Euler factorization}\label{sec:fact-intertw-oper}
Suppose for the moment that $F$ is global.  Intertwining operators may be factored, in a certain sense, over the places of $F$.  Recall that $\mathcal{I}(\chi)$ identifies with the restricted tensor product of the local principal series representations $\mathcal{I}_\mathfrak{p}(\chi_\mathfrak{p})$.  The restriction is with respect to the normalized spherical elements $f_\mathfrak{p}^0[\chi_\mathfrak{p}]$ defined for almost all $\mathfrak{p}$, with the normalization being $f_\mathfrak{p}^0[\chi_\mathfrak{p}](1) = 1$.  For a pure tensor $f = \otimes_\mathfrak{p} f_\mathfrak{p}$ in this restricted tensor product, we have
\begin{equation*}
  M_w f = \otimes_\mathfrak{p} ^* M_w f_\mathfrak{p},
\end{equation*}
where $M_w$ is the local intertwining operator, defined  for $\Re(\chi_\mathfrak{p}) \succ 0$ by integrating over $U(F_\mathfrak{p})$, and the meaning of $\otimes_\mathfrak{p}^*$ is that we regularize the tensor product using the ratio of $L$-functions describing the images $M_w f_\mathfrak{p}^0$ of the unramified vectors, as given in \cite[(4)]{MR0419366}.

Let $\psi$ be a nondegenerate unitary character of $N(\mathbb{A})$, trivial on $N(F)$.  Write $\psi_\mathfrak{p} : N(F_\mathfrak{p}) \rightarrow \U(1)$ for its local component at a place $\mathfrak{p}$.  As explained in \cite[\S4]{MR610479}, the functional equation for Eisenstein series (see, e.g., \eqref{eq:eisenstein-series-functional-equation} below) implies that
\begin{equation}\label{eq:m_w-f-=}
  M_w f = \otimes_{\mathfrak{p}}^* \mathcal{F}_{w} f_\mathfrak{p},
\end{equation}
where
\begin{equation*}
  \mathcal{F}_{w} f_\mathfrak{p} := \mathcal{F}_{w,\psi_\mathfrak{p}} f_\mathfrak{p}
\end{equation*}
and $\otimes^*_{\mathfrak{p}}$ is defined as above, but using a different ratio of $L$-values.  More precisely, suppose that $S$ is a finite set of places taken large enough that $f_\mathfrak{p} = f_\mathfrak{p}^0[\chi_\mathfrak{p}]$ and $\psi_\mathfrak{p}$ is unramified for all $\mathfrak{p} \notin S$.  Then
\begin{align*}
  \mathcal{F}_{w} \zeta_{F_\mathfrak{p}}(U, \chi_\mathfrak{p}) f_\mathfrak{p}^0[\chi_\mathfrak{p}] =
  \zeta_{F_\mathfrak{p}}(U,w \chi_\mathfrak{p})f_\mathfrak{p}^0[w \chi_\mathfrak{p}],
\end{align*}
for $\mathfrak{p} \notin S$.  The precise meaning of \eqref{eq:m_w-f-=} is then that
\begin{equation}\label{eq:m_w-f-=-1}
  M_w f
  =
  \frac{
    \zeta_F^{(S)}(U,w \chi)
  }{
    \zeta_F^{(S)}(U,\chi)
  }
  (\otimes_{\mathfrak{p} \in S}
  \mathcal{F}_w f_\mathfrak{p} ) \otimes (
  \otimes_{\mathfrak{p} \notin S} f_\mathfrak{p}^0).
\end{equation}

A more convenient formulation of \eqref{eq:m_w-f-=-1} may be given as follows.  Introduce the modified spherical elements
\begin{equation*}
  f_\mathfrak{p}^\sharp[\chi_\mathfrak{p}] := \zeta_{F_\mathfrak{p}}(U,\chi_\mathfrak{p}) f_\mathfrak{p}^0[\chi_\mathfrak{p}].
\end{equation*}
For $f \in \mathcal{I}(\chi)$, write
\begin{equation*}
  f = \otimes^\sharp f_\mathfrak{p}
\end{equation*}
to signify that for some large enough finite set of places $S$, we have $f_\mathfrak{p} = f_\mathfrak{p} ^\sharp [\chi_\mathfrak{p}]$ for all $\mathfrak{p} \notin S$, and
\begin{equation*}
  f = \zeta_F^{(S)}(U,\chi)
  (\otimes_{\mathfrak{p} \in S} f_\mathfrak{p} )
  \otimes
  (\otimes_{\mathfrak{p} \notin S} f_\mathfrak{p}^0 ).
\end{equation*}
Then for $\mathfrak{p} \notin S$, we have
\begin{equation}\label{eq:normalized-F-w-preserves-basic-vector}
  \mathcal{F}_w f_\mathfrak{p}^\sharp[\chi_\mathfrak{p}] = f_\mathfrak{p} ^\sharp [w \chi_\mathfrak{p}],
\end{equation}
and so we may rewrite \eqref{eq:m_w-f-=-1} as follows:
\begin{equation}\label{eq:m_w-f-=-2}
  M_w f = \otimes^\sharp \mathcal{F}_w f_\mathfrak{p}.
\end{equation}

\begin{remark}
  Such relations may be derived also from \cite[p373]{MR3468638}.
\end{remark}

\subsection{Representation theory of general linear groups}

\subsubsection{Kirillov model}\label{sec:local-prelim-kirillov-model}
Suppose now that $F$ is local and that $G$ is a nontrivial general linear group.

We recall that an irreducible (smooth admissible) representation $\pi$ of $G$ is \emph{generic} if it admits an equivariant embedding into the space $C^\infty(N \backslash G, \psi)$ for some (equivalently, any) nondegenerate unitary character $\psi$ of $N$; such an embedding is then uniquely determined up to a scalar (see \cite{MR348047}, \cite[\S1]{MR546599}), so its image, denoted $\mathcal{W}(\pi,\psi)$, is well-defined.

We extend $G$ to a general linear pair $(G,H)$, with $H$ possibly trivial.  The restriction $\psi|_{N_H}$ is a nondegenerate unitary character of $N_H$; for notational simplicity, we denote that restriction simply by $\psi$.  The restriction map from $\mathcal{W}(\pi,\psi)$ to the space $C^\infty(N \backslash H, \psi)$ is injective (see \cite[\S6.5]{MR748505}, \cite[\S10.2]{MR1999922}), with image called the \emph{Kirillov model} of $\pi$.  If $\pi$ is unitarizable, then an invariant norm on $\mathcal{W}(\pi,\psi)$ may be given by (see \cite[\S6.4]{MR748505}, \cite[\S10.2]{MR1999922})
\begin{equation*}
  \|W\|^2 := \int_{N_H \backslash H} |W|^2.
\end{equation*}
We discuss the finer details of the archimedean Kirillov model in \S\ref{sec:kirillov-model}.

\subsubsection{Essential vector}\label{sec:essential-vector}
We retain the setting of \S\ref{sec:local-prelim-kirillov-model}, but assume that $F$ is non-archimedean and $\psi$ is unramified.  Let $\pi$ be a generic irreducible representation of $G$, with Whittaker model $\mathcal{W}(\pi,\psi)$.

For a nonzero integral $\mathfrak{o}$-ideal $\mathfrak{q}$, we write $K_1(\mathfrak{q}) \leq K_0(\mathfrak{q}) \leq K$ for the standard congruence subgroups, consisting of those elements with bottom row lying in $(\mathfrak{q},\dotsc,\mathfrak{q},1 + \mathfrak{q})$ and $(\mathfrak{q},\dotsc,\mathfrak{q},\mathfrak{o})$, respectively.  Then $K_1(\mathfrak{q})$ is a normal subgroup of $K_0(\mathfrak{q})$, with quotient isomorphic to $\mathfrak{o}^\times / (\mathfrak{o}^\times \cap (1 + \mathfrak{q}))$.

For a subgroup $J$ of $G$, we denote by $\pi^J$ the subspace of $\pi$ consisting of vectors invariant by $J$.

We recall some results of \cite{MR620708, MR3138844, MR3001803}.
\begin{theorem}
  Let $\pi$ be a generic irreducible representation of $G$.  There is a unique nonzero integral $\mathfrak{o}$-ideal $\mathfrak{c} = \mathfrak{c}(\pi)$ so that $\pi^{K_1(\mathfrak{q})}$ is nontrivial precisely when $\mathfrak{q} \subseteq \mathfrak{c}$.  The space $\pi^{K_1(\mathfrak{c})}$ is one-dimensional.  Let $W \in \mathcal{W}(\pi,\psi)^{K_1(\mathfrak{c})}$ be any basis element.  Then $W(1) \neq 0$.
\end{theorem}
We refer to $\mathfrak{c}$ as the \emph{conductor ideal} of $\pi$ and to its absolute norm $C(\pi) := \# \mathfrak{o} / \mathfrak{c}$ as the \emph{conductor} of $\pi$.  We refer to any nonzero element of that space as an \emph{essential vector}.  We will often impose the normalization $W(1) = 1$.

If $\mathfrak{c} = \mathfrak{o}$, then $K_1(\mathfrak{c}) = K_0(\mathfrak{c}) = K$.  Otherwise, writing $\pi|_Z : F^\times \rightarrow \mathbb{C}^\times$ for the central character of $\pi$, the group $K_0(\mathfrak{c})$ acts on $\pi^{K_1(\mathfrak{c})}$ via the character $g \mapsto \pi|_Z(d_g)$ of $K_1(\mathfrak{c})$, where $d_g$ denotes the lower-right entry of $g$.  In particular, the conductor ideal of the central character $\pi|_Z$ divides that of the representation $\pi$.

\subsubsection{Standard $L$-factors}\label{sec:standard-l-factors}
Let $F$ be local, and let $G$ be a general linear group.  Let $\pi$ be a generic irreducible admissible representation of $G$.  The standard $L$-factor $L(\pi,s)$ was defined first by Godement--Jacquet \cite{MR0342495} as the greatest common divisor of a family of zeta integrals attached to the matrix coefficients of $\pi$.  It was later interpreted more broadly, by Jacquet--Piatetski-Shapiro--Shalika (see \cite[\S5.1]{MR701565}, \cite[Thm 4.3]{MR519356}) as the $\GL_n \times \GL_1$ case of Rankin--Selberg convolutions for $\GL_n \times \GL_m$.  We refer to the summaries \cite{MR546609} and \cite[\S2 and Appendix]{MR1395406} for further discussion.

Such $L$-factors may be expressed as finite products $L(\pi,s) = \prod_{j=1}^m \zeta_F(s + \mu_{j})$, where $0 \leq m \leq n$ and $\mu_j \in \mathbb{C}$.  If $F$ is archimedean, or if $F$ is non-archimedean and $\pi$ is unramified, then $m = n$.  We remark that if $\pi$ is unitary, then each $\Re(\mu_j) > -1/2$ \cite[\S2.5]{MR618323}, \cite[Prop 2.1]{MR1277950}, so $L(\pi,s)$ is holomorphic in the region $\Re(s) \geq 1/2-\eps$ for some $\eps > 0$.

\subsubsection{Zeta integrals}\label{sec:zeta-integrals}
Let $F$ be local and let $(G,H)$ be a general linear pair.  We summarize here some basic facts concerning local Rankin--Selberg integrals for $G \times H$.

As above, we choose a nondegenerate unitary character $\psi$ of $N$, and denote also by $\psi$ its restriction to $N_H$.

Let $\pi$ and $\sigma$ be generic irreducible representations of $G$ and $H$, respectively.  The following discussion applies also to
Whittaker-type representations in the sense of \cite[\S1.5]{MR3753910}, i.e., inductions of essentially square-integrable representations.  In particular, it applies to principal series representations in the sense of \S\ref{sec:induc-repr}.

Recall (\S\ref{sec:haar-measures}) that we have chosen Haar measures on $H$ and $N_H$ that, in the non-archimedean case, assign volume one to maximal compact subgroups.

For $(W,V,u) \in \mathcal{W}(\pi,\psi) \times \mathcal{W}(\sigma,\psi^{-1}) \times \mathbb{C}$, we define the zeta integral
\[
  Z(W, V, u) := \int_{h \in N_H \backslash H} W
  \begin{pmatrix}
    h &  \\
      & 1
  \end{pmatrix}
  V(h) |\det h|^{u-1/2} \, d h,
\]
provided that it converges. \index{local zeta integrals!$Z(W,V,u)$}

We recall some results of Jacquet, Piatetski-Shapiro and Shalika.
\begin{theorem}\label{thm:local-zeta-rankin-selberg-meromorphic-continuation}
  The integral $Z(W,V,u)$ converges absolutely for $\Re(u)$ sufficiently large in terms of $\pi$ and $\sigma$.  It extends meromorphically to the complex plane.  The normalized zeta integral
  \[
    \frac{Z(W,V,u)}{L(\pi \times \sigma, u)}
  \]
  extends to an entire function of $u$.
\end{theorem}
\begin{proof}
  See  Jacquet--Piatetski-Shapiro--Shalika \cite[\S2.7]{MR701565} and Jacquet \cite[\S5, \S8.1]{MR2533003}.
\end{proof}
The definition of $Z(W,V,u)$ is particularly simple in the special case that the restriction of the Whittaker function $W$ to the subgroup $H$ has compact support modulo $N_H$.  In that case, no meromorphic continuation is necessary - the zeta integral is entire in $u$.  This trivial case of Theorem \ref{thm:local-zeta-rankin-selberg-meromorphic-continuation} is in fact the only case required for the purposes of this paper.

We now specialize to the case that $\sigma$ is a principal series representation $\mathcal{I}(s)$ of $H$, as defined in \S\ref{sec:induc-repr} for $s \in \mathfrak{a}_{H,\mathbb{C}}^*$.  Write $(G,H) = (\GL_{n+1}(F), \GL_n(F))$.  Using standard diagonal coordinates, we may identify $\mathfrak{a}_{H,\mathbb{C}}^*$ with $\mathbb{C}^n$, hence $s$ with an $n$-tuple $(s_1,\dotsc,s_n)$ of complex numbers $s_j$.
\begin{theorem}\label{thm:jpss-essential-vector-factorized}
  Assume that $F$ is non-archimedean, that $\psi$ is unramified, and that $W \in \mathcal{W}(\pi,\psi)$ is the essential vector (\S\ref{sec:essential-vector}) with $W(1) = 1$.  Let $V \in \mathcal{W}(\mathcal{I}(s), \psi^{-1})$ be the spherical vector with $V(1) = 1$.  Then
  \begin{equation*}
    Z(W,V, u) = L(\pi \times \mathcal{I}(s), u) = \prod_{j=1}^{n} L(\pi, s_j + u),
  \end{equation*}
  where the RHS denotes a product of standard $L$-factors as in \S\ref{sec:standard-l-factors}.
\end{theorem}
\begin{proof}
  For the first identity, see Jacquet--Piatetski-Shapiro--Shalika \cite[p208]{MR620708}, Matringe \cite[Cor 3.3]{MR3138844} or the notes of Cogdell \cite[Thm 3.3]{MR2508768}.  For the second identity, see Jacquet--Piatetski-Shapiro--Shalika \cite[\S9.4, \S8.4, \S5.1]{MR701565}, \cite[Thm 4.3]{MR519356}.
\end{proof}

We abbreviate
\begin{equation*}
  Z(W,V) := Z(W,V,u)|_{u=1/2},
\end{equation*}
where the RHS is interpreted as the value (possibly infinite) of the meromorphic function $u \mapsto Z(W,V,u)$.  \index{local zeta integrals!$Z(W,V)$}

Given $f \in \mathcal{I}(s)$, we set
\begin{equation*}
  Z(W,f) := Z(W, W[f,\psi^{-1}]),
\end{equation*}
where $W[f,\psi^{-1}] \in \mathcal{W}(\mathcal{I}(s), \psi^{-1})$ denotes the image of $f$ under the Jacquet integral, as defined in \S\ref{sec:whittaker-functions-1}.  \index{local zeta integrals!$Z(W,f)$}

\subsubsection{Standard $L$-functions}\label{sec:standard-l-functions}
For simplicity, we assume here that $F = \mathbb{Q}$.

Let $G = \GL_n$ be a general linear group, $\pi$ a cuspidal automorphic representation of $\GL_n(\mathbb{A})$.  We recall some properties of the standard $L$-function $L(\pi,s)$ introduced by Godement--Jacquet \cite[Lemma 10.20]{MR0342495}.  We refer to Jacquet \cite{MR546609}, Rudnick--Sarnak \cite[\S2]{MR1395406}, and Iwaniec--Kowalski \cite[\S5]{MR2061214} for further background and details.

For $s \in \mathbb{C}$, we may form the twisted representation $\pi \otimes |.|^s$, consisting of all functions $g \mapsto \varphi(g) |\det g|^s$ with $\varphi \in \pi$.  There is a unique $\sigma \in \mathbb{R}$ for which the twist $\pi \otimes |.|^\sigma$ has unitary central character, and in that case, the Petersson inner product equips $\pi$ with the structure of a unitary representation.  There is thus little loss of generality in supposing $\pi$ to be unitary.

For each place $\mathfrak{p}$, the local component $\pi_\mathfrak{p}$ is a generic irreducible representation of $\GL_n(\mathbb{Q}_\mathfrak{p})$, unitary if $\pi$ is.  By the discussion of \S\ref{sec:standard-l-factors}, we obtain standard local $L$-factors $L(\pi_\mathfrak{p},s)$.  As shown by Godement--Jacquet \cite[Lemma 10.20]{MR0342495}, the Euler product
\begin{equation*}
  L(\pi,s) := \prod_p L(\pi_p,s)
\end{equation*}
converges for $\Re(s)$ sufficiently large (indeed, for $\Re(s) > 1$ if $\pi$ is unitary, as shown by Jacquet--Shalika \cite[I, Thm 5.3]{MR618323}).  The completed $L$-function $\Lambda(\pi,s) := L(\pi_\infty,s) L(\pi,s)$ satisfies a functional equation $\Lambda(\pi,s) = \eps(\pi,s) \Lambda(\tilde{\pi},1-s)$, where $\tilde{\pi}$ denotes the contragredient representation (for unitary $\pi$, the complex conjugate representation) and $\eps(\pi,s)$ has the form $\eps(\pi,s) = \eps(\pi) C(\pi_{\fin})^{1/2-s}$, with $|\eps(\pi)| = 1$ for unitary $\pi$.

The finite conductor $C(\pi_{\fin}) \in \mathbb{Z}_{\geq 1}$ of $\pi$ is defined to be the product $\prod_p C(\pi_p)$ of the local conductors, as defined in \S\ref{sec:essential-vector}; only finitely many of these local conductors are not equal to one, so the product is really a finite product.  Recall from \S\ref{sec:subc-bounds-assum} that if $\pi$ has archimedean $L$-function parameters $t_1,\dotsc,t_n$, i.e., if $L(\pi_\infty, s) = \prod_{j=1}^n \Gamma_{\mathbb{R}}(s + t_j)$, then its archimedean conductor is defined by $C(\pi_\infty) = \prod_{j=1}^n (3 + |t_j|)$, and its global analytic conductor by $C(\pi) = C(\pi_\infty) C(\pi_{\fin})$.

The completed $L$-function $\Lambda(\pi,s)$ is bounded in vertical strips, so $L(\pi,s)$ is of order one.  The absolute convergence of the Euler product, the functional equation, and the Phragmen--Lindel\"{o}f convexity principle then imply that $L(\pi,s)$ is polynomially bounded in vertical strips (see \cite[Lemma 5.2]{MR2061214}): for $\Re(s) \ll 1$,
\begin{equation}\label{eq:lpi-tfrac12-+-2}
  L(\pi, \tfrac{1}{2} + s) \ll C(\pi)^{\O(1)} (1 + |s|)^{\O(1)},
\end{equation}
provided that $n \geq 2$; when $n = 1$ and $L(\pi,s) = \zeta(s)$, this estimate must be modified to take into account the pole at $s=1$.

We recall the convexity bound:
\begin{lemma}\label{lem:convexity-bound}
  Let $\pi$ be a unitary cuspidal automorphic representation of $\GL_n$ over $\mathbb{Q}$, as above, with $n \geq 1$ fixed.  Then for
  $\sigma \in \mathbb{R}$ with $\sigma \geq -1/2$,
  \begin{equation}\label{eq:lpi-tfrac12-+sigma}
    L(\pi, \tfrac{1}{2} +\sigma)
    \ll
    C(\pi)^{\max(0,(1 - 2 \sigma)/4) + o(1)},
  \end{equation}
  except when $n=1$ and the LHS is of the form $\zeta(s)$ with $s$ in some fixed neighborhood of $1$.
\end{lemma}
\begin{proof}
  The case $\sigma \geq 0$ is recorded by Iwaniec--Kowalski \cite[Thm 5.41]{MR2061214}.  The general case follows from the functional equation and Phragmen--Lindel\"{o}f, as in \cite[\S5, Exercise 3]{MR2061214}.
\end{proof}
\begin{remark}
  We note that \eqref{eq:lpi-tfrac12-+sigma}, unlike \eqref{eq:lpi-tfrac12-+-2}, is uniform in $\pi$, and does not follow formally from the functional equation (without assuming strong hypotheses on $\pi$, such as the Ramanujan conjecture).  The key estimates underlying the proof of \eqref{eq:lpi-tfrac12-+sigma} are due to Iwaniec \cite{MR1067982} when $n=2$ and to Molteni \cite{MR1876443} for general $n$.  Despite the fact that the main results of this paper improve upon the estimate  \eqref{eq:lpi-tfrac12-+sigma} in many cases, we do not in any sense ``reprove'' that estimate: it provides an important input to our arguments.
\end{remark}

The twisted representation $\pi \otimes |.|^s$ has archimedean $L$-function parameters
\begin{equation*}
  t_1+s,\dotsc,t_n+s
\end{equation*}
and the same finite conductor as $\pi$, hence
\begin{equation*}
  C(\pi \otimes |.|^s) = C(\pi_{\fin}) \prod_{j=1}^n (3 + |t_j + s|).
\end{equation*}
Since $C(\pi \otimes |.|^s) \ll (1 + |s|)^{n} C(\pi)$, the convexity bound \eqref{eq:lpi-tfrac12-+sigma} contains the estimate: for $s \in \mathbb{C}$ with  $|\Re(s)| \leq 1/2$,
\begin{equation}\label{eqn:convexity-polynomially-in-s}
  L(\pi, \tfrac{1}{2} +s)
  \ll
  C(\pi)^{(1 - 2 s)/4 + o(1)}
  (1 + |s|)^{\O(1)}.
\end{equation}

\subsection{Eisenstein series}\label{sec:eisenstein-series}
Let $F$ be a number field.

\subsubsection{Definition and convergence}\label{sec:cuhjp7akc2}
We will be concerned primarily with ``minimal'' Eisenstein series coming from the Borel subgroup $Q$.  (We will eventually need to consider ``mirabolic'' or ``Siegel-type'' Eisenstein series, as well, but we treat those separately as they arise.)

For a function $f$ on $[G]_Q$, we denote by $\Eis(f) : [G] \rightarrow \mathbb{C}$ the usual series, provided that it converges absolutely:
\begin{equation*}
  \Eis(f)(g) := \sum_{\gamma \in Q(F) \backslash G(F)} f(\gamma g).
\end{equation*}
When we wish to explicate which groups are being considered, we write more verbosely $\Eis_Q^G$. \index{Eisenstein series!pseudo, $\Eis(f)$}

Standard convergence arguments as in \cite[Theoerm 12.4 and \S12.5]{MR2331343} give that for each $\eps > 0$, there exists $m \geq 0$ so that
\begin{equation}\label{eq:sup_g-in-gmathbba}
  \sup_{g \in G(\mathbb{A})} \|g\|^{-m}
  \sum_{\gamma \in Q(F) \backslash G(F)} |a(\gamma g)|^{(1+\eps) \rho_U}  < \infty,
\end{equation}
with $g = n a(g) k$ an Iwasawa decomposition.  It follows that $\Eis(f)$ converges absolutely if $f \in \mathcal{I}(\chi)$ for some $\chi \in \mathfrak{X}([A])$ with $\Re(\chi) \succ \rho_U$, in which case $\Eis(f) \in \mathcal{T}([G])$.  \index{Eisenstein series!spectral, $\Eis(f)$} The arguments of \cite[\S II.1.7, II.1.10]{MR1361168} show that $\Eis(f)$ converges for $f \in \mathcal{S}([G]_Q)$, in which case $\Eis(f) \in \mathcal{S}([G])$, and the induced map
\begin{equation*}
  \Eis : \mathcal{S}([G]_Q) \rightarrow \mathcal{S}([G])
\end{equation*}
is continuous.  (By the definition of the topologies, we can check the continuity after passing to invariants for some compact open subgroup $J$ of $G(\mathbb{A}_f)$.  We then appeal to \cite[Prop 1.13]{MR1001613}.)

\subsubsection{Meromorphic continuation and Mellin expansion}
For any holomorphic family $\{f[\chi]\}_{\chi \in \mathfrak{X}([A])}$, the assignment $\chi \mapsto \Eis(f[\chi])$, defined for $\Re(\chi) \succ \rho_U$ by the series, extends to a meromorphic vector-valued map $\mathfrak{X}([A]) \rightarrow \mathcal{T}([G])$.  We refer to \cite{MR2402686} for precise statements and proofs, and also to \cite[Thm 11.7]{MR3219530}.

Let $f \in \mathcal{S}([G]_Q)$, with corresponding holomorphic family of Mellin components $\{f[\chi]\}_{\chi \in \mathfrak{X}([A])}$ (necessarily of rapid vertical decay).  Let $\sigma \in \mathfrak{a}^*$ with $\sigma \succ \rho_U$.  By \eqref{eq:sup_g-in-gmathbba}, there exists $m \geq 0$ so that
\begin{equation*}
  \sup_{g \in G(\mathbb{A})} \|g\|^{-m}
  \int_{ (\sigma)  }
  \sum_{\gamma \in Q(F) \backslash G(F)} |f[\chi](\gamma g)| \, d \chi  < \infty,
\end{equation*}
and similarly for archimedean derivatives.  By interchanging integration with summation, it follows that
\begin{equation}\label{eq:eisf-=-int}
  \Eis(f) = \int_{(\sigma) } \Eis(f[\chi]) \, d \chi,
\end{equation}
where the integral converges in $\mathcal{T}([G])$.  In particular, for any $\varphi \in \mathcal{S}([G])$, we have
\begin{equation}\label{eq:int-_g-eisf}
  \int_{[G]} \Eis(f) \varphi
  = \int_{(\sigma)}
  \left( \int_{[G]} \Eis(f[\chi]) \varphi  \right) \, d \chi.
\end{equation}

\subsubsection{Whittaker functions}\label{sec:whittaker-functions}
Let $\psi : N(\mathbb{A}) \rightarrow \U(1)$ be a nondegenerate character, trivial on $N(F)$.  Provided that $\Eis(f)$ converges absolutely and locally uniformly, we have the standard unfolding calculation (see, e.g., \cite[Prop 7.1.3]{MR2683009})
\begin{equation}\label{eq:weisf-psig-=}
  W[\Eis(f),\psi](g) = W[f,\psi](g) := \int_{u \in  N(\mathbb{A})} f(w_G^{-1} u g) \psi^{-1}(u) \, d u.
\end{equation}

\subsubsection{Functional equation}
For $f \in \mathcal{I}(\chi)$, we have the functional equation
\begin{equation}\label{eq:eisenstein-series-functional-equation}
  \Eis(M_w f[\chi]) = \Eis(f[\chi])
\end{equation}
whenever both sides are defined.

\subsubsection{Constant terms}\label{sec:constant-terms}
Let $P$ be a parabolic subgroup of $G$ that contains $A$.  We factor $P = M U_P$, with $M$ the unique Levi factor containing $A$.

The constant term with respect to $P$ of a continuous function $\varphi$ on $[G]$ is the function $\varphi_P : [G]_P \rightarrow \mathbb{C}$ defined by $\varphi_P(g) := \int_{[U_P]} \varphi(u g) \, d u$.

For $f \in \mathcal{S}([G]_Q)$, or for $f \in \mathcal{I}(\chi)$ with $\Re(\chi) \succ \rho_U$, we have by standard unfolding calculations \cite[p101 and \S II.1.7]{MR1361168} that
\begin{equation}\label{eq:int_g-eisf-varphi}
  \int_{[G]} \Eis(f) \varphi = \int_{[G]_Q} f \varphi_Q \quad \text{ for } \varphi \in \mathcal{S}([G])
\end{equation}
and
\begin{equation}\label{eq:eisf_b-=-sum}
  \Eis(f)_Q = \sum_{w \in W} M_w f.
\end{equation}
Such identities extend meromorphically to $f \in \mathcal{I}(\chi)$ for almost all $\chi$.

For general $P$ as above, the intersection $Q_{M} := Q \cap M$ is a Borel subgroup of $M$, and for $f \in \mathcal{S}([G]_Q)$, we have
\begin{equation}\label{eq:eis_bgf_p-=-sum}
  \Eis_Q^G(f)_P = \sum_{w \in W(A,M)} \Eis_{Q_{M}}^{M} (M_w f),
\end{equation}
where $W(A,M)$ denotes the set of all $w \in W$ such that $w^{-1}(\alpha) > 0$ for all positive roots $\alpha$ for $(M,A)$, and where each Eisenstein series appearing on the RHS of \eqref{eq:eis_bgf_p-=-sum} converges absolutely.  This again follows from an unfolding calculation as in \emph{loc.\ cit.}, the main point being, as noted in \emph{loc.\ cit.}, that the set $W(A,M)$ defines a system of representatives for $Q(F) \backslash G(F) / P(F)$.

For a holomorphic family $\{f[\chi]\}_{\chi \in \mathfrak{X}([A])}$, the Eisenstein series $\Eis(f[\chi])$ is holomorphic in $\chi$ if and only its constant term $\Eis(f[\chi])_Q$ with respect to the Borel $Q$ has the same property (see \cite[\S I.4.10 and Remark p167]{MR1361168}).  By \eqref{eq:eisf_b-=-sum}, the latter condition holds provided that $M_w f[\chi]$ is holomorphic in $\chi$ for each $w \in W$.

\subsubsection{Inner product formulas}

We require the following consequence of the spectral theory of Eisenstein series.
\begin{proposition}\label{lem:let-f_1-f_2}
  Let $f_1, f_2 \in \mathcal{S}([G]_Q)$.  Assume that for each $j=1,2$ and $w \in W$, we have $M_w f_j \in \mathcal{S}([G]_Q)$.  Then
  \begin{align*}
    \left\langle \Eis(f_1), \Eis(f_2) \right\rangle
    &= \sum_{w \in W}
      \left\langle f_1, M_w f_2 \right\rangle \\
    &=
      \frac{1}{|W|}
      \left\langle \sum_{w \in W} M_w f_1, \sum_{w \in W } M_w f_2 \right\rangle.
  \end{align*}
\end{proposition}
\begin{proof}
  By \eqref{eq:int_g-eisf-varphi} and \eqref{eq:eis_bgf_p-=-sum}, we have
  \begin{equation*}
    \int_{[G]} \Eis(f_1) \Eis(f_2)
    =
    \int_{[G]_Q} f_1 \Eis(f_2)_Q
    = \sum_{w \in W}
    \int_{[G]_Q} f_1 \cdot  M_w f_2.
  \end{equation*}
  By applying complex conjugation to $f_2$, we obtain the first identity.

  We now study, for $w_1, w_2 \in W$, the integrals
  \begin{equation*}
    \int_{[G]_Q} M_{w_1} f_1 \cdot M_{w_2} f_2,
  \end{equation*}
  which converge in view of our assumption that each $M_w f_j$ lies in $\mathcal{S}([G]_Q)$.  We evaluate $M_{w_2} f_2$ by considering the Mellin expansion $f_2 = \int_{(\sigma)} f_2[\chi] \, d \chi$ (\S\ref{sec:mellin-paley-affine-quotients}), initially for sufficiently dominant $\sigma \in \mathfrak{a}^*$.  Since $f_2[\chi]$ has rapid vertical decay, the double integral defining $M_{w_2} f_2$ then converges absolutely, so the order of integration may be swapped, giving $M_{w_2} f_2 = \int_{(\sigma)} M_{w_2} f_2[\chi] \, d \chi$.  Our assumption that $M_{w_2} f_2$ lies in $\mathcal{S}([G]_Q)$ now implies that the family $M_{w_2} f_2[\chi]$ is holomorphic and has rapid vertical decay, so we may shift the contour to $\sigma = 0$.  Since also $M_{w_1} f_1 \in \mathcal{S}([G]_Q)$, it follows that
  \begin{equation*}
    \int_{g \in [G]_Q} |M_{w_1} f_1(g)| \int_{(0)} \left\lvert M_{w_2} f_2[\chi](g)  \right\rvert \, d \chi \, d g < \infty.
  \end{equation*}
  We may thus exchange integrals once more to arrive at
  \begin{equation*}
    \int_{[G]_Q} M_{w_1} f_1 \cdot M_{w_2} f_2
    = \int_{(0)} \left( \int_{[G]_Q} M_{w_1} f_1(g) \cdot M_{w_2} f_2[\chi](g) \, d g \right) \, d \chi.
  \end{equation*}
  The same identity holds with $f_2$ replaced by its complex conjugate.  Evaluating the resulting integral over $[G]_Q$ using the Iwasawa decomposition, we obtain
  \begin{equation*}
    \left\langle M_{w_1} f_1, M_{w_2} f_2 \right\rangle
    =
    \int_{(0)}
    \left\langle M_{w_1} f_1[\chi], M_{w_2} f_2[\chi] \right\rangle \, d \chi.
  \end{equation*}
  The unitarity of the intertwining operators $M_w$ on the unitary axis $\Re(\chi) = 0$ (see \eqref{eq:langle-f_1-f_2}) and their composition formula (see \eqref{eq:m_w_1w_2-chi-circ}) now give
  \begin{align*}
    \left\langle M_{w_1} f_1, M_{w_2} f_2 \right\rangle
    &=
      \left\langle f_1, M_{w_1}^{-1} M_{w_2} f_2 \right\rangle \\
    &=
      \left\langle f_1, M_{w_1^{-1} w_2} f_2 \right\rangle.
  \end{align*}
  Summing this relation over all $w_1$ and $w_2$ yields the second required identity in the statement of Proposition \ref{lem:let-f_1-f_2}.
\end{proof}
\begin{remark}
  Proposition \ref{lem:let-f_1-f_2} is essentially contained in \cite[\S3.6.6, p384]{MR3468638},  although not quite in the stated form.
\end{remark}

\subsection{Completed Eisenstein series}\label{sec:compl-eisenst-seri}
For simplicity, we now take $F = \mathbb{Q}$ and assume that $G$ arises from a split reductive group over $\mathbb{Z}$, which we continue to denote by $G$.  (The results of this section represent the special case $F = \mathbb{Q}$ and $S = \{\infty\}$, adequate for our immediate applications, of more general results concerning any number field $F$ and large enough finite set of places $S$.)  We assume that $K_p = G(\mathbb{Z}_p)$ for each finite prime $p$.

\subsubsection{Notation and basic definition}\label{sec:notat-basic-defin}
Recall from \S\ref{sec:induc-repr} that, for $s \in \mathfrak{a}_{\mathbb{C}}^*$, we denote by $\mathcal{I}(s) = \Ind_{Q(\mathbb{A})}^{G(\mathbb{A})} (|.|^s)$ the corresponding principal series representation of $G(\mathbb{A})$, which is the restricted tensor product, over places $\mathfrak{p}$ of $\mathbb{Q}$, of the local principal series representations $\mathcal{I}_\mathfrak{p}(s) = \Ind_{Q(\mathbb{Q}_\mathfrak{p})}^{G(\mathbb{Q}_\mathfrak{p})}(|.|_{\mathfrak{p}}^s)$.  For each finite prime $p$ and all $s \in \mathfrak{a}_{\mathbb{C}}^*$, we temporarily write
\begin{equation}\label{eq:holomorphic-families-f-p-of-s}
  f_p[s] \in \mathcal{I}_p(s) \subseteq C^\infty(U(\mathbb{Q}_p) \backslash G(\mathbb{Q}_p))
\end{equation}
for the spherical vector normalized by requiring that $f_p[s](1) = 1$.

As in \S\ref{sec:local-zeta-functions}, for $s \in \mathfrak{a}_{\mathbb{C}}^*$, we define $\zeta(U,s)$, initially for $\Re(s) \succ 0$ by the convergent Euler product \index{zeta functions!$\zeta(U,s)$}
\begin{equation}\label{eq:zeta-U-s-global-over-Z}
  \zeta(U,s) :=
  \prod_p
  \zeta_p(U,s)
  =
  \prod_{\alpha \in R_U^+}
  \zeta(1 + \alpha^\vee(s)),
\end{equation}
then in general via meromorphic continuation.  We also define the $W$-invariant polynomial \index{$W$-invariant polynomial $\mathcal{P}_G(s)$}
\begin{equation*}
  \mathcal{P}_G(s)
  :=
  \prod_{\alpha \in R}
  \alpha^\vee(s)
  (1 - \alpha^\vee(s))
  =
  \prod_{\alpha \in R_U^+}
  \alpha^\vee(s)^2
  (\alpha^\vee(s)^2 - 1).
\end{equation*}
The product $\mathcal{P}_G(s) \zeta(U,s)$ is then entire and of polynomial growth in vertical strips, hence acts as a multiplier on, e.g., $\mathcal{S}^e(U(\mathbb{R}) \backslash G(\mathbb{R}))$.

Let $f_\infty \in \mathcal{S}^e(U(\mathbb{R}) \backslash G(\mathbb{R}))$, with associated Mellin components $f_\infty[s] \in \mathcal{I}_\infty(s)$.

For $s \in \mathfrak{a}_{\mathbb{C}}^*$, define $\Phi[f_\infty][s] \in \mathcal{I}(s)$ by the formula \index{Eisenstein series!completed inducing data $\Phi[f]$}
\begin{equation}\label{eq:phis-:=-qn}
  \Phi[f_\infty][s] = \mathcal{P}_G(s) \zeta(U,s) \otimes_{\mathfrak{p}} f_\mathfrak{p}[s].
\end{equation}
This defines a holomorphic family of vectors of rapid vertical decay.  Per the discussion of \S\ref{sec:mellin-paley-affine-quotients}, this family arises by taking the Mellin components of some element $\Phi[f_\infty]$ of the Schwartz space $\mathcal{S}^e([G]_Q)$.  We obtain a $G(\mathbb{R})$-equivariant linear map
\begin{equation*}
  \Phi : \mathcal{S}^e(U(\mathbb{R}) \backslash G(\mathbb{R})) \rightarrow \mathcal{S}^e([G]_Q)
\end{equation*}
\begin{equation*}
  f_\infty \mapsto \Phi[f_\infty].
\end{equation*}
This map may be used to define the ``completed'' pseudo Eisenstein series
\begin{align*}
  \Eis[\Phi[f_\infty]]
  &= \int_{(\sigma)} \Eis[\Phi[f_\infty][s]] \, d \mu_{[A]}(s)
    \quad \quad (\sigma \succ \rho_U) \\
  &= \int_{(\sigma)} \mathcal{P}_G(s) \zeta(U,s)
    \Eis[\otimes_{\mathfrak{p}} f_\mathfrak{p}[s]] \, d \mu_{[A]}(s).
\end{align*}
In what follows, we drop the subscript $\infty$, writing simply $f$ for an element of $\mathcal{S}^e(U(\mathbb{R}) \backslash G(\mathbb{R}))$ and $\Phi[f]$ for the associated element of $\mathcal{S}^e([G]_Q)$.

We note that for $a \in A(\mathbb{R}) \hookrightarrow A(\mathbb{A})$, we have
\begin{equation*}
  \Phi[L(a) f] = L(a) \Phi[f].
\end{equation*}
Indeed, both sides have Mellin components
\begin{equation*}
  |a|^s \mathcal{P}_G(s) \zeta(U,s) f[s] \otimes (\otimes_p f_p[s]).
\end{equation*}

\begin{remark}\label{rmk:one-could-altern}
  Given $f \in \mathcal{S}^e(U(\mathbb{R}) \backslash G(\mathbb{R}))$, one could alternatively define $\Phi[f]$ as follows.  Set
  \begin{equation*}
    f^\sharp [s] := \mathcal{P}_G(s) \zeta(U,s) f[s].
  \end{equation*}
  This defines a holomorphic family of rapid vertical decay, corresponding to some element $f^\sharp \in \mathcal{S}^e(U(\mathbb{R}) \backslash G(\mathbb{R}))$.  Note here that $\zeta(U,s)$, defined initially by an Euler product over the primes, is being employed as a Mellin multiplier at the archimedean place.  At a finite place $p$, the holomorphic family \eqref{eq:holomorphic-families-f-p-of-s} corresponds to the element $f_p \in \mathcal{S}(U(\mathbb{Q}_p) \backslash G(\mathbb{Q}_p))$ given by the characteristic function of $U(\mathbb{Q}_p) K_p$.  A simple application of Mellin inversion gives
  \begin{equation*}
    \Phi[f](g) = \sum_{c \in A(\mathbb{Q})} f^\sharp (c g) \prod_p f_p(c g).
  \end{equation*}
  Indeed, the convergent integrals
  \begin{equation*}
    \int_{[A]} |a|^{-\rho_U - s}\Phi[f](a g)  \, d a
  \end{equation*}
  and
  \begin{align*}
    &\int_{[A]} |a|^{-\rho_U - s}\sum_{c \in A(\mathbb{Q})} f^\sharp (c a_\infty  g_\infty ) \prod_p f_p (c a_p g_p) \, d a \\
    &\quad =
      \int_{A(\mathbb{A}) } |a|^{-\rho_U - s} f^\sharp (a_\infty  g_\infty ) \prod_p f_p (a_p g_p) \, d a
  \end{align*}
  both evaluate to
  \begin{equation*}
    f^\sharp [s] (g_\infty ) \prod_p f_p [s] (g_p).
  \end{equation*}
  Here we note that the sum over $c$ is in fact a finite sum; for example, if $g \in G(\mathbb{R})$, then $\Phi[f](g) = \sum_{c \in A(\mathbb{Z})} f^\sharp (c g)$.  One derives different series representations, analogous to \eqref{eq:eg-=-sum}, by ``diverting'' the Mellin multiplier $\zeta(U,s)$ from the archimedean place to the finite place.  We explain this passage in greater detail, together with some variants, starting in \S\ref{sec:interl-basic-vect}.

  In the above optic, the associated Eisenstein series is given by
  \begin{equation}\label{eq:eisphif=-sum-_gamma}
    \Eis[\Phi[f]](g) = \sum_{\gamma \in Q(\mathbb{Q}) \backslash G(\mathbb{Q}) } \Phi[f](\gamma g)
    = \sum_{\gamma \in U(\mathbb{Q}) \backslash G(\mathbb{Q})} f^\sharp (\gamma g_\infty ) \prod_p f_p (\gamma g_p).
  \end{equation}
\end{remark}
\begin{example}\label{exa:example-sl2-completed-eisenstein-series}
  Take  $G = \SL_2$, and assume that $U$ is lower-triangular.  Then
  \begin{itemize}
  \item $U(\mathbb{R}) \backslash G(\mathbb{R})$ identifies with the punctured plane $\mathbb{R}^2 - \{0\}$ via the map associating to a matrix its top row,
  \item elements $f \in \mathcal{S}^e(U(\mathbb{R}) \backslash G(\mathbb{R}))$ identify with smooth even Schwartz functions $f : \mathbb{R}^2 \rightarrow \mathbb{C}$ that, together with all derivatives, decay rapidly near the origin, and
  \item the Eisenstein series $\Eis[\Phi[f]]$, defined above, identifies with the function $E : \SL_2(\mathbb{Z}) \backslash \SL_2(\mathbb{R}) \rightarrow \mathbb{C}$ given by
    \begin{equation}\label{eq:eg-=-sum}
      E(g) = \sideset{}{^*}\sum_{v \in \mathbb{Z}^2} f^\sharp( v g),
    \end{equation}
    where $f^\sharp$ is as in Remark \ref{rmk:one-could-altern} and $\sum^*$ denotes a sum over primitive vectors.
  \end{itemize}
  Indeed, consider the evaluation of \eqref{eq:eisphif=-sum-_gamma} at $g \in G(\mathbb{R})$:
  \begin{equation*}
    E(g) = \Eis[\Phi[f]](g) = \sum_{\gamma \in U(\mathbb{Q}) \backslash G(\mathbb{Q})} f^\sharp (\gamma g) \prod_p f_p(\gamma).
  \end{equation*}
  Identifying $U(\mathbb{Q}) \backslash G(\mathbb{Q})$ with $\mathbb{Q}^2 - \{0\}$ as above, the product $\prod_p f_p(\gamma)$ detects when, for each prime $p$, the image of $\gamma$ in $\mathbb{Q}_p^2 - \{0\}$ defines a primitive element of $\mathbb{Z}_p^2 - \{0\}$.  It is equivalent to ask that $\gamma$ itself defines a primitive element of $\mathbb{Z}^2$, whence \eqref{eq:eg-=-sum}.

  Alternatively, defining $f^{\flat} \in \mathcal{S}^e(U(\mathbb{R}) \backslash G(\mathbb{R}))$ by requiring that
  \begin{equation*}
    f^{\flat}[s] = \mathcal{P}_G(s) f[s],
    \quad
    f^\sharp [s] = \zeta(U,s) f^{\flat}[s],
  \end{equation*}
  we have
  \begin{equation}\label{eq:eg-=-sum-1}
    E(g) = \sum_{v \in \mathbb{Z}^2 - \{0\}} f^{\flat}(v g).
  \end{equation}
  Indeed, let $\alpha \in \Delta_U$ denote the unique simple root, and $\alpha^\vee$ the corresponding coroot.  The latter induces an identification of $\mathfrak{a}_{\mathbb{C}}^*$ with $\mathbb{C}$.  For $s \in \mathbb{C} \cong \mathfrak{a}_{\mathbb{C}}^*$, we have $\zeta(U,s) = \zeta(1 + s)$, while the character $|.|^s$ of $A(\mathbb{R})$ is given by
  \begin{equation}\label{eq:diagt-1ts-=}
    |\diag(t, 1/t)|^s = |t|^{-s}.
  \end{equation}
  By definition, $f[s](g) = \int_{a \in A(\mathbb{R})} |a|^{-\rho_U - s} f(a g) \, d a$ for $a \in \mathfrak{a}_{\mathbb{C}}^*$.  Identifying both $f[s]$ and $f$ with functions on $\mathbb{R}^2 - \{0\}$ as above and invoking \eqref{eq:diagt-1ts-=}, we see that this definition may be rewritten, for $s \in \mathbb{C} \cong \mathfrak{a}_{\mathbb{C}}^*$,
  \begin{equation*}
    f[s](v) = \int_{t \in \mathbb{R}^\times } |t|^{1+s} f(t v) \, \frac{d t}{|t|}
  \end{equation*}
  (for suitable normalization of Haar measure).  The same relation holds for $f^\sharp$ and $f^{\flat}$.  It follows readily that
  \begin{equation*}
    f^\sharp (v) = \sum_{c \in \mathbb{Z}_{\geq 1} } f^{\flat} (c v).
  \end{equation*}
  The claim \eqref{eq:eg-=-sum-1} then follows from the fact that every element of $\mathbb{Z} ^2 - \{0\}$ may be written uniquely as $c v$, where $c \in \mathbb{Z}_{\geq 1}$ and $v$ is primitive.

  One can check that $f^\flat$ is the image of $f$ under an explicit finite-order differential operator, coming from the scaling action of $\mathbb{R}^\times$ on $\mathbb{R}^2$.

  We note that the proofs of the basic properties of the associations $f \mapsto \Phi[f]$ and $f \mapsto \Eis[\Phi[f]]$ relevant for this paper can ultimately be reduced to analysis of this example (see \cite[\S3.5]{MR3468638}).
\end{example}

\subsubsection{Equivariance under Fourier transforms}
In what follows, we let $\psi$ be a nondegenerate unitary character of $N(\mathbb{A})$, trivial on $N(\mathbb{Q})$, that is unramified in the sense of \S\ref{sec:stand-nond-char}.  This condition determines $\psi$ ``up to signs.''  Recall the notation $\mathcal{S}^e(U(\mathbb{R}) \backslash G(\mathbb{R}))^w$ from \S\ref{sec:fourier-transforms}.  We define this notation with respect to the restriction $\psi_{\infty}$ of $\psi$ to $N(\mathbb{R})$, thus $\mathcal{S}^e(U(\mathbb{R}) \backslash G(\mathbb{R}))^w$ is the subspace of $\mathcal{S}^e(U(\mathbb{R}) \backslash G(\mathbb{R}))$  consisting of $\mathcal{F}_{w,\psi_\infty}$-invariant elements.

\begin{lemma}\label{lem:standard2:let-f_infty-in}
  Let $w \in W$ and $f \in \mathcal{S}^e(U(\mathbb{R}) \backslash G(\mathbb{R}))^w$.  Then
  \begin{equation*}
    M_w \Phi[f] =
    \Phi[f].
  \end{equation*}
\end{lemma}
\begin{proof}
  Abbreviate $\Phi := \Phi[f]$.  It admits the rapidly-convergent Mellin expansion
  \begin{equation}\label{eq:phi-=-int_sigma}
    \Phi = \int_{(\sigma)} \Phi[s] \, d s
  \end{equation}
  for each $\sigma \in \mathfrak{a}^*$.  Choose $\sigma$ such that $\sigma \succ \rho_{U}$.  Then for each $g \in G(\mathbb{A})$, the double integral
  \begin{equation*}
    M_w \Phi(g) = \int_{u \in
      (w U w^{-1} \cap U \backslash  U)(\mathbb{A})} \int_{(\sigma)}  \Phi[s](w^{-1} u g) \, d s \, d u
  \end{equation*}
  converges absolutely.  We may thus swap the integrals, giving
  \begin{equation*}
    M_w \Phi(g) = \int_{(\sigma) } M_w \Phi[s](g) \, d s.
  \end{equation*}
  If we can show that $M_w \Phi[s] = \Phi[w s]$ then it will follow from \eqref{eq:phi-=-int_sigma} that $M_w \Phi = \Phi$.  (Note that we have used crucially here that $\Phi \in \mathcal{S}([G]_{Q})$.)

  It remains to verify that $M_w \Phi[s] = \Phi[w s]$.  Recall the definition \eqref{eq:phis-:=-qn} of $\Phi[s]$.  Using the notation $\otimes^\sharp$ as in \S\ref{sec:fact-intertw-oper} to denote a tensor product regularized by $\zeta(U,s)$ and its local factors, we may rewrite that definition as
  \begin{equation*}
    \Phi[s] = \mathcal{P}_G(s) f[s] \otimes (\otimes_{p} ^\sharp f_p^\sharp [s]),
  \end{equation*}
  where $f_\mathfrak{p} ^\sharp [s]$ denotes the spherical element taking the value $\zeta_\mathfrak{p}(U,s)$ at the identity.  We have noted in \eqref{eq:normalized-F-w-preserves-basic-vector} that $\mathcal{F}_{w,\psi_p} f_p^\sharp[s] = f_p^\sharp [w s]$ (using here that $\psi_p$ is unramified).  The identity \eqref{eq:m_w-f-=-2} thus gives
  \begin{equation*}
    M_w \Phi[s] = \mathcal{P}_G(s) \mathcal{F}_{w,\psi_\infty} f[s] \otimes (\otimes_p^\sharp f_p^\sharp [w s]).
  \end{equation*}
  We have noted that the polynomial $\mathcal{P}_G(s)$ is $W$-invariant.  Our hypothesis is that $\mathcal{F}_{w,\psi_\infty} f[s] = f[ w s]$ for all $s \in \mathfrak{a}_{\mathbb{C}}^*$ (see Lemma \ref{lem:standard2:let-w-in}).  Therefore $M_w \Phi[s] = \Phi[w s]$, as required.
\end{proof}

\begin{remark}\label{rmk:one-could-likely}
  One could likely formulate and establish a more natural version of Lemma \ref{lem:standard2:let-f_infty-in}, as follows.  Let $X$ denote the affine closure of the quasi-affine variety $U \backslash G$.  Define the extended Schwartz space $\mathcal{S}(X(\mathbb{R}))$ as in \cite[Definition 5.5]{2021arXiv210310261G} or \cite[\S7.1]{MR1694894}; it should be a space intermediate between $\mathcal{S}(U(\mathbb{R}) \backslash G(\mathbb{R}))$ and $C^\infty(U(\mathbb{R}) \backslash G(\mathbb{R}))$.  The assignment $\Phi$ should extend to $\mathcal{S}(X(\mathbb{R}))$, hence in particular to its subspace $\mathcal{S}^e(X(\mathbb{R}))$.  The group $W$ should act on $\mathcal{S}^e(X(\mathbb{R}))$ (and $\mathcal{S}(X(\mathbb{R}))$) via the Fourier transforms $\mathcal{F}_{w,\psi_\infty}$, and we should have
  \begin{equation}\label{eq:Phi-intertwines-M-w-F-w}
    M_w \Phi[f] = \Phi[\mathcal{F}_{w,\psi_\infty} f].
  \end{equation}
  We were unable to extract such expectations from the literature on such extended Schwartz spaces, although the technique needed to formulate and establish identities such as \eqref{eq:Phi-intertwines-M-w-F-w} is certainly available (essentially already in \cite[\S3.5]{MR3468638}), and we suspect that such identities are well known to experts.  In any event, the simpler but less natural implication given by Lemma \ref{lem:standard2:let-f_infty-in} suffices for our applications.
\end{remark}

\subsubsection{Matrix coefficient bounds}\label{sec:matr-coeff-bounds}
Let $\mathcal{H} := \mathcal{H}_G$ denote the spherical Hecke algebra for $G(\mathbb{Z})$, i.e., the restricted product over all primes $p$ of the spherical Hecke algebras for $(G(\mathbb{Q}_p), G(\mathbb{Z}_p))$.

By definition, $\mathcal{H}$ consists of the smooth compactly-supported complex measures on $G(\mathbb{A}_f) = \prod ' G(\mathbb{Q}_p)$ that are bi-invariant under $\prod_p G(\mathbb{Z}_p)$.  Using the Haar measure on $G(\mathbb{A}_f)$ that assigns volume one to $\prod_{p} G(\mathbb{Z}_p)$, we may identify elements of $\mathcal{H}$ with functions on $G(\mathbb{A}_f)$.

We may speak of an element $t \in \mathcal{H}$ being nonnegative; this says that it defines a positive measure, or equivalently, that the corresponding function is nonnegative.

By the Satake isomorphism, characters (i.e., algebra homomorphisms) $\lambda : \mathcal{H} \rightarrow \mathbb{C}$ are in bijection with collections $A = (A_p)_p$, indexed by the primes $p$, of semisimple conjugacy classes $A_p$ in the dual group ${}^L G(\mathbb{C})$.  (In the case $G = \GL_n$ ultimately of interest to us, we have ${}^L G = \GL_n$, so such classes may be identified further with $n$-tuples of nonzero complex numbers.)

For $s \in \mathfrak{a}_{\mathbb{C}}^*$, let $\lambda_s$ denote the character of $\mathcal{H}$ attached to the principal series representation $\mathcal{I}(s)$.

The representation $\mathcal{I}(s)$ is tempered when $\Re(s) = 0$.  It follows from the results of \cite{MR946351} (or more elementarily) that if $t \in \mathcal{H}$ is nonnegative and $\Re(s) = 0$,  then
\begin{equation}\label{eq:lambd-leq-lambd}
  |\lambda_s(t)| \leq \lambda_0(t).
\end{equation}

\begin{lemma}\label{lem:standard2:matrix-coeff-bounds}
  Let $G$ be as above.  There exists $C_0 \geq 0$ and $d \in \mathbb{Z}_{\geq 0}$ with the following property.  Let $f \in \mathcal{S}^e(U(\mathbb{R}) \backslash G(\mathbb{R}))^W$.  For $Y \in A(\mathbb{R})^0$, set
  \begin{equation*}
    \Psi[Y] := \Eis[L(Y) \Phi[f]]  = \Eis[\Phi[L(Y) f]] \in \mathcal{S}([G]).
  \end{equation*}
  Let $t \in \mathcal{H}_G$ be nonnegative.  Let $g \in G(\mathbb{R})$ be arbitrary.  Write $\left\langle , \right\rangle$ for the Petersson inner product on $[G]$.  Then
  \begin{equation}\label{eq:leftlangle-g-t}
    |\left\langle g t \Psi[Y], \Psi[Y]  \right\rangle| \leq C_0 \lambda_0(t) \mathcal{S}_{d,0}(f)^2.
  \end{equation}
\end{lemma}
\begin{proof}
  We treat first the case $(Y,t,g) = (1,1,1)$.

  By Lemma \ref{lem:standard2:let-f_infty-in}, we have $M_w \Phi[f] = \Phi[f]$ for each $w \in W$.  In particular, $M_w \Phi[f] \in \mathcal{S}([G]_Q)$, so by Proposition \ref{lem:let-f_1-f_2}, we have
  \begin{equation*}
    \left\langle \Psi[1], \Psi[1]  \right\rangle
    = |W| \left\langle \Phi[f], \Phi[f] \right\rangle.
  \end{equation*}
  (The inner product on the LHS is in $L^2([G])$, that on the RHS in $L^2([G]_Q)$.)

  We may expand
  \begin{align*}
    \left\langle \Phi[f], \Phi[f] \right\rangle
    &=
      \int_{
      \substack{
      s \in \mathfrak{a}_{\mathbb{C}}^* :  \\
    \Re(s) = 0
    }
    }
    \left\langle \Phi[f][s], \Phi[f][s] \right\rangle \, d \mu_{[A]}(s) \\
    &=
      \int_{
      \substack{
      s \in \mathfrak{a}_{\mathbb{C}}^* :  \\
    \Re(s) = 0
    }
    }
    |\mathcal{P}_G(s) \zeta(U,s)|^2
    \left\langle f[s], f[s] \right\rangle \, d \mu_{[A]}(s).
  \end{align*}
  We choose $d \in \mathbb{Z}_{\geq 0}$ and $C_1 \geq 0$ so that the function $E_d$ on $\mathfrak{a}_{\mathbb{C}}^*$, as given by \eqref{eq:e_ds-:=-sum}, satisfies
  \begin{equation*}
    |\mathcal{P}_G(s) \zeta(U,s)|^2 \leq C_1 E_d(s)
  \end{equation*}
  whenever $\Re(s) = 0$.  This is possible in view of the convexity bound for the Riemann zeta function (or the well-known inequality $|t \zeta(1+it)| \leq C_2 \log(2 + |t|)$).  Using the strong approximation isomorphism $A(\mathbb{Q}) \backslash A(\mathbb{A}) / \prod_p A(\mathbb{Z}_p) \cong  A(\mathbb{Z})  \backslash A(\mathbb{R})$ for the split torus $A$ and our normalization of Haar measures, we verify that
  \begin{equation*}
    d \mu_{[A]}(s) = |A(\mathbb{Z})|\,  d \mu_{A(\mathbb{R})}(s).
  \end{equation*}
  Thus
  \begin{equation*}
    \left\langle \Phi[f], \Phi[f] \right\rangle \leq |A(\mathbb{Z})| C_1 J,
  \end{equation*}
  where
  \begin{equation*}
    J := \int_{
      \substack{
        s \in \mathfrak{a}_\mathbb{C}^* :  \\
        \Re(s)  = 0
      }
    }
    E_d(s) \left\langle f[s], f[s] \right\rangle d \mu_{A(\mathbb{R})}(s).
  \end{equation*}
  By \eqref{eq:sobolev-parseval-A-U-G}, we have $J = \mathcal{S}_{d,0}(f)^2$.  Thus the required bound \eqref{eq:leftlangle-g-t} holds in this case with $C_0 := |W| \cdot |A(\mathbb{Z})| \cdot  C_1$.

  In the general case, analogous arguments give the decomposition
  \begin{align*}
    &\left\langle g t \Psi[Y], \Psi[Y] \right\rangle \\
    &=
      |A(\mathbb{Z})|
      \sum_{w \in W}
      \int_{
      \substack{
      s \in \mathfrak{a}_{\mathbb{C}}^* :  \\
    \Re(s) = 0
    }
    }
    |\mathcal{P}_G(s) \zeta(U,s)|^2
    \lambda_s(t)
    Y^{s + \overline{w s}}
    \left\langle R(g) f[s], f[s] \right\rangle
    \, d \mu_{A(\mathbb{R})}(s).
  \end{align*}
  Using the inequality \eqref{eq:lambd-leq-lambd}, the identity $\left\lvert Y^{s+ \overline{w s}} \right\rvert = 1$ and Cauchy--Schwarz, we obtain now
  \begin{equation*}
    \left\lvert \langle g t \Psi[Y], \Psi[Y] \rangle \right\rvert \leq C_0 \lambda_0(t) \mathcal{S}_{d,0}(f) \mathcal{S}_{d,0}(R(g) f),
  \end{equation*}
  with the same $C_0$ as before.  We conclude by noting that $\mathcal{S}_{d,0}(R(g) f) = \mathcal{S}_{d,0}(f)$.
\end{proof}

\subsubsection{Crude growth bounds}

\begin{lemma}\label{lem:standard2:crude-growth-bounds-completed-eisenstein-series}
  For each $C_0 \geq 0$ and $x \in \mathfrak{U}(G(\mathbb{R}))$, there exists $C_1 \geq 0$ so that for each $C_2 \in \mathbb{R}$, there is a compact $\mathcal{D} \subseteq \mathfrak{a}^*$  and $\ell \in \mathbb{Z}_{\geq 0}$ with the following property.  For all $T \geq 1$, $f \in \mathcal{S}^e(U(\mathbb{R}) \backslash G(\mathbb{R}))$, $Y \in A(\mathbb{R})^0$ with each component $Y_j$ lying in $[T^{-1/2}, T^{1/2}]$, and $g \in [G]$, we have
  \begin{equation*}
    \|g\|_{[G]}^{C_0}
    \left\lvert     \det(g Y) \right\rvert^{C_2}
    \left\lvert \Eis[\Phi[ R(x) L(Y) f]](g) \right\rvert
    \leq C_1 T^{C_1} \nu_{\mathcal{D},\ell,T}(f).
  \end{equation*}
\end{lemma}
\begin{remark}
  The exponent $1/2$ in $[T^{-1/2}, T^{1/2}]$ is not important -- the same conclusion holds for any exponent.  The key feature here is that $C_2$ may be taken arbitrarily large, and either positive or negative, without increasing $C_1$.  The practical consequence is that, e.g., $\Eis[\Phi[L(Y) f]](g)$ is concentrated on $\det(g Y) \approx 1$.
\end{remark}
\begin{proof}
  By linearity, we my assume that $x = x_1 \dotsb x_k$ with each $x_j \in \mathcal{B}(G)$.  Then, using the inequality
  \begin{equation*}
    \nu_{\mathcal{D},\ell,T}(R(x) f) \leq
    T^k
    \nu_{\mathcal{D},\ell + k, T}(f),
  \end{equation*}
  we may and shall reduce to the case $x = 1$.

  The proof provides a convenient opportunity to illustrate the use of asymptotic notation, following the conventions of \S\ref{sec:asympt-notat-term}.

  By transfer, it is equivalent to prove the same statement, but with the quantification over the group $G$ and the scalars $C_0$, $C_1$ and $\ell$ restricted to fixed quantities.  The conclusion is filtered with respect to $C_1$ and $\ell$, so by idealization, we may commute them with the remaining quantifiers.  Our task is then to show the following:
  \begin{itemize}
  \item \emph{Let $G$ and $C_0 \geq 0$ be fixed.  Let $C_2 \in \mathbb{R}$, $T \geq 1$, $Y \in A(\mathbb{R})^0$ with each component $Y_j$ lying in $[T^{-1/2}, T^{1/2}]$, and $g \in [G]$.  Let $f \in \mathcal{S}^e(U(\mathbb{R}) \backslash G(\mathbb{R}))$.  There is a compact $\mathcal{D} \subseteq \mathfrak{a}^*$ and a fixed $\ell \in \mathbb{Z}_{\geq 0}$ so that}
    \begin{equation}\label{eq:g_gc_0-leftlv-detg}
      \|g\|_{[G]}^{C_0}
      \left\lvert     \det(g Y) \right\rvert^{C_2}
      \left\lvert \Eis[\Phi[ L(Y) f]](g) \right\rvert
      \ll T^{\O(1)} \nu_{\mathcal{D},\ell,T}(f).
    \end{equation}
  \end{itemize}
  The practical consequence of this formalism is that ``fixed'' means ``depending only upon $C_0$, and not upon $C_2$, $T$, etc.''

  Set $J := \prod_p G(\mathbb{Z}_p)$; it is an open subgroup of $G(\mathbb{A}_f)$ under which $\Phi[f][s]$ is invariant.  The continuity of $\Eis : \mathcal{S}([G]_Q) \rightarrow \mathcal{S}([G])$ implies that there is a fixed continuous seminorm $\nu$ on $\mathcal{S}([G]_Q)^J$ so that for all $\phi \in \mathcal{S}([G]_Q)^J$, we have
  \begin{equation}\label{eq:sup-_g-in}
    \sup_{g \in [G]} \|g\| ^{C_0} _{[G]} | \Eis(\phi) (g)| \leq \nu(\phi).
  \end{equation}
  Assume now that $\phi \in \mathcal{S}^e([G]_Q)^J$.  As usual, write $\phi[s] \in \mathcal{I}(s)$ for its Mellin components.  By the definition of the topology on $\mathcal{S}([G]_Q)^J$, we may find a fixed compact set $\mathcal{D}_0 \subseteq \mathfrak{a}^*$, fixed $y \in \mathfrak{U}(G(\mathbb{R}))$ and fixed $\ell_1 \geq 0$ so that
  \begin{equation}\label{eq:nuphi-ll-sup}
    \nu(\phi) \ll \sup_{
      \substack{
        s \in \mathfrak{a}_{\mathbb{C}}^* :
        \\
        \Re(s) \in \mathcal{D}_0
      }
    }
    \langle s \rangle^{\ell_1} \|R(y) \phi[s]\|,
  \end{equation}
  where we abbreviate $\langle s \rangle := 1 + |s|$.

  We will apply the estimate \ref{eq:sup-_g-in} to the following ``twist'' $\phi$ of $\Phi[L(Y) f]$:
  \begin{equation*}
    \phi(g) :=
    |\det(g Y)|^{C_2} \Phi[L(Y)f](g),
  \end{equation*}
  for which
  \begin{equation*}
    |\Eis(\phi)(g)| = |\det(g Y)|^{C_2} |\Eis[\Phi[L(Y) f]](g)|.
  \end{equation*}
  In view of \eqref{eq:sup-_g-in} and \eqref{eq:nuphi-ll-sup}, our task \eqref{eq:g_gc_0-leftlv-detg} reduces to verifying that for $\Re(s) \in \mathcal{D}_0$, we have
  \begin{equation}\label{eq:langle-s-rangleell_1}
    \langle s \rangle^{\ell_1} \|R(y) \phi[s]\| \ll T^{\O(1)} \nu_{\mathcal{D}, \ell, T}(f)
  \end{equation}
  for some compact $\mathcal{D} \subseteq \mathfrak{a}^*$ and some fixed $\ell \in \mathbb{Z}_{\geq 0}$.  The Mellin components of this twist are given by shifts of the original Mellin components.  Indeed, write $e \in \mathfrak{a}^*$ for the element corresponding to the diagonal matrix $(1,\dotsc,1)$, so that $|a|^{C_2 e} = |\det a|^{C_2}$.  Then
  \begin{align*}
    \phi[s](g)
    &=
      \int_{[A]} |a|^{-\rho_U - s} |\det(g Y)|^{C_2} \Phi[L(Y) f](a g) \, d a \\
    &=
      |\det(g)|^{C_2} \int_{[A]} |a|^{-\rho_U - s} Y^{C_2 e - \rho_U} \Phi[f](Y a g) \, d a.
  \end{align*}
  By the change of variables $a \mapsto Y^{-1} a$, it follows that
  \begin{equation}\label{eq:phis-=-left}
    \phi[s](g) = |\det(g)|^{C_2} Y^s \Phi[f][s - C_2 e].
  \end{equation}
  We have $\mathcal{P}_G(s - C_2 e) = \mathcal{P}_G(s)$ and $\zeta(U,s - C_2 e) = \zeta(U,s)$, so
  \begin{equation*}
    \Phi[f][s - C_2 e]
    = \mathcal{P}_G(s) \zeta(U,s) f[s - C_2 e].
  \end{equation*}
  For $\Re(s) \in \mathcal{D}_0$, we have $\mathcal{P}_G(s) \zeta(U,s) \ll \langle s \rangle^{\ell_2}$ for some fixed $\ell_2 \in \mathbb{Z}_{\geq 0}$.  We note finally that the factor $|\det(g)|^{C_2}$ has no effect on the norms $\|.\|$ obtained from $L^2(K)$.  Our task \eqref{eq:langle-s-rangleell_1} thereby reduces to verifying that for $\Re(s) \in \mathcal{D}_0$, we have
  \begin{equation}\label{eq:langle-s-rangleell_1-1}
    \langle s \rangle^{\ell_1 + \ell_2} Y^s \|R(y) f[s-C_2 e]\| \ll T^{\O(1)} \nu_{\mathcal{D},\ell,T}(f)
  \end{equation}
  for some compact $\mathcal{D} \subseteq \mathfrak{a}^*$ and some fixed $\ell \in \mathbb{Z}_{\geq 0}$.

  Let $\ell \in \mathbb{Z}_{\geq 0}$ be a fixed quantity larger than the maximum of $\ell_1 + \ell_2$ and the degree of $y$.  Define $\mathcal{D}$ to be the translated compact set $\mathcal{D}_0 - C_2 e$.  Then for $\Re(s) \in \mathcal{D}_0$, we have $\Re(s - C_2 e) \in \mathcal{D}$, and so from the definition of $\nu_{\mathcal{D},\ell,T}$,
  \begin{equation}\label{eq:ry-fs-c_2}
    \|R(y) f[s - C_2 e]\| \ll \langle s - C_2 e \rangle^{-\ell} T^{\langle \rho, s - C_2 e \rangle + \ell} \nu_{\mathcal{D},\ell,T}(f).
  \end{equation}
  Since $\Re(s) \in \mathcal{D}_0$, we have $\langle s - C_2 e \rangle \gg \langle s \rangle$.  Since $\langle \rho, s - C_2 e \rangle = \langle \rho, s \rangle$, we have $T^{\langle \rho, s - C_2 e \rangle + \ell} = T^{\O(1)}$.  By our hypotheses on $Y$, we have $Y^s = T^{\O(1)}$.   The required estimate \eqref{eq:langle-s-rangleell_1-1} follows.
\end{proof}

\subsection{Siegel domains}\label{sec:siegel-domains}
We take here for $G$ a general linear group over $\mathbb{Z}$, with $K_p = G(\mathbb{Z}_p)$.

For $c > 0$, we write $A(\mathbb{R})^0_{\geq c}$ for the set of all $a \in A(\mathbb{R})$ satisfying $|a|^{\alpha} \geq c$ for all $\alpha \in \Delta_N$, i.e., $a_i/a_{i+1} \geq c$ for $1 \leq i < \rank(G)$.

\begin{proposition}\label{prop:siegel-upper}
  There exists a compact subset $\Omega$ of $G(\mathbb{R})$ so that for each nonnegative integrable function $\phi$ on $[G]$,
  \begin{equation*}
    \int_{[G]} \phi \leq \int_{t \in A(\mathbb{R})_{\geq 1}^0} \int_{g \in \Omega \times \prod_p K_p} \phi(t g) \, d g \, \frac{d t}{\delta_N(t)}.
  \end{equation*}
\end{proposition}
\begin{proposition}\label{prop:siegel-lower}
  For each compact subset $\Omega$ of $G(\mathbb{R})$, there exists $C_0 \geq 0$ so that for each $t \in A(\mathbb{R})^0_{\geq 1}$ and each nonnegative integrable function $\phi$ on $[G]$,
  \begin{equation*}
    \int_{\Omega \times \prod_p K_p} \phi(t g) \, d g
    \leq C_0 \delta_N(t) \int_{[G]} \phi.
  \end{equation*}
\end{proposition}
\begin{proposition}\label{prop:siegel-lower-integrated}
  For each compact subset $\Omega$ of $G(\mathbb{R})$, there exists $C_1 \geq 0$ so that for each nonnegative integrable function $\phi$ on $[G]$,
  \begin{equation*}
    \int_{t \in A(\mathbb{R})_{\geq 1}^0} \int_{\Omega \times \prod_p K_p} \phi(t g) \, d g
    \, \frac{d t}{\delta_N(t)}
    \leq C_1 \int_{[G]} \phi.
  \end{equation*}
\end{proposition}
Each proposition follows readily from standard results on Siegel domains.  We could not quickly locate a suitable reference, so we record the short proofs.  First of all, it is known that $G(\mathbb{Q}) G(\mathbb{R}) \prod_p K_p = G(\mathbb{A})$ \cite[Prop 2.2]{MR202718} and $G(\mathbb{Z}) = G(\mathbb{Q}) \cap \prod_p K_p$, so we can reduce to verifying the corresponding assertions concerning $\phi$ on $G(\mathbb{Z}) \backslash G(\mathbb{R})$, where we integrate now over $\Omega$ rather than $\Omega \times \prod_p K_p$.  For $b, c > 0$, let $\Omega_N(b)$ denote the set of all $n \in N(\mathbb{R})$ satisfying $|n_{i j}| \leq b$ for all $i < j$, and denote by
\begin{equation*}
  \mathfrak{S}(b,c) := \Omega_N(b) A(\mathbb{R})^0_{\geq c} K_\infty  \subseteq G(\mathbb{R})
\end{equation*}
the corresponding standard Siegel domain.  A theorem of Siegel (see \cite[Prop 4.5]{MR147566}) implies that
\begin{enumerate}[(i)]
\item \label{itm:standard2:b-c-large-enough} if $b$ is large enough and $c$ is small enough (e.g., $b \geq 1/2$ and $c \leq \sqrt{3}/2$), then $G(\mathbb{Z}) \mathfrak{S}(b,c) = G(\mathbb{R})$, and
\item for all $b$ and $c$, we have
  \begin{equation}\label{eq:nubc}
    \nu(b,c) := \{\gamma \in G(\mathbb{Q}) : \gamma \mathfrak{S}(b,c) \cap \mathfrak{S}(b,c) \neq \emptyset \} < \infty.
  \end{equation}
\end{enumerate}

\begin{proof}[Proof of Proposition \ref{prop:siegel-upper}]
  Choosing $b$ and $c$ as in \eqref{itm:standard2:b-c-large-enough}, we have
  \begin{equation*}
    \int_{G(\mathbb{Z}) \backslash G(\mathbb{R})} \phi \leq \int_{\mathfrak{S}(b,c)} \phi
  \end{equation*}
  for all nonnegative $\phi$.  Let $f \in C_c^\infty(G(\mathbb{R}))$ be nonnegative and satisfy $\int f = 1$.  Write $R(f) \phi(x) := \int_{g \in G(\mathbb{R})} \phi(x g) f(g) \, d g$.  Applying the above inequality to $R(f) \phi$ gives
  \begin{align*}
    \int_{G(\mathbb{Z}) \backslash G(\mathbb{R})} \phi
    &\leq
      \int_{\mathfrak{S}(b,c)} R(f) \phi \\
    &=
      \int_{g \in G(\mathbb{R})} f(g)
      \int_{n \in \Omega_N(b)}
      \int_{t \in A(\mathbb{R})^0_{\geq c}}
      \int_{k \in K_\infty }
      \phi(n t k g)
      \, d k
      \, d n \, \frac{d t}{\delta_N(t)} \, d g.
  \end{align*}
  Let $a \in A(\mathbb{R})^0$ denote the element with $\det(a) = 1$ and $a_i/a_{i+1} = c$ for all $1 \leq i < \rank(G)$.  Then this last integral may be written
  \begin{equation*}
    \int_{t \in A(\mathbb{R})^0_{\geq 1}}
    \int_{g \in G(\mathbb{R})}
    \phi(t g)
    f_t(g) \, d g \, \frac{d t}{\delta_N(t)},
  \end{equation*}
  where
  \begin{equation*}
    f_t(g) :=
    \delta_N(a)^{-1}
    \int_{n \in \Omega_N(b)}
    \int_{k \in K_\infty }
    f((t^{-1} n t a k )^{-1} g)
    \, d k
    \, d n.
  \end{equation*}
  For each $n \in \Omega_N(b)$ and $t \in A(\mathbb{R})^0_{\geq 1}$, we have $t^{-1} n t \in \Omega_N$, hence $t^{-1} n t a k$ lies in the compact set $\Omega_N(b) a K_\infty$.  Therefore $f_t$ is supported on the compact set
  \begin{equation*}
    \Omega := (\Omega_N(b) a K_\infty)^{-1} \cdot \supp(f)
  \end{equation*}
  and bounded by
  \begin{equation*}
    C_2 := \delta_N(a)^{-1} \vol(\Omega_N(b), d n) \vol(K_\infty, d k) \|f\|_{\infty},
  \end{equation*}
  so
  \begin{equation*}
    \int_{G(\mathbb{Z}) \backslash G(\mathbb{R})} \phi \leq
    C_2 \int_{t \in A(\mathbb{R})^0_{\geq 1}}
    \int_{g \in \Omega } \phi(t g) \, d g \, \frac{d t}{\delta_N(t)},
  \end{equation*}
  as required.
\end{proof}
The proof of Proposition \ref{prop:siegel-lower} requires some further preliminaries.
\begin{lemma}\label{lem:siegel-upper-basic}
  For each $c > 0$, there exists $C_0 \geq 0$ so that for each nonnegative integrable function $\phi$ on $G(\mathbb{Z}) \backslash G(\mathbb{R})$, we have
  \begin{equation*}
    \int_{N(\mathbb{Z}) \backslash N(\mathbb{R})}
    \int_{A(\mathbb{R})^0_{\geq c}}
    \int_{K_\infty}
    \phi(n a k) \, d k \, \frac{d a}{\delta_N(a)} \, d n
    \leq C_0 \int_{G(\mathbb{Z}) \backslash G(\mathbb{R})} \phi.
  \end{equation*}
\end{lemma}
\begin{proof}
  Choose $b$ large enough that $\Omega_N(b)$ contains a fundamental domain for the compact quotient $N(\mathbb{Z}) \backslash N(\mathbb{R})$.  Then the conclusion holds with $C_0 := \nu(b,c)$ as defined in \eqref{eq:nubc}.
\end{proof}

\begin{lemma}\label{lem:scratch-research:each-c-Omega-N}
  For each $c > 0$ and compact subset $\Omega_N \subseteq N(\mathbb{R})$, there exists $C_0 \geq 0$ so that for each $t \in A(\mathbb{R})^0_{\geq c}$ and nonnegative integrable function $\phi$ on $N(\mathbb{Z}) \backslash N(\mathbb{R})$, we have
  \begin{equation*}
    \int_{\Omega_N}
    \phi(t n t^{-1}) \, d n \leq C_0 \int_{N(\mathbb{Z}) \backslash N(\mathbb{R})} \phi(n) \, d n.
  \end{equation*}
\end{lemma}
\begin{proof}
  The change of variables $n \mapsto t^{-1} n t$ gives
  \begin{equation*}
    \int_{\Omega_N}
    \phi(t n t^{-1}) \, d n
    =
    \frac{1}{\delta_N(t)}
    \int_{t \Omega_N t^{-1}}
    \phi(n) \, d n
    \leq
    \frac{C_1(t)}{\delta_N(t)} \int_{N(\mathbb{Z}) \backslash N(\mathbb{R}) }\phi(n) \, d n,
  \end{equation*}
  where
  \begin{equation*}
    C_1(t) := \sup_{n \in t \Omega_N t^{-1}}
    \# \left\{ \gamma \in N(\mathbb{Z}) : \gamma n \in t \Omega_N t^{-1} \right\}.
  \end{equation*}
  If $\gamma \in N(\mathbb{Z})$ satisfies the indicated condition for some $n \in t \Omega_N t^{-1}$, then $\gamma$ lies in $t \Omega_N t^{-1} n^{-1} \subseteq t(\Omega_N \cdot \Omega_N^{-1}) t^{-1}$.  It follows that the entries $\gamma_{i j} \in \mathbb{Z}$ for $i < j$ are bounded in magnitude by $C_2 t_i / t_j$ for some $C_2 = C_2(\Omega_N) > 0$, thus
  \begin{equation*}
    C_1(t) \leq \prod_{1 \leq i < j \leq \rank(G)} (1 + 2 C_2 t_i / t_j) \leq C_0 \delta(t)
  \end{equation*}
  for some $C_0 \geq 0$.  The required estimate then holds with this value of $C_0$.
\end{proof}

\begin{proof}[Proof of Proposition \ref{prop:siegel-lower}]
  By the Iwasawa decomposition, we have $\Omega \subseteq \Omega_N \Omega_A K_\infty$, where $\Omega_N \subseteq N(\mathbb{R})$ and $\Omega_A \subseteq A(\mathbb{R})^0$ are compact.  Then
  \begin{equation*}
    \int_{g \in \Omega } \phi(t g) \, d g
    \leq
    \int_{n \in \Omega_N}
    \int_{a \in \Omega_A}
    \int_{k \in K_\infty }
    \phi(t n a k)
    \, d k \, \frac{d a}{\delta_N(a)}  \, d n.
  \end{equation*}
  We write $t n a k = (t n t^{-1}) t a k$ and apply Lemma \ref{lem:scratch-research:each-c-Omega-N} to the $n$-integral, giving
  \begin{equation}\label{eq:int-_g-in-2}
    \int_{g \in \Omega } \phi(t g) \, d g \leq
    C_1 \int_{n \in N(\mathbb{Z}) \backslash N(\mathbb{R})}
    \int_{a \in \Omega_A } \int_{k \in K_\infty } \phi(n t a k)
    \, d k
    \, \frac{d a}{\delta_N(a)}
    \, d n,
  \end{equation}
  where $C_1$ depends only upon $\Omega$.  We substitute $a \mapsto t^{-1} a$, giving
  \begin{equation*}
    \int_{g \in \Omega } \phi(t g) \, d g \leq
    C_1 \delta_N(t) \int_{n \in N(\mathbb{Z}) \backslash N(\mathbb{R})}
    \int_{a \in t \Omega_A } \int_{k \in K_\infty } \phi(n a k)
    \, d k
    \, \frac{d a}{\delta_N(a)}
    \, d n.
  \end{equation*}
  We may find $c > 0$, depending only upon $\Omega$, so that $t \Omega_A \subseteq A(\mathbb{R})_{\geq c}^0$.  The conclusion then follows from Lemma \ref{lem:siegel-upper-basic}.
\end{proof}
\begin{proof}[Proof of Proposition \ref{prop:siegel-lower-integrated}]
  We argue as in the proof of Proposition \ref{prop:siegel-lower}, and integrate \eqref{eq:int-_g-in-2} over $t$:
  \begin{align*}
    &\int_{t \in A(\mathbb{R})^0_{\geq 1}} \int_{g \in \Omega } \phi(t g) \, d g \, \frac{d t}{\delta_N(t)} \\
    &\leq
      C_1 \int_{n \in N(\mathbb{Z}) \backslash N(\mathbb{R})}
      \int_{t \in A(\mathbb{R})^0_{\geq 1}}
      \int_{a \in \Omega_A } \int_{k \in K_\infty } \phi(n t a k)
      \, d k
      \, \frac{d a \, d t}{\delta_N(t a)}
      \, d n.
  \end{align*}
  We bring the $a$-integral outside and substitute $t \mapsto a^{-1} t$, giving
  \begin{equation*}
    C_1
    \int_{a \in \Omega_A }
    \int_{n \in N(\mathbb{Z}) \backslash N(\mathbb{R})}
    \int_{t \in a A(\mathbb{R})^0_{\geq 1}}
    \int_{k \in K_\infty } \phi(n t k)
    \, d k
    \, \frac{d t}{\delta_N(t)}
    \, d n \, d a.
  \end{equation*}
  We may find $c > 0$, depending only upon $\Omega_A$, so that that $a A(\mathbb{R})_{\geq 1}^0 \subseteq A(\mathbb{R})_{\geq c}^0$ for all $a \in \Omega_A$.  The conclusion then follows from Lemma \ref{lem:siegel-upper-basic} and the finiteness of the volume of $\Omega_A$.
\end{proof}

\subsection{Petersson norms vs.\ Whittaker norms}\label{sec:norms-cusp-forms}
Let $G$ be a general linear group over $\mathbb{Z}$.  Let $\psi$ be an unramified (\S\ref{sec:stand-nond-char}) nondegenerate unitary character of $N(\mathbb{A})$, trivial on $N(\mathbb{Q})$.  Let $\pi$ be a unitary cuspidal automorphic representation of $G(\mathbb{A})$.  For each prime number $p$, let $W_p \in \mathcal{W}(\pi_p, \psi_p)$ denote the essential vector (\S\ref{sec:essential-vector}), normalized by $W_p(1) = 1$.  Denote as usual by $C(\pi_{\fin}) \in \mathbb{Z}_{\geq 1}$ the finite conductor of $\pi$, i.e., the product of the local conductors $C(\pi_p)$ defined in \S\ref{sec:essential-vector}.

For each archimedean Whittaker function $W_\infty \in \mathcal{W}(\pi_\infty, \psi_\infty)$, there is a unique automorphic form $\varphi \in \pi$ whose Whittaker function
\begin{equation*}
  W_\varphi(g) = \int_{[N]} \varphi(u g) \psi^{-1}(u) \, d u
\end{equation*}
is the product of the $W_\mathfrak{p}$:
\begin{equation*}
  W_{\varphi}(g) = \prod_{\mathfrak{p}} W_{\mathfrak{p}}(g_\mathfrak{p}).
\end{equation*}
With respect to our fixed Haar measures, we define the Petersson norm
\begin{equation*}
  \|\varphi \|^2 := \int_{[G/Z]} |\varphi|^2
\end{equation*}
and the Whittaker norm
\begin{equation*}
  \|W_\infty \|^2 := \int_{N_H(\mathbb{R}) \backslash H(\mathbb{R})} |W_\infty|^2,
\end{equation*}
as in \S\ref{sec:local-prelim-kirillov-model}.  They may be compared via the following standard estimates, obtained by combining \cite[\S7]{MR2200849} and \cite[\S1.3]{MR2595006}.

\begin{lemma}\label{lem:nonstandard-norm-whittaker}
  Let $G$ and $\psi$ be as above, with $G$ fixed.  Let $\pi$ be a unitary cuspidal automorphic representation of $G$.  Let $W_\infty \in \mathcal{W}(\pi_\infty, \psi_\infty)$.  Define $\varphi \in \pi$ as in \S\ref{sec:norms-cusp-forms}.  Then
  \begin{equation*}
    \|\varphi \|^2
    \asymp  C(\pi_{\fin})^{o(1)}
    \res_{s=1}L(\pi \times \tilde{\pi}, s)
    \int_{N_H \backslash H}
    |W_\infty|^2
  \end{equation*}
  and
  \begin{equation*}
    \res_{s=1}L(\pi \times \tilde{\pi}, s) \ll C(\pi)^{o(1)},
  \end{equation*}
  hence
  \begin{equation*}
    \|\varphi \|^2
    \ll C(\pi) ^{o(1)}
    \int_{N_H \backslash H}
    |W_\infty|^2.
  \end{equation*}
\end{lemma}

\subsection{The distance function $d_H$}\label{sec:dist-funct-d_h}
We specialize the definition of  \cite[\S4.2]{2020arXiv201202187N} to the case that $(G,H)$ is a fixed general linear pair over $\mathbb{R}$, say $(G,H) = (\GL_{n+1}(\mathbb{R}), \GL_n(\mathbb{R}))$.  Set $\bar{G} := G/Z = \PGL_{n+1}(\mathbb{R})$.  We may regard $H$ as a subgroup of $\bar{G}$.  For a matrix $x$, we write $|x|$ for the maximum of the absolute values of its entries.  In particular, for $g \in \bar{G}$, we may regard $\Ad(g)$ as a matrix (using the standard basis for $\Lie(G)$) and define $|\Ad(g) - 1|$, which quantifies the distance from $g$ to the identity element.

\begin{definition}
  Let $g \in \bar{G}$.  We may write
  \begin{equation*}
    g = \begin{pmatrix}
      a & b \\
      c & d
    \end{pmatrix},
    \quad
    g^{-1} =
    \begin{pmatrix}
      a' & b' \\
      c' & d'
    \end{pmatrix},
  \end{equation*}
  where $a,b,c,d$ and $a',b',c',d'$ are matrices of respective dimensions $n \times n, n \times 1, 1 \times n, 1 \times 1$, well-defined up to simultaneous scaling.  We then define \index{distance function $d_H$}
  \begin{equation*}
    d_H(g) := \min\left(1, |b/d| + |b'/d'| + |c/d| + |c'/d'|\right),
  \end{equation*}
  with the convention that $d_H(g) := 1$ if either $d$ or $d'$ vanishes.
\end{definition}

Observe that $d_H(g) = 0$ if and only if $g$ lies in (the image of) $H$.  We may thus regard $d_H(g)$ as a quantification of the distance from $g$ to $H$.

\begin{remark}
  The definition given here differs mildly from that in \cite[\S4.2]{2020arXiv201202187N}.  Call the latter $d_H^0(g)$.  As explained in \emph{loc.\ cit.}, we have $d_H(g) \asymp d_H^0(g)$ for all $g \in \bar{G}$.  Since we will only use $d_H(g)$ for the purpose of quantifying estimates, the distinction between the two definitions is irrelevant.
\end{remark}

We recall \cite[Lemma 4.1]{2020arXiv201202187N}:
\begin{lemma}
  Suppose that $(G,H)$ is fixed.  Let $g \in \bar{G}$.
  \begin{enumerate}[(i)]
  \item We have
    \begin{equation}\label{eqn:d-H-Z-vs-dist-G-mod-Z}
      d_{H}(g) \ll |\Ad(g)-1|.
    \end{equation}
  \item
    Let $u_1, u_2 \in \bar{G}$.
    Suppose that $g$
    belongs to a fixed compact subset of $\bar{G}$
    and that
    \[
      |\Ad(u_j)-1| \lll
      d_{H}(g)
      \quad (j=1,2).
    \]
    Then
    \begin{equation}\label{eqn:d-H-Z-compare-u-g-g}
      d_{H}(u_1 g u_2) \asymp d_{H}(g).
    \end{equation}
  \item Let $h_1, h_2 \in H$.  Suppose that $g,h_1,h_2$
    each belong to some fixed compact subset of $\bar{G},H,H$, respectively.  Then
    \begin{equation}\label{eqn:d-H-h1gh2}
      d_{H}(h_1 g h_2)
      \asymp d_{H}(g).
    \end{equation}
  \end{enumerate}
\end{lemma}

\part{Division of the argument}\label{part:reduction-proof}

\section{Construction of test vectors: statement}\label{sec:20230516081722}
The following is the main result of Part \ref{part:constr-test-vect}, which relies in turn on Part \ref{part:asympt-analys-kirill}.  The statement is regrettably complicated, but has the advantage of being a purely local assertion concerning representations of real groups.  The sketch of \S\ref{sec:high-level-overview} may provide a useful roadmap for the relevance of many assertions in the following result to our global argument.  We note that, in the language of that sketch, the Whittaker function $W$ below describes the global cusp form $\varphi$.  The reader may also wish to compare with \cite[Thm 4.2]{2020arXiv201202187N}, the analogous result in the compact quotient case, which is simpler.

\begin{theorem}\label{thm:main-local-results}
  Suppose given a fixed general linear pair $(G,H)$ over $\mathbb{R}$, with accompanying notation $A$, $B$, $N$, $Q$, $U$, $Z$, $A_H$, $B_H$, $N_H$, $Q_H$, $U_H$, $Z_H$, as in \S\ref{sec:gener-line-groups-1}.  Let $\psi$ be a nondegenerate unitary character of $N$ equal to either the standard character or its inverse, as in \S\ref{sec:stand-nond-char}.  Let $T$ be a nonnegative real with $T \ggg 1$.  Let $\pi$ be a generic irreducible unitary representation of $G$, each of whose archimedean $L$-function parameters has magnitude $\asymp T$.  There exists
  \begin{itemize}
  \item $f, f_{\ast} \in \mathcal{S}^e(U_H \backslash H)^{W_H}$ (cf.\ \eqref{eq:mathc-backsl-gw-1}), a fixed element $\phi_{Z_H} \in C_c^\infty(Z_H)$,
  \item a smooth vector $W \in \mathcal{W}(\pi,\psi)$ of norm $\leq 1$ (with respect to the norm defined in \S\ref{sec:local-prelim-kirillov-model}), and
  \item $\omega_0 \in C_c^\infty(G)$, supported in each fixed neighborhood of the identity element,\footnote{For the sake of orientation, we note that in a paper without the nonstandard formalism, the support condition on $\omega_0$ would be phrased as follows: its support shrinks to the identity element, arbitrarily slowly, as the parameter $T$ tends off to $\infty$.}
  \end{itemize}
  with the following further properties.  Define
  \begin{equation*}
    \omega(g) := (\omega_0 \ast \omega_0^*)(g) := \int_{h \in G} \omega_0(g h^{-1}) \overline{\omega_0(h^{-1})} \, d h,
  \end{equation*}
  \begin{equation*}
    \omega^\sharp (g) := \int_{z \in Z}
    \pi|_Z (z) \omega(z g) \, d z,
  \end{equation*}
  where $\pi|_Z$ denotes the central character of $\pi$.
  \begin{enumerate}[(i)]
  \item \label{itm:sub-gln:5} $f_{\ast}(g) = \int_{z \in Z_H} f(g z^{-1}) \phi_{Z_H}(z) \, d z$.
  \item \label{itm:sub-gln:6} Each of the following assertions remains valid upon replacing $f$ with $f_{\ast}$.
  \item \label{itm:sub-gln:7} For each fixed compact subset $\mathcal{D}$ of $\mathfrak{a}^*$ and fixed $\ell \in \mathbb{Z}_{\geq 0}$, we have (see \S\ref{sec:local-sobolev-norms-prelims} for notation)
    \begin{equation*}
      \nu_{\mathcal{D},\ell,T}(f) \ll T^{o(1)}.
    \end{equation*}
  \item \label{itm:sub-gln:8} The local zeta integral $Z(W,f[s])$ (defined in \S\ref{sec:zeta-integrals}, and evaluated at the Mellin component $f[s]$ defined in \S\ref{sec:mellin-paley-affine-quotients}) is entire in $s$.  We have the lower bound
    \begin{equation*}
      Z(W,f[0]) \gg T^{- \dim(H)/4}.
    \end{equation*}
    For each $s \in \mathfrak{a}_{\mathbb{C}}^*$ with $\Re(s) \ll 1$, we have the upper bound
    \begin{equation*}
      Z(W,f[s]) \ll  T^{\rank(G) \trace(s)/2 - \dim(H)/4 + o(1) }
      (1 + |s|)^{-\infty }.
    \end{equation*}
  \item \label{itm:sub-gln:1} For each fixed $u \in \mathfrak{U}(G)$, we have
    \begin{equation*}
      \|\pi(u) ( \pi(\omega_0) W - W)\| \ll T^{-\infty}.
    \end{equation*}
  \item \label{itm:sub-gln:f-L-one-H}   We have $\int_{H} |\omega^\sharp| \leq T^{\rank(H)/2 + o(1)}$.
  \item \label{itm:sub-gln:3}  For each fixed compact subset $\Omega_H$ of $H$, there is a fixed compact subset $\Omega_H'$ of $H$ with the following property.  Let $\Psi_1, \Psi_2 : H \rightarrow \mathbb{C}$ be continuous functions.  For $j=1,2$, define the convolutions
    \begin{equation*}
      (\Psi_j \ast \phi_{Z_H})(x) := \int_{z \in Z_H} \Psi_j(x z^{-1}) \phi_{Z_H}(z) \, d z.
    \end{equation*}
    Let $\gamma \in \bar{G} - H$.  Then the integral
    \begin{equation}\label{eq:i-:=-int-1}
      I := \int_{x, y \in \Omega_H }
      \left\lvert
        (\Psi_1 \ast \phi_{Z_H})(x)
        (\Psi_2 \ast \phi_{Z_H})(y)
        \omega^\sharp(x^{-1} \gamma y)
      \right\rvert
      \, d x \, d y
    \end{equation}
    satisfies the estimate
    \begin{equation*}
      I \ll
      T^{\rank(H)/2 + o(1)}
      \min\left( 1,
        \frac{T^{-1/2}}{ d_H(\gamma)}
      \right)
      \|\Psi_1\|_{L^2(\Omega_H')}
      \|\Psi_2\|_{L^2(\Omega_H')}.
    \end{equation*}
  \end{enumerate}
\end{theorem}
We refer to \S\ref{sec:elem-form} for an elementary reformulation of this result, without the language of nonstandard analysis.
\begin{remark}\label{rmk:we-briefly-comment}
  We briefly comment on the significance of the numbered assertions, matching them with parts of the sketch given in \S\ref{sec:high-level-overview}.

  (i) and (ii) are purely technical.  As promised in \S\ref{sec:analyt-test-vect}, the functions on $U_{H} \backslash H$ that we construct are uniformly smooth under left translation by $A_H$, hence in particular approximately invariant under $Z_H$.  For later applications, the simplest way to ``remember'' this feature is to represent such elements explicitly as convolutions $f_{\ast}$ against some fixed function on $Z_H$.

  (iii) says that $f$ has the shape promised in \S\ref{sec:analyt-test-vect}.  The relation between bounds for $\nu_{\mathcal{D},\ell,T}(f)$ and that informal description will be made successively more explicit in Lemmas \ref{lem:standard:frakE-equivalences} and \ref{lem:standard:let-f-in-2}.

  (iv) gives the control over local zeta integrals promised around \eqref{eq:zvarphi-psi-approx}.

  (v) says that, for all practical purposes, $\pi(\omega_0) W \approx W$, as promised in \S\ref{sec:relat-trace-form} (recall that $W$ describes $\varphi$).

  (vi) is consistent with the shape of $\omega$ (or $\omega_0$) promised in \S\ref{sec:relat-trace-form}, in view of the stability of $\tau$ (see \S\ref{sec:stability}, Lemma \ref{lem:standard2:no-conductor-drop-vs-tau}, \cite[Thm 17.1]{nelson-venkatesh-1} and \cite[\S14.9]{2020arXiv201202187N} for details).

  (vii) asserts a precise form of the promised ``bilinear forms'' estimate \eqref{eq:t-n2-int-1}.
\end{remark}

\section{Growth bounds for Eisenstein series: statement}\label{sec:cnpy9ba61h}
We retain the setting and notation of \S\ref{sec:siegel-domains}, working with $G := \GL_n$ over $\mathbb{Z}$ for some natural number $n$.

Let $\Omega$ be a compact subset of $G(\mathbb{R})$.  For each $t \in A(\mathbb{R})^0_{\geq 1}$, we define a seminorm $\|.\|_{t,\Omega}$ on $L^2([G])$ by
\begin{equation}\label{eq:varphi-2_r-:=}
  \|\varphi \|^2_{t,\Omega} :=
  \int_{g \in \Omega \times \prod_p K_p }
  |\varphi(t g)|^2
  \, d g.
\end{equation}
As $t$ varies with $\Omega$ fixed, we may understand these seminorms as describing the local $L^2$-norm of $\varphi$ near a given point in $[G]$.

The following is the main result of Part \ref{part:local-l2-growth}.

\begin{theorem}\label{thm:growth-eisenstein-nonstandard}
  Let $n$ be a fixed natural number, $G := \GL_n$ over $\mathbb{Z}$.  Let $\Omega \subseteq G(\mathbb{R})$ be a fixed compact.  Let $T \geq 1$.  Let $Y \in A(\mathbb{R})^0$.  Assume that $T^{-\kappa} \leq |Y_j| \leq T^{\kappa}$ for all $j$ and some fixed $\kappa \in (0,1/2)$.  Let $t \in A(\mathbb{R})^0_{\geq 1}$.  Let $f \in \mathcal{S}^e(U(\mathbb{R}) \backslash G(\mathbb{R}))^W$ such that for each fixed compact subset $\mathcal{D}$ of $\mathfrak{a}^*$ and fixed $\ell \in \mathbb{Z}_{\geq 0}$, we have
  \begin{equation}\label{eq:nu_mathcald-ell-tf}
    \nu_{\mathcal{D},\ell,T}(f) \ll T^{o(1)},
  \end{equation}
  with the Sobolev norm $\nu_{\mathcal{D},\ell,T}$ as defined in \S\ref{sec:local-sobolev-norms-prelims}.  Then we have the following local $L^2$-norm bound for pseudo Eisenstein series, defined as in \S\ref{sec:notat-basic-defin}:
  \begin{equation}\label{eq:fraceis--phily}
    \frac{\|\Eis [ \Phi[L(Y) f] ] \|_{t,\Omega}^2}{\delta_N(t)}
    \ll T^{o(1)} \min \left( \det(Y)^{-1} t_1^{-n}, \det(Y) t_n^n \right).
  \end{equation}
\end{theorem}
We refer again to \S\ref{sec:elem-form} for an elementary reformulation of this result, without the language of nonstandard analysis.

\begin{remark}
  For the sake of comparison with the sketch of \S\ref{sec:high-level-overview}, we note that the hypothesis \eqref{eq:nu_mathcald-ell-tf} says that $f$ has the shape promised in \S\ref{sec:analyt-test-vect} (cf.\ Remark \ref{rmk:we-briefly-comment}, Lemmas \ref{lem:standard:frakE-equivalences} and \ref{lem:standard:let-f-in-2}), while the conclusion \eqref{eq:fraceis--phily} is essentially the promised pair of estimates \eqref{eq:psi-_l2t-mathcalh2} and \eqref{eq:psi-_l2t-mathcalh2-2}.  The reader might wish to pretend that $Y = 1$ on a first reading, since for the purpose of obtaining a qualitative subconvex bound (without an effective exponent), one can take $\kappa$ arbitrarily small.
\end{remark}

\section{Reduction of the proof}\label{sec:reduction-proof-0}
We aim now to deduce Theorem \ref{thm:main-result-reformulated} from Theorems \ref{thm:main-local-results} and \ref{thm:growth-eisenstein-nonstandard}.

\subsection{Setup and preliminary reductions}\label{sec:setup-prel-reduct}
Let $(G,H)$ be a fixed general linear pair over $\mathbb{Z}$, with $H$ nontrivial.  In this section, we label indices as follows:
\begin{equation*}
  (G,H) = (\GL_{n+1}, \GL_n).
\end{equation*}
Thus $n$ is a fixed natural number $\geq 1$.

Let $T \geq 1$.  Let $\pi$ be a unitary cuspidal automorphic representation of $G$, each of whose archimedean $L$-function parameters are $\asymp T$.  We must show, for suitable fixed $\delta > 0$, that
\begin{equation}\label{eq:lpi-tfrac12-ll}
  L(\pi,\tfrac{1}{2}) \ll C(\pi_{\fin})^{1/2} T^{(1 - \delta) (n+1)/4}.
\end{equation}

Consider first the case that $T \ll 1$.  Then $C(\pi_\infty) \ll 1$, so $C(\pi) \asymp C(\pi_{\fin})$.  In that case, the convexity bound (Lemma \ref{lem:convexity-bound}) implies that $L(\pi,\tfrac{1}{2}) \ll C(\pi_{\fin})^{1/4 + o(1)}$, while the RHS of \eqref{eq:lpi-tfrac12-ll} is $\asymp C(\pi_{\fin})^{1/2}$.  The required estimate \eqref{eq:lpi-tfrac12-ll} thus holds in that case.  We thereby reduce to the case that
\begin{equation}\label{eq:t-ggg-1}
  T \ggg 1.
\end{equation}

We may similarly reduce to the case that
\begin{equation}\label{eq:n_pi-=-to1}
  C(\pi_{\fin}) \leq T^{\O(1)}.
\end{equation}
Otherwise, the RHS of \eqref{eq:lpi-tfrac12-ll} is readily seen to be $\gg C(\pi)^{1/3}$, say, and the required bound follows again from convexity.

We appeal throughout this section to Theorem \ref{thm:main-local-results}.

\subsection{Construction of automorphic forms}\label{sec:constr-autom-forms}
We let $\psi$ denote the standard nondegenerate unitary character of $N(\mathbb{A})$, trivial on $N(\mathbb{Q})$ (\S\ref{sec:stand-nond-char}).  We denote also simply by $\psi$ its restriction to $N(\mathbb{Q}_\mathfrak{p})$ for each place $\mathfrak{p}$ of $\mathbb{Q}$, and also to $N_H(\mathbb{A}), N_H(\mathbb{Q}_\mathfrak{p})$.

Let $W_{\infty} \in \mathcal{W}(\pi_\infty, \psi)$ denote the smooth vector furnished by Theorem \ref{thm:main-local-results}.  Let $\varphi \in \pi$ denote the automorphic form attached to $W_{\infty}$, as in \S\ref{sec:norms-cusp-forms}, by requiring that its local components at finite primes be normalized essential vectors $W_{\pi,p} \in \mathcal{W}(\pi_p, \psi)$.  By the estimates of that section (cf.\ Lemma \ref{lem:nonstandard-norm-whittaker}) and our assumption \eqref{eq:n_pi-=-to1}, we then have the Petersson norm estimate
\begin{equation}\label{eq:varphi-2-asymp}
  \|\varphi \|
  \leq T^{o(1)}.
\end{equation}

Since $\varphi$ is cuspidal, its restriction to $[H]$ decays rapidly, so the integral $\int_{[H]} \varphi \Psi$ converges absolutely for any automorphic function $\Psi$ on $[H]$ of moderate growth.  We may unfold that integral using the Whittaker expansion (cf.\ \cite[Thm 5.9]{MR348047})
\begin{equation*}
  \varphi(g) = \sum_{\gamma \in N_H(\mathbb{Q}) \backslash H(\mathbb{Q})} W_\varphi(\gamma g).
\end{equation*}
Working formally for the moment, we obtain
\begin{equation}\label{eqn:unfolding-rs-1}
  \int_{[H]} \varphi \Psi
  =
  \int_{N_H(\mathbb{Q}) \backslash H(\mathbb{A})}
  W_\varphi \Psi
  =
  \int_{N_H(\mathbb{A}) \backslash H(\mathbb{A})}
  W_\varphi \tilde{W}_{\Psi},
\end{equation}
where the dual Whittaker function $\tilde{W}_{\Psi}$ is defined by
\begin{equation*}
  \tilde{W}_{\Psi}(g)
  = \int_{[N_H]}
  \Psi(u g) \psi(u ) \, d u.
\end{equation*}
Such calculations \eqref{eqn:unfolding-rs-1} are the starting point for the study of $\GL_{n+1} \times \GL_n$ zeta integrals due to Jacquet, Piatetski-Shapiro and Shalika.

For each $s \in \mathfrak{a}_{H,\mathbb{C}}^*$, we denote by $\mathcal{I}(s) = \Ind_{Q_H(\mathbb{A})}^{H(\mathbb{A})} (|.|^s)$ the corresponding principal series representation of $H(\mathbb{A})$, as defined in \S\ref{sec:induc-repr}; we recall that it is the restricted tensor product of local principal series representations $\mathcal{I}_\mathfrak{p}(s) = \Ind_{Q_H(\mathbb{Q}_\mathfrak{p})}^{H(\mathbb{Q}_\mathfrak{p})}(|.|_{\mathfrak{p}}^s)$.

Let $f, f_{\ast} \in \mathcal{S}^e(U_H(\mathbb{R}) \backslash H(\mathbb{R}))^{W_H}$ denote the vectors furnished by Theorem \ref{thm:main-local-results}, with associated Mellin components $f[s], f_{\ast}[s] \in \mathcal{I}_\infty(s)$.  Recall the elements
\begin{equation}\label{eq:phi-:=-phif}
  \Phi[f] \in \mathcal{S}^e([H]_{Q_H})
\end{equation}
constructed in \S\ref{sec:compl-eisenst-seri}, with Mellin components
\begin{equation*}
  \Phi[f][s] \in \mathcal{I}(s).
\end{equation*}
This notation and the following discussion applies also with $f$ replaced by $f_{\ast}$.

Suppose for the moment that $\Re(s) \succ \rho_U$, so that the series definition of, e.g., $\Eis[\Phi[f][s]] \in \mathcal{T}([H])$ applies.  By the unfolding identity \eqref{eq:weisf-psig-=}, the definition (\S\ref{sec:notat-basic-defin}) of $\Phi[f][s]$ and the Casselman--Shalika formula \eqref{eq:wf0chi-psi1-casselman-shalika}, we see that the global Whittaker function of $\Eis(\Phi[f][s])$ is given by
\begin{equation}\label{eq:weisphis-psi-1g}
  W[\Eis[\Phi[f][s]], \psi^{-1}](g)
  =
  \mathcal{P}_H(s)
  W[f[s], \psi^{-1}](g_\infty)
  \prod_p W_{s,p}(g_p)
\end{equation}
where $W_{s,p} \in \mathcal{W}(\mathcal{I}_p(s), \psi^{-1})$ denotes the spherical element of the local Whittaker model satisfying the normalization ${W}_{s,p}(1) = 1$.

Let $Y$ be an element of the identity component $A_H(\mathbb{R})^0$ of $A_H(\mathbb{R})$.  We identify $Y$ with an element of $A_H(\mathbb{A})$ via the standard inclusion.  We define the element
\begin{equation*}
  L(Y) \Phi[f] \in \mathcal{S}([H]_{B_H}),
\end{equation*}
as in \eqref{eq:normalized-left-translation}, by normalized left translation:
\begin{equation}\label{eq:ly-phig-=}
  L(Y) \Phi[f](g) = \delta_{U_H}^{-1/2}(Y) \Phi[f](Y g),
\end{equation}
or equivalently, via the Mellin integral
\begin{equation*}
  L(Y) \Phi[f] = \int_{(\sigma)} Y^{s} \Phi[f][s] \, d \mu_{[A_H]}(s),
\end{equation*}
convergent for all $\sigma \in \mathfrak{a}_H^*$.  (Here $Y^s$ abbreviates $|Y|^s$, noting that the entries of $Y$ are positive.)  We define the associated pseudo Eisenstein series \index{Eisenstein series!$\Psi = \Psi[Y]$}
\begin{equation*}
  \Psi[Y] := \Eis[L(Y) \Phi[f]] \in \mathcal{S}([H]),
\end{equation*}
which in turn admits the Mellin integral representation
\begin{align*}
  \Psi[Y] &=
         \int_{(\sigma)} Y^{s}
         \Eis[\Phi[s]]
         \, d \mu_{[A_H]}(s) \\
       &=
         \int_{(\sigma)}
         Y^{s} \mathcal{P}_H(s) \zeta(U_H,s)
         \Eis[f[s] \otimes (\otimes_{p} f_p[s])] \, d \mu_{[A_H]}(s).
\end{align*}
We analogously define
\begin{equation*}
  \Psi_{\ast}[Y] := \Eis[L(Y) \Phi[f_{\ast}]],
\end{equation*}
which admits an analogous integral representation.

\subsection{Reduction to period bounds}\label{sec:reduct-peri-bounds}
Since $\varphi|_{[H]}$ and $\Psi[Y]$ are (both) of rapid decay, the integral of their product
\begin{equation*}
  I(Y) := \int_{[H]} \varphi \Psi_{\ast}[Y]
\end{equation*}
converges absolutely.  The following lemma reduces our task of bounding $L(\pi,\tfrac{1}{2})$ to one of bounding the periods $I(Y)$ for $Y$ ``fairly close'' to the identity element $1$.

\begin{lemma}\label{lem:sub-gln:fix-kappa-in}
  Fix $\kappa \in (0,1/2)$ and $\delta > 0$.  Assume that whenever each component $Y_j$ of $Y$ satisfies $T^{-\kappa} \leq |Y_j| \leq T^{\kappa}$, we have
  \begin{equation}\label{eq:iy-ll-tn4}
    I(Y) \ll C(\pi_{\fin})^{n/2} T^{n/4 -\delta}.
  \end{equation}
  Then
  \begin{equation}\label{eq:lpi-tfrac12n-ll_eps}
    L(\pi,\tfrac{1}{2})^{n} \ll
    C(\pi_{\fin})^{n/2}
    T^{n(n+1)/4 - \min(\kappa/2,\delta) + o(1)}.
  \end{equation}
\end{lemma}
\begin{proof}
  We begin by expanding $I(Y)$ as a contour integral.  Since $\varphi$ is cuspidal, we have $\varphi|_{[H]} \in \mathcal{S}([H])$.  By \eqref{eq:int-_g-eisf}, it follows that for each $\sigma \in \mathfrak{a}_H^*$ with $\sigma \succ \rho_{U_H}$, we have
  \begin{equation}\label{eq:iy-=-int_c}
    I(Y)
    =
    \int_{(\sigma)}
    Y^{s}
    \tilde{I}(s)
    \, d s,
  \end{equation}
  where
  \[
    \tilde{I}(s) := \int_{[H]} \varphi \Eis[\Phi[f_{\ast}][s]].
  \]

  We now apply the unfolding identity \eqref{eqn:unfolding-rs-1}, the Whittaker function evaluation \eqref{eq:weisphis-psi-1g}, and the local unramified evaluation given by Theorem \ref{thm:jpss-essential-vector-factorized}.  Arguing formally for the moment, we obtain
  \begin{align*}
    \nonumber
    \tilde{I}(s)
    &= \int_{N_H(\mathbb{A}) \backslash H(\mathbb{A})} W_\varphi(g) W[\Eis[\Phi[f_{\ast}][s]], \psi^{-1}](g) \, d g \\
    \nonumber
    &=
      \mathcal{P}_H(s) Z_\infty(s)
      \prod_p
      \int_{N_H(\mathbb{Q}_p) \backslash H(\mathbb{Q}_p)} W_{\pi,p}(g) W_{s,p}(g) \, d g \\
    &=
      Z_\infty(s) \mathcal{L}(s),
  \end{align*}
  where
  \[
    Z_\infty(s)
    :=
    \int_{h \in N_H(\mathbb{R}) \backslash H(\mathbb{R})}
    W_{\infty}(h) W[f_{\ast}[s], \psi^{-1}](h) \, d h
  \]
  and
  \[
    \mathcal{L}(s) :=
    \mathcal{P}_H(s)
    \prod_{j=1}^n L(\pi, \tfrac{1}{2} + s_j).
  \]
  To make this calculation  rigorous, we first introduce a complex parameter $u \in \mathbb{C}$ of sufficiently large real part.  Then, by the same proof as in the case of a cuspidal Rankin--Selberg convolution (see \cite[\S2.2.2]{MR2508768}, \cite[\S13]{MR528964}), the integral
  \begin{equation*}
    \int_{N_H(\mathbb{A}) \backslash H(\mathbb{A})} W_\varphi(g) W[\Eis[\Phi[f_{\ast}][s]], \psi^{-1}](g) |\det g|^u \, d g
  \end{equation*}
  converges absolutely.  On the one hand, it folds up to the (convergent) global integral
  \begin{equation*}
    \tilde{I}(s,u) := \int_{[H]} \varphi \Eis[\Phi[f_{\ast}][s]] |\det|^u,
  \end{equation*}
  which is analytic in $u$ and satisfies $\tilde{I}(s,0) = \tilde{I}(s)$.  On the other hand, it factors as a product of local integrals $Z_\infty(s,u) \mathcal{L}(s + u e)$, where $e = (1,\dotsc,1)$ and $Z_\infty(s,u)$ is defined like $Z_\infty(s)$, but with an additional factor of $|\det h|^u$ in the integrand.  By specializing the resulting identity to $u=0$, we deduce the claimed factorization of $\tilde{I}(s)$.

  We abbreviate
  \begin{equation*}
    \langle s \rangle := 1 + |s|.
  \end{equation*}
  Recall from Theorem \ref{thm:main-local-results}, part \eqref{itm:sub-gln:8} that $Z_\infty(s)$ is entire and of rapid vertical decay, i.e., $Z_\infty(s) \ll \left\langle s \right\rangle^{-\infty}$ for $\Re(s) \ll 1$.  On the other hand, $\mathcal{L}(s)$ is entire and polynomially bounded in vertical strips, i.e., $\mathcal{L}(s) \ll C(\pi)^{\O(1)} \left\langle s \right\rangle^{\O(1)}$ for $\Re(s) \ll 1$ (see \eqref{eq:lpi-tfrac12-+-2}).  The product $\tilde{I}(s)$ is thus entire and of rapid decay in vertical strips, and so we may shift the contour integral in \eqref{eq:iy-=-int_c} to any $\sigma \in [-1/2, 1/2]^n$.

  By Theorem \ref{thm:main-local-results}, we have
  \[
    Z_\infty(s) \ll T^{\frac{n+1}{2} \sum s_j - n^2/4 + o(1)} \langle s \rangle^{-\infty},
    \quad
    Z_\infty(0) \gg T^{-{n^2}/{4}}.
  \]
  We appeal to the convexity bound \eqref{eqn:convexity-polynomially-in-s}, which gives for $|\Re(s_j)| \leq 1/2$ the estimate
  \begin{equation*}
    L(\pi, \tfrac{1}{2} +s_j)
    \ll
    C(\pi_{\fin})^{(1 - 2 s_j)/4}
    T^{(n+1)(1 - 2 s_j)/4 + o(1)}
    \langle s_j \rangle^{\O(1)}.
  \end{equation*}
  Estimating quite crudely with respect to the finite conductor $C(\pi_{\fin})$, we deduce that for $\max_j |\Re(s_j)| \leq 1/2$,
  \begin{equation*}
    \tilde{I}(s)
    \ll
    C(\pi_{\fin})^{n/2}
    T^{n/4 + o(1) }
    \langle s \rangle^{-\infty}.
  \end{equation*}
  By shifting the Mellin expansion \eqref{eq:iy-=-int_c} to the contour defined by the conditions $\Re(s_j) = \sigma_j$ with $\sigma_j = \pm 1/2$ (taking all possible signs), we deduce the following ``convexity'' bound for the period $I(Y)$:
  \begin{equation}\label{eq:iy-ll-t}
    I(Y) \ll
    C(\pi_{\fin})^{n/2}
    T^{n/4 + o(1) }
    (Y^\dagger )^{-1/2},
    \quad
    Y^\dagger :=
    \prod_{j=1}^n
    \max(Y_j, Y_j^{-1}).
  \end{equation}

  Let $\mathcal{Y}$ denote the set of all $Y$ for which $T^{-\kappa} \leq Y_j \leq T^{\kappa}$ for all $j$.  Write $d Y = d Y_1 \dotsb d Y_n /(Y_1 \dotsb Y_n)$ for the multiplicative Haar measure on $A_H(\mathbb{R})^0$.  Fix $\eps \in (0,1/2)$.  For $Y \notin \mathcal{Y}$, we have $Y_i \geq T^\kappa$ or $Y_i \leq T^{-\kappa}$ for some $i \in \{1, \dotsc, n\}$, hence $Y^\dagger \geq T^{\kappa} \prod_{j \neq i} \max(Y_j, Y_j^{-1})$, and therefore
  \begin{align*}
    \int_{Y \notin \mathcal{Y} } (Y^\dagger)^{\eps - 1/2} \, d Y
    &\leq
      n T^{(\eps - 1/2) \kappa} \left(
      \int_0^\infty
      \max(y, y^{-1})^{\eps - 1/2}
      \, \frac{d y}{y}
      \right)^{n - 1}
    \\
    &\ll T^{(\eps - 1/2) \kappa}.
  \end{align*}
  By \eqref{eq:iy-ll-t}, it follows that
  \begin{equation}\label{eq:int_y-notin-mathc}
    \int_{Y \notin \mathcal{Y} }
    |I(Y)| (Y^\dagger )^{\eps} \, d Y
    \ll
    C(\pi_{\fin})^{n/2}
    T^{n/4  + o(1) - (1/2-\eps) \kappa}.
  \end{equation}
  Recall from \eqref{eq:iy-ll-tn4} our assumption that $I(Y) \ll C(\pi_{\fin})^{n/2} T^{n/4 - \delta}$ for $Y \in \mathcal{Y}$.  Integrating this assumption over $\mathcal{Y}$ and combining with \eqref{eq:int_y-notin-mathc} gives the overall bound
  \begin{equation}
    \int_{Y} |I(Y)| (Y^\dagger )^\eps  \, d Y
    \ll
    C(\pi_{\fin})^{n/2}
    T^{n/4  - \min(\kappa/2,\delta)  + \eps (1 + \kappa) },
  \end{equation}
  say.  By Mellin inversion, we deduce that
  for
  $\max_j
  |\Re(s_j)| \leq \eps$,
  \begin{equation*}
    \tilde{I}(s)
    \ll
    C(\pi_{\fin})^{n/2}
    T^{n/4  - \min(\kappa/2,\delta)  + \eps (1 + \kappa) }.
  \end{equation*}
  By an $n$-fold application of Cauchy's residue theorem, we deduce the same bound for each fixed Taylor coefficient of $\tilde{I}(s)$ at $s=0$.  In particular,
  \begin{equation}
    Z_\infty(0)
    L(\pi,\tfrac{1}{2})^n \ll
    C(\pi_{\fin})^{n/2}
    T^{n/4  - \min(\kappa/2,\delta) + \eps (1 + \kappa) }.
  \end{equation}
  Since this estimate holds for each fixed $\eps > 0$, it holds also with $\eps (1 + \kappa)$ replaced by $o(1)$ (Lemma \ref{lem:overspill-A-vs-T-eps}).  By combining with the lower bound for $Z_\infty(0)$ noted previously, we deduce the claim \eqref{eq:lpi-tfrac12n-ll_eps}.
\end{proof}

By Lemma \ref{lem:sub-gln:fix-kappa-in}, our task reduces to estimating the period integrals
\begin{equation*}
  \int_{[H]} \varphi \Psi_{\ast}[Y]
\end{equation*}
for $Y \in A_H(\mathbb{R})^0$ with $T^{-\kappa} \leq Y_j \leq T^\kappa$, where $\kappa \in (0,1/2)$ is a fixed parameter.  Since we focus in what follows on individual values of $Y$, we abbreviate \index{Eisenstein series!$\Psi = \Psi[Y]$}
\begin{equation*}
  \Psi_{\ast} := \Psi_{\ast}[Y],
  \quad
  \Psi := \Psi[Y].
\end{equation*}


\subsection{Matrix coefficients of $\Psi$}\label{sec:matr-coeff-our}
Let $\mathcal{H}_H$ denote the spherical Hecke algebra for $H(\mathbb{Z})$, as in \S\ref{sec:matr-coeff-bounds}.
\begin{lemma}\label{lem:sub-gln:each-y-in}
  For each nonnegative $t \in \mathcal{H}_H$ and $h \in H(\mathbb{R})$, we have
  \begin{equation}
    \left\langle R(h) t \Psi, \Psi  \right\rangle \ll \lambda_0(t) T^{o(1)},
  \end{equation}
  and similarly for $\Psi_{\ast}$.
\end{lemma}
\begin{proof}
  Lemma \ref{lem:standard2:let-d_1-d_2} and assertion \eqref{itm:sub-gln:7} of Theorem \ref{thm:main-local-results} give, for fixed $d_1, d_2 \in \mathbb{Z}_{\geq 0}$,
  \begin{equation*}
    \mathcal{S}_{d_1,d_2}(f) \ll T^{o(1)}.
  \end{equation*}
  The conclusion then follows from Lemma \ref{lem:standard2:matrix-coeff-bounds}; the latter is applicable because, by construction, $f \in \mathcal{S}^e(U_H(\mathbb{R}) \backslash H(\mathbb{R}))^{W_H}$.
\end{proof}

\begin{corollary}\label{cor:standard2:we-have-psi}
  We have the $L^2([H])$ estimate $\|\Psi \| \ll T^{o(1)}$, and similarly for $\Psi_{\ast}$.
\end{corollary}
\begin{proof}
  Take $(t,h) = (1,1)$.
\end{proof}

\subsection{Introduction of an approximate idempotent}\label{sec:intr-an-appr}

\subsubsection{Preparatory lemmas}\label{sec:cuhjo7ygb4}
We record some preliminary lemmas, to be applied below in \S\ref{sec:cuhjo7z3jd}.

We define as in \cite[\S2.1.7]{michel-2009} a height function $\htt : [G] \rightarrow \mathbb{R}_{\geq 1}$.  It is essentially the norm $\|.\|_{[\bar{G}]}$, defined as in \S\ref{sec:schwartz-spaces} using the adjoint representation of $\bar{G}$.

\begin{lemma}\label{lem:standard:ht-neg-power-integrable-over-H}
  There exists $d_0 \geq 0$ so that $\int_{[H]} \htt(h)^{-d_0} \, d h < \infty$.
\end{lemma}
\begin{proof}
  We may apply Proposition \ref{prop:siegel-upper} to the general linear group $H$ to majorize the integral over $[H]$ by an integral over a Siegel domain for $H$.  On such a domain, the adjoint norm of $t \in A_H(\mathbb{R})^0_{\geq 1}$ with respect to $G$, and hence the height $\htt(t g)$ for $g$ in a fixed compact set, grows faster than some fixed positive power of $\lVert t \rVert$.  On the other hand, the Haar measure has density $d t/\delta_{N_H}(t)$.  Taking $d_0$ large enough yields the claimed integrability.
\end{proof}

We record a crude pointwise decay estimate for cuspidal automorphic functions, expressed using the adelic Sobolev norms of Michel--Venkatesh \cite[\S2]{michel-2009}.

Let \(\mathbf{X}\) denote either \([G]=G(\mathbb{Q})\backslash G(\mathbb{A})\), or the twisted space \(([G],\omega)\), where \(\omega:Z(\mathbb{Q})\backslash Z(\mathbb{A})\to \mathbb{C}^\times\) is a unitary character.  Define \(C^\infty(\mathbf{X})\) to be the space of functions \(f:[G]\to \mathbb{C}\) that are smooth in the \(G(\mathbb{R})\)-variable and right-invariant under some open compact subgroup of \(G(\mathbb{A}_{\fin})\) and, in the twisted case, satisfy \(f(zg)=\omega(z)f(g)\) for all \(z\in Z(\mathbb{A})\).  We write \(\|\cdot\|\) for the Petersson norm (possibly infinite) over $[G]$ or, in the twisted case, its quotient by the center.

For \(N\in \mathbb{Z}_{\ge 1}\), let
\begin{equation*}
  K(N) = \prod_{p} \ker \left( G(\mathbb{Z}_p) \rightarrow G(\mathbb{Z}/N\mathbb{Z}) \right) \leq  G(\mathbb{A}_{\fin})
\end{equation*}
denote the principal congruence subgroup of level \(N\).

\begin{lemma}\label{lem:standard:crude-height-decay}
  For each $d \geq 0$, there exist $u \in \mathfrak{U}(G(\mathbb{R}))$, $d' \geq 0$ and $C \geq 1$ such that for each $N \geq 1$ and cuspidal $f \in C^\infty(\mathbf{X})$ that is right-invariant under $K(N)$, we have for all $x \in [G]$ the bound
  \begin{equation}\label{eq:standard:crude-height-decay}
    |f(x)| \le C \htt(x)^{-d} N^{d'} \|R(u)f\|.
  \end{equation}
\end{lemma}
\begin{proof}
  By \cite[\S2.6.2, (2.14)]{michel-2009}, there exist $d_0 \geq 0$ and $C_0 \geq 1$ such that for every cuspidal $F \in C^\infty(\mathbf{X})$,
  \begin{equation}\label{eq:mv-pointwise-height}
    |F(x)| \le C_0\, \htt(x)^{-d}\, \mathcal{S}_{d_0}(F)
    \qquad \text{ for } x\in [G],
  \end{equation}
  where $\mathcal{S}_{d_0}$ denotes the $L^2$-Sobolev norm of \cite[\S2.3.4]{michel-2009} (possibly infinite).

  We apply \eqref{eq:mv-pointwise-height} to $F=f$.  By the construction of $\mathcal{S}_{d_0}$ in \cite[\S2.3.3--\S2.3.4]{michel-2009}, the finite-adelic part of the corresponding adelic operator acts on right $K(N)$-invariant functions with operator norm $\ll N$.  Therefore
  \begin{equation}\label{eq:sobolev-level-bound}
    \mathcal{S}_{d_0}(f) \ll N^{d_0}\,\|(1+\Delta_\infty)^{d_0} f\|,
  \end{equation}
  for a fixed elliptic operator $\Delta_\infty$ on $G(\mathbb{R})$, which may be given by $\sum_i (1 - X_i^2)$ for a basis $\{X_i\}_i$ of the Lie algebra of $G(\mathbb{R})$.  We may find $u\in \mathfrak{U}(G(\mathbb{R}))$ so that $(1+\Delta_\infty)^{d_0} f = R(u)f$.  Combining \eqref{eq:mv-pointwise-height} and \eqref{eq:sobolev-level-bound} then yields \eqref{eq:standard:crude-height-decay} (with $d' := d_0$).
\end{proof}

\subsubsection{Construction of the approximate idempotent}\label{sec:cuhjo7z3jd}

Let $\omega_{0,\infty} \in C_c^\infty(G(\mathbb{R}))$ be the element produced by Theorem \ref{thm:main-local-results} (denoted there by ``$\omega_0$'').  Recall that it is supported in each fixed neighborhood of the identity element.  In particular, it is supported in \emph{some} fixed compact subset  $E_{G(\mathbb{R})} \subseteq G(\mathbb{R})$.

Let $p$ be a prime number.  For each $q \in \{1, p, p^2, \dotsc\}$, we may define the congruence subgroup $K_0(q)_p \leq G(\mathbb{Z}_p)$ consisting of matrices with bottom row congruent to $(0,\dotsc,0,\ast)$ modulo $q$.  Let
\begin{equation*}
  d : K_0(q)_p \rightarrow \mathbb{Z}_p^\times / (\mathbb{Z}_p^\times \cap (1 + q \mathbb{Z}_p))
\end{equation*}
denote the homomorphism assigning to a matrix the class of its lower-right entry.  Let $C(\pi_p) \in \{1, p, p^2, \dotsc \}$ denote the conductor of the component $\pi_p$.  Let $\pi_p|_{Z} : \mathbb{Q}_p^\times \rightarrow \U(1)$ denote the central character of $\pi_p$.  Let $\omega_{0,p} \in C_c^\infty(G(\mathbb{Q}_p))$ denote the function supported on  $K_0(C(\pi_p))_p$ and given there by the following $L^2$-normalized character:
\begin{equation*}
  g \mapsto \vol(K_0(C(\pi_p))_p)^{-1/2} \pi_p|_{Z}(d(g))^{-1}.
\end{equation*}
Note that this character is well-defined due to the fact the conductor of $\pi_p|_{Z}$ divides that of $\pi_p$ (\S\ref{sec:essential-vector}).  The significance of this normalization is that the self-convolution of $\omega_{0,p}$ is bounded in magnitude by one (see Lemma \ref{lem:self-conv-omeg} below).

Let $\omega_0 = \otimes_{\mathfrak{p}} \omega_{0,\mathfrak{p}} \in C_c^\infty(G(\mathbb{A}))$ denote the tensor product of the functions defined above.

We have $C(\pi_{\fin}) = \prod_p C(\pi_p)$.  Set $K_0(\pi) := \prod_p K_0(C(\pi_p))_p$.  Then
\begin{equation}\label{eq:volk_0pi-=-cpiinfty}
  \vol(K_0(\pi)) = C(\pi_{\fin})^{-n+o(1)}.
\end{equation}

\begin{lemma}\label{lem:standard:we-have-begin}
  We have
  \begin{equation}\label{eq:volk_0pi-int_h-varph}
    \vol(K_0(\pi))^{1/2}
    \int_{[H]} \varphi \cdot \Psi_{\ast}
    =
    \int_{[H]} \pi(\omega_0) \varphi \cdot \Psi_{\ast}
    + \O(T^{-\infty}).
  \end{equation}
\end{lemma}
\begin{proof}
  Abbreviate $c := \vol(K_0(\pi))^{1/2} = \prod_p c_p$, with $c_p := \vol(K_0(C(\pi_p))_p)^{1/2}$.  By Corollary \ref{cor:standard2:we-have-psi}, we have $\|\Psi\| \ll T^{o(1)}$.  By Cauchy--Schwarz, our task reduces to verifying that
  \begin{equation*}
    \int_{[H]} | c \varphi - \pi(\omega_0) \varphi |^2 \ll T^{-\infty}.
  \end{equation*}

  Choose $d$ large enough that $\htt(\cdot)^{-2d}$ is integrable over $[H]$ (Lemma \ref{lem:standard:ht-neg-power-integrable-over-H}).  Apply Lemma \ref{lem:standard:crude-height-decay} to the cuspidal function
  \begin{equation*}
    \varphi' := c \varphi - \pi(\omega_0) \varphi,
  \end{equation*}
  which is right-invariant under $K(C(\pi_{\fin}))$ (see \S\ref{sec:essential-vector}).  We obtain an estimate of the form
  \begin{equation*}
    |\varphi'(x)| \ll \htt(x)^{-d} C(\pi_{\fin})^{d'} \|\pi(u)\varphi'\|
  \end{equation*}
  for all $x \in [G]$.  Recalling our assumption \eqref{eq:n_pi-=-to1} concerning $C(\pi_{\fin})$, we reduce to verifying that for each fixed $u \in \mathfrak{U}(G(\mathbb{R}))$, we have
  \begin{equation}\label{eq:c-varphi-piomega_0}
    \|\pi(u) (c \varphi - \pi(\omega_0) \varphi) \| \ll T^{-\infty}.
  \end{equation}

  By construction, we have $\pi(\omega_{0,p}) W_p = c_p W_p$ for each finite prime $p$.  Thus $c \varphi - \pi(\omega_0) \varphi$ is the automorphic form attached, as in \S\ref{sec:norms-cusp-forms}, to the archimedean Whittaker function $W_\infty - \pi(\omega_{0,\infty}) W_\infty$.  By the results of \S\ref{sec:norms-cusp-forms} (applied as in \S\ref{sec:constr-autom-forms}, using our assumption \eqref{eq:n_pi-=-to1}), we reduce further to verifying that for each fixed $u \in \mathfrak{U}(G(\mathbb{R}))$, we have
  \begin{equation*}
    \|\pi(u) (W_\infty - \pi(\omega_{0,\infty}) W_\infty) \| \ll T^{-\infty},
  \end{equation*}
  where now $\|.\|$ denotes the norm defined using the Kirillov model.  This last estimate is the content of part \eqref{itm:sub-gln:1} of Theorem \ref{thm:main-local-results}.
\end{proof}

We thereby reduce to estimating the integrals
\begin{equation*}
  \int_{[H]} \pi(\omega_0) \varphi \cdot \Psi_{\ast}.
\end{equation*}

\subsection{Introduction of an amplifier}\label{sec:intr-an-ampl}
Let $S$ be the finite set of places of $\mathbb{Q}$ consisting of the infinite place $\infty$ together with any finite primes at which $\pi$ is ramified.  Let $\mathcal{H}_G^S$ denote the restricted product over $p \notin S$ of the spherical Hecke algebras for $(G(\mathbb{Q}_p), G(\mathbb{Z}_p))$.

The algebra $\mathcal{H}_G^S$ is generated by the elements $T_p(a)$, defined for primes $p \notin S$ and $(n+1)$-tuples of integers $a = (a_1,\dotsc,a_{n+1}) \in \mathbb{Z}^{n+1}$, given by the characteristic functions of the double cosets
\begin{equation*}
  G(\mathbb{Z}_p) \diag(p^{a_1}, \dotsc, p^{a_{n+1}}) G(\mathbb{Z}_p) \subseteq G(\mathbb{Q}_p)
\end{equation*}
(see for instance \cite[\S5.2.1]{2020arXiv201202187N}).  For $j \geq 0$, we introduce the shorthand $T_p[j] := T_p ((j,0,\dotsc,0))$ as well as its normalized variant
\begin{equation*}
  t_{p,j} := p^{ \frac{n j}{2} } T_p[j].
\end{equation*}

Let $L > 1$ with $\log L / \log T$ fixed and positive (to be determined later, in \S\ref{sec:optimization}).  Let $\mathbb{P}$ denote the set of all primes $p \notin S$ that lie in the interval $[L,2 L]$.  By the divisor bound and the prime number theorem, we then have
\begin{equation*}
  T^{-o(1)} L \leq |\mathbb{P}| \leq L.
\end{equation*}

Let $\lambda_\pi : \mathcal{H}_G^S \rightarrow \mathbb{C}$ denote the algebra homomorphism \index{representations!Hecke algebra homomorphism $\lambda_\pi$} describing the eigenvalues for the action of $\mathcal{H}_G^S$ on the $\prod_{p \notin S} G(\mathbb{Z}_p)$-invariant subspace of $\pi$.  By \cite[(5.6)]{2020arXiv201202187N} (which relies in turn on \cite[Cor 4.3]{2014arXiv1405.6691B}), we have $\sum_{j=1}^{n+1} \lvert \lambda_{\pi}(t_{p,j}) \rvert \gg 1$ for each $p \notin S$.  By the pigeonhole principle, we may thus find $j_0 \in \{1, \dotsc, n+1\}$ so that
\begin{equation*}
  \sum_{p \in \mathbb{P}} \lvert \lambda_{\pi}(t_{p,j_0}) \rvert \gg |\mathbb{P}|.
\end{equation*}

Write $\omega_{0,S} := \otimes_{p \in S} \omega_{0,p}$.  Then $\omega_0 = \omega_{0,S} \otimes (\otimes_{p \notin S} \omega_{0,p})$, with $\omega_{0,p}$ unramified for each $p \notin S$.  For $p \in \mathbb{P}$, set $x_p := \overline{\sgn(\lambda_\pi(t_{p,j_0}))}$, where $\sgn(0) := 0$ and $\sgn(z) := z/|z|$ for $z \neq 0$.  Then $|x_p| \leq 1$, and
\begin{equation*}
  \pi \left( \omega_{0,S} \otimes \sum_{p \in \mathbb{P} } x_p t_{p,j_0} \right) \varphi
  =
  \left( \sum_{p \in \mathbb{P} } |\lambda_\pi(t_{p,j_0}) | \right)
  \pi(\omega_0) \varphi.
\end{equation*}

By combining the above identities and estimates with \eqref{eq:volk_0pi-=-cpiinfty} and \eqref{eq:volk_0pi-int_h-varph}, we deduce the following implication:
\begin{align}
  &\mathcal{R} := \left\lvert
    \int_{[H]}
    \pi \left( \omega_{0,S} \otimes \sum_{p \in \mathbb{P} } x_p t_{p,j_0} \right) \varphi
    \cdot \Psi_{\ast}
    \right\rvert^2
    \ll
    T^{n/2 - 2 \delta } L^2 \label{eq:leftlvert-int-_h}
  \\
  &\quad \implies
    \int_{[H]} \varphi \Psi_{\ast} \ll
    C(\pi_{\fin})^{n/2}
    T^{n/4 - \delta + o(1)}. \nonumber
\end{align}
Our task thereby reduces to suitably estimating the quantity $\mathcal{R}$.

\subsection{Pretrace inequality and multiplication of Hecke operators}\label{sec:pretr-ineq-mult-1}
Following \cite[\S6.3]{2020arXiv201202187N} line-by-line, we apply the pretrace inequality given in \cite[Lemma 5.5]{2020arXiv201202187N} (which explicitly avoids requiring compactness for the quotients $[G]$ and $[H]$) and multiply the Hecke operators arising in the geometric expansion.  We summarize the estimate obtained in this way in the lemma stated below.
\begin{notation}
  Let $Z$ denote the center of $G$ and $\bar{G} = G/Z$ the adjoint group.  (We will consider $\bar{G}(R)$ only when $R$ is either $\mathbb{Q}$, $\mathbb{A}$ or $\mathbb{Z}_p$; in each case, $\bar{G}(R) = G(R) / Z(R)$.)  Let $\pi|_Z : Z(\mathbb{A}) \rightarrow \U(1)$ denote the central character of $\pi$.  For $\omega \in C_c^\infty(G(\mathbb{A}))$, define
  \begin{equation*}
    \omega^\sharp(g) := \int_{z \in Z(\mathbb{A})} \pi|_Z(z) \omega(z g) \, d z
  \end{equation*}
  and
  \begin{equation*}
    \mathcal{R}(\omega) := \int_{x, y \in [H]} \overline{\Psi_{\ast}(x)} \Psi_{\ast}(y)
    \sum_{\gamma \in \bar{G}(\mathbb{Q})} \omega^\sharp(x^{-1} \gamma y) \, d x \, d y.
  \end{equation*}
\end{notation}
\begin{notation}\label{notation:convolve-adjoint-local-global-pretrace}
  For elements $\omega, \omega_1, \omega_2$ of $C_c^\infty(G(\mathbb{A}))$ or of $C_c^\infty(G(\mathbb{Q}_\mathfrak{p}))$, write $\omega^*(g) := \overline{\omega(g^{-1})}$ for the adjoint and $\omega_1 \ast \omega_2$ for the convolution.
\end{notation}
\begin{notation}\label{sec:notation-pretr-ineq-mult-3}
  For $p_1, p_2 \in \mathbb{P}$ and $\mathfrak{j} \in \mathbb{Z}_{\geq 0}$, let $\omega[p_1, p_2, \mathfrak{j}]$ denote the factorizable test function $\omega = \otimes \omega_\mathfrak{p} \in C_c^\infty(G(\mathbb{A}))$ with local components described as follows.
  \begin{itemize}
  \item For $p \in S$, take $\omega_p := \omega_{0,p} \ast \omega_{0,p}^*$.
  \item For $p \notin S \cup  \{p_1, p_2\}$, take for $\omega_p$ the characteristic function of $G(\mathbb{Z}_p)$.
  \item For $p_1 \neq p_2$, take
    \begin{equation}\label{eq:omega_p_1-:=-t_p_1}
      \omega_{p_1} :=  t_{p_1,\mathfrak{j}}
      =
      p_1 ^{- \frac{n \mathfrak{j}}{2}} T_{p_1} (\mathfrak{j},0,\dotsc,0),
    \end{equation}
    \begin{equation}\label{eq:omega_p_2-:=-t_p_2}
      \omega_{p_2} :=  t_{p_2,\mathfrak{j}}^* = p_2 ^{- \frac{n \mathfrak{j}}{2}} T_{p_2} (0, \dotsc, 0, -\mathfrak{j}).
    \end{equation}
  \item For $p_1 = p_2$, take
    \begin{equation}\label{eq:omeg-=-omeg}
      \omega_{p_1} = \omega_{p_2} := p_1^{- n \mathfrak{j}} T_{p_1}(\mathfrak{j},0,\dotsc,0,-\mathfrak{j}).
    \end{equation}
  \end{itemize}
  By analogy to the adelic case, we write $\pi|_{Z_\mathfrak{p}}$ for the central character of $\pi_\mathfrak{p}$ and set
  \begin{equation}\label{eq:omeg-=-int}
    \omega_\mathfrak{p}^\sharp(g) = \int_{z \in Z(\mathbb{Q}_\mathfrak{p})} \pi|_{Z_\mathfrak{p}}(z) \omega_\mathfrak{p}(z g) \, d z.
  \end{equation}

\end{notation}
\begin{lemma}
  The quantity $\mathcal{R}$ defined in \eqref{eq:leftlvert-int-_h} satisfies the estimate
  \begin{equation}\label{eq:mathcalr-ll-sum}
    \mathcal{R} \ll \sum_{p_1 \neq p_2 \in \mathbb{P} } \sum_{\mathfrak{j}=1}^{n+1}
    \left\lvert \mathcal{R}(\omega[p_1, p_2, \mathfrak{j}]) \right\rvert
    + \sum_{p_1 = p_2 \in \mathbb{P} } \sum_{\mathfrak{j}=0}^{n+1}
    \left\lvert \mathcal{R} (\omega[p_1, p_2, \mathfrak{j}]) \right\rvert.
  \end{equation}
\end{lemma}
\begin{proof}
  This follows from the arguments of \cite[\S6.3]{2020arXiv201202187N}, applied verbatim.
\end{proof}

We henceforth focus on individual values of $\mathfrak{j}, p_1, p_2$ appearing on the RHS of \eqref{eq:mathcalr-ll-sum}, abbreviate
\begin{equation*}
  \omega := \omega[p_1,p_2,\mathfrak{j}],
\end{equation*}
and aim to estimate $\mathcal{R}(\omega)$.  As a first step, we decompose
\begin{equation}\label{eq:mathc-=-mathc}
  \mathcal{R}(\omega) = \mathcal{M} + \mathcal{E},
\end{equation}
where $\mathcal{M}$ denotes the contribution from $\gamma$ lying in $H(\mathbb{Q}) \hookrightarrow \bar{G}(\mathbb{Q})$ and $\mathcal{E}$ denotes the remaining contribution.

We will estimate $\mathcal{M}$ in \S\ref{sec:main-term-estimates} and $\mathcal{E}$ in \S\ref{sec:error-term-estimates}.  We then combine the estimates so obtained in \S\ref{sec:optimization}.

We conclude this subsection with some basic estimates for the factors $\omega_\mathfrak{p} ^\sharp$, to be applied below.
\begin{lemma}\label{lem:let-p-notin}
  Let $p \notin S \cup \{p_1, p_2\}$.  Then $\left\lvert \omega_p^\sharp \right\rvert$ is the pullback under $G(\mathbb{Q}_p) \rightarrow \bar{G}(\mathbb{Q}_p)$ of the characteristic function of $\bar{G}(\mathbb{Z}_p)$; in particular,
  \begin{equation}\label{eq:omega_p-sharp-_infty}
    \|\omega_p^\sharp \|_{\infty} = 1, \quad \int_{H(\mathbb{Q}_p)} |\omega_p^\sharp| = 1.
  \end{equation}
\end{lemma}
\begin{proof}
  Clear by construction.  Indeed, $\bar{G}(\mathbb{Z}_p)$ is the image of $G(\mathbb{Z}_p)$, and for each $g \in G(\mathbb{Q}_p)$, the set $\{z \in Z(\mathbb{Q}_p) : z g \in G(\mathbb{Z}_p)\}$ is either empty or a coset of $Z(\mathbb{Z}_p)$.  Because $\pi_p$ is unramified and unitary, its central character takes a constant value of magnitude one on each such coset.  Moreover, each such coset has volume one, and the set in question is nonempty precisely when $g \in \bar{G}(\mathbb{Z}_p)$.  The first identity in \eqref{eq:omega_p-sharp-_infty} is then clear.  Since $\bar{G}(\mathbb{Z}_p) \cap H(\mathbb{Q}_p) = H(\mathbb{Z}_p)$ and the Haar measure on $H(\mathbb{Q}_p)$ assigns volume one to $H(\mathbb{Z}_p)$, the second identity in \eqref{eq:omega_p-sharp-_infty} follows.
\end{proof}
\begin{lemma}\label{lem:self-conv-omeg}
  Let $p \in S - \{\infty \}$.
  \begin{enumerate}[(i)]
  \item \label{item:self-conv-omeg-2}  The support of $\omega_p^\sharp$ is contained in $\bar{G}(\mathbb{Z}_p)$.
  \item \label{item:self-conv-omeg-1} We have $\|\omega_{p} ^\sharp\|_{\infty} \leq 1$.
  \item \label{item:self-conv-omeg-3}  We have $\int_{H(\mathbb{Q}_p)} |\omega_p^\sharp| \leq 1$.
  \end{enumerate}
\end{lemma}
\begin{proof}
  We recall from \S\ref{sec:intr-an-appr} the definition of $\omega_{0,p}$.  Set $J := K_0(C(\pi_p))_p \leq K_p \leq G(\mathbb{Q}_p)$.  Let $\chi : J \rightarrow \U(1)$ denote the character of $J$ given by $g \mapsto \pi_p|_{Z}(d(g))^{-1}$.  Let $e \in C_c^\infty(G(\mathbb{Q}_p))$ denote the function supported on $J$ and given there by $\vol(J)^{-1} \chi$.  Then $e^* = e$ and $e \ast e = e$.  Since $\omega_{0,p} = \vol(J)^{1/2} e$, it follows that $\omega_{0,p} \ast \omega_{0,p}^* = \vol(J) e$ is the function $1_J \chi$ supported on $J$ and given there by $\chi$.  Since $J \leq G(\mathbb{Z}_p)$, we deduce that $|\omega_{0,p} \ast \omega_{0,p}^*|$ is bounded by the characteristic function of $G(\mathbb{Z}_p)$.  We may thus bound $|\omega_p^\sharp|$ exactly as in Lemma \ref{lem:let-p-notin}.
\end{proof}
\begin{lemma}\label{lem:we-have-begin}
  We have
  \begin{equation*}
    \prod_{p \in \{p_1, p_2\}}
    \|\omega_p^\sharp \|_{\infty} \ll L^{-n \mathfrak{j}}.
  \end{equation*}
\end{lemma}
\begin{proof}
  The basis elements $T_p(a)$ for the Hecke algebra have the property that for each $g \in G(\mathbb{Q}_p)$, the set of all $z \in Z(\mathbb{Q}_p)$ for which $z g$ lies in the support of $T_p(a)$ is either empty or a coset of $Z(\mathbb{Z}_\mathfrak{p})$.  Thus for $p \in \{p_1, p_2\}$, we have $\|\omega_p^\sharp \|_{\infty } \leq \|\omega_p\|_{\infty}$.  Since also $\|T_p(a)\|_{\infty} = 1$, we deduce from the definitions \eqref{eq:omega_p_1-:=-t_p_1}, \eqref{eq:omega_p_2-:=-t_p_2} and \eqref{eq:omeg-=-omeg} that
  \begin{itemize}
  \item if $p_1 \neq p_2$, then $\|\omega_{p_1} ^\sharp \|_{\infty} = p_1^{- \frac{n \mathfrak{j}}{2}}$ and $\|\omega_{p_2} ^\sharp \|_{\infty} = p_2^{- \frac{n \mathfrak{j}}{2}}$, while
  \item if $p_1 = p_2 =: p$, then $\|\omega_{p} ^\sharp \|_{\infty} = p^{- n \mathfrak{j}}$.
  \end{itemize}
  In either case, the required bound follows from the fact that $p \asymp L$.
\end{proof}

\subsection{Main term estimates}\label{sec:main-term-estimates}
\begin{lemma}\label{lem:sub-gln:main-term-estimates}
  We have
  \begin{equation*}
    \mathcal{M} \ll T^{n/2 + o(1)} L^{-\mathfrak{j}}.
  \end{equation*}
\end{lemma}
\begin{proof}
  The main term $\mathcal{M}$ defined in the previous section unfolds, as follows:
  \begin{align*}
    \mathcal{M}
    &:= \int_{x, y \in [H]} \overline{\Psi_{\ast}(x)} \Psi_{\ast}(y)
      \sum_{\gamma \in H(\mathbb{Q})} \omega^\sharp(x^{-1} \gamma y) \, d x \, d y \\
    &= \int_{x \in [H]} \int_{y \in H(\mathbb{A})}\overline{\Psi_{\ast}(x)} \Psi_{\ast}(y)  \omega^\sharp(x^{-1} y) \, d x \, d y \\
    &= \int_{h \in H(\mathbb{A})}  \omega^\sharp(h) \langle R(h) \Psi_{\ast}, \Psi_{\ast} \rangle \, d h.
  \end{align*}
  Recall that $\omega = \otimes_{\mathfrak{p}} \omega_\mathfrak{p}$.  For $p \notin S$, write $t_p$ for the restriction of $\omega_p^\sharp$ to $H(\mathbb{Q}_p)$, weighted by the Haar measure on $H(\mathbb{Q}_p)$ (which, by convention, assigns volume one to $H(\mathbb{Z}_p)$), so that $t := \otimes_{p \notin S} t_p$ defines a nonnegative element of the spherical Hecke algebra $\mathcal{H}_H$ (see \S\ref{sec:matr-coeff-our}).  Writing $\mathbb{Q}_S := \prod_{\mathfrak{p} \in S} \mathbb{Q}_\mathfrak{p}$, we have
  \begin{equation*}
    \mathcal{M} = \int_{h \in H(\mathbb{Q}_S) } \omega_\infty^\sharp(h) \langle R(h) t \Psi_{\ast}, \Psi_{\ast} \rangle \, d h.
  \end{equation*}
  We recall that, by construction, $\Psi_{\ast}$ is invariant under $\prod_p H(\mathbb{Z}_p)$.  By Lemmas \ref{lem:sub-gln:each-y-in} and \ref{lem:let-p-notin}, it follows that
  \begin{equation*}
    \mathcal{M} \ll \lambda_0(t) T^{o(1)} \int_{H(\mathbb{R})} |\omega_\infty^\sharp|
    \prod_{p \in S - \{\infty\}}
    \int_{H(\mathbb{Q}_p)} |\omega_p^\sharp|.
  \end{equation*}

  By construction (i.e., Theorem \ref{thm:main-local-results}), we have $\int_{H(\mathbb{R})} |\omega_\infty^\sharp| \ll T^{n/2+o(1)}$.

  By Lemma \ref{lem:self-conv-omeg}, we have
  \begin{equation*}
    \prod_{p \in S} \int_{H(\mathbb{Q}_p)} |\omega_p^\sharp| \leq 1.
  \end{equation*}

  To bound $\lambda_0(t)$, we factor it as $\prod_p \lambda_0(t_p)$ and bound each factor individually.  Recall that $\omega_p = \omega_p[p_1, p_2, \mathfrak{j}]$ and that $p_1, p_2 \in [L,2 L]$.  The discussion of \cite[\S6.4]{2020arXiv201202187N}, specialized to $\vartheta = 0$, applies verbatim,\footnote{
    There is a typo in the cited reference, where the index $k \in \{1, 2\}$ is used in place of $\mathfrak{j}$ in an exponent.
  } giving
  \begin{align*}
    p \notin \{p_1, p_2 \} &\implies \lambda_0(t_p) = 1, \\
    p_1 \neq p_2 &\implies \lambda_0(t_{p_1}), \lambda_0(t_{p_2}) \ll L^{-\mathfrak{j}/2}, \\
    p_1 = p_2 &\implies \lambda_0(t_{p_1}) \ll L^{-\mathfrak{j}}.
  \end{align*}
  The claimed bound follows in either case.
\end{proof}

\subsection{Matrix counting}\label{sec:matrix-counting}
We pause to record some estimates to be applied in the following subsection.  We note that the quantitative details of this subsection are relevant only for explicating the precise numerical savings $\delta_n^\sharp$ obtained in our results.  The reader who is primarily interested in understanding the proof of a qualitative subconvex bound is thus encouraged to skim this subsection.

Recall that $(G,H) = (\GL_{n+1},\GL_n)$ over $\mathbb{Z}$, with $n \geq 1$ fixed.  For $t \in A_H(\mathbb{R})^0$, we set
\begin{equation*}
  t^\dagger  := \prod_{i=1}^n \max(t_i, t_i^{-1}).
\end{equation*}
Recall that $T \ggg 1$ and $\bar{G} = G/Z$.

\begin{definition}\label{defn:standard2:let-ell-ell}
  Let $\ell$ and $\ell '$ be natural numbers with $\ell \asymp \ell '$ and $\log \ell / \log T \asymp 1$.  In particular, $\ell, \ell ' \ggg 1$.  Let $t, u \in A_H(\mathbb{R})^0_{\geq 1}$.  We regard these as included in $G(\mathbb{R})$, with $t_{n+1}, u_{n+1} := 1$; we also write simply $t$ and $u$ for their images in $\bar{G}(\mathbb{R})$ when doing so introduces no confusion.  Let $0 < X \leq 1$.  Let $\mathcal{D}$ be a subset of $\bar{G}(\mathbb{R})$.  Let
  \begin{equation*}
    \Sigma(t, u, \ell, \ell', X, \mathcal{D})
  \end{equation*}
  denote the set of all $\gamma \in \bar{G}(\mathbb{Q})$ with the following properties:
  \begin{enumerate}[(i)]
  \item \label{itm:scratch-research:Sigma-count-1}
    There is a lift $\tilde{\gamma} \in G(\mathbb{Q})$ of $\gamma$ with integral entries and determinant $\pm \ell$.
  \item \label{itm:scratch-research:Sigma-count-2}
    There is a lift $\tilde{\gamma}^{-1} \in G(\mathbb{Q})$ of $\gamma^{-1}$ with integral entries and determinant $\pm \ell '$.
  \item \label{itm:scratch-research:Sigma-count-3} We have $t^{-1} \gamma u \in \mathcal{D}$.
  \item \label{itm:scratch-research:Sigma-count-4}  We have $d_H(t^{-1} \gamma u) \leq X$.
  \item \label{itm:scratch-research:Sigma-count-5} We have $\gamma \notin H(\mathbb{Q})$,  regarding the latter as a subgroup of $\bar{G}(\mathbb{Q})$.
  \end{enumerate}
\end{definition}

For the following lemmas, we retain the setting of Definition \ref{defn:standard2:let-ell-ell}, we fix a compact symmetric (i.e., inversion-invariant) neighborhood $E \subseteq \bar{G}(\mathbb{R})$ of the identity element, we fix a compact symmetric subset $\Omega \subseteq H(\mathbb{R})$, and we abbreviate
\begin{equation*}
  \Sigma := \Sigma(t,u,\ell, \ell', X, \Omega E \Omega).
\end{equation*}

\begin{lemma}[Crude \emph{a priori} bound]\label{lem:standard2:crude-a-priori-counting}
  We have
  \begin{equation*}
    |\Sigma| \ll
    \left( \ell t^\dagger u^\dagger \max(|\det(t^{-1} u)|, |\det(u^{-1} t)|) \right)^{n+1}.
  \end{equation*}
\end{lemma}
\begin{proof}
  We consider the case $|\det(t^{-1} u)| \geq 1$; the other case can be treated by applying the same argument to the inverses of all involved matrices and recalling that $\ell ' \asymp \ell$.

  Set $D := |\det(t)^{-1} \ell \det(u)|$.  By hypothesis,
  \begin{equation}\label{eq:d-geq-ell}
    D \geq \ell \geq 1.
  \end{equation}
  Our task is to verify that
  \begin{equation*}
    |\Sigma| \ll (D t^\dagger u^\dagger)^{n+1}.
  \end{equation*}

  Let $\gamma \in \Sigma$.  Choose a lift $\tilde{\gamma}$ as in condition \eqref{itm:scratch-research:Sigma-count-1}.  The rescaled matrix
  \begin{equation*}
    \gamma^{(1)} := D^{- \frac{1}{n+1}} \tilde{\gamma }
  \end{equation*}
  then satisfies
  \begin{equation*}
    \det(t^{-1} \gamma^{(1)} u) = \pm 1.
  \end{equation*}
  The image of $t^{-1} \gamma^{(1)} u$ in $\bar{G}(\mathbb{R})$ is $t^{-1} \gamma u$, which, by condition \eqref{itm:scratch-research:Sigma-count-3}, lies in a fixed compact subset of $\bar{G}(\mathbb{R})$.  Thus $t^{-1} \gamma^{(1)} u$ lies in a fixed compact subset of $G(\mathbb{R})$.  Let $a_{i j}$ ($1 \leq i, j \leq n+1$) denote the matrix entries of $\gamma^{(1)}$.  We have $a_{i j} \in D^{- \frac{1}{n+1}} \mathbb{Z}$ and
  \begin{equation}\label{eq:r_i-1-s_j}
    t_i^{-1} u_j a_{i j} \ll 1
    \quad
    (1 \leq i, j \leq n+1),
  \end{equation}
  where, as noted above, $t_{n+1}, u_{n+1} := 1$.

  For each $\gamma$, there are exactly two choices for $\tilde{\gamma}$, differing by $\pm 1$.  For this reason, it suffices to bound the number of possibilities for $\gamma^{(1)}$, or equivalently, for the coefficients $a_{i j}$.  Thus
  \begin{equation}\label{eq:sigma-ll-prod_i}
    |\Sigma| \ll
    \prod_{i,j= 1}^{n+1} \left( 1 + D^{\frac{1}{ n + 1}} t_i u_j^{-1} \right).
  \end{equation}
  We conclude by estimating crudely via the inequality
  \begin{equation*}
    1 + D^{\frac{1}{ n + 1}} t_i u_j^{-1}
    \leq
    D^{\frac{1}{ n + 1}} ( 1 + t_i ) ( 1 + u_j^{-1} ),
  \end{equation*}
  which follows from \eqref{eq:d-geq-ell}, and the inequalities
  \begin{equation*}
    \prod_{i=1}^{n+1} (1 + t_i) \leq 2^{n+1} t^\dagger ,
    \quad
    \prod_{j=1}^{n+1} (1 + u_j^{-1}) \leq 2^{n+1} u^\dagger.
  \end{equation*}
\end{proof}

\begin{lemma}\label{lem:scratch-research:if-sigma-nonempty}
  Assume that
  \begin{equation}\label{eq:dett-1-u}
    \det(t^{-1} u) = T^{o(1)}.
  \end{equation}
  If $\Sigma$ is nonempty, then
  \begin{equation}\label{eq:d--frac2n+1-0}
    \ell^{-\frac{2}{n+1} - o(1)} \ll t_i / u_i \ll \ell^{\frac{2}{n+1}+o(1)},
  \end{equation}
\end{lemma}
\begin{proof}
  The proof will show more precisely that if there exists $\gamma \in \bar{G}(\mathbb{Q})$ satisfying conditions \eqref{itm:scratch-research:Sigma-count-1} and \eqref{itm:scratch-research:Sigma-count-2}, then the indicated estimate holds.

  We define $D$, $\tilde{\gamma}$, $\gamma^{(1)}$ and $a_{i j}$ as in the proof of Lemma \ref{lem:standard2:crude-a-priori-counting}.

  By \eqref{eq:dett-1-u}, we have
  \begin{equation}\label{eq:d-=-ell1}
    D = \ell^{1 + o(1)}.
  \end{equation}
  In particular, $D \geq 1$.

  Write $r_1 \geq \dotsb \geq r_{n+1}$ (resp.\ $s_1 \geq \dotsb \geq s_{n+1}$) for the numbers $t_1,\dotsc,t_n,1$ (resp.\ $u_1,\dotsc,u_n,1$) arranged in non-increasing order.  Since $\gamma^{(1)}$ has nonzero determinant, we may find a permutation $\sigma$ of $\{1, \dotsc, n+1\}$ so that $a_{i,\sigma(i)}$ is nonzero for each $i \in \{1, \dotsc, n+1\}$.  Since $a_{i j} \in D^{- \frac{1}{n+1}} \mathbb{Z}$, it follows then from \eqref{eq:r_i-1-s_j} that
  \begin{equation}\label{eq:r_i-1-s_sigmai}
    r_i^{-1} s_{\sigma(i)} \ll D^{\frac{1}{n+1}}.
  \end{equation}
  Since $r$ and $s$ are non-increasing, the same estimate holds with $\sigma(i) = i$, i.e.,
  \begin{equation}\label{eq:r_i-1-s_i}
    r_i^{-1} s_i \ll D^{\frac{1}{n+1}}.
  \end{equation}

  Similarly, setting $D' := |\det(u)^{-1} \ell ' \det(t)| \asymp \ell^{1+o(1)}$ and arguing as above with a lift $\tilde{\gamma}^{-1}$ of $\gamma^{-1}$ as in condition \eqref{itm:scratch-research:Sigma-count-2}, we obtain $s_i^{-1} r_i \ll (D') ^{\frac{1}{n+1}}$.  Combining the estimates derived thus far, we deduce that
  \begin{equation}\label{eq:ell--frac1n+1}
    \ell^{- \frac{1}{n+1} + o(1)} \ll s_i^{-1} r_i \ll \ell^{\frac{1}{n+1} + o(1)}.
  \end{equation}

  We pass from \eqref{eq:ell--frac1n+1} to \eqref{eq:d--frac2n+1-0} by a simple combinatorial argument, which we record for completeness.  For $i \in \{1,\dotsc,n+1\}$, define $a_i,b_i \in \mathbb{R}$ by writing
  \begin{equation*}
    t_i = \ell^{a_i}, \quad u_i = \ell^{b_i}.
  \end{equation*}
  For example, $a_{n+1} = b_{n+1} = 0$.  The desired conclusion \eqref{eq:d--frac2n+1-0} is that for each fixed $c > \frac{1}{n+1}$ and all $i = \{1, \dotsc, n\}$, we have $|a_i - b_i| \leq 2 c$.  On the other hand, our hypothesis \eqref{eq:ell--frac1n+1} implies there is a permutation $\sigma$ of $\{1,\dotsc,n+1\}$ so that $|a_i - b_{\sigma(i)}| \leq c$ for all $i \in \{1, \dotsc, n+1\}$.  The passage from our hypothesis to the desired conclusion is the content of
  Lemma \ref{lem:scratch-research:let-n-geq}, stated at the end of this subsection.
\end{proof}

\begin{lemma}\label{lem:standard2:we-have-begin}
  Assume that \eqref{eq:dett-1-u} holds.  Then
  \begin{equation}\label{eq:sigma_r-sx-ll-0}
    |\Sigma| \ll X \ell^{3 (n+1) + o(1)}
    \min(\delta_H(t) t^\dagger , \delta_H(u) u^\dagger),
  \end{equation}
  where
  \begin{equation*}
    \delta_H(t) := \delta_{N_H}(t) = t_1^{n-1} t_2^{n-3} \dotsb t_n^{1-n}.
  \end{equation*}
\end{lemma}
\begin{proof}
  We may assume that $\Sigma$ is nonempty, so that \eqref{eq:d--frac2n+1-0} holds.  By the symmetry of the problem with respect to $t$ and $u$, it will suffice to show that
  \begin{equation}\label{eq:sigma_r-sx-ll-0-b}
    |\Sigma| \ll X \ell^{3 (n+1) + o(1)}
    \delta_H(t) t^\dagger.
  \end{equation}

  Let $\gamma \in \Sigma$.  We define $D$, $\tilde{\gamma}$, $\gamma^{(1)}$ and $a_{i j}$ as in the proof of Lemma \ref{lem:standard2:crude-a-priori-counting}, and reduce again to counting the possibilities for the $a_{i j}$.

  We consider first the case $X = 1$, in which the condition \eqref{itm:scratch-research:Sigma-count-4} on $d_H$ in the definition of $\Sigma$ is tautological.  We appeal in this case to the estimate \eqref{eq:sigma-ll-prod_i} for $|\Sigma|$, established under weaker hypotheses than the present ones.  Recall from \eqref{eq:d-=-ell1} that $D = \ell^{1 + o(1)}$.  By \eqref{eq:d--frac2n+1-0}, we have $t_i u_j^{-1} \ll \ell^{\frac{2}{n+1} + o(1)} t_i t_j^{-1}$.  Thus, estimating a bit crudely with respect to $\ell$, we obtain
  \begin{equation*}
    1 + D^{\frac{1}{n + 1}} t_i u_j^{-1}
    \ll
    1 +
    \ell^{\frac{3}{n+1} + o(1)} t_i t_j^{-1}
    \leq
    \ell^{\frac{3}{n+1} + o(1)}(1 +
    t_i t_j^{-1}),
  \end{equation*}
  and so
  \begin{equation*}
    |\Sigma| \ll \ell^{3 (n+1) + o(1)} \prod_{i,j=1}^{n+1} (1 + t_i t_j^{-1}).
  \end{equation*}
  By estimating the contribution to this last product from indices grouped according to whether they are equal to $n+1$ or not, we see that it is $\ll \delta_H(t) t^\dagger$.  Indeed:
  \begin{itemize}
  \item Since $t_{n+1} = 1$, the contribution from when $i=n+1$ or $j=n+1$ is $2 \prod_{j=1}^{n} (1 + t_j^{-1}) (1 + t_j) \asymp t^\dagger$.
  \item For $i,j \leq n$ with $i \geq j$, we have $1 + t_i t_j^{-1} \leq 2$ in view of the $N$-dominance of $t$.
  \item For $i,j \leq n$ with $i < j$, we have $t_i t_j^{-1} \geq 1$, hence $1 + t_i t_j^{-1} \asymp t_i t_j^{-1}$; we conclude by noting that $\prod_{i < j \leq n} t_i t_j^{-1} = \delta_H(t)$.
  \end{itemize}
  Thus
  \begin{equation}\label{eq:sigma-ll-ell3}
    |\Sigma| \ll \ell^{3 (n+1) + o(1)}
    \delta_H(t) t^\dagger,
  \end{equation}
  giving the required estimate \eqref{eq:sigma_r-sx-ll-0-b}.

  It remains to address the case $X < 1$.  By the definition of $d_H(t^{-1} \gamma u)$, the lower-right entries of $t^{-1} \gamma u$ and its inverse are nonzero.  In particular, $a_{n+1,n+1} \neq 0$, and
  \begin{equation}\label{eq:d_hr-1-gamma}
    X \geq d_H(t^{-1} \gamma u) \asymp
    \max_{1 \leq i \leq n}
    t_i^{-1} \frac{|a_{i,n+1}|}{ | a_{n+1, n+1}|}
    +
    \max_{1 \leq j \leq n}
    u_j \frac{|a_{n+1, j}|}{ | a_{n+1, n+1}|}
    + \dotsb,
  \end{equation}
  where $\dotsb$ denotes the contribution of the analogous quantities defined in terms of the inverse matrix $u^{-1} \gamma^{-1} t$.  By \eqref{eq:r_i-1-s_j}, we have $a_{n+1,n+1} \ll 1$, and similarly for the corresponding entry of the inverse.  Condition \eqref{itm:scratch-research:Sigma-count-4} thus gives
  \begin{equation}\label{eq:a_j-n+1-ll}
    a_{j,n+1} \ll X t_j, \quad
    a_{n+1,j} \ll X u_j^{-1}.
  \end{equation}
  Since
  \begin{itemize}
  \item each of the entries $a_{j,n+1}$ or $a_{n+1,j}$ with $j \in \{1,\dotsc,n\}$ lies in $D^{\frac{-1}{n+1}} \mathbb{Z}$, and
  \item not all such entries vanish (thanks to condition \eqref{itm:scratch-research:Sigma-count-5}),
  \end{itemize}
  we deduce from \eqref{eq:a_j-n+1-ll} that
  \begin{equation}\label{eq:x-max_1-leq}
    X \ell^{\frac{1}{n+1}} \max_{1 \leq j \leq n} \max(t_j, u_j^{-1}) \gg \ell^{ - o(1)}.
  \end{equation}
  We now count as in the case $X = 1$.  In that case, the number of possibilities for $a_{j,n+1}$ and $a_{n+1,j}$ ($1 \leq j \leq n$) was majorized, up to a factor of $\ell^{o(1)}$, by
  \begin{equation*}
    C_1 := \prod_{1 \leq j \leq n } (1 + \ell^{\frac{1}{n+1}} t_j)  (1 + \ell^{\frac{1}{n+1}} u_j^{-1}).
  \end{equation*}
  In the present case $X < 1$, we obtain instead the modified count
  \begin{equation}\label{eq:c_x-:=-prod_1}
    C_X := \prod_{1 \leq j \leq n } (1 + X \ell^{\frac{1}{n+1}} t_j)  (1 + X \ell^{\frac{1}{n+1}} u_j^{-1}).
  \end{equation}

  Using \eqref{eq:x-max_1-leq},  we see that
  \begin{equation}\label{eq:c_x-ll-x}
    C_X \ll X  C_1.
  \end{equation}
  Indeed, let $j_0$ be an index at which the maximum in \eqref{eq:x-max_1-leq} is attained.  For the factor in \eqref{eq:c_x-:=-prod_1} indexed by $j = j_0$, we have by \eqref{eq:x-max_1-leq} that
  \begin{equation*}
    \frac{
      (1 + X \ell^{\frac{1}{n+1}} t_j)  (1 + X \ell^{\frac{1}{n+1}} u_j^{-1})
    }{
      (1 + \ell^{\frac{1}{n+1}} t_j)  (1 + \ell^{\frac{1}{n+1}} u_j^{-1})
    }
    \ll X,
  \end{equation*}
  while for the remaining indices $j \neq j_0$, we estimate the corresponding factor in \eqref{eq:c_x-:=-prod_1} trivially using that $X < 1$.

  Counting the remaining entries ($a_{i j}$ with $1 \leq i, j \leq n$ and $a_{n+1,n+1}$) as in the case $X = 1$ yields a bound for $|\Sigma|$ that saves an additional factor of $X$ over the bound \eqref{eq:sigma-ll-ell3} obtained in the case $X = 1$.  We thereby obtain the required bound \eqref{eq:sigma_r-sx-ll-0-b}.
\end{proof}

We conclude \S\ref{sec:matrix-counting} by recording the combinatorial argument required to complete the proof of Lemma \ref{lem:scratch-research:if-sigma-nonempty}.  The following discussion will not otherwise be applied in this paper.
\begin{definition}
  Given $n \in \mathbb{Z}_{\geq 0}$ and $n$-tuples $a = (a_1,\dotsc,a_n) \in \mathbb{R}^n$ and $b = (b_1,\dotsc,b_n)\in \mathbb{R}^n$, we define
  \begin{equation*}
    d(a,b) := \min_{\sigma \in S(n)}
    \max_{1 \leq i \leq n}
    |a_i - b_{\sigma(i)}|.
  \end{equation*}
\end{definition}
\begin{lemma}\label{lem:manhattan-greedy}
  If $a,b \in \mathbb{R}^n$ are non-decreasing, then
  \begin{equation*}
    d(a,b) = \max_{1 \leq i \leq n} |a_i - b_i|.
  \end{equation*}
\end{lemma}
\begin{proof}
  Write $d_\sigma(a, b) := \max_{1 \leq i \leq n} \lvert a_i - b_{\sigma(i)} \rvert$.  We must show that $d_{\mathrm{id}}(a,b) \leq d_\sigma(a,b)$ for each $\sigma \in S(n)$.  This is clear if $\sigma = \mathrm{id}$, so we may assume otherwise.  We may then find $i < j$ with $\sigma(i) > \sigma(j)$.  Let $\sigma' \in S(n)$ be obtained from $\sigma$ by swapping $\sigma(i)$ and $\sigma(j)$, but leaving all other values unchanged.  The length of $\sigma'$ is less than that of $\sigma$, so by induction on that length, it will be enough to show that $d_{\sigma'}(a,b) \leq d_{\sigma}(a,b)$.  In verifying this, we may reduce to the case of just two indices $i$ and $j$, since all other indices are unaffected by the swap from $\sigma$ to $\sigma'$.  Suppose then that $n=2$ and that $\sigma$ is the transposition swapping $1$ and $2$.  We must show that
  \begin{equation}\label{eq:cuhjotetbn}
    \max(\lvert a_1 - b_1 \rvert, \lvert a_2 - b_2 \rvert) \leq \max(\lvert a_1 - b_2 \rvert, \lvert a_2 - b_1 \rvert).
  \end{equation}
  Suppose for instance that the maximum on the left hand side is attained by the first quantity (the other case is similar).  If $a_1 \leq b_1$, then, since $b_2 \geq b_1$, we have $\lvert a_1 - b_1 \rvert \leq \lvert a_1 - b_2 \rvert$, so \eqref{eq:cuhjotetbn} holds.  Similarly, if $a_1 \geq b_1$, then, since $a_2 \geq a_1$, we have $\lvert a_1 - b_1 \rvert \leq \lvert a_2 - b_1 \rvert$, so \eqref{eq:cuhjotetbn} holds again.  This completes the proof.
\end{proof}
\begin{definition}
  For $a \in \mathbb{R}^{n+1}$ with $n \in \mathbb{Z}_{\geq 0}$, and $k \in \{1,\dotsc, n+1\}$, define
  \begin{equation*}
    a^{(k)} := (a_1,\dotsc,a_{k-1},a_{k+1},\dotsc,a_{n}) \in \mathbb{R}^{n}.
  \end{equation*}
\end{definition}
\begin{lemma}\label{lem:scratch-research:let-n-geq}
  Let $n \in \mathbb{Z}_{\geq 0}$, $a,b \in \mathbb{R}^{n+1}$, $k,\ell \in \{1,\dotsc,n+1\}$, $c \geq 0$.  Suppose that
  \begin{equation*}
    d (a,b) \leq c, \quad
    a_k = b_{\ell}.
  \end{equation*}
  Then
  \begin{equation*}
    d (a^{(k)}, b^{(\ell)}) \leq 2 c.
  \end{equation*}
\end{lemma}
\begin{proof}
  We may assume that $a$ and $b$ are non-decreasing, so that by Lemma \ref{lem:manhattan-greedy}, we have
  \begin{equation}\label{eq:a_i-b_i-leq}
    |a_i - b_i| \leq c
  \end{equation}
  for all $i \in \{1,\dotsc,n+1\}$.  It is enough to show then that
  \begin{equation}\label{eq:ak_i-aell_i-leq}
    |a^{(k)}_i - b^{(\ell)}_i| \leq 2 c
  \end{equation}
  for all $i \in \{1,\dotsc,n\}$.  Without loss of generality, suppose $k \leq \ell$.  For $i < k$ (resp.\ $i \geq \ell$), we have $a^{(k)}_i = a_i$ and $b^{(\ell)}_i = b_i$ (resp.\ $a^{(k)}_i = a_{i+1}$ and $b^{(\ell)}_i = b_{i+1}$), so \eqref{eq:ak_i-aell_i-leq} follows from \eqref{eq:a_i-b_i-leq}.  It remains to treat the range $k \leq i < \ell$, in which $a_i^{(k)} = a_{i+1}$ and $b_i^{(\ell)} = b_i$.  By assumption, we have
  \begin{equation*}
    a_k - c \leq b_k \leq \dotsb \leq b_{\ell} = a_k \leq \dotsb \leq a_{\ell} \leq b_{\ell} + c = a_k + c.
  \end{equation*}
  Therefore all of the elements $a_k,\dotsc,a_{\ell}, b_k,\dotsc,b_{\ell}$ lie in the interval $[a_k-c, a_k+c]$ of length $2 c$, and so
  \begin{equation*}
    |a_i^{(k)} - b_i^{(\ell)}| = |a_{i+1} - b_i| \leq 2 c,
  \end{equation*}
  as required.
\end{proof}

\subsection{Error term estimates}\label{sec:error-term-estimates}
By definition,
\begin{equation}\label{eq:mathcale-:=E--int}
  \mathcal{E}
  := \int_{x, y \in [H]} \overline{\Psi_{\ast}(x)} \Psi_{\ast}(y)
  \sum_{\gamma \in \bar{G}(\mathbb{Q}) - H(\mathbb{Q})} \omega^\sharp(x^{-1} \gamma y) \, d x \, d y.
\end{equation}
We recall the notation:
\begin{itemize}
\item $\Psi_{\ast}$ is shorthand for $\Psi_{\ast}[Y]$, the pseudo Eisenstein series attached in \S\ref{sec:constr-autom-forms} to some $Y \in A_H(\mathbb{R})^0$ satisfying $T^{-\kappa} \leq Y_j \leq T^{\kappa}$ for some fixed $\kappa \in (0,1/2)$.
\item $\omega^\sharp \in C_c^\infty(\bar{G}(\mathbb{A}))$, $\bar{G} = G/Z$, is the integral over $Z(\mathbb{A})$, weighted by the central character of $\pi$, of the test function $\omega \in C_c(G(\mathbb{A}))$ attached to some triple $(p_1,p_2,\mathfrak{j})$ (see \S\ref{sec:pretr-ineq-mult-1}), which is in turn an ``amplified'' self-convolution of the test function $\omega_{0,\infty}$ produced by Theorem \ref{thm:main-local-results}.
\end{itemize}
The main result of \S\ref{sec:error-term-estimates} is as follows.
\begin{proposition}\label{prop:standard2:with-notat-assumpt}
  With notation and assumptions as above,
  \begin{equation*}
    \mathcal{E} \ll
    T^{n/2 - 1/2 + n \kappa + o(1)}
    L^{(3 (n+1)^2 + n) \mathfrak{j}}.
  \end{equation*}
\end{proposition}
The proof will be broken up into several steps, and occupies the remainder of \S\ref{sec:error-term-estimates}.  The crucial ingredients are
\begin{itemize}
\item the ``bilinear forms'' or ``transversality'' estimate supplied by part \eqref{itm:sub-gln:3} of Theorem \ref{thm:main-local-results},
\item the ``growth bounds for Eisenstein series'' estimate supplied by Theorem \ref{thm:growth-eisenstein-nonstandard}, and
\item the simple estimates for matrix counts given above.
\end{itemize}

\subsubsection{Siegel domains, triangle inequality}
We begin by bounding the integrand in absolute value and applying the Siegel domain majorization afforded by Proposition \ref{prop:siegel-upper} to both variables.  Recall that $A_H(\mathbb{R})^0_{\geq 1}$ denotes the set of $N$-dominant elements of the connected component of the diagonal subgroup $A_H(\mathbb{R})$.  Each factor in the integrand of the integral \eqref{eq:mathcale-:=E--int} defining $\mathcal{E}$ is $\prod_{p} H(\mathbb{Z}_p)$-invariant in both variables.  We deduce that there is a fixed compact subset $\Omega$ of $H(\mathbb{R})$ so that
\begin{equation}\label{eq:mathcale-ll-int}
  \mathcal{E} \ll \int_{t, u \in A_H(\mathbb{R})^0_{\geq 1} } \sum_{\gamma \in \bar{G}(\mathbb{Q}) - H(\mathbb{Q})} I_{t, u}(\gamma)
  \, \frac{d t \, d u}{\delta_H(t) \delta_H(u)},
\end{equation}
where
\begin{equation*}
  I_{t,u}(\gamma) :=
  \int_{x \in t \Omega}
  \int_{y \in u \Omega}
  \left\lvert \Psi_{\ast}(x) \Psi_{\ast}(y)
    \omega^\sharp (x^{-1} \gamma y)
  \right\rvert
  \,d x \,d y,
\end{equation*}
where $d x$ and $d y$ denote our chosen Haar measure on $H(\mathbb{R})$.

In what follows, and unless otherwise stated, $t$ and $u$ denote elements of $A_H(\mathbb{R})^0_{\geq 1}$.

\subsubsection{Factorization}
By part \eqref{itm:sub-gln:5} of Theorem \ref{thm:main-local-results} and the $H(\mathbb{R})$-equivariance of the association $f \mapsto \Psi$, we may write
\begin{equation}\label{eq:psiy-=-psi_0}
  \Psi_{\ast}(y) = (\Psi \ast \phi_{Z_H})(y) := \int_{z \in Z_H(\mathbb{R})} \Psi(y z^{-1}) \phi_{Z_H}(z) \, d z,
\end{equation}
where
\begin{itemize}
\item $\phi_{Z_H} \in C_c^\infty(Z_H(\mathbb{R}))$ is fixed,
\item $\Psi : [H] \rightarrow \mathbb{C}$, attached as in \S\ref{sec:compl-eisenst-seri} to the element $f \in \mathcal{S}^e(U_H(\mathbb{R}) \backslash H(\mathbb{R}))^{W_H}$ furnished by Theorem \ref{thm:main-local-results} via the recipe
  \begin{equation*}
    \Psi := \Eis[L(Y) \Phi[f]] = \Eis[\Phi[L(Y) f]] \in \mathcal{S}([H]),
  \end{equation*}
  has the same essential properties as $\Psi_{\ast}$ (indeed, all except the explicit factorization property \eqref{eq:psiy-=-psi_0}); in particular, it is smooth and of rapid decay.
\end{itemize}

\subsubsection{Application of transversality}
Here we apply the ``bilinear forms'' estimate afforded by part \eqref{itm:sub-gln:3} of Theorem \ref{thm:main-local-results}, whose proof relies in turn on the volume estimates established in \cite[\S15--16]{2020arXiv201202187N}.

Recall that part \eqref{itm:sub-gln:3} of Theorem \ref{thm:main-local-results} asserts ``for each fixed compact $\Omega_{H} \subseteq H(\mathbb{R})$, there exists a compact $\Omega_{H}' \subseteq H(\mathbb{R})$ so that $\dotsb$.''  Let $\Omega ' \subseteq H(\mathbb{R})$ denote the set obtained as ``$\Omega_{H}'$'' by taking for ``$\Omega_H$'' our given set $\Omega$.  By enlarging $\Omega '$ if needed, we may and shall assume that $\Omega '$ is symmetric (i.e., invariant by inversion).  For $t \in A_H(\mathbb{R})^0_{\geq 1}$, we define the seminorms $\|.\|_t$ on $L^2(H(\mathbb{Z}) \backslash H(\mathbb{R}))$ as in \eqref{eq:varphi-2_r-:=}, but using $\Omega '$, i.e.,
\begin{equation*}
  \|\phi\|^2_t := \int_{t \Omega' } |\phi|^2
  =
  \int_{h \in \Omega' \subseteq H(\mathbb{R})}
  |\phi(t h)|^2
  \, d h,
\end{equation*}
where the integral is taken with respect to the Haar measure on $H(\mathbb{R})$.

\begin{lemma}\label{lem:sub-gln:each-gamma-in}
  For each $\gamma \in \bar{G}(\mathbb{Q}) - H(\mathbb{Q})$, we have
  \begin{equation*}
    I_{t, u}(\gamma) \ll T^{n/2 + o(1)} L^{-n \mathfrak{j}}
    \min \left( 1, \frac{T^{-1/2}}{d_H(t^{-1} \gamma u)} \right)
    \|\Psi \|_t \|\Psi \|_{u}.
  \end{equation*}
\end{lemma}
\begin{proof}
  For $x \in t \Omega$ and $y \in u \Omega$, we have
  \begin{equation*}
    \omega^\sharp (x^{-1} \gamma y)
    =
    \omega_\infty^\sharp (x^{-1} \gamma y)
    \prod_p \omega_p^\sharp (\gamma).
  \end{equation*}
  By Lemmas \ref{lem:let-p-notin}, \ref{lem:self-conv-omeg} and \ref{lem:we-have-begin}, we have
  \begin{equation*}
    \prod_p \omega_p^\sharp (\gamma) \ll L^{-n \mathfrak{j}}.
  \end{equation*}
  Thus
  \begin{equation*}
    I_{t, u}(\gamma) \ll L^{-n \mathfrak{j}} \int_{x \in t \Omega } \int_{y \in u \Omega } \left\lvert \Psi(x) \Psi(y)
      \omega_{\infty}^{\sharp}(x^{-1} \gamma y) \right\rvert \,d x \,d y.
  \end{equation*}
  Setting $\Psi_1(x) := \Psi(t x)$ and $\Psi_2(y) := \Psi(u y)$, we may rewrite this last double integral as
  \begin{equation*}
    \int_{x ,y \in \Omega } \left\lvert (\Psi_1 \ast \phi_{Z_H})(x) (\Psi_2 \ast \phi_{Z_H})(y) \omega_\infty^\sharp(x^{-1} t^{-1} \gamma u y) \right\rvert \,d x \,d y
  \end{equation*}
  and then estimate it, using part \eqref{itm:sub-gln:3} of Theorem \ref{thm:main-local-results}, as
  \begin{equation*}
    \ll
    T^{n/2 + o(1)}
    \min \left( 1, \frac{T^{-1/2}}{d_H(t^{-1} \gamma u)} \right)
    \|\Psi_1\|_{L^2(\Omega')}
    \|\Psi_2\|_{L^2(\Omega')}.
  \end{equation*}
  By construction of the norms $\|.\|_t, \|.\|_{u}$, we have
  \begin{equation*}
    \|\Psi_1\|_{L^2(\Omega ')} = \|\Psi\|_t, \quad
    \|\Psi_2\|_{L^2(\Omega ')} = \|\Psi\|_{u}.
  \end{equation*}
  The required bound follows.
\end{proof}

For $t, u \in A_H(\mathbb{R})^0_{\geq 1}$, set
\begin{equation*}
  \Sigma_{t,u} :=
  \left\{ \gamma \in \bar{G}(\mathbb{Q}) - H(\mathbb{Q}) : \text{$\omega^\sharp(x^{-1} \gamma y) \neq 0$ for some $(x,y) \in t \Omega \times u \Omega $} \right\}.
\end{equation*}
For $0 < X \leq 1$, let $\Sigma_{t,u}^X \subseteq \Sigma_{t,u}$ denote the subset defined by the further condition $d_H(t^{-1} \gamma u) \leq X$.  For example, $\Sigma_{t,u}^1 = \Sigma_{t,u}$.  Define
\begin{equation}\label{eq:i-:=-int}
  I := \int_{t, u \in A_H(\mathbb{R})^0_{\geq 1}} \|\Psi\|_t \|\Psi\|_u
  \sum_{\gamma \in \Sigma_{t,u}}
  \min \left( T^{1/2}, \frac{1}{d_H(t^{-1} \gamma u)} \right) \, \frac{d t \, d u}{\delta_H(t) \delta_H(u)}.
\end{equation}
From \eqref{eq:mathcale-ll-int} and Lemma \ref{lem:sub-gln:each-gamma-in}, we see that
\begin{equation*}
  \mathcal{E} \ll T^{n/2 - 1/2 + o(1)} L^{- n \mathfrak{j}} I.
\end{equation*}
The proof of Proposition \ref{prop:standard2:with-notat-assumpt} thereby reduces to that of the estimate
\begin{equation}\label{eq:i-ll-tnkappa}
  I \ll
  T^{n \kappa + o(1)}
  L^{(3 (n+1)^2  + 2 n) \mathfrak{j}}.
\end{equation}

\subsubsection{Application of matrix counting}
The sets $\Sigma_{t,u}^X$ are contained in those considered in \S\ref{sec:matrix-counting}.  More precisely, let $E \subseteq \bar{G}(\mathbb{R})$ denote the image of the set $\{x y^{-1} : x,y \in E_{G(\mathbb{R})}\}$, where $E_{G(\mathbb{R})} \subseteq G(\mathbb{R})$ is any fixed symmetric compact subset that contains the support of $\omega_0$.
\begin{lemma}
  There are natural numbers $\ell, \ell' \in \mathbb{Z}_{\geq 1}$ with
  \begin{equation}\label{eq:ell-asymp-ell}
    \ell \asymp \ell ' \asymp L^{(n+1) \mathfrak{j}}
  \end{equation}
  so that
  \begin{equation}\label{eq:sigma_t-ux-subseteq}
    \Sigma_{t,u}^X \subseteq \Sigma(t,u,\ell, \ell', X, \Omega E \Omega).
  \end{equation}
\end{lemma}
\begin{proof}
  We will see that the required conclusion holds with
  \begin{equation*}
    \ell := \prod_p \ell_p, \quad
    \ell_p := \begin{cases}
      1 & \text{ if } p \notin \{p_1, p_2\}, \\
      p^{\mathfrak{j}} & \text{ if } p = p_1 \neq p_2, \\
      p^{n \mathfrak{j} } &  \text{ if } p = p_2 \neq p_1, \\
      p^{(n + 1) \mathfrak{j} } &  \text{ if } p = p_1 = p_2,
    \end{cases}
  \end{equation*}
  \begin{equation*}
    \ell ' := \prod_p \ell '_p, \quad
    \ell'_p := \begin{cases}
      1 & \text{ if } p \notin \{p_1, p_2\}, \\
      p^{n\mathfrak{j}} & \text{ if } p = p_1 \neq p_2, \\
      p^{\mathfrak{j} } &  \text{ if } p = p_2 \neq p_1, \\
      p^{(n + 1) \mathfrak{j} } &  \text{ if } p = p_1 = p_2.
    \end{cases}
  \end{equation*}
  Since $p_1, p_2 \asymp L$, the estimate \eqref{eq:ell-asymp-ell} holds with this choice.  It remains to verify the containment \eqref{eq:sigma_t-ux-subseteq}.

  Let $\gamma \in \Sigma_{t,u}^X$.  The conditions \eqref{itm:scratch-research:Sigma-count-4} and \eqref{itm:scratch-research:Sigma-count-5} in the definition of $\Sigma(t,u,\ell, \ell', X, \Omega E \Omega)$ follow immediately from the definition of $\Sigma_{t,u}^X$, so our task is to verify the remaining conditions, namely, that
  \begin{equation}\label{eq:t-1-gamma}
    t^{-1} \gamma u \in \Omega  E \Omega
  \end{equation}
  and that $\gamma$ (resp.\ $\gamma^{-1}$) admits a lift to $G(\mathbb{Q})$ with integral entries and determinant $\pm \ell$ (resp.\ $\pm \ell '$).

  By the definition of $\Sigma_{t,u}^X$, we may find $(x,y) \in t \Omega \times u \Omega$ with $\omega^\sharp(x^{-1} \gamma y) \neq 0$, i.e.,
  \begin{equation}\label{eq:omega_infty-sharp-x}
    \omega_\infty^\sharp (x^{-1} \gamma y) \neq 0
  \end{equation}
  and
  \begin{equation}\label{eq:omega_p-sharp-gamma}
    \omega_p^\sharp (\gamma ) \neq 0
  \end{equation}
  for all primes $p$.  The support of $\omega_\mathfrak{p} ^\sharp$ is contained in the image of the support of $\omega_\mathfrak{p}$, which we now discuss case-by-case (recalling from \S\ref{sec:pretr-ineq-mult-1} the construction of $\omega$).

  Since $\omega_{0,\infty}$ is supported in $E_{G(\mathbb{R})}$ and $\omega_{\infty} = \omega_{0,\infty} \ast \omega_{0,\infty}^*$, we deduce that $\omega_\infty^\sharp$ is supported in the set $E$ defined above.  From \eqref{eq:omega_infty-sharp-x} and the symmetry of $\Omega$, the required membership \eqref{eq:t-1-gamma} follows.

  We observe that $\supp(\omega_{p}^\sharp)$ is contained in the image $\Xi_p \subseteq \bar{G}(\mathbb{Q}_p)$ of the set
  \begin{equation*}
    \tilde{\Xi}_p := \left\{ g \in M_{n+1}(\mathbb{Z}_p) : \det(g) \in \ell_p \mathbb{Z}_p^\times  \right\}.
  \end{equation*}
  Indeed:
  \begin{itemize}
  \item If $p \notin \{p_1, p_2\}$, then $\supp(\omega_{\mathfrak{p}}) = \bar{G}(\mathbb{Z}_p)$.
  \item If $p_1 \neq p_2$ and $p = p_1$, then $\supp(\omega_{p_1}^\sharp)$ is the image of the support of $\omega_{p_1} = T_{p_1}(\mathfrak{j},0,\dotsc,0)$.
  \item If $p_1 \neq p_2$ and $p = p_2$, then $\supp(\omega_{p_2}^\sharp)$ is the image of the support of $\omega_{p_2} = T_{p_2}(0,\dotsc,0,-\mathfrak{j})$, or equivalently, that of the central shift $T_{p_2}(\mathfrak{j},\dotsc,\mathfrak{j},0)$.
  \item If $p_1 = p_2 = p$, then $\supp(\omega_p^\sharp)$ is the image of the support of $\omega_p = T_p(\mathfrak{j},0,\dotsc,0,-\mathfrak{j})$, or equivalently, that of the central shift $T_p(2 \mathfrak{j}, \mathfrak{j}, \dotsc, \mathfrak{j}, 0)$.
  \end{itemize}
  To conclude, we need only note that for $a \in \mathbb{Z}_{\geq 0}^{n+1}$, the Hecke operator $T_p(a)$ (see \S\ref{sec:intr-an-ampl}), regarded as measure on $G(\mathbb{Q}_p) = \GL_{n+1}(\mathbb{Q}_p)$, is supported on
  \begin{equation*}
    \{g \in M_{n+1}(\mathbb{Z}_p) : \det(g) \in p^{\sum a_j} \mathbb{Z}_p^\times \}.
  \end{equation*}

  From \eqref{eq:omega_p-sharp-gamma}, we deduce that $\gamma \in \Xi_p$ for all $p$.  Since $\mathbb{Q}$ has class number one, we may find a lift $\tilde{\gamma} \in G(\mathbb{Q})$ of $\gamma$ that lies in $\tilde{\Xi}_p$ for all $p$.  This lift has integral entries and determinant $\pm \ell$.

  By applying the above argument to the inverses of all involved matrices, we see that $\gamma^{-1}$ admits a lift $\tilde{\gamma}^{-1} \in G(\mathbb{Q})$ with integral entries and determinant $\pm \ell '$.

  The proof is now complete.
\end{proof}

By construction, the sets $E$ and $\Omega$ satisfy the assumptions required by \S\ref{sec:matrix-counting}.  The results of that section thus apply here, giving the following.

\begin{lemma}\label{lem:standard2:cool-counts-for-SIgma-t-u}
  Let $t, u \in A_H(\mathbb{R})^0_{\geq 1}$ and $0 < X \leq 1$.  Define $t^\dagger = \prod_{i=1}^n \max(t_i, t_i^{-1})$, as in \S\ref{sec:matrix-counting}.  We have
  \begin{equation}\label{eq:sigma-ll-ln+12-a-priori}
    |\Sigma_{t,u}| \ll
    L^{(n+1)^2 \mathfrak{j}}
    \left( t^\dagger u^\dagger \max(|\det(t^{-1} u)|, |\det(u^{-1} t)|) \right)^{n+1}.
  \end{equation}
  Assume further that
  \begin{equation}\label{eq:dett-1-u-b}
    \det(t^{-1} u) = T^{o(1)}.
  \end{equation}
  Then $\Sigma_{t,u}^X$ is nonempty only if
  \begin{equation}\label{eq:d--frac2n+1}
    L^{-2 \mathfrak{j}  - o(1)} \ll t_i / u_i \ll L^{2 \mathfrak{j} +o(1)},
  \end{equation}
  in which case
  \begin{equation}\label{eq:sigma_r-sx-ll}
    |\Sigma_{t,u}^X| \ll X L^{3(n+1)^2 \mathfrak{j} + o(1)} \min(\delta_H(t) t^\dagger, \delta_H(u) u^\dagger),
  \end{equation}
\end{lemma}

\subsubsection{Application of crude bounds for Eisenstein series}\label{sec:appl-crude-bounds}
We now introduce some \emph{a priori} bounds concerning $\|\Psi\|_t$.  These will permit us to reduce the range of pairs $(t,u)$ to consider.

\begin{lemma}\label{lem:each-fixed-c_0-Psi-crude}
  For each fixed $C_0 \geq 0$, there is a fixed $C_1 \geq 0$ so that for each fixed $C_2 \in \mathbb{R}$, we have for all $g \in [H]$
  \begin{equation*}
    \Psi(g)
    \ll
    \|g\|_{[H]}^{-C_0}
    T^{C_1}
    \left\lvert \det(g Y) \right\rvert^{C_2}.
  \end{equation*}
\end{lemma}
\begin{proof}
  We appeal to Lemma \ref{lem:standard2:crude-growth-bounds-completed-eisenstein-series} (using transfer to restrict quantification over $C_0$ through $\mathcal{D}$ to fixed quantities, and with $-C_2$ playing the role of $C_2$).  Since $\Psi = \Eis[\Phi[f]]$ with $f \in \mathcal{S}^e(U(\mathbb{R}) \backslash G(\mathbb{R}))^W$, we obtain (with quantifiers as indicated)
  \begin{equation*}
    \|g\|_{[H]}^{C_0} \left\lvert \det(g Y) \right\rvert^{-C_2} \left\lvert \Psi(g) \right\rvert \ll T^{C_1} \nu_{\mathcal{D},\ell,T}(f).
  \end{equation*}
  By assertion \eqref{itm:sub-gln:7} of Theorem \ref{thm:main-local-results}, we have $\nu_{\mathcal{D},\ell,T}(f) \ll T^{o(1)}$.  We conclude by replacing $C_1$ by $C_1 + 1$, say, to absorb the factor $T^{o(1)}$.
\end{proof}

\begin{lemma}\label{lem:standard2:each-fixed-c_0}
  For each fixed $C_0 \geq 0$, there is a fixed $C_1 \geq 0$ so that for each fixed $C_2 \in \mathbb{R}$, we have for all $t \in A_H(\mathbb{R})^0_{\geq 1}$
  \begin{equation}\label{eq:psi-_t-ll}
    \|\Psi \|_t \ll (t^\dagger )^{-C_0} T^{C_1} |\det(t Y)|^{C_2}.
  \end{equation}
\end{lemma}
\begin{proof}
  This is a consequence of Lemma \ref{lem:each-fixed-c_0-Psi-crude} and the estimate $\log(3 + \|t u\|_{[H]}) \asymp \log(3 + t^\dagger)$ for $u \in \Omega$, which follows (in sharper forms) from a well-known theorem of Siegel; see for instance \cite[Lemma 2.1]{MR1001613} or \cite[A.1.1 (viii)]{2016arXiv160206538B}.
\end{proof}

\begin{lemma}\label{lem:standard2:each-fixed-c}
  For each fixed $C \geq 0$, there is a subset $\mathcal{C}$ of $A_H(\mathbb{R})^0_{\geq 1}$ with the following properties.
  \begin{enumerate}[(i)]
  \item For each $t \in \mathcal{C}$, we have
    \begin{equation}\label{eq:detr-=-dety}
      \det(t Y)  = T^{o(1)}, \quad t^\dagger \ll T^{\O(1)}.
    \end{equation}
  \item The integral $I$ defined in \eqref{eq:i-:=-int} satisfies the estimate
    \begin{equation}\label{eq:mathcale-ll-tn2-restricted-to-C}
      I \ll
      \int_{t, u  \in \mathcal{C} }
      \|\Psi\|_{t}
      \|\Psi\|_{u}
      \sum_{\gamma \in \Sigma_{t,u}}
      \min \left( T^{1/2}, \frac{1}{d_H(t^{-1} \gamma u)} \right)
      \, \frac{d t \, d u}{\delta_H(t) \delta_H(u)}
      + T^{-C}.
    \end{equation}
  \end{enumerate}
\end{lemma}
\begin{proof}
  By \eqref{eq:sigma-ll-ln+12-a-priori}, we have
  \begin{equation}\label{eq:sum-_gamma-in-2}
    \sum_{\gamma \in \Sigma_{t,u}}
    \min \left( T^{1/2}, \frac{1}{d_H(t^{-1} \gamma u)} \right)
    \leq T^{1/2} |\Sigma_{t,u}| \ll
    T^{1/2}
    L^{(n+1)^2 \mathfrak{j}} \mathcal{G}(t,u),
  \end{equation}
  where we abbreviate
  \begin{equation*}
    \mathcal{G}(t,u)  := \left( t^\dagger u^\dagger \max(|\det(t^{-1} u)|, |\det(u^{-1} t)|) \right)^{n+1}.
  \end{equation*}
  We may fix $C_0 \geq 0$ (large enough in terms of $n$) for which
  \begin{equation}\label{eq:int-_t-u}
    \int_{t, u \in A_H(\mathbb{R})^0_{\geq 1}}
    \mathcal{G}(t,u) (t^\dagger u^\dagger )^{1- C_0} \, \frac{d t \, d u}{\delta_H(t) \delta_H(u)} < \infty.
  \end{equation}
  We then fix $C_1$ according to the conclusion of Lemma \ref{lem:standard2:each-fixed-c_0}.

  For $\eps \in (0,1)$, define
  \begin{equation*}
    \mathcal{C}(\eps) := \left\{ t \in A_H(\mathbb{R})^0_{\geq 1} : \left\lvert \det (t Y ) \right\rvert \in [T^{- \eps }, T^\eps ]\right\}.
  \end{equation*}
  When $\eps$ is fixed, the conclusion of Lemma \ref{lem:standard2:each-fixed-c_0} implies that
  \begin{equation}\label{eq:t-notin-mathcalceps}
    t \notin \mathcal{C}(\eps)
    \implies
    \|\Psi\|_t
    \leq
    (t^\dagger )^{-C_0} T^{-1/\eps}.
  \end{equation}
  To see this, we apply \eqref{eq:psi-_t-ll} with $\pm C_2$ fixed but sufficiently large in terms of $C_1$ and $\eps $, using that $T \ggg 1$ to clean up the implied constant.  By overspill, we may find some $\eps \lll 1$ for which this implication \eqref{eq:t-notin-mathcalceps} remains valid.  We work henceforth with such a value of $\eps$.  For all $u \in A_H(\mathbb{R})^0_{\geq 1}$, we see by taking $C_2 = 0$ in  \eqref{eq:psi-_t-ll} that
  \begin{equation}\label{eq:psi_0-_u-ll}
    \|\Psi\|_u \ll (u^\dagger )^{-C_0} T^{C_1}.
  \end{equation}
  Combining \eqref{eq:sum-_gamma-in-2}, \eqref{eq:int-_t-u}, \eqref{eq:t-notin-mathcalceps} and \eqref{eq:psi_0-_u-ll}, we obtain
  \begin{equation*}
    \int_{
      \substack{
        t, u \in A_H(\mathbb{R})^0_{\geq 1} :  \\
        t \notin \mathcal{C}(\eps)
      }
    }
    \|\Psi\|_t \|\Psi\|_u
    \sum_{\gamma \in \Sigma_{t,u}}
    \min \left( T^{1/2}, \frac{1}{d_H(t^{-1} \gamma u)} \right) \, \frac{d t \, d u}{\delta_H(t) \delta_H(u)}
    \ll T^{-\infty}.
  \end{equation*}
  A similar estimate applies with the roles of $t$ and $u$ reversed.  Thus by restricting the integral in the definition \eqref{eq:i-:=-int} to $t,u \in \mathcal{C}(\eps)$, we introduce an additive error of $\O(T^{-\infty})$.

  Next, for each $m \geq 0$, set $\mathcal{C}(\eps,m) := \left\{ t \in \mathcal{C}(\eps) : t^\dagger \leq T^m \right\}$.  If $t \notin \mathcal{C}(\eps,m)$, then it follows again from Lemma \ref{lem:standard2:each-fixed-c_0} that
  \begin{equation*}
    \|\Psi\|_t \ll (t^\dagger )^{-C_0} T^{C_1} \leq (t^\dagger ) ^{1 - C_0} T^{C_1 - m}.
  \end{equation*}
  By fixing $m$ large enough, we see as in the previous paragraph that the contribution to $I$ from $t \notin \mathcal{C}(\eps,m)$ or $u \notin \mathcal{C}(\eps,m)$ is $\O(T^{-C})$.  The conclusion of the lemma then holds with $\mathcal{C} := \mathcal{C}(\eps,m)$.
\end{proof}

\subsubsection{Clean-up}
Since an additive error $T^{-C}$, with $C$ large enough but fixed, is acceptable for the purposes of proving \eqref{eq:i-ll-tnkappa}, we reduce to showing that the estimate \eqref{eq:i-ll-tnkappa} holds for the modification $I_{\mathcal{C}}$ of $I$ obtained by restricting the integral over $t,u$ to some set $\mathcal{C}$ as in the conclusion of Lemma \ref{lem:standard2:each-fixed-c}:
\begin{equation*}
  I_{\mathcal{C}} :=
  \int_{
    t, u \in \mathcal{C}
  }
  \|\Psi\|_t \|\Psi\|_u
  \sum_{\gamma \in \Sigma_{t,u}}
  \min \left( T^{1/2}, \frac{1}{d_H(t^{-1} \gamma u)} \right) \, \frac{d t \, d u}{\delta_H(t) \delta_H(u)}.
\end{equation*}
We henceforth focus on one such set $\mathcal{C}$.  By Cauchy--Schwarz, we have $I \leq I_0$, where
\begin{equation*}
  I_0
  := \int_{t, u \in \mathcal{C} } \|\Psi\|_t^2
  \sum_{\gamma \in \Sigma_{t,u}}
  \min \left( T^{1/2}, \frac{1}{d_H(t^{-1} \gamma u)} \right)
  \, \frac{d t \, d u}{\delta_H(t) \delta_H(u)}.
\end{equation*}

\begin{lemma}
  We have
  \begin{equation}\label{eq:i_1-ll-l}
    I_0 \ll
    L^{(3(n+1)^2 + 2n) \mathfrak{j} + o(1)}
    \int_{t \in \mathcal{C} } \|\Psi\|_t^2
    t^\dagger  \, \frac{d t}{\delta_H(t)}.
  \end{equation}
\end{lemma}
\begin{proof}
  For $t,u \in \mathcal{C}$, the hypothesis \eqref{eq:dett-1-u-b} is satisfied (indeed, $\det(t)$ and $\det(u)$ are both of the form $\det(Y)^{-1} T^{o(1)}$).  We may thus appeal to the refined conclusions of Lemma \ref{lem:standard2:cool-counts-for-SIgma-t-u}.  First of all, set
  \begin{equation*}
    \mathcal{C}_2 := \left\{ (t,u) \in \mathcal{C} \times \mathcal{C} : \Sigma_{t,u} \neq \emptyset
    \right\}.
  \end{equation*}
  Then, by \eqref{eq:d--frac2n+1}, we have
  \begin{equation*}
    (t,u) \in \mathcal{C}_2 \implies
    L^{-2 \mathfrak{j} - o(1)} \ll t_i / u_i \ll L^{2 \mathfrak{j} + o(1)} \text{ for all } i \in \{1,\dotsc,n\}.
  \end{equation*}
  In particular,
  \begin{itemize}
  \item for each $t \in \mathcal{C}$, the set of all $u \in \mathcal{C}$ with $(t,u) \in \mathcal{C}_2$ has volume $\ll (\log L)^{\O(1)} \ll L^{o(1)}$, and
  \item we have $u^\dagger \leq L^{2 n \mathfrak{j} + o(1)} t^\dagger$.
  \end{itemize}
  It follows that
  \begin{equation}\label{eq:int-_u-in}
    \int_{u \in \mathcal{C} : (t,u) \in \mathcal{C}_2 }
    u^\dagger
    \, d u
    \ll  L^{2 n \mathfrak{j}  + o(1)} t^\dagger.
  \end{equation}

  Recall that the subset $\Sigma_{t,u}^X \subseteq \Sigma_{t,u}$ is cut out by the condition $d_H(t^{-1} \gamma u) \leq X$.  By dyadic decomposition, we have
  \begin{equation*}
    \sum_{\gamma \in \Sigma_{t,u}}
    \min \left( T^{1/2}, \frac{1}{d_H(t^{-1} \gamma u)} \right) \ll  \sum_{
      \substack{
        X \in \exp(\mathbb{Z}) :  \\
        T^{-1/2} < X \leq 1
      }
    }
    \frac{| \Sigma_{t,u}^X|}{X}
    +
    \frac{|\Sigma_{t,u}^{X_0}|}{X_0}
  \end{equation*}
  with $X_0 := 100 T^{-1/2}$, say.  We estimate the cardinalities $|\Sigma_{t,u}^X|$ using \eqref{eq:sigma_r-sx-ll}, giving
  \begin{equation*}
    \frac{|\Sigma_{t,u}^X|}{X} \ll L^{3(n+1)^2 \mathfrak{j} + o(1)} \min(\delta_H(t) t^\dagger, \delta_H(u) u^\dagger).
  \end{equation*}
  There are $\O(\log T)$ summands, giving an overall estimate
  \begin{equation*}
    \sum_{\gamma \in \Sigma_{t,u}}
    \min \left( T^{1/2}, \frac{1}{d_H(t^{-1} \gamma u)} \right)
    \ll L^{3(n+1)^2\mathfrak{j} + o(1)} \min(\delta_H(t) t^\dagger, \delta_H(u) u^\dagger).
  \end{equation*}
  We take the second alternative in the minimum and insert the resulting estimate into the definition of $I_0$, giving
  \begin{equation*}
    I_0
    \ll
    L^{3(n+1)^2 \mathfrak{j} + o(1)}
    \int_{(t, u) \in \mathcal{C}_2 } \|\Psi\|_t^2
    \delta_H(u) u^\dagger \, \frac{d t \, d u}{\delta_H(t) \, \delta_H(u)}.
  \end{equation*}
  Estimating the $u$-integral using \eqref{eq:int-_u-in}, we arrive at the desired estimate \eqref{eq:i_1-ll-l}.
\end{proof}

The proof of our earlier reduction \eqref{eq:i-ll-tnkappa} of  Proposition \ref{prop:standard2:with-notat-assumpt} thereby reduces to that of the estimate
\begin{equation}\label{eq:int-_t-in}
  \int_{t \in \mathcal{C} } \|\Psi\|_t^2
  t^\dagger  \, \frac{d t}{\delta_H(t)}
  \ll
  T^{n \kappa + o(1)}.
\end{equation}

\subsubsection{Application of growth bounds for Eisenstein series}\label{sec:appl-growth-bounds}
We now apply Theorem \ref{thm:growth-eisenstein-nonstandard} to deduce the key estimate \eqref{eq:int-_t-in}, following some preliminary lemmas.

\begin{lemma}\label{lem:standard2:we-have-begin-1-growth}
  For each $t \in A_H(\mathbb{R})^0_{\geq 1}$, we have
  \begin{equation}\label{eq:fracpsi-_r2d-ll}
    \frac{\|\Psi\|_t^2}{\delta_H(t)}  \ll T^{o(1)} \min \left( \det(Y)^{-1} t_1^{-n}, \det(Y) t_n^n \right).
  \end{equation}
\end{lemma}
\begin{proof}
  We appeal to Theorem \ref{thm:growth-eisenstein-nonstandard}.  Recall that $\Psi = \Eis[L(Y) \Phi[f]] = \Eis[\Phi[L(Y) f]]$.  The required hypotheses concerning $f$ follow from assertion \eqref{itm:sub-gln:7} of Theorem \ref{thm:main-local-results}.
\end{proof}

\begin{lemma}\label{lem:standard:any-dominant-element}
  For each $t \in A_H(\mathbb{R})^0_{\geq 1}$, we have
  \begin{equation*}
    t^\dagger \leq \max\left( \frac{1}{\det(t)}, \det(t), \frac{t_1^n}{\det(t)}, \frac{\det(t)}{ t_n^{n}} \right).
  \end{equation*}
\end{lemma}
\begin{proof}
  Consider first the case that $t_1 \leq 1$.  Then $t_j \leq 1$ for all $j$, so $t^\dagger = \det(t)^{-1}$.  Similarly, if $t_n \geq 1$, then $t^\dagger = \det(t)$.  In either case, the required inequality holds.  We may thus suppose that $t_1 \geq 1 \geq t_n$.  We may then choose $m \in \{1, \dotsc, n-1\}$ so that $t_1 \geq \dotsb \geq t_m \geq 1 \gg t_{m+1} \geq \dotsb \geq t_n$, in which case
  \begin{equation*}
    t^\dagger = \frac{t_1 \dotsb t_m}{t_{m+1} \dotsb t_n}.
  \end{equation*}
  Since $t_1 \dotsb t_m \cdot t_{m+1} \dotsb t_n = \det(t)$, it follows that
  \begin{equation*}
    t^\dagger = \det(t)^{-1} (t_1 \dotsb t_m ) ^2 = \det(t) (t_{m + 1 } \dotsb t_n ) ^{- 2},
  \end{equation*}
  and so
  \begin{equation*}
    t^\dagger \leq \min(\det(t)^{-1} t_1^{2 m},  \det(t) t_n^{-2(n-m)}).
  \end{equation*}
  If $m \leq n/2$, then we may conclude by taking the first alternative in the minimum, recalling that $t_1 \geq 1$.  If $m \geq n/2$, so that $n-m \leq n/2$, then we may conclude by taking the second alternative, recalling that $t_n \leq 1$.
\end{proof}

In view of our assumption that $Y_j \in [T^{-\kappa}, T^{\kappa}]$, we have $\det(Y) \in [T^{- n \kappa}, T^{n \kappa}]$.  The desired estimate \eqref{eq:int-_t-in} is thus an immediate consequence of the following.
\begin{proposition}
  We have
  \begin{equation*}
    \int_{t \in \mathcal{C} } \|\Psi\|_t^2 t^\dagger \, \frac{d t}{\delta_H(t)}
    \ll
    \max(\det(Y), \det(Y)^{-1}) T^{o(1)}.
  \end{equation*}
\end{proposition}
\begin{proof}
  Lemma \ref{lem:standard:any-dominant-element} gives
  \begin{equation*}
    \int_{t \in \mathcal{C} } \|\Psi\|_t^2 t^\dagger \, \frac{d t}{\delta_H(t)} \ll
    J_+ + J_- + J_0,
  \end{equation*}
  \begin{equation*}
    J_{\pm} := \int_{t \in \mathcal{C} } \|\Psi\|_t^2 \det(t)^{\pm 1} \, \frac{d t}{d_H(t)},
  \end{equation*}
  \begin{equation*}
    J_0 := \int_{t \in \mathcal{C} } \|\Psi\|_t^2 \max(\det(t)^{-1} t_1^n, \det(t) t_n^{-n}) \, \frac{d t}{ \delta_H(t) }.
  \end{equation*}

  Recall from \eqref{eq:detr-=-dety} that for $t \in \mathcal{C}$, we have
  \begin{equation}\label{eq:dett-=-dety}
    \det(t) = \det(Y)^{-1} T^{o(1)}
  \end{equation}
  and
  \begin{equation}\label{eq:t-dagger-ll}
    t^\dagger \ll T^{\O(1)}.
  \end{equation}
  Recall also from Corollary \ref{cor:standard2:we-have-psi}  that $\|\Psi\| \ll T^{o(1)}$.  By Proposition \ref{prop:siegel-lower-integrated}, it follows that
  \begin{equation*}
    \int_{t \in \mathcal{C} } \|\Psi\|_t^2 \, \frac{d t}{ \delta_H(t) } \ll T^{o(1)}.
  \end{equation*}
  Thus
  \begin{equation*}
    J_+ + J_- \ll
    \max(\det(Y), \det(Y)^{-1}) T^{o(1)}.
  \end{equation*}
  It remains only to estimate $J_0$.  By Lemma \ref{lem:standard2:we-have-begin-1-growth} and the estimate \eqref{eq:dett-=-dety}, we have
  \begin{equation*}
    J_0 \ll T^{o(1)} \int_{t \in \mathcal{C}} \, d t.
  \end{equation*}
  By \eqref{eq:t-dagger-ll}, we see that $\int_{t \in \mathcal{C}} \, d t \ll (\log T)^{\O(1)} \ll T^{o(1)}$.  The proof is now complete.
\end{proof}

\subsection{Optimization}\label{sec:optimization}
We now optimize (as in \cite[\S6.6]{2020arXiv201202187N}).  The results of \S\ref{sec:main-term-estimates} and \S\ref{sec:error-term-estimates} combine to give
\begin{equation*}
  \mathcal{M} + \mathcal{E} \ll
  T^{n/2 + o(1)} (
  L^{-\mathfrak{j}}+ T^{-1/2 + n \kappa} L^{(3(n+1)^2 + n) \mathfrak{j}})
\end{equation*}
We insert this bound into the decompositions \eqref{eq:mathcalr-ll-sum} and \eqref{eq:mathc-=-mathc} to see that the quantity $\mathcal{R}$ defined in \eqref{eq:leftlvert-int-_h} satisfies, for fixed $\eps > 0$,
\begin{align*}
  \frac{\mathcal{R}}{L^2 T^{n/2}}
  &\ll T^{o(1)}  \sum_{\mathfrak{j}=1}^{n+1} \left( L^{-\mathfrak{j}} + T^{-1/2 + n \kappa} L^{(3(n+1)^2  + n) \mathfrak{j}}
    \right) \\
  &\quad + T^{o(1)} L^{-1} \sum_{\mathfrak{j} = 0}^{n+1} \left(
    L^{- \mathfrak{j}} + T^{-1/2 + n \kappa}  L^{(3(n+1)^2  + n) \mathfrak{j}}
    \right).
\end{align*}
By the results of \S\ref{sec:reduct-peri-bounds} and \S\ref{sec:intr-an-ampl}, we have the following implication, for fixed $\kappa \in (0, 1/2)$ and $\delta > 0$:
\begin{equation*}
  \frac{\mathcal{R}}{L^2 T^{n/2}} \ll T^{- 2 \delta + o(1)}
  \implies
  L(\pi,\tfrac{1}{2})^n \ll C(\pi_{\fin})^{n/2} T^{n(n+1)/4 - \min(\kappa/2,\delta) + o(1)}.
\end{equation*}
To optimize, we take $\kappa = 2 \delta$, which forces $\delta < 1/4$, and $L = T^{2 \delta}$.  Our estimate for $\mathcal{R}$ then reads
\[
  \frac{\mathcal{R}}{L^2 T^{n/2}}
  \ll
  T^{- 2 \delta + o(1)} ( 1 + T^{\alpha} )
\]
with
\begin{align*}
  \alpha
  &:=
    2 \delta -1/2 + n \kappa
    +
    (3 (n+1)^2 + n)(n+1) \delta  \\
  &=
    (2 + 2 n + (3 (n+1)^2  + n)(n+1)) \delta - 1/2.
\end{align*}
We have $\alpha \leq 0$ provided that
\begin{equation*}
  \delta \leq \frac{1}{2 (2 + 2 n + (3 (n+1)^2 + n)(n+1))},
\end{equation*}
in which case we obtain
\begin{equation*}
  L(\pi,\tfrac{1}{2})^n \ll C(\pi_{\fin})^{n/2} T^{n(n+1)/4 - \delta + o(1)}.
\end{equation*}
Equivalently, for any fixed
\begin{equation*}
  \delta < \delta_{n+1} ^\sharp  := \frac{2}{n(n+1) (2 + 2 n + (3 (n+1)^2 + n)(n+1))}.
\end{equation*}
we have
\begin{equation*}
  L(\pi,\tfrac{1}{2}) \ll C(\pi_{\fin})^{1/2} T^{\frac{n+1}{4}(1 - \delta)}.
\end{equation*}
We verify readily that $\delta_n^\sharp$ simplifies to the fraction indicated in the introduction, i.e., that,
\begin{equation*}
  n(n+1) (2 + 2 n + (3 (n+1)^2 + n)(n+1))
  = 3 ( n+1)^5 - 2 (n+1)^4 - (n+1)^2.
\end{equation*}
\qed

This completes the deduction of Theorem \ref{thm:main-result-reformulated} from Theorems \ref{thm:main-local-results} and \ref{thm:growth-eisenstein-nonstandard}.  The remainder of the paper is devoted to proofs of the latter results.  Theorem \ref{thm:main-local-results} is the subject of Parts \ref{part:constr-test-vect} and \ref{part:asympt-analys-kirill}, while Theorem \ref{thm:growth-eisenstein-nonstandard} is taken up in Part \ref{part:local-l2-growth}.

\part{Construction of test vectors}\label{part:constr-test-vect}

Here we develop the proof of Theorem \ref{thm:main-local-results}.  That result asserts the existence of, among other things, vectors $f \in \mathcal{S}^e(U_H \backslash H)^{W_H}$ and $W \in \mathcal{W}(\pi,\psi)$ having favorable properties.  The main point of Part \ref{part:constr-test-vect} is to construct $f$ and $W$.  The construction of $W$ is given in \S\ref{sec:asympt-analys-local}, that of $f$ in \S\ref{sec:some-analysis-basic}.  Some ingredients and language common to both constructions are given in \S\ref{sec:modul-bump-funct}.

\section{Lie-theoretic preliminaries}\label{sec:modul-bump-funct}

Let $F$ be a fixed completion of $\mathbb{Q}$, thus $F$ is either $\mathbb{R}$ or $\mathbb{Q}_p$ for some fixed prime $p$.  This section concerns Lie theory over $F$.  We remark that any algebraic group over a local field of characteristic zero may be regarded as a Lie group over some such $F$, so there is little loss of generality in our assumptions.  Our aim is to record some basic properties of ``microlocalized'' functions on Lie groups over $F$.   It seemed clearest to develop the main arguments first in the category of analytic manifolds, and then to specialize at the end to the setting of Lie groups.

The $p$-adic case of the present discussion is not applied in this paper.  It could be useful for generalizing the results of this paper to the depth aspect, although simpler approaches are likely possible for that.  Our main motivation for discussing the $p$-adic case is that many statements and proofs are simpler in that case, and we feel that including them helps motivate the ideal that we are striving for with our treatment of the real case.

The main exports of this section are Definition \ref{defn:we-henceforth-denote} and Propositions \ref{lem:standard2:class-mathfrakcm-p}, \ref{lem:assume-that-rho}, \ref{lem:archimedean-vector-field-act-on-C-2}, \ref{lem:sub-gln:let-gamma-in} and \ref{lem:standard:subcl-mathfr-theta}.  These give a well-defined notion of, roughly, a ``smooth modulated bump concentrated on $p + \O(T^{-1/2})$ and oscillating with frequency $\ll T$'' with good functorial properties (under pushforward, the action of a Lie group, convolution, etc.).  The proofs are direct applications of Taylor's theorem, using that the quadratic and higher terms do not contribute significantly in the indicated range.

\subsection{Preliminaries}

\subsubsection{Pontryagin duality for vector spaces}\label{sec:pontry-dual-vect}
In this section, ``vector space'' always means ``finite-dimensional vector space over $F$.''

Let $V$ be a vector space.  We denote by $V^*$ the dual space, consisting of $F$-linear maps $V \rightarrow F$.  We denote by $V^\wedge$ the Pontryagin dual, consisting of continuous group homomorphisms $V \rightarrow \U(1)$, with $\U(1) \subseteq \mathbb{C}^\times$ the circle group.  Using (say) the ``standard'' nontrivial unitary character $\psi_F : F \rightarrow \U(1)$ (\S\ref{sec:stand-nond-char}), we may identify $V^*$ with $V^\wedge$: to $\ell \in V^*$ corresponds $[V \ni v \mapsto \psi_F(\ell(v)) \in \U(1)] \in V^\wedge$.  In particular, $V^\wedge$ is naturally a vector space of the same dimension as $V$.  We will generally avoid making such identifications here.

As a matter of notation, given $\xi \in V^\wedge$, we denote by \index{Lie algebra!$\psi_{\xi}$}
\begin{equation*}
  \psi_\xi : V \rightarrow \U(1)
\end{equation*}
the homomorphism canonically associated to it.  We emphasize that this device is purely notational, reflecting that we view $\xi$ as a vector under addition and $\psi_\xi$ as a function under multiplication, thus $\psi_{\xi_1 + \xi_2} = \psi_{\xi_1} \psi_{\xi_2}$.

In the archimedean case $F = \mathbb{R}$, we define the bilinear pairing
\begin{equation*}
  \langle ,  \rangle : V \otimes V^\wedge \rightarrow i \mathbb{R}
\end{equation*}
to be the differential of the canonical pairing $V \times V^\wedge \rightarrow \U(1)$, i.e.,
\begin{equation}\label{eq:crb2bvh8nl}
  \langle x, \xi  \rangle := \partial_{t=0} \psi_{t \xi}(x) = \partial_{t=0} \psi_{\xi}(t x),
\end{equation}
so that $\psi_{\xi}(x) = \exp(\langle x, \xi  \rangle)$.  In this way, $V^\wedge$ identifies naturally with the imaginary dual of $V$.

For the sake of orientation and future reference, we record the following:
\begin{lemma}\label{lemma:crb2bvqn54}
  For each vector field $X$ on $V$, we have
  \begin{equation*}
    (X \psi_{\xi}) = \langle \xi, X \rangle \psi_{\xi}.
  \end{equation*}
\end{lemma}
\begin{proof}
  For each $x \in V$, we have
  \begin{align*}
    (X \psi_{\xi})(x)
    &=
      \partial_{t = 0} \psi_{\xi}(x + t X(x)) \\
    &=
      \psi_{\xi}(x)
      \partial_{t = 0} \psi_{\xi}(t X(x)) \\
    &=
      \langle \xi, X(x) \rangle
      \psi_{\xi}(x),
  \end{align*}
  where
  \begin{itemize}
  \item the first identity is the definition of the action of the vector field $X$,
  \item the second is the homomorphism property of $\psi_{\xi}$, and
  \item the third follows from the definition \eqref{eq:crb2bvh8nl} of the pairing $\langle , \rangle$.
  \end{itemize}
\end{proof}

\subsubsection{Local maps}\label{sec:local-maps}
Given topological spaces $X$ and $Y$, a \emph{local map} $f : X \xdashrightarrow{} Y$ is a pair $(\dom(f), f)$, where
\begin{itemize}
\item $\dom(f)$, called the \emph{domain} of $f$, is an open subset of $X$, and
\item $f : \dom(f) \rightarrow Y$ is a continuous map.
\end{itemize}
We say that $f$ is \emph{defined} at an element $p$ of $X$ if $p \in \dom(f)$.

Given two local maps $f : X \xdashrightarrow{} Y$ and $g : Y \xdashrightarrow{} Z$, we may define their composition $g \circ f : X \xdashrightarrow{} Z$, which is a local map with (possibly empty) domain
\begin{equation*}
  \dom(g \circ f) := \dom(f) \cap f^{-1}(\dom(g)).
\end{equation*}

\subsubsection{Multi-index notation}\label{sec:multi-index-notation}
Let $V$ be a vector space.  By choosing a basis for $V$, we may identify it with $F^n$ for some $n$.  By a \emph{multi-index} for $V$, we mean an $n$-tuple $\alpha = (\alpha_1, \dotsc, \alpha_n ) \in \mathbb{Z} _{\geq 0} ^n$ of nonnegative integers.  For a multi-index $\alpha$ and $x \in V$, we set $x^\alpha := \prod_i x_i^{\alpha_i}$ and $|\alpha| := \sum_i \alpha_i$.  In the archimedean case (thus $V = \mathbb{R}^n$), we denote by $\partial^\alpha := \prod_i \frac{\partial^{\alpha_i}}{\partial^{\alpha_i} x_i}$ the translation-invariant differential operator on $C^\infty(V)$ associated to a multi-index $\alpha$.

\subsubsection{Analytic maps}\label{sec:analytic-functions}
Let $V$ be a vector space.  By choosing a basis, we may identify it with $F^n$ for some $n \in \mathbb{Z}_{\geq 0}$.  We denote by $|.|$ the norm on $F^n$ given by $|x| := \max_{i} |x_i|$.

Let $V_1$ and $V_2$ be vector spaces, and let $f : V_1 \xdashrightarrow{} V_2$ be a local map.  We recall that
\begin{itemize}
\item $f$ is \emph{analytic} at a point $p \in \dom(f)$ if it may be expressed in some neighborhood of $p$ by a convergent power series, and
\item $f$ is \emph{analytic} if it is analytic at every point of its domain.
\end{itemize}
By choosing bases, we may identify $V_i$ with $F^{n_i}$ for some $n_i \in \mathbb{Z}_{\geq 0}$.  For analytic $f$, the power series representation of $f$ at $p \in \dom(f)$ may be written
\begin{equation}\label{eq:phip-+-x}
  f(p + x) = \sum_{\alpha } c_\alpha x^\alpha,
\end{equation}
where $\alpha$ runs over multi-indices for $V_1$ and $c_\alpha \in V_2 = F^{n_2}$.  The analyticity condition implies that there exist $C \geq 0$ and $r > 0$ so that for all $\alpha$,
\begin{equation}\label{eq:c_alpha-leq-c_0}
  |c_\alpha| \leq C r^{-|\alpha|}.
\end{equation}
We say then that $f$ is \emph{$(C,r)$-controlled at $p$}.

The following elementary results may be verified by effectivizing the proofs of \cite[p69, Theorem]{MR2179691} and \cite[p70, Theorem]{MR2179691}, respectively.
\begin{lemma}\label{lem:standard2:elementary-analytic-translate}
  Let $V_1$ and $V_2$ be vector spaces equipped with bases, hence with norms $|.|$ as above.  Let $C \geq 0$ and $r, r', r'' > 0$ with $r' + r'' < r$.  There exists $C' \geq 0$ with the following property.  Let $p \in V_1$.  Let $f : V_1 \xdashrightarrow{} V_2$ be analytic, defined at $p$, and $(C,r)$-controlled at $p$.  Let $p' \in V_1$ with $|p' - p| \leq r'$.  Then $f$ is $(C',r'')$-controlled at $p'$.
\end{lemma}
\begin{lemma}\label{lem:standard2:elementary-analytic-compose}
  Let $V_1, V_2, V_3$ be vector spaces equipped with bases, hence with norms $|.|$ as above.  Let $C_1, C_2 \geq 0$ and $r_1, r_2 > 0$.  There exists $C_3 \geq 0$ and $r_3 > 0$ with the following property.  Let $p_1 \in V_1$ and $p_2 \in V_2$.  For $i=1,2$, let $f_i : V_i \xdashrightarrow{} V_{i+1}$ be analytic, defined at $p_i$, and $(C_i,r_i)$-controlled at $p_i$.  Then $f_2 \circ f_1 : V_1 \xdashrightarrow{} V_3$ is analytic, defined at $p_1$, and $(C_3,r_3)$-controlled at $p_1$.
\end{lemma}

\subsubsection{Manifolds}\label{sec:manifolds}
We recall that an \emph{analytic manifold} $M$ over $F$ is a topological space equiped with a maximal atlas of coordinate charts, taking values in finite-dimensional vector spaces over $F$, for which the  transition maps are analytic (see \cite[Part II, Chapter III]{MR2179691} for details).  In what follows, ``manifold'' always means ``analytic manifold over $F$,'' and ``chart'' refers to the analytic structure.

A \emph{morphism} $f : M_1 \rightarrow M_2$ of manifolds is an analytic map, i.e., a map that is described in local coordinates, with respect to every pair (equivalently, sufficiently many pairs) of charts $c_1$ on $M_1$ and $c_2$ on $M_2$, by an analytic map of vector spaces as in \S\ref{sec:analytic-functions}.  We then define \emph{isomorphisms} in the usual way as the morphisms that admit two-sided inverses.

For each $p \in M$, we denote by $T_p(M)$ the tangent space, $T_p^*(M)$ the cotangent space, and $T_p^\wedge(M)$ the Pontryagin dual of the tangent space.

\subsubsection{Pointed manifolds}\label{sec:pointed-manifolds}
By a \emph{pointed manifold} $(M,p)$, we mean a pair consisting of a manifold $M$ and a point $p \in M$.  A \emph{morphism} of pointed manifolds $(M_1, p_1) \rightarrow (M_2, p_2)$ is a morphism of manifolds $M_1 \rightarrow M_2$ that maps $p_1$ to $p_2$.  A \emph{local morphism} of pointed manifolds $f : (M_1,p_1) \xdashrightarrow{} (M_2,p_2)$ is a local map for which
\begin{itemize}
\item $p_1$ belongs to the domain of $f$,
\item $f$ is analytic on its domain, and
\item $f(p_1) = p_2$.
\end{itemize}
It would be cleaner for some purposes to work systematically with germs, but we do not do so here.  In particular, the term ``local'' is meant to emphasize that our primary concern is the behavior of $f$ in some neighborhood of $p_1$ (for instance, when speaking about whether $f$ defines a ``local isomorphism'' in the sense specified below).  On the other hand, we do not rule out the possibility that $f$ happens to be globally defined.

A composition of local morphisms of pointed manifolds, defined as in \S\ref{sec:local-maps}, is likewise a local morphism.

Let $f : (M_1,p_1) \xdashrightarrow{} (M_2,p_2)$ and $g : (M_2,p_2) \xdashrightarrow{} (M_1,p_1)$ be local morphisms of pointed manifolds.  We say that $g$ is an \emph{inverse} for $f$ if $\dom(g) = \image(f)$, $\image(g) = \dom(f)$, and $f \circ g$ (resp.\ $g \circ f$) coincides with the identity map on $\dom(g)$ (resp.\ $\dom(f)$); in that case, $g$ is uniquely determined, and we denote it by $f^{-1}$.  We say that $f$ is a \emph{local isomorphism} if it admits an inverse.  We note that if $f$ satisfies the weaker condition that there exists $g$ for which $f \circ g$ (resp.\ $g \circ f$) coincides with the identity map on some neighborhood of $p_2$ (resp.\ $p_1$), then by shrinking the domain of $f$ to a suitable neighborhood of $p$, we obtain a local isomorphism in the indicated sense.

In this terminology, a chart $c$ for $M$ at $p$ is a local isomorphism $(M,p) \xdashrightarrow{} (V,q)$ for some vector space $V$.  By composing with a translation, we may always arrange that $q$ is the origin $0 \in V$.  The derivative of $c'(p)$ of $c$ at $p$ induces a linear isomorphism between $T_p(M)$ and $T_{c(p)} V = V$.  There is thus little loss of generality in taking $c$ of the form
\begin{equation}\label{eq:c-:-m}
  c : (M,p) \xdashrightarrow{} (T_p(M), 0),
  \quad
  c'(p) = \text{ the identity map on } T_p(M).
\end{equation}
For example, if $M$ is a Lie group and $p$ the identity element, then the logarithm defines a chart of this form.

\subsubsection{Local submersion theorem}\label{sec:inverse-funct-theor}
Let $\rho : (M_1, p_1) \xdashrightarrow{} (M_2, p_2)$ be a morphism of pointed manifolds.  We recall that $\rho$ is \emph{submersive} at $p_1$ if the map $(d \rho)_{p_1} : T_{p_1}(M_1) \rightarrow T_{p_2}(M_2)$ is surjective.  In that case, the implicit function theorem implies that there are vector spaces $V_j$, coordinate charts $c_j : (M_j, p_j) \xdashrightarrow{} (V_j,0)$ and a surjective linear map $\ell : V_1 \rightarrow V_2$ so that the diagram
\begin{equation*}
  \begin{CD}
    M_1 @> \rho >>M_2\\
    @Vc_1VV  @VVc_2V \\
    V_1 @>>\ell>V_2\\
  \end{CD}
\end{equation*}
is defined and commutes in some neighborhood of $p_1$.

We recall that a morphism $\rho : M_1 \rightarrow M_2$ of manifolds is \emph{submersive} if for each $p \in M$, the induced morphism $(M_1, p) \rightarrow (M_2, \rho(p))$ is submersive in the above sense.

\subsubsection{Smooth measures}
Let $V$ be a vector space.  We may speak of smooth functions $V \rightarrow \mathbb{C}$.  In the archimedean case, this carries the usual meaning (``smooth over $\mathbb{R}$'').  In the non-archimedean case, it means ``locally constant;'' we will primarily consider compactly-supported measures, for which it is equivalent to ask for invariance under some open subgroup of $V$.  By a \emph{smooth measure} on $V$, we mean a complex-valued measure obtained by multiplying a Haar measure by a smooth function.

Let $M$ be a manifold.  By a \emph{smooth measure} $\gamma$ on $M$, we mean a complex-valued measure with the following property: for each $p \in M$, there is a coordinate chart $c : (M,p) \xdashrightarrow{} (V, 0)$ and a smooth measure $\gamma_c$ on $V$ so that
\begin{equation*}
  c_*(\gamma|_{\dom(c)}) = \gamma_c |_{\image(c)}.
\end{equation*}
We denote by $\mathcal{M}(M)$ (resp.\ $\mathcal{M}_c(M)$) the space of smooth measures (resp.\ the space of smooth measures of compact support) on $M$.

Let $\rho : M_1 \rightarrow M_2$ be a submersive morphism of manifolds.  Using the implicit function theorem, we see that the pushforward map $\rho_*$ on measures induces a surjective linear map
\begin{equation*}
  \rho_* : \mathcal{M}_c(M_1) \rightarrow \mathcal{M}_c(M_2).
\end{equation*}

Let $\phi : M_1 \rightarrow M_2$ be an isomorphism of manifolds.  The pushforward map $\phi_*$ on measures induces isomorphisms $\phi_* : \mathcal{M}(M_1) \rightarrow \mathcal{M}(M_2)$, mapping $\mathcal{M}_c(M_1)$ to $\mathcal{M}_c(M_2)$.  We write
\begin{equation*}
  \phi^* := \phi_*^{-1} = (\phi^{-1})_* : \mathcal{M}(M_2) \rightarrow \mathcal{M}(M_1)
\end{equation*}
for the inverse map.  It will be convenient to describe such maps in local coordinates.  In doing so, we may assume that each $M_i$ is an open subset of some vector space $V_i$.  Choose Haar measures $d x$ on $V_1$ and $d y$ on $V_2$.  Then, writing $y = \phi(x)$, we have
\begin{equation*}
  d y = |\phi '(x)| \, d x,
\end{equation*}
where $|\phi '(x)|$ denotes the normalized absolute value of the Jacobian of $\phi$ with respect to $d x$ and $d y$; if we identify both $V_1$ and $V_2$ with $F^n$ and choose $d x$ and $d y$ compatibly, then $|\phi '(x)|$ is the absolute value of the determinant of the matrix of partial derivatives.  Let $\gamma \in \mathcal{M}(M_2)$.  We identify $\gamma$ (resp.\ $\phi^* \gamma$) with a smooth function on $M_2$ (resp.\ $M_1$) by dividing by the Haar measure $d y$ (resp.\ $d x$).  Then
\begin{equation}\label{eq:phi-gammax-=jacobian}
  \phi^* \gamma(x) = |\phi '(x)| \gamma(\phi(x)).
\end{equation}

In the archimedean case, given a vector field $X$ on $M$ and a smooth measure $\gamma$, we may define the smooth measure $X \gamma$ via duality.  We explicate this definition in local coordinates.  Suppose that $M$ is an open subset of $\mathbb{R}^n$.  Identify $\gamma$ with a smooth function on $M$ by dividing by Lebesgue measure.  Write $\nabla_X \gamma$ for the ordinary derivative of the smooth function $\gamma$ with respect to $X$.  Let $\div X : M \rightarrow \mathbb{R}$ denote the divergence of $X$ with respect to the given coordinates.  Then the smooth measure $X \gamma$, identified with a smooth function by dividing by Lebesgue measure, is given by
\begin{equation*}
  X \gamma = \nabla_X \gamma + (\div X) \gamma.
\end{equation*}
The verification of this identity is an application of integration by parts.

\subsection{Some nonstandard formalism}

\subsubsection{Infinitesimal neighborhoods}
Let $X$ be a fixed topological space (e.g., a fixed manifold $M$).
\begin{notation}
  For $x,y \in X$, we write $x \simeq y$ to signify that there is a fixed element $z \in X$ so that $x$ and $y$ are both contained in each fixed neighborhood of $z$.
\end{notation}
\begin{example}
  Suppose that $X$ is a fixed normed space $(X,|.|)$, so that notation such as $x \ll y$ and $x = o(y)$ is defined for $x \in X$ and $y \in X \cup \mathbb{R}_{\geq 0}$.  Then for $x,y \in X$ with $x,y \ll 1$, the notation $x \simeq y$ has the same meaning as the notation $x = y + o(1)$ defined in \S\ref{sec:asymptotic-notation-overview} and \S\ref{sec:furth-notat-conv}.
\end{example}

\begin{notation}
  Let $p$ belong to some fixed compact subset of $X$.  We write $p + o(1)$ for the subclass of $X$ consisting of all $x \in X$ for which $x \simeq p$.  We write ``$x = p + o(1)$'' and ``$x \in p + o(1)$'' synonymously; each of these has the same meaning as ``$x \simeq p$.''
\end{notation}

\begin{lemma}
  The following conditions on an element $p \in X$ are equivalent:
  \begin{enumerate}[(i)]
  \item $p$ belongs to some fixed compact subset of $X$.
  \item There is a fixed element $q \in X$ with $p \simeq q$.
  \end{enumerate}
  In that case, the element $q$ is unique provided that $X$ is Hausdorff.
\end{lemma}
\begin{proof}
  This is \cite[Thm 6.1]{MR469763}.
\end{proof}

We refer to the element $q$, as in the conclusion of the lemma, as the \emph{fixed part} of $p$.  This element has the following property: a fixed open subset $E$ of $X$ contains $q$ if and only if it contains $p + o(1)$.\footnote{ The informal content of ``fixed part'', following the motivating discussion of \S\ref{sec:asymptotic-notation-overview}, may be understood as follows.  Suppose given a compact Hausdorff space $X$ and an element $p = p_T \in X$ that depends upon some asymptotic parameter $T \rightarrow \infty$.  Then, possibly after passing to a subsequence of parameters $T_j \rightarrow \infty$, there is a unique limit $q = \lim_{j \rightarrow \infty} p_{T_j}$.
  In particular, for each neighborhood $U$ of $q$, there exists $j_0 = j_0(U)$ so that $p_{T_{j}}$ lies in $U$ whenever $j \geq j_0$.  The nonstandard terminology expresses an equivalent feature but in language that obviates the need to have most quantities depend upon a parameter or to repeatedly pass to subsequences.
}

\subsubsection{Control}\label{sec:control}
We may prepend the adjective ``fixed'' to any of the objects defined above: ``fixed manifold,'' ``fixed chart,'' etc.  To allow more flexibility, we define ``controlled'' analogues of such objects.

\begin{definition}\label{defn:standard2:emphc-point-manif}
  Let $(M,p)$ be a pointed manifold.  We say that $(M,p)$ is \emph{controlled} if
  \begin{itemize}
  \item $M$ is fixed, and
  \item $p$ belongs to some fixed compact subset of $M$.
  \end{itemize}
\end{definition}
The terminology applies in particular to pointed vector spaces $(V,p)$.

\begin{definition}\label{defn:standard2:let-f-controlled-at-p}
  Let $f : V_1 \xdashrightarrow{} V_2$ be an analytic map between vector spaces.  Fix bases for $V_1, V_2$, hence norms $|.|$ as in \S\ref{sec:analytic-functions}.  Let $p \in \dom(f)$.  We say that $f$ is \emph{controlled at $p$} if $f$ is $(C,r)$-controlled at $p$, in the sense of \S\ref{sec:analytic-functions}, for some $C \ll 1$ and $r \gg 1$.
\end{definition}

\begin{lemma}\label{lem:standard2:sett-defin-refd}
  In the setting of Definition \ref{defn:standard2:let-f-controlled-at-p}, if $f$ is controlled at $p$, then it is controlled at every element of the form $p + o(1)$.
\end{lemma}
\begin{proof}
  By Lemma \ref{lem:standard2:elementary-analytic-translate} (with quantification over $V_1, V_2, |.|, C, r, r', r'', C'$ restricted to fixed quantities).
\end{proof}

\begin{example}\label{example:standard2:if-v_1-v_2}
  If $V_1$, $V_2$, $f$ and $p$ are fixed, then $f$ is controlled at $p$.  Lemma \ref{lem:standard2:sett-defin-refd} implies that the same conclusion holds more generally if $p$ belongs to some fixed compact set, since in that case we have $p = q + o(1)$ for some fixed $q$.
\end{example}

\begin{lemma}\label{lem:standard2:controlled-composition}
  Let $f_1 : V_1 \xdashrightarrow{} V_2$ and $f_2 : V_2 \xdashrightarrow{} V_3$ be defined and controlled at $p_1$ and $p_2$, respectively.  Then $f_2 \circ f_1$ is defined and controlled at $p_1$.
\end{lemma}
\begin{proof}
  By Lemma \ref{lem:standard2:elementary-analytic-compose}.
\end{proof}

\begin{definition}\label{defn:let-v_1-v_2}
  We say that a local morphism $f : (V_1, p_1) \xdashrightarrow{} (V_2, p_2)$ of controlled pointed vector spaces is \emph{controlled} if $f$ is defined and controlled at every element of the form $p_1 + o(1)$.
\end{definition}

\begin{remark}
  It is clear that Definition \ref{defn:let-v_1-v_2} depends only upon the class $p_1 + o(1)$ and not the basepoint $p_1$ itself.  The crucial condition is that $f$ be controlled at \emph{some} element of the form $p_1 + o(1)$; in that case, the Taylor series for $f$ at such a point permits us to extend the definition of $f$ to a region containing every element of the form $p_1 + o(1)$.
\end{remark}

The basic axioms of a category are satisfied (although for set-theoretic reasons, we will not introduce that language explicitly):
\begin{lemma}\label{lem:identity-map--m-0}
  The identity map $(V,p) \xdashrightarrow{} (V,p)$ on a fixed pointed vector space is controlled.  The composition $(V_1,p_1) \xdashrightarrow{} (V_2, p_2) \xdashrightarrow{} (V_3, p_3)$ of any pair of controlled local morphisms is controlled.
\end{lemma}
\begin{proof}
  The assertion concerning the identity map is clear, while that concerning compositions follows from Lemma \ref{lem:standard2:controlled-composition}.
\end{proof}

We now extend the above discussion to manifolds.  Let $f : M_1 \rightarrow M_2$ be a morphism (i.e., analytic map) between manifolds $M_1$ and $M_2$.  Thus, for each $p_1 \in M_1$, each chart $c_1 : (M_1, p_1) \xdashrightarrow{} (V_1, 0)$ and each chart $c_2 : (M_2, f(p_2)) \xdashrightarrow{} (V_2,0)$, the composition
\begin{equation*}
  f_{c_1,c_2} : (V_1,0) \xdashrightarrow{c_1^{-1}} (M_1, p_1) \xdashrightarrow{f} (M_2, p_2) \xdashrightarrow{c_2} (V_2,0)
\end{equation*}
is analytic in the sense of \S\ref{sec:analytic-functions}.

\begin{definition}\label{defn:standard2:f-controlled-M1-M2}
  Let $f : (M_1,p_1) \xdashrightarrow{} (M_2,p_2)$ be a local morphism between controlled pointed manifolds.  Let $q_1 \in M_1$ and $q_2 \in M_2$ denote the respective fixed parts of $p_1$ and $p_2$.  We say that $f$ is \emph{controlled} if
  \begin{enumerate}[(i)]
  \item the domain of $f$ contains every element of the form $p_1 + o(1)$,
  \item for each $x = p_1 + o(1)$, we have $f(x) = p_2 + o(1)$, and
  \item for all fixed charts $c_1 : (M_1,q_1) \xdashrightarrow{} (V_1,0)$ and $c_2 : (M_2, q_2) \xdashrightarrow{} (V_2,0)$, the map $f_{c_1,c_2} : (V_1,0) \xdashrightarrow{} (V_2,0)$ as above is controlled at the origin (in the sense of Definition \ref{defn:standard2:let-f-controlled-at-p}).
  \end{enumerate}
\end{definition}

\begin{lemma}\label{lem:identity-map-m}
  The identity map $M \rightarrow M$ on a fixed manifold is controlled.  Any fixed morphism between fixed manifolds is controlled.  The composition of any pair of controlled morphisms is controlled.
\end{lemma}
\begin{proof}
  Immediate from the definitions and Lemma \ref{lem:identity-map--m-0}.
\end{proof}

\begin{definition}
  We say that a controlled local morphism $f : (M_1, p_1) \xdashrightarrow{} (M_2, p_2)$ is a \emph{controlled local isomorphism} if it admits an inverse (resp.\ $g : (M_2, p_2) \xdashrightarrow{} (M_1, M_2)$ that is a controlled local morphism; such an inverse is then unique and denoted $f^{-1}$.
\end{definition}

\begin{definition}
  Given a controlled pointed manifold $(M,p)$ and a fixed vector space $V$, a \emph{controlled chart} $c : (M,p) \xdashrightarrow{} (V,0)$ is a chart that defines a controlled local isomorphism.
\end{definition}

\begin{lemma}\label{lem:families-controlled-maps}
  Let $(X,x_0),(Y,y_0),(Z,z_0)$ be controlled pointed manifolds.  Let $\phi : X \times Z \rightarrow Y$ be a fixed analytic map with $\phi(x_0,z_0) = y_0$.  Then for each $z \in Z$ with $z \simeq z_0$, we have $\phi(x_0,z) \simeq y_0$, and the map
  \begin{equation*}
    \phi_z : (X,x_0) \rightarrow (Y,\phi(x_0,z))
  \end{equation*}
  \begin{equation*}
    x \mapsto \phi(x,z)
  \end{equation*}
  defines a controlled local morphism.
\end{lemma}
\begin{proof}
  We can reduce to the case that $X,Y,Z$ are vector spaces with $x_0,y_0,z_0$ of the form $o(1)$, say.  Since $\phi$ is fixed and analytic, we may write, for $(x,z) \in X \times Z$ with $x,z = o(1)$,
  \begin{equation*}
    \phi_z(x) = \phi(x,z) = \sum_{\alpha,\beta} c_{\alpha,\beta} x^\alpha z^\beta
  \end{equation*}
  where $|c_{\alpha, \beta}| \leq C r^{-|\alpha| - |\beta|}$ for some $C \ll 1$ and $r \gg 1$.  Thus
  \begin{equation*}
    \phi_z(x) = \sum_{\alpha} c_\alpha(z) x^\alpha, \quad
    c_\alpha(z) := \sum_{\beta} c_{\alpha,\beta} z^\beta.
  \end{equation*}
  Since $z = o(1)$, we have $\sum_{\beta} r^{-|\beta|} |z^\beta| \ll 1$, hence $c_\alpha(z) \ll r^{-|\alpha|}$.  Thus $\phi_z$ is controlled at $x$, as required.
\end{proof}

\subsubsection{Distances}\label{sec:distances}
\begin{definition}\label{defn:standard2:let-m-p-distances}
  Let $(M,p)$ be a controlled pointed manifold, with $p$ having fixed part $q$.  Let $c : (M,q) \xdashrightarrow{} (V,0)$ be a fixed chart.  Let $|.|$ be a fixed norm on $V$.  We then define, for any $x,y \in M$ of the form $x,y = p + o(1)$ (equivalently, $q + o(1)$) the distance $\dist_{c,|.|}(x,y)$ to be $|c(x) - c(y)|$.
\end{definition}
\begin{lemma}
  In the context of Definition \ref{defn:standard2:let-m-p-distances}, for any two fixed choices $c_i, |.|_i$ ($i=1,2$) of fixed charts $c_i : (M, q) \xdashrightarrow{} (V_i,0)$ and fixed norms $|.|_i$ on $V$, we have $\dist_{c_1,|.|_1}(x,y) \asymp \dist_{c_2,|.|_2}(x,y)$.
\end{lemma}
\begin{proof}
  The proof is carried out in a more general setting, below, in the proof of Lemma \ref{lem:controlled-uniform-continuity-origin}.
\end{proof}

We henceforth write $\dist_M$, or $\dist$ for short, for some function of the form $\dist_{c,|.|}$ as above, which is thus ``essentially well-defined.''  In particular, for $C \ll 1$ or $\eps \lll 1$, the conditions ``$\dist(x,y) \lll C$'' and ``$\dist(x,y) \ll \eps$'' are well-defined, i.e., independent of the choice of $(c,|.|)$.  The same conditions will often be denoted respectively by
\begin{equation*}
  x = y + o(C), \quad x = y + \O(\eps).
\end{equation*}

Similar definitions apply to elements of $T_p^*(M)$ and $T_p^\wedge(M)$.  For example, given $C \geq 0$ and $\theta \in T_p^\wedge(M)$, we say that $\theta \ll C$ if for some (equivalently, any) controlled chart $c : (M,p) \xdashrightarrow{} (V,0)$, the following condition holds.  Let $c_* \theta \in T_{c(p)}^\wedge(V)$ denote the pushforward of $\theta$.  Using the standard identification of $T_{c(p)} V$ with $V$, we may identify $c_* \theta$ with an element of $V^\wedge$.  Then, with respect to some (equivalently, any) fixed norm $|.|$ on $V^\wedge$, we have $|c_* \theta| \ll C$.

Controlled morphisms satisfy a uniform Lipschitz condition.
\begin{lemma}\label{lem:controlled-uniform-continuity-origin}
  Let $f : (M_1, p_1) \xdashrightarrow{} (M_2, p_2)$ be a controlled local morphism between controlled pointed manifolds.  Then for all $x,y \in M_1$ with $x,y \simeq p_1$, we have
  \begin{equation*}
    \dist(f(x), f(y)) \ll \dist(x,y).
  \end{equation*}
\end{lemma}
\begin{proof}
  Passing to local coordinates, we reduce to the case that each $M_i$ is a vector space $V_i$.  It is enough to show then that for all $x, y \in V_1$ with $x, y = o(1)$, we have
  \begin{equation}\label{eq:fp_1-+-x}
    f(p_1 + x) - f(p_1 + y) \ll x - y.
  \end{equation}
  We note that, per the above discussion, the precise meaning of \eqref{eq:fp_1-+-x} is that the estimate
  \begin{equation*}
    |f(p_1 + x) - f(p_1 + y)|_2 \ll |x - y|_1
  \end{equation*}
  holds for some (equivalently, any) fixed norms $|.|_1$ and $|.|_2$ on $V_1$ and $V_2$.

  Writing $c_\alpha \in V_2$ ($\alpha \in \mathbb{Z}_{\geq 0}^{\dim(V_1)}$) for the Taylor coefficients of $f$ at $p$, we have
  \begin{equation*}
    f(p_1 + x) - f(p_2 + y) = \sum_{|\alpha| \geq 1} c_\alpha (x^\alpha - y^\alpha).
  \end{equation*}
  Fix a norm $|.|$ on $V_1$, and choose $\eps > 0$ with $\eps = o(1)$ so that $|x|, |y| \leq \eps$.  Any difference of monomials $x^\alpha - y^\alpha$ with $|\alpha| \geq 1$ is divisible by some difference of coordinates $x_i - y_i$.  The resulting quotient polynomial is a sum of $\O(|\alpha|)$ many monomials, each of degree $|\alpha| -1$, and with coefficients of size $\O(1)$.  Thus
  \begin{equation*}
    x^\alpha - y^\alpha \ll |\alpha| |x-y| \eps^{|\alpha| - 1}.
  \end{equation*}
  Using that $c_\alpha \ll \exp(C_1 |\alpha|)$ with $C_1 \ll 1$, we deduce that
  \begin{equation*}
    f (p_1 + x) - f (p_2 + y) \ll |x-y| \sum_{|\alpha| \geq 1} |\alpha| (\eps \exp(C_1))^{|\alpha|-1}.
  \end{equation*}
  This last sum may be bounded by a product of convergent geometric series, of size $\O(1)$.
\end{proof}

\begin{corollary}\label{lem:standard:let-phi-controlled-dist-preserve}
  Let $f : (M_1, p_1) \xdashrightarrow{} (M_2, p_2)$ be a controlled local isomorphism between controlled pointed manifolds.  Then for all $x,y \in M_1$ with $x,y \simeq p_1$, we have
  \begin{equation*}
    \dist(f(x), f(y)) \asymp \dist(x,y).
  \end{equation*}
\end{corollary}

\subsubsection{Bumps}

\begin{definition}
  Let $V$ be a fixed vector space.  Let $p \in V$.  Let $\eps > 0$ with $\eps \lll 1$.  We say that a smooth function $\gamma$ on $V$ is an \emph{$\eps$-bump at $p$} if it satisfies the following conditions.
  \begin{enumerate}[(i)]
  \item In the non-archimedean case, we require that
    \begin{equation*}
      \|\gamma\|_\infty \ll \eps^{-\dim(V)}
    \end{equation*}
    and that, for all $x, u \in V$,
    \begin{align*}
      \gamma(x) \neq 0 &\implies x - p \ll  \eps, \\
      u \lll \eps &\implies \gamma(x+u) = \gamma(x).
    \end{align*}
  \item \label{item:epsilon-bump-3} In the archimedean case, we require that $\gamma$ be supported on $p + o(1)$ and that for all fixed $k,m \in \mathbb{Z}_{\geq 0}$, fixed vector fields $X_1,\dotsc,X_m$ on $V$, and all $x \in V$ (with $x = p + o(1)$), we have
    \begin{equation}\label{eq:y_1-dotsb-y_m}
      X_1 \dotsb X_m \gamma(x)
      \ll \eps^{-(\dim(V) + m)} \left( 1 + \frac{|x-p|}{\eps} \right)^{-k}.
    \end{equation}
  \end{enumerate}
  In particular, in either case, $\gamma$ is supported on $p + o(1)$.

  Similarly, we say that a smooth measure $\gamma$ on $V$ is an \emph{$\eps$-bump at $p$} if the smooth function obtained by dividing $\gamma$ by a fixed Haar measure is an $\eps$-bump in the above sense.
\end{definition}
\begin{remark}
  The formulation of condition \eqref{item:epsilon-bump-3} is convenient because when we pass to the setting of (fixed) manifolds, ``(fixed) vector fields'' has an intrinsic meaning, independent of coordinate charts.  For computations, it is convenient to note the following reformulation of that condition: for each fixed $k \in \mathbb{Z}_{\geq 0}$ and multi-index $\alpha$ for $V_1$, we have
  \begin{equation}\label{eq:part-gamm-ll}
    \partial^\alpha \gamma(x) \ll \eps^{ -(\dim (V) + |\alpha|)}  \left( 1 + \frac{|x-p|}{\eps} \right)^{-k}.
  \end{equation}
  The equivalence between the two formulations follows from the product rule for derivatives.
\end{remark}

\begin{lemma}\label{lem:let-v_1-v_2}
  Let $V_1$ and $V_2$ be fixed vector spaces.  Let $p_1 \in V_1$ and $p_2 \in V_2$.  Let $\phi : (V_1,p_1) \xdashrightarrow{} (V_2,p_2)$ be a controlled local isomorphism.  Then pushforward and pullback under $\phi$ induce mutually inverse bijective class maps
  \begin{equation}\label{eq:lefttexteps-bumps-at}
    \left\{\text{$\eps$-bumps at $p_1$} \right\}
    \leftrightarrow
    \left\{\text{$\eps$-bumps at $p_2$} \right\},
  \end{equation}
  either in the sense of functions or measures.
\end{lemma}
\begin{proof}
  Consider first the case of functions.  Let $\gamma_2$ be an $\eps$-bump function at $p_2$.  Let $\gamma_1$ denote the pullback of $\gamma_2$ under $\phi$, i.e., the function supported on $p_1 + o(1)$ and given there by $\gamma_1(x) := \gamma_2(\phi(x))$.  Since $\phi$ is controlled, all of its fixed partial derivatives are of size $\O(1)$.  By the chain rule, it follows that condition \eqref{eq:part-gamm-ll} for $\gamma_2$ implies the same condition for $\gamma_1$.  Therefore $\gamma_1$ is an $\eps$-bump at $p_1$.  The same argument applied to $\phi^{-1}$ yields the required bijection.

  For the case of measures, we use the formula \eqref{eq:phi-gammax-=jacobian} expressing the pullback in terms of the Jacobian $|\phi '(x)|$.  Since $\phi$ is controlled, that Jacobian and its first partial derivatives are of size $\O(1)$.  The proof is otherwise as in the case of functions.
\end{proof}

We now extend to manifolds:

\begin{definition}\label{defn:let-m-p-1}
  Let $(M,p)$ be a controlled pointed manifold.  Let $\eps \lll 1$.  By an \emph{$\eps$-bump at $p$}, we mean a smooth measure $\gamma$ on $M$, supported on $p + o(1)$, such that for some controlled chart $c : (M,p) \xdashrightarrow{} (V,0)$, the pushforward $c_* \gamma$ is an $\eps$-bump at $0$.
\end{definition}
\begin{remark}
  We could just a well make the analogous definition for functions, or more generally for sections of fixed bundles, but in what follows, we require only the given definition concerning measures.
\end{remark}
\begin{lemma}
  If the hypotheses of Definition \ref{defn:let-m-p-1} hold for some $c$, then they hold for all $c$.
\end{lemma}
\begin{proof}
  By Lemma \ref{lem:let-v_1-v_2}.
\end{proof}
\begin{lemma}\label{lem:standard2:any-controlled-local}
  Any controlled local isomorphism of pointed manifolds $\rho : (M_1, p_1) \xdashrightarrow{} (M_2, p_2)$ induces a bijective class map of $\eps$-bumps as in \eqref{eq:lefttexteps-bumps-at}.
\end{lemma}
\begin{proof}
  Immediate from the definitions, since (by Lemma \ref{lem:identity-map-m})  controlled charts for $(M_1,p_1)$ and $(M_2,p_2)$ correspond to one another via $\rho$.
\end{proof}

\begin{lemma}\label{lem:standard2:bump-preserved-under-vector-field}
  Assume that $F = \mathbb{R}$.  Let $X$ be a fixed vector field on $M$.  Let $\gamma$ be an $\eps$-bump on $M$ at $p$.  Then $\eps X \gamma$ is also an $\eps$-bump on $M$ at $p$.
\end{lemma}
\begin{proof}
  Immediate from the definition (especially \eqref{eq:y_1-dotsb-y_m}).
\end{proof}

\begin{lemma}\label{lem:standard2:let-f-vanish-at-p-eps-bump}
  Assume that $F = \mathbb{R}$.  Let $f : M \rightarrow \mathbb{C}$ be a smooth function that vanishes at $p$ and has the property that its image under any fixed differential operator has size $\O(1)$ on $p + o(1)$.  Let $\gamma$ be an $\eps$-bump on $M$ at $p$.  Then $\eps^{-1} f \gamma$ is also an $\eps$-bump on $M$ at $p$.
\end{lemma}
\begin{proof}
  We may assume that $M$ is a vector space $V$.  After dividing by a fixed Haar measure, it will suffice to consider the variant problem in which $\gamma$ is an $\eps$-bump function (rather than a measure).  We must check then that for each fixed multi-index $\alpha$, fixed $k \in \mathbb{Z}_{\geq 0}$ and all $x = p + o(1)$,
  \begin{equation}\label{eq:partialalpha-eps-1}
    \partial^\alpha (\eps^{-1} f  \gamma)(x) \ll \eps^{ -(\dim (V) + |\alpha|)}  \left( 1 + \frac{|x-p|}{\eps} \right)^{-k}.
  \end{equation}
  Our hypotheses on $f$ imply that
  \begin{equation*}
    f(x) \ll |x - p| \ll \eps \left( 1 + \frac{|x-p|}{\eps } \right)
  \end{equation*}
  and, for $|\alpha| \geq 1$,
  \begin{equation*}
    \partial^\alpha f(x) \ll 1.
  \end{equation*}
  Using these estimates, the product rule for differentiation and \eqref{eq:part-gamm-ll}, the required estimate \eqref{eq:partialalpha-eps-1} follows.
\end{proof}

\subsection{Modulated bumps on manifolds}\label{sec:class-bumps}

\subsubsection{Basic definition}

\begin{definition}\label{defn:standard2:given-point-manif}
  Given a pointed manifold $(M,p)$, a chart $c$ as in \eqref{eq:c-:-m} and an element $\theta \in T_p^\wedge(M)$, we denote by\index{Lie algebra!$\psi_{\theta,c}$}
  \begin{equation*}
    \psi_{\theta,c} : M \xdashrightarrow{} \U(1)
  \end{equation*}
  the local map, defined at $p$, given by $x \mapsto \psi_\theta(c(x))$.
\end{definition}

\begin{lemma}\label{lem:standard2:bumps-modulated-become-bumps}
  Let $(M, p)$ be a controlled pointed manifold.  Let $\theta_1, \theta_2 \in T_p^\wedge(M)$ and $\eps > 0$ with $\eps \lll 1$ and
  \begin{equation}\label{eq:epsth-thet-ll}
    \eps |\theta_1 - \theta_2| \ll 1, \quad
    \eps^2 |\theta_1| \ll 1,
    \quad
    \eps^2 |\theta_2| \ll 1.
  \end{equation}
  Let $\gamma$ be an $\eps$-bump at $p$.  Let $c_1$ and $c_2$ be controlled charts $c_i : (M,p) \xdashrightarrow{} (T_p(M), 0)$.  Then
  \begin{equation}\label{eq:psi_theta_1-c_1-1}
    \psi_{\theta_1,c_1}^{-1} \psi_{\theta_2, c_2} \gamma
    \quad
    \text{ and }
    \quad
    \psi_{\theta_2,c_2}^{-1} \psi_{\theta_1, c_1} \gamma
  \end{equation}
  are $\eps$-bumps at $p$.
\end{lemma}
\begin{proof}
  By symmetry, it suffices to show that the first quantity in \eqref{eq:psi_theta_1-c_1-1} is an $\eps$-bump.  By applying a controlled chart, we may reduce to the case that $(M,p) = (V,0)$ for a fixed vector space $V$.  We may assume then that $c_1$ is the identity, while $c_2$ is some controlled morphism $\phi : (V,0) \xdashrightarrow{} (V,0)$ with
  \begin{equation}\label{eq:phi0-=-0}
    \phi(0) = 0, \quad \phi '(0) = \text{identity}.
  \end{equation}
  Then
  \begin{equation*}
    \varphi(x) := \psi_{\theta_1,c_1}^{-1}(x) \psi_{\theta_2,c_2}(x) = \psi_{\theta_2}(\phi(x)) \psi_{\theta_1}^{-1}(x).
  \end{equation*}
  It is clear that $\varphi \gamma$ satisfies the same support condition as $\gamma$.  Our task reduces to verifying the following:
  \begin{itemize}
  \item In the non-archimedean case, for $x \in V$ with $x \lll \eps$, we have
    \begin{equation*}
      \varphi(x) = 1.
    \end{equation*}
  \item In the archimedean case, for each fixed multi-index $\alpha$ and each $x \in V$ with $x \lll 1$, we have
    \begin{equation}\label{eq:part-varphi-x}
      \partial^{\alpha} \varphi (x)
      \ll \eps^{-|\alpha|} \left( 1 + \frac{|x|}{\eps} \right)^{\O(1)}.
    \end{equation}
  \end{itemize}


  Consider first the non-archimedean case.  Let $x \lll 1$.  Since $\phi$ is controlled at $0$, we have $\phi(x) = \phi(0) + \phi '(0) x + \O(|x|^2)$.  Assuming further that $x \lll \eps$, we deduce from \eqref{eq:phi0-=-0} that $\phi(x) = x + o(\eps^2)$.  Since $\eps^2 \theta_2 \ll 1$, it follows that $\psi_{\theta_2}(\phi(x)) = \psi_{\theta_2}(x)$.  Therefore $\varphi(x) = \psi_{\theta_2 - \theta_1}(x)$.  On the other hand, since $\eps |\theta_1 - \theta_2| \ll 1$, we have $|x| \cdot |\theta_1 - \theta_2| \lll 1$.  Therefore $\varphi(x) = 1$, as required.

  We turn to the archimedean case.  By definition,
  \begin{equation*}
    \varphi(x) = \exp(\Phi(x)), \quad
    \Phi(x) := \left\langle \theta_2, \phi(x) \right\rangle - \left\langle \theta_1, x \right\rangle.
  \end{equation*}
  Since $\Phi$ is imaginary-valued, we may reduce via the chain rule to verifying the analogue of \eqref{eq:part-varphi-x} for $\Phi$.  Let $x \lll 1$.  Since $\phi$ is controlled, we have
  \begin{equation*}
    \Phi(x) = \langle \theta_2 - \theta_1, x \rangle + \O(\max(|\theta_1|,|\theta_2|) |x|^2),
  \end{equation*}
  thus
  \begin{align*}
    \Phi(x)
    &\lll
      | \theta_1 - \theta_2| \, |x| +
      \max(|\theta_1|, |\theta_2|) |x|^2 \\
    &\leq
      \eps | \theta_1 - \theta_2|  \left( 1 + \frac{|x|}{\eps } \right)
      +
      \eps^2 \max(|\theta_1|, |\theta_2|)
      \left( 1 + \frac{|x|}{\eps } \right)^2.
  \end{align*}
  The required bound \eqref{eq:part-varphi-x} in the case $\alpha = 0$ thus follows from our hypothesis \eqref {eq:epsth-thet-ll}.  Suppose next that $|\alpha| = 1$.  Then, using again that $\phi$ is controlled, we obtain
  \begin{align*}
    \partial^\alpha \Phi(x)
    &\ll |\theta_1 - \theta_2| + \max(|\theta_1|,|\theta_2|) |x| \\
    &\ll \eps^{-1} \left(
      \eps |\theta_1 - \theta_2| +
      \eps^2 \max(|\theta_1|,|\theta_2|) \left( 1 + \frac{x}{|\eps|} \right) \right),
  \end{align*}
  and we may conclude as in the $\alpha = 0$ case.  Finally, for fixed $\alpha$ with $|\alpha| \geq 2$, we have
  \begin{equation*}
    \partial^\alpha \Phi(x)
    \ll \max(|\theta_1|, |\theta_2|) \ll \eps^{-2},
  \end{equation*}
  which is at least as good as the required estimate.
\end{proof}

\begin{definition}\label{defn:standard:let-v-be}
  Let $(M,p)$ be a controlled pointed manifold.  Given $\theta \in T_p^\wedge(M)$ with $\theta \ll 1$, $T \in F$ with $|T| \ggg 1$, $\delta \in (0,1/2]$ fixed, and a controlled chart $c : (M,p) \xdashrightarrow{} (T_p(M),0)$, we define $\mathfrak{C}(M,p,\theta,T,\delta)_c$ to consist of all smooth measures $\gamma$ on $M$, supported on $p + o(1)$, for which $\psi_{T \theta,c}^{-1} \gamma$ is a $|T|^{\delta-1}$-bump at $p$, and, in the archimedean case, for which $\gamma$ is in fact supported on $p + \O(|T|^{\delta_+-1})$ for some fixed element
  \begin{equation}\label{eq:delt-in-beginc-1}
    \delta_+ \in \begin{cases}
      (\delta,1/2) & \text{ if } \delta < 1/2, \\
      \{1/2\} & \text{ if } \delta = 1/2.
    \end{cases}
  \end{equation}
\end{definition}

\begin{remark}\label{remark:crb2bpyvfj}
  Informally, an element of the class $\mathfrak{C}(M,p,\theta,T,\delta)$ is an approximate modulated delta mass, concentrated on elements of the form $p + \O(|T|^{\delta-1})$ and oscillating there with frequency $T \theta$.

  The support condition in the archimedean case may seem artificial.  It is motivated by two considerations.  On the one hand, we want to retain some absolute control over the support of $\gamma$, which is more convenient in practice than the decay condition \eqref{eq:y_1-dotsb-y_m} alone.  On the other hand, we will eventually produce and study $\gamma$  via truncated Fourier transforms (\S\ref{sec:quant-prel}); it will be useful then to truncate to $p + \O(|T|^{\delta_+-1})$, rather than to the smallest possible range $p + \O(|T|^{\delta-1})$, so that the truncation has a truly negligible effect rather than merely a small effect.
\end{remark}

\begin{proposition}\label{lem:standard2:class-mathfrakcm-p}
  The class $\mathfrak{C}(M,p,\theta,T,\delta)_c$ given by Definition \ref{defn:standard:let-v-be} is independent of the choice of $c$.
\end{proposition}
\begin{proof}
  Let $c_1, c_2 : (M, p) \xdashrightarrow{} (T_p(M),0)$ be controlled charts.  Let $\gamma$ be an element of $\mathfrak{C}(M,p,\theta,T,\delta)_{c_1}$, so that $\gamma$ satisfies the stated support conditions and has the property that its multiple $\psi_{T \theta,c_1}^{-1} \gamma$ is a $|T|^{\delta-1}$-bump at $p$.  We apply Lemma \ref{lem:standard2:bumps-modulated-become-bumps} with $\eps := |T|^{\delta-1}$ and $\theta_1 = \theta_2 = T \theta$.  In view of our assumption $\delta \leq 1/2$, we have
  \begin{equation*}
    \eps^2 |\theta|  \ll |T|^{2 (\delta - 1) + 1} \leq 1,
  \end{equation*}
  so condition \ref{eq:epsth-thet-ll} is satisfied.  The conclusion of that lemma implies that $\psi_{T \theta, c_2}^{-1} \gamma$ is a $|T|^{\delta-1}$-bump at $p$.  Therefore $\gamma \in \mathfrak{C}(M,p,\theta,T,\delta)_{c_2}$.  The reverse inclusion is proved in the same way.
\end{proof}

\begin{definition}\label{defn:we-henceforth-denote}
  We henceforth denote by $\mathfrak{C}(M,p,\theta,T,\delta)$ \index{classes!modulated bumps on a manifold, $\mathfrak{C}(M,p,\theta,T,\delta)$}the common value of the classes attached as in Definition \ref{defn:standard:let-v-be} to controlled charts $c$.
\end{definition}

\begin{example}\label{example:standard2:smooth-modulated-bump-in-frak-C}
  In the setting of Definition \ref{defn:standard:let-v-be}, suppose given $t \in F$ with $|t| \asymp |T|^{\delta-1}$.  Let $c : (M, p ) \xdashrightarrow{} (V,0)$ be a controlled chart.  Fix a smooth compactly-supported measure $\nu$ on $V$.  Define the measure $\gamma$ on $M$ by requiring that, for $f \in C_c(M)$,
  \begin{equation*}
    \int_M f \gamma
    = \int_{x \in V}
    f(c^{-1}(x/t))
    \psi_{T \theta}(x)
    \, d \nu(x).
  \end{equation*}
  Then we may verify readily that $\gamma$ defines an element of $\mathfrak{C}(M,p,\theta,T,\delta)$.
\end{example}

\begin{lemma}\label{lem:if-rho-local}
  Let $\rho : (M_1, p_1) \xdashrightarrow{} (M_2,p_2)$ be a controlled local isomorphism, and let $\theta_1 \in T_{p_1}^\wedge(M_1)$ and $\theta_2 \in T_{p_2}^\wedge(M_2)$ with $\theta_1, \theta_2 \ll 1$ and such that $\theta_1$ is the pullback of $\theta_2$ under the differential of $\rho$.  Then the pushforward map $\rho_*$ induces a bijective class map $\mathfrak{C}(M_1,p_1,\theta_1,T,\delta) \rightarrow \mathfrak{C}(M_2,p_2,\theta_2,T,\delta)$.
\end{lemma}
\begin{proof}
  Immediate from the definitions (as in the proof of Lemma \ref{lem:standard2:any-controlled-local}).
\end{proof}

\subsubsection{Explication in local coordinates}

\begin{example}\label{example:standard2:let-v-p}
  Let $(V,p)$ be a controlled pointed vector space.  Let $\theta$ be an element of the Pontryagin dual $V^\wedge$ with $\theta \ll 1$.  Then, expanding all definitions, we see that $\mathfrak{C}(V, p, \theta, T, \delta)$ is the class of smooth measures $\gamma$ on $V$ with the following properties.

  Fix a Haar measure $d x$ on $V$.  Using $d x$, we may and shall identify smooth measures on $V$ with smooth functions on $V$.  Since any two choices for $d x$ are related by a positive real $C \asymp 1$, the conditions that follow are clearly independent of the choice of $d x$.

  In the non-archimedean case,
  \begin{equation*}
    \|\gamma\|_{\infty} \ll |T|^{(1 - \delta) \dim(V)}
  \end{equation*}
  and, for all $x,u \in V$,
  \begin{align*}
    \gamma(x) \neq 0 &\implies x - p \ll  |T|^{\delta-1}, \\
    u \lll |T|^{\delta-1} &\implies \gamma(x+u) = \psi_{T \theta}(u) \gamma(x).
  \end{align*}

  In the archimedean case, we require the following slightly more technical conditions:
  \begin{enumerate}[(i)]
  \item \label{item:class-CGthetaTdelta:1} There is a fixed element
    \begin{equation*}
      \delta_+ \in \begin{cases}
        (\delta,1/2) & \text{ if } \delta < 1/2, \\
        \{1/2\} & \text{ if } \delta = 1/2
      \end{cases}
    \end{equation*}
    so that $\gamma$ is supported on $p + \O(|T|^{-1 + \delta_+})$.
  \item \label{item:class-CGthetaTdelta:2} Let $\gamma^{\flat}$ denote the smooth measure given by
    \begin{equation}\label{eq:gammaflatx-:=-psi_t}
      \gamma^{\flat}(x) := \psi_{T \theta}^{-1}(x) \gamma(x).
    \end{equation}
    Identify it with a smooth function on $V$ by dividing by a fixed Haar measure.  Then for all fixed $k \in \mathbb{Z}_{\geq 0}$, fixed multi-indices $\alpha \in \mathbb{Z}_{\geq 0}^{\dim(V)}$ and all $x \in V$,
    \begin{equation}\label{eq:y_1-dotsb-y_m-2}
      \partial^\alpha  \gamma^{\flat}(x)
      \ll |T|^{(1 - \delta) (\dim(V) + |\alpha|)} \left( 1 + \frac{|x-p|}{|T|^{\delta -1}} \right)^{-k}.
    \end{equation}
  \end{enumerate}
\end{example}

\subsubsection{Insensitivity}

\begin{lemma}\label{lem:standard2:retain-setting-above}
  Retain the setting of the above definitions.  Let $\theta_1, \theta_2 \in T_p^\wedge(M)$ with
  \begin{equation*}
    \theta_1 = \theta_2 + \O(|T|^{-\delta}).
  \end{equation*}
  Then $\mathfrak{C}(M,p,\theta_1,T,\delta) = \mathfrak{C}(M,p,\theta_2,T,\delta)$.
\end{lemma}
\begin{proof}
  Let $c : (M, p) \xdashrightarrow{} (T_p(M), 0)$ be a controlled chart.  Let $\gamma \in \mathfrak{C}(M,p,\theta_1,T,\delta)$; in particular, $\psi_{T \theta_1, c}^{-1} \gamma$ is a $|T|^{\delta-1}$-bump at $p$.  Lemma \ref{lem:standard2:bumps-modulated-become-bumps}, applied exactly as in the proof of Proposition \ref{lem:standard2:class-mathfrakcm-p}, implies that $\psi_{T \theta_2, c}^{-1} \gamma$ is likewise a $|T|^{\delta-1}$-bump at $p$.  The support conditions concerning $\gamma$ do not depend upon $\theta$, so we deduce the containment $\mathfrak{C}(M,p,\theta_1,T,\delta) \subseteq \mathfrak{C}(M,p,\theta_2,T,\delta)$.  The reverse containment is verified in the same way.
\end{proof}

\begin{lemma}\label{lem:insensitivity}
  Retain the setting of Example \ref{example:standard2:let-v-p}.  Let $p_1, p_2 \in V$ and $\theta_1, \theta_2 \in V^\wedge$ with $\theta_1, \theta_2 \ll 1$.  Assume that
  \begin{equation*}
    p_1 = p_2 + o(|T|^{\delta-1}),
  \end{equation*}
  \begin{equation*}
    \theta_1 = \theta_2 + O(|T|^{-\delta}).
  \end{equation*}
  Then $\mathfrak{C}(V,p_1,\theta_1,T,\delta) = \mathfrak{C}(V,p_2,\theta_2,T,\delta)$.
\end{lemma}
\begin{proof}
  It is easy to see that none of the conditions stated in Example  \ref{example:standard2:let-v-p} change if we swap $p_1$ for $p_2$, so we may reduce to the case that $p_1 = p_2 = p$.  We then appeal to Lemma \ref{lem:standard2:retain-setting-above}.
\end{proof}

For the remainder of this section, we retain the setting of the above definitions concerning $M, p, \theta, T, \delta$.  We focus moreover on individual values of $T$ and $\delta$, and accordingly abbreviate
\begin{equation*}
  \mathfrak{C}(M,p,\theta) := \mathfrak{C}(M,p,\theta,T,\delta).
\end{equation*}
We will also write, e.g., $p_1, p_2$ or $\theta_1, \theta_2$ for quantities satisfying the same conditions as $p$ and $\theta$.

\subsubsection{Total variation norm}

\begin{lemma}\label{lem:standard:each-gamma-in}
  For each $\gamma \in \mathfrak{C}(M,p,\theta)$, the total variation $\|\gamma \|_1$ is $\ll 1$.
\end{lemma}
\begin{proof}
  We can check this in local coordinates, where it is immediate from the concentration property \eqref{eq:y_1-dotsb-y_m-2} applied with $\alpha =0$.
\end{proof}

\subsubsection{Modulation}\label{sec:modulation}

\begin{lemma}\label{lem:modulation}
  Let $\delta \in T_p^\wedge(M)$ with $\delta \ll 1$.  Let $c : (M,p) \xdashrightarrow{} (T_p(M), 0)$ be a controlled chart.  Then multiplication by $\psi_{T \delta,c}$ (Definition \ref{defn:standard2:given-point-manif}) induces a bijective class map $\mathfrak{C}(M,p,\theta) \rightarrow \mathfrak{C}(M,p,\theta+\delta)$.
\end{lemma}
\begin{proof}
  Immediate from Definition \ref{defn:standard:let-v-be}.
\end{proof}

\subsubsection{Pushforward}
Let $\rho : (M_1, p_1) \xdashrightarrow{} (M_2, p_2)$ be a fixed local morphism of fixed pointed manifolds.  (The discussion to follow applies with ``fixed'' replaced by ``controlled,'' but is slightly more involved, and we do not require such generality.)  Then, if $\rho$ is submersive, the charts $c_1, c_2$ and the linear map $\ell$ arising from application of the inverse function theorem (\S\ref{sec:inverse-funct-theor}) may also be taken fixed.

Let $\theta_2 \in T_{p_2}(M_2)^\wedge$ with $\theta_2 \ll 1$.  Set $\theta_1 := \rho^* \theta_2 \in T_{p_1}(M_1)^\wedge$ denote the pullback of $\theta_2$ under the differential of $\rho$.  As noted earlier, we then have $\theta_1 \ll 1$.

\begin{proposition}\label{lem:assume-that-rho}
  Assume that $\rho$ is submersive at $p_1$.  Then the pushforward map $\rho_*$ on spaces of measures induces a surjective class map
  \begin{equation}\label{eq:phi_-:-mathcalcm_1}
    \rho_* : \mathfrak{C}(M_1,p_1,\theta_1) \rightarrow \mathfrak{C}(M_2,p_2,\theta_2).
  \end{equation}
\end{proposition}
\begin{proof}
  We may reduce to the case that $(M_1,p_1) = (V_1,0)$ and $(M_2,p_2) = (V_2,0)$, where $V_1, V_2$ are vector spaces, and $\rho : V_1 \rightarrow V_2$ is a fixed surjective linear map.  Set $U := \ker(\rho)$.  By fixing a splitting of $\rho$, we may identify $V_1$ with $V_2 \oplus U$.  We write elements of $V_1$ as $(x,y)$, with $x \in V_2$ and $y \in U$, so that $\rho(x,y) = x$.  By fixing Haar measures on all vector spaces, we may identify smooth measures with smooth functions.  For $\gamma \in \mathfrak{C}(V_1, 0, \theta_1)$, we have
  \begin{equation*}
    \rho_*(\gamma)(x) = \int_{y \in U}
    \gamma(x, y) \, d y.
  \end{equation*}
  We are given $\theta_2 \in V_2^\wedge$ with $\theta_2 \ll 1$.  The element $\theta_1 = \rho^* \theta_2 \in V_1^\wedge$ is given by $\langle \theta_1, (x,y) \rangle = \langle \theta_2, x \rangle$.

  Let us show first that if $\gamma \in \mathfrak{C}(V_1, 0, \theta_1)$, then $\rho_*(\gamma) \in \mathfrak{C}(V_2,0,\theta_2)$.  We consider the non-archimedean case; the archimedean case may be treated similarly, substituting discussion of invariance with estimates involving derivatives.  We are given that $\gamma$ has $L^\infty$-norm $\ll |T|^{(1 - \delta) \dim(V_1)}$, is supported on elements of size $\O(|T|^{\delta-1})$, and, for all $u \in V_1$ of size $o(|T|^{\delta-1})$, satisfies $\gamma(v + u) = \psi_{T \theta_1}(u) \gamma(v)$ for all $v \in V_1$.  For each $x \in V_2$, the volume of $\{y \in U : \gamma(x,y) \neq 0\}$ is $\ll |T|^{\delta-1}$, so the triangle inequality gives $\|\rho_*(\gamma)\|_{\infty} \ll |T|^{(1 - \delta) \dim(V_2)}$.  The support condition on $\gamma$ implies that $\rho_*(\gamma)$ is supported on elements of the form $\O(|T|^{\delta-1})$.  The equivariance property of $\gamma$ implies that $\rho_*(\gamma)(x + u) = \psi_{T \theta_2}(u) \rho_*(\gamma)(x)$ for all $u \in V_2$.  Thus $\rho_*(\gamma) \in \mathfrak{C}(V_2,0,\theta_2)$, as claimed.

  It remains to check that the class map \eqref{eq:phi_-:-mathcalcm_1} is surjective.  Let $\gamma_1 \in \mathfrak{C}(V_2,0,\theta_2)$ be given.  We may find $\gamma_2 \in \mathfrak{C}(U,0,0)$ with $\int_U \gamma_2 = 1$.  (For instance, take the rescaling by $t \in F^\times$, with $|t| \asymp |T|^{1-\delta}$, of a fixed bump function of total mass one.)  Define $\gamma(x,y) := \gamma_1(x) \gamma_2(y)$.  It is then easy to see that $\gamma \in \mathfrak{C}(V_1,0,\theta_1)$ and $\rho_*(\gamma) = \gamma_1$.
\end{proof}

\begin{remark}\label{remark:crb2bq4vf9}
  Proposition \ref{lem:assume-that-rho} describes the image of $\mathfrak{C}(M_1, \rho_1, \theta_1)$ under pushforward when $\theta_1$ lies in the image of the pullback map $T_{p_2}(M_2)^\wedge \rightarrow T_{p_1}(M_1)^\vee$.  On the other hand, if $\theta_1$ is far enough away from that image, then one can show that the induced class map is negligible in a certain sense; special cases of this feature follow from the star product asymptotics recalled in \S\ref{sec:star-prod-asymptotics}.
\end{remark}

\subsubsection{Negligible bumps}
Given any complex scalar $t$, we may define the rescaled class $t \mathfrak{C}(M,p,\theta,T,\delta)$, consisting of smooth measures on $M$ of the form $t \gamma$ for some $\gamma \in \mathfrak{C}(M,p,\theta,T,\delta)$.  This notation applies in particular when $t$ is a power of $|T|$.  We denote by
\begin{equation*}
  |T|^{-\infty} \mathfrak{C}(M,p,\theta,T,\delta)
\end{equation*}
the intersection over all fixed $m \geq 0$ of the classes $|T|^{-m} \mathfrak{C}(M,p,\theta,T,\delta)$ (cf.\ \S\ref{sec:furth-notat-conv}).  Passing to local coordinates as in Example \ref{example:standard2:let-v-p}, we readily obtain the following explicit description of that class:
\begin{lemma}\label{lem:standard2:let-v-p}
  Let $(V,p)$ be a controlled point vector space, and let $(\theta,T,\delta)$ be as in Example \ref{example:standard2:let-v-p}.  We identify smooth measures on $V$ with smooth functions via some fixed Haar measure.

  In the non-archimedean case, $|T|^{-\infty} \mathfrak{C}(V,p,\theta,T,\delta)$ identifies with the space of smooth functions $\gamma$ on $V$ such that
  \begin{equation}\label{eq:gamma-_infty-ll}
    \|\gamma \|_{\infty} \ll |T|^{-\infty},
  \end{equation}
  (i.e., following the conventions of \S\ref{sec:asympt-notat-term}, $\|\gamma \|_{\infty} \ll |T|^{-m}$ for each fixed $m$), and, for all $x,u \in V$,
  \begin{align*}
    \gamma(x) \neq 0 &\implies x - p \ll  |T|^{\delta-1}, \\
    u \lll |T|^{\delta-1} &\implies \gamma(x+u) = \psi_{T \theta}(u) \gamma(x).
  \end{align*}

  In the archimedean case,  $|T|^{-\infty} \mathfrak{C}(V,p,\theta,T,\delta)$ identifies with the space of smooth functions $\gamma$ on $V$ that are supported on  $p + \O(|T|^{-1 + \delta_+})$ for some fixed $\delta_+$ as in  \eqref{item:class-CGthetaTdelta:1} and such that for all fixed $\ell \in \mathbb{Z}_{\geq 0}$ and $\alpha \in \mathbb{Z}_{\geq 0}^{\dim(V)}$,
  \begin{equation*}
    \|\partial^\alpha \gamma\|_{\infty}
    \ll |T|^{-\ell}.
  \end{equation*}
\end{lemma}
\begin{proof}
  Abbreviate $\mathfrak{C} := \mathfrak{C}(V,p,\theta,T,\delta)$.  It is easy to see that any $\gamma$ satisfying the above conditions lies in $|T|^{-\ell} \mathfrak{C} $ for each fixed $\ell$, hence lies in $|T|^{-\infty} \mathfrak{C}$.  Conversely, suppose $\gamma$ lies in $|T|^{-\ell} \mathfrak{C}$ for each fixed $\ell$.  Consider first the non-archimedean case.  Then our hypothesis, applied with $\ell = 0$, gives the required support and invariance conditions.  Applying our hypothesis again but with $\ell$ arbitary then gives the required $L^\infty$ bound \eqref{eq:gamma-_infty-ll}.  Consider next the archimedean case.  Our hypothesis, applied with $\ell = 0$, gives the required support condition.  Applying our hypothesis again but with $\ell$ arbitrary gives, with $\gamma^{\flat}$ defined as in \eqref{eq:gammaflatx-:=-psi_t}, the following assertion: for fixed $k,\ell \in \mathbb{Z}_{\geq 0}$ and $\alpha \in \mathbb{Z}_{\geq 0}^{\dim(V)}$,
  \begin{equation*}
    \partial^\alpha \gamma^{\flat}(x) \ll
    |T|^{-\ell + (1 - \delta)(\dim(V) + |\alpha|)} \left( 1 + \frac{|x - p|}{|T|^{\delta - 1}} \right)^{-k}.
  \end{equation*}
  In view of the support condition that $x - p \ll |T|^{\delta_+ - 1}$, it is equivalent to ask simply that for all fixed $\ell, m \in \mathbb{Z}_{\geq 0}$ and fixed $\alpha$, we have
  \begin{equation}\label{eq:part-gamm-ll-1}
    \|\partial^\alpha \gamma^{\flat}\|_{\infty} \ll |T|^{-\ell}.
  \end{equation}
  Finally, since the estimate $\partial^\alpha \psi_{T \theta}(x) \ll |T|^{|\alpha|}$ holds for fixed $\alpha$, we deduce that the estimates \eqref{eq:gamma-_infty-ll} and \eqref{eq:part-gamm-ll-1}, taken over all fixed $\ell \geq 0$, are equivalent.
\end{proof}

The following technical lemma will be needed at one point, in the proof of Proposition \ref{prop:standard2:each-fixed-x}.
\begin{lemma}\label{lem:standard2:overspill-frak-C-M}
  Suppose that $\gamma \in |T|^{-\infty} \mathfrak{C}(M,p,\theta,T,\delta)$.  Then there exists a natural number $\ell \ggg 1$ such that $|T|^{\ell} \gamma \in |T|^{-\infty} \mathfrak{C}(M,p,\theta,T,\delta)$.
\end{lemma}
\begin{proof}
  Abbreviate $\mathfrak{C} := \mathfrak{C}(M,p,\theta,T,\delta)$.  We can assume that $(M,p) = (V,p)$, so that the description of Lemma \ref{lem:standard2:let-v-p} applies.  Our hypothesis is that $\gamma \in |T|^{-\infty} \mathfrak{C}$.  The desired conclusion is that $|T|^{\ell} \gamma \in |T|^{-\infty} \mathfrak{C}$ for some $\ell \geq 0$ with $\ell \ggg 1$.  As we now explain, this is a simple application of overspill.

  In the non-archimedean case, the required support and invariance conditions for $|T|^{\ell} \gamma$, for any $\ell$, follow immediately from our hypothesis.  Our hypothesis \eqref{eq:gamma-_infty-ll} implies also that, for each fixed $\ell$, we have $\||T|^{\ell} \gamma \|_{\infty} \leq |T|^{-\ell}$.  By overspill, we may find $\ell \ggg 1$ for which $\||T|^{\ell} \gamma \|_{\infty} \leq |T|^{-\ell}$, hence $\||T|^{\ell} \gamma \|_{\infty} \ll |T|^{-\infty}$, as required.  (For this last implication, we use that if $c \geq 0$ satisfies $c \leq \lvert T \rvert^{- \ell}$, then for each fixed $m$, we have $m \leq \ell$, hence $\lvert T \rvert^{- \ell} \leq \lvert T \rvert^{- m}$, hence $c \leq \lvert T \rvert^{- m}$; since $m$ was arbitrary, it follows that $c \ll \lvert T \rvert^{- \infty}$.)

  In the archimedean case, the support condition again follows immediately from our hypothesis.  We also know that for each fixed $\ell \in \mathbb{Z}_{\geq 0}$, we have $\|\partial^\alpha |T|^{\ell} \gamma(x)\|_{\infty} \leq |T|^{-\ell}$ for all multi-indices $\alpha$ with $|\alpha| \leq \ell$, say.  By overspill, the same conclusion holds for some $\ell \ggg 1$.  With this choice of $\ell$, we then have $\|\partial^\alpha |T|^{\ell} \gamma (x)\|_{\infty} \ll |T|^{-\infty}$, whence $|T|^{\ell} \gamma \in \mathfrak{C}$, as required.
\end{proof}

\subsubsection{Fourier-analytic description}

For a vector space $V$ over $F$, we denote by $\mathcal{S}(V)$ the Schwartz space of $V$, by $\mathcal{M}_{\mathcal{S}}(V)$ the space of Schwartz measures on $V$ (i.e., measures of the form $\varphi(x) \, d x$, where $\varphi \in \mathcal{S}(V)$ and $d x$ is a Haar measure on $V$), and by $V^\wedge$ the Pontryagin dual of $V$.  We denote by $f \mapsto f^\wedge$ the Fourier transform $\mathcal{M}_{\mathcal{S}}(V) \rightarrow \mathcal{S}(V^\wedge)$ given by
\begin{equation*}
  f^\wedge(\xi) = \int_{x \in V} f(x) \psi_\xi(x),
\end{equation*}
and by $f \mapsto f^\vee$ its inverse.  Given $T \in F^\times$, we define the rescaled Fourier transform
\begin{equation*}
  \mathcal{F}_T : \mathcal{M}_{\mathcal{S}}(V) \rightarrow \mathcal{S}(V^\wedge)
\end{equation*}
by the formula $\mathcal{F}_T(f)(\xi) := f^\wedge(T \xi)$.  We denote by $\mathcal{F}_T^{-1}$ its inverse.

\begin{definition}\label{defn:let-mathfrakav-theta}
  Let $\mathfrak{A}(V,\theta,T,\delta)$ denote the class of all smooth functions $a$ on $V^\wedge$ with the following properties.
  \begin{itemize}
  \item In the archimedean case, each derivative of $a$ is of rapid decay.  Moreover, for each fixed $m, \ell \in \mathbb{Z}_{\geq 0}$ and fixed vector fields $X_1,\dotsc,X_m$ on $V^\wedge$, we have
    \begin{equation*}
      X_1 \dotsb X_m a(\xi) \ll |T|^{m \delta} \left\langle |T|^{\delta} (\xi + \theta ) \right\rangle^{-\ell}.
    \end{equation*}
  \item In the non-archimedean case, for $\xi, \eta \in V^\wedge$,
    \begin{align*}
      a(\xi ) \neq 0 &\implies |\xi + \theta| \ll |T|^{-\delta}, \\
      |\eta| \lll |T|^{-\delta} &\implies a(\xi + \eta  ) = a(\xi).
    \end{align*}
  \end{itemize}
\end{definition}

\begin{lemma}\label{lem:let-c-frakC-vs-frakA}
  Let $c : (M,p) \xdashrightarrow{} (T_p(M), 0)$ be a controlled chart.
  \begin{enumerate}[(i)]
  \item Let $\gamma \in \mathfrak{C}(M,p,\theta,T,\delta)$.  Then the rescaled Fourier transform $\mathcal{F}_T(c_*\gamma) \in \mathcal{S}(T_p^\wedge(M))$ lies in $\mathfrak{A}(T_p(M),\theta,T,\delta)$.
  \item Suppose $F$ is non-archimedean.   Let $a \in \mathfrak{A}(T_p(M),\theta,T,\delta)$.  Then $c^* \mathcal{F}_T^{-1}(a) \in \mathfrak{C}(M,p,\theta,T,\delta)$.  Thus the maps
    \begin{equation*}
      \mathcal{F}_T \circ c_* : \mathfrak{C}(M,p,\theta,T,\delta) \rightarrow \mathfrak{A}(T_p(M),\theta,T,\delta),
    \end{equation*}
    \begin{equation*}
      c^* \circ \mathcal{F}_T^{-1}  : \mathfrak{A}(T_p(M),\theta,T,\delta) \rightarrow \mathfrak{C}(M,p,\theta,T,\delta)
    \end{equation*}
    are mutually inverse class bijections.
  \item \label{itm:standard:fourier-description-of-bumps-arch} Suppose $F$ is archimedean.  Let $a \in \mathfrak{A}(T_p(M),\theta,T,\delta)$.  Let $\chi \in C_c^\infty(M)$ be supported on $p + \O(|T|^{\delta_+-1})$ for some fixed $\delta_+$ as in \eqref{eq:delt-in-beginc-1}, and suppose that $\partial^\alpha \chi(x) \ll |T|^{(1 - \delta) |\alpha|}$ for all fixed $\alpha$ and all $x$.  Then $\chi \cdot c^* \mathcal{F}_T^{-1}(a)$ lies in the class $\mathfrak{C}(M,p,\theta,T,\delta)$.  Conversely, every element of that class arises in this way.
  \end{enumerate}
\end{lemma}
\begin{proof}
  Assertions (i) and (ii), and all but the final assertion of (iii), are essentially a multi-dimensional variant of Lemma \ref{lem:four-transf-begin}; they follow again by repeated integration by parts in the definition of the Fourier integral.  For the final assertion, let $\gamma \in \mathfrak{C}(M,p,\theta,T,\delta)$.  We may construct $\chi$, as indicated, with the additional property that $\chi$ takes the value $1$ on all elements of the form $p + o(|T|^{\delta_+-1})$.  Then $\chi \gamma = \gamma$, hence $\gamma = \chi \cdot c^* \mathcal{F}_T^{-1} (a)$ with $a := \mathcal{F}_T (c_* \gamma)$.
\end{proof}

\begin{remark}
  In the non-archimedean case, we thus have a clear Fourier-analytic interpretation of $\mathfrak{C}(M,p,\theta,T,\delta)$.  In the archimedean case, for $\chi$ as in assertion \eqref{itm:standard:fourier-description-of-bumps-arch} and $\gamma \in \mathfrak{C}(M,p,\theta,T,\delta)$, it is not hard to see that $\chi \cdot \gamma - \gamma \in |T|^{-\infty} \mathfrak{C}(M,p,\theta,T,\delta)$.  Thus the maps
  \begin{equation*}
    \mathcal{F}_T \circ c_* : \mathfrak{C}(M,p,\theta,T,\delta) \rightarrow \mathfrak{A}(T_p(M),\theta,T,\delta),
  \end{equation*}
  \begin{equation*}
    \chi \cdot c^* \circ \mathcal{F}_T^{-1}  : \mathfrak{A}(T_p(M),\theta,T,\delta) \rightarrow \mathfrak{C}(M,p,\theta,T,\delta)
  \end{equation*}
  are, in a certain sense, asymptotically inverse to one another.
\end{remark}

\subsubsection{Action by vector fields}

\begin{lemma}\label{lem:archimedean-vector-field-act-on-C}
  Assume that $F = \mathbb{R}$.  Let $X$ be a fixed vector field on $M$.  Let $\gamma \in \mathfrak{C}(M,p,\theta)$.  Then
  \begin{equation*}
    X \gamma - \langle T \theta, X(p) \rangle \gamma
    \in
    |T|^{1-\delta}
    \mathfrak{C}(M,p,\theta).
  \end{equation*}
\end{lemma}
\begin{proof}
  We may reduce to the case $(M,p) = (V,0)$, $\theta \in V^\wedge$.  We write $\gamma = \psi_{T \theta} \gamma^{\flat}$, where $\gamma^{\flat}$ is a $|T|^{\delta-1}$-bump.  By the product rule, $X \gamma = (X \psi_{T \theta}) \gamma^{\flat} + \psi_{T \theta} (X \gamma^{\flat})$.  By Lemma \ref{lem:standard2:bump-preserved-under-vector-field}, we know that $|T|^{\delta - 1} X  \gamma^{\flat}$ is a $|T|^{\delta-1}$-bump, hence that $\psi_{T \theta} (X \gamma^{\flat}) \in |T|^{1- \delta} \mathfrak{C}(M,p,\theta)$.  We thereby reduce to verifying that
  \begin{equation}\label{eq:x-psi_t-theta}
    (X \psi_{T \theta}) \gamma^{\flat} - \langle T \theta, X(0) \rangle \gamma \in |T|^{1-\delta} \mathfrak{C}(M,p,\theta).
  \end{equation}
  To that end, observe first (by Lemma~\ref{lemma:crb2bvqn54}) that
  \begin{equation*}
    (X \psi_{T \theta}) = \langle T \theta, X \rangle \psi_{T \theta}.
  \end{equation*}
  It follows that
  \begin{equation*}
    (X \psi_{T \theta}) \gamma^{\flat} - \langle T \theta, X(0) \rangle \gamma
    =
    \langle
    T \theta, X - X(0)
    \rangle
    \gamma.
  \end{equation*}
  Now, write
  \begin{equation*}
    f(x) :=
    \frac{
      \langle
      T \theta, X(x) - X(0)
      \rangle
    }{
      |T|
    }.
  \end{equation*}
  Then $f$ is a smooth function that vanishes at the origin and has each fixed derivative of size $\O(1)$.  By Lemma \ref{lem:standard2:let-f-vanish-at-p-eps-bump}, we deduce that $|T|^{1 - \delta} f \gamma^{\flat}$ is a $|T|^{\delta-1}$ bump, thus $|T|^{1- \delta} f \gamma \in \mathfrak{C}(M,p,\theta)$, which implies the required membership \eqref{eq:x-psi_t-theta} in view of the inequality $\delta \leq 1 - \delta$ following from our assumption that $\delta \leq 1/2$.
\end{proof}

\subsubsection{Action by Lie groups}\label{sec:translation}
Let $G$ be a fixed analytic group over $F$, i.e., a topological group that is also a manifold (analytic, by our conventions) for which the multiplication and inversion maps are analytic \cite[\S II.IV]{MR2179691}.  We denote by $1 \in G$ the identity element and by $\Lie(G)$ the Lie algebra, which is a vector space over $F$.  The group logarithm defines a (fixed) chart $(G,1) \xdashrightarrow{} (\Lie(G), 0)$.

Assume given a fixed action $G \times M \rightarrow M$ on our manifold $M$.  We denote this action simply by juxtaposition.  This action induces one of $G$ on functions on $M$, measures on $M$, and so on.

\begin{lemma}\label{lem:standard2:G-defines-controlled-local-isoms}
  Let $p$ belong to some fixed compact subset of $M$.  Let $g \in G$ with $g = 1 + o(1)$.  Then
  \begin{enumerate}[(i)]
  \item $g p$ likewise belongs to some fixed compact subset of $M$, so $(M,p)$ and $(M, gp)$ are controlled pointed manifolds in the sense of \S\ref{sec:control}, and
  \item the induced map
    \begin{equation}\label{eq:g-:-m}
      g : (M,p) \rightarrow  (M, g p)
    \end{equation}
    defines a controlled local isomorphism.
  \end{enumerate}
\end{lemma}
\begin{proof}
  Lemma \ref{lem:families-controlled-maps} (applied with $X, Y = M$ and $Z = G$) implies that \eqref{eq:g-:-m} is defines controlled local morphism.  The required conclusion then follows from the same argument applied to $g^{-1}$.
\end{proof}

\begin{definition}
  Let $p$ belong to a fixed compact subset of $M$.  Let
  \begin{equation*}
    \mu : T_p^\wedge(M) \rightarrow \Lie(G)^\wedge
  \end{equation*}
  denote the dual of the differential of the orbit map $G \rightarrow M$, $g \mapsto g p$.
\end{definition}

We have the following consequence of Lemma \ref{lem:archimedean-vector-field-act-on-C}:
\begin{proposition}\label{lem:archimedean-vector-field-act-on-C-2}
  Assume that $F$ is archimedean.  Let $x \in \Lie(G)$.  Let $X$ denote the vector field on $M$ associated to $x$ and the given action of $G$.  Let $\gamma \in \mathfrak{C}(M,p,\theta)$.  Then
  \begin{equation*}
    X \gamma - \langle \mu(T \theta), x \rangle \gamma
    \in
    |T|^{1-\delta}    |x|    \mathfrak{C}(M,p,\theta).
  \end{equation*}
\end{proposition}
\begin{proof}
  Fix a basis $\{X_1,\dotsc,X_n\}$ of $\Lie(G)$, and write $X = \sum c_j X_j$.   Then $|x| \asymp \max_j |c_j|$.  By definition, $\langle \mu(T \theta), x \rangle = \langle T \theta, X(p) \rangle$.  The claim follows by summing up the application of Lemma \ref{lem:archimedean-vector-field-act-on-C} to each basis element $X_j$.
\end{proof}

We now record the effect of the group action on the classes $\mathfrak{C}(M,p,\theta)$.  Note that for each $g \in G$ with $g = 1 + o(1)$, we may define
\begin{equation*}
  \psi_{\mu(T \theta)}(\log(g)) \in \U(1).
\end{equation*}
\begin{proposition}\label{lem:sub-gln:let-gamma-in}
  Let $\gamma \in \mathfrak{C}(M,p,\theta)$.  Let $g \in G$ with $g = 1 + o(|T|^{\delta-1})$.  Then the pushforward $g_* \gamma$ of $\gamma$ under the action map belongs to $\mathfrak{C}(M,p,\theta)$.  More precisely,
  \begin{itemize}
  \item in the non-archimedean case, we have
    \begin{equation}\label{eq:g_-gamma-equivariance-non-arch}
      g_* \gamma = \psi_{\mu(T \theta)} (\log (g)) \gamma,
    \end{equation}
    while
  \item in the archimedean case,
    \begin{equation*}
      g_* \gamma - \psi_{\mu(T \theta)}(\log(g)) \gamma \in |T|^{1 - \delta} \dist_G(g,1) \mathfrak{C}(M,p,\theta).
    \end{equation*}
  \end{itemize}
\end{proposition}
\begin{proof}
  We show first that $g_* \gamma \in \mathfrak{C}(M,p,\theta)$.  We choose a controlled chart $c : (M,p) \xdashrightarrow{} (V,0)$ with $V = T_p(M)$.  Let $g \in G$ with $g = 1 + o(|T|^{\delta-1})$.  Let $\phi_g : V \xdashrightarrow{} V$ denote the transport, via $c$, of the action map for $g$.  By Lemma \ref{lem:standard2:G-defines-controlled-local-isoms}, it defines a controlled local isomorphism $(V,0) \xdashrightarrow{} (V, \phi_g(0))$.  Using $c$, we may identify $\mathfrak{C}(M,p,\theta)$ with $\mathfrak{C}(V,0,\eta)$ for some $\eta \in V^\wedge$ with $\eta \ll 1$.  Define $\eta '$ by requiring that $\phi '_g(0)^* \eta ' = \eta$.  By Lemma \ref{lem:if-rho-local}, we have
  \begin{equation*}
    (\phi_g)_* c_* \gamma \in \mathfrak{C}(V, \phi_g(0), \eta ').
  \end{equation*}
  Since $\|g - 1\| \lll |T|^{\delta-1} \leq |T|^{-1/2}$, we have $\phi_g(0) \lll |T|^{\delta-1}$ and  $\eta ' = \eta + o(|T|^{-1/2})$, so by Lemma \ref{lem:insensitivity},
  \begin{equation*}
    \mathfrak{C}(V, \phi_g(0), \eta ')
    =
    \mathfrak{C}(V,0,\eta ') = \mathfrak{C}(V,0,\eta).
  \end{equation*}
  Thus $(\phi_g)_* c_* \gamma  \in \mathfrak{C}(V,0,\eta)$ and so $g_* \gamma \in \mathfrak{C}(M,p,\theta)$, as required.

  We now establish the refined assertions.  In the non-archimedean case, we see from
  \begin{itemize}
  \item the hypothesis $g = 1 + o(|T|^{\delta-1})$,  and
  \item Lemma \ref{lem:controlled-uniform-continuity-origin}, applied to the given action map $G \times M \rightarrow  M$,
  \end{itemize}
  that, for all $v \in V$ with $v = o(1)$, we have
  \begin{equation*}
    \phi_g(v) = v + o(|T|^{\delta-1}).
  \end{equation*}
  The claimed equivariance property \eqref{eq:g_-gamma-equivariance-non-arch} then follows immediately from Definition \ref{defn:standard:let-v-be}.

  Turning to the archimedean case, write $g = \exp(x)$ with $x \in \Lie(G)$.  Then each of the elements $\exp(t x)$, for $t \in [0,1]$, is likewise of the form $1 + o(|T|^{\delta-1})$, so by what was just shown, each of the elements
  \begin{equation*}
    \gamma_t := \psi_{\mu(T \theta)}^{-1}(t x) (\phi_{\exp(t x)})_* \gamma,
  \end{equation*}
  for $t \in [0,1]$, lies in the class $\mathfrak{C}(V,0,\theta,T,\delta)$.  Our aim is to verify that
  \begin{equation*}
    \gamma_1 - \gamma_0 \in |T|^{1-\delta} \|g - 1\| \mathfrak{C}(V,0,\theta).
  \end{equation*}
  By the fundamental theorem of calculus, it is enough to show that, for all $t \in [0,1]$,
  \begin{equation}\label{eq:partial_t-gamma_t-in}
    \partial_t \gamma_t  \in |T|^{1-\delta} \|g - 1\| \mathfrak{C}(V,0,\theta).
  \end{equation}
  Let $X$ denotes vector field on $V$, defined near the origin, associated to $x$.  Then
  \begin{equation*}
    \partial_t \gamma_t
    =
    X
    \gamma_t
    -
    \langle \mu(T \theta), x \rangle \gamma_t.
  \end{equation*}
  By Proposition \ref{lem:archimedean-vector-field-act-on-C-2}, it follows that
  \begin{equation*}
    \partial_t \gamma_t \in |T|^{1 - \delta} |x| \mathfrak{C}(V,0,\theta).
  \end{equation*}
  Since $|x| \asymp \|g-1\|$, the required membership \eqref{eq:partial_t-gamma_t-in} follows.
\end{proof}

\subsection{Modulated bumps on Lie groups}\label{sec:spec-group-sett}
We may specialize Definition \ref{defn:we-henceforth-denote} to the case that $(M,p) = (G,1)$ for some fixed analytic group $G$ with identity element $1$.  The group $G$ acts on itself in three standard ways: by left translation, right translation, and conjugation.  The discussion of \S\ref{sec:translation} applies to all three actions (although in the case of conjugation, more precise estimates are possible).

We often focus on the case that $p$ is the identity element $1$, and accordingly abbreviate \index{classes!modulated bumps on a Lie group, $\mathfrak{C}(G,\theta,T,\delta)$}
\begin{equation}\label{eq:mathfrakcg-theta-t}
  \mathfrak{C}(G,\theta,T,\delta) := \mathfrak{C}(G,1,\theta,T,\delta)
\end{equation}
We further abbreviate $\mathfrak{C}(G,\theta) := \mathfrak{C}(G, \theta, T,\delta)$ when $T$, $\delta$ are clear from context.

Recall that $\mathcal{M}_c(G)$ denotes the space of compactly-supported smooth measures on $G$, of which $\mathfrak{C}(G,\theta)$ is a subclass.  For $\gamma \in \mathcal{M}_c(G)$, we denote by $\gamma^*$ its adjoint, i.e., the complex conjugate of the pushforward under the inversion map.  For $\gamma_1, \gamma_2 \in \mathcal{M}_c(G)$, we denote by $\gamma_1 \ast \gamma_2$ their convolution, i.e., the pushforward of $\gamma_1 \times \gamma_2$ under the multiplication map $G \times G \rightarrow G$.
\begin{proposition}\label{lem:standard:subcl-mathfr-theta}
  The subclass $\mathfrak{C}(G,\theta)$ of $\mathcal{M}_c(G)$ is closed under taking adjoints and convolution.
\end{proposition}
\begin{proof}
  Let $\iota : G \rightarrow G$ denote the inversion map.  Then $\theta$ is the pullback of $- \theta$ under the differential of $\iota$, so by Lemma \ref{lem:if-rho-local}, we obtain the class map $\iota_* : \mathfrak{C}(G,\theta) \rightarrow \mathfrak{C}(G,-\theta)$.  We verify readily that complex conjugation induces $\mathfrak{C}(G,-\theta) \rightarrow \mathfrak{C}(G,\theta)$.  We conclude that the adjoint map $\gamma \mapsto \gamma^*$ induces $\mathfrak{C}(G,\theta) \rightarrow \mathfrak{C}(G,\theta)$.

  The multiplication map $\mu : G \times G \rightarrow G$ is submersive at the identity and satisfies $\mu^*(\theta) = (\theta,\theta)$, so closure under convolution follows from Proposition \ref{lem:assume-that-rho}.
\end{proof}

\section{Parameters and stability}\label{sec:parameters-stability}
In this section, given a Lie group $G$ over $\mathbb{R}$, we often denote simply by $\mathfrak{g} := \Lie(G)$ its Lie algebra and by $\mathfrak{g}^\wedge := \Lie(G)^\wedge$ the Pontryagin dual of the Lie algebra.  As explained in \S\ref{sec:pontry-dual-vect}, we may identify $\mathfrak{g}^\wedge$ with the imaginary dual $i \mathfrak{g}^* = \Hom(\mathfrak{g}, i \mathbb{R})$ of $\mathfrak{g}$.  We write $\mathfrak{g}_\mathbb{C} := \mathfrak{g}\otimes_{\mathbb{R}} \mathbb{C}$ for the complexification of the Lie algebra and $\mathfrak{g}_\mathbb{C}^*$ for the complex dual, which identifies naturally with the complexification of $\mathfrak{g}^\wedge$.

We recall that $\mathfrak{U}(G)$ and $\mathfrak{Z}(G)$ denote the universal enveloping algebra and center for $\mathfrak{g}_{\mathbb{C}}$.  We similarly define \index{Lie algebra!enveloping algebra $\mathfrak{U}$, center $\mathfrak{Z}$} $\mathfrak{U}(\mathfrak{h})$, $\mathfrak{U}(\mathfrak{h}_{\mathbb{C}})$ for any Lie algebra $\mathfrak{h}$.

\subsection{Preliminaries on real reductive groups}\label{sec:prel-real-reduct}
Let $G$ be a connected reductive group over $\mathbb{R}$.  In particular, we may regard $G$ as a Lie group over $\mathbb{R}$.

\subsubsection{Regular elements}\label{sec:regular-elements}
A subscripted $\reg$, as in $\mathfrak{g}_{\reg}$ or $\mathfrak{g}^\wedge_{\reg}$, denotes the subset of \emph{regular} elements, i.e., those with $\mathfrak{g}$-centralizer of minimal dimension.\index{Lie algebra!regular subset $\mathfrak{g}^\wedge_{\reg}$}

\subsubsection{Invariant-theoretic quotients}
We denote by $[\mathfrak{g}_\mathbb{C}^*]$ the geometric invariant theory \index{Lie algebra!GIT quotient $[\mathfrak{g}^\wedge]$} quotient $\mathfrak{g}_\mathbb{C}^* \git G$, i.e., the spectrum of the ring of $G$-invariant polynomials $\mathfrak{g}_\mathbb{C}^* \rightarrow \mathbb{C}$, and $[\mathfrak{g}^\wedge]$ for its real form corresponding to polynomials that are real-valued on $\mathfrak{g}^\wedge$.  These come with natural maps
\[
  \mathfrak{g}_\mathbb{C}^* \rightarrow [\mathfrak{g}_\mathbb{C}^*], \quad \mathfrak{g}^\wedge \rightarrow [\mathfrak{g}^\wedge],
\]
which we denote by $\xi \mapsto [\xi]$.

There is a natural scaling action of $\mathbb{R}$ on $[\mathfrak{g}^\wedge]$ and of $\mathbb{C}$ on $[\mathfrak{g}_{\mathbb{C}}^*]$, compatible with the vector space operations of $\mathbb{R}$ on $\mathfrak{g}^\wedge$ and
of
$\mathbb{C}$
on
$\mathfrak{g}_{\mathbb{C}}^*$.

\subsubsection{Harish--Chandra isomorphism}\label{sec:cui1fnha7b}
The Harish--Chandra isomorphism (see, e.g., \cite[\S9.4]{nelson-venkatesh-1}) identifies the center $\mathfrak{Z}(\mathfrak{g}_\mathbb{C})$ of the universal enveloping algebra with the ring $\Sym(\mathfrak{g}_{\mathbb{C}})^G$ of regular functions on $[\mathfrak{g}_\mathbb{C}^*]$.  We recall its construction.  Choose a Cartan subalgebra $\mathfrak{t} \leq \mathfrak{g}_\mathbb{C}$ and a decomposition $\mathfrak{g}_\mathbb{C} = \mathfrak{n}_{\mathbb{C}}^- \oplus \mathfrak{t}_\mathbb{C} \oplus \mathfrak{n}_{\mathbb{C}}$ into positive and negative root spaces.  We may identify $\mathfrak{U}(\mathfrak{t}_\mathbb{C})$ with the symmetric algebra $\Sym(\mathfrak{t}_\mathbb{C})$.  Then
\begin{equation}\label{eqn:decomposition-for-HC}
  \mathfrak{U}(\mathfrak{g}_\mathbb{C})
  = (\mathfrak{n}_{\mathbb{C}}^- \mathfrak{U}(\mathfrak{g}_{\mathbb{C}})
  + \mathfrak{U}(\mathfrak{g}_{\mathbb{C}})
  \mathfrak{n}_{\mathbb{C}})
  \oplus \Sym(\mathfrak{t}_\mathbb{C}).
\end{equation}
Let $\rho \in \mathfrak{t}_\mathbb{C}^*$ denote the half-sum of positive roots.  Let $\sigma$ denote the algebra automorphism of $\Sym(\mathfrak{t}_\mathbb{C})$ given on $\mathfrak{t}_\mathbb{C}$ by $t \mapsto t - \rho(t)$.  The Harish--Chandra isomorphism is then the composition
\[
  \mathfrak{Z}(\mathfrak{g}_\mathbb{C}) \xrightarrow{\HACH} \Sym(\mathfrak{t}_\mathbb{C})^W
  \cong
  \Sym(\mathfrak{g}_{\mathbb{C}})^G \cong \{\text{regular functions on } [\mathfrak{g}_\mathbb{C}^*]\}
\]
\[
  \HACH(z) := \sigma(z_T),
\]
where $z_T \in \Sym(\mathfrak{t}_\mathbb{C})$ denotes the component of $z$ with respect to the decomposition \eqref{eqn:decomposition-for-HC}.

\subsubsection{Infinitesimal characters}\label{sec:infin-char}
Let $\pi$ be an irreducible admissible representation of $G$.  Then $\mathfrak{Z}(\mathfrak{g}_{\mathbb{C}})$ acts on $\pi$ by scalars.  Using the Harish--Chandra isomorphism, we may thus associate to $\pi$ an \index{representations!infinitesimal character $\lambda_\pi$} \emph{infinitesimal character} $\lambda_\pi \in [\mathfrak{g}^*_{\mathbb{C}}]$, characterized as follows: $z \in \mathfrak{Z}(\mathfrak{g}_{\mathbb{C}})$ acts on $\pi$ as multiplication by $\HACH(z)(\lambda_\pi)$, where we identify $\HACH(z)$ with a regular function on $[\mathfrak{g}^*_{\mathbb{C}}]$ as above.  If $\pi$ is unitary, then $\lambda_\pi \in [\mathfrak{g}^\wedge]$ (see \cite[\S9.5]{nelson-venkatesh-1}).

\subsubsection{Dual groups}
We write $G^\vee$ for the complex dual group and $\mathfrak{g}^\vee$ for its complex Lie algebra.  For any Cartan subalgebras $\mathfrak{t}_\mathbb{C} \leq \mathfrak{g}_\mathbb{C}$ and $\mathfrak{t}^\vee \leq \mathfrak{g}^\vee$, we may canonically identify
\begin{equation}\label{eqn:inf-chars-via-langlands-params}
  \mathfrak{g}^\vee \git G^\vee
  \cong \mathfrak{t}^\vee / W
  \cong \mathfrak{t}_\mathbb{C}^*/W \cong [\mathfrak{g}_\mathbb{C}^*],
\end{equation}
where $W$ denotes the Weyl group (see \cite[\S15.1]{nelson-venkatesh-1}).

\subsection{General linear groups}\label{sec:gener-line-groups}
We may specialize the above discussion to the case of a general linear group $G$ over the reals, say $G = \GL_n(\mathbb{R})$.  Then $\mathfrak{g}$ is the space $\glLie_n(\mathbb{R})$ of $n \times n$ real matrices.  We have $G^\vee = \GL_n(\mathbb{C})$ and $\mathfrak{g}^\vee = \glLie_n(\mathbb{C})$.

\subsubsection{Characteristic polynomials}\label{sec:char-polyn}
We denote by $e : \glLie_n(\mathbb{R}) \rightarrow \Sym(\mathfrak{g})$ the ``identity map'' coming from the identity $\glLie_n(\mathbb{R}) = \mathfrak{g}$ and the inclusion $\mathfrak{g} \hookrightarrow \Sym(\mathfrak{g})$.  We regard $e$ as a matrix $(e_{i j})$, with entries $e_{i j} \in \mathfrak{g} \hookrightarrow \Sym(\mathfrak{g})$, and form the characteristic polynomial
\begin{align}\label{eq:detx-+-e}
  \det(X + e)
  &:=
    \sum_{\sigma \in S(n)}
    (-1)^\sigma
    \prod_{i=1}^n
    (1_{i=\sigma(i)} X
    + e_{i,\sigma(i)})
  \\ \nonumber
  &=
    X^n +
    \mathfrak{c}_1 X^{n-1}
    + \dotsb + \mathfrak{c}_n,
\end{align}
which defines a polynomial in the variable $X$ with coefficients $\mathfrak{c}_j$ in the ring $\Sym(\mathfrak{g})^G$ of $G$-invariant elements of $\Sym(\mathfrak{g})$.  We may identify each $\mathfrak{c}_j$ with a $G$-invariant polynomial $\mathfrak{g}_\mathbb{C}^* \rightarrow \mathbb{C}$, hence with a regular function $[\mathfrak{g}_\mathbb{C}^*] \rightarrow \mathbb{C}$.  The resulting map
\[
  \det(X + e) : [\mathfrak{g}_\mathbb{C}^*] \rightarrow \{ \text{monic degree $n$ polynomials over $\mathbb{C}$} \}
\]
is an isomorphism.  For $\lambda \in [\mathfrak{g}_{\mathbb{C}}^*]$, we denote by $\mathcal{P}_\lambda(X) \in \mathbb{C}[X]$ the corresponding polynomial, i.e., \index{Lie algebra!characteristic polynomial $\mathcal{P}_\lambda(X)$}
\begin{equation}\label{eq:polynomial-P-lambda}
  \mathcal{P}_\lambda(X) := \det(X+e)(\lambda) = X^n + \mathfrak{c}_1(\lambda) X^{n-1} + \dotsb + \mathfrak{c}_n(\lambda).
\end{equation}

We have noted that the space $\mathfrak{g}^\wedge$ identifies (canonically) with the imaginary dual of $\mathfrak{g}$,  hence, via the trace pairing, with the space $i \glLie_n(\mathbb{R})$ of $n \times n$ imaginary matrices:
\[
  i \glLie_n(\mathbb{R}) \cong \mathfrak{g}^\wedge
\]
\[
  \xi \mapsto [x \mapsto \trace(x \xi)].
\]
With respect to this identification, we see that for $\xi \in \mathfrak{g}^\wedge$ with image $[\xi] \in [\mathfrak{g}^\wedge]$,
\begin{equation}\label{eq:cui1fo922g}
  \det(X + \xi) = \mathcal{P}_{[\xi]}(X) = \det(X + e)(\xi).
\end{equation}
In other words, the map $\xi \mapsto \mathcal{P}_{[\xi]}(X)$ associates to $\xi$ its characteristic polynomial.

\subsubsection{Eigenvalues}
Let $\lambda \in [\mathfrak{g}_\mathbb{C}^*]$.  Under \eqref{eqn:inf-chars-via-langlands-params}, we may identify $\lambda$ with an element of $\mathfrak{t}^\vee/W$.  We write $t_\lambda \in \mathfrak{t}^\vee$ for any lift of that element, well-defined up to $W$.  We write
\begin{equation*}
  \ev(\lambda) = \{\lambda_1,\dotsc,\lambda_n\}
\end{equation*}
for the multiset \index{Lie algebra!eigenvalue multiset $\ev(\lambda)$} of eigenvalues of $t_\lambda$.  (This definition is consistent with that of \cite[\S13.4.1]{nelson-venkatesh-1}.)  The characteristic polynomial $\det(X + t_\lambda) \in \mathbb{C}[X]$ then coincides with $\mathcal{P}_\lambda(X)$, as defined above:
\begin{equation}\label{eq:detx-+-t_lambda}
  \det(X + t_\lambda) = \mathcal{P}_{\lambda}(X).
\end{equation}

\subsubsection{Characterizations of $[\mathfrak{g}^\wedge] \subseteq [\mathfrak{g}_{\mathbb{C}}^*]$}
We observe that the following conditions on $\lambda \in [\mathfrak{g}_\mathbb{C}^*]$ are equivalent to each other.
\begin{itemize}
\item $\lambda \in [\mathfrak{g}^\wedge]$.
\item $- \overline{t_\lambda} \in W \cdot t_\lambda$.
\item $\mathfrak{c}_1(\lambda), \mathfrak{c}_3(\lambda), \mathfrak{c}_5(\lambda), \dotsc \in i \mathbb{R}$, $\mathfrak{c}_2(\lambda), \mathfrak{c}_4(\lambda), \mathfrak{c}_6(\lambda), \dotsc \in \mathbb{R}$.
\end{itemize}
Thus the map $\lambda \mapsto \mathcal{P}_\lambda$ defines a bijection between $[\mathfrak{g}^\wedge]$ and the space of polynomials $X^n + c_1 X^{n-1} + \dotsb + c_n$ with $c_{2 j} \in \mathbb{R}$ and $c_{2j - 1} \in i \mathbb{R}$.

\subsubsection{Scaling}
The scaling action of $\mathbb{R}$ on $[\mathfrak{g}^\wedge]$ (or of $\mathbb{C}$ on $[\mathfrak{g}_{\mathbb{C}}^*]$) may be described explicitly as follows.  For $(t,\lambda) \in \mathbb{R} \times [\mathfrak{g}^\wedge]$, the scaled element $t \lambda \in [\mathfrak{g}^\wedge]$ is characterized by
\begin{equation}\label{eq:cal-P-t-lambda}
  \mathcal{P}_{t \lambda}(X) = t^{n} \mathcal{P}_{\lambda}(t^{-1} X)
  =
  X^n + t \mathfrak{c}_1(\lambda) X^{n-1} + \dotsb + t^n \mathfrak{c}_n(\lambda).
\end{equation}
Equivalently, the roots of $\mathcal{P}_{t \lambda}$ are obtained by multiplying those of $\mathcal{P}_{\lambda}$ by $t$.

\subsubsection{Infinitesimal characters and parameters}
Let $\pi$ be an irreducible admissible representation of $G$.  It has an infinitesimal character $\lambda_\pi \in [\mathfrak{g}_{\mathbb{C}}^*]$, which, as we have noted (\S\ref{sec:infin-char}), belongs to $[\mathfrak{g}^\wedge]$ when $\pi$ is unitary.  We may in any case associate to $\pi$ a semisimple element $t_{\lambda_\pi}$ and a multiset $\ev(\lambda_{\pi}) = \{\lambda_{\pi,1},\dotsc, \lambda_{\pi,n}\}$ of complex numbers.  We refer to these numbers as the ``archimedean Satake parameters'' of $\pi$. \index{representations!archimedean Satake parameters $\lambda_{\pi,i}$}

By the discussion of \cite[\S15.4]{nelson-venkatesh-1}, the standard local $L$-factor for $\pi$ is given in terms of the infinitesimal character by
\[
  L(\pi,s) =
  \prod_j \Gamma_\mathbb{R}(s + \lambda_{\pi,j}^+ + a_j),
\]
where the elements $a_j \in \{0, 1\}$ depend upon $\pi$ and we set
\[
  (x + i y)^+ := |x| + i y, \quad \Gamma_\mathbb{R}(s) := \pi^{-s/2} \Gamma(s/2), \quad \Gamma_{\mathbb{C}}(s) := 2 (2 \pi)^{-s} \Gamma(s).
\]
(Strictly speaking, it is assumed in parts of \emph{loc.\ cit.}\ that $\pi$ is tempered, but the same argument applies in general, using the
Langlands decomposition \cite[Theorem 8.54]{MR855239} and the compatibility of infinitesimal characters with induction \cite[Proposition 8.22]{MR855239} to reduce to the tempered case.  See also \cite[Appendix]{MR1395406}.)

In \S\ref{sec:subc-bounds-assum}, we referred to the numbers $t_j = \lambda_{\pi,j}^+ + a_j$ as the ``archimedean $L$-function parameters'' of $\pi$ (or more precisely, of a global automorphic representation having $\pi$ as its local component at $\infty$).  These are not quite the same as the ``archimedean Satake parameters'' $\lambda_{\pi,j}$ described above, but they have very similar asymptotic behavior.  For instance:
\begin{lemma}\label{lem:standard:let-t-geq}
  Let $T \geq 1$.  Let $\pi$ be an irreducible admissible representation of $G$.  The following are equivalent:
  \begin{enumerate}[(i)]
  \item \label{itm:standard2:compare-parameters-1} The rescaled infinitesimal character $T^{-1} \lambda_\pi$ lies in some fixed compact subset of $[\mathfrak{g}_{\mathbb{C}}^*]$.
  \item \label{itm:standard2:compare-parameters-2} Each archimedean Satake parameter  $\lambda_{\pi,j}$ is of the form $\O(T)$.
  \item  \label{itm:standard2:compare-parameters-3} Each archimedean $L$-function parameter $\lambda_{\pi,j}^+ + a_j$ is of the form $\O(T)$.
  \end{enumerate}
\end{lemma}
\begin{proof}
  By \eqref{eq:detx-+-t_lambda} and \eqref{eq:cal-P-t-lambda},
  \begin{equation*}
    \mathcal{P}_{T^{-1} \lambda_\pi}(X) = \prod_j (X + T^{-1} \lambda_{\pi,j}).
  \end{equation*}
  The equivalence between \eqref{itm:standard2:compare-parameters-1} and \eqref{itm:standard2:compare-parameters-2}
  follows from the fact that a degree $n$ monic polynomial has coefficients of size $\O(1)$ if and only if its (negated) roots are of size $\O(1)$.
  The equivalence between \eqref{itm:standard2:compare-parameters-2} and \eqref{itm:standard2:compare-parameters-3} follows from the fact that $a_j \in \{0,1\}$ and $|z| = |z^+|$.
\end{proof}

\subsection{General linear pairs}
Let $(G,H)$ be a general linear pair (\S\ref{sec:gener-line-groups-1}) over $\mathbb{R}$, say $(G,H) = (\GL_n(\mathbb{R}), \GL_{n-1}(\mathbb{R}))$.  We recall that $H$ is embedded in $G$ via the upper-left corner map $h \mapsto \left(
  \begin{smallmatrix}
    h&0\\
    0&1 \\
  \end{smallmatrix}
\right)$.  The above discussion applies both to $\mathfrak{g} = \Lie(G)$ and to $\mathfrak{h} = \Lie(H)$.

For $\xi \in \mathfrak{g}^\wedge$, we let $\xi_H \in \mathfrak{h}^\wedge$ denote its restriction.  We note that, if we identify each of these Lie algebras and their duals with spaces of matrices via the trace pairing, then $\xi_H$ is the upper-left corner of $\xi$, i.e., $\xi = \left(
  \begin{smallmatrix}
    \xi_H&\ast\\
    \ast       &\ast \\
  \end{smallmatrix}
\right)$.

\subsubsection{Stability}\label{sec:stability}

\begin{definition}
  We say that a pair $(\lambda, \mu ) \in [\mathfrak{g} ^\wedge ] \times [\mathfrak{h} ^\wedge ]$ is \emph{stable} if the eigenvalue multisets $\eval(\lambda)$ and $\eval(\mu)$ have no common element.
\end{definition}

\begin{definition}
  We say that $\xi \in \mathfrak{g}^\wedge$ is \emph{stable}   \index{Lie algebra!stable subset $\mathfrak{g}^\wedge_{\stab}$} if it is stable with respect to $H$ in the sense of geometric invariant theory, i.e., if it has finite $H$-stabilizer and (Zariski) closed $H$-orbit.  We denote by $\mathfrak{g}^\wedge_{\stab} \subseteq \mathfrak{g}^\wedge$ the subset consisting of stable elements.
\end{definition}

The two definitions are related as follows (see, e.g., \cite[Lemma 11.4]{2020arXiv201202187N}):
\begin{lemma}
  The stable subset $\mathfrak{g}^\wedge_{\stab}$ of $\mathfrak{g}^\wedge$ is the preimage of $\{\text{stable } (\lambda, \mu ) \in [\mathfrak{g}^\wedge] \times [\mathfrak{h}^\wedge]\}$ under the map $\xi \mapsto ([\xi], [\xi_H])$.
\end{lemma}
In particular, $\mathfrak{g}^\wedge_{\stab}$ is the dense open subset of $\mathfrak{g}^\wedge$ given by the nonvanishing of the polynomial
\begin{equation*}
  \mathfrak{g}^\wedge  \ni \xi \mapsto \mathcal{R}(\xi) := \operatorname{Resultant} (\mathcal{P}_{[\xi]}, \mathcal{P}_{[\xi_H]}).
\end{equation*}
Here $[\xi] \in [\mathfrak{g}^\wedge] \subseteq [\mathfrak{g}_{\mathbb{C}}^*]$ and $[\xi_H] \in [\mathfrak{h}^\wedge] \subseteq [\mathfrak{h}_{\mathbb{C}}^*]$ give rise to polynomials $\mathcal{P}_{[\xi]}, \mathcal{P}_{[\xi_H]}$ via \eqref{eq:polynomial-P-lambda}.

We obtain the following equivalence:
\begin{lemma}\label{lem:standard2:stability-vs-resultant}
  Let $\tau$ belong to some fixed compact subset of $\mathfrak{g}^\wedge$.  Then the following are equivalent:
  \begin{enumerate}[(i)]
  \item $\tau$ belongs to some fixed compact subset of $\mathfrak{g}^\wedge_{\stab}$.
  \item $\mathcal{R}(\xi) \gg 1$.
  \end{enumerate}
\end{lemma}

\subsubsection{The parameter $\theta$}\label{sec:parameter-theta}
Let $P$ denote the mirabolic subgroup of $G$, consisting of elements with bottom row $(0,\dotsc,0,1)$.  Recall that $\bar{B}_H$ denotes the lower-triangular Borel subgroup of $H$.  The Lie algebra of the mirabolic subgroup then decomposes as
\begin{equation*}
  \Lie(P) = \Lie(N) \oplus \Lie(\bar{B}_H).
\end{equation*}

Let $\psi$ be a nondegenerate unitary character of $N$.  We denote by \index{Lie algebra!$\theta_P(\psi)$}
\begin{equation*}
  \theta_P(\psi)  \in \Lie(P)^\wedge
\end{equation*}
the unique element with trivial restriction to $\Lie(\bar{B}_H)$ that satisfies, for all $x \in \Lie(N)$,
\begin{equation*}
  \psi(\exp(x)) = \exp (\langle x, \theta_P(\psi)  \rangle).
\end{equation*}
We recall that, by Definition \ref{defn:standard2:given-point-manif}, we may associate to each $\xi \in \Lie(P)^\wedge$ and chart $c : (P,1) \xdashrightarrow{} (\Lie(P),0)$ a function $\psi_{\xi, c} : P \xdashrightarrow{} \U(1)$, defined near the identity element.  (The use of ``$\psi$'' in this notation is purely symbolic, unrelated to the character $\psi$ under consideration.)  In particular, we may define $\psi_{\theta_P(\psi),\log}$  with respect to the chart $\log :  (P, 1) \xdashrightarrow{} (\Lie(P), 0)$.  Then
\begin{equation*}
  \psi(x) = \psi_{\theta_P(\psi), \log}(x) \quad \text{ for } x \in N,
\end{equation*}
while $\psi_{\theta_P(\psi), \log}(x) = 1$ for $x \in \bar{B}_H$.

Using the trace form, we may identify $\Lie(P)^\wedge$ with a quotient of $\mathfrak{g}^\wedge = i \glLie_n(\mathbb{R})$, e.g., when $n=4$,
\begin{equation*}
  \Lie(P)^\wedge \leftrightarrow
  \begin{pmatrix}
    \ast & \ast & \ast & ? \\
    \ast & \ast & \ast & ? \\
    \ast & \ast & \ast & ? \\
    \ast & \ast & \ast & ?
  \end{pmatrix}.
\end{equation*}
The meaning of this notation is that $\Lie(P)^\wedge$ identifies with the space of equivalence clases of $n \times n$ imaginary matrices, where two such matrices are declared equivalent if they have the same entries except possibly in their rightmost columns.  We may write
\begin{equation}\label{eq:psiu-=-sum_i=1n}
  \psi(u) = \sum_{i=1}^{n-1} \exp(\eta_i u_{i,i+1})
\end{equation}
for some nonzero imaginary numbers $\eta_1,\dotsc,\eta_{n-1}$ (e.g., each $\eta_j = 2 \pi \sqrt{-1}$ for a suitable sign).  Then, under the above identification,
\begin{equation}\label{eq:theta-leftr-beginpm}
  \theta_P(\psi)
  \leftrightarrow
  \begin{pmatrix}
    0 & 0 & 0 & ? \\
    \eta_1  & 0 & 0 & ? \\
    0 & \eta_2  & 0 & ? \\
    0 & 0 & \eta_3  & ?
  \end{pmatrix}.
\end{equation}

\subsubsection{The parameter $\tau$}\label{sec:parameter-tau}
For $\xi \in \mathfrak{g}^\wedge$, we denote by $\xi_P \in \Lie(P)^\wedge$ its restriction.  We have the following simple linear algebra lemma, closely related to ``rational canonical form:''
\begin{lemma}\label{lem:construction-tau-from-theta}
  For each nondegenerate unitary character $\psi$ of $N$ and $\lambda \in [\mathfrak{g}^\wedge]$, there is a unique element $\tau = \tau(\psi,\lambda) \in \mathfrak{g}^\wedge$ such that \index{Lie algebra!$\tau(\psi,\lambda)$}
  \[
    [\tau] = \lambda, \quad \tau_P = \theta_P(\psi).
  \]
  Moreover, $\tau$ may be given explicitly with respect to the identification $\mathfrak{g}^\wedge = i \glLie_n(\mathbb{R})$ by
  \[
    \tau_{i j} =
    \begin{cases}
      \eta_j   & \text{ if } i = j+1, \\
      (-\eta_n)^{-1} \dotsb (-\eta_{i+1})^{-1} \mathfrak{c}_{n+1-i}
            & \text{ if } j = n, \\
      0 & \text{ otherwise}.
    \end{cases}
  \]
\end{lemma}
For example, when $n=4$ and $\psi$ is as in \eqref{eq:psiu-=-sum_i=1n}, the element $\tau$ is given by
\begin{equation}\label{eqn:tau-three-by-three-theta-lambda}
  \tau
  =
  \begin{pmatrix}
    0 & 0 & 0 &  - \mathfrak{c}_4(\lambda)  / \eta_1 \eta_2 \eta_3 \\
    \eta_1  & 0 & 0 & \mathfrak{c}_3(\lambda)  / \eta_2 \eta_3\\
    0 & \eta_2 & 0 & -\mathfrak{c}_2(\lambda)/ \eta_3 \\
    0 & 0 & \eta_3 & \mathfrak{c}_1(\lambda)
  \end{pmatrix}.
\end{equation}
\begin{proof}
  It is clear that $\tau$ as in \eqref{eqn:tau-three-by-three-theta-lambda} satisfies $\tau_P = \theta_P(\psi)$, and a direct calculation gives $\det(X + \tau) = \mathcal{P}_\lambda(X)$.

  As for uniqueness, the condition $\tau_P = \theta_P(\psi)$ implies that $\theta_P(\psi)$ looks like \eqref{eqn:tau-three-by-three-theta-lambda}, except possibly with a different rightmost column; to determine that column, we evaluate the characteristic polynomial $\det(X + \tau)$ and compare its coefficients with those of $\mathcal{P}_\lambda(X)$.
\end{proof}

Lemma \ref{lem:construction-tau-from-theta} allows us to make the following definition.
\begin{definition}\label{defn:parameter-tau-pi-T}
  Let $\psi$ be a nondegenerate unitary character of $N$.  Let $T \geq 1$.  Let $\pi$ be an irreducible unitary representation of $G$.  We denote by $\tau(\pi,\psi,T) \in \mathfrak{g}^\wedge$ the unique element with \index{Lie algebra!$\tau(\pi,\psi,T)$}
  \begin{equation}\label{eq:tau-characterize-via-pi}
    [\tau(\pi,\psi,T)] = T^{-1} \lambda_\pi,
    \quad
    \tau(\pi,\psi,T)_P = \theta_P(\psi).
  \end{equation}
\end{definition}

\begin{lemma}
  In the setting of Definition \ref{defn:parameter-tau-pi-T}, assume that $\psi$ is fixed and that $T^{-1} \lambda_\pi$ belongs to some fixed compact subset of $[\mathfrak{g}^\wedge]$.  Then $\tau(\pi,\psi,T)$ belongs to some fixed compact subset of $\mathfrak{g}^\wedge_{\reg}$.
\end{lemma}
\begin{proof}
  Fix $C \geq 0$ so that each $|\mathfrak{c}_j (T^{-1} \lambda_\pi)| \leq C$.  Let $\mathcal{D} \subseteq \mathfrak{g}^\wedge$ denote the set of matrices of the form, e.g., for $n=4$,
  \begin{equation*}
    \begin{pmatrix}
      0 & 0 & 0 & \ast \\
      \eta_1  & 0 & 0 & \ast \\
      0 & \eta_2  & 0 & \ast \\
      0 & 0 & \eta_3  &  \ast
    \end{pmatrix},
  \end{equation*}
  where each entry $\ast$ has magnitude bounded by $C$.  Then $\tau(\pi,\psi,T) \in \mathcal{D}$.   Every element of $\mathcal{D}$ is regular (e.g., because the standard basis element $e_1$ is a cyclic vector for such elements), and $\mathcal{D}$ is fixed and compact, thus $\mathcal{D}$ is a fixed compact subset of $\mathfrak{g}^\wedge_{\reg}$.
\end{proof}

The following result refines the equivalence between upper bounds given in Lemma \ref{lem:standard:let-t-geq}:
\begin{lemma}\label{lem:standard2:no-conductor-drop-vs-tau}
  Fix a nondegenerate unitary character $\psi$ of $N$.  Let $T \ggg 1$.  Let $\pi$ be an irreducible unitary representation of $G$.  Assume that $T^{-1} \lambda_\pi$ belongs to some fixed compact subset of $[\mathfrak{g}^\wedge]$.  The following are then equivalent.
  \begin{enumerate}[(i)]
  \item
    \label{itm:standard2:characterize-arch-parameters-conductor-drop-1} $\tau(\pi,\psi,T)$ belongs to some fixed compact subset of $\mathfrak{g}^\wedge_{\stab}$.
  \item \label{itm:standard2:characterize-arch-parameters-conductor-drop-2} Each archimedean Satake parameter of $\pi$ satisfies $\lambda_{\pi,j} \asymp T$.
  \item \label{itm:standard2:characterize-arch-parameters-conductor-drop-3} Each archimedean $L$-function parameter of $\pi$ satisfies $\lambda_{\pi,j}^+ + a_j \asymp T$.
  \end{enumerate}
\end{lemma}
\begin{proof}
  The equivalence of \eqref{itm:standard2:characterize-arch-parameters-conductor-drop-2} and \eqref{itm:standard2:characterize-arch-parameters-conductor-drop-3} follows from our assumption $T \ggg 1$ and the fact that $a_j \in \{0,1\}$.  To see the equivalence between \eqref{itm:standard2:characterize-arch-parameters-conductor-drop-1} and \eqref{itm:standard2:characterize-arch-parameters-conductor-drop-2}, we note that, writing $\tau = \tau(\pi,\psi,T)$,
  \begin{equation*}
    \mathcal{P}_{[\tau]} = \prod_j (X + T^{-1} \lambda_{\pi,j}),
    \quad
    \mathcal{P}_{[\tau_H]}(X) = X^{n-1},
  \end{equation*}
  where the second identity follows from the fact that $\tau_H$ is lower-triangular nilpotent, e.g., for $n = 4$ and $\tau$ as in \eqref{eqn:tau-three-by-three-theta-lambda},
  \begin{equation*}
    \tau_H =
    \begin{pmatrix}
      0      & 0 & 0 \\
      \eta_1 & 0 & 0 \\
      0 & \eta_2 & 0 \\
    \end{pmatrix}.
  \end{equation*}
  It follows that
  \begin{equation*}
    \mathcal{R}(\tau) = \pm \prod_{j=1}^n (T^{-1} \lambda_{\pi,j})^{n-1}.
  \end{equation*}
  In view of our assumption that each $T^{-1} \lambda_{\pi, j} = \O(1)$, we deduce that
  \begin{equation*}
    \mathcal{R}(\tau) \gg 1 \iff \text{ each } T^{-1} \lambda_{\pi,j} \asymp 1.
  \end{equation*}
  The required equivalence then follows from Lemma \ref{lem:standard2:stability-vs-resultant}.
\end{proof}

\section{Localized elements of the Kirillov model}\label{sec:asympt-analys-local}
Let $(G,H)$ be a fixed general linear pair over $\mathbb{R}$.  Recall (\S\ref{sec:gener-line-groups-1}) that we write $N \leq G$ and $N_H \leq H$ for the standard upper-triangular maximal unipotent subgroups.  Recall also that $B \geq N$ is the Borel subgroup containing $N$, that $B_H$ and $N_H$ denote the subgroups of $H$ obtained via intersection, and that $Q_H$ is a Borel subgroup of $H$, chosen according to convenience.  In this section, we take for $Q_H$ the opposite of $B_H$:
\begin{equation*}
  Q_H := \bar{B}_H.
\end{equation*}

We denote as in \S\ref{sec:parameter-theta} by $P \leq G$ the mirabolic subgroup, and by $U_P \leq P$ its unipotent radical.  For illustration, if $(G, H) =(\mathrm{GL}_3(\mathbb{R}), \mathrm{GL}_2(\mathbb{R}))$, then
\begin{equation*}
  N =
  \begin{pmatrix}
    1 & \ast & \ast \\
    0 & 1 & \ast \\
    0 & 0 & 1 \\
  \end{pmatrix},
  \quad
  B =
  \begin{pmatrix}
    \ast & \ast & \ast \\
    0 & \ast & \ast \\
    0 & 0 & \ast \\
  \end{pmatrix},
  \quad
  Q_H =
  \begin{pmatrix}
    \ast & 0 & 0 \\
    \ast & \ast & 0 \\
    0 & 0 & 1 \\
  \end{pmatrix},
\end{equation*}
\begin{equation*}
  N_H
  =
  \begin{pmatrix}
    1 & \ast & 0 \\
    0 & 1 & 0 \\
    0 & 0 & 1 \\
  \end{pmatrix},
  \quad
  P =
  \begin{pmatrix}
    \ast & \ast & \ast \\
    \ast & \ast & \ast \\
    0 & 0 & 1 \\
  \end{pmatrix},
  \quad
  U_P =
  \begin{pmatrix}
    1 & 0 & \ast \\
    0 & 1 & \ast \\
    0 & 0 & 1 \\
  \end{pmatrix}.
\end{equation*}
The motivation for choosing $Q_H$ as we have is that the maps
\begin{equation}\label{eq:crb6fpkcd0}
  N_H \times Q_H \rightarrow H
  \quad \text{ and } \quad
  N \times Q_H \rightarrow P,
\end{equation}
given in either case by $(n,g) \mapsto g^{-1} n$, define open immersions with dense image, and in particular, local analytic isomorphisms in some neighborhood of the identity.  We may thus regard $Q_H$ as an approximate fundamental domain for the quotients $N_H \backslash H$ and $N \backslash P$.


Let $\pi$ be a generic irreducible unitary representation of $G$.

The main result of this section is Theorem \ref{thm:big-group-whittaker-behavior-under-G}, stated below in \S\ref{sec:behavior-under-g}.  Informally speaking,
that result supplies us with a class of microlocalized vectors $v$ in $\pi$ that concentrate in the Kirillov model near a specific
point, together with some quantitative control over how the group $G$ acts on such vectors.  A central role is played by convolution operators $\pi(\omega)$ attached to test functions $\omega$ whose properties were described heuristically in \S\ref{sec:relat-trace-form}.

\subsection{The Kirillov model}\label{sec:kirillov-model}
Let $\psi$ be a nondegenerate unitary character of $N$.  We denote also simply by $\psi$ the restriction of $\psi$ to $N_H$, which is again nondegenerate.

We denote in what follows by $\pi^0$ the underlying Hilbert space, by $\pi^\infty$ the space of smooth vectors, and by $\pi^{-\infty}$ the space of distributional vectors, i.e., the continuous dual of $\pi^\infty$; these spaces are related via  natural inclusions
\begin{equation*}
  \pi^\infty \subseteq \pi^0 \subseteq \pi^{-\infty}.
\end{equation*}

By the definition of ``generic,'' there is a nonzero distributional vector \index{representations!Whittaker vector $\Theta_{\pi,\psi}$} $\Theta_{\pi,\psi} \in \pi^{-\infty}$, unique up to scalars \cite{MR348047}, for which $\pi(n) \Theta_{\pi,\psi} = \psi(n) \Theta_{\pi,\psi}$ for all $n \in N$.  The map $v \mapsto [g \mapsto \langle \pi(g) v, \Theta_{\pi,\psi} \rangle]$ defines the $\psi$-Whittaker model $\mathcal{W}(\pi,\psi)$ of $\pi$, which consists of smooth $\pi$-isotypic $W : G \rightarrow \mathbb{C}$ satisfying $W(n g) = \psi(n) W(g)$ for all $(n,g) \in N \times G$.  It follows from the results of \cite{MR1999922} (see also \cite{MR748505} and \cite[Ch. III]{MR0579172}) that by restricting elements of the model $\mathcal{W}(\pi,\psi)$ to the subgroup $H$, we may realize $\pi^0$ as the Hilbert space of all measurable functions $W : H \rightarrow \mathbb{C}$ satisfying
\begin{equation}\label{eqn:kirillov-model-defn}
  W(n h) = \psi(n) W(h)
  \text{ for all } (n,h) \in N_H \times H,
\end{equation}
\begin{equation}\label{eq:int_n_h-backslash-h}
  \int_{N_H \backslash H} |W|^2 < \infty,
\end{equation}
with an invariant inner product on $\pi$ given by integration over $N_H \backslash H$.  We recall (see \S\ref{sec:local-prelim-kirillov-model}) that the Kirillov model of $\pi$ is the image of $\mathcal{W}(\pi, \psi)$ under restriction to $H$.

In what follow, we focus primarily on the subspace $C_c^\infty(N_H \backslash H,\psi)$ of compactly-supported smooth functions in the Kirillov model, defined (by analogy to \S\ref{sec:whittaker-functions-1}) to consist of smooth functions $H \rightarrow \mathbb{C}$ satisfying \eqref{eqn:kirillov-model-defn} for which $\lvert W \rvert : N_H \backslash H \rightarrow \mathbb{R}_{\geq 0}$ has compact support.  The motivation for restricting in this way is that we are ultimately interested in the invariance properties of certain modulated bump functions.  We topologize this space (and analogously, other $C_c^\infty$ spaces: $C_c^\infty(H)$, $C_c^\infty(P)$, etc.) with its standard LF-topology, as the inductive limit over compact subsets $E \subseteq N_H \backslash H$ of the subspaces $C_c^\infty(E, \psi)$ consisting of elements supported on the preimage of $E$, where each $C_c^\infty(E, \psi)$ is a Fréchet space when equipped with the family of seminorms given by $f \mapsto \sup_{h \in H} \lvert  D f(h)  \rvert$, with $D$ ranging over the left-invariant differential operators on $H$.

It is known, and relatively straightforward to see, that the space $C_c^\infty(N_H \backslash H, \psi)$ is contained in the Hilbert space $\pi^0$ and is smooth under the action of $H$.  Jacquet \cite[Prop 5]{MR2733072} established the deeper fact that it is smooth under the action of $G$, i.e., that
\begin{equation*}
  C_c^\infty(N_H \backslash H, \psi) \subseteq \pi^\infty.
\end{equation*}
Jacquet's proof shows in particular that the inclusion is continuous and that for any bounded open subset $E$ of $N_H \backslash H$, the subspace $C_c^\infty(E, \psi)$ is invariant by $\mathfrak{U}(G)$.

In more detail, by the proof of \cite[Prop 5]{MR2733072}, we may represent each element of $C_c^\infty(N_H \backslash H, \psi)$ as the convolution
\begin{equation*}
  \pi(f) \Theta_{\pi,\psi} = \int_{h \in H} f(h) \pi(h) \Theta_{\pi, \psi} \, d h
\end{equation*}
of the distributional vector $\Theta_{\pi, \psi}$ by some $f \in C_c^\infty(H)$.  We recall the formula given in \cite[Lem 7]{MR2733072} for that convolution: for $h \in H$,
\begin{equation}\label{eq:pif-theta_pi-h}
  \pi(f) \Theta_{\pi,\psi} (h)
  =
  \int_{n \in N_H}
  f(h^{-1} n)
  \psi(n) \, d n.
\end{equation}
The following result describes the action of $\mathfrak{U}(G)$ on such convolutions, showing in particular that each subspace $C_c^\infty(E, \psi)$ as above is indeed invariant by $\mathfrak{U}(G)$.

\begin{theorem}[Jacquet]\label{thm:jacquet-on-kirillov-model}
  Let $f \in C_c^\infty(H)$.  Then $\pi(f) \Theta_{\pi,\psi} \in \pi^\infty$.  More precisely, for each $x \in \mathfrak{U}(G)$, we have
  \begin{equation}\label{eqn:pi-x-pi-f-Theta-pi-f-prime}
    \pi(x) \pi(f) \Theta_{\pi,\psi} = \pi(f') \Theta_{\pi,\psi},
  \end{equation}
  where $f' \in C_c^\infty(H)$ with $\supp(f') \subseteq \supp(f)$, and $f \mapsto f'$ is continuous.
\end{theorem}
\begin{proof}
  See \cite[\S3.2]{MR2733072}.
\end{proof}

We note that the above discussion may be reformulated using the bijection $N_H \backslash H = N \backslash P$.  In particular:
\begin{theorem}[Jacquet]\label{thm:jacquet-on-kirillov-model-P}
  Let $f \in C_c^\infty(P)$.  Then $\pi(f) \Theta_{\pi,\psi}$, defined using (say) a left Haar measure on $P$, lies in $\pi^\infty$, and the resulting map $C_c^\infty(P) \rightarrow \pi^\infty$ is continuous.
\end{theorem}
\begin{proof}
  \cite[Prop 4]{MR2733072} gives the same result, but with $P$ replaced by $H$.  Recall that $U_P$ denotes the unipotent radical of $P$.  The group $P$ is the semidirect product of its subgroups $H$ and $U_P$.  Using our given Haar measures $d u$ on $U_P$ and $d h$ on $H$, we define a left Haar measure $d p$ on $P$ by
  \begin{equation*}
    \int_{P} F(p) \, d p
    =
    \int_{H} \int_{U_P} F(h u) \, d u \, d h
  \end{equation*}
  for $F \in C_c(P)$.  Then, since $\pi(u)\Theta_{\pi,\psi}=\psi(u)\Theta_{\pi,\psi}$, we have
  \begin{align*}
    \pi(f) \Theta_{\pi, \psi}
    &=
      \int_{H} \int_{U_P}
      f(h u)
      \pi(h u)
      \Theta_{\pi, \psi}
      \, d u \, d h
    \\
    &=
      \int_{H} \int_{U_P}
      f(h u)
      \psi(u)
      \pi(h)
      \Theta_{\pi, \psi}
      \, d u \, d h \\
    &=
      \int_{H}
      f'(h)
      \pi(h)
      \Theta_{\pi, \psi} \, d h
      =
      \pi(f ') \Theta_{\pi, \psi},
  \end{align*}
  where $f' \in C_c^\infty(H)$ is defined by
  \begin{equation*}
    f'(h) := \int_{U_P} f(h u) \psi(u) \, d u.
  \end{equation*}
  The map $C_c^\infty(P) \rightarrow C_c^\infty(H)$, $f \mapsto f'$ is continuous.  The cited result gives the continuity of $f' \mapsto \pi(f') \Theta_{\pi,\psi}$, so the required continuity of $f \mapsto \pi(f) \Theta_{\pi, \psi}$ follows.
\end{proof}
We refer to \cite{MR3352025} for an analogous result concerning $\GL_n(\mathbb{C})$, which may be useful in generalizing the results of this paper.

\subsection{A class of bump functions in the Kirillov model}\label{sec:class-bump-functions}
We now take the nondegenerate unitary character $\psi$ of $N$ to be \emph{fixed}.  We define $\theta_P(\psi) \in \Lie(P)^\wedge$ as in \S\ref{sec:parameter-theta}.  We denote by $\theta_H(\psi) \in \Lie(H)^\wedge$ the restriction of $\theta_P(\psi)$.  \index{Lie algebra!$\theta_P(\psi)$} \index{Lie algebra!$\theta_H(\psi)$}

Let $T$ be a positive parameter with $T \ggg 1$.  We define the rescaled nondegenerate character $\psi_T$ of $N$ by requiring that, for $x \in \Lie(N)$,
\begin{equation*}
  \psi_T(\exp(x)) = \psi(\exp(T x)).
\end{equation*}
Thus $\psi_T(\exp(x)) = \psi_{T \theta_P(\psi), \log}(x)$.  Equivalently, $\psi_T(n) = \psi (T^{\rho_N^\vee} n T^{-\rho_N^\vee})$.  In the notation of \S\ref{sec:parameter-theta}, we pass from $\psi$ to $\psi_T$ by replacing each $\eta_i$ by its multiple $T \eta_i$.

In this section, we identify $\pi$ with its $\psi_T$-Whittaker model $\mathcal{W}(\pi,\psi_T)$.  We identify the latter, in turn, with the $\psi_T$-Kirillov model, the subspace of $C^\infty(N_H \backslash H, \psi_T)$ given by the image of $\mathcal{W}(\pi, \psi_T)$ under restriction.  By restricting further to $Q_H$, we may identify the $\psi_T$-Kirillov model with a subspace of $C^\infty(Q_H)$.  By multiplying by some fixed right Haar measure on $Q_H$, we may identify $C^\infty(Q_H)$, in turn, with the space of smooth measures $\mathcal{M}(Q_H)$.  We denote by $\iota : \pi^\infty \hookrightarrow \mathcal{M}(Q_H)$ the resulting embedding.  By Jacquet's results, the maps
\begin{equation}\label{eq:mathc-right-mathc}
  \mathcal{M}_c(H) \rightarrow \mathcal{M}(Q_H)
\end{equation}
and
\begin{equation}\label{eq:crb6fqwwre}
  \mathcal{M}_c(P) \rightarrow \mathcal{M}(Q_H),
\end{equation}
both given by
\begin{equation*}
  f \mapsto \iota(\pi(f) \Theta_{\pi,\psi_T}),
\end{equation*}
are defined.

Recall from \S\ref{sec:class-bumps} and \S\ref{sec:spec-group-sett} that for Lie groups $G'$, we have defined certain subclasses $\mathfrak{C}(G', \dotsc) \subseteq \mathcal{M}_c(G')$.  The following lemma concerns the restriction to such subclasses of the maps \eqref{eq:mathc-right-mathc} and \eqref{eq:crb6fqwwre}.
\begin{lemma}\label{lem:fix-delta-in}
  Fix $\delta \in (0,1/2]$.  The above maps induce surjective class maps
  \begin{equation*}
    \mathfrak{C}(H,-\theta_H(\psi),T,\delta) \rightarrow \mathfrak{C}(Q_H,0,T,\delta),
  \end{equation*}
  \begin{equation*}
    \mathfrak{C}(P,-\theta_P(\psi),T,\delta) \rightarrow \mathfrak{C}(Q_H,0,T,\delta).
  \end{equation*}
  In particular, $\mathfrak{C}(H,-\theta_H(\psi),T,\delta)$ and $\mathfrak{C}(P,-\theta_P(\psi),T,\delta)$ have the same images under the above maps.
\end{lemma}
\begin{proof}
  We start by giving the proof for the first map, involving $H$.  The proof will make use of the fact, noted after \eqref{eq:crb6fpkcd0}, that the map
  \begin{equation*}
    j: N_H \times Q_H \rightarrow H,
    \qquad
    (n, g) \mapsto g^{-1} n
  \end{equation*}
  defines a local analytic isomorphism in some neighborhood of the identity.  An essentially identical argument works for $P$, using the analogous map $N \times Q_H \rightarrow P$.

  To simplify notation in what follows, we write simply $\theta$ for $\theta_H(\psi)$.

  The conclusion is independent of our normalization of Haar measure, provided that it is fixed, since any two fixed Haar measures are related by a scalar $\asymp 1$.  We may thus suppose given Haar measures on $H$ and $N_H$, and a right Haar measure on $Q_H$, that are related by the identity
  \begin{equation}\label{eq:int-_h-phih}
    \int_{H} \phi(h) \, d h
    =
    \int_{n \in N_H}
    \int_{g \in Q_H} \phi(g^{-1} n) \, d n \, d g,
  \end{equation}
  using here that $N_H$ is unimodular.

  Let us explicate the indicated maps.  Let $f \in \mathcal{M}_c(H)$.  We may write $f(h) = f_0(h) \, d h$ with $f_0 \in C_c^\infty(H)$.  By \eqref{eq:pif-theta_pi-h}, we have
  \begin{equation*}
    \pi(f) \Theta_{\pi,\psi_T}(g) =  \int_{n \in N_H} f_0(g^{-1} n ) \psi_T(n) \, d n.
  \end{equation*}
  Thus $\iota(\pi(f) \Theta_{\pi,\psi_T})$ is the measure on $Q_H$ assigning to a function $\phi \in C_c(Q_H)$ the value
  \begin{equation}\label{eq:int-_q_h-phi}
    \int_{Q_H} \phi \cdot \iota(\pi(f) \Theta_{\pi,\psi_T})
    =
    \int_{g \in Q_H}
    \int_{n \in N_H}
    f_0(g^{-1} n) \psi_T(n) \phi(g)
    \, d n
    \, d g.
  \end{equation}

  Consider the chart
  \begin{equation*}
    c : H \xdashrightarrow{} \Lie(N_H) \oplus  \Lie(Q_H) \cong \Lie(H),
  \end{equation*}
  defined near the identity, given by $j^{-1}$ composed with the logarithm.  Let $\psi_{T \theta, c} : H \xdashrightarrow{} \U(1)$ denote the function given, as in Definition \ref{defn:standard2:given-point-manif}, by $h \mapsto \psi_{T \theta}(c(h))$.  Thus for $(n,g) \in N_H \times Q_H$, we have $\psi_{T \theta, c}(g^{-1} n) = \psi_T(n)$.  We may rewrite the RHS of \eqref{eq:int-_q_h-phi} as
  \begin{equation}\label{eq:int_g-in-q_h}
    \int_{g \in Q_H} \phi(g) \int_{n \in N_H} \psi_{T \theta, c}(g^{-1} n) f_0(g^{-1} n) \, d n \, d g.
  \end{equation}
  Consider the following analytic map, defined near the origin and submersive there:
  \begin{equation*}
    \rho : H \xdashrightarrow{} Q_H
  \end{equation*}
  \begin{equation*}
    g^{-1} n \mapsto g \text{ for } (n,g) \in N_H \times Q_H.
  \end{equation*}
  Then, in view of \eqref{eq:int-_h-phih},  we may rewrite \eqref{eq:int_g-in-q_h} further as
  \begin{equation*}
    \int_{g \in H} \phi(\rho(g)) \psi_{T \theta , c}(g) f_0(g) \, d g.
  \end{equation*}
  Let $m : \mathcal{M}_c(H) \rightarrow \mathcal{M}_c(H)$ denote multiplication by $\psi_{T \theta, c}$.  The map \eqref{eq:mathc-right-mathc} may be written as the composition
  \begin{equation*}
    \mathcal{M}_c(H) \xrightarrow{m} \mathcal{M}_c(H) \xrightarrow{\rho_*} \mathcal{M}(Q_H).
  \end{equation*}
  By Lemma \ref{lem:modulation} and Proposition \ref{lem:assume-that-rho}, this composition induces class maps
  \begin{equation*}
    \mathfrak{C}(H,-\theta,T,\delta) \xrightarrow{m} \mathfrak{C}(H,0,T,\delta) \xrightarrow{\rho_*} \mathfrak{C}(Q_H,0,T,\delta),
  \end{equation*}
  with the first map bijective and the second map surjective.  The required conclusion follows.
\end{proof}

\begin{definition}\label{defn:we-denote-mathfr-frakW}
  We denote by \index{classes!bump functions in the Kirillov model, $\mathfrak{W}(\pi,\psi,T,\delta)$} $\mathfrak{W}(\pi,\psi,T,\delta)$ the subclass of $\mathcal{W}(\pi,\psi_T)$ consisting of all elements $W$ satisfying the following conditions, which, in view of Lemma \ref{lem:fix-delta-in}, are equivalent to one another.
  \begin{itemize}
  \item $W = \pi(\gamma) \Theta_{\pi,\psi_T}$ for some $\gamma \in \mathfrak{C}(H,-\theta_H(\psi),T,\delta)$.
  \item $W = \pi(\gamma) \Theta_{\pi,\psi_T}$ for some $\gamma \in \mathfrak{C}(P,-\theta_P(\psi),T,\delta)$.
  \item $\iota(W) \in \mathfrak{C}(Q_H, 0, T, \delta)$.
  \end{itemize}
\end{definition}
Informally, $\mathfrak{W}(\pi,\psi,T,\delta)$ consists of $W \in \mathcal{W}(\pi,\psi_T)$ for which $W|_{Q_H}$ is a smooth $L^1$-normalized bump function concentrated on $1 + \O(T^{\delta-1})$.

\begin{lemma}\label{lem:let-w-in-square-integral}
  Let $W \in \mathfrak{W}(\pi,\psi,T,\delta)$.  Then
  \begin{equation*}
    \int_{N_H \backslash H} |W|^2 \ll T^{(1- \delta) \dim(Q_H)}.
  \end{equation*}
\end{lemma}
\begin{proof}
  With suitable measure normalization, we have $\int_{N_H \backslash H} |W|^2 = \int_{Q_H} |W|^2$, and $W|_{Q_H} \in \mathfrak{C}(Q_H,0,T,\delta)$.  We can compute the integral in local coordinates.  The $L^1$-norm of $W|_{Q_H}$ is $\ll 1$, while the $L^\infty$-norm is $\ll T^{(1 - \delta) \dim(Q_H)}$.  The claim follows.
\end{proof}

\begin{lemma}\label{lem:there-exists-w-big-group-big}
  There exists $W \in \mathfrak{W}(\pi,\psi,T,\delta)$ so that $W|_{Q_H} \geq 0$, $\int_{Q_H} W = 1$.
\end{lemma}
\begin{proof}
  It is easy to see that there exists $\gamma \in \mathfrak{C}(Q_H,0,T,\delta)$ which is positive and has total mass one (take the rescaling of a fixed bump function, as in Example \ref{example:standard2:smooth-modulated-bump-in-frak-C}).  By the surjectivity assertion of Lemma \ref{lem:fix-delta-in}, there exists $W \in \mathfrak{W}(\pi,\psi,T,\delta)$ so that $\iota(W) = \gamma$.  The claim follows.
\end{proof}

\subsection{Behavior under $G$}\label{sec:behavior-under-g}
We will show that elements $W$ of the class $\mathfrak{W}(\pi,\psi,T,\delta)$ are approximately invariant under integral operators defined by test functions $\omega_0 \in C_c^\infty(G)$ with favorable properties.

\begin{theorem}\label{thm:big-group-whittaker-behavior-under-G}
  Fix $\eps > 0$.  Fix $\delta \in (0,1/2)$ with $1/2 - \delta$ sufficiently small in terms of $\eps$.  Let $T \ggg 1$.  Let $\pi$ be a generic irreducible unitary representation of $G$.  Assume that the rescaled infinitesimal character $T^{-1} \lambda_\pi$ belongs to a fixed compact set.

  There exists $\omega_0 \in C_c^\infty(G)$ with the following properties.  Set $\bar{G} = G/Z$.  Define $\omega := \omega_0 \ast \omega_0^\ast \in C_c^\infty(G)$ and $\omega^\sharp \in C_c^\infty(\bar{G})$, as previously, by
  \begin{equation*}
    \omega(g) := \int_{h \in G} \omega_0(g h^{-1}) \overline{\omega_0(h^{-1})} \, d h,
  \end{equation*}
  \begin{equation*}
    \omega^\sharp (g) := \int_{z \in Z}
    \pi|_Z (z) \omega(z g) \, d z,
  \end{equation*}
  where $\pi|_Z$ denotes the central character of $\pi$.  Then:
  \begin{enumerate}[(i)]
  \item\label{enumerate:cnpstqit4d} $\omega(g) \neq 0$ only if $g = 1 + o(1)$ and $\Ad^*(g) \tau = \tau + o(T^{-1/2})$, where $\tau = \tau(\pi,\psi,T)$ is given by \eqref{eq:tau-characterize-via-pi}.
  \item\label{enumerate:cnpstqivhi} $\|\omega^\sharp\|_{\infty} \ll T^{\dim(N)+\eps}$.
  \item\label{enumerate:cnpstqiwp2} If $\tau$ lies in some fixed compact subset of $\mathfrak{g}^\wedge_{\stab}$, then $\int_H |\omega^\sharp| \ll T^{\rank(H)/2 + \eps}$.
  \item\label{enumerate:cnpstqix13} For each $W \in \mathfrak{W}(\pi,\psi,T,\delta)$ and fixed $x \in \mathfrak{U}(G)$, we have
    \begin{equation*}
      \|\pi(x) (\pi(\omega_0)W - W)\| \ll T^{-\infty}.
    \end{equation*}
  \end{enumerate}
\end{theorem}
The proof is given in Part \ref{part:asympt-analys-kirill}.

We refer again to \S\ref{sec:elem-form} for an elementary reformulation of this result, without the language of nonstandard analysis.


\section{Analysis on $U \backslash G$: qualitative}\label{sec:some-analysis-basic}

This section supplies the qualitative local analysis underlying the construction of ``$f$'' sketched in \S\ref{sec:analyt-test-vect}.  We introduce certain functions $\mathcal{J}[\psi,\beta,\gamma]$ on $U \backslash G$, obtained from the open-cell distribution $A N \ni a n \mapsto \psi(n)$ by smoothing with respect to the two-sided action of $A \times G$, using Schwartz data $\beta$ on $A$ and $\gamma$ on $G$.  A key feature of this construction is its compatibility with intertwining operators.  The central result of this section is Proposition \ref{lem:sub-gln:each-psi-beta}, which shows that such functions are Schwartz and depend continuously on the smoothing data $(\beta,\gamma)$.  In verifying this, a key technical input is a uniform analytic-continuation estimate for Jacquet integrals.  We further study the Mellin components and Whittaker transforms attached to such functions.  In \S\ref{sec:cunt6n2c9p}, these qualitative statements are upgraded to quantitative estimates for a class of functions $\mathfrak{F}(U \backslash G,\psi,T)$ depending upon a large parameter $T$.  In \S\ref{sec:constr-test-vect-proofs}, such functions supply the local test vector data needed to prove Theorem~\ref{thm:main-local-results}.  Specifically, \S\ref{sec:some-analysis-basic} and \S\ref{sec:cunt6n2c9p} serve to establish assertions \eqref{item:lemma-there-exists-gamma-W-6} and \eqref{item:lemma-there-exists-gamma-W-7} of Lemma~\ref{lem:there-exists-begin} (via Propositions \ref{lem:all-f-in} and \ref{prop:there-exists-gamma}, respectively), which in turn reformulate parts of Theorem~\ref{thm:main-local-results}.

We note that the group denoted $G$ here plays the role of the smaller general linear group ``$H$'' in other parts of the paper.

\subsection{Preliminaries}
Let $G$ be a split connected reductive group over an archimedean local field $F$, thus $F = \mathbb{R}$ or $F = \mathbb{C}$.  In what follows, we view $G$ primarily as a Lie group over $\mathbb{R}$ (even if $F = \mathbb{C}$), so that the discussion of \S\ref{sec:modul-bump-funct} applies.  We expect that the discussion to follow applies more generally to any quasi-split real reductive group (cf.\ Lemma \ref{lem:assume-that-g} below).

We retain the general notation of \S\ref{sec:gener-prel}, thus $B = A N$ and $Q = AU$ are Borel subgroups of $G$, with common maximal torus $A$.  For computational convenience, we assume that these choices are opposite to one another, i.e.,
\begin{equation*}
  N = \bar{U},
\end{equation*}
so that the image of $A N$ in $U \backslash G$ is open.

Given a measure $\gamma$ on $G$ and $g \in G$, we write $g \ast \gamma$ and $\gamma \ast g$ for the convolutions of $\gamma$ with the $\delta$-mass at $g$.  We define in the same way $x \ast \gamma$ and $\gamma \ast x$ for $x \in \mathfrak{U}(G)$, regarding $x$ as a distribution supported at the identity element.

Recall (\S\ref{sec:schwartz-spaces}) that we have defined the spaces $\mathcal{S}(G)$ and $\mathcal{S}(U \backslash G)$ of Schwartz functions.  These spaces are related as follows:
\begin{lemma}\label{lem:cal-S-G-to-cal-S-U-G}
  For $f \in \mathcal{S}(G)$ and $g \in G$, the integral $f_U(g) := \int_{u \in U} f(u g) \, d u$ converges absolutely and defines an element $f_U$ of $\mathcal{S}(U \backslash G)$.  The resulting map $\mathcal{S}(G) \rightarrow \mathcal{S}(U \backslash G)$ is continuous.
\end{lemma}
\begin{proof}
  This follows readily the definitions and the fact that some negative power of the norm $\|.\|$ on $G$ is integrable over $U$ \cite[Cor 1.3] {MR1001613}.
\end{proof}

\subsection{Smoothened Whittaker distributions}

\begin{definition}
  Let $\psi$ be a nondegenerate unitary character of $N$.  Let $\beta \in \mathcal{S}(A)$.  We define a function \index{Jacquet integrals and distributions!$\mathcal{J}[\psi,\beta]$}
  \begin{equation*}
    \mathcal{J}[\psi,\beta] : U \backslash G \rightarrow \mathbb{C}
  \end{equation*}
  by requiring that $\mathcal{J}[\psi,\beta]$ be supported on the open cell represented by elements of the form $a n$ with $(a,n) \in A \times N$, and given on such elements by
  \begin{equation*}
    \mathcal{J}[\psi,\beta](a n) = \delta_U^{1/2}(a) \beta(a) \psi(n).
  \end{equation*}
\end{definition}
It is clear from construction that $\mathcal{J}[\psi,\beta]$ defines a bounded function on $U \backslash G$; more precisely,
\begin{equation}\label{eq:mathc-beta_-leq}
  \|\mathcal{J}[\psi,\beta]\|_{\infty} \leq \|\delta_U^{1/2} \beta \|_{\infty}.
\end{equation}
It is smooth under the left action of $A$.  For $a \in A$, we recall from \eqref{eq:normalized-left-translation} the normalized left translation action $L(a) f(g) := \delta_U^{-1/2}(a) f(a g)$, which induces an action of $\mathfrak{U}(A)$ that we denote in the same way.  For $x \in A$ or $x \in \mathfrak{U}(A)$, we have
\begin{equation}\label{eq:lx-mathcalj_psi-beta}
  L(x) \mathcal{J}[\psi,\beta] = \mathcal{J}[\psi,L(x) \beta].
\end{equation}

\begin{definition}
  Let $(\psi,\beta)$ be as above, and let $\gamma \in \mathcal{S}(G)$.  We define  \index{Jacquet integrals and distributions!$\mathcal{J}[\psi,\beta,\gamma]$}
  \begin{equation*}
    \mathcal{J}[\psi,\beta,\gamma] : = R(\gamma) \mathcal{J}[\psi,\beta] : U \backslash G \rightarrow \mathbb{C},
  \end{equation*}
  i.e.,
  \begin{equation}\label{eq:cuml2ory7z}
    \mathcal{J}[\psi,\beta,\gamma](h) :=
    \int_{g \in G}
    \gamma(g) \mathcal{J}[\psi,\beta](h g) \, d g.
  \end{equation}
\end{definition}
Since $\mathcal{J}[\psi,\beta]$ is bounded and $\gamma$ is integrable, such integrals converge absolutely.  For $\gamma_1, \gamma_2 \in \mathcal{S}(G)$, we have
\begin{equation}\label{eq:rgamm-mathc-beta}
  R(\gamma_1) \mathcal{J}[\psi,\beta,\gamma_2]
  = \mathcal{J}[\psi,\beta,\gamma_1 \ast \gamma_2].
\end{equation}
The functions $\mathcal{J}[\psi,\beta,\gamma]$ are smooth; more precisely, for $x \in \mathfrak{U}(G)$, we have
\begin{equation}\label{eq:rx-mathcald_psi-beta}
  R(x) \mathcal{J}[\psi,\beta,\gamma] = \mathcal{J}[\psi,\beta,x \ast  \gamma].
\end{equation}
We note the corresponding identity for the left action by $A$ on $U \backslash G$: for $x \in A$ or $x \in \mathfrak{U}(A)$, we have (by \eqref{eq:lx-mathcalj_psi-beta})
\begin{equation}\label{eq:lx-mathcald_psi-beta}
  L(x) \mathcal{J}[\psi,\beta,\gamma] = \mathcal{J}[\psi,L(x) \beta, \gamma].
\end{equation}
We note also that, in view of \eqref{eq:mathc-beta_-leq},
\begin{equation}\label{eq:mathc-beta-gamm}
  \|\mathcal{J}[\psi,\beta,\gamma]\|_{\infty} \leq
  \|\delta_U^{1/2} \beta \|_{\infty} \|\gamma \|_{1},
\end{equation}
where $\|\gamma \|_1$ denotes the $L^1$-norm (i.e., the total variation norm of the corresponding Schwartz measure).

The Mellin transform of $\mathcal{J}[\psi,\beta,\gamma]$ admits the following direct description.  It motivates the analytic-continuation estimate recorded next and will be used in the proof of Proposition \ref{lem:sub-gln:each-psi-beta}.
\begin{lemma}\label{lem:mellin-expansion-cal-J}
  Let $\psi$ be a nondegenerate unitary character of $N$ and $(\beta,\gamma) \in \mathcal{S}^e(A) \times \mathcal{S}(G)$.  Set $f := \mathcal{J}[\psi,\beta,\gamma]$.  For $h \in G$, set
  \begin{equation*}
    \alpha[h](g) := \int_{u \in U} \gamma(h^{-1} u g)\,d u.
  \end{equation*}
  For $\chi \in \mathfrak{X}(A)$, define
  \begin{equation*}
    \alpha[h,\chi](g)
    :=
    \int_{a \in A} \delta_U^{-1/2}(a)\chi(a)\alpha[h](a g)\,d a
  \end{equation*}
  whenever this integral converges, and write $W[h,\chi,\psi]$ as shorthand for $W[\alpha[h,\chi],\psi^{-1}]$, with $W[\,\cdot\,,\,\cdot\,]$ denoting the Jacquet integral of \S\ref{sec:whittaker-functions-1} (note that this shorthand suppresses the dependence of these functions upon $\gamma$).  Assume that Haar measures are normalized so that the multiplication map $A \times U \times N \rightarrow G$ is measure-preserving.  For any strictly dominant $\sigma \in \mathfrak{a}^*$, we have the convergent Mellin expansion
  \begin{equation}\label{eq:mellin-expansion-cal-J}
    f(a k)
    =
    \int_{(\sigma)}
    f[s](k)\,
    \delta_U^{1/2}(a) |a|^s
    \, d \mu_A(s),
  \end{equation}
  where $f[s] \in \mathcal{I}(s)$ is given by
  \begin{equation}\label{eq:mellin-component-cal-J}
    f[s](h) = \tilde{\beta}(|.|^s)\,W[h,|.|^s,\psi](1).
  \end{equation}
\end{lemma}
\begin{proof}
  Substituting $g \mapsto h^{-1} g$ in \eqref{eq:cuml2ory7z} gives
  \begin{equation*}
    f(h) = \int_{g \in G} \gamma(h^{-1} g)\mathcal{J}[\psi,\beta](g)\,d g.
  \end{equation*}
  With our normalization of Haar measures, we have
  \begin{equation*}
    \int_{G} \phi (g) \,d g
    =
    \int_N \int_A \int_U \phi (u a n ) \, d u \, \frac{d a}{ \delta_U(a)} \, d n.
  \end{equation*}
  Inserting the definition of $\mathcal{J}[\psi,\beta]$ gives
  \begin{equation*}
    f(h) =
    \int_{n \in N}
    \int_{a \in A}
    \int_{u \in U}
    \gamma(h^{-1} u a n)
    \beta(a) \psi(n)
    \, d u \, \frac{d a}{\delta_U^{1/2}(a)} \, d n.
  \end{equation*}
  With the notation introduced in the statement, this says
  \begin{equation}\label{eq:fh-int-a-n}
    f(h)
    =
    \int_{n \in N} \psi(n)
    \left(
      \int_{a \in A}
      \beta(a)\,\alpha[h](a n)\,\frac{d a}{\delta_U^{1/2}(a)}
    \right)\,d n.
  \end{equation}
  We insert the Mellin expansion $\beta(a)=\int_{(\sigma)} \tilde{\beta}(|.|^s)|a|^s\,d\mu_A(s)$, which converges for all $\sigma \in \mathfrak{a}^*$, and take $\sigma$ strictly dominant.  By the absolute convergence of the Jacquet integral in the strictly dominant case (\S\ref{sec:whittaker-functions-1}), we may rearrange the resulting double integral to obtain
  \begin{equation*}
    f(h)
    =
    \int_{(\sigma)} \tilde{\beta}(|.|^s)\,W[h,|.|^s,\psi](1)\,d\mu_A(s),
  \end{equation*}
  The assignment $h \mapsto \alpha[h, \chi]$, and hence also $h \mapsto W[h,\chi,\psi](1)$, transforms on the left under $A$ by $\delta_U^{1/2}\chi$.  Writing $h = a k$ with $(a,k) \in A \times K$ then gives \eqref{eq:mellin-expansion-cal-J} with $f[s]$ as in \eqref{eq:mellin-component-cal-J}.
\end{proof}

\subsection{Bounds for the Jacquet integral}\label{sec:analyt-cont-jacq}
The Mellin formula \eqref{eq:mellin-component-cal-J} expresses the Mellin components $f[s]$ of $\mathcal{J}[\psi,\beta,\gamma]$ in terms of Jacquet integrals.  Since such integrals are known to be analytic in the parameter $s$, we deduce the analytic continuation of $f[s]$ to all $s$.  We require the finer property that $f[s]$ decays rapidly in vertical strips.  In view of the rapid decay of $\tilde{\beta}(\lvert . \rvert^s)$, we require a uniform version of the analytic continuation of Jacquet integrals, with polynomial dependence on the imaginary part of the inducing character.  We record such a result here.

Since the notation used in this subsection is a bit different from that in the rest of the paper (where we typically take $G$ split, for simplicity), we recap all necessary notation.

Let $G$ be a quasi-split reductive group over $\mathbb{R}$.  We assume moreover, for technical reasons (see Lemma \ref{lem:assume-that-g}), that $G$ arises via restriction of scalars from a split group over $\mathbb{R}$ or $\mathbb{C}$.  In particular, the discussion to follow applies to the groups $\GL_n(\mathbb{R})$ and $\GL_n(\mathbb{C})$.

We assume given a pair of Borel subgroups $B = M N$ and $Q = M U$ of $G$, with maximal torus $M$ and respective unipotent radicals $N$ and $U$.  (We will use the latter for induced representations, the former for Whittaker functions.)

Let $X(M)$ denote the group of algebraic characters of $M$.  Write
\begin{equation*}
  \mathfrak{a}^* := X(M) \otimes_{\mathbb{Z}} \mathbb{R} \hookrightarrow \mathfrak{a}_{\mathbb{C}}^* := X(M) \otimes_{\mathbb{Z}} \mathbb{C}.
\end{equation*}
(We may identify $\mathfrak{a}^*$ with the dual of $\Lie(A)$, where $A \leq M$ denotes the maximal $\mathbb{R}$-split subtorus.)  Let $\mathfrak{X}(M)$ denote the group of complex characters $\chi : M \rightarrow \mathbb{C}^\times$.  To each $s \in \mathfrak{a}_{\mathbb{C}}^*$, we may associate a character $|.|^s \in \mathfrak{X}(M)$, denoted $m \mapsto |m|^{s}$.  For each $\chi \in \mathfrak{X}(M)$, we may define $\Re(\chi) \in \mathfrak{a}^*$ by writing $|\chi(m)| = |m|^{\Re(\chi)}$.  For $\chi \in \mathfrak{X}(M)$ and $s \in \mathfrak{a}^*_{\mathbb{C}}$, we abbreviate $\chi_s := \chi |.|^s \in \mathfrak{X}(M)$.  We denote by $\delta_U : M \rightarrow \mathbb{R}^\times_+$ the modulus character.

Let $K \leq G$ be a maximal compact subgroup such that $K \cap M$ is a maximal compact subgroup of $M$.  Let $W$ denote the Weyl group for $(G,A)$.  Let $\psi$ be a nondegenerate unitary character of $N$.  We choose Haar measures on $N$ and $G$.  We choose a $W$-invariant norm $|.|$ on $\mathfrak{a}_{\mathbb{C}}^*$.

Let $\Delta \subseteq \mathfrak{a}^*$ denote the set of simple roots attached to $U$.  In this section, ``dominant'' always refers to the sense defined by $U$ (rather than $N$).

For each $\chi \in \mathfrak{X}(M)$, we denote by $\mathcal{I}(\chi)$ the corresponding principal series representation of $G$, consisting of all smooth functions $f : U \backslash G \rightarrow \mathbb{C}$ satisfying $f(m g) = \delta_U^{1/2}(m) \chi(m) f(g)$ for all $(m,g) \in M \times G$.

Let $w_0 \in G$ be a representative for the Weyl group element such that $w_0 N = \bar{U} w_0$, where $\bar{U}$ denotes the maximal unipotent subgroup opposite to $U$.

Let $\chi \in \mathfrak{X}(M)$ and $f \in \mathcal{I}(\chi)$.  We define $W[f,\psi] : G \rightarrow \mathbb{C}$, first for strictly dominant $\chi$ by the absolutely convergent Jacquet integral
\begin{equation}\label{eq:wf-psig-:whittaker-analytic}
  W[f,\psi](g) := \int_{n \in N } f (w_0 n g) \psi^{-1}(n) \, d n,
\end{equation}
then in general by analytic continuation along flat sections (see \cite{MR563369}, \cite[\S15.4]{MR1170566}).  It defines an equivariant map
\begin{equation*}
  W[\cdot, \psi] : \mathcal{I}(\chi) \rightarrow C^\infty(N \backslash G, \psi).
\end{equation*}
We have seen some references in which it is noted that such maps are ``locally uniform'' provided that $\chi$ is constrained to a given compact subset of $\mathfrak{X}(M)$  (see, e.g., \cite[\S3.5]{MR2533003}).  We require the following slightly finer uniformity condition, in which only the real part of $\chi$ is constrained.
\begin{proposition}\label{prop:jacquet-integral-uniformity-question:let-eta-be}
  Let $\chi \in \mathfrak{X}(M)$.  Let $\mathcal{D}$ be a compact subset of $\mathfrak{a}^*$.  There exists $C_0 \geq 0$, $\ell \in \mathbb{Z}_{\geq 0}$ and a seminorm $\mu$ on $C^\infty(K)$ so that for all $s \in \mathfrak{a}_{\mathbb{C}}^*$ with $\Re(s) \in \mathcal{D}$, we have for all $f \in \mathcal{I}(\chi_s)$ the estimate
  \begin{equation*}
    |W[f,\psi](1)|
    \leq
    C_0
    (1 + |s|)^{\ell}
    \mu(f).
  \end{equation*}
\end{proposition}
We expect that such an estimate should follow (in sharper forms) by a careful quantification of the proofs given in \cite{MR563369} or \cite[\S15.4]{MR1170566}.  For the sake of completeness, we record an alternative proof below, roughly in the spirit of \cite[Proof of Lemma 5.6]{MR1800349}, that instead assumes the analytic continuation of the Jacquet integral and deduces the desired uniformity using standard bounds for intertwining operators and the functional equation.  The proof will follow a series of lemmas.

For $w \in W$ and $\chi \in \mathfrak{X}(M)$, we define the meromorphic families of intertwining operators
\begin{equation*}
  M_w = M_w[\chi] : \mathcal{I}(\chi) \rightarrow \mathcal{I}(w \chi)
\end{equation*}
and their normalized variants
\begin{equation*}
  \mathcal{F}_w = \mathcal{F}_{w,\psi}[\chi] : \mathcal{I}(\chi) \rightarrow \mathcal{I}(w \chi)
\end{equation*}
exactly as in \S\ref{sec:local-intertw-oper-four}.  We retain the notational conventions of those sections.  We denote by $c(w,\chi)$ the proportionality constant (``local coefficient'') relating the two:
\begin{equation*}
  \mathcal{F}_w[\chi] = c(w,\chi) M_w[\chi].
\end{equation*}

By a \emph{root polynomial} $p$ on $\mathfrak{a}_{\mathbb{C}}^*$, we mean a polynomial of the form
\begin{equation*}
  p(s) = \prod_{(\alpha,c) \in \Sigma}
  ( \alpha^\vee(s) - c  )
\end{equation*}
for some finite set $\Sigma$ consisting of pairs of roots $\alpha$ and complex numbers $c$.

\begin{lemma}\label{lem:jacquet-integral-uniformity-question:let-chi-in}
  Let $w \in W$.  Let $\chi \in \mathfrak{X}(M)$.  Let $\mathcal{D}$ be a compact subset of $\mathfrak{a}^*$.  There is a root polynomial $p$ on $\mathfrak{a}_{\mathbb{C}}^*$ and scalars $C_0 \geq 0$, $\ell \in \mathbb{Z}_{\geq 0}$ so that for all $s \in \mathfrak{a}_{\mathbb{C}}^*$ with $\Re(s) \in \mathcal{D}$, we have
  \begin{equation*}
    |p(s) c(w,\chi_s)| \leq C_0 (1 + |s|)^{\ell}.
  \end{equation*}
\end{lemma}
\begin{proof}
  We may express $c(w,\chi)$ as a product of local Tate $\gamma$-factors (see \cite[Thm 1]{MR816396}).  It follows from this expression that $c(w,\chi_s)$ has at most finitely many poles with $\Re(s) \in \mathcal{D}$.  We construct $p$ to kill each such pole.  The required upper bound then follows from Stirling's formula (as, for instance, in \cite[\S2.2]{2021arXiv210112106B}).
\end{proof}
\begin{lemma}\label{lem:jacquet-integral-uniformity-question:let-chi-in-1}
  Let $w \in W$.  Let $\chi \in \mathfrak{X}(M)$.  Let $d \in \mathbb{R}$.  There exists $\ell \in \mathbb{Z}_{\geq 0}$, a seminorm $\mu$ on $C^\infty(K)$ and a root polynomial $p$ on $\mathfrak{a}_{\mathbb{C}}^*$ with the following property.  Let $s \in \mathfrak{a}_{\mathbb{C}}^*$ such that $\alpha^\vee(\Re(s)) \geq d$ for all $\alpha \in \Delta$ and $p(s) \neq 0$.  Then $s$ is not a pole of $M_w : \mathcal{I}(\chi_s) \rightarrow \mathcal{I}(w \chi_s)$, and for all $f \in \mathcal{I}(\chi_s)$,
  \begin{equation*}
    \|p(s) M_w f\|_{L^\infty(K)} \leq  (1 + |s|)^{\ell} \mu(f).
  \end{equation*}
\end{lemma}
\begin{proof}
  This is a special case of a result of Lapid \cite[Cor 2]{MR2402686}, whose proof refers in turns to Wallach \cite[\S10.1.11]{MR1170566}.
\end{proof}

\begin{lemma}\label{lem:standard:let-alpha-F-w-alpha-holom}
  Let $\alpha \in \Delta$.  Let $w_\alpha$ denote the corresponding simple reflection, thus $w_\alpha^2 = 1$.  The meromorphic family of operators
  \begin{equation*}
    \mathcal{F}_{w_\alpha} : \mathcal{I}(w_\alpha \chi ) \rightarrow \mathcal{I}(\chi)
  \end{equation*}
  is holomorphic on $\{\chi :  \alpha^\vee(\Re(\chi))   \geq 0\}$, or indeed, on $\{\chi :  \alpha^\vee(\Re(\chi)) > -1\}$.
\end{lemma}
\begin{proof}
  We can reduce to the case $G = \SL_2$ and appeal to the theory of one-dimensional local Tate integrals.  For instance, in the case $F = \mathbb{R}$, the conclusion of the lemma is noted by Shahidi \cite{MR563369}, just before the statement of \cite[Lemma 2.3.2]{MR563369}, to follow from \cite[Prop 2.2.1]{MR563369}.  The general case may be treated in the same way.
\end{proof}

\begin{lemma}\label{lem:jacquet-integral-uniformity-question:let-w-in}
  Let $w \in W$.  Then
  \begin{equation*}
    \mathcal{F}_w : \mathcal{I}(w^{-1} \chi) \rightarrow \mathcal{I}(\chi)
  \end{equation*}
  is holomorphic for $\chi$ with dominant real part.
\end{lemma}
\begin{proof}
  We factor $w$ as a minimal product $w_{\alpha_1} \dotsb w_{\alpha_{\ell}}$ of simple reflections, as in \cite[Thm 2.1.1]{MR610479}; by \cite[\S3.2]{MR610479}, $\mathcal{F}_w$ then factors as a composition of operators $\mathcal{F}_{w_{\alpha_j}}$ to which Lemma \ref{lem:standard:let-alpha-F-w-alpha-holom} applies.
\end{proof}

\begin{lemma}\label{lem:assume-that-g}
  Assume that $G$ arises, possibly via restriction of scalars, from a split reductive group over $\mathbb{R}$ or $\mathbb{C}$.  Let $\chi \in \mathfrak{X}(M)$.  Let $f \in \mathcal{I}(\chi)$ be $K$-finite.  For $s \in \mathfrak{a}_{\mathbb{C}}^*$, let $f_s \in \mathcal{I}(\chi_s)$ be the element with $f_s|_K = f$.  Then the analytic function $s \mapsto W[f_s,\psi](1)$ is of finite order $\leq 1$.
\end{lemma}
\begin{proof}
  This is a special case of a result of McKee \cite[Thm 5.2]{MR3043499}.
\end{proof}
\begin{remark}
  For the argument that follows, we require the weaker conclusion that the restriction of $s \mapsto W[f_s,\psi](1)$ to $\{s : \Re(s) \in \mathcal{D} \}$ is of finite order for each compact subset $\mathcal{D}$ of $\mathfrak{a}^*$.  This conclusion may be available for more general groups $G$ (and without the assumption of $K$-finiteness) via the technique of \cite[Proof of Lemma 5.6]{MR1800349}, but we were unable to verify all details.
\end{remark}

\begin{proof}[Proof of Proposition \ref{prop:jacquet-integral-uniformity-question:let-eta-be}]
  For notational convenience, we will assume that $N = \bar{U}$ is opposite to $U$, so that $w_0$ may be taken to be the identity element and the definition \ref{eq:wf-psig-:whittaker-analytic} simplifies to
  \begin{equation}\label{eq:wf-psi1-=whittaker-integral-simplified}
    W[f, \psi](1) = \int_{n \in N } f (n) \psi^{-1}(n) \, d n.
  \end{equation}

  Let us agree to write $f_s \in \mathcal{I}(\chi_s)$ for the section denoted ``$f$'' in the statement of Proposition \ref{prop:jacquet-integral-uniformity-question:let-eta-be}, and $f$ for the unique section in $\mathcal{I}(\chi)$ whose restriction to $K$ is $f_s$.  It then suffices to show the following:

  \emph{Claim.}  For each finite subset $\mathcal{X}$ of  $\mathfrak{X}(M)$ and compact $\mathcal{D} \subseteq \mathfrak{a}^*$ as above, there exists $C_0 \geq 0$, $\ell \in \mathbb{Z}_{\geq 0}$ and a seminorm $\mu$ on $C^\infty(K)$ so that for all $\chi \in \mathcal{X}$ and $f \in \mathcal{I}(\chi)$,
  \begin{equation}\label{eq:wfchi-psi1-leq}
    \Re(s) \in \mathcal{D} \, \implies \, |W[f_s, \psi](1)| \leq C_0 (1 + |s|)^{\ell} \mu(f).
  \end{equation}

  In verifying this, it will be convenient to assume that $\Re(\chi) = 0$ for each $\chi \in \mathcal{X}$ (otherwise, replace $\chi$ with $\chi |.|^{-\Re(\chi)}$ and enlarge $\mathcal{D}$ suitably).  We may assume also -- enlarging $\mathcal{X}$ as necessary -- that $w \chi \in \mathcal{X}$ for all $(w,\chi) \in W \times \mathcal{X}$.

  By the continuity of both sides of the desired inequality, we may assume that $f$ is $K$-finite (but, of course, the quantities $C_0$, $\ell$ and $\mu$ must be independent of the $K$-type of $f$).

  We may find a compact subset $\mathcal{D}_0$ of $\mathfrak{a}^*$ with the following properties:
  \begin{enumerate}[(i)]
  \item Every element of $\mathcal{D}_0$ is strictly dominant.
  \item $\mathcal{D}$ is contained in the convex hull of $\cup_{w \in W} w \mathcal{D}_0$.
  \end{enumerate}
  For $w \in W$, we have $w(\chi_s) = (w \chi)_{w s}$.  By hypothesis, we have $w \chi \in \mathcal{X}$.  Since the norm $|.|$ on $\mathfrak{a}_{\mathbb{C}}^*$ is $W$-invariant, we have $|s| = |w s|$ for all $(s,w) \in \mathfrak{a}_{\mathbb{C}}^* \times W$.  In view of Lemma \ref{lem:assume-that-g} and our assumption that $f$ is $K$-finite, the Phragmen--Lindel\"{o}f principle is applicable.  The required estimate \eqref{eq:wfchi-psi1-leq} will thus follow if we can show, with quantifiers as above, that
  \begin{equation}\label{eq:wfchi-psi1-leq-2}
    (\Re(s),w) \in \mathcal{D}_0 \times W \, \implies \,
    |W[f_{w s}, \psi](1)| \leq C_0 (1 + |s|)^{\ell} \mu(f).
  \end{equation}

  We consider first the case $w = 1$.  For $\sigma \in \mathfrak{a}^*$, let $\phi_{\sigma } \in \mathcal{I}(|.|^{\sigma})$ denote the spherical vector taking the value $1$ at the identity.  In particular, $\phi_\sigma \geq 0$.  If $\sigma$ is strictly dominant, then the integral $\int_{N} \phi_\sigma$ converges and is bounded by a continuous function of $\sigma$ (see \cite[\S10.1.2]{MR1170566}).  In particular,
  \begin{equation*}
    \sup_{\sigma \in \mathcal{D}_0} \int_N \phi_\sigma(n) \, d n < \infty.
  \end{equation*}
  On the other hand, for $f \in \mathcal{I}(\chi)$, $s \in \mathfrak{a}_\mathbb{C}^*$ and $n \in N$, we have
  \begin{equation*}
    |f_s(n)| \leq \|f\|_{L^\infty(K)} \phi_{\Re(s)}(n).
  \end{equation*}
  Thus, recalling \eqref{eq:wf-psi1-=whittaker-integral-simplified}, we obtain
  \begin{equation*}
    \Re(s) \in \mathcal{D}_0
    \,
    \implies
    \,
    |W[f_s,\psi](1)|
    \leq
    C_0
    \|f\|_{L^\infty(K)},
  \end{equation*}
  where $C_0$ depends only upon $\mathcal{D}_0$.
  This completes the proof of the ``$w = 1$ case'' of the implication \eqref{eq:wfchi-psi1-leq-2}.

  For the general case, we appeal to the functional equation \eqref{eq:wmathcalf_w-fchi-psi} (i.e., the defining property of $\mathcal{F}_{w^{-1}}$), which gives
  \begin{equation*}
    W[f_{w s}, \psi](1)
    =
    W[\mathcal{F}_{w^{-1}} f_{w s}, \psi](1).
  \end{equation*}
  We have $\mathcal{F}_{w^{-1}} f_{w s} \in \mathcal{I}((w^{-1}\chi)_s)$, so for $\Re(s) \in \mathcal{D}_0$, the $w = 1$ case gives
  \begin{equation*}
    |W[\mathcal{F}_{w^{-1}} f_{w s}, \psi](1)| \leq C_0 \|\mathcal{F}_{w^{-1}} f_{w s}\|_{L^\infty(K)}.
  \end{equation*}
  Our task thereby reduces to verifying that
  \begin{equation}\label{eq:res-w-in}
    (\Re(s), w) \in \mathcal{D}_0 \times W \implies \|\mathcal{F}_{w^{-1}} f_{w s}\|_{L^\infty(K)} \leq C_0 (1 + |s|)^{\ell} \mu(f).
  \end{equation}
  To that end, we write $\mathcal{F}_{w^{-1}} = \mathcal{F}_{w^{-1}}[\chi_{w s}] : \mathcal{I}(\chi_{w s}) \rightarrow \mathcal{I}((w^{-1} \chi)_{s})$ as
  \begin{equation*}
    \mathcal{F}_{w^{-1}}[\chi_{w s}] = c(w^{-1},\chi_{w s}) M_{w^{-1}}[\chi_{w s}],
  \end{equation*}
  as above.  Combining Lemmas \ref{lem:jacquet-integral-uniformity-question:let-chi-in} and \ref{lem:jacquet-integral-uniformity-question:let-chi-in-1} over the finite set $W \times \mathcal{X}$, we see that there is
  \begin{itemize}
  \item a root polynomial $p$ on $\mathfrak{a}_{\mathbb{C}}^*$,
  \item $C_0 \geq 0$,
  \item $\ell \in \mathbb{Z}_{\geq 0}$ and
  \item a seminorm $\mu$ on $C^\infty(K)$
  \end{itemize}
  so that
  \begin{equation}\label{eq:res-w-in-1}
    (\Re(s), w) \in \mathcal{D}_0 \times W \implies \|p(s) \mathcal{F}_{w^{-1}} f_{w s}\|_{L^\infty(K)} \leq C_0 (1 + |s|)^{\ell} \mu(f).
  \end{equation}
  This is nearly the required implication, except for the factor $p(s)$.  To address that, we note that the real part of $w^{-1}(\chi_{w s})$ is the real part of $s$, which is strictly dominant, so by Lemma \ref{lem:jacquet-integral-uniformity-question:let-w-in}, the operator $\mathcal{F}_{w^{-1}}[\chi_{w s}]$ is holomorphic in the relevant range.  We may thus pass from \eqref{eq:res-w-in-1} to \eqref{eq:res-w-in} via Cauchy's integral formula.
\end{proof}

\begin{remark}\label{rmk:standard2:jacquet-analytic-refine}
  It would be natural to further quantify Proposition \ref{prop:jacquet-integral-uniformity-question:let-eta-be} by explicating the dependence upon $\chi$.  This would require explicating the analogous dependence for \cite[Cor 2]{MR2402686}, hence in turn for \cite[\S10.1.11]{MR1170566}.
\end{remark}

\subsection{Schwartz property}\label{sec:qual-jacquet-schwartz}
We now return to the general setting of \S\ref{sec:some-analysis-basic}.  We have already noted that the functions $\mathcal{J}[\psi,\beta,\gamma]$ are smooth and bounded.  That they decay rapidly is the content of the following proposition, which is the main technical output of this section.  Its continuity assertion gives bounds for fixed Schwartz seminorms in terms of seminorms of the smoothing data.  After the rescaling introduced in \S\ref{sec:analysis-U-G-quantitative-setup}, these become weak, polynomial-in-$T$ bounds for elements of $\mathfrak{F}(U \backslash G,\psi,T)$ (see first Lemma \ref{lem:standard2:mathfr-backsl-g}, then Lemma \ref{lem:scratch-research:each-fixed-seminorm}).  Such bounds provide preliminary control at infinity; the sharper $T$-uniform estimates needed in the main quantitative arguments are proved later in \S\ref{sec:analysis-U-G-quantitative-setup}.
\begin{proposition}\label{lem:sub-gln:each-psi-beta}
  For each nondegenerate unitary character $\psi$ of $N$ and $(\beta,\gamma) \in \mathcal{S}^e(A) \times \mathcal{S}(G)$, we have $\mathcal{J}[\psi,\beta,\gamma] \in \mathcal{S}^e(U \backslash G)$.  For each $\psi$, the map
  \begin{equation*}
    \mathcal{S}^e(A) \times \mathcal{S}(G) \rightarrow \mathcal{S}^e(U \backslash G)
  \end{equation*}
  \begin{equation}\label{eq:beta-gamma-mapsto}
    (\beta,\gamma) \mapsto \mathcal{J}[\psi,\beta,\gamma]
  \end{equation}
  is continuous.
\end{proposition}
\begin{proof}
  It is clear from the equivariance property \eqref{eq:lx-mathcald_psi-beta} that $\mathcal{J}[\psi,\beta,\gamma]$ is left-invariant under the maximal compact subgroup of $A$, hence lies in $\mathcal{S}^e(U \backslash G)$ provided that it lies in $\mathcal{S}(U \backslash G)$.  To establish membership in the latter and verify the required continuity, we must show that for each $(y,m) \in \mathfrak{U}(G) \times \mathbb{Z}_{\geq 0}$, there are
  continuous seminorms $\mu_1$ on $\mathcal{S}^e(A)$ and $\mu_2$ on $\mathcal{S}(G)$ so that for all $(\beta,\gamma) \in \mathcal{S}^e(A) \times \mathcal{S}(G)$, we have
  \begin{equation*}
    \sup_{h \in U \backslash G} \|h\|_{U \backslash G}^m \left\lvert R(y) \mathcal{J}[\psi,\beta,\gamma](h) \right\rvert
    \leq
    \mu_1(\beta) \mu_2(\gamma).
  \end{equation*}
  Using the equivariance property \eqref{eq:rx-mathcald_psi-beta}, we may reduce to the case $y=1$.  Let us abbreviate
  \begin{equation*}
    f := \mathcal{J}[\psi,\beta,\gamma].
  \end{equation*}
  By Paley--Wiener theory for $\mathcal{S}^e(U \backslash G)$, it will suffice to show that there is a holomorphic family $f[s] \in \mathcal{I}(s)$, indexed by $s \in \mathfrak{a}_{\mathbb{C}}^*$, with the following properties.
  \begin{itemize}
  \item The family $\{f[s]\}$ is of rapid vertical decay (\S\ref{sec:mellin-paley-affine-quotients}), so that its Mellin inverse may be defined.  That Mellin inverse equals $f$.
  \item For each compact $\mathcal{D} \subseteq \mathfrak{a}^*$ and $m \in \mathbb{Z}_{\geq 0}$, there are seminorms $(\mu_1,\mu_2)$ as above so that
    \begin{equation}\label{eq:sup_-substack-s-1}
      \sup_{
        \substack{
          s \in \mathfrak{a}_{\mathbb{C}}^* :   \\
          \Re(s) \in \mathcal{D}
        }
      }
      (1 + |s|)^{m} \sup_{k \in K} |f[s](k)|  \leq \mu_1(\beta) \mu_2(\gamma).
    \end{equation}
  \end{itemize}
  By Lemma \ref{lem:mellin-expansion-cal-J}, we may take $f[s]$ as in \eqref{eq:mellin-component-cal-J}; that lemma also supplies a Mellin expansion of $f$ (initially for strictly dominant $\sigma$).  It remains to verify the rapid vertical decay and the uniform estimates \eqref{eq:sup_-substack-s-1}.

  We may understand $W[k,|.|^s, \psi](1)$ as the application of the Jacquet integral to the holomorphic family of vectors $\alpha[k,|.|^s]$.  Since the Jacquet integral varies analytically (see \S\ref{sec:whittaker-functions-1} and \S\ref{sec:analyt-cont-jacq}), we see that $\{f[s]\}$ defines a holomorphic family.  We claim that this family is of rapid vertical decay.  In checking this claim, we may once again invoke the equivariance property \eqref{eq:rx-mathcald_psi-beta} to reduce to verifying the following qualitative form of \eqref{eq:sup_-substack-s-1}: for each compact subset $\mathcal{D}$ of $\mathfrak{a}^*$ and $m \in \mathbb{Z}_{\geq 0}$,
  \begin{equation}\label{eq:sup_-substack-s}
    \sup_{
      \substack{
        s \in \mathfrak{a}_{\mathbb{C}}^* :  \\
        \Re(s) \in \mathcal{D}
      }
    }
    (1 + |s|)^{m}
    \sup_{k \in K}
    \left\lvert
      f[s](k)
    \right\rvert < \infty.
  \end{equation}
  We verify this below, and in doing so, see that the LHS is in fact bounded by continuous seminorms of $\beta$ and $\gamma$, giving the required continuity per the discussion preceding \eqref{eq:sup_-substack-s-1}.

  To that end, Proposition \ref{prop:jacquet-integral-uniformity-question:let-eta-be} gives $d \in \mathbb{Z}_{\geq 0}$ and, after absorbing a scalar, a seminorm $\mu_2$ on $C^\infty(K)$, depending only upon $\mathcal{D}$, so that for all $k \in K$ and $s \in \mathfrak{a}_{\mathbb{C}}^*$ with $\Re(s) \in \mathcal{D}$,
  \begin{equation*}
    |W[k,|.|^s,\psi](1)| \leq (1 + |s|)^d \mu_2(\alpha[k,|.|^s]).
  \end{equation*}
  We choose a continuous seminorm $\mu_1$ on $\mathcal{S}(A)$ so that
  \begin{equation*}
    \sup_{s : \Re(s) \in \mathcal{D}} (1 + |s|)^{m+d} \left\lvert \tilde{\beta}(|.|^s) \right\rvert \leq \mu_1(\beta).
  \end{equation*}
  Therefore \eqref{eq:sup_-substack-s} is bounded by
  \begin{equation*}
    \mu_1(\beta) \sup_{ s : \Re(s) \in \mathcal{D} } \sup_{k \in K} \mu_2(\alpha[k,|.|^s]).
  \end{equation*}
  On the other hand:
  \begin{itemize}
  \item For each compact subset $\mathcal{D}$ of $\mathfrak{a}^*$, the map
    \begin{equation*}
      \mathcal{S}(U \backslash G)
      \ni \eta
      \mapsto \sup_{
        \substack{
          \chi \in \mathfrak{X}(A):  \\
          \Re(\chi) \in \mathcal{D}
        }
      }
      \mu_2(\eta[\chi])
    \end{equation*}
    defines a continuous seminorm on $\mathcal{S}(U \backslash G)$ (by \S\ref{sec:mellin-paley-affine-quotients}).
  \item
    The family of maps
    \begin{equation*}
      \mathcal{S}(G) \ni \gamma \mapsto \alpha[k] \in \mathcal{S}(U \backslash G)
      \quad (k \in K)
    \end{equation*}
    is equicontinuous (by Lemma \ref{lem:cal-S-G-to-cal-S-U-G}).
  \end{itemize}
  Combining these, we deduce that there is a continuous seminorm $\mu_3$ on $\mathcal{S}(G)$ so that
  \begin{equation*}
    \sup_{ s : \Re(s) \in \mathcal{D} } \sup_{k \in K} \mu_2(\alpha[k,|.|^s]) \leq \mu_3(\gamma).
  \end{equation*}
  This completes the proof of \eqref{eq:sup_-substack-s} in the required quantitative form.
\end{proof}
\begin{remark}
  The conclusion of the proposition is closely related to arguments employed in some \emph{proofs} of the analytic continuation of Jacquet integrals (cf.\ \cite[Thm 2.1]{MR563369}), and could likely be deduced by quantifying those arguments.  We have instead argued the other way, deducing the conclusion from (a quantitative form of) the analytic continuation of Jacquet integrals.
\end{remark}
\begin{remark}
  In \cite{MR2058615}, Jacquet gave some results for $\GL_n$ that are similar in spirit.
\end{remark}
\begin{remark}\label{rmk:standard2:we-expect-it-extend-S}
  We expect it to be true and provable that the same assignment defines a continuous map $\mathcal{S}(A) \times \mathcal{S}(G) \rightarrow \mathcal{S}(U \backslash G)$.  We do not require such an extension, whose verification would require additional effort (see Remark \ref{rmk:standard2:jacquet-analytic-refine}).
\end{remark}

\subsection{Bounds for Mellin components}\label{sec:qual-mellin-component-bounds}
Having established the Schwartz property for the functions constructed above, we now record a couple bounds for the Mellin components and Whittaker functions attached to general Schwartz functions $f \in \mathcal{S}(U \backslash G)$.  Recall from \S\ref{sec:mellin-paley-affine-quotients} that the Mellin components $f[\chi]$, for $\chi \in \mathfrak{X}(A)$, are defined by
\begin{equation*}
  f[\chi](g) = \int_{a \in A} \delta_U^{1/2}(a)\chi(a)f(a^{-1} g)\,\frac{d a}{\delta_U(a)}.
\end{equation*}
We record two bounds for such components.  The first gives their rapid decay in vertical strips.  The second, obtained by combining the first with the bounds for the Jacquet integral established in \S\ref{sec:analyt-cont-jacq}, gives a corresponding bound for the attached Whittaker functional.  In \S\ref{sec:analysis-U-G-quantitative-setup}, they will be specialized to elements of the class $\mathfrak{F}(U \backslash G,\psi,T)$ (see Lemmas \ref{lem:standard:let-chi-in} and \ref{lem:scratch-research:let-f-in} and Proposition \ref{lem:all-f-in}) en route to the Whittaker estimates proved there.

\begin{lemma}\label{lem:standard:each-comp-mathc}
  For each compact $\mathcal{D} \subseteq \mathfrak{a}^*$, $m \in \mathbb{Z}_{\geq 0}$ and seminorm $\kappa_0$ on $C^\infty(K)$, there is a seminorm $\kappa$ on $\mathcal{S}(U \backslash G)$ so that for all $f \in \mathcal{S}(U \backslash G)$ and all $\chi \in \mathfrak{X}^e(A)$ with $\Re(\chi) \in \mathcal{D}$, we have
  \begin{equation*}
    \kappa_0(f[\chi]) \leq \kappa(f) C(\chi)^{-m}.
  \end{equation*}
\end{lemma}
\begin{proof}
  We may find $x \in \mathfrak{U}(A)$ so that $L(x) f[\chi] = c f[\chi]$ with $|c| \geq C(\chi)^{m}$.  Since $L(x) : \mathcal{S}(U \backslash G) \rightarrow \mathcal{S}(U \backslash G)$ is continuous, we may reduce to the case $m = 0$.

  We may assume that $\kappa_0(v) = \|R(y) v\|_{L^\infty(K)}$ for some $y \in \mathfrak{U}(K)$, so that
  \begin{equation*}
    \kappa_0(f[\chi]) =
    \sup_{k \in K}
    \left\lvert
      \int_{a \in A}
      \delta_U^{1/2}(a) \chi(a)
      R(y) f(a^{-1} k)
      \, \frac{d a}{\delta_U(a)}
    \right\rvert.
  \end{equation*}

  We may find $\ell \geq 0$ so that
  \begin{equation*}
    C_0 := \sup_{
      \substack{
        \chi \in \mathfrak{X}^e(A) :  \\
        \Re(\chi) \in \mathcal{D}
      }
    }
    \int_{a \in A} \|a\|^{-\ell} \delta_U^{1/2}(a) |\chi(a)|  \, \frac{d a}{\delta_U(a)} < \infty.
  \end{equation*}
  We may find $C_1 \geq 0$ so that $\|a k\|_{U \backslash G}^{\ell} \leq C_1 \|a\|^{\ell}$ for all $(a,k) \in A \times K$.  The conclusion then holds with
  \begin{equation*}
    \kappa(f) := C_0 C_1 \sup_{g \in U \backslash G} \|g\|_{U \backslash G}^{\ell} |R(y) f(g)|.
  \end{equation*}
\end{proof}

\begin{lemma}\label{lem:standard:each-comp-mathc-1}
  For each compact $\mathcal{D} \subseteq \mathfrak{a}^*$ and $m \in \mathbb{Z}_{\geq 0}$, there is a seminorm $\kappa$ on $\mathcal{S}(U \backslash G)$ so that for all $f \in \mathcal{S}(U \backslash G)$ and all unramified $\chi \in \mathfrak{X}(A)$ with $\Re(\chi) \in \mathcal{D}$, we have
  \begin{equation*}
    |W[f[\chi], \psi](1)| \leq \kappa(f) C(\chi)^{-m}.
  \end{equation*}
\end{lemma}
\begin{proof}
  We combine Proposition \ref{prop:jacquet-integral-uniformity-question:let-eta-be} and Lemma \ref{lem:standard:each-comp-mathc}.
\end{proof}

\subsection{Mellin components and intertwining operators}\label{sec:qual-jacquet-mellin-intertwining}
We now specialize the preceding Mellin-component formalism to the functions $\mathcal{J}[\psi,\beta,\gamma]$.  The first goal is to compute their Mellin components by testing against principal-series vectors, in a form compatible with Jacquet integrals.  The main consequence is that the assignment $(\beta,\gamma) \mapsto \mathcal{J}[\psi,\beta,\gamma]$ intertwines the normalized Fourier transforms on $\mathcal{S}(U \backslash G)$ with the natural Weyl action on the $A$-variable.  We also record a pairing formula for Whittaker values of Mellin components, used later in the lower-bound argument.

\begin{lemma}\label{lem:scratch-research:let-chi-in-1-v-J}
  Let $\chi \in \mathfrak{X}(A)$ with strictly dominant real part, and let $v \in \mathcal{I}(\chi)$.  Then for all $\beta \in \mathcal{S}(A)$, the integral
  \begin{equation*}
    \int_{U \backslash G} v \cdot \overline{\mathcal{J}[\psi, \beta ]}
  \end{equation*}
  converges absolutely and is continuous in $v$ and $\beta$.  If the Haar measure on $U \backslash G$ is normalized so that the multiplication map $N \times A \rightarrow U \backslash G$ is measure-preserving, then the above integral is equal to
  \begin{equation*}
    \overline{\tilde{\beta}(\bar{\chi}^{-1})}
    W[v, \psi](1).
  \end{equation*}
\end{lemma}
\begin{proof}
  For $(a,n) \in A \times N$, we have $v(a n) = \delta_U^{1/2} \chi(a) v(n)$, $\mathcal{J}[\psi,\beta](a n) = \delta_U^{1/2}(a) \beta(a) \psi(n)$ and $d(a n) = \delta_U^{-1}(a) \, d a \, d n$, so
  \begin{equation*}
    \int_{U \backslash G} v \cdot \overline{\mathcal{J}[\psi, \beta ]}
    =
    \int_{a \in A}
    \int_{u \in N}
    \chi(a) \overline{\beta}(a)
    v(n)
    \psi^{-1}(n) \, d n \, d a.
  \end{equation*}
  This integral factors as a product of absolutely convergent integrals over $A$ and $N$, which evaluate as indicated.  These absolutely convergent integrals are moreover bounded by continuous seminorms, whence the continuity.
\end{proof}

For the following discussion, we take $(\beta,\gamma) \in \mathcal{S}^e(A) \times \mathcal{S}(G)$.

We have seen that $\mathcal{J}[\psi,\beta,\gamma]$ lies in $\mathcal{S}(U \backslash G)$, so we may define its Mellin components for $\chi \in \mathfrak{X}(A)$ (which, since $\beta \in \mathcal{S}^e(A)$, vanish for $\chi$ outside $\mathfrak{X}^e(A)$).  For typographical reasons, we often denote these by
\begin{equation*}
  \mathcal{J}[\psi,\beta,\gamma,\chi] := \mathcal{J}[\psi,\beta,\gamma][\chi] \in \mathcal{I}(\chi).
\end{equation*}

Write, as before, $\gamma^*(g) = \overline{\gamma(g^{-1})}$.  The following technical lemma justifies an exchange of integrals used below in Lemmas \ref{lem:standard:chi-in-mathfrakxa}, \ref{lem:sub-gln:each-psi-beta-1} and \ref{lem:scratch-research:let-beta_1-beta_2}.
\begin{lemma}\label{lem:scratch-research:let-chi-in}
  Let $\chi \in \mathfrak{X}(A)$ and $v \in \mathcal{I}(\chi)$.   We have

  \begin{equation}\label{eq:int-_n-v}
    \overline{\tilde{\beta}(\bar{\chi}^{-1})}
    W[R(\gamma^*) v,\psi](g)
    =
    \int_{N} v \cdot \overline{\mathcal{J}[\psi, \beta, \gamma \ast g^{-1}, \bar{\chi}^{-1}]}.
  \end{equation}
\end{lemma}
\begin{proof}
  We observe first that the RHS of the stated identity makes sense.  Indeed, we have seen that $\mathcal{J}[\psi,\beta,\gamma \ast g^{-1}]$ lies in $\mathcal{S}(U \backslash G)$, so its Mellin components define a holomorphic family of vectors.  Up to normalization of Haar measure, the integral on the RHS of \eqref{eq:int-_n-v} over $N$ may be replaced with the same integral over $K$.  In particular, that integral converges absolutely for all $\chi$.

  We may reduce to the case $g = 1$, as follows.  We have
  \begin{equation*}
    W[ R(\gamma^* ) v, \psi](g) = W[R(g) R(\gamma^* ) v, \psi](1),
  \end{equation*}
  $R(g ) R(\gamma^*) = R(g \ast \gamma^*)$ and $g \ast \gamma^* = (\gamma \ast g^{-1})^*$.  Thus if the identity \eqref{eq:int-_n-v} holds for $g=1$ and all $\gamma$, then it holds for all $g$.  Our task is thus to verify that
  \begin{equation}\label{eq:overl-1-wgamma}
    \overline{\tilde{\beta}(\bar{\chi}^{-1})}
    W[R(\gamma^*) v,\psi](1)
    =
    \int_{N} v \cdot \overline{\mathcal{J}[\psi, \beta, \gamma, \bar{\chi}^{-1}]}.
  \end{equation}

  Since $\gamma$ and its derivatives are integrable, the assignment $v \mapsto R(\gamma^*) v$ defines a continuous map from the representation $\mathcal{I}(\chi)$, equipped with the $L^2(K)$-topology, to the same representation equipped with the usual smooth topology.  It follows that both sides of \eqref{eq:overl-1-wgamma}  depend continuously upon $v$ with respect to the $L^2(K)$-topology.  By a limiting argument, we may thus reduce to the case that the restriction of $v$ to $N$ is compactly-supported.

  We define in that case a holomorphic family $v[s] \in \mathcal{I}(\chi |.|^s)$ by requiring that $v[s]|_{N} = v|_N$, so that $v[0] = v$.  The LHS of \eqref{eq:overl-1-wgamma} is, by definition, the value at $s = 0$ of the analytic continuation of
  \begin{equation}\label{eq:overl-1.-bars}
    \overline{\tilde{\beta}(\bar{\chi}^{-1}|.|^{-\bar{s}})}
    W[R(\gamma^*) v[s], \psi](1),
  \end{equation}
  defined initially for $\Re(s)$ sufficiently dominant by an absolutely convergent integral.  For such $s$, we may apply Lemma \ref{lem:scratch-research:let-chi-in-1-v-J} to rewrite the above as the absolutely convergent integral
  \begin{equation*}
    \int_{U \backslash G}
    R(\gamma^*) v[s] \cdot
    \overline{\mathcal{J}[\psi,\beta]},
  \end{equation*}
  This last integral may be expanded as the double integral (absolutely convergent, again by Lemma \ref{lem:scratch-research:let-chi-in-1-v-J})
  \begin{equation*}
    \int_{g \in U \backslash G}
    \int_{h \in G} \gamma^*(h)
    v[s](g h) \overline{\mathcal{J}[\psi,\beta](g)} \, d h \, d g.
  \end{equation*}
  The substitutions $g \mapsto g h^{-1}$  and $h \mapsto h^{-1}$ allow us to further rewrite the above as
  \begin{equation*}
    \int_{g \in U \backslash G}
    \int_{h \in G} \gamma^*(h^{-1})
    v[s](g) \overline{\mathcal{J}[\psi,\beta](g h)} \, d h \, d g,
  \end{equation*}
  which, recalling the definition $R(\gamma) \mathcal{J}[\psi,\beta] = \mathcal{J}[\psi,\beta,\gamma]$, may in turn be rewritten
  \begin{equation}\label{eq:int-_g-in-gu}
    \int_{g \in U \backslash G}
    v[s](g) \overline{\mathcal{J}[\psi,\beta,\gamma](g)} \, d g.
  \end{equation}
  Since (by Proposition \ref{lem:sub-gln:each-psi-beta}) we have $\mathcal{J}[\psi,\beta,\gamma] \in \mathcal{S}(U \backslash G)$, this last integral converges absolutely for all $s$, hence defines an analytic function of $s$ given by the same formula.  Specializing the equality between \eqref{eq:overl-1.-bars} and \eqref{eq:int-_g-in-gu} to the case $s = 0$ gives the desired identity \eqref{eq:overl-1-wgamma}.
\end{proof}

For the following result, we denote by $\check{\gamma}(g) := \gamma(g^{-1})$ the ``real adjoint'' of $\gamma$.
\begin{lemma}\label{lem:standard:chi-in-mathfrakxa}
  For $\chi \in \mathfrak{X}(A)$ and $v \in \mathcal{I}(\chi^{-1})$, we have
  \begin{equation}\label{eq:mathc-beta-gamm-1}
    \int_{N}
    v \cdot
    \mathcal{J}[\psi,\beta,\gamma \ast g^{-1}][\chi]
    =
    \tilde{\beta}(\chi)
    W[R(\check{\gamma}) v, \psi^{-1}](g).
  \end{equation}
\end{lemma}
\begin{proof}
  We apply complex conjugation to Lemma \ref{lem:scratch-research:let-chi-in}.
\end{proof}

Recall from \S\ref{sec:fourier-transforms} the unitary operators $\mathcal{F}_{w,\psi}$ on $\mathcal{S}(U \backslash G)$ attached to Weyl group elements $w \in W$.
\begin{lemma}\label{lem:sub-gln:each-psi-beta-1}
  For $w \in W$, we have
  \begin{equation*}
    \mathcal{F}_{w,\psi^{-1}} \mathcal{J}[\psi,\beta,\gamma] = \mathcal{J}[{\psi,{}^w \beta,\gamma}],
  \end{equation*}
  where ${}^w \beta(a) := \beta(w^{-1} a)$.
\end{lemma}
\begin{proof}
  For the proof, we abbreviate $\mathcal{F}_w := \mathcal{F}_{w,\psi^{-1}}$.  It is enough to check that for each $\chi \in [A]^\wedge$, we have
  \begin{equation*}
    (\mathcal{F}_w \mathcal{J}[\psi,\beta,\gamma])[ w \chi] = \mathcal{J}[{\psi,{}^w \beta, \gamma}][w \chi].
  \end{equation*}
  By definition,
  \begin{equation*}
    (\mathcal{F}_w \mathcal{J}[\psi,\beta,\gamma])[ w \chi] = \mathcal{F}_w \mathcal{J}[\psi,\beta,\gamma][\chi]
    :=
    \mathcal{F}_w (\mathcal{J}[\psi,\beta,\gamma][\chi]).
  \end{equation*}
  It is enough to check that for each $v \in \mathcal{I}(w \chi^{-1})$, we have
  \begin{equation*}
    (\mathcal{F}_w \mathcal{J}[\psi,\beta,\gamma][\chi], v)
    =
    (\mathcal{J}[{\psi,{}^w \beta, \gamma}][w \chi], v),
  \end{equation*}
  where $(,)$ denotes the invariant pairing given by integrating the product over $N$.  (We recall from \S\ref{sec:induc-repr} that this coincides up to scaling with the pairing given by integration over $K$.)  By \eqref{eq:mathc-1-mathc} and \eqref{eq:f_1chi-f_2chi-1} followed by Lemma \ref{lem:standard:chi-in-mathfrakxa} (specialized to $g=1$), we have
  \begin{equation*}
    (\mathcal{F}_w \mathcal{J}[\psi,\beta,\gamma][\chi], v) = (\mathcal{J}[\psi,\beta,\gamma][\chi], \mathcal{F}_{w^{-1}} v)
    = \tilde{\beta}(\chi) W[R(\check{\gamma}) \mathcal{F}_{w^{-1}} v, \psi^{-1}](1).
  \end{equation*}
  By another application of Lemma \ref{lem:standard:chi-in-mathfrakxa}, we have
  \begin{equation*}
    (\mathcal{J}[\psi, {}^{w} \beta, \gamma][w \chi], v)
    = {}^{w} \tilde{\beta}(w \chi)
    W[R(\check{\gamma}) v, \psi^{-1}](1).
  \end{equation*}
  In view of the identity $\tilde{\beta}(\chi) = {}^w \tilde{\beta}(w \chi)$, we thereby reduce to checking that
  \begin{equation*}
    W[R(\check{\gamma}) \mathcal{F}_{w^{-1}} v, \psi^{-1}](1)
    =
    W[R(\check{\gamma}) v, \psi^{-1}](1).
  \end{equation*}
  This follows from the equivariance property $R(\check{\gamma}) \mathcal{F}_{w^{-1}} v = \mathcal{F}_{w^{-1}} R(\check{\gamma}) v$ of $\mathcal{F}_{w^{-1}}$ and the defining property \eqref{eq:wmathcalf_w-fchi-psi} of $\mathcal{F}_{w^{-1}}$.
\end{proof}

The following technical lemma is used later in the proof of Proposition \ref{prop:there-exists-gamma}.  It expresses the Whittaker value $W[f[\chi], \psi](g)$ as a pairing given by integration over $N$.  In the application, we will take $\gamma = \gamma_0^* \ast \gamma_0$ and $\chi$ trivial, so that at $g = 1$, the pairing is a nonnegative integral, a feature that enables us to bound it from below.
\begin{lemma}\label{lem:scratch-research:let-beta_1-beta_2}
  Let $\beta_1, \beta_2 \in \mathcal{S}^e(A)$, $\gamma_1, \gamma_2 \in \mathcal{S}(G)$.  Set $\gamma := \gamma_2^* \ast \gamma_1$ and
  \begin{equation*}
    f := \mathcal{J}[\psi, \beta_1, \gamma] \in \mathcal{S}^e(U \backslash G).
  \end{equation*}
  Then for each $\chi \in \mathfrak{X}(A)$ and $g \in G$,
  \begin{equation*}
    \overline{\tilde{\beta}_2(\bar{\chi}^{-1})}
    W[f[\chi], \psi](g)
    =
    \int_{N}
    \mathcal{J}[\psi, \beta_1, \gamma_1, \chi]
    \overline{\mathcal{J}[\psi, \beta_2, \gamma_2 \ast g^{-1}, \bar{\chi}^{-1}]}.
  \end{equation*}
\end{lemma}
\begin{proof}
  Unwinding the definitions and applying \eqref{eq:rgamm-mathc-beta}, we see that
  \begin{equation*}
    f[\chi] = R(\gamma_2^*) \mathcal{J}[\psi,\beta_1,\gamma_1][\chi].
  \end{equation*}
  By Lemma \ref{lem:scratch-research:let-chi-in} applied to $v := \mathcal{J}[\psi,\beta_1,\gamma_1][\chi]$,
  we then have
  \begin{equation*}
    \overline{\tilde{\beta}_2(\bar{\chi}^{-1})}
    W[f[\chi], \psi](g)
    =
    \int_N v \cdot \overline{\mathcal{J}[\psi,\beta_2, \gamma_2 \ast g^{-1}, \bar{\chi}^{-1}]},
  \end{equation*}
  as required.
\end{proof}

\section{Analysis on $U \backslash G$: quantitative}\label{sec:cunt6n2c9p}

The goal of this section is to construct the local functions $f$ anticipated in the sketch of \S\ref{sec:analyt-test-vect}.  We package them in Definition \ref{defn:standard2:we-denote-mathfrakfu} as a class $\mathfrak{F}(U \backslash G,\psi,T)$ depending upon a large parameter $T$.  Elements of this class are obtained from the qualitative functions of \S\ref{sec:some-analysis-basic} by replacing $\psi$ with the rescaled character $\psi_T$, choosing the $A$-cutoff $\beta$ from a fixed bounded class, and smoothing on the right using a modulated bump measure $\gamma \in \mathfrak{C}(G,-\theta_G(\psi),T,\tfrac{1}{2})$.  The class is designed to have three useful features: compatibility with normalized intertwining operators, rapid decay of Mellin components, and microlocalization inherited from $\gamma$.  After developing these properties, we prove the two quantitative outputs needed later: upper and lower bounds for Whittaker functions of Mellin components (Proposition \ref{lem:all-f-in} and Proposition \ref{prop:there-exists-gamma}).  As noted earlier, these results are applied in \S\ref{sec:constr-test-vect-proofs} (with the role of the present section's $G$ played by $H$) to establish assertions \eqref{item:lemma-there-exists-gamma-W-6} and \eqref{item:lemma-there-exists-gamma-W-7} of Lemma~\ref{lem:there-exists-begin}, which feeds into the proof of Theorem~\ref{thm:main-local-results}.

\subsection{Setup}\label{sec:analysis-U-G-quantitative-setup}
We retain the setting of \S\ref{sec:some-analysis-basic} (but not the more general setting of \S\ref{sec:analyt-cont-jacq}); in particular, $G$ is a split reductive group over an archimedean local field.  We assume, of course, that the group $G$ is fixed, and continue to assume that $N$ is opposite to $U$.

We fix a nondegenerate unitary character $\psi$ of $N$.  Its pullback $\exp^* \psi$ to the Lie algebra defines an element of the Pontryagin dual $\Lie(N)^\wedge$.  By analogy to \S\ref{sec:parameter-theta}, we denote by
\begin{equation*}
  \theta_G(\psi) \in \Lie(G)^\wedge
\end{equation*}
the element whose restriction to $\Lie(N)$ is $\exp^*\psi$ and whose restriction to $\Lie(Q)$ is trivial.  Thus, with notation as in \S\ref{sec:pontry-dual-vect}, we have $\psi(\exp(x)) = \exp (\langle x, \theta_G(\psi)  \rangle)$ for $x \in \Lie(N)$, while $\langle x, \theta_G(\psi) \rangle = 0$ for $x \in \Lie(Q)$.

Let $T$ be a positive parameter with $T \ggg 1$.  We define the class
\begin{equation*}
  \mathfrak{C}(G,-\theta_G(\psi),T,\tfrac{1}{2})
\end{equation*}
as in \S\ref{sec:spec-group-sett}.  We recall (Proposition \ref{lem:standard:subcl-mathfr-theta}) that this class is closed under taking adjoints and convolution.

\subsection{The class $\mathfrak{B}(A)$}\label{sec:class-mathfrakba}

\begin{definition}
  We denote by $\mathfrak{B}(A)$ \index{classes!$\mathfrak{B}(A)$} the class of functions $\beta \in \mathcal{S}^e(A)$ that lie in some fixed bounded subset of $\mathcal{S}(A)$.  We denote by $\mathfrak{B}(A)^W \leq \mathfrak{B}(A)$ the subclass consisting of $W$-invariant elements.
\end{definition}

It is clear that if $x \in \mathfrak{U}(A)$ is fixed and $\beta \in \mathfrak{B}(A)$, then $L(x) \beta \in \mathfrak{B}(A)$.

\begin{lemma}\label{lem:standard:each-beta-in}
  Each $\beta \in \mathfrak{B}(A)$ may be decomposed as $\beta = \beta^\sharp + \beta^\flat$, where $\beta^\sharp, \beta^\flat \in \mathfrak{B}(A)$ have the following properties:
  \begin{enumerate}[(i)]
  \item $\beta^\sharp(a) \neq 0$ only if $\|a\| \ll T^{o(1)}$.
  \item $\beta^\flat \in T^{-\infty} \mathfrak{B}(A)$.
  \end{enumerate}
  The same conclusion holds with $\mathfrak{B}(A)$ replaced by $\mathfrak{B}(A)^W$.
\end{lemma}
\begin{proof}
  We may and shall assume that the norm $\|.\|$ on $A$ is defined so as to be $W$-invariant.

  Fix $\phi \in C_c^\infty(\mathbb{R})$ taking the value $1$ in a neighborhood of zero.  For each $t \geq 1$, define $\beta_t \in C_c^\infty(A)$ by $\beta_t(x) := \phi(\|x\|/t) \beta(x)$.  Then the family $\{\beta_t\}_{t \geq 1}$ lies in a fixed bounded subset of $\mathcal{S}(A)$ (and of $\mathcal{S}(A)^W$, provided that $\beta \in \mathcal{S}(A)^W$), and so $\beta_t$ lies in $\mathfrak{B}(A)$ (resp.\ $\mathfrak{B}(A)^W$) for each $t \geq 1$.

  We specialize to $t = t(\eps) := T^{\eps}$ with $\eps > 0$.  We claim that for each fixed $\eps > 0$ and fixed seminorm $\kappa$ for $\mathcal{S}(A)$, we have
  \begin{equation}\label{eq:cump8aw0t2}
    \kappa(\beta - \beta_{t(\eps)}) \ll T^{- \infty}.
  \end{equation}
  Indeed, we have $(\beta - \beta_{t(\eps)})(x) = \Phi_\eps(x) \beta(x)$ with $\Phi_\eps(x) := (1 - \phi(\|x\| / T^\eps))$.  Define the fixed quantity $c := \inf_{x \in \supp(\phi)} \lVert x \rVert$.  The function $\Phi_\eps$, together with each of its fixed derivatives, has magnitude $\ll 1$ and is supported on elements $x$ with $\|x\| \geq c T^{\eps}$.  Since each fixed derivative $\beta '$ of $\beta$ (of any order) satisfies $\|x\|^m |\beta '(x)| \ll 1$ for each fixed $m \geq 0$, the claim follows.

  The estimate \eqref{eq:cump8aw0t2} implies that, if we fix a countable collection $\{\kappa_j\}_{j \geq 1}$ of defining seminorms for $\mathcal{S}(A)$, then for each fixed $\eps > 0$, we have $\kappa_j(\beta - \beta_{t(\eps)}) \leq T^{-1/\eps}$ for all $j \leq 1/\eps$.  By the overspill principle, we may thus find $\eps > 0$ with $\eps \lll 1$ so that $\kappa_j(\beta - \beta_{t(\eps)}) \ll T^{-\infty}$ for each fixed $j$.  The same then holds for each fixed seminorm $\kappa$.  With this value of $\eps$, we take $\beta^\sharp := \beta_{t(\eps)}$ and $\beta^\flat := \beta - \beta_{t(\eps)}$.  Then $\beta^\flat \in T^{-\infty} \mathfrak{B}(A)$, while $\beta^\sharp$ has the required support properties.
\end{proof}

\subsection{The class $\mathfrak{F}(U \backslash G, T, \psi)$}

As in \S\ref{sec:class-bump-functions}, we define the rescaled nondegenerate unitary character $\psi_T$ of $N$ by requiring that for $x \in \Lie(N)$, we have $\psi_T(\exp(x)) = \psi(\exp(T x))$.

\begin{definition}\label{defn:standard2:we-denote-mathfrakfu}
  We denote by $\mathfrak{F}(U \backslash G, \psi, T)$ the class \index{classes!$\mathfrak{F}(U \backslash G, \psi, T)$} of functions $U \backslash G \rightarrow \mathbb{C}$ of the form
  \begin{equation*}
    f =
    \mathcal{J}_T[\psi,\beta,\gamma]
    :=
    L(T^{\rho_U^\vee })
    \mathcal{J}[\psi_T,\beta,\gamma]
    =
    L(T^{\rho_U^\vee })
    R(\gamma)
    \mathcal{J}[\psi_T,\beta]
  \end{equation*}
  for some $\beta \in \mathfrak{B}(A)$ and $\gamma \in \mathfrak{C}(G,-\theta_G(\psi),T,\tfrac{1}{2})$.
\end{definition}

\begin{example}\label{example:cump8dyq7s}
  We abbreviate $\theta := \theta_G(\psi)$.  Fix $\phi \in C_c^\infty(\mathfrak{g})$.  Let $\gamma_0$ denote the smooth measure on $\mathfrak{g}$ given by
  \begin{equation*}
    \gamma_0(x) = T^{\dim(G)/2} \phi(T^{1/2} x) \psi_{- T \theta}(x) \, d x.
  \end{equation*}
  The measure $\gamma_0$ is supported near the origin, hence identifies via pushforward under the exponential map with a smooth measure on $G$, characterized as follows: for any locally integrable function $\Phi$ on $G$,
  \begin{equation}\label{eq:int-_g-gamma_0}
    \int_{G} \gamma_0 \cdot \Phi
    = \int_{x \in \mathfrak{g} }
    \phi(x)
    \Phi (\exp (\tfrac{x}{T^{1/2} }))
    \exp (- \langle T \theta, \tfrac{x}{ T^{1/2} } \rangle)
    \, d x.
  \end{equation}
  It is easy to see, as in Example \ref{example:standard2:smooth-modulated-bump-in-frak-C}, that $\gamma_0$ lies in the class $\mathfrak{C}(G,-\theta,T,\tfrac{1}{2})$.  With this choice, the corresponding function $f = \mathcal{J}_T[\psi, \beta, \gamma_0]$ lies in $\mathfrak{F}(U \backslash G, \psi, T)$ for any $\beta \in \mathfrak{B}(A)$.  A variant of this construction (with $\gamma_0$ replaced by $\gamma_0^* \ast \gamma_0$) will be used below in the proof of the key Whittaker function lower bound (Proposition~\ref{prop:there-exists-gamma}).
\end{example}

\begin{remark}\label{rmk:we-discuss-informal-cal-F-U-G-T}
  We discuss the informal content of the above definition.
  \begin{enumerate}[(i)]
  \item
    As we will see below (see, e.g,. Proposition \ref{prop:standard:each-f-in}), elements of $\mathfrak{F}(U \backslash G, \psi, T)$ concentrate near the image of $T^{-\rho^\vee}$ in $U \backslash G$.  Informally speaking, they are uniformly smooth under left translation by $A$, while under right translation by $G$, they are approximate eigenfunctions of group elements of size $1 + \O(T^{-1/2})$, with eigenvalue described by $\theta_G(\psi)$.
  \item The uniform smoothness under left translation by $A$ implies that the Mellin components $f[s] \in \mathcal{I}(s)$ are small unless $s$ is small.  Since $f$ concentrates near $T^{-\rho^\vee}$, the components $f[s]$ will have magnitude roughly $T^{- \langle \rho^\vee, s \rangle}$ and satisfy the same equivariance condition under small elements of $G$ as $f$.  The analogous condition for the $p$-adic analogue of our discussion (see~\cite{2025arXiv2503.12310}) turns out to describe a one-dimensional subspace of $\mathcal{I}(s)$ consisting of \emph{Howe vectors} (see \cite[\S5]{MR2016587}, \cite[\S7.1]{MR2192818}, \cite[\S5]{MR3683106}).  One can thus understand the Mellin components $f[s]$ as analogues, over $\mathbb{R}$, of Howe vectors.
  \item The reader might wonder why we construct the class $\mathfrak{F}(U \backslash G, \psi, T)$ in terms of the distributions $\mathcal{J}[\psi,\beta]$, rather than in some more direct way.  The reason is that by doing so, we trivialize the problem of computing the images of elements of $\mathfrak{F}(U \backslash G, \psi, T)$ under intertwining operators (see Lemma \ref{lem:each-f-in}), which will be very convenient for our global study of the Eisenstein series associated to such elements (Part \ref{part:local-l2-growth}).
  \item
    The preceding properties enter the two Whittaker estimates below in different ways.  For the upper bound in Proposition \ref{lem:all-f-in}, after passing to the essential part \(f^\sharp\), we use only support and size information for its Mellin component \(f[\chi]\): its restriction to \(N\) is supported on \(1+\O(T^{-1/2+o(1)})\) and has the expected pointwise size.  For the lower bound in Proposition \ref{prop:there-exists-gamma}, the oscillation is essential: the modulated bump \(\gamma_0\) is chosen so that, near \(n=1\), its phase cancels the oscillation of \(\psi_T((n\exp(x/T^{1/2}))_N)\), making the relevant value of the trivial Mellin component roughly \(\psi_T(n)\int \phi\).  Taking \(\gamma=\gamma_0^*\ast\gamma_0\) then turns the Whittaker value at \(g=1\) into a nonnegative integral.
  \item
    We record an approximate reference picture for what elements $f$ of $\mathfrak{F}(U \backslash G, \psi, T)$ look like in the special case $G = \SL_2(\mathbb{R})$.    Taking for $U$ the lower-triangular unipotent subgroup, we can identify $U \backslash G$ with $\mathbb{R}^2 - \{0\}$ via the map taking a matrix to its top row.  Under this identification, it turns out that
    \begin{equation*}
      f(x,y) \approx
      T^{-1/4}
      1_{x \asymp T^{1/2}}
      1_{y \ll 1}
      \psi(T y/x)
    \end{equation*}
    for some fixed unitary character $\psi$ of $\mathbb{R}$.  Equivalently, in polar coordinates,
    \begin{equation*}
      f(r \cos \alpha, r \sin \alpha) \approx
      T^{-1/4} 1_{r \asymp T^{1/2}}
      1 _{\alpha \ll T^{-1/2}} \psi(T \alpha).
    \end{equation*}
  \end{enumerate}
\end{remark}

\begin{lemma}\label{lem:standard2:mathfr-backsl-g}
  We have $\mathfrak{F}(U \backslash G, \psi, T) \subseteq \mathcal{S}(U \backslash G)$.
\end{lemma}
\begin{proof}
  This is a special case of Proposition \ref{lem:sub-gln:each-psi-beta}.
\end{proof}

\begin{lemma}\label{lem:cal-J-T-vs-cal-J}
  We have
  \begin{equation*}
    \mathcal{J}_T[\psi,\beta,\gamma] =  \delta_U^{-1/2}(T^{\rho_U^\vee}) R(T^{-\rho_N^\vee}) \mathcal{J}[\psi,\beta, T^{\rho_N^\vee} \ast \gamma \ast T^{- \rho_N^\vee }].
  \end{equation*}
\end{lemma}
\begin{proof}
  Set $f := \mathcal{J}_T[\psi,\beta,\gamma] = R(\gamma) L(T^{\rho_U^\vee }) \mathcal{J}[\psi_T,\beta]$.  We have
  \begin{equation*}
    L(T^{\rho_U^\vee}) \mathcal{J}[\psi_T,\beta] = \delta_U^{-1/2}(T^{\rho_U^\vee}) R(T^{-\rho_N^\vee}) \mathcal{J}[\psi,\beta].
  \end{equation*}
  Indeed, evaluating at $a n$ for $(a,n) \in A \times N$ gives for the LHS
  \begin{align*}
    L(T^{\rho_U^\vee}) \mathcal{J}[\psi_T,\beta](a n) &=
                                                        \delta_U^{-1/2}(T^{\rho_U^\vee})  \mathcal{J}[\psi_T,\beta](T^{\rho_U^\vee} a n) \\
                                                      &=
                                                        \delta_U^{1/2}(a) \beta(T^{\rho_U^\vee} a) \psi_T(n),
  \end{align*}
  and for the RHS
  \begin{align*}
    \delta_U^{-1/2}(T^{\rho_U^\vee})  R(T^{-\rho_N^\vee}) \mathcal{J}[\psi,\beta](a n)
    &=
      \delta_U^{-1/2}(T^{\rho_U^\vee})  \mathcal{J}[{\psi,\beta}](a n T^{-\rho_N^\vee}) \\
    &=
      \delta_U^{-1/2}(T^{\rho_U^\vee})   \mathcal{J}[{\psi,\beta}](T^{\rho_U^\vee} a (T^{\rho_N^\vee} n T^{-\rho_N^\vee})) \\
    &=
      \delta_U^{1/2}(a) \beta(T^{\rho_U^\vee} a) \psi_T(n),
  \end{align*}
  as claimed.  (We have used here that $N$ is opposite to $U$, so that $\rho_U^\vee = - \rho_N^\vee$.)  Thus
  \begin{align*}
    f &=
        R(\gamma)     L(T^{\rho_U^\vee}) \mathcal{J}[\psi_T,\beta] \\
      &= \delta_U^{-1/2}(T^{\rho_U^\vee}) R(\gamma) R(T^{-\rho_N^\vee}) \mathcal{J}[\psi,\beta] \\
      &= \delta_U^{-1/2}(T^{\rho_U^\vee}) R(T^{-\rho_N^\vee}) R(T^{\rho_N^\vee} \ast \gamma \ast T^{- \rho_N^\vee }) \mathcal{J}[\psi,\beta] \\
      &= \delta_U^{-1/2}(T^{\rho_U^\vee}) R(T^{-\rho_N^\vee}) \mathcal{J}[\psi,\beta,  T^{\rho_N^\vee} \ast \gamma \ast T^{- \rho_N^\vee }].
  \end{align*}
\end{proof}

\begin{lemma}\label{lem:scratch-research:each-fixed-seminorm}
  For each fixed seminorm $\kappa$ on $\mathcal{S}(U \backslash G)$ and each $f \in \mathfrak{F}(U \backslash G, \psi, T)$, we have $\kappa(f) \ll T^{\O(1)}$.
\end{lemma}
\begin{proof}
  Write $f = \mathcal{J}_T[\psi,\beta,\gamma]$.  By Lemma \ref{lem:cal-J-T-vs-cal-J}, we have
  \begin{equation*}
    f = \delta_U^{-1/2}(T^{\rho_U^\vee}) \mathcal{J}[\psi,\beta, \gamma \ast T^{-\rho_N^\vee}].
  \end{equation*}
  By Proposition \ref{lem:sub-gln:each-psi-beta}, there are fixed seminorms $\mu_1$ on $\mathcal{S}(A)$ and $\mu_2$ on $\mathcal{S}(G)$ so that
  \begin{equation*}
    \kappa(f) \leq \delta_U^{-1/2}(T^{\rho_U^\vee}) \mu_1(\beta) \mu_2(\gamma \ast T^{-\rho_N^\vee}).
  \end{equation*}
  Clearly $\delta_U^{-1/2}(T^{\rho_U^\vee}) = T^{\O(1)}$.  By the definition of $\mathfrak{B}(A) \ni \beta$, we have $\mu_1(\beta) \ll 1$.  It remains to check that $\mu_2(\gamma \ast T^{-\rho_N^\vee}) \ll T^{\O(1)}$.  Indeed, it follows readily from the definitions that
  \begin{itemize}
  \item  $\mu(\gamma) \ll T^{\O(1)}$ for each fixed seminorm $\mu$ on $\mathcal{S}(G)$ and all $\gamma \in \mathfrak{C}(G,-\theta_G(\psi),T,\tfrac{1}{2})$, and that
  \item  translation by $T^{-\rho_N^\vee}$ at most polynomially distorts seminorms on $\mathcal{S}(G)$.
  \end{itemize}
\end{proof}
Lemma \ref{lem:scratch-research:each-fixed-seminorm} implies that each $f \in \mathfrak{F}(U \backslash G, \psi, T)$ is of rapid decay, but in a quantitatively weak sense with respect to $T$.  We will return to the issue of quantitative decay shortly.

\subsection{Equivariance under Fourier transforms}

\begin{lemma}\label{lem:each-f-in}
  Let $(\beta,\gamma) \in \mathfrak{B}(A) \times \mathfrak{C}(G,-\theta_G(\psi),T,\tfrac{1}{2})$ and $w \in W$.  Then
  \begin{equation*}
    \mathcal{F}_{w,\psi^{-1} } \mathcal{J}_T[\psi,\beta,\gamma]
    = \mathcal{J}_T[\psi,{}^w \beta, \gamma].
  \end{equation*}
  In particular, if $\beta \in \mathfrak{B}(A)^W$, then
  \begin{equation*}
    \mathcal{F}_{w,\psi^{-1} } \mathcal{J}_T[\psi,\beta,\gamma]
    = \mathcal{J}_T[\psi,\beta,\gamma]
  \end{equation*}
  for all $w \in W$.
\end{lemma}
\begin{proof}
  By Lemma \ref{lem:cal-J-T-vs-cal-J} and the $G$-equivariance of $\mathcal{F}_{w,\psi}$, it is enough to check that
  \begin{equation*}
    \mathcal{F}_{w,\psi^{-1}} \mathcal{J}[\psi,\beta, \gamma' ]
    =
    \mathcal{J}[{\psi,{}^w \beta, \gamma' }]
  \end{equation*}
  for all $\gamma ' \in \mathcal{S}(G)$, which follows in turn from Lemma \ref{lem:sub-gln:each-psi-beta-1}.
\end{proof}

\subsection{Behavior with respect to differentiation}

\begin{lemma}\label{lem:standard:each-fixed-x}
  For each fixed $x \in \mathfrak{U}(A)$ and each $f \in \mathfrak{F}(U \backslash G, \psi, T)$, we have $L(x) f \in \mathfrak{F}(U \backslash G, \psi, T)$.
\end{lemma}
\begin{proof}
  Write $f = \mathcal{J}_T[\psi,\beta,\gamma] = L(T^{\rho_U^\vee}) \mathcal{J}[\psi_T,\beta,\gamma]$.  By \eqref{eq:lx-mathcald_psi-beta} and the commutativity of $A$, we have
  \begin{equation*}
    L(x) f = L(T^{\rho_U^\vee}) \mathcal{J}[\psi_T, L(x) \beta, \gamma] = \mathcal{J}_T[\psi,L(x) \beta,\gamma].
  \end{equation*}
  Since $L(x) \beta \in \mathfrak{B}(A)$, it follows that $L(x) f \in \mathfrak{F}(U \backslash G, \psi, T)$.
\end{proof}

\begin{lemma}\label{lem:standard:each-fixed-x-1}
  For each fixed $x \in \mathfrak{g}$ and each $f \in \mathfrak{F}(U \backslash G, \psi, T)$, we have
  \begin{equation}\label{eq:crb6g0q5u7}
    \frac{R(x) - \langle T \theta_G(\psi)  , x  \rangle }{T^{1/2}} f \in \mathfrak{F}(U \backslash G, \psi, T).
  \end{equation}
  In particular,
  \begin{equation}\label{eq:crb6g0rgwi}
    R(x) f \in T \mathfrak{F}(U \backslash G, \psi, T).
  \end{equation}
\end{lemma}
\begin{proof}
  We first establish \eqref{eq:crb6g0q5u7}.  Write $f = L(T^{\rho_U^\vee}) \mathcal{J}[\psi_T,\beta,\gamma]$.  By \eqref{eq:rx-mathcald_psi-beta}, we have
  \begin{equation*}
    R(x) f = L(T^{\rho_U^\vee}) \mathcal{J}[\psi_T,\beta,x \ast \gamma],
  \end{equation*}
  so it is enough to check that
  \begin{equation*}
    x \ast \gamma -  \langle T \theta_G(\psi), x \rangle \gamma
    \in T^{1/2} \mathfrak{C}(G,-\theta_G(\psi),T,\tfrac{1}{2}),
  \end{equation*}
  which is the content of Lemma \ref{lem:archimedean-vector-field-act-on-C} applied to $x$, regarded as a right-invariant vector field on $G$.

  To deduce \eqref{eq:crb6g0rgwi} from \eqref{eq:crb6g0q5u7}, we use that $\langle \theta_G(\psi), x \rangle$ is fixed, which gives
  \begin{equation*}
    \langle T \theta_G(\psi), x \rangle f \in T \mathfrak{F}(U \backslash G, \psi, T),
  \end{equation*}
  and that $T \geq 1$, so that
  \begin{equation*}
    T^{1/2} \mathfrak{F}(U \backslash G, \psi, T) \subseteq T \mathfrak{F}(U \backslash G, \psi, T).
  \end{equation*}
\end{proof}

\subsection{Quantitative decay}\label{sec:cnjgzs9fmn}
Lemma \ref{lem:standard2:mathfr-backsl-g} implies that each element of $\mathfrak{F}(U \backslash G, \psi, T)$ lies in $\mathcal{S}(U \backslash G)$, hence decays at infinity.  However, it gives little quantitative information concerning the dependence of that decay upon $T$ besides the soft polynomial control offered by Lemma \ref{lem:scratch-research:each-fixed-seminorm}.  In this section, we show that such elements concentrate near the image of $T^{-\rho_U^\vee}$ in $U \backslash G$.

\subsubsection{Distance functions and restriction norms}\label{sec:dist-funct-restr}
In what follows, we use the standard nontrivial additive character $\psi_F$ of $F$ to identify Pontryagin duals of vector spaces over $F$ with linear duals.

Recall the notation $\psi$ and $\theta_G(\psi)$ from \S\ref{sec:analysis-U-G-quantitative-setup}.

We fix on the compact real manifold $Q \backslash G$ a locally Euclidean metric.  We denote by $\dist_{Q \backslash G}(x,y)$ the associated distance between a pair of points $x, y \in Q \backslash G$.  In particular, $\dist_{Q \backslash G}(x,y) \ll 1$ for all $x$ and $y$.  We use this notation also when $x$ and/or $y$ belongs to $G$, referring in that case to their images in $Q \backslash G$.

We fix a Euclidean norm $|.|$ on $\Lie(G)$, hence a dual norm $|.|$ on $\Lie(G)^\wedge$.  We denote by $\nu : Q \backslash G \rightarrow \mathbb{R}_{\geq 0}$ the function sending $Q g$ to the norm of the restriction of $\theta$ to the subspace $\Ad(g)^{-1} \Lie(Q)$, i.e.,
\begin{equation*}
  \nu(g) := \sup_{0 \neq y \in \Ad(g)^{-1} \Lie(Q)}
  \frac{|\langle \theta, y \rangle|}{ |y|}.
\end{equation*}
We may also think of $\nu$ as a left-$Q$-invariant function on $G$.

The following estimate will play an important role in the ``integration by parts'' arguments of \S\ref{sec:integration-parts}.
\begin{lemma}\label{lem:sub-gln:g-in-q}
  For $g \in Q \backslash G$, we have
  \begin{equation}\label{eq:dist_q-backslash-gg}
    \dist_{Q \backslash G}(g, 1) \asymp \nu(g).
  \end{equation}
\end{lemma}
\begin{proof}
  Since $\psi$ is fixed, we abbreviate $\theta := \theta_G(\psi)$ here.

  We observe first that $\nu(g) = 0$ only if $g \in Q$.  Indeed, if $\nu(g) = 0$, then $\theta$ has trivial restriction to $\Ad(g)^{-1} \Lie(Q)$, or equivalently, $\Ad^*(g) \theta$ has trivial restriction to $\Lie(Q)$.  Let us fix a faithful representation of $\mathfrak{g}$.  The corresponding trace pairing then induces an equivariant identification of $\Lie(G)^\wedge$ with $\Lie(G)$.  With respect to this identification, the orthogonal complement of $\Lie(Q)$ is its maximal nilpotent subalgebra $\Lie(U)$.  Thus $\theta$ and $\Ad^*(g) \theta$ both identify with regular nilpotent elements of $\Lie(U)$.  By Lemma \ref{lem:scratch-research:let-mathfrakb-be} below, we deduce that $g \in Q$, as required.


  The function $\nu$ is continuous, vanishes only at the trivial coset, and is defined on the compact set $Q \backslash G$.  It follows that for each fixed neighborhood $E$ of the trivial coset in $Q \backslash G$, we have $\nu(g) \asymp 1$ for all $g \notin E$.  The claimed estimate \eqref{eq:dist_q-backslash-gg} thus holds except possibly when $\dist_{Q \backslash G}(g,1) \lll 1$.  It remains to address the latter case.

  By the Bruhat decomposition, the natural map $N \rightarrow Q \backslash G$ is an injective immersion whose image contains an open neighborhood of the identity coset.  It will thus suffice to check that if $n \in N$ with $n \simeq 1$, then $\dist_{N}(n,1) \asymp \nu(n)$ (with $\dist_N$ some fixed Euclidean distance, as in \S\ref{sec:distances}).

  To that end, let $\mu : N \rightarrow \Lie(Q)^\wedge$ denote the composition
  \begin{equation*}
    N \xrightarrow{n \mapsto \Ad^*(n) \theta } \Lie(G)^\wedge \xrightarrow{\text{restriction}} \Lie(Q)^\wedge.
  \end{equation*}
  This map has the following properties.
  \begin{enumerate}[(i)]
  \item $\mu(1) = 0$.  Indeed, $\theta|_{\Lie(Q)} = 0$.
  \item For each fixed compact subset $\Omega$ of $N$, we have $\nu(n) \asymp \mu(n)$ for all $n \in \Omega$.  Indeed, since $n$ lies in a fixed compact set, the norm of the restriction of $\theta$ to $\Ad(n)^{-1} \Lie(Q)$ is comparable to the norm of the restriction of $\Ad^*(n) \theta$ to $\Lie(Q)$, i.e., to $|\mu(n)|$.
  \end{enumerate}
  It will thus suffice to show that if $n \in N$ with $n \simeq 1$, then $\dist_N(n,1) \asymp \mu(n)$.

  To establish this last estimate, it suffices to verify that $\mu$ is a local immersion at the identity, i.e., that the derivative $d \mu$ of $\mu$ at $1$ is injective.  That derivative is the composition
  \begin{equation*}
    \Lie(N) \xrightarrow{x \mapsto [x,\theta]} \Lie(G)^\wedge \xrightarrow{\text{restriction}} \Lie(Q)^\wedge,
  \end{equation*}
  where we write $[x,\theta] := \ad^*(x)\theta$.
  Its kernel consists of all $x \in \Lie(N)$ for which $\langle [x,\theta], y \rangle = 0$ for all $y \in \Lie(Q)$.  Since $\langle [x,\theta], y \rangle = -\langle x, [y,\theta] \rangle$, we have in fact
  \begin{equation*}
    \ker(d \mu) = \Lie(N) \cap [\Lie(Q), \theta]^\perp.
  \end{equation*}
  Using the trace pairing, we identify $\theta$ with an element of $\Lie(U)$.

  We claim that $[\Lie(Q), \theta] = \Lie(U)$.  Indeed, let $n$ denote the rank of $G$.  (Here we define ``rank'' and ``dimension'' either over $F$ or $\mathbb{R}$; as long as we are consistent, it does not matter which.)  Since $\theta$ is regular, its $\Lie(G)$-centralizer is $n$-dimensional; since $[\Lie(Q), \theta] \subseteq \Lie(U)$ and $[\Lie(Q):\Lie(U)] = n$, the claim follows.

  The trace pairing induces a perfect pairing between $\Lie(N)$ and $\Lie(U)$, so the claim implies that
  \begin{equation*}
    \Lie(N) \cap [\Lie(Q), \theta]^\perp = \Lie(N) \cap \Lie(U)^\perp
  \end{equation*}
  is trivial, hence that $d \mu$ is injective, as required.
\end{proof}

We record below some linear algebraic facts that were needed above.  We first recall
\cite[Cor 5.6]{MR0114875}:
\begin{lemma}\label{lem:kostant-regular-nilpotent}
  Let $x$ be a nilpotent element of a complex semisimple Lie algebra $\mathfrak{g}$.  Then $x$ is regular if and only if it is contained in exactly one maximal nilpotent subalgebra.
\end{lemma}
As a consequence, we deduce:
\begin{lemma}\label{lem:scratch-research:let-mathfrakb-be}
  Let $G$ be a quasi-split real reductive group, with Lie algebra $\mathfrak{g}$.
  \begin{enumerate}[(i)]
  \item \label{itm:scratch-research:lemma-general-linear-remedial} Let $x$ be a regular nilpotent element of $\mathfrak{g}$.  Then $x$ is contained in at most one Borel subalgebra of $\mathfrak{g}$.
  \item Let $B$ be a Borel subgroup of $G$, with corresponding Borel subalgebra $\mathfrak{b}$ of $\mathfrak{g}$.  Let $x$ be a regular nilpotent of $\mathfrak{g}$ that is contained in $\mathfrak{b}$.  Let $g \in G$ be such that $\Ad(g) x \in \mathfrak{b}$.  Then $g \in B$.
  \end{enumerate}
\end{lemma}
\begin{proof}
  (i): It suffices to prove this after passing to the semisimple part of $\mathfrak{g}$ and extending scalars to $\mathbb{C}$, after which we can apply Lemma \ref{lem:kostant-regular-nilpotent}.  (We remark that for the cases $G = \GL_n$ of primary interest, the required conclusion is an easy exercise with Jordan form.)


  (ii): Our hypotheses imply that $\Ad(g)^{-1} \mathfrak{b}$ and $\mathfrak{b}$ are both Borel subalgebras containing $x$, so by part (i), we have $\Ad(g)^{-1} \mathfrak{b} = \mathfrak{b}$.  Thus $g$ normalizes the connected group $B$, which forces $g \in B$ \cite[Thm 11.16]{MR1102012}.
\end{proof}

\subsubsection{Integration by parts}\label{sec:integration-parts}

\begin{lemma}\label{lem:scratch-research:let-h-in}
  Let $f \in \mathfrak{F}(U \backslash G, \psi, T)$.  If $h \in U \backslash G$ satisfies $\dist_{Q \backslash G}(h,1) \gg T^{-1/2+\eps}$ for some fixed $\eps > 0$, then for each fixed $x \in \mathfrak{U}(A)$ and $y \in \mathfrak{U}(G)$,
  \begin{equation*}
    L(x) R(y) f(h) \ll
    \|h\|_{U \backslash G}^{-\infty}
    T^{-\infty}.
  \end{equation*}
\end{lemma}
\begin{proof}
  Using Lemmas \ref{lem:standard:each-fixed-x} and \ref{lem:standard:each-fixed-x-1}, we may and shall reduce to the case $(x,y) = (1,1)$.  Our task is then to show that for each fixed $m \geq 0$ and each $f \in \mathfrak{F}(U \backslash G, \psi, T)$, we have
  \begin{equation}\label{eq:fg-ll-g_u}
    f(h) \ll \|h\|_{U \backslash G}^{-m} T^{-m}.
  \end{equation}

  By Lemma \ref{lem:scratch-research:each-fixed-seminorm} applied to $\kappa(f) := \sup_{g \in U \backslash G} \|g\|_{U \backslash G}^{m+1} |f(g)|$, we have
  \begin{equation}\label{eq:fg-ll-g_u-1}
    f(h) \ll \|h\|_{U \backslash G}^{-m - 1} T^{C_m}
  \end{equation}
  for some fixed $C_m \geq 0$.

  If $\|h\|_{U \backslash G} \geq T^{m + C_m}$, then \eqref{eq:fg-ll-g_u-1} implies \eqref{eq:fg-ll-g_u}.  We focus henceforth on the remaining range $\|h\|_{U \backslash G} < T^{m + C_m}$.  We then have $T^{-m(m+ C_m) - m} \leq \|h\|_{U \backslash G}^{-m} T^{-m}$, so it will suffice to show the following:
  \begin{itemize}
  \item for each fixed $\eps > 0$, $C \geq 0$ and $m \in \mathbb{Z}$, we have $f(h) \ll T^{-m}$ whenever $\|h\|_{U \backslash G} < T^{C}$ and $\dist_{Q \backslash G}(h,1) \gg T^{-1/2+\eps}$.
  \end{itemize}

  Write $f = \mathcal{J}_T[\psi,\beta,\gamma] = L(T^{\rho_U^\vee}) \mathcal{J}[\psi_T, \beta, \gamma]$ with $(\beta,\gamma) \in \mathfrak{B}(A) \times  \mathfrak{C}(G,-\theta_G(\psi),T,\tfrac{1}{2})$.  By definition,
  \begin{equation*}
    \delta_U^{1/2}(T^{\rho_U^\vee})
    f(h)
    =
    \int_{g \in G}
    \gamma(g)
    \phi_{h,\beta}(g),
    \quad
    \phi_{h,\beta}(g) :=
    \mathcal{J}[\psi_T, \beta](T^{\rho_U^\vee} h g) \, d g.
  \end{equation*}
  Let us henceforth fix $C \geq 0$ and $\eps > 0$.  It will then suffice to show that for each fixed $m \in \mathbb{Z}_{\geq 0}$, the following assertion $S(m)$ holds:
  \begin{itemize}
  \item for all $\beta \in \mathfrak{B}(A)$, $\gamma \in \mathfrak{C}(G,-\theta_G(\psi),T,\tfrac{1}{2})$ and $h \in G$ (or $h \in U \backslash G$) with $\|h\|_{U \backslash G} < T^C$ and $\dist_{Q \backslash G}(h,1) \gg T^{-1/2+\eps}$, we have
    \begin{equation*}
      \int_{g \in G} \gamma(g) \phi_{h,\beta}(g) \ll T^{-m \eps}
    \end{equation*}
  \end{itemize}
  We will do so by induction on $m$.  The idea is that integration by parts will save one factor of $T^{-\eps}$ per inductive step.  The base case $m = 0$ follows from Lemma \ref{lem:standard:each-gamma-in} and \eqref{eq:mathc-beta_-leq}, which give $\|\gamma \|_1 \ll 1$ and $\|\phi_{h,\beta}\|_{\infty} \ll 1$.  Our task is thus to verify the inductive step: fixing $m \in \mathbb{Z}_{\geq 0}$ and assuming $S(m)$, we will deduce $S(m+1)$.

  Retain the hypotheses of $S(m+1)$.  By Lemma \ref{lem:sub-gln:g-in-q}, we have $\nu(h) \gg T^{-1/2+\eps}$.  By the definition of $\nu$, we may thus find
  \begin{equation*}
    0 \neq y \in \Ad(h)^{-1} \Lie(Q) \subseteq \Lie(G)
  \end{equation*}
  with
  \begin{equation*}
    \frac{\langle T \theta_G(\psi), y \rangle}{\lvert y \rvert} \gg T^{1/2 + \eps}.
  \end{equation*}
  By scaling $y$ suitably, we may arrange that $y \ll T^{-1/2}$ and $\langle T \theta_G(\psi), y \rangle \gg T^{\eps}$.  We integrate by parts with respect to the right-invariant vector field $Y$ on $G$ corresponding to $y$, namely $Y(g) := \partial_{t=0} \exp(t y) g$.  This gives
  \begin{equation*}
    \int_{g \in G} Y \gamma(g) \phi_{h, \beta } (g)
    +
    \int_{g \in G} \gamma(g) Y \phi_{h, \beta } (g) = 0.
  \end{equation*}
  By Proposition \ref{lem:archimedean-vector-field-act-on-C-2}, we have
  \begin{equation*}
    Y \gamma = -\langle T \theta_G(\psi), y \rangle \gamma + \gamma_0
  \end{equation*}
  for some
  \begin{equation*}
    \gamma_0 \in
    \lvert y \rvert T^{-1/2}
    \mathfrak{C}(G,-\theta_G(\psi),T,\tfrac{1}{2})
    \subseteq
    \mathfrak{C}(G,-\theta_G(\psi),T,\tfrac{1}{2}).
  \end{equation*}
  Thus
  \begin{equation}\label{eq:langle-t-theta-1}
    \langle T \theta_G(\psi) , y \rangle \int_{g \in G} \gamma(g) \phi_{h,\beta}(g)
    =
    \int_{g \in G} \gamma_0(g) \phi_{h,\beta}(g) + \int_{g \in G} \gamma(g) Y \phi_{h,\beta}(g).
  \end{equation}
  Since $\langle T \theta_G(\psi), y \rangle \gg T^{\eps}$, our task reduces to verifying that each term on the RHS of \eqref{eq:langle-t-theta-1} is $\ll T^{-m \eps}$.  For the first term, the required bound is the content of our inductive hypothesis $S(m)$, applied with $\gamma_0$ in place of $\gamma$.  It thus remains only to show that
  \begin{equation}\label{eq:int-_g-in-1-9}
    \int_{g \in G} \gamma(g) Y \phi_{h,\beta}(g) \ll T^{-m \eps}.
  \end{equation}
  To that end, we write
  \begin{equation*}
    y = \Ad(h)^{-1} x,
  \end{equation*}
  with $x \in \Lie(Q)$, and decompose
  \begin{equation*}
    x = x_U + x_A, \quad
    (x_U, x_A) \in \Lie(U)  \times \Lie(A).
  \end{equation*}
  Then
  \begin{align*}
    Y \phi_{h,\beta}(g)
    &= \partial_{t=0} \mathcal{J}[\psi_T,\beta](T^{\rho_U^\vee} h \exp(t y) g) \\
    &=
      \partial_{t=0} \mathcal{J}[\psi_T, \beta](T^{\rho_U^\vee} \exp(t x) h g)
    \\
    &=
      \partial_{t=0} \mathcal{J}[\psi_T, \beta](\exp(t x_A) T^{\rho_U^\vee} h g)
    \\
    &=
      \mathcal{J}[\psi_T, L(x_A) \beta](T^{\rho_U^\vee} h g),
  \end{align*}
  where
  \begin{itemize}
  \item in the third step, we have used that $\mathcal{J}[\psi_T,\beta]$ is left-$U$-invariant (and that $A$ is commutative), and
  \item in the fourth step, we have appealed to \eqref{eq:lx-mathcalj_psi-beta}.
  \end{itemize}
  We claim that $x_A \ll T^{-1/2}$.  To see this, it is equivalent to show that each coefficient of the characteristic polynomial of $T^{1/2} x_A$ (with respect to the given linear embedding of $\Lie(G)$) is $\O(1)$.  The characteristic polynomials of the elements $x_A, x$ and $y$ are all the same, so the conclusion follows from the fact that $y \ll T^{-1/2}$.  In particular, $L(x_A) \beta = T^{-1/2} \beta_0$ for some $\beta_0 \in \mathfrak{B}(A)$, and so by our inductive hypothesis $S(m)$,
  \begin{equation*}
    \int_{g \in G} Y \phi_{h,\beta}(g) =
    T^{-1/2}
    \int_{g \in G} \phi_{h,\beta_0}(g) \ll T^{-1/2 - m \eps}.
  \end{equation*}
  We thereby deduce the required estimate \eqref{eq:int-_g-in-1-9} (in stronger form).
\end{proof}

\subsubsection{Decomposition into essential and negligible parts}

The next proposition is an intermediate estimate.  It decomposes elements of $\mathfrak{F}(U \backslash G,\psi,T)$ into an essential part, localized near the expected region in $U \backslash G$, and a negligible remainder.  We will use it below to derive bounds for Mellin components and hence the Whittaker bounds that are the main outputs of this section; it will not be exported beyond that role.

\begin{proposition}\label{prop:standard:each-f-in}
  Each $f \in \mathfrak{F}(U \backslash G, \psi, T)$ may be decomposed as $f = f^\sharp + f^\flat$, where $f^\sharp, f^{\flat} \in \mathcal{S}(U \backslash G)$ have the following properties:
  \begin{enumerate}[(i)]
  \item \label{itm:standard:prop-decompose-frakF-4} For $g \in G$, we have $f^\sharp (g) \neq 0$ only if $\dist_{Q \backslash G}(g,1) \ll T^{-1/2 + o(1)}$.  More precisely, for $(a,n) \in A \times N$, we have
    \begin{equation*}
      f^\sharp (T^{-\rho_U^\vee} a n) \neq 0 \implies \|a\| \ll T^{o(1)}, \quad n = 1 + \O(T^{-1/2+o(1)}).
    \end{equation*}
  \item $\|f^\sharp \|_{\infty} \ll \delta_U^{1/2}(T^{-\rho_U^\vee})$.
  \item $\kappa(f^{\flat}) \ll T^{-\infty}$ for each fixed seminorm $\kappa$ on $\mathcal{S}(U \backslash G)$.
  \end{enumerate}
\end{proposition}
\begin{proof}
  Write $f = L(T^{\rho_U^\vee }) \mathcal{J}[\psi_T, \beta, \gamma]$, as usual.  We first decompose $\beta = \beta^\sharp + \beta^\flat$ according to Lemma \ref{lem:standard:each-beta-in}.  By combining the conclusion of that lemma with the \emph{a priori} bounds of Lemma \ref{lem:scratch-research:each-fixed-seminorm}, we see that
  \begin{equation*}
    \kappa(L(T^{\rho_U^\vee}) \mathcal{J}[\psi_T,\beta^\flat,\gamma]) \ll T^{-\infty}
  \end{equation*}
  for each fixed seminorm $\kappa$ on $\mathcal{S}(U \backslash G)$.  The contribution of $\beta^{\flat}$ may thus be absorbed into the term $f^{\flat}$ that we will produce below, and so we may assume that $\beta = \beta^\sharp$.  In particular, we have reduced to the case that $\beta(a) \neq 0$ only if $\|a\| \ll T^{o(1)}$.

  Fix $\phi \in C_c^\infty(\mathbb{R})$ taking the value $1$ in a neighborhood of zero.  For small $t > 0$, define $f_t : U \backslash G \rightarrow \mathbb{C}$ by
  \begin{equation*}
    f_t(g) :=
    f(g) \phi (\dist_{Q \backslash G}(g,1) / t).
  \end{equation*}
  Then $f_t$ is smooth and of compact support, hence lies in $\mathcal{S}(U \backslash G)$, as does the difference $f - f_t$.  We specialize to $t = t(\eps) := T^{-1/2+\eps}$, with $\eps > 0$ to be determined.
  The integration-by-parts estimate of Lemma \ref{lem:scratch-research:let-h-in}, applied on the support of $1 - \phi(\dist_{Q \backslash G}(\cdot,1)/t)$ and combined with the product rule, shows that, for each fixed $\eps > 0$, each fixed Schwartz seminorm of $f - f_t$ is $\ll T^{-\infty}$.  Fixing a defining countable collection $\{\kappa_j\}_j$ of seminorms on $\mathcal{S}(U \backslash G)$, we deduce in particular that, for each fixed $\eps > 0$, we have $\kappa_j(f-f_t) \leq T^{-1/\eps}$ for each $j \leq 1/\eps$.  By overspill, we may thus choose $\eps > 0$ with $\eps \lll 1$ so that, for each fixed seminorm $\kappa$ on $\mathcal{S}(U \backslash G)$, we have $\kappa(f-f_t) \ll T^{-\infty}$.  With this choice, set
  \begin{equation*}
    f^\sharp := f_t,
    \quad
    f^\flat := f - f_t.
  \end{equation*}
  By construction, assertion (iii) holds; it remains to verify assertions (i) and (ii).
  \begin{enumerate}[(i)]
  \item The first claimed support property of $f^\sharp $ follows immediately from its construction and the compact support of $\phi$.  For the second assertion, note that, by construction, for $(a,n) \in A \times N$,
    \begin{equation}\label{eq:f-sharp-t}
      f^\sharp (T^{-\rho_U^\vee} a n) = f( T^{-\rho_U^\vee} a n) \phi\left(\frac{\dist_{Q \backslash G}(n,1) }{ T^{-1/2+\eps}}\right).
    \end{equation}
    Since the natural map $N \rightarrow Q \backslash G$ is an injective immersion (at the identity), we have
    \begin{equation*}
      \dist_{Q \backslash G} (n,1) \ll T^{-1/2+\eps} \implies \|n - 1\| \ll T^{-1/2 + \eps}.
    \end{equation*}
    Thus the fixed compact support of $\phi$ forces $f^\sharp (T^{- \rho_U^\vee } a n )$ to vanish unless $n = 1 + \O(T^{-1/2+\eps})$.

    It remains to constrain the support of $f^\sharp (T^{- \rho_U^\vee } a n)$ with respect to the variable $a$.
    Our task is to show that for $(a,n) \in A \times N$,
    \begin{equation}\label{eq:f-t-}
      f^\sharp (T^{- \rho_U^\vee } a n) \neq 0
      \implies
      \|a\| \ll T^{o(1)}.
    \end{equation}
    The support condition that we imposed on $\beta$ implies that
    \begin{equation}\label{eq:mathc-betag-neq}
      (u,b,m) \in U \times A \times N,\;
      \mathcal{J}[\psi_T,\beta](u b m) \neq 0
      \implies
      \|b\| \ll T^{o(1)}.
    \end{equation}

    We observe that if $(a,g) \in A \times G$ with $g = 1 + o(1)$, then $(a g)_A$ is defined and satisfies
    \begin{equation}\label{eq:mathfraka-a-g}
      a^{-1} (a g)_A = 1 + o(1)
    \end{equation}
    To see this, we use that the multiplication map $U \times A \times N \rightarrow G$ is an injective immersion with open image to write $g = u b n$ with $(u,b,n) \in U \times A \times N$.  Then, since $A$ normalizes $U$ and $N$, we have $(a g)_A = a b$, hence $a^{-1} (a g)_A = b = 1 + o(1)$.

    We observe next that for $(a,n) \in A \times N$,
    \begin{align*}
      f(T^{- \rho_U^\vee } a n)
      &= L(T^{\rho_U^\vee}) R(\gamma) \mathcal{J}[\psi_T,\beta](T^{- \rho_U^\vee } a n)
      \\
      &=
        R(\gamma) L(T^{\rho_U^\vee})  \mathcal{J}[\psi_T,\beta](T^{- \rho_U^\vee } a n)
      \\
      &=
        \delta_U^{-1/2}(T^{\rho_U^\vee})
        \int_{g \in G}  \mathcal{J}[\psi_T,\beta]( a n g) \gamma(g),
    \end{align*}
    hence
    \begin{equation}\label{eq:f-t--1}
      f(T^{- \rho_U^\vee } a n)
      =
      \delta_U^{-1/2}(T^{\rho_U^\vee})
      \int_{g \in G}
      \mathcal{J}[\psi_T,\beta](a n g) \gamma(g).
    \end{equation}
    Suppose now that $f^\sharp(T^{-\rho_U^\vee} a n) \neq 0$.  We saw above that then $n = 1 + \O(T^{-1/2+\eps})$, hence in particular $n = 1 + o(1)$.  Since $f^\sharp = f \cdot \phi(\cdots)$, this also implies $f(T^{-\rho_U^\vee} a n) \neq 0$, so by \eqref{eq:f-t--1}, there exists $g \in \supp(\gamma)$ with $\mathcal{J}[\psi_T,\beta](a n g) \neq 0$.  By the definition of $\mathfrak{C}(G,-\theta_G(\psi),T,1/2)$ (Definition \ref{defn:standard:let-v-be}, \eqref{eq:mathfrakcg-theta-t}), we have $g = 1 + \O(T^{-1/2}) = 1 + o(1)$, hence $n g = 1 + o(1)$.
    Applying \eqref{eq:mathfraka-a-g} with $g$ replaced by $n g$, we obtain
    \begin{equation*}
      a^{-1} (a n g)_A = 1 + o(1).
    \end{equation*}
    On the other hand, \eqref{eq:mathc-betag-neq} gives $\|(a n g)_A\| \ll T^{o(1)}$.  It follows that $\|a\| \ll T^{o(1)}$, proving \eqref{eq:f-t-}.

  \item
    Since $\phi$ is fixed, we have $\|\phi \|_\infty \ll 1$, hence $\|f^\sharp \| _\infty \ll \|f\| _\infty$.  We recall that $f =  R(\gamma) L(T^{\rho_U^\vee }) \mathcal{J}[\psi_T, \beta]$, hence
    \begin{equation*}
      \|f\|_\infty \leq \|\gamma\|_1 \|L(T^{\rho_U^\vee }) \mathcal{J}[\psi_T, \beta]\|_{\infty}.
    \end{equation*}
    By Lemma \ref{lem:standard:each-gamma-in}, we have $\|\gamma \|_1 \ll 1$.  For $(a,n) \in A \times N$, we have
    \begin{equation*}
      L(T^{\rho_U^\vee}) \mathcal{J}[\psi_T,\beta](a n)
      =
      \delta_U^{-1/2}(T^{\rho_U^\vee})  \delta_U^{1/2}(T^{\rho_U^\vee} a) \beta(T^{\rho_U^\vee} a) \psi_T(n).
    \end{equation*}
    Since $\|\delta_U^{1/2} \beta \|_{\infty} \ll 1$, it follows that
    \begin{equation*}
      \|L(T^{\rho_U^\vee}) \mathcal{J}[\psi_T,\beta]\|_{\infty} \ll \delta_U^{-1/2}(T^{\rho_U^\vee}).
    \end{equation*}
    The required bound for $\|f\|_\infty$ follows.
  \end{enumerate}

\end{proof}

\subsection{$L^2$-estimates}

The following bound records the natural $L^2$-scale suggested by Proposition \ref{prop:standard:each-f-in}.  We include it only for completeness, as our later arguments use sharper Mellin and Whittaker variants of this estimate.
\begin{lemma}
  For $f \in \mathfrak{F}(U \backslash G, \psi, T)$, we have $\int_{U \backslash G} |f|^2 \ll T^{-\dim(N)/2 + o(1)}$.
\end{lemma}
\begin{proof}
  We decompose $f = f^\sharp + f^\flat$ as in Proposition \ref{prop:standard:each-f-in}.  Writing $\|.\|_2$ for the $L^2$-norm on $U \backslash G$, the triangle inequality gives $\|f\|_2 \leq \|f^\sharp \|_2 + \|f^\flat \|_2$.  The conclusion of the cited proposition gives $\|f^\flat \|_2 \ll T^{-\infty}$.  By our Haar measure convention \eqref{eq:int-_g-in-U-bakslash-G}, we have
  \begin{equation*}
    \|f^\sharp \|_2^2
    =
    \int_{a \in A}
    \int_{n \in N}
    |f^\sharp|^2 (a n) \, \frac{d a }{ \delta_U(a)} \, d n.
  \end{equation*}
  The change of variables $a \mapsto T^{-\rho_U^\vee} a$ gives
  \begin{equation*}
    \|f^\sharp \|_2^2
    =
    \delta_U(T^{\rho_U^\vee})
    \int_{a \in A}
    \int_{n \in N}
    |f^\sharp|^2 (T^{-\rho_U^\vee} a n) \, \frac{d a }{ \delta_U(a)} \, d n.
  \end{equation*}
  By construction, $f^\sharp (T^{- \rho_U^\vee } a n ) \neq 0$ only if $\|a\| \ll T^{o(1)}$ and $n = 1 + \O(T^{-1/2+o(1)})$.  For such $a$, we have $\delta_U(a) = T^{o(1)}$.  On the other hand, $\|f^\sharp \|  _\infty \ll \delta_U(T^{\rho_U^\vee})^{-1/2}$.  These size and support estimates yield the required bound.
\end{proof}

\subsection{Basic Mellin component bounds}

The following two lemmas specialize to $f \in \mathfrak{F}(U \backslash G,\psi,T)$ the general vertical-decay bound for Mellin components recorded in Lemma \ref{lem:standard:each-comp-mathc}: they convert the support and size estimates from Proposition \ref{prop:standard:each-f-in}, expressed in Bruhat coordinates, into pointwise and $L^2(K)$ bounds for the Mellin components $f[\chi]$.  These lemmas will serve as the basis for our subsequent analysis of the Whittaker functions attached to $f[\chi]$ (see \S\ref{sec:whittaker-upper-bounds} and \S\ref{sec:whitt-funct-lower}).

For what follows, we note that the condition $\Re(\chi) \ll 1$ (i.e., $\lvert \Re(\chi) \rvert \leq C$ for some fixed norm $\lvert . \rvert$ on $\mathfrak{a}^\ast$ and some $C \geq 0$) is equivalent to $\Re(\chi)$ belonging to some fixed compact subset of $\mathfrak{a}^*$; these two formulations of that condition are used interchangeably.

\begin{lemma}\label{lem:standard:let-chi-in}
  Let $\chi \in \mathfrak{X}(A)$ with $\Re(\chi) \ll 1$.  Let $f \in \mathfrak{F}(U \backslash G, \psi, T)$.  Let $k \in K$.  Then
  \begin{equation*}
    f[\chi](k) \ll T^{-\infty} C(\chi)^{-\infty}
  \end{equation*}
  unless $\dist_{Q \backslash G}(k,1) \ll T^{-1/2+o(1)}$, in which case
  \begin{equation*}
    f[\chi](k) \ll T^{\langle \rho_U^\vee , \Re(\chi) \rangle + o(1)} C(\chi)^{-\infty}.
  \end{equation*}
\end{lemma}
\begin{proof}
  First, for each fixed $m \in \mathbb{Z}_{\geq 0}$, we may find a fixed element $x \in \mathfrak{U}(A)$ so that $d \chi(x) \gg C(\chi)^m$, hence
  \begin{equation*}
    f[\chi](k) \ll C(\chi)^{-m} L(x) f[\chi](k).
  \end{equation*}
  By Lemma \ref{lem:standard:each-fixed-x}, we may thus reduce to verifying the ostensibly weaker estimates obtained by omitting the ``factor'' $C(\chi)^{-\infty}$.  Next, write $f = f^\sharp + f^\flat$ as in Proposition \ref{prop:standard:each-f-in}.  Using Lemma \ref{lem:standard:each-comp-mathc}, we see that $f^\flat$ satisfies the required conclusions (in stronger form).  It remains to treat the contribution of $f^\sharp$.  To that end, we change variables in the definition to see that
  \begin{equation*}
    f^\sharp[\chi](k) = \delta_U^{1/2} \chi  (T^{\rho_U^\vee}) \int_{a \in A} \delta_U^{1/2} \chi(a) f^\sharp(T^{- \rho_U^\vee } a^{-1} k) \, d a.
  \end{equation*}
  By construction, we have $f^\sharp  (T^{- \rho_U^\vee } a^{-1} k) \neq 0$ only if $\dist_{Q \backslash G}(k,1) \ll T^{-1/2 + o(1)}$.  In that case, we may write $k = u b n$ with $(u,b,n) \in U \times A \times N$, $\|b\| \ll T^{o(1)}$ and $n = 1 + \O(T^{-1/2+o(1)})$.  Then
  \begin{equation*}
    f^\sharp (T^{- \rho_U^\vee } a^{-1} k)
    =
    f^\sharp (T^{- \rho_U^\vee } a^{-1} b n).
  \end{equation*}
  This vanishes unless $\|a^{-1} b \| \ll T^{o(1)}$, in which case $\|a\| \ll T^{o(1)}$.  Since $\Re(\chi) \ll 1$ and $\|f^\sharp \| _\infty \ll \delta_U^{-1/2}(T^{\rho_U^\vee}) T^{o(1)}$, we obtain
  \begin{equation}\label{eq:f-sharp-chik}
    f^\sharp [\chi](k) \ll \chi (T^{\rho_U^\vee })
    T^{o(1)}
    \ll
    T^{\langle \rho_U^\vee , \Re(\chi) \rangle + o(1)},
  \end{equation}
  as required.
\end{proof}

\begin{lemma}\label{lem:scratch-research:let-f-in}
  Let $f \in \mathfrak{F}(U \backslash G, \psi, T)$.  Let $\chi \in \mathfrak{X}(A)$ with $\Re(\chi) \ll 1$.  Then for fixed $x \in \mathfrak{U}(A)$, $m \in \mathbb{Z}_{\geq 0}$ and $y_1,\dotsc,y_m \in \Lie(G)$,
  \begin{equation*}
    \int_K \left\lvert L(x) R(y) f[\chi] \right\rvert^2 \ll T^{2 \langle \rho_U^\vee, \Re(\chi) \rangle -\dim(N)/2 + 2 m + o(1)} C(\chi)^{-\infty}.
  \end{equation*}
  In particular,
  \begin{equation}\label{eq:int_k-fchi2-ll}
    \int_K \left\lvert f[\chi] \right\rvert^2 \ll T^{2 \langle \rho_U^\vee, \Re(\chi) \rangle -\dim(N)/2 + o(1)} C(\chi)^{-\infty}.
  \end{equation}
\end{lemma}
\begin{proof}
  Using Lemmas \ref{lem:standard:each-fixed-x} and \ref{lem:standard:each-fixed-x-1}, we may reduce to the case $x = 1$, $m = 0$.  We then apply Lemma \ref{lem:standard:let-chi-in} and integrate over $K$, using that for $r \lll 1$ (in particular, $r = T^{-1/2+\eps}$),
  \begin{equation*}
    \vol \{k \in K : \dist_{Q \backslash G}(k,1) \leq r \} \asymp r^{\dim(N)}.
  \end{equation*}
\end{proof}

\subsection{Whittaker function upper bounds}\label{sec:whittaker-upper-bounds}

We now record the first of the two main exports of this section: a uniform upper bound that will be carried into \S\ref{sec:constr-test-vect-proofs} (where the role of the present section's $G$ is played by $H$).
\begin{proposition}\label{lem:all-f-in}
  For all $f \in \mathfrak{F}(U \backslash G, \psi, T)$, all $g \in Q$ with $g = 1 + o(T^{-1/2})$, and all $\chi \in \mathfrak{X}(A)$ with $\Re(\chi) \ll 1$, we have
  \begin{equation}\label{eq:wfchi-psi_tg-ll}
    W[f[\chi], \psi_T](g)
    \ll T^{\langle \rho_U^\vee, \Re(\chi) \rangle - \dim(N)/2 + o(1)} C(\chi)^{-\infty}.
  \end{equation}
\end{proposition}
\begin{proof}
  As in the proof of Lemma \ref{lem:standard:let-chi-in}, we may apply suitable fixed elements of $\mathfrak{U}(A)$ and use Lemma \ref{lem:standard:each-fixed-x} to reduce to proving the estimate without the ``factor'' $C(\chi)^{-\infty}$.

  Write $f = \mathcal{J}_T[\psi, \beta, \gamma]$.  By \eqref{eq:rgamm-mathc-beta}, we have $R(g) f = \mathcal{J}_T[\psi,\beta, g \ast \gamma]$.  By Proposition \ref{lem:sub-gln:let-gamma-in} and our assumption $g = 1 + o(T^{-1/2})$, we have
  \begin{equation*}
    g \ast \gamma \in \mathfrak{C}(G,-\theta_G(\psi),T,\tfrac{1}{2}),
  \end{equation*}
  hence likewise
  \begin{equation*}
    R(g) f \in \mathfrak{F}(U \backslash G, \psi, T).
  \end{equation*}
  We may thus reduce further to the case $g = 1$.  Our task is then to show that
  \begin{equation*}
    W[f[\chi], \psi_T](1)
    \ll T^{\langle \rho_U^\vee, \Re(\chi) \rangle - \dim(N)/2 + o(1)}.
  \end{equation*}
  To that end, we write $f = f^\sharp + f^\flat$ as in Proposition \ref{prop:standard:each-f-in}.  By Lemma \ref{lem:standard:each-comp-mathc-1}, we see that $f^\flat$ gives an acceptable contribution.  It remains to estimate $W [f^\sharp [\chi ], \psi_T ] (1)$.  Using the support condition on $f^\sharp$ stated in Proposition \ref{prop:standard:each-f-in}, we see that $f^\sharp [\chi ]|_{N}$ is supported on $1 + \O(T^{-1/2+o(1)})$.
  In particular, $f^\sharp [\chi ]|_{N}$ has compact support, so the Jacquet integral for $W [f^\sharp [\chi ], \psi_T ] (1)$ converges absolutely for all $\chi$, without need for analytic continuation.  Thus
  \begin{equation*}
    |W[f^\sharp [\chi],\psi_T](1)|
    =
    \left\lvert
      \int_N f^\sharp [\chi](n) \psi_T^{-1}(n) \, d n
    \right\rvert
    \leq
    \int_{N} |f^\sharp [\chi](n)| \, d n.
  \end{equation*}
  By the proof of Lemma \ref{lem:standard:let-chi-in} (namely, \eqref{eq:f-sharp-chik}), we see moreover that
  \begin{equation*}
    f^\sharp [\chi](n) \ll  T^{\langle \rho_U^\vee , \Re(\chi) \rangle + o(1)}.
  \end{equation*}
  Combining this pointwise upper bound with the noted support condition yields the reduced form of the required estimate \eqref{eq:wfchi-psi_tg-ll} for $W[f^\sharp [\chi],\psi_T](1)$.
\end{proof}

\subsection{Whittaker function lower bounds}\label{sec:whitt-funct-lower}

We now develop the second of the two main exports of this section: a lower bound (Proposition~\ref{prop:there-exists-gamma}) for the Whittaker functions attached to suitable elements of the class $\mathfrak{F}(U \backslash G, \psi, T)$.  Like the upper bound, this will be applied in \S\ref{sec:constr-test-vect-proofs}, with the role of the present section's $G$ played by $H$.

\begin{notation}
  The following notation is active only for \S\ref{sec:whitt-funct-lower} (note that it differs from that employed in the proof of Proposition \ref{prop:standard:each-f-in}).  Let $g \in G$ (or $U \backslash G$).  If $g$ is of the form $u a n$ for some $(u,a,n) \in U \times A \times N$, then we define
  \begin{equation*}
    g_A := a,
    \quad
    g_N := n.
  \end{equation*}
  Otherwise, we leave $g_A$ and $g_N$ undefined.
\end{notation}

\begin{lemma}
  Let $T \ggg 1$.  Fix $\phi \in C_c^\infty(\mathfrak{g})$ with $\int \phi \neq 0$.  Let $t > 0$ with $t \lll T^{-1/2}$.  Introduce the temporary notation: for $n \in N$,
  \begin{equation*}
    I(n) :=
    \int_{x \in \mathfrak{g}}
    \phi(x) \delta_U^{1/2}((n \exp(t x))_A)
    \psi_T ( (n \exp(t x))_N)
    \exp (- \langle T \theta_G(\psi), t x \rangle) \, d x.
  \end{equation*}
  Then
  \begin{equation*}
    \int_{n \in N}
    \left\lvert
      I(n)
    \right\rvert^2
    \, d n
    \gg T^{-\dim(N)/2}.
  \end{equation*}
\end{lemma}
\begin{proof}
  It will suffice to show that $I(n) \gg 1$ whenever $n = 1 + \O(T^{-1/2})$.  For such $n$, we may write $n = \exp(y)$ with $y \in \Lie(N)$.  Using Baker--Campbell--Hausdorff, we see that for $x \in \supp(\phi)$,
  \begin{equation*}
    (n \exp(tx))_N = \exp(y + t x_N + \O(T^{-1/2} t)),
  \end{equation*}
  where $x_N \in \Lie(N)$ denotes the component of $x$ with respect to the decomposition $\Lie(G) = \Lie(N) \oplus \Lie(Q)$.  Since $T^{1/2} t \lll 1$, it follows that
  \begin{equation*}
    \psi_T((n \exp(t x))_N) = \psi_T(n) \exp (\langle T \theta_G(\psi), t x \rangle) + o(1).
  \end{equation*}
  We derive similarly that
  \begin{equation*}
    \delta_U^{1/2}((n \exp(t x))_A) = 1 + o(T^{-1/2} t) = 1 + o(T^{-1}) = 1 + o(1).
  \end{equation*}
  Thus the integrand in the definition of $I(n)$ is $\psi_T(n) \phi(x) + o(|\phi(x)|)$ on its support.  Integrating over $x$, we obtain
  \begin{equation*}
    I(n) = \psi_T(n) \int_{\mathfrak{g} } \phi(x) \, d x + o( \int_{x \in \mathfrak{g} } |\phi(x)| \, d x)
    \asymp 1.
  \end{equation*}
\end{proof}

\begin{corollary}\label{cor:sub-gln:there-exists-fixed}
  There exists a fixed $\phi \in C_c^\infty(\mathfrak{g})$ with the following property.  Set
  \begin{align*}
    &I(n) \\
    &:=
      \int_{x \in \mathfrak{g}}
      \phi(x) \delta_U^{1/2} \left( (n \exp(\tfrac{x}{T^{1/2}}))_A  \right)
      \psi_T ( (n \exp(\tfrac{x}{T^{1/2}}))_N )
      \exp (- \langle T \theta_G(\psi), \tfrac{x}{T^{1/2}} \rangle) \, d x.
  \end{align*}
  Then
  \begin{equation*}
    \int_{n \in N}
    \left\lvert
      I(n)
    \right\rvert^2
    \, d n
    \gg T^{-\dim(N)/2}.
  \end{equation*}
\end{corollary}
\begin{proof}
  Call the LHS of the desired estimate $I(\phi)$.  Fix $\phi_1 \in C_c^\infty(\mathfrak{g})$ with $\int \phi_1 \neq 0$.  For $\eps> 0$, set $\phi_\eps(x) := \eps^{-\dim(\mathfrak{g})} \phi_1(x/\eps)$.  The previous lemma implies that whenever $\eps \lll 1$, we have $I(\phi_\eps) \gg T^{-\dim(N)/2}$.  (Indeed, after the change of variables \(x\mapsto \varepsilon x\), we see that \(I(n)\) is as in Lemma 10.24 with \(t=\varepsilon/T^{1/2}\) and $\phi = \phi_1$.)  We deduce via overspill that the same estimate holds for some fixed $\eps > 0$.  (Indeed, the set $\{\eps > 0 : I(\phi_{\eps}) \geq \eps T^{-\dim(N)/2} \}$ contains all positive reals $\eps \lll 1$, hence must contain some fixed positive real.)  We conclude by taking $\phi = \phi_\eps$.
\end{proof}

For the following proposition and its proof, we denote by
\begin{equation*}
  \mathbf{1} := |.|^0 \in \mathfrak{X}(A)
\end{equation*}
the trivial character.  Thus, for instance,
\begin{equation*}
  \mathcal{J}_T[\psi,\beta,\gamma][\mathbf{1}] =
  \mathcal{J}_T[\psi,\beta,\gamma,\mathbf{1}]
\end{equation*}
denotes the Mellin component at the trivial character of the elements $\mathcal{J}_T[\psi,\beta,\gamma] \in \mathcal{S}^e(U \backslash G)$ constructed previously.
\begin{proposition}\label{prop:there-exists-gamma}
  There exists $\gamma \in \mathfrak{C}(G,-\theta_G(\psi),T,\tfrac{1}{2})$ so that for each $\beta \in \mathfrak{B}(A)$ with $\int_{A} \beta \asymp 1$, the element
  \begin{equation}\label{eq:f-:=-mathcalj_tpsi}
    f := \mathcal{J}_T[\psi,\beta,\gamma] = L (T^{\rho_U^\vee }) \mathcal{J}[\psi_T,\beta,\gamma] \in \mathfrak{F}(U \backslash G, \psi, T)
  \end{equation}
  has the following property.  There is a fixed $C_0 > 0$ so that for all $g \in Q$ for which $g = 1 + \O(T^{-1/2-\eps})$ for some fixed $\eps > 0$, we have
  \begin{equation}\label{eq:real-part-W-lower-bound}
    \Re(W[f[\mathbf{1}], \psi_T](g)) \geq C_0 T^{-\dim(N)/2}.
  \end{equation}
\end{proposition}
\begin{proof}
  We abbreviate $\theta := \theta_G(\psi)$.  Let $\phi$ be as in the conclusion of Corollary \ref{cor:sub-gln:there-exists-fixed}.  Let $\gamma_0 \in \mathfrak{C}(G, - \theta, T, \tfrac{1}{2})$ be constructed from $\phi$ as in Example \ref{example:cump8dyq7s}.  By Proposition \ref{lem:standard:subcl-mathfr-theta}, we have
  \begin{equation*}
    \gamma := \gamma_0^* \ast \gamma_0 \in \mathfrak{C}(G,-\theta,T,\tfrac{1}{2}).
  \end{equation*}
  We claim that this element $\gamma$ satisfies the required conclusions.

  Indeed, let $\beta$ be as in our hypotheses, and define $f$ as in \eqref{eq:f-:=-mathcalj_tpsi}.  Without loss of generality, suppose that $\int_{A} \beta = 1$.  Then by Lemma \ref{lem:scratch-research:let-beta_1-beta_2} (applied with ``$\psi$'' the rescaled character $\psi_T$) and the identity
  \begin{align*}
    \mathcal{J}_T[\psi,\beta,\gamma_0,\mathbf{1}](h)
    &=
      \delta_U^{-1/2} (T^{\rho_U^\vee}) \mathcal{J}[\psi_T,\beta,\gamma_0, \mathbf{1}](T^{\rho_U^\vee} h)
    \\
    &=
      \mathcal{J}[\psi_T,\beta,\gamma_0,\mathbf{1}](h),
  \end{align*}
  we have
  \begin{equation}\label{eq:whittaker-evaluate-at-1-via-integral-distributions-convolved-appl}
    W[f[\mathbf{1}], \psi_T](g) = \int_N \mathcal{J}[\psi_T, \beta, \gamma_0,\mathbf{1}] \overline{\mathcal{J}[\psi_T, \beta, \gamma_0 \ast g^{-1},\mathbf{1}]}.
  \end{equation}
  We must bound the real part of this from below.

  Consider first the case of the identity element $g = 1$.  An important feature in this case is that the integrand in \eqref{eq:whittaker-evaluate-at-1-via-integral-distributions-convolved-appl} is nonnegative.  For $n \in N$, we have
  \begin{align*}
    \mathcal{J}[\psi_T,\beta,\gamma_0,\mathbf{1}](n)
    &=
      \int_{a \in A} \delta_U^{-1/2}(a) \mathcal{J}[\psi_T, \beta, \gamma_0](a n) \, d a \\
    &=
      \int_{a \in A} \delta_U^{-1/2}(a)
      \int_{h \in G }  \gamma_0(h) \mathcal{J}[\psi_T,\beta](a n h) \, d a.\
  \end{align*}
  By definition, $\mathcal{J}[\psi_T,\beta](g) =  \delta_U^{1/2}(g_A) \beta(g_A) \psi_T(g_N)$.  We have
  \begin{equation*}
    (a n h)_A = a (n h)_A, \quad (a n h)_N = (n h)_N,
  \end{equation*}
  so the above simplifies to
  \begin{equation*}
    \int_{a \in A}
    \int_{h \in G} \gamma_0(h)
    \delta_U^{1/2}((n h)_A)
    \beta(a (n h)_A)
    \psi_T((n h)_N) \, d a.
  \end{equation*}
  This double integral converges absolutely.  Indeed, the integrand is compactly-supported in $h$ due to the factor $\gamma_0(h)$, then compactly-supported in $a$ due to the factor $\beta(a(n h)_A)$.  After interchanging the order of integration and substituting $a \mapsto a (n h)_A^{-1}$, the remaining $a$-integral evaluates to $\int_A \beta  = 1$, and we are left with
  \begin{equation*}
    \int_{h \in G} \gamma_0(h)
    \delta_U^{1/2}((n h)_A)
    \psi_T((n h)_N).
  \end{equation*}
  Appealing to the definition \eqref{eq:int-_g-gamma_0} of $\gamma_0$, we conclude that
  \begin{align*}
    &\mathcal{J}[\psi_T,\beta,\gamma_0, \mathbf{1}](n)   \\
    &\quad
      =
      \int_{x \in \mathfrak{g} }
      \phi (x) \delta_U^{1/2}((n \exp(\tfrac{x}{T^{1/2}}))_A) \psi_T ((n \exp(\tfrac{x}{T^{1/2} }))_N)
      \psi ( - \langle T \theta, \tfrac{x}{ T^{1/2} } \rangle) \, d x.
  \end{align*}
  By the conclusion of Corollary \ref{cor:sub-gln:there-exists-fixed}, it follows that
  \begin{equation*}
    \int_{n \in N} \left\lvert \mathcal{J}[\psi_T,\beta,\gamma_0, \mathbf{1}](n) \right\rvert^2 \, d n \gg T^{-\dim(N)/2}.
  \end{equation*}
  Via \eqref{eq:whittaker-evaluate-at-1-via-integral-distributions-convolved-appl}, we obtain the required lower bound \eqref{eq:real-part-W-lower-bound} in the case $g = 1$.

  It remains to treat the more general case in which $\dist(g,1) \ll T^{-1/2-\eps}$.  We will deduce this case from the $g = 1$ case via ``continuity,'' as follows.  Set
  \begin{equation*}
    v_g := \mathcal{J}[\psi_T,\beta,\gamma_0 \ast g^{-1}, \mathbf{1}] = \mathcal{J}_T[\psi,\beta,\gamma_0 \ast g^{-1}, \mathbf{1}].
  \end{equation*}
  Let $\langle ,  \rangle$ denote the inner product on the unitary representation $\mathcal{I}(1)$ defined by integration over $N$, and $\|v\| := \langle v, v \rangle^{1/2}$ the associated norm.  We may rewrite \eqref{eq:whittaker-evaluate-at-1-via-integral-distributions-convolved-appl} as
  \begin{equation*}
    W[f[\mathbf{1}], \psi_T](g) = \langle v_1, v_g \rangle.
  \end{equation*}
  We have seen already that $\langle v_1, v_1 \rangle \gg T^{-\dim(N)/2}$.  To obtain the required estimate, it is enough to check that $\langle v_1, v_g - v_1 \rangle \lll T^{-\dim(N)/2}$ whenever $\dist(g,1) \ll T^{-1/2-\eps}$.  By Cauchy--Schwarz, we reduce further to checking that
  \begin{equation}\label{eq:v_12-ll-t}
    \|v_1\|^2 \ll T^{-\dim(N)/2 + o(1)},
  \end{equation}
  \begin{equation}\label{eq:v_g-v_12-ll}
    \|v_g - v_1\|^2 \ll T^{-\dim(N)/2 - 2 \eps + o(1)}.
  \end{equation}

  The first estimate \eqref{eq:v_12-ll-t} is a consequence of Lemma \ref{lem:scratch-research:let-f-in}, and in particular the $L^2$ estimate \eqref{eq:int_k-fchi2-ll}, applied to $\mathcal{J}_T[\psi,\beta,\gamma_0] \in \mathfrak{F}(U \backslash G, \psi, T)$.

  For the second estimate, we recall (from \S\ref{sec:analysis-U-G-quantitative-setup}) that $\theta = \theta_G(\psi)$ was defined to have trivial restriction to $\Lie(Q)$.  Since $g \in Q$, it follows that $\psi_{- T \theta, \log}(g) = 1$, so by Proposition \ref{lem:sub-gln:let-gamma-in}, we obtain
  \begin{equation*}
    \gamma_0 \ast g^{-1} - \gamma_0 \in T^{1/2} \dist_G(g,1) \mathfrak{C}(G,-\theta,T,\tfrac{1}{2})
    \subseteq
    T^{-\eps} \mathfrak{C}(G,-\theta,T,\tfrac{1}{2}).
  \end{equation*}
  It follows that $v_g - v_1 = f_g[\mathbf{1}]$, where
  \begin{equation*}
    f_g := \mathcal{J}_T[\psi, \beta, \gamma_0 \ast g^{-1} - \gamma_0] \in T^{-\eps} \mathfrak{F}(U \backslash G, \psi, T).
  \end{equation*}
  We deduce \eqref{eq:v_g-v_12-ll} once again via Lemma \ref{lem:scratch-research:let-f-in}.
\end{proof}

\section{Construction of test vectors: proofs}\label{sec:constr-test-vect-proofs}

\subsection{Setup}\label{sec:constr-test-vec-setup}
We retain the setting of Theorem \ref{thm:main-local-results}: $(G,H) = (\GL_{n+1}(\mathbb{R}),\GL_n(\mathbb{R}))$ for some fixed $n$ (with the usual accompanying notation), $\psi$ is a fixed nondegenerate unitary character of $N$, $T \ggg 1$ is a large positive real, and $\pi$ is a generic irreducible unitary representation of $G$ with archimedean $L$-function parameters of magnitude $\asymp T$.  We will construct the required vectors and verify all claimed properties.

We continue to denote also by $\psi$ its restriction to $N_H$.  We continue to denote by $P$ the mirabolic subgroup of $G$, consisting of elements with bottom row $(0,\dotsc,0,1)$.  We define $\theta_P(\psi)$ as in \S\ref{sec:parameter-theta} and $\theta_H(\psi^{-1}) = - \theta_H(\psi)$ as in \S\ref{sec:analysis-U-G-quantitative-setup}.  We define the rescaled character $\psi_T$ as in \S\ref{sec:class-bump-functions}.  We define the classes $\mathfrak{W}(\pi,\psi,T,\delta)$ (for fixed $\delta \in (0,1/2)$) as in Definition \ref{defn:we-denote-mathfr-frakW}, $\mathfrak{F}(U_H \backslash H, \psi^{-1}, T)$ as in Definition \ref{defn:standard2:we-denote-mathfrakfu}, and $\mathfrak{C}(\dotsb)$ as in \S\ref{sec:spec-group-sett}.  We set $\tau := \tau(\pi,\psi,T)$ as in Definition \ref{defn:parameter-tau-pi-T}.  Then, as noted in Lemma \ref{lem:standard2:no-conductor-drop-vs-tau}, $\tau$ lies in a fixed compact subset of $\mathfrak{g}^\wedge_{\stab}$.

\subsection{Estimates for rescaled local zeta integrals}\label{sec:estim-resc-local}
Let $f \in \mathcal{S}(U_H \backslash H)$.  For $s \in \mathfrak{a}_{H,\mathbb{C}}^*$, we define the Mellin components $f[s] \in \mathcal{I}(s)$ as in \S\ref{sec:mellin-paley-affine-quotients}, following the ``unramified character'' notation (i.e., $f[s] = f[|.|^s]$, $\mathcal{I}(s) = \mathcal{I}(|.|^s)$).  In particular, $\mathcal{I}(0)$ denotes the induction of the trivial character, and $f[0]$ the corresponding component of $f$; these are the only components relevant for the following discussion.  (By contrast, in the statement of Proposition \ref{prop:there-exists-gamma}, we referred to the trivial character as $\mathbf{1}$.  This minor adjustment of notation should introduce no confusion.)

Let $W \in \mathcal{W}(\pi,\psi_T)$.  As in \S\ref{sec:zeta-integrals}, we may define, for $s \in \mathfrak{a}_{H,\mathbb{C}}^*$, the local zeta integral
\begin{equation*}
  Z_T(W,f[s])
  :=
  \int_{h \in N_H \backslash H}
  W(h) \cdot W[f[s], \psi_T^{-1}](h) \, d h.
\end{equation*}
The subscripted $T$ reflects that we are working with $\psi_T$-Whittaker models rather than $\psi$-Whittaker models, as in \S\ref{sec:zeta-integrals}.  Normally, this definition requires some interpretation via meromorphic continuation.  However, if the restriction of $W$ to $Q_H$ is compactly-supported, then the above formula converges absolutely for all $s$.  Indeed, it is then the integral of a continuous function over a compact set.  In such cases, there is no need to appeal to meromorphic continuation: the same formula works for all $s$, and the zeta integral is entire in $s$.  This happens if
\begin{equation}\label{eq:cuiqe676f7}
  W = \pi(\gamma) \Theta_{\pi,\psi_T} \qquad \text{ for some } \gamma \in \mathcal{M}_c(H),
\end{equation}
with the distribution $\Theta_{\pi, \psi_T}$ as defined in \S\ref{sec:kirillov-model}.  For instance, \eqref{eq:cuiqe676f7} holds if $W$ lies in the class $\mathfrak{W}(\pi,\psi,T,\delta)$, which we recall (see Definition \ref{defn:we-denote-mathfr-frakW}) consists of $W$ given by \eqref{eq:cuiqe676f7} for some $\gamma \in \mathfrak{C}(H,-\theta_H(\psi),T,\delta) \subseteq \mathcal{M}_c(H)$.  This is the only case we need to consider here.

In the definition of $Z_T$, we can replace the integral over $N_H \backslash H$ with one over $Q_H$ with respect to a right Haar measure.

We begin with the main calculation:

\begin{lemma}\label{lem:there-exists-begin}
  There exist
  \begin{equation*}
    \gamma \in \mathfrak{C}(H, -\theta_H(\psi^{-1}), T, \tfrac{1}{2}),
    \quad
    W \in \mathcal{W}(\pi,\psi_T),
    \quad
    \omega_0 \in C_c^\infty(G)
  \end{equation*}
  with the following properties.  Let $\beta \in \mathfrak{B}(A_H)^{W_H}$ with $\int_{A_H} \beta  = 1$.  Define
  \begin{equation*}
    f := \mathcal{J}_T[\psi^{-1},\beta,\gamma] \in \mathfrak{F}(U_H \backslash H, T, \psi^{-1}) \subseteq \mathcal{S}(U_H \backslash H).
  \end{equation*}
  Define $\omega$ and $\omega^\sharp$ as in the statement of Theorem \ref{thm:big-group-whittaker-behavior-under-G}.  Then:
  \begin{enumerate}[(i)]
  \item \label{item:lemma-there-exists-gamma-W-1} $\int_{N_H \backslash H} |W|^2 \ll T^{\dim(Q_H) / 2 + o(1)}$.
  \item \label{itm:lemma-exists-omega-2}  $\omega(g) \neq 0$ only if $g = 1 + o(1)$ and $\Ad^*(g) \tau = \tau + o(T^{-1/2})$.
  \item \label{itm:lemma-exists-omega-3} $\|\omega^\sharp\|_{\infty} \ll T^{\dim(N)+o(1)}$
  \item $\int_H |\omega^\sharp| \ll T^{\rank(H)/2 + o(1)}$.
  \item \label{item:lemma-there-exists-gamma-W-5} For each fixed $x \in \mathfrak{U}(G)$, we have $\|\pi(x) (\pi(\omega)W - W)\| \ll T^{-\infty}$.
  \item \label{item:lemma-there-exists-gamma-W-6}  $Z_T(W,f[0]) \gg T^{-\dim(N_H)/2}$.
  \item \label{item:lemma-there-exists-gamma-W-7} For all $s \in \mathfrak{a}_{H,\mathbb{C}}^*$ with $\Re(s) \ll 1$, we have
    \begin{equation*}
      Z_T(W,f[s]) \ll T^{\langle \rho_{U_H}^\vee, s \rangle - \dim(N_H)/2 + o(1)} \langle s \rangle^{-\infty}.
    \end{equation*}
  \end{enumerate}
\end{lemma}
\begin{proof}
  The proof begins with a few technical arguments that consolidate earlier results and perform some ``$\eps$ clean-up.''  Our main task is then to estimate the local zeta integrals $Z_T(W,f[s])$.

  We let $\gamma$ be as given by Proposition \ref{prop:there-exists-gamma}, and let $\beta$ and $f$ be as in the statement above.  The conclusion of that proposition provides a fixed $C_0 > 0$ so that for each fixed $\eps_0 > 0$,
  \begin{equation}\label{eq:rewf0-psi_tg-gg}
    g \in Q_H, \quad |g - 1| \leq T^{-1/2-\eps_0}
    \quad
    \implies
    \quad
    \Re(W[f[0], \psi_T^{-1}](g)) \geq C_0 T^{-\dim(N_H)/2}.
  \end{equation}
  By overspill, we deduce that there exists $\eps_0 > 0$ with $\eps_0 \lll 1$ for which the implication \eqref{eq:rewf0-psi_tg-gg} remains valid.  We henceforth work with such a value of $\eps_0$.

  Fix a basis $\mathcal{B}$ of $\mathfrak{g}$.  We observe the following consequence of Lemma \ref{lem:there-exists-w-big-group-big} and Theorem \ref{thm:big-group-whittaker-behavior-under-G}: for each fixed $\eps > 0$, there is a fixed $\delta \in (0,1/2)$ and elements $W \in \mathfrak{W}(\pi,\psi,T,\delta)$ and $\omega \in C_c^\infty(G)$ so that
  \begin{equation}\label{eq:int_n_h-backslash-h-1}
    \int_{N_H \backslash H} |W|^2 \leq   T^{\dim(Q_H)/2 + \eps},
  \end{equation}
  \begin{equation}\label{eq:w_q_h-geq-0}
    W|_{Q_H} \geq 0, \quad \int_{Q_H} W =  1,
  \end{equation}
  \begin{equation}\label{eq:g-in-q_h}
    g \in Q_H, W(g) \neq 0 \implies |g - 1| \leq \eps T^{-1/2-\eps_0},
  \end{equation}
  \begin{equation}\label{eq:g-in-g}
    g \in G, \omega(g) \neq 0 \implies |g - 1| \leq \eps, |\Ad^*(g) \tau - \tau| \leq \eps T^{-1/2},
  \end{equation}
  \begin{equation}\label{eq:omeg-leq-tnn+12}
    \|\omega^{\sharp}\|_{\infty} \leq T^{\dim(N) + \eps}, \quad
    \int_H \lvert \omega^\sharp \rvert \leq T^{\rank(H)/2+\eps},
  \end{equation}
  and, for all $m \in \mathbb{Z}_{\geq 0}$ with $m \leq 1/\eps$ and all $x_1,\dotsc,x_m \in \mathcal{B}$,
  \begin{equation}\label{eq:pix_1-dotsb-x_m}
    \|\pi(x_1 \dotsb x_m) (\pi(\omega_0) W - W)\| \leq T^{-1/\eps}.
  \end{equation}
  Indeed, given such an $\eps$, we fix $\delta \in (0,1/2)$, with $1/2 - \delta$ small enough in terms of $\eps$, and choose $W \in \mathfrak{W}(\pi,\psi,T,\delta)$ according to Lemma \ref{lem:there-exists-w-big-group-big}.  By Lemma \ref{lem:fix-delta-in} and the definition of $\mathfrak{C}(Q_H,0,T,\delta)$, we have
  \begin{equation*}
    W|_{Q_H}(g) \neq 0 \implies |g - 1| \ll T^{-1+\delta_+}
  \end{equation*}
  for some fixed $\delta_+ \in (\delta,1/2)$.  Since $\eps_0 \lll 1$ and $T \ggg 1$, we have $T^{-1+\delta_+} \lll T^{-1/2-\eps_0}$, so the support condition \eqref{eq:g-in-q_h} holds.  The $L^2$-estimate \eqref{eq:int_n_h-backslash-h-1} follows from Lemma \ref{lem:let-w-in-square-integral}, while the assertions in \eqref{eq:w_q_h-geq-0} are the content of the conclusion of Lemma \ref{lem:there-exists-w-big-group-big}.  We then take $\omega_0$ as in Theorem \ref{thm:big-group-whittaker-behavior-under-G} (applied with $\eps/2$ in place of $\eps$, so that we may absorb any implied constants $\ll 1$ into the factors $T^{\eps/2} \ggg 1$).  The remaining assertions \eqref{eq:g-in-g}, \eqref{eq:omeg-leq-tnn+12} and \eqref{eq:pix_1-dotsb-x_m} are then the content of the conclusion of Theorem \ref{thm:big-group-whittaker-behavior-under-G}.

  By overspill, we may find $\eps > \eps_0$ with $\eps \lll 1$ for which there exist $W \in \mathcal{W}(\pi,\psi_T)$ and $\omega_0 \in C_c^\infty(G)$ satisfying the stated assertions \eqref{eq:int_n_h-backslash-h-1}--\eqref{eq:pix_1-dotsb-x_m} .  We henceforth work such values of $\eps, W$ and $\omega_0$.  Assertions \eqref{item:lemma-there-exists-gamma-W-1}--\eqref{item:lemma-there-exists-gamma-W-5} then follow from \eqref{eq:int_n_h-backslash-h-1}, \eqref{eq:g-in-g}, \eqref{eq:omeg-leq-tnn+12} and \eqref{eq:pix_1-dotsb-x_m}.  It remains to verify the stated bounds \eqref{item:lemma-there-exists-gamma-W-6} and \eqref{item:lemma-there-exists-gamma-W-7} for $Z_T$.

  We observe by the support condition \eqref{eq:g-in-q_h} for $W|_{Q_H}$, together with the conclusion of Proposition \ref{prop:there-exists-gamma}, that the lower bound \eqref{eq:rewf0-psi_tg-gg} for the real part of $W[f[0],\psi_T^{-1}]$ holds on the support of $W|_{Q_H}$.  By \eqref{eq:w_q_h-geq-0}, it follows that
  \begin{equation*}
    \Re(Z_T(W,f[0]))
    =
    \int_{Q_H}
    W \cdot \Re(W[f[0],\psi_T^{-1}])
    \geq C_0 T^{-\dim(N_H) / 2}.
  \end{equation*}
  Thus assertion \eqref{item:lemma-there-exists-gamma-W-6} holds.

  Let $g \in \supp(W|_{Q_H})$.  By \eqref{eq:g-in-q_h}, we have $g = 1 + o(T^{-1/2})$.  In that case, Proposition \ref{lem:all-f-in} gives that for $\Re(s) \ll 1$,
  \begin{equation}
    W[f[s], \psi_T^{-1}](g) \ll T^{\langle \rho_{U_H}^\vee, s \rangle - \dim(N_H)/2 + o(1)} \langle s \rangle^{-\infty }.
  \end{equation}
  Assertion \eqref{item:lemma-there-exists-gamma-W-7} then follows by integrating over $Q_H$ against $W$ and using that $\int_{Q_H} \lvert W \rvert \ll 1$ (by \eqref{eq:w_q_h-geq-0}).
\end{proof}

\subsection{Transversality}\label{sec:cnpv1hjjro}
In this section, we verify that the $\omega_0 \in C_c^\infty(G)$ furnished by Lemma \ref{lem:there-exists-begin} satisfies the ``transversality'' or ``bilinear forms'' assertion \eqref{itm:sub-gln:3} of Theorem \ref{thm:main-local-results} for a wide class of $\phi_{Z_H} \in C_c^\infty(Z_H)$.

We continue to write $\bar{G} = G / Z$, where $Z < G$ and $Z_H \leq H$ denote the centers.  Recall from \S\ref{sec:dist-funct-d_h} the definition of $d_H : \bar{G} \rightarrow [0,1]$.
\begin{theorem}\label{thm:bilinear-forms-cited}
  Let $\mathfrak{S}$ be a compact subset of $\mathfrak{g}^\wedge_{\stab}$.  There are bounded symmetric open neighborhoods $\mathcal{Z} \subseteq Z_H$ and $\mathcal{G} \subseteq \bar{G}$ of the respective identity elements, which may be chosen arbitrarily small, for which the following assertions hold.  Let $\Omega$ be a compact subset of $H$.  There exists $C \geq 0$ with the following property.  Let $u_1 ^0 , u_2 ^0 : H \rightarrow \mathbb{R}_{\geq 0}$ be nonnegative measurable functions.  Define the convolutions
  \begin{equation*}
    u_i(x) := \int_{z \in \mathcal{Z} } u_i^0 (x z^{-1}) \, d z.
  \end{equation*}
  Let $r \geq 0$ and $\gamma \in \bar{G} - H$.  Assume that $\tau \in \mathfrak{S}$.  Define
  \begin{equation*}
    \mathcal{C}_\tau(r) = \left\{ g \in \bar{G} : D_\tau(g) \leq r \right\},
    \quad
    D_\tau(g) := \max \paren{ \abs{\Ad \paren{g} \tau-\tau}, \abs{\Ad \paren{g^{-1}} \tau, \tau} }.
  \end{equation*}
  Then the integral
  \begin{equation*}
    I := \int_{x, y \in \Omega } u_1(x) u_2(y) 1_{\mathcal{G} \cap \mathcal{C}_\tau(r)}  (x^{-1} \gamma y) \,d x \, d y
  \end{equation*}
  satisfies the estimate
  \begin{equation}\label{eq:i-leq-c}
    I \leq C r^{\dim_F(H)} \min \left( 1, \frac{r}{ d_H(\gamma)} \right) \|u_1^0 \|_{L^2(\Omega \mathcal{Z} )} \|u_2^0 \|_{L^2(\Omega \mathcal{Z})}.
  \end{equation}
\end{theorem}
\begin{proof}
  This is \cite[Theorem 15.1]{2020arXiv201202187N}, specialized here to $(\GL_{n+1}(\mathbb{R}), \GL_n(\mathbb{R}))$.
\end{proof}
\begin{theorem}[Reformulation of Theorem \ref{thm:bilinear-forms-cited}]\label{thm:bilinear-forms-reformulated}
  Assume that $\tau$ belongs to a fixed compact subset of $\mathfrak{g}^\wedge_{\stab}$.  There are fixed bounded symmetric open neighborhoods $\mathcal{Z} \subseteq Z_H$ and $\mathcal{G} \subseteq \bar{G}$ of the respective identity elements, which may be chosen arbitrarily small, for which the following assertions hold.  Let $\Omega$ be a fixed compact subset of $H$.  Let $u_1 ^0 , u_2 ^0 : H \rightarrow \mathbb{R}_{\geq 0}$ be nonnegative measurable functions.  Let $r \geq 0$ and $\gamma \in \bar{G} - H$.  Then
  \begin{equation}\label{eq:i-ll-rdim_fh}
    I \ll r^{\dim(H)} \min \left( 1, \frac{r}{ d_H(\gamma)} \right) \|u_1^0 \|_{L^2(\Omega \mathcal{Z} )} \|u_2^0 \|_{L^2(\Omega \mathcal{Z})}.
  \end{equation}
\end{theorem}
\begin{proof}
  This is an exercise in applying the axioms of IST (\S\ref{sec:axioms-ist}).  By transfer, Theorem \ref{thm:bilinear-forms-cited} holds upon restricting quantification over the quantities $\mathfrak{S}$ through $C$ to be fixed.  By idealization, we can then commute quantification over $C$ with the remaining quantifiers.  We arrive in this way at the stated reformulation.
\end{proof}

\begin{proposition}\label{prop:standard2:let-tau-belong}
  Retain the setting of \S\ref{sec:constr-test-vec-setup} -- in particular, $\pi$, $T$ and $\tau = \tau(\pi,\psi,T)$ are defined.  Let $\omega_0$ be as furnished by Lemma \ref{lem:there-exists-begin}; define $\omega$ and $\omega^\sharp$ as in the statement of Theorem \ref{thm:big-group-whittaker-behavior-under-G}.  Suppose that $\tau$ belongs to a fixed compact subset of $\mathfrak{g}^\wedge_{\stab}$.  There is a fixed bounded neighborhood $\mathcal{Z}$ of the identity in $Z_H$ with the following property.  Let $\phi_{Z_H} \in C_c^\infty(\mathcal{Z})$ with $\|\phi_{Z_H}\|_{L^\infty} \ll 1$.  For each fixed compact $\Omega \subseteq H$, there is a fixed compact $\Omega ' \subseteq H$ with the following property.  Let $\Psi_1, \Psi_2 : H \rightarrow \mathbb{C}$ be continuous functions.  Then for each $\gamma \in \bar{G} - H$, the integral
  \begin{equation*}
    I := \int_{x, y \in \Omega }
    \left\lvert
      (\Psi_1 \ast \phi_{Z_H})(x)
      (\Psi_2 \ast \phi_{Z_H})(y)
      \omega^\sharp(x^{-1} \gamma y)
    \right\rvert
    \, d x \, d y
  \end{equation*}
  satisfies the estimate
  \begin{equation*}
    I \ll
    T^{\rank(H)/2 + o(1)}
    \min\left( 1,
      \frac{T^{-1/2}}{ d_H(\gamma)}
    \right)
    \|\Psi_1\|_{L^2(\Omega')}
    \|\Psi_2\|_{L^2(\Omega')}.
  \end{equation*}
\end{proposition}
\begin{proof}
  We apply Theorem \ref{thm:bilinear-forms-reformulated} to any fixed compact subset of $\mathfrak{g}^\wedge_{\stab}$ that contains $\tau$.  This application yields $\mathcal{Z}$ and $\mathcal{G}$.  Let $\phi_{Z_H}$ and $\Omega$ as in our hypotheses.  Take for $\Omega '$ any fixed compact set that contains $\Omega \mathcal{Z}$.  Let $\Psi_1, \Psi_2$ and $\gamma$ be as in our hypotheses.  We will verify that the stated estimate holds.

  To that end, we first observe that for all $g \in \bar{G}$
  \begin{equation*}
    \omega^\sharp(g) \ll
    1_{\mathcal{G} \cap \mathcal{C}_\tau(T^{-1/2})}     T^{\dim(N) + o(1)}.
  \end{equation*}
  This follows from parts of Lemma \ref{lem:there-exists-begin}.  Indeed, part \eqref{itm:lemma-exists-omega-2} implies that $\omega^\sharp(g)$ vanishes unless $g = 1 + o(1)$ and $\Ad^*(g) \tau = \tau + o(T^{-1/2})$, in which case $g \in \mathcal{G}$ and $D_\tau(g) \lll T^{-1/2}$, hence $g \in \mathcal{C}_\tau\paren{T^{-1/2}}$, while part \eqref{itm:lemma-exists-omega-3} says that $\|\omega^\sharp\|_{\infty} \ll T^{\dim(N)+o(1)}$.

  We now set $u_j^0 := \Psi_j$ and define $u_j$ as in Theorem \ref{thm:bilinear-forms-reformulated}.  Since $\lVert \phi_{Z_H} \rVert_{L^\infty} \ll 1$, we then have $\abs{ \paren{\Psi_j \ast \phi_{Z_H}}(x) } \ll u_j(x)$.  Moreover,
  \begin{equation*}
    \lVert u_j^0 \rVert_{L^2(\Omega \mathcal{Z})} \leq \lVert \Psi_j \rVert_{L^2(\Omega ')}.
  \end{equation*}
  We conclude via \eqref{eq:i-ll-rdim_fh} and the identity $\dim(N) - \dim(H)/2 = \rank(H)/2$.  (For example, if $(G,H) = (\GL_{n+1}(\mathbb{R}), \GL_n(\mathbb{R}))$, then this reads $n(n+1)/2 - n^2/2 = n/2$.)
\end{proof}

\subsection{Book-keeping}\label{sec:book-keeping}
Here we complete the proof of Theorem \ref{thm:main-local-results}.  Our main task is to collect and normalize estimates given previously.

We retain the setup of \S\ref{sec:constr-test-vec-setup}.

Let
\begin{equation*}
  W_T \in \mathfrak{W}(\pi,\psi,T,\delta) \subseteq \mathcal{W}(\pi,\psi_T),
\end{equation*}
\begin{equation*}
  \gamma \in \mathfrak{C}(H,-\theta_H(\psi^{-1}),T,\tfrac{1}{2})  \text{ and }  \omega_0 \in C_c^\infty(G)
\end{equation*}
be as furnished by Lemma \ref{lem:there-exists-begin}.

Let $\tau := \tau(\pi,T,\psi)$ as in \S\ref{sec:parameter-tau}.  Let $\mathcal{Z} \subseteq Z_H$ be the fixed bounded subset of the identity furnished by Proposition \ref{prop:standard2:let-tau-belong}.  We fix $\phi_{Z_H} \in C_c^\infty(\mathcal{Z}) \subseteq C_c^\infty(Z_H)$ with $\int_{Z_H} \phi_{Z_H} = 1$ with respect to the our fixed Haar measure on $Z_H$.  We fix $\beta_0 \in \mathfrak{B}(A_H)^{W_H}$ (\S\ref{sec:class-mathfrakba}) with $\int_{A_H} \beta_0 = 1$.  We define $\beta$ to be the convolution
\begin{equation}\label{eq:hx-:=-int}
  \beta (x) := \int_{z \in Z_H} \beta_0(x z^{-1}) \phi_{Z_H}(z) \, d z.
\end{equation}
Then $\beta$ likewise lies in $\mathfrak{B}(A_H)^{W_H}$ and satisfies $\int_{A_H} \beta = 1$.

We set
\begin{equation*}
  c_2 := T^{\dim(N_H)/4}
\end{equation*}
and
\begin{equation*}
  f := c_2 \mathcal{J}_T[\psi^{-1},\beta_0,\gamma],
  \quad
  f_{\ast} := c_2 \mathcal{J}_T[\psi^{-1},\beta,\gamma].
\end{equation*}

We define $a_T \in A \leq G$ by the formula
\begin{equation}\label{eq:a_t-=-trho_n_hvee}
  a_T := T^{\rho_{N_H}^\vee} z_H(T^{\rank(G)/2}),
  \quad
  z_H(t) := \diag(t,t,\dotsc,t,1).
\end{equation}
Explicitly, if $(G,H) = (\GL_{n+1}(\mathbb{R}), \GL_n(\mathbb{R}))$, then
$a_T = \diag(T^{n}, T^{n-1}, \dotsc, T, 1)$; the equivalence of this formula with \eqref{eq:a_t-=-trho_n_hvee} follows from the identity of $(n+1)$-tuples
\begin{align*}
  &\left( \frac{n-1}{2},
    \frac{n-3}{2},
    \dotsc,
    \frac{1-n}{2}, 0
    \right) \\
  &\quad +
    \left( \frac{n + 1}{2}, \frac{n + 1}{2}, \dotsc,
    \frac{n + 1}{2}, 0\right)
  \\
  &\quad =
    \left( n, n - 1, \dotsc, 1, 0 \right).
\end{align*}

Set
\begin{equation*}
  W_1(g) := W_T(a_T^{-1} g).
\end{equation*}
Then $W_1 \in \mathcal{W}(\pi,\psi)$.  Set
\begin{equation*}
  e(H) := \langle \rho_{N_H}, \rho_{N_H}^\vee  \rangle.
\end{equation*}
Then $\delta_{N_H}(T^{\rho_N^\vee}) = T^{e(H)}$, so the change of variables $h \mapsto a_T h$ on $N_H \backslash H$ has Jacobian $T^{-e(H)}$, giving
\begin{equation}\label{eq:int-_n_h-backslash}
  \int_{N_H \backslash H} |W_1|^2
  =
  T^{-e(H)}
  \int_{N_H \backslash H} |W_T|^2.
\end{equation}
By \eqref{eq:int-_n_h-backslash} and assertion \eqref{item:lemma-there-exists-gamma-W-1} of Lemma \ref{lem:there-exists-begin}, there exists
\begin{equation*}
  c_1 \asymp  T^{e(H)/2 - \dim(Q_H)/4}.
\end{equation*}
so that $(\int_{N_H \backslash H} |W_1|^2)^{1/2} \leq c_1^{-1}$.  Thus
\begin{equation*}
  W := c_1 W_1 \in \mathcal{W}(\pi,\psi)
\end{equation*}
satisfies $\int_{N_H \backslash H} |W|^2 \leq 1$.

We claim that the datum $(f, f_{\ast}, \phi_{Z_H}, W, \omega_0)$ that we have constructed has the required properties.

We verify first that each element lies in the required space.
\begin{itemize}
\item That $f, f_{\ast}$ lie in $\mathcal{S}^e(U_H \backslash H)^{W_H}$ is a consequence of Proposition \ref{lem:sub-gln:each-psi-beta} and Lemma \ref{lem:each-f-in}.
\item It is clear by construction that $W$ is a smooth vector of norm $\leq 1$.
\item It is clear from \eqref{eq:g-in-g} that $\omega_0$ is supported in any fixed neighborhood of the identity in $G$.
\end{itemize}
We now check each numbered assertion from Theorem \ref{thm:main-local-results} in turn:
\begin{enumerate}[(i)]
\item The identity $f_{\ast}(g) = \int_{z \in Z_H} f(g z^{-1}) \phi_{Z_H}(z) \, d z$ is immediate from the construction, using the equivariance property \eqref{eq:lx-mathcald_psi-beta}.
\item The arguments below apply just as well to $f$ and to $f_{\ast}$.
\item We must check that for each fixed compact subset $\mathcal{D}$ of $\mathfrak{a}^*$ and fixed $\ell \in \mathbb{Z}_{\geq 0}$, we have the seminorm bound $\nu_{\mathcal{D},\ell,T}(f) \ll T^{o(1)}$.  Recalling from \S\ref{sec:local-sobolev-norms-prelims} the definitions of such seminorms, the required bound is a consequence of Lemma \ref{lem:scratch-research:let-f-in}, which gives, for fixed $x_1,\dotsc,x_m \in \Lie(H)$,
  \begin{align*}
    \|R(x_1 \dotsb x_m) f[s]\|
    &\asymp c_2 \left(\int_K \left\lvert R(x_1 \dotsb x_m) \mathcal{J}_T[\psi^{-1},\beta,\gamma][s] \right\rvert^2\right)^{1/2} \\
    &\ll c_2 T^{\langle \rho_U^\vee, \Re(s) \rangle - \dim(N)/4 + m +  o(1)} \langle s \rangle^{-\infty} \\
    &= T^{\langle \rho_U^\vee, \Re(s) \rangle + m + o(1)} \langle s \rangle^{-\infty}.
  \end{align*}
\item We must now verify the claimed upper and lower bounds for local zeta integrals.  It is straightforward but tedious to deduce those bounds from Lemma \ref{lem:there-exists-begin}, as follows.  The change of variables $g \mapsto a_T g$ on $N_H \backslash H$, with Jacobian $T^{-e(H)}$, gives
  \begin{equation*}
    Z(W_1, f[s]) = T^{-e(H)} \int_{g \in N_H \backslash H} W_1(a_T g) W[f[s], \psi^{-1}](a_T g) \, d g.
  \end{equation*}
  We aim now to compare $W[f[s], \psi^{-1}](a_T g)$ with $W[f[s], \psi_T^{-1}](g)$.
  Since $f[s] : H \rightarrow \mathbb{C}$ has central character given by $\trace(s)$, we have
  \begin{equation*}
    W[f[s], \psi^{-1}](a_T g) = T^{\rank(G) \trace(s)/2} W[f[s], \psi^{-1}](T^{\rho_{N_H}^\vee} g).
  \end{equation*}

  We pause to record a general identity relating Jacquet integrals for the character $\psi$ and its rescaling $\psi_T$ (valid more generally for any $f \in \mathcal{I}(s)$):
  \begin{equation}\label{eq:cunvjucq6r}
    W[f, \psi](T^{-\rho_U^\vee} g) = T^{\langle \rho_U, \rho_U^\vee  \rangle - \rho_U^\vee(s)} W[f, \psi_T](g).
  \end{equation}
  Indeed, by definition, the LHS is
  \[
    \int_{n \in N} f(n T^{-\rho_U^\vee} g) \psi^{-1}(n) \,  d n.
  \]
  Recall that $N$ is opposite to $U$.  We apply the change of variables $N \ni n \mapsto T^{-\rho_U^\vee} n T^{\rho_U^\vee}$, with Jacobian $T^{2 \langle \rho_U, \rho_U^\vee  \rangle}$, to transform the above integral to
  \[
    T^{2 \langle \rho_U, \rho_U^\vee  \rangle}
    \int_{n \in N} f(T^{-\rho_U^\vee} n g) \psi_T^{-1}(n) \,  d n.
  \]
  We conclude via the transformation property
  \begin{equation*}
    f(T^{-\rho_U^\vee} n g) =
    \delta_U^{1/2}(T^{- \rho_U^\vee})
    \lvert T^{- \rho_U^\vee} \rvert^s
    f(n g)
    =
    T^{-\langle \rho_U + s,  \rho_U^\vee \rangle } f(n g).
  \end{equation*}
  (In verifying this, we recall from \S\ref{sec:cnjeor2mwy} that $\delta_U^{1/2} = \lvert . \rvert^{\rho_U}$.)

  By \eqref{eq:cunvjucq6r} and the identity $\rho_{N_H}^\vee = - \rho_{U_H}^\vee$, we obtain
  \begin{equation*}
    W[f[s], \psi^{-1}](T^{\rho_{N_H}^\vee} g) = T^{e(H)/2 - \rho_{U_H}^\vee(s)} W[f[s], \psi_T^{-1}](g).
  \end{equation*}
  Using also $W_1(a_T g)=W_T(g)$, we obtain the desired comparison
  \begin{equation}\label{eq:z_1w_1-f-s}
    Z(W_1,f[s]) =
    T^{\rank(G) \trace(s)/2 - e(H)/2 - \rho_{U_H}^\vee(s)}
    Z_T(W_T, f[s]).
  \end{equation}
  By Lemma \ref{lem:there-exists-begin} and \eqref{eq:z_1w_1-f-s}, we obtain
  \begin{equation}\label{eq:z_1w-f-0}
    Z(W,f[0]) \gg Z(c_1 W_1, f, 0) \gg c_1 c_2 T^{-e(H)/2 -\dim(N_H)/2}
  \end{equation}
  and, for $\Re(s) \ll 1$,
  \begin{equation}\label{eq:z_1w-f-s}
    Z(W,f[s]) \ll c_1 c_2 T^{-e(H)/2  - \dim(N_H)/2 + \rank(G) \trace(s)/2 + o(1)} \langle s \rangle^{-\infty}.
  \end{equation}
  We calculate
  \begin{equation*}
    c_1 c_2 T^{-e(H)/2 - \dim(N_H)/2}
    \asymp
    T^{ - \dim(Q_H)/4 - \dim(N_H)/4}
    =
    T^{-\dim(H)/4}.
  \end{equation*}
  With this calculation, assertion \eqref{itm:sub-gln:8} of Theorem \ref{thm:main-local-results}, namely that
  \begin{equation*}
    Z(W,f[0]) \gg T^{- \dim(H)/4},
  \end{equation*}
  \begin{equation*}
    \Re(s) \ll 1 \implies Z(W,f[s]) \ll  T^{\rank(G) \trace(s)/2 - \dim(H)/4 + o(1) }
    (1 + |s|)^{-\infty },
  \end{equation*}
  follows from \eqref{eq:z_1w-f-0} and \eqref{eq:z_1w-f-s}.
\item The required estimate $\|\pi(u) ( \pi(\omega_0) W - W)\| \ll T^{-\infty}$ for fixed $u \in \mathfrak{U}(G)$ is immediate from the corresponding assertion of Lemma \ref{lem:there-exists-begin}.
\item The required estimate $\int_{H} |\omega^\sharp| \leq T^{\rank(H)/2 + o(1)}$ is likewise immediate from Lemma \ref{lem:there-exists-begin}.
\item We must check that for each compact subset $\Omega_H$ of $H$, there is a fixed compact subset $\Omega_H'$ of $H$ with a certain property.  We appeal to Proposition \ref{prop:standard2:let-tau-belong}.  To see that the latter is applicable, we observe that
  \begin{itemize}
  \item $\tau$ belongs to a fixed compact subset of $\mathfrak{g}^\wedge_{\stab}$, in view of Lemma \ref{lem:standard2:no-conductor-drop-vs-tau} and our assumption that the archimedean $L$-function parameters of $\pi$ have magnitude $\asymp T$, and
  \item $\phi_{Z_H} \in C_c^\infty(Z_H)$ is fixed, hence in particular  $\|\phi_{Z_H}\|_{L^\infty} \ll 1$.
  \end{itemize}
  Taking $\Omega_H'$ as in the conclusion of Proposition \ref{prop:standard2:let-tau-belong}, the required estimate is then precisely the content of that proposition.
\end{enumerate}

The proof of Theorem  \ref{thm:main-local-results} is now complete modulo the proof of Theorem \ref{thm:big-group-whittaker-behavior-under-G}, which is the subject of Part \ref{part:asympt-analys-kirill}.

\part{Asymptotics of the Kirillov model}\label{part:asympt-analys-kirill}
The purpose of this part is to supply the proof of Theorem \ref{thm:big-group-whittaker-behavior-under-G}.  The proof is given at the end, in \S\ref{sec:proof-theor-refthm:b}, following several preliminaries.

\section{Overview}\label{sec:overview-asymp-kirillov}
Here we summarize the main points of Part \ref{part:asympt-analys-kirill}, without technicalities.  We retain the notation of \S\ref{sec:asympt-analys-local}.  In particular, $G = \GL_n(\mathbb{R})$, and $P \leq G$ denotes the mirabolic subgroup.  We are given a Whittaker function $W$ that is compactly-supported in the Kirillov model and transforms nicely under $P$ (e.g., is a bump function with given support); we would like to understand how $G$ acts on $W$.  Jacquet's results, recalled in \S\ref{sec:kirillov-model}, tell us that $W$ is smooth under $G$.  We would like to understand this smoothness quantitatively.

For example, if we differentiate $W$ some finite number $m$ times with respect to basis elements of $\Lie(G)$, to what extent does this increase the norm of $W$?  It turns out that for the $W$ of interest to us (elements of the class $\mathfrak{W}(\pi,\psi,T,\delta)$), the norm goes up by roughly a factor of $T^m$.  A direct quantification of Jacquet's arguments (\S\ref{sec:empha-priori-bounds}) gives instead an upper bound of the form $T^{C_m}$ for some $C_m \geq 0$.  For different classes of $W$ than those considered here (e.g., those relevant for ``analytic newvector theory'' \cite{JN19a, 2020arXiv200109640J}), it is necessary to take $C_m$ much larger than $m$ (e.g., $C_m = n m$ in the newvector setting), so the conclusion is not completely formal.  The $p$-adic analogue of this phenomenon is that the depth of a generic irreducible representation can be (and often is) significantly smaller than the conductor.

We will see, in fact, that the $W$ of interest to us behave like approximate eigenfunctions under $\mathfrak{U}(G)$, with eigenvalue described by the rescaling $T \tau$ of the element $\tau = \tau(\pi,\psi,T) \in \mathfrak{g}^\wedge$ as in Theorem \ref{thm:big-group-whittaker-behavior-under-G}.  This feature is ultimately responsible for our proof that the function $\omega_0 \in C_c^\infty(G)$ that we construct approximately preserves $W$ in the required sense.

The general theme of understanding $G$-behavior from $P$-behavior is a recurring one in the theory of Whittaker functions for $\GL_n$, fundamental to the theory of newvectors and that of the Kirillov model itself, see for instance \cite{MR620708, MR3138844, MR3001803} and \cite{MR0425030, MR0579172, MR748505,MR1999922}.

The basic idea underlying our approach to the study of $W$ is as in Jacquet's work.  Let $E_{i j}$ ($1 \leq i, j \leq n$) denote the standard basis elements of $\Lie(G)$, regarded as degree one elements of $\mathfrak{U}(G)$.  The elements indexed by $i < n$ lie in the Lie algebra of the mirabolic subgroup $P$, so their action is very explicit; the problem is to understand how the ``bottom row'' acts.  On the other hand, we know that $\mathfrak{Z}(G)$ acts via scalars, described by the infinitesimal character $\lambda_\pi$.  We can try to combine what we know about the actions of $\Lie(P)$ and $\mathfrak{Z}(G)$ to say something about that of $\Lie(G)$.

The $n=2$ case is particularly simple, and has appeared in some form in the subconvexity literature (see \cite[Lemma 8.4]{venkatesh-2005}, \cite[\S3.2.5]{michel-2009}).  We recall that case in detail, since it motivates what we do below.  The action of the central element $E_{1 1} + E_{2 2}$ is described by the central character of $\pi$, so we focus on the action of $\Lie(\SL_2(\mathbb{R}))$.  The main point is then to understand the action of $E_{2 1}$ in terms of that of $E_{12}$ and $h := E_{1 1} - E_{2 2}$.  The Casimir element is given (up to normalization) by
\begin{equation}\label{eq:u-=-tfrac12}
  u = \tfrac{1}{2} h^2 + E_{12} E_{21} + E_{2 1} E_{1 2}
  =\tfrac{1}{2}  h^2 - h + 2 E_{12} E_{21}.
\end{equation}
In this case, we can identify $N_H \backslash H = H$ with $\GL_1(\mathbb{R})$ via the map $y \mapsto \begin{pmatrix}
  y & 0 \\
  0 & 1
\end{pmatrix}$.  The Kirillov model of $\pi$ is thus an embedding $\pi \hookrightarrow L^2(\GL_1(\mathbb{R}))$.  Suppose that we have defined this model using the character $\begin{pmatrix}
  1 & x \\
  0 & 1
\end{pmatrix} \mapsto e^{i x}$.  Then the action of $\Lie(P)$ in the Kirillov model is given by
\begin{equation*}
  \pi(E_{12}) W(y) = i y W(y),
  \quad
  \pi(h) W(y) = 2 y \partial_y W(y).
\end{equation*}
Suppose that the Casimir element $u$ acts on $\pi$ via the eigenvalue $\lambda$.  Then, using \eqref{eq:u-=-tfrac12}, we obtain an explicit formula for the action of $E_{2 1}$,
\begin{equation}\label{eq:pie_2-1-wy}
  \pi(E_{2 1}) W(y) =
  \frac{1}{ 2 i y} \left(
    \lambda - \frac{1}{2} (2 y \partial_y)^2 + 2 y \partial_y
  \right) W(y),
\end{equation}
which we can use to study $W$ further.  We remark that in this case, the $W$ of interest to us are essentially bump functions on $y = T + \O(T^{1/2})$, in representations with Casimir eigenvalue $\lambda$ of size roughly $T^2$.  We see from \eqref{eq:pie_2-1-wy} that applying $\pi(E_{21 })$ to such an element $W$ increases its norm by roughly a factor of $T$.

The simplifying feature of the $n=2$ case is that the operator $\pi(E_{12})$ is explicitly invertible.  For $n \geq 3$, it seems unlikely to us that the operators $\pi(E_{n j})$ may be described quite as explicitly as in \eqref{eq:pie_2-1-wy}, e.g., as differential operators with rational coefficients.  Instead, we will approximate the $\pi(E_{n j})$, or more precisely their action on the class of Whittaker functions of interest, by convolution operators on $P$ attached to well-behaved symbols on $\Lie(P)^\wedge$.  This is the content of Theorem \ref{thm:g-via-p} below.  By a ``quantum analogue'' of the proof of the uniqueness of $\tau = \tau(\pi, \psi, T)$ (Lemma \ref{lem:construction-tau-from-theta}), we will eventually deduce the required $G$-behavior of $W$ from the known $P$-behavior.  Jacquet's results play an important role in getting this argument started: they tell us that the $G$-behavior of $W$ is ``polynomial in $T$'' \S\ref{sec:empha-priori-bounds}.

After understanding the behavior of $W$ under $\mathfrak{U}(G)$, we will want to construct the approximate idempotent $\omega_0 \in C_c^\infty(G)$ for $W$ required by Theorem \ref{thm:big-group-whittaker-behavior-under-G}.  We do this using convolution operators on $G$ attached to carefully-chosen symbols on $\Lie(G)^\wedge$.  We show first, using the approximate eigenvector property of $W$ under $\mathfrak{U}(G)$, that for any symbol taking the value $1$ on $T \tau + \O(T^{1/2})$, the corresponding convolution operator nearly preserves $W$.  We show next, using the action of the center of the universal enveloping algebra, that such a symbol may be truncated with negligible effect to within roughly $\O(1)$ of the locus of the infinitesimal character $\lambda_\pi$; this is the point of  \S\ref{sec:local-refin-symb}.  (This locus plays the role of the coadjoint orbit in our setup -- recall that $\pi$ is not assumed tempered.)  The resulting truncation then yields the desired $\omega_0$.

\section{Microlocal analysis on Lie group representations}
Here we recall some of the results of \cite[Parts 1 and 2]{nelson-venkatesh-1} and \cite[\S8-\S10]{2020arXiv201202187N}.  This review is non-exhaustive; we focus on what is needed for the applications of this paper, adding a couple minor new points along the way.  In particular, we do not recall any results related to the Kirillov character formula or trace norm estimates, as those do not play a role in this paper.

\subsection{Preliminaries}
Fix a Lie group $G$.  We retain the general notation of \S\ref{sec:parameters-stability}.  In particular, $\mathfrak{g}$ denotes the Lie algebra and $\mathfrak{g}^\wedge = i \mathfrak{g}^* = \Hom(\mathfrak{g},i \mathbb{R})$ its Pontryagin (equivalently, imaginary) dual.  We denote typical elements of $\mathfrak{g}$ (resp.\ $\mathfrak{g}^\wedge$) by $x,y,z$ (resp.\ $\xi,\eta,\zeta,\tau$).  We write $\langle x, \xi \rangle \in i \mathbb{R}$ for the result of the canonical pairing.  In this section, we often abbreviate
\begin{equation*}
  x \xi := \langle x, \xi \rangle,
\end{equation*}
so that the $\U(1)$-valued pairing between $\mathfrak{g}$ and $\mathfrak{g}^\wedge$ is given by $(x,\xi) \mapsto e^{ x \xi}$.

The group $G$ acts on $\mathfrak{g}$ via the adjoint representation $x \mapsto \Ad(g) x$ and on $\mathfrak{g}^\wedge$ via the coadjoint representation $\xi \mapsto \Ad^*(g) \xi$.  We often abbreviate $g \xi = \Ad^*(g) \xi$.  The Lie algebra $\mathfrak{g}$ acts on $\mathfrak{g}$ via the adjoint representation $\ad_x y = [x,y]$ and on $\mathfrak{g}^\wedge$ via the coadjoint representation $\ad_x^* \xi$.  We often abbreviate the latter by $[x,\xi]$.

For a Lie subgroup $H$ of $G$, we obtain a natural inclusion $\mathfrak{h} \leq \mathfrak{g}$ and restriction map $\mathfrak{g}^\wedge \twoheadrightarrow \mathfrak{h}^\wedge$.  For $\xi \in \mathfrak{g}^\wedge$, we write $\xi_H \in \mathfrak{h}^\wedge$ for its restriction.

We fix Haar measures $d x$ on $\mathfrak{g}$ and $d \xi$ on $\mathfrak{g}^\wedge$ that are Fourier dual in the sense that for Schwartz functions $a$ on $\mathfrak{g}^\wedge$ and $\phi$ on $\mathfrak{g}$, the Fourier transforms given by \index{Lie algebra!Fourier transforms $a^\vee, \phi^\wedge$}
\begin{equation*}
  a^\vee(x) := \int_{\xi \in \mathfrak{g} ^\wedge }
  a(\xi) e^{- x \xi} \, d \xi,
  \quad
  \phi^\wedge(\xi)
  = \int_{x \in \mathfrak{g} }
  \phi(x)
  e^{x \xi} \, d x
\end{equation*}
define mutually inverse transforms.

\subsection{Wavelength parameters}
Let $\h \in (0,1]$.  We regard \index{wavelength parameter $\h$} $\h$ as a ``wavelength parameter.''  For any function $a$ on $\mathfrak{g}^\wedge$, we define its rescaling \index{symbols!rescaling $a_{\h}(\xi) := a(\h \xi)$}
\begin{equation*}
  a_{\h}(\xi) := a(\h \xi).
\end{equation*}

\subsection{Symbol classes}
The following definitions come from \cite[\S9]{2020arXiv201202187N}.  We simply recall the definitions here, referring to the cited reference for further discussion, intuition and motivation.

In what follows, we repeatedly use the notation
\begin{equation*}
  \langle \xi \rangle := (1 + |\xi|^2)^{1/2}.
\end{equation*}

\subsubsection{Underlying spaces}
For each integer $m \in \mathbb{Z}$ and multi-index $\alpha \in \mathbb{Z}_{\geq 0}^{\dim(\mathfrak{g})}$ (\S\ref{sec:multi-index-notation}), we define a seminorm $\nu_{m,\alpha}$ on $C^\infty(\mathfrak{g}^\wedge)$ by
\[\nu_{m,\alpha}(a) := \sup_{\xi \in \mathfrak{g}^\wedge}
  \langle \xi \rangle^{|\alpha| - m} \left\lvert \partial^\alpha a(\xi) \right\rvert,
\]
where
\begin{itemize}
\item $|\alpha| := \sum \alpha_j$ denotes the order of $\alpha$, and
\item $\partial^\alpha$ denotes the differential operator defined using some fixed basis of $\mathfrak{g}^\wedge$.
\end{itemize}

We denote by \index{symbols!underlying space $\underline{S}^m$} $\underline{S}^m \subseteq C^\infty(\mathfrak{g}^\wedge)$ the subspace on which $\nu_{m,\alpha}$ is finite-valued for all $\alpha$.  We extend this definition to $m = \pm \infty$ by taking the union or intersection over all integers $m$.  Then $\underline{S}^m$ is naturally a Frechet space for $m < \infty$ and an inductive limit of such spaces for $m = \infty$.  We note that $\underline{S}^{-\infty}$ is the Schwartz space.

\subsubsection{Basic classes}
For fixed $m \in \mathbb{Z}$, we denote by \index{symbols!basic class $S_\delta^m$} $S^m_\delta \subseteq \underline{S}^{\infty}$ the subclass consisting of $a \in \underline{S}^\infty$ for which
\begin{equation*}
  \nu_{m,\alpha}(a) \ll \h^{-\delta |\alpha|}
\end{equation*}
for each fixed $\alpha$.  We write more verbosely $S_\delta^m(\mathfrak{g}^\wedge)$ or $S_\delta^m[\h]$ or $S_\delta^m(\mathfrak{g}^\wedge)[\h]$ when we wish to indicate the group and/or wavelength parameter under consideration.

We write $S_{\delta}^{-\infty}$ for the intersection of the classes $S^m_\delta$ taken over all fixed $m$.

For a scalar $c$, we write $c S_\delta^m$ for the class $\{c v : v \in S_\delta^m\}$.  We write $\h^\infty S_\delta^m$ for the intersection over all fixed $N$ of the classes $\h^N S_\delta^m$.

\subsubsection{The negligible class}
The class $\h^\infty S_\delta^{-\infty}$ is independent of $\delta$; it consists of $a \in \underline{S}^\infty$ satisfying $\partial^\alpha a(\xi) \ll \h^N \langle \xi  \rangle^{-N}$ for all fixed $\alpha$ and $N$.  We denote this class simply by $\h^\infty S^{-\infty}$.   \index{symbols!negligible class $\h^\infty S^{-\infty}$}

\begin{remark}
  We sometimes refer to the intersection $\h^\infty S^{-\infty} \cap \underline{S}^{-\infty}$, which is not redundant.  The class $\h^\infty S^{-\infty}$ consists of symbols whose \emph{fixed} symbol norms decay rapidly as $\h \rightarrow 0$, while the set $\underline{S}^{-\infty}$ consists of symbols for which \emph{all} symbol norms are finite, hence which define Schwartz functions.  In practice, the primary content is carried by $\h^\infty S^{-\infty}$ --- we intersect with $\underline{S}^{-\infty}$ when we need error terms that are genuinely Schwartz, as required by some results in the literature.
\end{remark}

\subsubsection{Refined classes}\label{sec:refined-classes}
Suppose now that $G$ is connected reductive, so that the notation $\mathfrak{g}^\wedge_{\reg}$ from \S\ref{sec:regular-elements} applies.  Assume also that $\h \lll 1$.

Let $\tau \in \mathfrak{g}^\wedge_{\reg}$.  We denote by $\mathfrak{g}_{\tau} \subseteq \mathfrak{g}$ the $\mathfrak{g}$-centralizer of $\tau$.  We denote by $\mathfrak{g}_\tau^{\perp} \subseteq \mathfrak{g}^\wedge$ the complement of $\mathfrak{g}_\tau$ with respect to the canonical bilinear pairing between $\mathfrak{g}$ and $\mathfrak{g}^\wedge$.  We fix an inner product on $\mathfrak{g}^\wedge$, hence a dual inner product on $\mathfrak{g}$.  This choice defines orthogonal complements $\mathfrak{g}_{\tau}^{\flat}$ (resp.\ $\mathfrak{g}_\tau^{\perp \flat}$) to $\mathfrak{g}_\tau$ (resp.\ $\mathfrak{g}_\tau^{\perp}$), hence orthogonal direct sum decompositions \index{Lie algebra!$\mathfrak{g}_\tau, \mathfrak{g}_\tau^\perp,
  \mathfrak{g}_\tau^{\flat}, \mathfrak{g}_\tau^{\perp \flat}$}
\begin{equation}\label{eq:mathfr-=-mathfr}
  \mathfrak{g}^\wedge = \mathfrak{g}_\tau^\perp \oplus \mathfrak{g}_\tau^{\perp \flat},
  \quad
  \mathfrak{g} = \mathfrak{g}_\tau^{\flat} \oplus \mathfrak{g}_\tau.
\end{equation}
Given $\xi \in \mathfrak{g}^\wedge$ (resp.\ $x \in \mathfrak{g}^\wedge$), we write
\begin{equation*}
  (\xi ', \xi '' ) \in \mathfrak{g}_\tau^{\perp} \times \mathfrak{g}_{\tau}^{\perp \flat}
  \quad
  \text{
    (resp.
    $(x ', x '' ) \in \mathfrak{g}_\tau^{\flat} \times  \mathfrak{g}_\tau$)
  }
\end{equation*}
for the coordinates with respect to \eqref{eq:mathfr-=-mathfr}.  We call these the \emph{$\tau$-coordinates} of $\xi$ (resp.\ $x$).  \index{Lie algebra!$\tau$-coordinates}

The inner products on $\mathfrak{g}$ and $\mathfrak{g}^\wedge$ induce inner products on each of the subspaces in \eqref{eq:mathfr-=-mathfr}.  We choose orthonormal bases for each of these subspaces, hence Euclidean coordinates with respect to which we employ multi-index notation (\S\ref{sec:multi-index-notation}) in what follows.

Let $\delta ', \delta ''$ be fixed quantities with $0 < \delta ' < \delta '' < 2 \delta ' < 1$, and let $\tau$ be an element of $\mathfrak{g}^\wedge$ that belongs to some fixed compact subset of $\mathfrak{g}^\wedge_{\reg}$.  We define $S^{\tau}_{\delta ', \delta ''} \subseteq \underline{S}^{-\infty}$ (more verbosely, $S_{\delta ', \delta ''}^{\tau}(\mathfrak{g}^\wedge)$, etc.) to be the subclass consisting of elements that are   \index{symbols!refined class $S^\tau_{\delta ', \delta ''}$}
\begin{itemize}
\item supported on $\tau + \O(\h^{ \delta '' - \delta '})$, and
\item satisfy, in $\tau$-coordinates, the derivative bounds
  \begin{equation}\label{eq:partial_xi--alpha}
    \partial_{\xi '} ^{\alpha } \partial_{\xi ''} ^{\beta } a (\xi) \ll \h^{- \delta ' |\alpha| - \delta '' |\beta|}
  \end{equation}
  for all fixed multi-indices $(\alpha,\beta) \in \mathbb{Z}_{\geq 0}^{\dim(\mathfrak{g}) - \rank(\mathfrak{g})} \times \mathbb{Z}_{\geq 0}^{\rank(\mathfrak{g})}$.
\end{itemize}
A convenient criterion for checking membership in the classes $S^\tau_{\delta ', \delta ''}$ is given as follows.
\begin{lemma}\label{lem:basepoint-insensitivity}
  With notation as above, let $\tau_1, \tau_2$ be elements of $\mathfrak{g}^\wedge$ that belong to some fixed compact subset of $\mathfrak{g}^\wedge_{reg}$.  If
  \begin{equation}\label{eq:tau_1-=-tau_2}
    \tau_1 = \tau_2 + \O (\h^{\delta '' - \delta '}),
  \end{equation}
  then $S_{\delta ', \delta ''}^{\tau_1} = S_{\delta ', \delta ''}^{\tau_2}$; otherwise, $S^{\tau_1}_{\delta', \delta ''} \cap S^{\tau_2}_{\delta', \delta ''} = \{0\}$.  More precisely, if \eqref{eq:tau_1-=-tau_2} holds, then for each smooth function $a : \mathfrak{g}^\wedge \rightarrow \mathbb{C}$ and each $\xi = \tau_j + \O (\h^{\delta '' - \delta '})$, the estimate \eqref{eq:partial_xi--alpha} holds with respect to $\tau_1$-coordinates if and only if it holds with respect to $\tau_2$-coordinates.
\end{lemma}
\begin{proof}
  This is \cite[Lemma 9.9]{2020arXiv201202187N}.
\end{proof}

\subsection{Cutoffs}

\subsubsection{Basic definition}
We fix a small enough symmetric neighborhood $\mathcal{G}_0$ of the origin $\mathfrak{g}$.  We assume in particular that
\begin{equation*}
  x \ast y := \log(\exp(x) \exp(y))
\end{equation*}
is defined and analytic for $(x,y) \in \mathcal{G}_0 \times \mathcal{G}_0$.  By a \index{Lie algebra!cutoff $\chi$} \emph{cutoff}, we mean a smooth function $\chi : \mathfrak{g} \rightarrow [0,1]$ that is supported on $\mathcal{G}_0$, satisfies $\chi(x) = \chi(-x)$, and is identically $1$ in some (not necessarily fixed) neighborhood of the origin.

It will always be possible, and often convenient, to choose the cutoff supported arbitrarily close to the origin.  In practice, we will work with cutoffs supported on elements of size $o(1)$, or much smaller.

Fix a slightly smaller fixed symmetric bounded neighborhood $\mathcal{G}_1$ with the property that the closure of $\mathcal{G}_1 \mathcal{G}_1 := G_1 \ast G_1$ is contained in $\mathcal{G}_0$.  Similarly construct $\mathcal{G}_2, \mathcal{G}_3, \dotsc$, with each bearing the same relationship to the previous as $\mathcal{G}_1$ does to $\mathcal{G}_0$ (i.e., the closure of $G_n \ast G_n$ is contained in $G_{n - 1}$).  Then for any cutoffs $\chi_1, \chi_2$ supported in $\mathcal{G}_{n+1}$, we may find a cutoff $\chi '$ supported on $\mathcal{G}_n$ with the property that $\chi ' (x \ast y) = 1$ whenever $\chi_1(x) \chi_2(y) \neq 0$.


\subsubsection{Basic quantitative variant}
For fixed $\delta \in (0, 1]$, we denote by $\mathcal{X}_{\delta}$, or $\mathcal{X}_\delta(\mathfrak{g})$ when we wish to explicate the group under consideration, the class of cutoffs $\chi$ on $\mathfrak{g}$  for which
\begin{itemize}
\item $\chi(x) = 1$ for $x = o(\h^{1 - \delta})$, and
\item $\partial^{\alpha} \chi(x) \ll \h^{-\O(1)}$ for each fixed $\alpha$.
\end{itemize}
The first condition may be reformulated as follows: there is a fixed $r > 0$ so that $\chi(x) = 1$ for $|x| \leq r \h^{1 - \delta}$.  It is easy to see that for each $\chi_1, \chi_2 \in \mathcal{X}_{\delta}$ supported in $\mathcal{G}_{n+1}$ for some fixed $n$, there exists $\chi ' \in \mathcal{X}_{\delta}$ supported in $\mathcal{G}_n$ so that $\chi '(x \ast y) = 1$ whenever $\chi_1(x) \chi_2(y) \neq 0$.

For fixed $\delta \in [0,1)$, we say that a cutoff $\chi$ is \emph{$\delta$-admissible} if it lies in $\mathcal{X}_{\delta_+}$ for some fixed $\delta_+ > \delta$.  For example, any fixed cutoff is $\delta$-admissible.  The point of this definition is that it captures the essential support properties of the Fourier transforms of symbols in $S_\delta^\infty$.
\begin{lemma}\label{lem:fix-0-leq}
  Fix $0 \leq \delta < 1$.  Let $a \in S_\delta^{\infty}$.  Let $\chi$ be a $\delta$-admissible cutoff.  Define
  \begin{equation}\label{eq:a-:=-chi}
    a' := (\chi a_{\h}^\vee)^\wedge.
  \end{equation}
  Then
  \begin{equation}\label{eq:a-equiv-pmodhinfty}
    a' \equiv a_{\h} \mod{\h^\infty S^{-\infty} \cap \underline{S}^{-\infty}}.
  \end{equation}
\end{lemma}
\begin{proof}
  As noted in \cite[\S8.2]{nelson-venkatesh-1}, one sees by repeated partial integration that the symbol $a$ admits a distributional Fourier transform, smooth away from the origin, that decays rapidly, together with all derivatives, outside the support of $\chi$.  More precisely,
  \begin{equation*}
    \partial^{\alpha} ((\chi - 1) a_{\h}^\vee) (x)
    \ll
    \h^{N(\delta_+ - \delta)} \langle x \rangle^{- N}
  \end{equation*}
  for each fixed $\alpha$ and $N$ and some fixed $\delta_+ > \delta$.    The required bound then follows by repeated integration by parts in the definition of the Fourier integral.
\end{proof}

\subsubsection{Refined quantitative variant}
An analogous discussion applies to the refined symbol classes.  Let $\tau$ belong to some fixed compact subset of $\mathfrak{g}^\wedge_{\reg}$, and let $0 < \delta ' < \delta '' < 2 \delta '' < 1$ be fixed.  Using $\tau$-coordinates on $\mathfrak{g}$ (\S\ref{sec:refined-classes}), we define for $r > 0$ the set
\begin{equation*}
  \mathcal{D} (r) := \left\{ x \in \mathfrak{g} : |x'| \leq r \h^{1 - \delta '}, |x''| \leq r \h^{1 - \delta ''} \right\}.
\end{equation*}
Informally, the scaling in the definition of $\mathcal{D}(r)$ is Fourier-dual to the ``coin-shaped'' frequency regions described in \S\ref{sec:relat-trace-form}: localization of a function to $|x'| \ll \h^{1-\delta'}$ and $|x''| \ll \h^{1-\delta''}$ corresponds to its Fourier transform being essentially supported on a region with $|\xi'| \ll \h^{-1+\delta'}$ and $|\xi''| \ll \h^{-1+\delta''}$, with $\delta''>\delta'$ producing a ``thicker'' direction and a ``thinner'' direction in frequency space.

We denote by $\mathcal{X}_{\tau,\delta ', \delta ''}$ (more verbosely, $\mathcal{X}_{\tau, \delta ', \delta ''}(\mathfrak{g})$) the class of cutoffs $\chi$ on $\mathfrak{g}$ such that
\begin{itemize}
\item there is a fixed $r > 0$ so that $\chi(x) = 1$ for $x \in \mathcal{D}(r)$, and
\item $\partial^{\alpha} \chi(x) \ll \h^{-\O(1)}$ for each fixed $\alpha$.
\end{itemize}
We say that a cutoff $\chi$ is \emph{$(\tau, \delta ', \delta '')$-admissible} if it lies in $\mathcal{X}_{\tau, \delta '_+, \delta ''_+}$ for some fixed $\delta'_+ > \delta '$ and $\delta ''_+ > \delta ''$.  We have the following analogue of Lemma \ref{lem:fix-0-leq}, which follows again by partial integration (compare with \cite[(14.6)]{2020arXiv201202187N}).
\begin{lemma}\label{lem:let-tau-belong}
  Let $\tau$ belong to a fixed compact subset of $\mathfrak{g}^\wedge_{\reg}$, and let $0 < \delta ' < \delta '' < 2 \delta ' < 1$ be fixed.  Let $a \in S^\tau_{\delta ', \delta ''}$.  Let $\chi$ be a $(\tau,\delta ', \delta '')$-admissible cutoff.  Define $a'$ as in \eqref{eq:a-:=-chi}.  Then \eqref{eq:a-equiv-pmodhinfty} holds.
\end{lemma}

\subsection{Star product}\label{sec:star-product}

\subsubsection{Fourier transform and basic definition}\label{sec:four-transf-basic}
Recall that for any function $a : \mathfrak{g}^\wedge \rightarrow \mathbb{C}$, we define its rescaling $a_{\h}$ by the formula $a_{\h}(\xi) := a(\h \xi)$.  If $a$ is integrable, then its rescaled inverse Fourier transform is given explicitly by
\begin{equation*}
  a_{\h}^\vee(x)
  :=(a_{\h})^\vee(x)
  = \h^{-\dim(\mathfrak{g})} a^\vee(x/\h).
\end{equation*}
Multiplying the latter by a cutoff $\chi$ defines a compactly-supported function $\chi a_{\h}^\vee$ on $\mathfrak{g}$.  More canonically, we obtain a compactly-supported measure $\chi a_{\h}^\vee \, d x$ on $\mathfrak{g}$.  We denote by \index{operators!$\widetilde{\Opp}$}
\begin{equation*}
  \widetilde{\Opp}_{\h}(a:\chi)
\end{equation*}
the compactly-supported measure on $G$ obtained via pushforward under the exponential map of $\chi a_{\h}^\vee \, d x$.  This measure assigns to any test function $\phi$ on $G$ the integral
\begin{equation*}
  \int_{x \in \mathfrak{g} } \chi(x) a_{\h}^\vee(x) \phi(\exp(x)) \, d x.
\end{equation*}
When $G$ is unimodular, we will often identify $\widetilde{\Opp}_{\h}(a:\chi)$ with the function on $G$ obtained by dividing by some fixed Haar measure $d g$ on $G$.    The choice of wavelength parameter $\h$ and cutoff $\chi$ is often clear by context; in such cases, we abbreviate $\widetilde{\Opp}(a) := \widetilde{\Opp}_{\h}(a:\chi)$.

The convolution of such measures is described by a certain bilinear operator $\star_{\h}$, called the (rescaled) \emph{star product}, for which \index{symbols!rescaled star product $\star_{\h}$}
\begin{equation}\label{eq:widet-ast-widet}
  \widetilde{\Opp}_{\h}(a:\chi_1) \ast
  \widetilde{\Opp}_{\h}(b:\chi_2)
  =
  \widetilde{\Opp}_{\h}(a \star_{\h} b:\chi')
\end{equation}
if $\chi_1, \chi_2, \chi '$ are cutoffs for which $\chi ' (x \ast y) = 1$ when $\chi_1(x) \chi_2(y) \neq 0$.  Explicitly, the star product is characterized by
\begin{equation*}
  (a \star_{\h} b)_{\h} = \chi_1 a_{\h}^\vee \star \chi_2 b_{\h}^\vee,
\end{equation*}
where on the RHS, $\star$ denotes the result of taking the measures $\chi_{1} a_{\h}^\vee \, d x$ and $\chi_{2} b_{\h}^\vee \, d y$ on $\mathfrak{g}$, pushing them forward to $G$, convolving them, and then pulling the result back to $\mathfrak{g}$ (see \cite[\S8.7]{2020arXiv201202187N} or \cite[\S2.5]{nelson-venkatesh-1} for details).  It is clear that $\star_{\h}$ defines a continuous bilinear operation on the Schwartz space $\underline{S}^{- \infty}$.  It depends upon $\chi_1$ and $\chi_2$; we suppress this dependence from our notation because it is mild (e.g., in the senses quantified by Lemmas \ref{lem:fix-0-leq}, \ref{lem:let-tau-belong}, \ref{lem:fix-delta-in-1}, \ref{lem:fix-delta-in-2}).

\subsubsection{Formal expansion}\label{sec:formal-expansion}
The star product admits a formal asymptotic expansion  \index{symbols!star product component $\star^j$}
\begin{equation*}
  a \star_{\h} b \sim \sum_{j \geq 0} \h^j a \star^j b,
\end{equation*}
discussed in detail in \cite[\S4.2, \S4.6]{nelson-venkatesh-1} and \cite[\S9.2]{2020arXiv201202187N}.  We recall parts of that discussion.  We have $a \star^0 b = ab$.  We note that $\star_{\h}$ depends on a choice of cutoffs, but the bidifferential operators $\star^j$ are determined solely by the Taylor expansion of the Baker--Campbell--Hausdorff series at the origin, hence are independent of that choice.  For $(x,y,\zeta) \in \mathfrak{g} \times \mathfrak{g} \times \mathfrak{g}^\wedge$ with $x,y$ sufficiently small, we set
\begin{equation*}
  \{x,y\} := \log(\exp(x) \exp(y)) - x - y.
\end{equation*}
Then the function $(x,y,\zeta) \mapsto \exp (\langle \{x, y\}, \zeta  \rangle)$ is analytic near the origin, hence admits a Taylor expansion, which we denote in multi-index notation by
\begin{equation*}
  \exp (\langle \{x, y\}, \zeta  \rangle)
  = \sum_{\alpha, \beta , \gamma } c_{\alpha \beta \gamma } x^\alpha y^\beta \zeta^\gamma.
\end{equation*}
The $\star^j$ are then the finite bidifferential operators, with polynomial coefficients, given by
\begin{equation*}
  a \star^j b(\zeta) = \sum  _{|\alpha| + |\beta| - |\gamma| = j} c_{\alpha \beta \gamma} \zeta^\gamma \partial^\alpha a(\zeta) \partial^\beta b(\zeta).
\end{equation*}
The coefficients $c_{\alpha \beta \gamma}$ satisfy the basic support conditions
\begin{equation}\label{eq:cuhmwyqdas}
  |\gamma| \leq \min(|\alpha|,|\beta|), \qquad  \max(\lvert \alpha \rvert, \lvert \beta \rvert) \leq j
\end{equation}
(see \cite[(4.5)]{nelson-venkatesh-1}, \cite[(9.3)]{2020arXiv201202187N}; in fact, the second condition in \eqref{eq:cuhmwyqdas} follows from the first and the identity $\lvert \alpha \rvert + \lvert \beta \rvert - \lvert \gamma \rvert = j$).  A finer support condition may be described by writing $\exp (\langle \{x, y\}, \tau  \rangle)$ in $\tau$-coordinates as
\begin{equation*}
  \sum_{\alpha ', \alpha '', \beta ', \beta '' , \gamma }
  c_{\alpha ' \alpha '' \beta ' \beta '' \gamma }
  (x ' ) ^{\alpha ' }
  (x '') ^{\alpha ''}
  (y ') ^{\beta '}
  (y '' ) ^{\beta ''}
  \tau^\gamma.
\end{equation*}
The coefficients in this expansion are then nonvanishing only if (see \cite[Lemma 9.13]{2020arXiv201202187N})
\begin{equation*}
  |\alpha '| + |\beta '| \geq 2 |\gamma|.
\end{equation*}
Setting $j = |\alpha'| + |\alpha ''| + |\beta '| + |\beta ''| - |\gamma|$, this last inequality may be rewritten
\begin{equation}\label{eq:alpha-refined-support-condition}
  |\alpha '| + |\beta ' | + 2 | \alpha '' | + 2 |\beta ''| \leq 2 j.
\end{equation}

\begin{example}\label{exa:star1-structure}
  Let us explicate $\star^j$ when $j=1$.  Using that $\{x,0\}=0=\{0,y\}$, we see that $c_{\alpha\beta\gamma}=0$ whenever $|\alpha|=0$ or $|\beta|=0$, except when $(\alpha,\beta,\gamma)=(0,0,0)$.  In particular, for $j\ge 1$ the operator $\star^j$ involves only terms with $|\alpha|\ge 1$ and $|\beta|\ge 1$.  Specializing to $j=1$, the conditions $|\alpha|+|\beta|-|\gamma|=1$ and $|\gamma|\le \min(|\alpha|,|\beta|)$ force $|\alpha|=|\beta|=|\gamma|=1$.  Thus $\star^1$ has bidegree $(1,1)$, with homogeneous linear coefficients (in fact, it is a constant multiple of the Poisson bracket \cite{gutt1983explicit}).  This is consistent with the fact that the BCHD law is given to first order by a multiple of the Lie bracket.  Moreover, in $\tau$-coordinates, we see that the terms contributing to $\star^1$ have $(\alpha '', \beta '') =(0,0)$.
\end{example}

\subsubsection{Extended star product}

\begin{theorem}\label{thm:extended-star-product}
  Let $\chi_1$ and $\chi_2$ be cutoffs.  Then the star product, defined initially for Schwartz functions using $\chi_1$ and $\chi_2$, extends uniquely to a compatible family of continuous bilinear maps
  \begin{equation*}
    \star_{\h} : \underline{S}^{m_1} \times \underline{S}^{m_2} \rightarrow \underline{S}^{m_1 + m_2},
  \end{equation*}
  with the convention that $\infty + (-\infty) := -\infty$.  For each $(a,b) \in \underline{S}^{m_1} \times \underline{S}^{m_2}$ and $J \in \mathbb{Z}_{\geq 0}$, we have
  \begin{equation*}
    a \star_{\h} b \equiv \sum_{0 \leq j < J} \h^j a \star^j b \mod{ \underline{S} ^{m_1 + m_2 - J} },
  \end{equation*}
  where the remainder term defines a continuous map $\underline{S} ^{m_1 } \times \underline{S} ^{m_2 } \rightarrow \underline{S} ^{m_1 + m_2 - J}$.
\end{theorem}
\begin{proof}
  This is \cite[Thm 1]{nelson-venkatesh-1}.
\end{proof}

\subsubsection{Asymptotics}\label{sec:star-prod-asymptotics}

\begin{theorem}\label{thm:basic-star-prod}
  Fix $\delta \in [0,1/2)$ and $m_1, m_2 \in \mathbb{Z} \cup \{\pm \infty \}$.  Let $\chi_{1}$ and $\chi_2$ be fixed cutoffs.  Let $\star_{\h}$ denote the star product, as defined using $\chi_1$ and $\chi_2$.  Then the extended rescaled star product $\star_{\h} : \underline{S}^\infty \times \underline{S}^\infty \rightarrow \underline{S}^\infty$ induces class maps
  \begin{equation*}
    \star_{\h} : S_{\delta}^{m_1} \times S_{\delta}^{m_2} \rightarrow S_{\delta}^{m_1 + m_2}.
  \end{equation*}
  For $(a,b) \in S_{\delta}^{m_1} \times S_{\delta}^{m_2}$ and fixed $j \in \mathbb{Z}_{\geq 0}$, we have
  \begin{equation*}
    a \star^j b \in \h^{-2 \delta j} S_\delta^{m_1 + m_2 - j}.
  \end{equation*}
  Moreover, for fixed $J \in \mathbb{Z}_{\geq 0}$,
  \begin{equation*}
    a \star_{\h} b \equiv \sum_{0 \leq j < J} \h^j a \star^j b \mod{
      \h^{(1 - 2 \delta) J} S_\delta  ^{m_1 + m_2 - J}.
    }
  \end{equation*}
\end{theorem}
\begin{proof}
  The required conclusion for infinite values of $m$ follows formally from the case of finite values of $m$.  In that case, the claimed result is given by \cite[Thm 9.5]{2020arXiv201202187N}.
\end{proof}
\begin{lemma}\label{lem:concl-theor-refthm:b}
  The conclusion of Theorem \ref{thm:basic-star-prod} remains valid if we let $\chi_1$ and $\chi_2$ be $\delta$-admissible cutoffs supported in $\mathcal{G}_1$, not necessarily fixed.
\end{lemma}
\begin{proof}
  Let $\chi_1', \chi_2'$ be fixed cutoffs supported in $\mathcal{G}_1$; write $\star_{\h}'$ for the star product defined using them, to which Theorem \ref{thm:basic-star-prod} applies.  Write $\star_{\h}$ for that defined using $\chi_1$ and $\chi_2$.  Assume first that $\chi_2 = \chi_2'$.
  By definition,
  \begin{equation*}
    (a \star_{\h} b - a \star_{\h} ' b)_{\h}^\vee
    = (\chi_1 - \chi_1 ') a_{\h}^\vee \star \chi_2 ' b_{\h}^\vee.
  \end{equation*}
  Define $a'' \in \underline{S}^\infty$ by requiring that $(a_{\h}'')^\vee  := (\chi_1 - \chi_1 ') a_{\h}^\vee$.  By Lemma \ref{lem:fix-0-leq}, we have $a'' \in \h^\infty S^{-\infty} \cap \underline{S}^{-\infty}$.  We may choose a fixed cutoff $\chi_1 ''$, supported in $\mathcal{G}_0$, with $\chi_1'' \equiv 1$ on the supports of $\chi_1$ and $\chi_1 '$.  Then $\chi_1 '' (a_{\h}'')^\vee  = (a_{\h}'')^\vee$.  Define $\star_{\h}''$ using the fixed cutoffs $\chi_1''$ and $\chi_2$.  Then
  \begin{equation*}
    a \star_{\h} b - a \star_{\h}' b = a'' \star_{\h}'' b,
  \end{equation*}
  which, by Theorems \ref{thm:extended-star-product} and \ref{thm:basic-star-prod}, lies in $\h^\infty S^{-\infty} \cap \underline{S}^{-\infty}$.

  We argue similarly if $\chi_1 = \chi_1'$, and then conclude via bilinearity.
\end{proof}

\begin{theorem}\label{thm:refined-star-prod}
  Assume that $\h \lll 1$.  Fix $\delta ', \delta ''$ satisfying $0 < \delta ' < \delta '' < 2 \delta ' < 1$.  Let $\chi_1$ and $\chi_2$ be fixed cutoffs; define $\star_{\h}$ using these.  Write $j, J$ for fixed elements of $\mathbb{Z}_{\geq 0}$ and $\tau, \tau_1, \tau_2$ for
  elements of some fixed compact subset of $\mathfrak{g}^\wedge_{\reg}$.
  \begin{enumerate}[(i)]
  \item \label{item:star-prod-1} $\star_{\h}$ enjoys the mapping property
    \begin{equation}\label{eqn:star-h-new-mapping-property-super-localized}
      \star_{\h} :
      S^{\tau_1}_{\delta', \delta''}
      \times
      S^{\tau_2}_{\delta', \delta''}
      \rightarrow
      S^{\tau_1}_{\delta', \delta ''}
      \cap S^{\tau_2}_{\delta', \delta ''}
      + \h^{\infty} S^{-\infty}.
    \end{equation}
  \item \label{item:star-prod-2} If $\tau_1 - \tau_2 \ggg \h^{\delta '' - \delta '}$, then $a \star^j b = 0$ for all $a \in S^{\tau_1}_{\delta',\delta''}$, $b \in S^{\tau_2}_{\delta ', \delta ''}$.
  \item \label{item:star-prod-3} $\star^j$ enjoys the mapping property
    \begin{equation}\label{eqn:star-j-new-mapping-property}
      \star^j :
      S_{\delta ', \delta ''}^{\tau }
      \times
      S_{\delta ', \delta ''}^{\tau }
      \rightarrow \h^{- 2 \delta' j}
      S^{\tau}_{\delta', \delta ''}.
    \end{equation}
    For any $a,b \in S_{\delta ', \delta ''}^\tau$, we have the asymptotic expansion
    \begin{equation}\label{eqn:star-j-asymp-expn-new}
      a \star_{\h} b
      \equiv  \sum_{0 \leq j < J}
      \h^j a \star^j b
      \mod{
        \h^{(1- 2 \delta') J} S^{\tau}_{\delta', \delta''}
        +
        \h^{\infty}
        S^{-\infty}}.
    \end{equation}
  \item \label{item:star-prod-4} Fix $\delta \in [0,\delta']$.  We have the mapping properties
    \begin{equation*}
      \star_{\h} : S^{\infty}_{\delta} \times S^{\tau}_{ \delta
        ', \delta ''}
      \rightarrow S^{\tau}_{ \delta ', \delta ''}
      + \h^\infty S^{-\infty},
    \end{equation*}
    \begin{equation*}
      \star^j : S^{\infty}_{\delta} \times S^{\tau}_{ \delta
        ', \delta ''}      \rightarrow \h^{-(\delta +  \delta ')
        j}
      S^{\tau}_{ \delta ', \delta ''}
    \end{equation*}
    and, for $(a,b) \in S^{\infty}_{\delta} \times S^{\tau}_{ \delta ', \delta ''}$, the asymptotic expansion
    \begin{equation*}
      a \star_{\h} b
      \equiv  \sum_{0 \leq j < J}
      \h^j a \star^j b
      \mod{
        \h^{(1- \delta - \delta') J}  S^{\tau}_{ \delta ', \delta ''}
        +
        \h^{\infty} S^{-\infty}}.
    \end{equation*}
    Analogous results hold with $S^{\infty}_{\delta} \times S^{\tau}_{ \delta ', \delta ''}$ replaced by $S^{\tau}_{ \delta ', \delta ''} \times S^{\infty}_{\delta}$.
  \end{enumerate}
\end{theorem}
\begin{proof}
  This is \cite[Thm 9.12]{2020arXiv201202187N}.
\end{proof}

The analogue of Lemma \ref{lem:concl-theor-refthm:b} holds, with the same proof: for each assertion involving $\star_{\h}$ applied to an element $a$ of the class $S_\delta^\infty$ (resp.\ $S_{\delta ', \delta ''}^{\tau}$), we can take the cutoff $\chi$ assigned to $a$ to be any $\delta$-admissible (resp.\ $(\tau, \delta ', \delta '')$-admissible) cutoff.

From now on, we do not write explicitly which cutoffs are used to define $\star_{\h}$, which should be clear by context in each case.

\subsection{Operator assignment}
Let $\pi$ be a unitary representation of $G$.

\subsubsection{Basic definition}\label{sec:basic-definition}
We denote as in \S\ref{sec:kirillov-model} by
\begin{equation*}
  \pi^\infty \subseteq \pi^0 \subseteq \pi^{-\infty}
\end{equation*}
the spaces of smooth, norm-bounded and distributional vectors.  By an \emph{operator} on $\pi$, we mean a linear map $\pi^{-\infty} \rightarrow \pi^{-\infty}$.  (The operators of interest to us will have more specific mapping properties than this.)

Let $a \in \underline{S}^{- \infty}$.  Let $\chi$ be a cutoff.  Recall, from \S\ref{sec:four-transf-basic}, the definition of the smooth compactly-supported measure $\widetilde{\Opp}_{\h}(a:\chi)$ on $G$.  We obtain an operator $\Opp_{\h}(a:\pi,\chi)$ on $\pi$ by integrating: \index{operators!$\Opp$}
\begin{equation*}
  \Opp_{\h}(a:\pi,\chi) := \int_{g \in G} \widetilde{\Opp}_{\h}(a:\chi) \pi(g)
  =
  \int_{x \in \mathfrak{g} } \chi(x) a_{\h}^\vee(x) \pi(\exp(x)) \, d x.
\end{equation*}
The convolution formula \eqref{eq:widet-ast-widet} translates to the composition formula
\begin{equation}\label{eq:composition-formula-most-basic}
  {\Opp}_{\h}(a:\pi,\chi_1) \ast
  {\Opp}_{\h}(b:\pi,\chi_2)
  =
  {\Opp}_{\h}(a \star_{\h} b:\pi,\chi')
\end{equation}

In what follows, the representation $\pi$ will always be clear from context, so we often abbreviate
\begin{equation*}
  \Opp_{\h}(a:\chi) := \Opp_{\h}(a:\pi,\chi).
\end{equation*}

\subsubsection{Some operator classes}\label{sec:negligible-operators}
The discussion here is highly non-exhaustive, since we require far less precision in our discussion of operators than in \cite{nelson-venkatesh-1, 2020arXiv201202187N}.
\begin{definition}
  Let $\underline{\Psi }^{-\infty}$ denote the space of operators $T$ on $\pi$ such that for all $x,y \in \mathfrak{U}(G)$, the operator $\pi(x) T \pi(y)$ induces a bounded operator on $\pi$.\index{operators!underlying space $\underline{\Psi}^{-\infty}$ of smoothing operators}
\end{definition}
\begin{definition}
  Let $\underline{\Psi }^{\infty}$ denote the space \index{operators!underlying space $\underline{\Psi}^{\infty}$ of finite-order operators} of operators $T$ on $\pi$ such that there exists $m \in \mathbb{Z}_{\geq 0}$ so that for all $k \in \mathbb{Z}_{\geq 0}$ and $x_1,\dotsc,x_k \in \mathfrak{g}$, there exists $C_0 \geq 0$ so that for all $v \in \pi^\infty$, we have
  \begin{equation*}
    \|[\pi(x_1), \dotsc, [\pi(x_k), T] ] v\| \leq C_0 \sum_{y_1,\dotsc,y_m \in \mathcal{B}(G) \cup \{1\}} \|\pi(y_1 \dotsb y_m) v\|,
  \end{equation*}
  where we recall from \S\ref{sec:cumos34zyo} that $\mathcal{B}(G)$ is a fixed basis for $\Lie(G)$.
\end{definition}
\begin{example}
  For $f \in C_c^\infty(G)$, the integral operator $\pi(f)$ lies in $\underline{\Psi }^{-\infty}$.  For $x \in \mathfrak{U}(G)$, the differential operator $\pi(x)$ lies in $\underline{\Psi }^{\infty}$.
\end{example}
\begin{remark}
  By \cite[(3.1), Lem 3.5, Prop 8.2]{nelson-venkatesh-1}, the abbreviated definitions of $\underline{\Psi }^{\pm \infty}$ given here coincide with the those in \cite[\S10]{2020arXiv201202187N}.
\end{remark}
\begin{lemma}\label{lem:foll-incl-hold}
  The following inclusions hold, where juxtaposition denotes either composition of operators or application of operators to vectors.
  \begin{enumerate}[(i)]
  \item $\underline{\Psi}^{-\infty} \subseteq \underline{\Psi}^{\infty}$.
  \item $\underline{\Psi}^{\infty} \underline{\Psi}^{-\infty}, \underline{\Psi}^{-\infty} \underline{\Psi}^{\infty} \subseteq \underline{\Psi}^{\infty}$.
  \item $\underline{\Psi }^{\infty} \pi^\infty \subseteq \pi^\infty$.
  \item $\underline{\Psi }^{-\infty} \pi^{-\infty} \subseteq \pi^\infty$.
  \end{enumerate}
\end{lemma}
\begin{proof}
  These are all essentially immediate from the definitions.  We refer to \cite[\S3]{nelson-venkatesh-1} for further discussion.
\end{proof}
\begin{definition}
  Let $\Psi^{-\infty}$ denote the class of all $T \in \underline{\Psi }^\infty$ such that for all fixed $x,y \in \mathfrak{U}(G)$, the operator $\pi(x) T \pi(y)$ induces a bounded operator on $\pi$ with operator norm $\O(1)$.

  We may define the class $c \Psi^{-\infty}$ for any scalar $c$.  We denote by $\h^\infty \Psi^{-\infty}$ the intersection of the classes $\h^N \Psi^{-\infty}$ over all fixed $N$.    \index{operators!operator class $\Psi^{-\infty}$}
\end{definition}

\begin{lemma}
  If $T_1$ and $T_2$ lie in $\Psi^{-\infty}$, then so does their composition $T_1 T_2$.
\end{lemma}
\begin{proof}
  This is a special case of \cite[Lemma 10.2]{2020arXiv201202187N} (see also \cite[\S3.4]{nelson-venkatesh-1}).
\end{proof}



\subsubsection{Qualitative theory}

\begin{theorem}\label{thm:qualitative-extended-operator-assignment}
  Let $\chi$ be a cutoff.  The assignment
  \begin{equation*}
    \underline{S}^{-\infty} \ni a \mapsto \Opp_{\h}(a:\chi) =: \Opp(a)
  \end{equation*}
  extends to a linear map
  \begin{equation*}
    \Opp : \underline{S}^{\infty} \rightarrow
    \{\text{operators on } \pi\}
  \end{equation*}
  with the following properties.
  \begin{enumerate}[(i)]
  \item $\Opp(\underline{S}^\infty) \subseteq \underline{\Psi}^{\infty}$.
  \item $\Opp(\underline{S}^{-\infty}) \subseteq \underline{\Psi}^{-\infty}$.
  \item $\Opp(\underline{S}^{\infty}) \pi^\infty \subseteq \pi^\infty$.
  \item\label{enumerate:cuhmv6fu1d} $\Opp(\underline{S}^{-\infty}) \pi^{-\infty} \subseteq \pi^\infty$.
  \item\label{enumerate:cuhmv6fwiz} $\Opp (\underline{S}^0) \pi \subseteq \pi$.
  \item The composition formula \eqref{eq:composition-formula-most-basic} remains valid, i.e., for $\chi_1, \chi_2, \chi '$ as before and $a, b \in \underline{S}^{\infty}$,
    \begin{equation}\label{eq:qualitative-composition-formula}
      \Opp_{\h}(a:\chi_1)
      \Opp_{\h}(b:\chi_2)
      =
      \Opp_{\h}(a \star_{\h} b:\chi').
    \end{equation}
  \item\label{enumerate:cui1f04tkm} For any polynomial $p \in \underline{S}^\infty$ (i.e., $p \in \Sym(\mathfrak{g}_{\mathbb{C}})$), we have
    \begin{equation}\label{eq:oppp-=-pisymp_h}
      \Opp(p) = \pi(\sym(p_{\h})),
    \end{equation}
    where $p_{\h}$ denotes the rescaling $p_{\h}(\xi) := p(\h \xi)$ and $\sym : \Sym(\mathfrak{g}_\mathbb{C}) \rightarrow \mathfrak{U}(\mathfrak{g}_\mathbb{C})$ denotes the symmetrization map, i.e., the linear isomorphism sending each monomial to the average of its permutations.
  \item\label{enumerate:cui1f0676u} For two polynomial symbols $p, q \in \Sym(\mathfrak{g}_{\mathbb{C}})$, the composition formula holds exactly: $\Opp(p) \Opp(q) = \Opp(p \star_{\h} q)$.
  \end{enumerate}
\end{theorem}
\begin{proof}
  The extension, inclusions and composition formula are special cases of \cite[Thm 10.6]{2020arXiv201202187N}.  For example, \emph{loc.\ cit.}\ asserts that \(\Opp(\underline{S}^m) \subseteq \underline{\Psi}^m\) for all \(m\), which gives \(\Opp(\underline{S}^m) \pi^s \subseteq \pi^{s - m}\) for the Sobolev spaces \(\pi^s\) (\(s \in \mathbb{Z}\)) defined in \cite[\S10.1]{2020arXiv201202187N}.  Taking \(m = - \infty\) gives \(\Opp(\underline{S}^{-\infty}) \pi^s \subseteq \pi^\infty\) for each \(s\); taking the union over \(s\) then gives \eqref{enumerate:cuhmv6fu1d}.  Taking \((m, s) =(0, 0)\) gives \eqref{enumerate:cuhmv6fwiz}.

  Assertion \eqref{enumerate:cui1f04tkm} is given by \cite[\S5.2, \S8.1]{nelson-venkatesh-1}.

  For assertion \eqref{enumerate:cui1f0676u}, we specialize \eqref{eq:qualitative-composition-formula} to the case that $a$ and $b$ are polynomials, and must check then that the right-hand side is independent of the choice of cutoff $\chi '$.  Indeed, given any such choice, we may find $\chi_1$ and $\chi_2$ that are compatible with $\chi '$ in the usual sense ($\chi '(x \ast y) = 1$ whenever $\chi_1(x) \chi_2(y) \neq 0$), which follows from the continuity of $(x,y) \mapsto x \ast y$.  Then, applying \eqref{eq:qualitative-composition-formula} in reverse and \eqref{eq:oppp-=-pisymp_h} to each factor, we obtain $\Opp(a \star_{\h} b, \chi ') = \pi(\sym(a_{\h})) \pi(\sym(b_{\h}))$, which is indeed independent of $\chi '$.
\end{proof}

\begin{lemma}\label{lem:let-chi-be}
  Let $\chi$ be any $0$-admissible cutoff.  Abbreviate $\Opp(a) := \Opp_{\h}(a:\chi)$.  Then
  \begin{equation*}
    \Opp(\h^\infty S^{-\infty}) \subseteq \h^\infty \Psi^{-\infty}.
  \end{equation*}
\end{lemma}
\begin{proof}
  This follows in much sharper forms from the results of \cite[\S10]{2020arXiv201202187N}.
\end{proof}

\subsubsection{Irrelevance of cutoffs}
\begin{lemma}\label{lem:fix-delta-in-1}
  Fix $\delta \in [0,1)$.  Let $a \in S_\delta^\infty$.  Let $\chi_1$ and $\chi_2$ be $\delta$-admissible cutoffs, supported in $\mathcal{G}_1$.  Then
  \begin{equation}\label{eq:opp_h-opp_h-in}
    \Opp_{\h}(a:\chi_1) - \Opp_{\h}(a:\chi_2) \in \h^\infty \Psi^{-\infty} \cap \underline{\Psi }^{-\infty}.
  \end{equation}
  More precisely, if $\chi '$ is a $\delta$-admissible cutoff satisfying $\chi ' (x \ast y) = 1$ whenever $\chi_1(x) \chi_2(y) \neq 0$, then
  \begin{equation}\label{eq:opp_h-opp_h-=}
    \Opp_{\h}(a:\chi_1) - \Opp_{\h}(a:\chi_2) = \Opp_{\h}(a' : \chi ') \quad \text{for some $a' \in \h^\infty S^{-\infty} \cap \underline{S}^{-\infty}$.}
  \end{equation}
\end{lemma}
\begin{proof}
  By Lemma \ref{lem:fix-0-leq}, the difference $\Opp_{\h}(a:\chi_1) - \Opp_{\h}(a:\chi_2)$ lies in the class $\Opp_{\h}(\h^\infty S^{-\infty} \cap \underline{S}^{-\infty})$.  By Theorem \ref{thm:qualitative-extended-operator-assignment} and Lemma \ref{lem:let-chi-be}, that class is contained in $\h^\infty \Psi^{-\infty} \cap \underline{\Psi}^{-\infty}$.  The same argument gives the more precise assertion.  (Compare with \cite[Lemma 10.7]{2020arXiv201202187N} and \cite[\S5.4]{nelson-venkatesh-1}.)
\end{proof}

\begin{lemma}\label{lem:fix-delta-in-2}
  Let $\tau$ belong to a fixed compact subset of $\mathfrak{g}^\wedge_{\reg}$.  Fix $0 < \delta ' < \delta '' < 2 \delta ' < 1$.  Let $a \in S_{\delta', \delta ''}^{\tau}$.  Let $\chi_1$ and $\chi_2$ be $(\tau, \delta ', \delta '')$-admissible cutoffs.  Then \eqref{eq:opp_h-opp_h-in} holds.  More precisely, if $\chi '$ is any $(\tau, \delta ', \delta '')$-admissible cutoff with $\chi ' (x \ast y) = 1$ whenever $\chi_1(x) \chi_2(y) \neq 0$, then \eqref{eq:opp_h-opp_h-=} holds.
\end{lemma}
\begin{proof}
  The proof is similar to that of Lemma \ref{lem:fix-delta-in-1}, using Lemma \ref{lem:let-tau-belong} instead of Lemma \ref{lem:fix-0-leq}.
\end{proof}

\subsubsection{Composition}

\begin{theorem}\label{thm:composition-formula-quantitative}
  Fix $\delta \in [0,1/2)$.  For \(a \in S_\delta^\infty\), abbreviate $\Opp(a) := \Opp_{\h}(a:\chi)$ for some $\delta$-admissible cutoff $\chi$, supported in $\mathcal{G}_2$.  Then for $a,b \in S_\delta^\infty$, we have
  \begin{equation}\label{eq:oppa-oppb-equiv}
    \Opp(a) \Opp(b) \equiv \Opp(a \star_{\h} b) \mod{ \h^\infty \Psi^{-\infty} \cap \underline{\Psi }^{-\infty}}.
  \end{equation}
\end{theorem}
\begin{proof}
  By \eqref{eq:qualitative-composition-formula} and Lemma \ref{lem:fix-delta-in-1}.  (Compare with \cite[Thm 10.8, (10.8)]{2020arXiv201202187N}.)
\end{proof}

Similarly:
\begin{theorem}
  Suppose $G$ is connected reductive and $\h \lll 1$.  Fix $0 < \delta ' < \delta '' < 2 \delta ' < 1$, and let $\tau$ belong to a fixed compact subset of $\mathfrak{g}^\wedge_{\reg}$.  Set $\Opp(a) := \Opp_{\h}(a:\chi)$ for some $(\tau, \delta ', \delta'')$-admissible cutoff $\chi$, supported in $\mathcal{G}_2$.  Then the composition formula \eqref{eq:oppa-oppb-equiv} remains valid for $a,b \in S_{\delta '}^{\infty} \cup S^{\tau}_{\delta ', \delta ''}$.
\end{theorem}

Many variants on these results are possible, and we do not state them exhaustively.  The point is just that for all practical purposes, when working with a symbol $a \in S^\infty_\delta$ (resp.\ $a \in S^\tau_{\delta ', \delta ''}$), we can define $\Opp_{\h}(a:\chi)$ using a $\delta$-admissible (resp.\ $(\tau, \delta ', \delta '')$-admissible) cutoff $\chi$ supported on $o(1)$.  We can then compose such operators using \eqref{eq:oppa-oppb-equiv}.  The ambiguity arising from the choice of cutoff always lies in the negligible class $\h^\infty \Psi^{-\infty} \cap \underline{\Psi }^{-\infty}$.

\subsection{Borel summation}

\begin{lemma}\label{lem:borel-summation}
  Fix $0 \leq \delta < 1$.  Suppose given for each $j \in \mathbb{Z}_{\geq 0}$ a smooth function $a_j$ on $\mathfrak{g}^\wedge$, a real number $e_j$ and an integer $m_j$ so that
  \begin{itemize}
  \item for each fixed $j$, the quantities $e_j$ and $m_j$ are fixed, and we have $a_j \in \h^{e_j} S^{m_j}_\delta$, and
  \item $e_j \rightarrow \infty$ and $m_j \rightarrow -\infty$ as $j \rightarrow \infty$.
  \end{itemize}
  Then for each $L_0 \ggg 1$, there exists $L \in \mathbb{Z}_{\geq 0}$ with $1 \lll L  \leq L_0$ so that the sum
  \begin{equation}\label{eq::=-sum-_0}
    a := \sum_{0 \leq j < L} a_j
  \end{equation}
  belongs to $S^{\infty}_{\delta}$ and satisfies, for fixed $N \in \mathbb{Z}_{\geq 0}$ and large enough fixed $J \in \mathbb{Z}_{\geq 0}$, the asymptotic expansion
  \begin{equation}\label{eq:cui1g9lsvk}
    a \equiv \sum_{0 \leq j < J} a_j \mod{\h^N S^{-N}_\delta}.
  \end{equation}
\end{lemma}
\begin{proof}
  This is a simple overspill exercise, essentially identical to the standard argument for establishing Borel summation in the pseudodifferential calculus.  Define $J_0(N)$ to be the smallest $J$ for which $j \geq J \implies a_j \in \h^N S_\delta^{-N}$.  Temporarily suppose $L$ is fixed.  Then for all $N \leq L$ and all $J$ with $J_0(N) \leq J \leq L$, the sum $b := \sum_{J \leq j < L} a_j$ lies in the symbol class $\h^{N} S_\delta^{-N}$.  This last membership means that $b$ lies in the space of symbols $\underline{S}^{\infty}$ and satisfies $\nu_{-N,\alpha}(b) \ll \h^{N - \delta |\alpha|}$ for each fixed multi-index $\alpha$.  In particular, since $\h \lll 1$, we have $\nu_{-N,\alpha}(b) \leq \h^{N - 1 - \delta |\alpha|}$ for all multi-indices $\alpha$ with $|\alpha| \leq L$.  Thus the set (cf.\ \eqref{eq:cumotvxqcc})
  \begin{equation*}
    \left\{
      L \in \mathbb{Z}_{\geq 0} :
      \begin{gathered}
        b \in \underline{S}^{\infty} \text{ and }
        \nu_{-N,\alpha}(b) \leq \h^{N - 1 - \delta |\alpha|} \\
        \text{ for all } N \leq L, \, J_0(N) \leq J \leq L, \,  |\alpha| \leq L
      \end{gathered}
    \right\},
  \end{equation*}
  with $b$ defined as above, contains all fixed $L$.  By overspill, it contains some $L \ggg 1$ with $L \leq L_0$.  For this $L$, we have $N \leq L$ all fixed $N$.  Thus for each fixed $N$ and all fixed $J \geq J_0(N)$, we have
  \begin{equation*}
    \sum_{J \leq j < L } a_j \in \h^{N-1} S_{\delta}^{-N}.
  \end{equation*}
  The required memberships follow readily from this.\end{proof}

\begin{lemma}
  Let $\tau \in \mathfrak{g}^\wedge$ belong to a fixed compact set of regular elements.  Fix $0 < \delta ' < \delta '' < 2 \delta ' < 1$.  Suppose given for each $j \in \mathbb{Z}_{\geq 0}$ a smooth function $a_j$ on $\mathfrak{g}^\wedge$ and a real number $e_j$ so that
  \begin{itemize}
  \item for each fixed $j$, the quantities $e_j$ and $m_j$ are fixed, and we have $a_j \in \h^{e_j} S^{\tau}_{\delta', \delta ''}$, and
  \item $e_j \rightarrow \infty$ as $j \rightarrow \infty$.
  \end{itemize}
  Assume moreover that there is a fixed $R > 0$ so that each $a_j$ is supported on $\{\xi : |\xi - \tau| \leq R \h^{\delta '' - \delta '}\}$.  Then for each $L_0 \ggg 1$, there exists $L \in \mathbb{Z}_{\geq 0}$ with $1 \lll L \leq L_0$ so that the sum $a := \sum_{0 \leq j < L} a_j$ belongs to $S^{\tau}_{\delta', \delta ''}$ and satisfies, for fixed $N \in \mathbb{Z}_{\geq 0}$ and large enough fixed $J \in \mathbb{Z}_{\geq 0}$, the asymptotic expansion
  \[
    a \equiv \sum_{0 \leq j < J} a_j \mod{\h^N S^{\tau }_{\delta', \delta ''}}.
  \]
\end{lemma}
\begin{proof}
  The proof is very similar to that of Lemma \ref{lem:borel-summation}: in place of $\nu_{-N,\alpha}$, we use the seminorms relevant for the definition of $S^\tau_{\delta ', \delta ''}$.  The important point is that since we have fixed $R$, the support condition in the definition of $S_{\delta ', \delta ''}^{\tau}$ holds automatically for the partial sum $a$.
\end{proof}

\section{Abstract study of localized vectors}\label{sec:abstr-study-local}
In this section, we work over any fixed Lie group $G$.  We assume given a wavelength parameter $\h \in (0,1]$ and a fixed element $0 \leq \delta < 1/2$, so that we may define the symbol class $S_\delta^\infty = S_\delta^\infty(\mathfrak{g}^\wedge)$.  Many definitions in this section, like that of the symbol class, depend upon the parameter $\h$, but this dependence is suppressed to simplify the notation.

Our aim is to develop some formal consequences of the symbol calculus in a general setting, where the necessary inputs to our arguments can be presented most clearly.  Example \ref{exa:cal-W-delta-mirabolic}, below, suggests our primary motivation doing so.


\subsection{$\O(1)$-modules}
Recall from \S\ref{sec:classes}  our convention on ``class.''  We introduce one further specialized piece of terminology.
\begin{definition}
  An \emph{$\O(1)$-module} is a subclass $M$ of a complex vector space for which $0 \in M$,
  \[
    m_1, m_2 \in M \implies m_1 + m_2 \in M
  \]
  and
  \[
    c \in \mathbb{C}, \quad c \ll 1, \quad m \in M \implies c m \in M.
  \]
\end{definition}
The conditions defining an $\O(1)$-module $M$ may be rephrased as follows: for all fixed $n \in \mathbb{Z}_{\geq 0}$ and all scalars (i.e., complex numbers) $c_1, \dotsc,c_n = \O(1)$, we have
\begin{equation*}
  m_1,\dotsc,m_n \in M \implies c_1 m_1 + \dotsb + c_n m_n \in M.
\end{equation*}
In particular, an $\O(1)$-module need not be a linear subspace: it is only required to be stable under \emph{bounded} linear combinations.
\begin{example}
  Given a vector space $V$ equipped with some collection of seminorms $\{\nu_i\}_{i \in I}$, the class of all $v \in V$ for which $\nu_i(v) = \O(1)$ for each fixed $i$ is an $\O(1)$-module.  Note that neither $V$ nor the collection of seminorms need be fixed.  For instance, $V$ could be the space $C_c^\infty(E)$ for a non-fixed subspace $E$ of some fixed Euclidean space.  Examples like the latter arise naturally for us and may motivate why we have not phrased matters in terms of, say, bounded subsets of Fréchet spaces.
\end{example}
The following lemma will not be used at all; it is included just for the sake of illustration.
\begin{lemma}
  Let $M$ be an $\O(1)$-module.  Then $M$ is a subspace if and only if it is a set.
\end{lemma}
\begin{proof}
  That subspaces are sets follows from the definitions.  Conversely, suppose $M$ is a set.  By definition, $M$ contains zero, the sum of any two elements, and the negative of any element (because $-1 \ll 1$), so it is an abelian subgroup.  It remains only to verify that it is closed under scalar multiplication.  Let $m \in M$.  Consider the set $S_m := \{c \in \mathbb{C} : c m \in M\}$; we must check that $S_m = \mathbb{C}$.  Since $M$ is an abelian subgroup, so is $S_m$.  By hypothesis, $S_m$ contains all scalars of the form $\O(1)$, and in particular, the unit disc, which generates $\mathbb{C}$, giving $S_m = \mathbb{C}$, as required.
\end{proof}

We note that if $M$ is an $\O(1)$-module, then so is $c M$ for any scalar $c$.  For an $\O(1)$-module $M$, we write
\begin{equation*}
  \h^\infty M
\end{equation*}
for the intersection of $\h^N M$ taken over fixed integers $N$; it is an $\O(1)$-module.

\subsection{Localized vectors: definition}\label{sec:local-vect-defin}
First, some general notation.  For $\tau \in \mathfrak{g}^\wedge$ with $\tau \ll 1$ and fixed $r \in \mathbb{Z}_{\geq 0}$, we write
\[I_{0,\tau}^r \subseteq S_0^\infty\] for the subclass consisting of symbols that vanish to order at least $r$ at $\tau$.  The polynomial elements of $I_{0,\tau}^r$ are then precisely those polynomials of degree $\O(1)$ and coefficients of size $\O(1)$ that vanish to order at least $r$ at $\tau$.

For any polynomial function $p : \mathfrak{g}^\wedge \rightarrow \mathbb{C}$, we define $\Opp(p) \in \mathfrak{U}(\mathfrak{g}_\mathbb{C})$ by
\begin{equation}\label{eq:oppp-:=-symp_h}
  \Opp(p) := \sym(p_{\h}),
\end{equation}
where, as in Theorem \ref{thm:qualitative-extended-operator-assignment}, we denote by $\sym : \Sym(\mathfrak{g}_\mathbb{C}) \rightarrow \mathfrak{U}(\mathfrak{g}_\mathbb{C})$ the symmetrization map and $p_{\h}(\xi) := p(\h \xi)$ the $\h$-rescaling of $p$.  As noted in the statement of that theorem, the formula \eqref{eq:oppp-:=-symp_h} is satisfied by our ``usual'' operator assignment.

\begin{definition}\label{defn:localized-vectors}
  Let $M$ be an $\O(1)$-submodule of a $\mathfrak{U}(\mathfrak{g}_\mathbb{C})$-module.  Fix $\nu \in (0,1/2]$.  Let $\tau \in \mathfrak{g}^\wedge$ with $\tau \ll 1$.  We say that $v \in M$ is \emph{$(\mathfrak{g},\nu)$-localized at $\tau$ inside $M$} if for each fixed $r \in \mathbb{Z}_{\geq 0}$ and each polynomial $p \in I_{0,\tau}^r$, we have
  \[
    \Opp(p) v \in \h^{r \nu} M.
  \]
  When $M$ is clear from context, we omit ``inside $M$.''  We say that $M$ is $(\mathfrak{g},\nu)$-localized at $\tau$ if every $v \in M$ has this property.  When $\mathfrak{g}$ is clear from context, we write simply ``$\nu$-localized.''
\end{definition}
\begin{remark}
  Definition \ref{defn:localized-vectors} depends upon the $\O(1)$-module $M$ to which we regard $v$ as belonging.
\end{remark}
\begin{remark}\label{rmk:r-=-0}
  The $r = 0$ case of the condition in Definition \ref{defn:localized-vectors} is equivalent to the following: for fixed $m \in \mathbb{Z}_{\geq 0}$ and fixed $x_1,\dotsc,x_m \in \mathfrak{g}$, we have
  \begin{equation*}
    \pi(x_1 \dotsb x_m) v \in \h^{-m} M.
  \end{equation*}
\end{remark}

We record a simple criterion.
\begin{lemma}\label{lem:localization-criterion-degree-1}
  Let $M$ be an $O(1)$-submodule of a $\mathfrak{U}(\mathfrak{g}_\mathbb{C})$-module.  Suppose that for each $v \in M$ and each fixed $x \in \mathfrak{g}$, we have $\Opp(x) v \equiv x(\tau) v \mod{\h^\nu M}$.  Then $M$ is $\nu$-localized at $\tau$.
\end{lemma}
\begin{proof}
  We must show that for each fixed $r \in \mathbb{Z}_{\geq 0}$ and each polynomial $p \in I^r_{0,\tau}$, we have
  \begin{equation}\label{eqn:Opp-p-M-h-r-nu-M}
    \Opp(p) M \subseteq  \h^{r \nu} M.
  \end{equation}
  By linearity, we may assume that $p$ is homogeneous and vanishes at $\tau$ to order exactly $r$, and then further that $p(\xi) = \prod_{j=1}^r (x_j(\xi) - x_j(\tau))$ for some fixed $x_1,\dotsc,x_r \in \mathfrak{g}$.  By hypothesis, $\Opp(x_j - x_j(\tau)) M \subseteq \h^{\nu} M$.  Iterating this hypothesis yields the required containment \eqref{eqn:Opp-p-M-h-r-nu-M}.
\end{proof}

\subsection{$S^\infty_\delta$-modules}\label{sec:sinfty_delta-modules}

\begin{definition}\label{defn:s-infinity-delta-module}
  By an \emph{$S_\delta^\infty$-module}, we mean a triple $(\underline{M}, M, \Opp)$, where
  \begin{itemize}
  \item $\underline{M}$ is a $\mathfrak{U}(\mathfrak{g}_\mathbb{C})$-module,
  \item $M \subseteq \underline{M}$ is an $\O(1)$-submodule, and
  \item $\Opp$ is a linear map $\Opp : \underline{S}^\infty \rightarrow \End(\underline{M})$
  \end{itemize}
  satisfying the following conditions.
  \begin{enumerate}[(i)]
  \item  \label{item:Sinfdelmod-1} For each $p \in \Sym(\mathfrak{g})$ and $v \in \underline{M}$, we have $\Opp(p) v = \sym(p_{\h}) v$.
  \item \label{item:Sinfdelmod-2} For each $a \in S_\delta^\infty$ and $v \in M$, we have $\Opp(a) v \in M$.  Thus $\Opp(a)$ defines a class map $M \rightarrow M$.
  \item \label{item:Sinfdelmod-3} For each $a,b \in S^\infty_\delta$ and $v \in M$, we have $\Opp(a) \Opp(b) v \equiv \Opp(a \star_{\h} b) v \pmod{\h^\infty M}$.
  \end{enumerate}
  When we wish to emphasize which group is being considered, we write more verbosely ``$S^\infty_\delta(\mathfrak{g}^\wedge)$-module.''
\end{definition}
We often abbreviate the tuple $(\underline{M}, M, \Opp)$ by simply $M$, and write, e.g., ``Let $M$ be an $S_\delta^\infty$-module.''  We then refer to $\underline{M}$ as the underlying $\mathfrak{U}(\mathfrak{g}_\mathbb{C})$-module for $M$.

\begin{example}
  For any unitary representation $\pi$ of $G$, the space $M := \pi^\infty$ of smooth vectors is naturally an $S_\delta^\infty$-module, but not a very interesting one.  The issue is that $c M = M$ for \emph{all} nonzero scalars $c$.  In particular, $\h^\ell M = M$ for all $\ell$.  The more interesting examples noted below are obtained by restricting to proper subclasses of $\pi^{\infty}$ consisting of vectors that are ``controlled'' in a suitable sense.
\end{example}

\begin{example}\label{example:class-V-wavelength-h}
  Let $\pi$ be a unitary representation of $G$.  Let $\mathcal{V}$ \index{classes!$\mathcal{V}$} denote the class consisting of all smooth vectors $v \in \pi$ such that for (all fixed $m \in \mathbb{Z}_{\geq 0}$  and) all fixed $x_1,\dotsc,x_m \in \mathfrak{g}$, we have
  \[
    \|\pi(x_1 \dotsb x_m) v\| \ll \h^{-m}.
  \]
  Then $\mathcal{V}$ is naturally an $S^\infty_\delta$-module for each fixed $\delta \in [0,1/2)$.  Explicitly, we may take as the underlying $\mathfrak{U}(\mathfrak{g}_\mathbb{C})$-module the space of all smooth vectors in $\pi$, and define $\Opp(a) v := \Opp_{\h}(a:\chi) v$ for some $\chi \in \mathcal{X}_{\delta_+}$, with $\delta_+ \in (\delta, 1/2)$ fixed.

  Informally, $\mathcal{V}$ consists of vectors having norm $\ll 1$ and wavelength $\gg \h$.
\end{example}

\begin{example}
  If $M$ is an $S_\delta^\infty$-module, then so is $c M$ for any scalar $c$.  If $c \ll 1$, then $c M$ is a $S_\delta^\infty$-submodule of $M$.
\end{example}

\begin{example}\label{exa:cal-W-delta-mirabolic}
  Take $G = \GL_n(\mathbb{R})$, with mirabolic subgroup $P$.  With $T := 1/\h$, the class $\mathcal{W}_\delta = \mathfrak{W}(\pi,\psi,T,\delta)$ (Definition \ref{defn:we-denote-mathfr-frakW} and \S\ref{sec:quant-prel}), consisting of Whittaker functions $W \in \mathcal{W}(\pi,\psi_T)$ having certain smoothness and support properties, is naturally an $S_\delta^\infty(\mathfrak{p}^\wedge)$-module (Lemma \ref{lem:class-map-defined-cal-W-delta-module}).  It is not hard to see (Proposition \ref{prop:local-with-resp-p}) that, regarding $\mathcal{W}_\delta$ as an $\O(1)$-submodule of a $\mathfrak{U}(\mathfrak{p}_\mathbb{C})$-module, it is $(\mathfrak{p},\delta)$-localized at $\theta = \theta_P(\psi)$.  We will show eventually (Theorem \ref{thm:cal-W-delta-g-acts-localized}) that, regarding $\mathcal{W}_\delta$ as an $\O(1)$-submodule of a $\mathfrak{U}(\mathfrak{g}_\mathbb{C})$-module, it is $(\mathfrak{g},\delta)$-localized at $\tau = \tau(\pi,\psi,T)$.  This passage from $\mathfrak{p}$ to $\mathfrak{g}$ is the main step in the proof of the main result (Theorem \ref{thm:big-group-whittaker-behavior-under-G}) of Part \ref{part:asympt-analys-kirill}.  The main point of \S\ref{sec:abstr-study-local} is to encapsulate some of the technical arguments required to justify that passage.
\end{example}

\begin{example}\label{example:delta-to-zero-module}
  If $\delta_0 < \delta$, then any $S^\infty_\delta$-module is naturally an $S^\infty_{\delta_0}$-module, hence in particular an $S^\infty_0$-module.
\end{example}

The following result describes some basic properties of localized vectors:
\begin{enumerate}[(i)]
\item they are approximately annihilated by symbols that vanish in a large enough neighborhood of the point at which they localize, and
\item the condition of Definition \ref{defn:localized-vectors}, formulated using polynomials, extends to more general symbols.
\end{enumerate}

\begin{theorem}\label{thm:S-infinity-modules-localized-stuff}
  Fix $\delta \in [0,1/2)$ and $\nu \in (0,1/2]$.  Let $M$ be an $S^\infty_\delta$-module.  Let $v \in M$ be $\nu$-localized at $\tau$.
  \begin{enumerate}[(i)]
  \item \label{item:1-thm:S-infinity-modules-localized-stuff} If $\nu > \delta$, then $\Opp(a) v \in \h^\infty M$ for each $a \in S^\infty_\delta$ that vanishes on $\tau + o(\h^\delta)$.
  \item \label{item:2-thm:S-infinity-modules-localized-stuff} If $\delta > 0$, then $\Opp(a) v \in \h^{r \nu} M$ for each $a \in I_{0,\tau}^r$ (not necessarily polynomial).
  \end{enumerate}
\end{theorem}
\begin{proof}[Proof of part (i)]
  We define $d \in S^2_0$ by $d(\xi) = |\xi - \tau|^2$.  Then $d^r \in S^{2 r}_0$ is a polynomial symbol that vanishes to order $2 r$ at $\tau$, i.e., $d^r \in I_{0,\tau}^{2 r}$, so by hypothesis,
  \begin{equation}\label{eq:oppdr-v-in}
    \Opp(d^r) v \in \h^{2 r \nu } M.
  \end{equation}

  The proof consists of a few steps.

  \emph{Step 1: discarding large frequencies.}  First, we fix a radius $R > 0$ with $R \asymp 1$, taken large enough that $|\tau| \leq R/2$.  We will verify that $\Opp(a) v \in \h^\infty M$ for any $a \in S^0_0$ for which $a(\xi) = 0$ whenever $|\xi| \leq R$.  To that end, we construct for fixed $\ell \geq 0$ a symbol $q \in S^{-2 \ell}_0$ so that
  \[
    q \star_{\h} d^{\ell} \equiv a\pmod{ \h^\infty S^{-\infty} }.
  \]
  We do this by setting
  \begin{equation}\label{eqn:divide-star-prod-by-d-N}
    q_0 :=
    \frac{
      a
    }{
      d^\ell
    },
    \quad
    q_1 := \frac{- q_0 \star^1 d^\ell }{d^\ell},
    \quad
    q_2 :=
    \frac{- q_0 \star^2 d^\ell
      + q_1 \star^1 d^\ell
    }{
      d^\ell
    },
  \end{equation}
  and so on, so that, formally,
  \[
    (\sum_{k \geq 0} \h^k q_k) \star_{\h} d^\ell \sim a.
  \]
  The point is that division by $d^\ell$ is defined because all symbols occurring in the numerators of \eqref{eqn:divide-star-prod-by-d-N} vanish near the origin.  By iterated application of the quotient rule for derivatives and the mapping properties \eqref{eqn:star-j-new-mapping-property} for $\star^j$, we see that
  \begin{equation*}
    q_j \in S^{-2 \ell-j}_0.
  \end{equation*}
  Setting $q := \sum_{j \leq J} \h^j q_j$ for some small enough $J \ggg 1$ (in the sense afforded by Lemma \ref{lem:borel-summation}), we have $q \in S^{-2 \ell}_0$ and
  \begin{equation}\label{eqn:q-star-d-N-equiv-1-minus-a}
    q \star_{\h} d^\ell \equiv a \pmod{\h^\infty S^{-\infty}}.
  \end{equation}
  We now apply the axioms of Definition \ref{defn:s-infinity-delta-module}: the axiom \eqref{item:Sinfdelmod-2} gives \(\Opp(\h^\infty S^{-\infty}) M \subseteq \h^\infty M\) (allowing us to discard the error term above), and then the composition formula axiom \eqref{item:Sinfdelmod-3} gives
  \begin{equation}\label{eqn:composition-formula-appl-to-v-minus-opp-a-v}
    \Opp(a) v
    \equiv
    \Opp(q) \Opp(d^\ell) v \mod{\h^\infty M}.
  \end{equation}
  By another application of axiom \eqref{item:Sinfdelmod-2}, we have $\Opp(q) M \subseteq M$.  By \eqref{eq:oppdr-v-in}, we have $\Opp(d^\ell) v \in \h^{2 \ell \nu} M$.  Therefore $\Opp(a) v \in \h^{2 \ell \nu} M$.  Since $\ell$ was arbitrary, we conclude as required that $\Opp(a) v \in \h^\infty M$.

  \emph{Step 2: reduction to symbols with support $\O(1)$.}  Assuming now that the lemma holds in the special case that $a$ is supported on $\O(1)$, we deduce the general case.  To that end, we fix $R > 0$ large enough, as in Step 1, and choose a symbol $c \in S^{-\infty}_0$ taking the value $1$ on $\{\xi : |\xi| \leq R\}$ and vanishing outside $\{\xi : |\xi| < 2 R\}$.  Here and henceforth, let $\equiv$ denote equivalence modulo $\h^\infty M$.  The symbol $c-1 \in S^0_0$ vanishes on $\{ \xi : |\xi| \leq R \}$, so by Step 1, we have $v \equiv \Opp(c) v$.  It follows that for any $a \in S_\delta^\infty$, we have
  \[
    \Opp(a) v \equiv \Opp(a) \Opp(c) v \equiv \Opp(a \star_{\h} c) v.
  \]
  (In the first step, we used that $\Opp(a) \h^\infty M \subseteq \h^\infty M$, which is a consequence of the more precise containment $\Opp(a) M \subseteq M$ following from axiom  \eqref{item:Sinfdelmod-2}.  In the second step, we used axiom \eqref{item:Sinfdelmod-3}.)
  Setting $b := \sum_{j \leq J} \h^j a \star^j c$ for any small enough $J \ggg 1$ (Lemma \ref{lem:borel-summation}), we have $b \in S^{-\infty}_\delta$ and $\Opp(a) v \equiv \Opp(b) v$.  Suppose now that $a$ vanishes on $\tau + o(\h^\delta)$.  Since $\supp(b) \subseteq \supp(a) \cap \supp(c)$, we see that $b$ vanishes on $\tau + o(\h^\delta)$ and is supported on $\O(1)$.  By our hypothesis that the conclusion holds for symbols of support $\O(1)$ (such as $b$), we have $\Opp(b) v \in \h^\infty M$.  The required conclusion $\Opp(a) v \in \h^\infty M$ then follows.  This completes the proof of the claimed reduction.

  \emph{Step 3: completion of the proof.}  We suppose now that $a \in S^{\infty}_\delta$ is supported on $\O(1)$ (hence lies in $S^{-\infty}_\delta$) and vanishes on $\tau + o(\h^\delta)$.  By overspill, we may find a fixed $r_0 > 0$ so that $a$ vanishes on $\{\tau + \xi : |\xi| \leq r_0 \h^{\delta}\}$.  We repeat the division argument described in Step 1 to obtain a symbol $q \in \h^{-2 \ell \delta} S^{-\infty}_\delta$ for which \eqref{eqn:q-star-d-N-equiv-1-minus-a} holds. 
  As before, $\Opp(a) v \equiv \Opp(q) \Opp(d^\ell ) v \mod {\h^\infty M}$.  By hypothesis, $\Opp(d^\ell) v \in \h^{2 \ell \nu} M$ and $\Opp(q) M \subseteq \h^{-2 \ell \delta} M$, hence $\Opp(a) v \in \h^{2 \ell (\nu - \delta)} M$.  Since $\nu > \delta$, we may conclude by taking $\ell$ sufficiently large.
\end{proof}
\begin{proof}[Proof of part (ii)]
  We may assume -- after shrinking $\delta$ as necessary (cf.\ Example \ref{example:delta-to-zero-module}) -- that $0 < \delta < \nu$.  In particular, the conclusion of part (i) applies.  Let $a \in I_{0,\tau}^r$.  Fix $\ell \geq r$, and let $p$ denote the Taylor polynomial approximation of degree $\ell-1$ for $a$ at $\tau$.  Then $p \in S^{\ell-1}_0 \cap I_{\tau,0}^r$, so by hypothesis, $\Opp(p) v \in \h^{r \nu} M$.  Set $q := a - p$.  Then $q \in S^\infty_0$.  Moreover, for $\xi \ll 1$ and each fixed multi-index $\alpha$, we have by Taylor's theorem
  \begin{equation}\label{eq:cuhmwsku3k}
    \partial^\alpha q(\tau + \xi) \ll \min(1, |\xi|^{\ell-|\alpha|}).
  \end{equation}
  To estimate $\Opp(q) v$, we choose an ``envelope'' $c \in S^{-\infty}_\delta$ for which
  \begin{equation*}
    c(\xi) = \begin{cases}
      1 & \text{ if } \xi = \tau + o(\h^\delta), \\
      0 &  \text{ if } \xi \neq \tau + \O(\h^\delta).
    \end{cases}
  \end{equation*}
  Then $c-1$ vanishes on $\tau + o(\h^\delta)$, so part (i) gives $v \equiv \Opp(c) v \mod{\h^\infty M}$.

  On the other hand, the noted derivative bounds for $q$, together with the support condition on $c$, imply that
  \begin{equation}\label{eq:cuhmwswmbr}
    q \star_{\h} c \in \h^{\ell \delta} S_\delta^{-\infty}.
  \end{equation}
  To see this, we apply the star product asymptotics \ref{thm:basic-star-prod} for the symbol class map \(\star_{\h} : S_\delta^\infty \times S_\delta^{- \infty} \rightarrow S_\delta^{- \infty}\) (using that \(S_0^\infty \subseteq S_\delta^\infty\)), with the index truncation \(J\) taken large enough that \((1 - 2 \delta) J \geq \ell \delta\).  In this way, we reduce to checking that \(\h^j q \star^j c \in \h^{\ell \delta} S_\delta^{-\infty}\) for each fixed \(j\).  As noted in \cite[\S4.2]{nelson-venkatesh-1}, the bidifferential operator \(\star^j\) is of degree \(\leq j\) in each variable, with polynomial coefficients, hence lies in \(S^\infty\); since multiplication induces \(S^\infty \times S_\delta^{- \infty} \rightarrow S_\delta^{- \infty}\), we reduce to checking that
  \begin{equation*}
    \h^j (\partial^\alpha q)(\partial^\beta c) \in \h^{\ell \delta} S_\delta^{- \infty}
  \end{equation*}
  for any fixed multi-indices \(\alpha\) and \(\beta\) with \(\lvert \alpha \rvert, \lvert \beta \rvert \leq j\).  We now open the definition of \(S_\delta^{- \infty}\) and use again \(c\) is supported on a fixed compact set to reduce to checking that \(\partial^\gamma \left( \h^j(\partial^\alpha q)(\partial^\beta c) \right)\) has \(L^\infty\)-norm \(\ll \h^{\ell \delta} \h^{- \delta \lvert \gamma \rvert}\) for each fixed multi-index \(\gamma\).  We may write \(\partial^\gamma \left( (\partial^\alpha q)(\partial^\beta c) \right)\) as a linear combination, with fixed coefficients, of \((\partial^{\alpha + \mu} q)(\partial^{\beta + \nu} c)\), where \(\mu\) and \(\nu\) are fixed multi-indices with \(\mu + \nu = \gamma\), hence \(\lvert \mu \rvert + \lvert \nu \rvert = \lvert \gamma \rvert\).  Since \(c \in S_\delta^{-\infty}\), we have \(\partial^{\beta + \nu}(\tau + \xi) \ll \h^{- \delta(\lvert \beta \rvert + \lvert \nu \rvert)}\).  By \eqref{eq:cuhmwsku3k}, we have \(\partial^{\alpha + \mu} q(\tau + \xi) \ll \lvert \xi \rvert^{\ell - \lvert \alpha \rvert - \lvert \mu \rvert}\), which is \(\ll \h^{\delta \left( \ell - \lvert \alpha \rvert - \lvert \mu \rvert \right)}\) on the support of \(c\) (or its derivatives).  Multiplying these bounds together, we reduce to checking that
  \begin{equation*}
    \h^j \h^{\delta(\ell - \lvert \alpha \rvert - \lvert \mu \rvert)} \h^{- \delta(\lvert \beta \rvert + \lvert \nu \rvert)} \ll \h^{\ell \delta - \delta(\lvert \mu \rvert + \lvert \nu \rvert)}.
  \end{equation*}
  Canceling the common factors involving \(\mu\), \(\nu\) and \(\ell\), we reduce to checking that \(\h^j \h^{- \delta(\lvert \alpha \rvert + \lvert \beta \rvert)} \ll 1\), which follows from the fact that \(\lvert \delta \rvert \leq 1/2\) and \(\lvert \alpha \rvert, \lvert \beta \rvert \leq j\).  This completes the verification of \eqref{eq:cuhmwswmbr}.

  Applying \(\Opp\) to \eqref{eq:cuhmwswmbr} gives
  \[
    \Opp(q) v \equiv \Opp(q) \Opp(c) v \equiv \Opp(q \star_{\h} c) v \in \h^{\ell \delta } M,
  \]
  and so
  \[
    \Opp(a) v = \Opp(p) v + \Opp(q) v \in \h^{r \nu} M + \h^{\ell \delta} M.
  \]
  Since $\delta > 0$, we may choose $\ell$ large enough that $\ell \delta > r \nu$.  The proof is then complete.
\end{proof}

\section{Localizing refined symbols near an infinitesimal character
}\label{sec:local-refin-symb}
In this section we record an estimate (Theorem \ref{thm:local-refin-symb}) to be applied below in \S\ref{sec:appr-idemp}.  While we do not apply that estimate further in this paper, it may be more broadly useful, so we develop it here a bit more generally than required.

Throughout \S\ref{sec:local-refin-symb}, we adopt the following setting.  Fix a connected reductive group $G$ over the reals.  Let $\h > 0$ with $\h \lll 1$.  Let
\begin{equation}\label{eq:cuhnuoy8vh}
  \tau \in(\text{some fixed compact subset of } \mathfrak{g}^\wedge_{\reg}).
\end{equation}
Fix $0 < \delta ' < \delta '' < 2 \delta ' < 1$.  We may then define the symbol class $S^\tau_{\delta ', \delta ''}$.

\subsection{Statement of result}
\begin{theorem}\label{thm:local-refin-symb}
  Let $\pi$ be an irreducible unitary representation of $G$.
  Let $a \in S^\tau_{\delta', \delta ''}$.  Suppose that $[\tau] = \h \lambda_\pi$ and that, in $\tau$-coordinates,
  \begin{equation}\label{eqn:xi-delta-primez-tau}
    \supp(a) \subseteq
    \left\{\xi : \xi' - \tau ' \ll \h^{\delta '},
      \xi '' - \tau '' \gg \h^{\delta ''}
    \right\}
  \end{equation}
  Let $\chi$ be any $(\tau, \delta ', \delta '')$-admissible cutoff.  Then $\Opp(a) := \Opp_{\h}(a:\pi, \chi )$ belongs to the operator class $\h^\infty \Psi^{-\infty}$.
\end{theorem}

The proof is given in \S\ref{sec:completion-proof} after some preliminaries.  The informal content is that if $a$ is sufficiently separated from the locus described by the infinitesimal character of $\pi$, then $\Opp(a)$ is negligibly small.  The result is ``sharp up to epsilons'' for reasons explained in the final paragraph of \cite[\S2.4]{2020arXiv201202187N} .

We note that the special case of Theorem \ref{thm:local-refin-symb} in which $\pi$ is \emph{tempered} follows from the proof of \cite[Lemma 14.3]{2020arXiv201202187N}, which invoked the Kirillov character formula.  The main point of this section is to give a more general argument that applies also in the non-tempered case.  As a substitute for the Kirillov formula, we integrate by parts with respect to the center of the universal enveloping algebra.


\subsection{Using invariant polynomials to separate a symbol from a coadjoint orbit}
\label{sec:}
\begin{proposition}\label{lem:produce-inv-polynomials-large-on-symbol-support}
  Let $a \in S^\tau_{\delta ', \delta ''}$ satisfy the support condition \eqref{eqn:xi-delta-primez-tau}.  Set $n := \rank(\mathfrak{g})$.  There exist $a_1,\dotsc,a_n \in S^\tau_{\delta ', \delta''}$ and $G$-invariant polynomial symbols $f_1,\dotsc,f_n \in S^\infty_0$ with the following properties.
  \begin{enumerate}[(i)]
  \item $a = a_1 + \dotsb + a_n$.
  \item $f_1(\tau) = \dotsb = f_n(\tau) = 0$.
  \item $f_j(\xi) \gg \h^{\delta ''}$ for all $\xi \in \supp(a_j)$.
  \end{enumerate}
\end{proposition}
The proof is given below after a few preliminaries.

As noted in \cite[\S9.2]{nelson-venkatesh-1}, the ring $\Sym(i \mathfrak{g})^G$ of real-valued $G$-invariant polynomials $\mathfrak{g}^\wedge \rightarrow \mathbb{R}$ admits a finite set of algebraically independent homogeneous generators $p_1,\dotsc,p_n$.  Fixing such generators yields an identification of $[\mathfrak{g}^\wedge]$ with the Euclidean space $\mathbb{R}^n$.  We may speak in particular of the Euclidean distance, denoted $\dist(\lambda,\mu)$, between elements of $[\mathfrak{g}^\wedge]$.  Recall that we have equipped the summands in the decomposition
\begin{equation}\label{eqn:g-wedge-tau-decmop-std}
  \mathfrak{g}^\wedge = \mathfrak{g}_\tau^\perp \oplus \mathfrak{g}_{\tau}^{\perp \flat}
\end{equation}
with Euclidean coordinates defined using orthonormal bases.  We refer to \cite[\S9.4.1]{2020arXiv201202187N} for a detailed discussion and illustration of such coordinates, which may aid the reader in following the arguments below.  We recall in particular Figure~\ref{fig:tau-coordinates}, taken from \emph{loc.\ cit.}.
\setlength{\unitlength}{1.5cm}
\begin{figure}
  \begin{picture}(4,3)
    \put(-1,0){\vector(1,0){6}}
    \put(-1,0){\vector(0,1){1.5}}
    {%
      \thicklines
      \color{black}%
      \multiput(2.1,0.2)(0,0.1){10}{\line(0,1){0.05}}
      \put(2.1,1.3){$\tau + \mathfrak{g}_\tau^{\perp \flat}$}
    }
    {%
      \thicklines
      \color{black}%
      \multiput(0,0.5)(0.1,0){40}{\line(1,0){0.05}}
      \put(4.2,0.4){$\tau + \mathfrak{g}_\tau^\perp = \tau + [\mathfrak{g},\tau]$
      }
      {%
        \thicklines
        \color{black}%
        \qbezier(0,1)(2,0)(4,1)
        \put(-0.6,1.1){$G \cdot \tau$}
      }

      \color{black}
      \put(-1.4,1.4){$\xi''$}
      \put(4.9,-0.25){$\xi'$}
      \put(2,0.3){$\tau$}
      \put(1.95,0.5){\circle*{0.1}}
    }
  \end{picture}
  \caption{ The coadjoint orbit $G \cdot \tau$ near $\tau$ in $\tau$-coordinates.  } \label{fig:tau-coordinates}
\end{figure}

We consider the map
\[
  \gamma : \mathfrak{g}^\wedge \rightarrow \mathfrak{g}_\tau^\perp \times [\mathfrak{g}^\wedge]\]
\[
  \gamma(\xi) := (\xi ', [\xi]).
\]
For instance, in the above figure, this map resembles ``$(x, y) \mapsto(x, y - x^2)$''.  It yields a coordinate chart near $\tau$.  Indeed, since the derivative at $\tau$ of $\mathfrak{g}^\wedge \rightarrow [\mathfrak{g}^\wedge]$ is surjective (see \cite[Thm 0.1]{MR0158024} or \cite[Theorem 13.5]{2020arXiv201202187N}) and has kernel $\mathfrak{g}_\tau^\perp$, we see that the derivative $(d \gamma)_\tau$ at $\tau$ of $\gamma$ is given with respect to the decomposition \eqref{eqn:g-wedge-tau-decmop-std} by a block diagonal matrix $\left(
  \begin{smallmatrix}
    1&0\\
    0 &M_\tau
  \end{smallmatrix}
\right)$ consisting of the identity transformation $1$ of $\mathfrak{g}_\tau^\perp$ and an isomorphism
\[M_\tau : \mathfrak{g}_\tau^{\perp \flat} \rightarrow T_{[\tau]}([\mathfrak{g}^\wedge]) \cong \mathbb{R}^n.\]

The operator norms of $M_\tau$ and $M_\tau^{-1}$ are $\O(1)$ for $\tau$ as indicated (indeed, both $M_\tau$ and $M_\tau^{-1}$ depend continuously on $\tau$, so by \eqref{eq:cuhnuoy8vh}, they lie in some fixed compact set of linear transformations).  For $\xi = \tau + o(1)$, the derivative $(d \gamma)_\xi$ varies smoothly with $\xi$, so by the mean value theorem, we have the operator norm estimate
\begin{equation}\label{eqn:estimate-d-gamma-xi-near-tau}
  (d \gamma)_{\xi}
  = \left(
    \begin{smallmatrix}
      1&\\
       &M_\tau
    \end{smallmatrix}
  \right) + \O(|\xi - \tau|).
\end{equation}
In particular, $(d \gamma)_{\xi}$ and its inverse have operator norms $\O(1)$.  Since the second partial derivatives of $\gamma$ on arguments $\xi = \tau + o(1)$ are of size $\O(1)$, we deduce that for $r \lll 1$, the restriction of $\gamma$ to the ball $\{\xi : |\xi - \tau| < r\}$ defines a diffeomorphism onto its image that approximately preserves distances in the sense that $|\xi-\eta| \asymp \dist(\gamma(\xi), \gamma(\eta))$.  (Alternatively, we could have appealed here to Lemma \ref{lem:controlled-uniform-continuity-origin}, noting that the map $\gamma$ is controlled in view of our hypothesis \eqref{eq:cuhnuoy8vh} on $\tau$).

For any smooth function $a : \mathfrak{g}^\wedge \rightarrow \mathbb{C}$ supported on $\tau + o(1)$, we obtain a smooth function
\[
  \tilde{a} : \mathfrak{g}_\tau^\perp \times [\mathfrak{g}^\wedge] \rightarrow \mathbb{C},
\]
supported on $\gamma(\tau) + o(1)$ and given there by
\begin{equation}\label{eq:tild-:=-axi}
  \tilde{a}(\gamma(\xi)) := a(\xi).
\end{equation}

Let $\tilde{S}^{\tau}_{\delta ', \delta ''}$ denote the class consisting of smooth functions $\tilde{a}$ as above that are supported on $\gamma(\tau) + \O(\h^{ \delta '' - \delta'})$ and satisfy the derivative bounds
\begin{equation}\label{eqn:tilde-a-deriv-bounds}
  \partial_{\xi '}^{\alpha'}
  \partial_{\lambda}^{\alpha ''}
  \tilde{a}(\xi ', \lambda )
  \ll
  \h^{- \delta ' |\alpha '|
    - \delta ''| \alpha ''|}.
\end{equation}

\begin{lemma}
  Let $a : \mathfrak{g}^\wedge \rightarrow \mathbb{C}$ and $\tilde{a} : \mathfrak{g}_\tau^{\perp} \times [\mathfrak{g}^\wedge] \rightarrow \mathbb{C}$ be smooth functions, supported respectively on $\tau + o(1)$ and $\gamma(\tau) + o(1)$, and related as in \eqref{eq:tild-:=-axi}.  The following conditions are equivalent.
  \begin{enumerate}[(i)]
  \item $a \in S_{\delta ', \delta ''}^{\tau}$.
  \item $\tilde{a} \in \tilde{S}^{\tau}_{\delta ', \delta ''}$.
  \end{enumerate}
\end{lemma}
\begin{proof}
  We follow the proof of \cite[Lemma 9.9]{2020arXiv201202187N} line-by-line.  The crucial estimate is \eqref{eqn:estimate-d-gamma-xi-near-tau}.  
\end{proof}
We note also that a smooth function $f : \mathfrak{g}^\wedge \rightarrow \mathbb{C}$ is a $G$-invariant polynomial if and only if the corresponding function $\tilde{f}$, defined via \eqref{eq:tild-:=-axi}, is given by a polynomial function of the second coordinate.

With these preliminaries settled, we now proceed:
\begin{proof}[Proof of Proposition \ref{lem:produce-inv-polynomials-large-on-symbol-support}]
  Let $\xi$ be an element of the support of $a$.  We consider the second order Taylor expansion of $[\xi]$ based at $\tau$: 
  \[
    [\xi] = [\tau] + M_\tau (\xi '' - \tau '') + \O(|\xi - \tau |^2).
  \]
  Since $|\xi - \tau |^2 \ll \h^{2 \delta '} \lll \h^{\delta ''}$ and $|\xi '' - \tau ''| \gg \h^{\delta ''}$, we deduce from the noted operator norm bound $M_\tau^{-1} = \O(1)$ that
  \[
    |[\xi] - [\tau]| \asymp |\xi '' - \tau ''| \gg \h^{ \delta ''}.
  \]
  Moreover, since $\xi = \tau + \O(\h^{\delta '})$, we have
  \[
    \gamma(\xi) = \gamma(\tau) + \O(\h^{\delta '}).
  \]

  Define \(\tilde{a} \in \tilde{S}^{\tau}_{\delta ', \delta ''}\) as above.  Then
  \begin{equation}\label{eq:cuhn9dsj1l}
    \text{\(\tilde{a}(\xi ', \lambda)\) is supported on \(\lambda - [\tau] \gg \h^{ \delta ''}\)}.
  \end{equation}

  Writing \(\lambda = (\lambda_1,\dotsc,\lambda_n)\) and \([\tau] = ([\tau]_1,\dotsc,[\tau]_n)\) with respect to the identification of \([\mathfrak{g}^\wedge]\) with \(\mathbb{R}^n\), we claim that we may smoothly decompose
  \begin{equation}\label{eq:cuhn9ar8ty}
    \tilde{a} = \tilde{a}_1 + \dotsb + \tilde{a}_n,
  \end{equation}
  where each \(\tilde{a}_j \in \tilde{S}^{\tau}_{\delta ', \delta ''}\) is supported on \(|\lambda_j - [\tau]_j| \gg \h^{\delta ''}\).  To see this, first use \eqref{eq:cuhn9dsj1l} to fix \(c>0\) so that
  \begin{equation}\label{eq:cuhn9dzoi5}
    \supp(\tilde{a}) \subseteq \{(\xi ',\lambda) : \max_{1 \leq j \leq n} |\lambda_j - [\tau]_j| \geq c \h^{\delta ''}\}.
  \end{equation}
  Fix a smooth partition of unity \(1=\phi_0+\phi_1\) of \(\mathbb{R}\), with \(\phi_0\) supported on \([-c,c]\) and \(\phi_1\) supported on \(\{x : |x| \geq c/2\}\), and define a smooth partition of unity $1 = \Phi_0 + \Phi_1 + \dotsb + \Phi_n$ of \(\mathbb{R}^n\) by
  \[
    \Phi_0(x) := \prod_{k=1}^n \phi_0(x_k), \qquad
    \Phi_j(x) := \phi_1(x_j) \prod_{k<j} \phi_0(x_k) \quad (1 \leq j \leq n).
  \]
  On the region \(\max_j |x_j| \geq c\), we have \(\Phi_0(x)=0\) and hence \(\sum_{j=1}^n \Phi_j(x)=1\).  For \(1 \leq j \leq n\), set
  \[
    \tilde{a}_j(\xi ',\lambda) := \tilde{a}(\xi ',\lambda)\,\Phi_j\!\left(\frac{\lambda-[\tau]}{\h^{\delta ''}}\right).
  \]
  By the support condition \eqref{eq:cuhn9dzoi5}, we see that the claimed decomposition \eqref{eq:cuhn9ar8ty} holds.  Moreover, \(\tilde{a}_j\) is supported where \(|\lambda_j-[\tau]_j| \geq (c/2)\h^{\delta ''} \gg \h^{\delta ''}\).  Finally, for each multi-index \(\alpha''\), we have
  \[
    \partial_\lambda^{\alpha''}\Phi_j\!\left(\frac{\lambda-[\tau]}{\h^{\delta ''}}\right) \ll \h^{-\delta ''|\alpha''|},
  \]
  so by the product rule for differentiation and \eqref{eqn:tilde-a-deriv-bounds}, we see that \(\tilde{a}_j\) satisfies the same bounds \eqref{eqn:tilde-a-deriv-bounds}, hence \(\tilde{a}_j \in \tilde{S}^{\tau}_{\delta ', \delta ''}\).

  The required conclusion holds with $a_j$ the element of $S^{\tau}_{\delta ', \delta ''}$ corresponding to $\tilde{a}_j$ and $f_j(\xi) := [\xi]_j - [\tau]_j$.
\end{proof}

\subsection{Harish--Chandra vs.\ symmetrization}\label{sec:harish-chandra-vs}
Consider the composition
\begin{equation}\label{eqn:HACH-sym}
  \Sym(\mathfrak{g}_\mathbb{C})^G
  \xrightarrow{\sym} \mathfrak{Z}(\mathfrak{g}_\mathbb{C})
  \xrightarrow{\HACH}
  \Sym(\mathfrak{g}_\mathbb{C})^G.
\end{equation}
Recall that $\sym$ is a vector space isomorphism and $\HACH$ is an algebra isomorphism.

This composition is typically not the identity map, but converges to the identity in the ``classical limit,'' e.g., in the sense of the lemma recorded below.

Let $f \in S^\infty_0$ be a polynomial symbol.  Let $f'$ denote the polynomial for which $f_{\h}' := [\xi \mapsto f'(\h \xi)]$ is the inverse image of $f_{\h} := [\xi \mapsto f(\h \xi)]$ under the above composition.  Thus $\HACH(\sym(f_{\h}')) = f_{\h}$.  By \cite[\S5.2]{nelson-venkatesh-1}, we have
\[
  \Opp(f') = f(\h \lambda_\pi),
\]
i.e., $\Opp(f')$ is the scalar operator on $\pi$ given by multiplication by $f(\h \lambda_\pi)$.
\begin{lemma}
  We have $f' \in S^\infty_0$ and $f ' - f \in  \h S^\infty_0$.
\end{lemma}
\begin{proof}
  As noted in \cite[\S9.9]{nelson-venkatesh-1}, the composition \eqref{eqn:HACH-sym} has the property that if $p \in \Sym(\mathfrak{g}_\mathbb{C})^G$ has order $\leq n$, then $\HACH(\sym(p)) - p$ has order $\leq n-1$.  The conclusion follows readily from this (see \cite[\S10.3]{nelson-venkatesh-1} for further details).
\end{proof}

\subsection{Division by an invariant symbol}
The technique of approximately dividing one symbol by another is pervasive in microlocal analysis.  Some applications of that technique were recorded in \cite[\S10]{nelson-venkatesh-1}.  We record here a further application.
\begin{proposition}\label{prop:divide-by-G-inv-symb}
  Let $a \in S_{\delta ', \delta ''}^{\tau}$.  Let $f \in S^\infty_0$ be $G$-invariant.  Suppose that $f(\xi) \gg \h^{\delta ''}$ for all $\xi \in \supp(a)$.  Then there exist $b, b' \in \h^{-\delta ''} S^\tau_{\delta', \delta ''}$ so that
  \[
    a \equiv f \star_{\h} b \mod{\h^\infty S^{-\infty}},
  \]
  \[
    a \equiv b' \star_{\h} f \mod{\h^\infty S^{-\infty}}.
  \]
\end{proposition}
The proof requires a few preliminaries.

Working in $\tau$-coordinates $\xi = (\xi ', \xi '')$, we write $\partial_1^{\alpha '} \partial_2^{\alpha ''}$ for the differential operator given by $\alpha '$ and $\alpha ''$ partial derivatives in the $\xi'$ and $\xi ''$ variables, respectively.
\begin{lemma}\label{lem:derivatives-along-orbits-inv-fn}
  Let $\tau \in \mathfrak{g}^\wedge_{\reg}$.  Let $f \in C^\infty(\mathfrak{g}^\wedge)$ be $G$-invariant.  Then in $\tau$-coordinates, we have
  \[
    |\alpha '| =1 \quad \implies \quad \partial_{1}^{\alpha '} f(\tau) = 0.
  \]
\end{lemma}
\begin{proof}
  This follows from the fact that $\mathfrak{g}_\tau^\perp$ is the tangent plane at $\tau$ to the coadjoint orbit containing $\tau$ and that $f$ is constant on each such orbit.
\end{proof}

\begin{lemma}\label{lem:star-j-against-invariant-function}
  Let $f \in S_0^\infty$ be $G$-invariant.  Let $a \in S^\tau_{\delta ', \delta ''}$.
  \begin{enumerate}[(i)]
  \item We have $f \star^1 a = 0$.
  \item For fixed $j \in \mathbb{Z}_{\geq 2}$, we have $f \star^j a \in \h^{\delta '' - 2 j \delta '} S_{\delta ', \delta ''}^\tau$.
  \end{enumerate}
  Analogous results hold for $a \star^j f$.
\end{lemma}
\begin{proof}
  (i): The restriction of $f \star^1 a$ to a coadjoint orbit is a multiple of the Poisson bracket of the restrictions of $f$ and $a$ (cf.\ \cite{gutt1983explicit}).  Since $f$ is $G$-invariant, its restriction to each such orbit is constant, and so the Poisson bracket in question vanishes.  This completes the proof.

  As motivation for the proof of (ii), we pause to reprove (i) using ``bare hands''.

  We must show that $f \star^1 a(\xi) = 0$ for all $\xi \in \mathfrak{g}^\wedge$.

  We first reduce to the case $\xi = \tau$, as follows.  By hypothesis, the symbol $a$ and hence also $f \star^1 a$ is supported on $\tau + \O(\h^{\delta '' - \delta '})$, so we may reduce to the case that $\xi$ is of this form.  By Lemma \ref{lem:basepoint-insensitivity}, we then have $S_{\delta ', \delta''}^\tau = S_{\delta ', \delta''}^\xi$.  Since in particular $\xi = \tau + o(1)$, our assumption that $\tau$ lies in a fixed compact subset of $\mathfrak{g}^\wedge_{\reg}$ likewise holds for $\xi$.  Our hypotheses thus hold with $\tau$ replaced by $\xi$, so we obtain the required reduction after renaming $\xi$ to $\tau$.

  We must now verify that the quantity $f \star^1 a(\tau)$ vanishes.  As discussed in Example \ref{exa:star1-structure}, we may write this quantity as a finite linear combination of terms
  \[
    \left( \partial_1^{\alpha '} f(\tau) \right) \left( \partial_1^{\beta '} a(\tau) \right) \left( \tau^\gamma \right), \quad |\alpha '| = |\beta '| = |\gamma| = 1,
  \]
  where the differential operators are defined using $\tau$-coordinates.  By Lemma \ref{lem:derivatives-along-orbits-inv-fn}, each of these terms vanish.  This completes our motivating ``bare hands'' reproof of (i).

  (ii): By Lemma \ref{lem:basepoint-insensitivity} and arguments similar to those above, we reduce to verifying that $f \star^j a$, for fixed $j \geq 2$, satisfies the required derivative bounds at $\tau$.  By the support condition \eqref{eq:alpha-refined-support-condition} for the coefficients defining the star product and the product rule for differentiation, we may write $\partial_1^{\kappa '} \partial_2^{\kappa ''} (f \star^j a)(\tau)$ as a finite linear combination of terms
  \[
    \left( \partial_{1}^{\alpha ' + \mu '} \partial_2^{\alpha '' + \mu ''} f(\tau) \right) \left( \partial_{1}^{\beta ' + \nu '} \partial_2^{\beta '' + \nu ''} a(\tau) \right) \left( \partial_1^{\lambda '} \partial_2^{\lambda ''} \tau^{\gamma} \right)
  \]
  where
  \[
    \lambda ' + \mu ' + \nu ' = \kappa ', \quad \lambda '' + \mu '' + \nu '' = \kappa ''
  \]
  and
  \begin{equation}\label{eqn:key-alpha-beta-ineqs}
    |\alpha '| + |\beta '|
    +
    2 |\alpha ''| + 2 |\beta ''|
    \leq 2 j,
    \quad
    |\alpha '| + |\alpha ''| \geq 1,
    \quad
    |\beta '| + |\beta ''| \geq 1.
  \end{equation}
  We must check that any such term is
  \[
    \ll \h^{- \delta ' |\kappa '| - \delta '' | \kappa ''| + \delta '' - 2 j \delta '}.
  \]
  Since $a \in S^\tau_{\delta ', \delta ''}$, we have
  \begin{equation}\label{eqn:estimates-for-a-for-f-star-a}
    \partial_{1}^{\beta ' + \nu '} \partial_2^{\beta '' + \nu ''} a(\tau)
    \ll
    \h^{- \delta '( |\beta '| +  |\nu '|) - \delta ''(|\beta ''| + |\nu ''|)},
  \end{equation}
  so it is enough to check that
  \begin{equation}\label{eqn:f-star-a-desiderata-f}
    \partial_{1}^{\alpha ' + \mu '} \partial_2^{\alpha '' + \mu ''} f(\tau)
    \h^{ - \delta ' |\beta '|-  \delta ''|\beta ''|}
    \ll
    \h^{- \delta ' |\mu '| - \delta '' | \mu ''| + \delta '' - 2 j \delta '}.
  \end{equation}

  The function $f$ enjoys the estimates
  \begin{equation}\label{eqn:f-has-mild-alpha-mu-derivs}
    \partial_{1}^{\alpha ' + \mu '} \partial_2^{\alpha '' + \mu ''} f(\tau)
    \ll 1
  \end{equation}
  and, by Lemma \ref{lem:derivatives-along-orbits-inv-fn},
  \begin{equation}\label{eqn:f-star-a-f-vanishing-condition}
    \partial_{1}^{\alpha ' + \mu '} \partial_2^{\alpha '' + \mu ''} f(\tau)
    = 0 \text{ if } |\alpha ' + \mu '| = 1, \, |\alpha '' + \mu ''| = 0.
  \end{equation}
  In particular, the analogue of \eqref{eqn:estimates-for-a-for-f-star-a} for $f$
  \begin{equation}\label{eqn:estimates-for-a-for-f-star-f}
    \partial_{1}^{\alpha ' + \mu '} \partial_2^{\alpha '' + \mu ''} f(\tau)
    \ll
    \h^{- \delta '( |\alpha '| +  |\mu '|) - \delta ''(|\alpha ''| + |\mu ''|)}
  \end{equation}
  is valid, but typically much weaker than reality.  It is nevertheless instructive to note that this last estimate quite nearly yields \eqref{eqn:f-star-a-desiderata-f}: it gives that the LHS of \eqref{eqn:f-star-a-desiderata-f} is
  \[
    \h^{- \delta ' |\mu '| - \delta '' | \mu ''| - \delta ' |\alpha '| - \delta ' |\beta '| - \delta '' | \alpha '' | - \delta '' |\beta ''| },
  \]
  which, thanks to \eqref{eqn:key-alpha-beta-ineqs}, is merely a factor $\h^{\delta ''}$ shy of the required estimate \eqref{eqn:f-star-a-desiderata-f}.  It will thus suffice to improve the estimate \eqref{eqn:estimates-for-a-for-f-star-f} by at least a factor of $\h^{\delta ''}$.

  If $|\alpha '' + \mu ''| \geq 1$, then the required improvement is immediate from \eqref{eqn:f-has-mild-alpha-mu-derivs}, so we may suppose that $\alpha '' = \mu '' = 0$.  If $|\alpha ' + \mu '| \geq 2$, then we may again apply \eqref{eqn:f-has-mild-alpha-mu-derivs} to improve \eqref{eqn:estimates-for-a-for-f-star-f} by a factor $\h^{2 \delta '} \leq \h^{\delta ''}$.  If $|\alpha ' + \mu '| = 1$, then the LHS of \eqref{eqn:estimates-for-a-for-f-star-f} vanishes in view of \eqref{eqn:f-star-a-f-vanishing-condition}.  The remaining case $\alpha ' = \mu ' = 0$ does not arise due to the inequality $|\alpha '| + |\alpha ''| \geq 1$ noted in \eqref{eqn:key-alpha-beta-ineqs}.  The proof is now complete.
\end{proof}

\begin{lemma}\label{lem:b0-in-h-neg-delta-dubprime}
  Under the hypotheses of Proposition \ref{prop:divide-by-G-inv-symb}, we have $a/f \in \h^{- \delta ''} S_{\delta', \delta ''}^{\tau}$.
\end{lemma}
\begin{proof}
  The support condition is clear.  We must check the derivative bounds.  By Lemma \ref{lem:basepoint-insensitivity} and arguments similar to those in the proof of Lemma \ref{lem:star-j-against-invariant-function}, it will suffice to check those bounds at the point $\tau$.  We may assume that $\tau \in \supp(a)$, so that, by hypothesis, $f(\tau) \gg \h^{\delta ''}$.  Since $a \in S^\tau_{\delta ', \delta ''}$, we may reduce via the product rule for derivatives to verifying that $1/f$ satisfies, at the point $\tau$, the derivative bounds
  \begin{equation}\label{eqn:required-partials-for-1-over-f}
    \partial_{1}^{\alpha '}
    \partial_{2}^{\alpha ''} \frac{1}{f}(\tau) \ll
    \h^{- \delta '' - \delta ' |\alpha'| - \delta '' |\alpha ''|}.
  \end{equation}

  To that end, observe first for a smooth function $\phi$ of one variable, the $n$th derivative $(1/\phi)^{(n)}$ of $1/\phi$ may be written using the quotient rule for derivatives as a linear combination of terms
  \[
    \frac{1}{\phi } \left( \frac{\phi^{(1)} }{\phi } \right)^{e_1} \left( \frac{\phi^{(2)} }{\phi } \right)^{e_2} \dotsb \left( \frac{\phi^{(n)} }{\phi } \right)^{e_n}
  \]
  where the exponents $e_1,\dotsc,e_n$ are nonnegative integers satisfying $\sum j e_j = n$.  Note in particular that if at some point the first derivative $\phi^{(1)}$ vanishes but $\phi$ does not vanish, then the above expression vanishes unless $e_1 = 0$, in which case $\sum e_j \leq n/2$, and so the total number of denominator factors $\phi$ is at most $1 + n/2$ .

  Similar observations apply to partial derivatives in higher dimensions.  Writing $m := |\alpha '|$ and $n := |\alpha ''|$, we may express the LHS of \eqref{eqn:required-partials-for-1-over-f} as
  a linear combination of terms
  \begin{equation}\label{eqn:terms-for-partial-derivs-of-f}
    \frac{1}{f}
    \prod_{i=0}^{m}
    \prod_{j=0}^{n}
    \frac{
      (\partial_1^i \partial_2^j f)^{a_{i j}}
    }{f^{a_{i j}}}
    (\tau),
  \end{equation}
  with notation as follows.
  \begin{itemize}
  \item $\partial_i^j$ is shorthand for an operator of the form $\partial_i^{\beta}$ with $|\beta| = j$, not necessarily with the same $\beta$ in each occurrence.  For instance, if $a_{ij } = 3$, then $(\partial_1^i \partial_2^j f)^{a_{i j}}$ denotes some product of the form $(\partial_1^{\beta} \partial_2^{\gamma } f) (\partial_1^{\beta'} \partial_2^{\gamma'} f) (\partial_1^{\beta''} \partial_2^{\gamma'' } f)$ with $|\beta| = |\beta '| = |\beta ''| = i$ and $|\gamma| = |\gamma '| = |\gamma ''| = j$.
  \item The $a_{i j}$ are nonnegative integers satisfying
    \[
      \sum_{i} \sum_j i a_{i j} = m, \quad \sum_{i} \sum_j j a_{i j} = n.
    \]
  \end{itemize}
  We may and shall assume that $a_{0 0} = 0$, since the corresponding factor is $1$.  By Lemma \ref{lem:derivatives-along-orbits-inv-fn}, we have $\partial_1^1 \partial_2^0 f (\tau) = 0$, so we may assume moreover that $a_{1 0} = 0$.  We then rewrite \eqref{eqn:terms-for-partial-derivs-of-f} as $\frac{1}{f}(\tau) P Q$, where $P$ and $Q$ denote the respective contributions from $j = 0$ and from $j \geq 1$.  Since $f(\tau) \gg \h^{\delta ''}$, our goal bound \eqref{eqn:required-partials-for-1-over-f} will follow if we can show that
  \[
    P \ll \h^{- \delta ' m}, \quad Q \ll \h^{ - \delta '' n}.
  \]
  Using the estimates $\partial^\beta f(\tau) \ll 1$ and $f(\tau) \gg \h^{\delta ''}$, we see that
  \[
    P \ll \h^{- \sum_i a_{i 0} \delta ''}, \quad Q \ll \h^{- \sum_i \sum_{j \geq 1} a_{i j} \delta ''}.
  \]
  By the stated vanishing conditions on the $a_{i j}$, we have
  \[
    \sum_i a_{i 0} \leq \frac{1}{2} \sum_i i a_{i 0} \leq \frac{m}{2}, \quad \sum_i \sum_{j \geq 1} a_{i j} \leq \sum_i \sum_j j a_{i j} = n.
  \]
  Our goal bound follows upon recalling that $\delta ''/2 \leq \delta '$.
\end{proof}

We turn finally to the proof of Proposition \ref{prop:divide-by-G-inv-symb}.  Retain its hypotheses.  We will construct $b$; the construction of $b'$ is similar.  We define $b_j$ ($j \in \mathbb{Z}_{\geq 0}$) inductively by requiring that $f \star_{\h} \sum_{j \geq 0} \h^j b_j \sim a$ in a formal sense.  Thus
\[
  b_0 := \frac{a}{f}, \quad b_1 := \frac{- f \star^1 b_0}{f}, \quad b_2 := \frac{- f \star^2 b_0 + f \star^1 b_1}{f},
\]
and so on.
\begin{lemma}
  For fixed $j$, we have $b_j \in \h^{-\delta '' - 2 j \delta '} S^\tau_{\delta ', \delta ''}$.
\end{lemma}
\begin{proof}
  We induct on $j$.  The case $j=0$ is the content of Lemma \ref{lem:b0-in-h-neg-delta-dubprime}, so suppose $j \geq 1$.  By our inductive hypothesis and Lemma \ref{lem:star-j-against-invariant-function}, we see that for $\ell \in \{1, \dotsc, j\}$, we have $f \star^{\ell} b_{j - \ell} \in \h^{- 2 j \delta ' } S^\tau_{\delta ', \delta''}$.  By Lemma \ref{lem:b0-in-h-neg-delta-dubprime} (applied with $a$ replaced by $\h^{2 j \delta '} f \star^{\ell} b_{j - \ell}$, which satisfies the required hypotheses), we deduce that
  \[
    \frac{f \star^{\ell} b_{j - \ell}}{f} \in \h^{- \delta '' - 2 j \delta '} S^\tau_{\delta ', \delta ''}.
  \]
  By construction, $b_j$ is an alternating sum of such terms.  The proof is complete.
\end{proof}

We define
\[
  b := \sum_{0 \leq j < L} \h^j b_j
\]
for some $L \ggg 1$ as in the conclusion of Lemma \ref{lem:borel-summation}.  Then $b \in \h^{- \delta ''} S^\tau_{\delta ', \delta ''}$, and the required approximation $a \equiv f \star_{\h} b \mod{\h^\infty S^{-\infty}}$ holds.  The proof of Proposition \ref{prop:divide-by-G-inv-symb} is now complete.

\subsection{Completion of the proof}\label{sec:completion-proof}
We finally prove Theorem \ref{thm:local-refin-symb}.  We retain its setting.  We decompose $a = a_1 + \dotsb + a_n$ and construct $f_1,\dotsc,f_n$ according to Proposition \ref{lem:produce-inv-polynomials-large-on-symbol-support}.  Let $f_j'\in S^\infty_0$ be attached to $f_j$ as in \S\ref{sec:harish-chandra-vs}.  Then $f_j' - f_j \in \h S^\infty_0$, so for $\xi \in \supp(a_j)$, we have $f_j'(\xi) = f_j(\xi) + \O(\h)$.  In particular, since $\delta '' < 1$, the lower bound $f_j'(\xi) \gg \h^{\delta ''}$ persists for $f_j'$.  Using Proposition \ref{prop:divide-by-G-inv-symb}, we may thus produce $b_1,\dotsc,b_n \in \h^{- \delta ''} S^\tau_{\delta ', \delta ''}$ so that $a_j \equiv f_j' \star_{\h} b_j \mod{\h^\infty S^{-\infty}}$.  By Lemma \ref{lem:let-chi-be} and the composition formula (Theorem \ref{thm:composition-formula-quantitative}), we have
\[
  \Opp(a_j) \equiv \Opp(f_j') \Opp(b_j) \mod{\h^\infty \Psi^{-\infty}}.
\]
Since $[\tau] = \h \lambda_\pi$, we have
\[
  \Opp(f_j') = f_j(\h \lambda_\pi) = f_j(\tau) = 0.
\]
Thus each $\Opp(a_j)$ and hence also $\Opp(a)$ belongs to $\h^\infty \Psi^{-\infty}$, as required.

\section{Localization with respect to $P$}\label{sec:cuhn9f8y8q}
For the remainder of Part \ref{part:asympt-analys-kirill}, we specialize to
\begin{equation*}
  G := \GL_n(\mathbb{R}),
\end{equation*}
with accompanying notation $H, N, N_H, P, A$ as in \S\ref{sec:asympt-analys-local}.

\subsection{Qualitative preliminaries}\label{sec:qual-prel}

Let $\psi$ be a nondegenerate unitary character of $N$; as usual, denote also by $\psi$ its restriction to $N_H$.  Let $\pi$ be a generic irreducible unitary representation of $G$.  We identify $\pi$ with its Whittaker model:
\begin{equation*}
  \pi = \mathcal{W}(\pi,\psi).
\end{equation*}
Recall from \S\ref{sec:basic-definition} that we denote by $\pi^\infty \subseteq \pi^0 \subseteq \pi^{-\infty}$ the respective spaces of smooth, norm-bounded and distributional vectors for $\pi$.  We realize $\pi$ in its Kirillov model, thereby identifying its underlying Hilbert space $\pi^0$ with $L^2(N_H \backslash H, \psi)$.  We recall that $\Theta_{\pi,\psi} \in \pi^{-\infty}$ denotes the distributional vector with the property that for all $W \in \mathcal{W}(\pi,\psi)$, we have $\langle W, \Theta_{\pi,\psi} \rangle = W(1)$.

The smooth and distributional parts of the restrictions of $\pi$ to the subgroups $P$ and $H$ are related via natural inclusions
\begin{equation*}
  \pi^\infty \subseteq \pi|_P^\infty \subseteq \pi|_H^\infty \subseteq \pi^0  \subseteq \pi|_H^{-\infty}
  \subseteq \pi|_P^{-\infty} \subseteq \pi^{-\infty}.
\end{equation*}
\begin{lemma}
  We have $\Theta_{\pi,\psi} \in \pi|_H^{-\infty}$.
\end{lemma}
\begin{proof}
  The Sobolev lemma for $(N_H \backslash H,\psi)$ implies that we may find $d \in \mathbb{Z}_{\geq 0}$ and $C_0 \geq 0$ so that for all $f \in C^\infty(N_H \backslash H, \psi)$, we have
  \begin{equation*}
    |f(1)| \leq C_0 \sum_{x_1, \dotsc, x_d \in \mathcal{B}(H) \cup \{1\}} \|R(x_1 \dotsb x_d) f\|_{L^2(N_H \backslash H)}.
  \end{equation*}
  In particular, for $W \in \pi$,
  \begin{equation*}
    |\langle W, \Theta_{\pi,\psi} \rangle| \leq C_0  \sum_{x_1, \dotsc, x_d \in \mathcal{B}(H) \cup \{1\}} \|\pi(x_1 \dotsb x_d) W\|.
  \end{equation*}
  This implies that $\Theta_{\pi,\psi}$ lies in $\pi|_H^{-\infty}$ (or more specifically, in the negative index Sobolev space $\pi|_H^{-d}$ as defined in \cite[\S3.2]{nelson-venkatesh-1}).
\end{proof}

Recall that $\h \in (0,1]$ denotes a ``wavelength parameter.''

\begin{lemma}
  For each cutoff $\chi$ and $a \in \underline{S}^{-\infty}(\mathfrak{p}^\wedge)$, we have
  \begin{equation*}
    \Opp_{\h}(a:\chi) \Theta_{\pi,\psi} \in \pi^\infty.
  \end{equation*}
\end{lemma}
\begin{proof}
  By definition, $\Opp_{\h}(a:\chi) = \pi(\omega)$ with $\omega = \widetilde{\Opp}_{\h}(a:\chi)$ a smooth compactly-supported measure on $P$, or, after dividing by some fixed (left) Haar measure, a smooth compactly-supported function.  The conclusion is thus a special case of Jacquet's results (Lemma \ref{thm:jacquet-on-kirillov-model-P}), which give $\pi(\omega) \Theta_{\pi,\psi} \in \pi^\infty$ for all smooth compactly-supported measures $\omega$ on $P$.
\end{proof}

\begin{lemma}\label{lem:action-opp-on-cal-W-delta}
  Suppose that $\chi_1, \chi_2, \chi '$ are cutoffs with $\chi '(x \ast y) = 1$ whenever $\chi_1(x) \chi_2(y) \neq 0$.  Then for all
  \begin{equation*}
    (a,b) \in \underline{S} ^\infty (\mathfrak{p}^\wedge ) \times \underline{S} ^{\infty } (\mathfrak{p} ^\wedge ),
  \end{equation*}
  we have
  \begin{equation*}
    \Opp_{\h}(a:\chi_1)
    \Opp_{\h}(b:\chi_2)  \Theta_{\pi,\psi}
    =
    \Opp_{\h}(a \star_{\h} b:\chi ') \Theta_{\pi,\psi},
  \end{equation*}
  with the star product $\star_{\h}$ defined relative to $\chi_1$ and $\chi_2$ as in \S\ref{sec:star-product}.
\end{lemma}
\begin{proof}
  We evaluate the general composition formula (Theorem \ref{thm:qualitative-extended-operator-assignment}) on $\Theta_{\pi,\psi}$.
\end{proof}

\subsection{Quantitative preliminaries}\label{sec:quant-prel}

\subsubsection{Scaling parameters}
We now
set
\begin{equation*}
  T := 1/\h
\end{equation*}
and normalize the wavelength parameter $\h \lll 1$ by requiring that the rescaled infinitesimal character $\h \lambda_\pi = T^{-1} \lambda_\pi$ belongs to a fixed compact set (as we may, in the setting of Theorem \ref{thm:big-group-whittaker-behavior-under-G}).

\subsubsection{Characters}
We take $\psi$ fixed, and define the rescaled character $\psi_T$ as in \S\ref{sec:class-bump-functions}.  The discussion of \S\ref{sec:qual-prel} then applies, with $\psi_T$ playing the role of ``$\psi$.''

For the remainder of Part \ref{part:asympt-analys-kirill}, we identity $\pi$ with its $\psi_T$-Whittaker model:
\begin{equation*}
  \pi = \mathcal{W}(\pi, \psi_T).
\end{equation*}

\subsubsection{A class of localized Whittaker functions}
Recall (Definition \ref{defn:we-denote-mathfr-frakW}) the classes $\mathfrak{W}(\pi,\psi,T,\delta) \subseteq \mathcal{W}(\pi,\psi_T)$.  We abbreviate  \index{classes!bump functions in the Kirillov model, $\mathcal{W}_\delta$}
\begin{equation*}
  \mathcal{W}_\delta := \mathfrak{W}(\pi, \psi, T, \delta).
\end{equation*}
Informally, $\mathcal{W}_\delta$ consists of vectors given in the $\psi_T$-Kirillov model by smooth bumps on $1 + \O(\h^{1-\delta})$.

\subsubsection{The parameter $\theta$}
In what follows, we abbreviate simply by $\theta$ the quantities $\theta_P(\psi)$ or $\theta_H(\psi)$, defined as in \S\ref{sec:class-bump-functions}; it should be clear from context to which element we are referring.

\subsubsection{A class of $\theta$-localized symbols}
We abbreviate, with notation as in Definition \ref{defn:let-mathfrakav-theta},   \index{symbols!$\theta$-localized, $\mathcal{A}_\delta(\mathfrak{p}^\wedge)$}
\begin{equation*}
  \mathcal{A}_\delta(\mathfrak{p}^\wedge) := \mathfrak{A}(\mathfrak{p}^\wedge, -\theta, T, \delta),
\end{equation*}
Expanding that definition, $\mathcal{A}_\delta(\mathfrak{p}^\wedge)$ is the subclass of $\underline{S}^{-\infty}(\mathfrak{p}^\wedge)$ consisting of symbols $a$ satisfying
\begin{equation}\label{eq:partialalpha-axi-ll}
  \partial^\alpha a(\theta + \xi) \ll \h^{-\delta |\alpha|} \langle \h^{-\delta} \xi \rangle^{-m}
\end{equation}
for all fixed $\alpha,m$, where $\langle \xi \rangle := (1 + |\xi|^2)^{1/2}$.  Informally, such symbols are concentrated on the ball $\theta + \O(\h^\delta)$ and vary mildly there.  Clearly
\begin{equation*}
  \mathcal{A}_{\delta}(\mathfrak{p}^\wedge) \subseteq S_{\delta}^{-\infty}(\mathfrak{p}^\wedge) \cap \underline{S}^{-\infty}(\mathfrak{p}^\wedge).
\end{equation*}

\begin{lemma}\label{lem:star-product-asymptotics-applied-to-cal-A-delta}
  Fix $r \in \mathbb{Z}_{\geq 0}$.  Let $a \in \mathcal{A}_\delta(\mathfrak{p}^\wedge)$.  Suppose that $b \in S^\infty_0(\mathfrak{p}^\wedge)$ vanishes to order at least $r$ at $\theta$ (i.e., $\partial^\alpha b(\theta) = 0$ whenever $\lvert \alpha \rvert < r$).  Define $\star_{\h}$ with respect to any pair of $\delta$-admissible cutoffs $\chi_1, \chi_2$.  Then
  \begin{equation*}
    a \star_{\h} b, \, \, b \star_{\h} a \in \h^{r \delta} \mathcal{A}_\delta(\mathfrak{p}^\wedge).
  \end{equation*}
\end{lemma}
\begin{proof}
  Since $a \in \underline{S}^{-\infty}(\mathfrak{p}^\wedge)$ and $b \in \underline{S}^\infty(\mathfrak{p}^\wedge)$, we know by Theorem \ref{thm:extended-star-product} that $a \star_{\h} b \in \underline{S}^{-\infty}(\mathfrak{p}^\wedge)$.  It remains to verify the analogue of \eqref{eq:partialalpha-axi-ll}, i.e., that for fixed $\alpha$ and $m$, the rescaled star product $a \star_{\h} b$ (say) satisfies the estimate (with $c = a \star_{\h} b$)
  \begin{equation}\label{eq:partial-c-theta-xi-lemma}
    \partial^\alpha c(\theta + \xi) \ll \h^{\delta (r - |\alpha|)} \langle \h^{-\delta} \xi \rangle^{-m}.
  \end{equation}
  It is clear that if $\ell$ is fixed large enough in terms of $\alpha$ and $m$, then any element of $\h^\ell S_{\delta}^{-\ell}$ satisfies \eqref{eq:partial-c-theta-xi-lemma}.  By the star product asymptotics (Theorem \ref{thm:basic-star-prod}), we thereby reduce to verifying that for each fixed $j \in \mathbb{Z}_{\geq 0}$, the homogeneous star product component $\h^j a \star^j b$ satisfies \eqref{eq:partial-c-theta-xi-lemma}, i.e.,
  \begin{equation}\label{eq:cuhmv904yq}
    \h^j \partial^\alpha (a \star^j b)(\theta + \xi) \ll \h^{\delta (r - |\alpha|)} \langle \h^{-\delta} \xi \rangle^{-m}.
  \end{equation}

  Fix $\ell \in \mathbb{Z}_{\geq 0}$ so that $b \in S^\ell_{0}(\mathfrak{p}^\wedge)$.  We claim that for each fixed $\alpha$ multi-index, we have
  \begin{equation}\label{eq:cuhmv966ci}
    \partial^\alpha b(\theta + \xi ) \ll |\xi|^{\max(0, r - |\alpha|)} \langle \xi  \rangle^{\ell - |\alpha|}.
  \end{equation}
  When $\lvert \alpha \rvert \geq r$ or $\lvert \xi \rvert \geq 1$, the claim follows from the definition of $S_0^{\ell}(\mathfrak{p}^\wedge)$, which gives
  \begin{equation*}
    \partial^\alpha b(\theta + \xi) \ll \langle \theta + \xi \rangle^{\ell - \lvert \alpha \rvert}
    \ll |\xi|^{\max(0, r - |\alpha|)} \langle \xi  \rangle^{\ell - |\alpha|},
  \end{equation*}
  where in the last step, we used that $\langle \theta + \xi \rangle \asymp \langle \xi \rangle$ (since $\theta \ll 1$) and that either $\max(0, r - \lvert \alpha \rvert) = 0$ (in the case $\lvert \alpha \rvert \geq r$) or that $\lvert \xi \rvert^{\max(0, r - \lvert \alpha \rvert)} \geq 1$ (in the case $\lvert \xi \rvert \geq 1$).  We have reduced to the case that $\lvert \alpha \rvert < r$ and $\lvert \xi \rvert \leq 1$.  Our task is then to show that
  \begin{equation}\label{eq:cuhmv9y0kx}
    \partial^\alpha b(\theta + \xi) \ll \lvert \xi \rvert^{r - \lvert \alpha \rvert}.
  \end{equation}
  We apply Taylor's theorem to $\partial^\alpha b$ taken to order $r - \lvert \alpha \rvert$, which gives
  \begin{equation*}
    \partial^\alpha b(\theta + \xi)
    = \sum_{\lvert \beta \rvert < r - \lvert \alpha \rvert}
    \frac{\partial^{\alpha + \beta} b(\theta)}{\beta !} \xi^\beta
    +
    \sum_{\lvert \beta \rvert = r - \lvert \alpha \rvert}
    \frac{\partial^{\alpha + \beta} b(\theta + t_\xi \xi)}{\beta!} \xi^\beta
  \end{equation*}
  for some $t_\xi \in(0, 1)$.  Our hypotheses on $b$ imply that the first sum above vanishes identically.  For the second sum, we use the definition of $S_0^{\ell}(\mathfrak{p}^\wedge)$ to see that $\partial^{\alpha + \beta} b(\theta + t_\xi \xi) \ll 1$.  Since $\xi^\beta \ll \lvert \xi \rvert^{r - \lvert \alpha \rvert}$, we obtain \eqref{eq:cuhmv9y0kx}, hence the claim.

  We now verify \eqref{eq:cuhmv904yq}.  We consider first the case $j = 0$, in which $\h^j a \star^j b = a b$.  By the product rule for differentiation, we deduce from \eqref{eq:cuhmv966ci} that for each fixed $k \geq 0$,
  \begin{equation*}
    \partial^\alpha (a b)(\theta + \xi) \ll
    \max_{
      \substack{
        d_1, d_2 \geq 0:   \\
        d_1 + d_2 = |\alpha|
      }
    }
    |\xi|^{\max(0,r - d_1)}
    \h^{- \delta d_2}
    \langle \xi  \rangle^{\ell - d_1}
    \langle \h^{-\delta} \xi  \rangle^{-k}.
  \end{equation*}
  To deduce \eqref{eq:cuhmv904yq} from this, we need to verify that for all $d_1$ and $d_2$ as above, there is a fixed $k \geq 0$ so that
  \begin{equation}\label{eq:ximax0-r-d_1}
    |\xi|^{\max(0,r - d_1)}
    \langle \xi  \rangle^{\ell - d_1}
    \ll
    \h^{\delta (r - d_1)}
    \langle \h^{- \delta} \xi  \rangle^k.
  \end{equation}
  This follows easily, for any $k \geq r + \ell$, by some case analysis:
  \begin{itemize}
  \item Suppose $|\xi| \leq \h^\delta$.  Then $\langle \xi \rangle \asymp 1$ and $\langle \h^{-\delta} \xi \rangle \asymp 1$.  If $r \geq d_1$, then $|\xi|^{\max(0,r - d_1)} \leq \h^{\delta (r - d_1)}$, while if $r < d_1$, then $|\xi|^{\max(0, r - d_1)} \leq 1 \leq \h^{\delta(r-d_1)}$.
  \item Suppose $|\xi| > \h^\delta$ and $|\xi| \ll 1$.  Then $\langle \xi \rangle \asymp 1$ and $\langle \h^{-\delta} \xi \rangle \asymp \h^{-\delta}$.  The LHS of \eqref{eq:ximax0-r-d_1} is $\ll 1$, while the RHS is $\gg 1$ provided that $k \geq r$.
  \item Suppose $|\xi| \ggg 1$.  Then the LHS is
    \begin{equation*}
      \ll |\xi|^{r + \ell - 2 d_1} \leq |\xi|^{r + \ell},
    \end{equation*}
    while the RHS is
    \begin{equation*}
      \gg \h^{\delta (r - d_1) - \delta k} |\xi|^k
      \geq |\xi|^{r+\ell}
    \end{equation*}
    provided that $k \geq r + \ell$.
  \end{itemize}

  The case $j \geq 1$ may be reduced to the case $j=0$.  To see this, we recall from \S\ref{sec:formal-expansion} that $\h^j a \star^j b(\zeta)$ is given by a linear combination of expressions
  \begin{equation*}
    \h^j \zeta^\gamma \partial^\alpha a(\zeta) \partial^\beta b(\zeta)
  \end{equation*}
  indexed by multi-indices $\alpha, \beta, \gamma$ satisfying the conditions
  \begin{equation*}
    |\alpha| + |\beta| - |\gamma| = j,
    \quad
    |\gamma| \leq \min(|\alpha|,|\beta|),
  \end{equation*}
  which imply that $\max(|\alpha|,|\beta|) \leq j$.  It follows from the definitions that
  \begin{equation*}
    \partial^\alpha a \in \h^{- \delta |\alpha|} \mathcal{A}_{\delta}(\mathfrak{p}^\wedge), \quad
    \partial^{\beta} b \in S_0^{\ell-|\beta|}(\mathfrak{p}^\wedge).
  \end{equation*}
  Moreover, $\partial^{\beta} b$ vanishes to order $\geq \max(0,r - |\beta|)$ at $\theta$.  The $j = 0$ analysis gives
  \begin{equation*}
    \h^j \partial^\alpha a \cdot \partial^\beta b
    \in
    \h^{j - \delta |\alpha| + \delta ( r - |\beta|)} \mathcal{A}_{\delta}(\mathfrak{p}^\wedge).
  \end{equation*}
  Since $\max(|\alpha|,|\beta|) \leq j$ and $\delta \leq 1/2$, we have
  \begin{equation*}
    j - \delta |\alpha| + \delta ( r - |\beta|)
    \geq \delta r.
  \end{equation*}
  Therefore $\h^j \partial^\alpha a \cdot \partial^\beta b \in \h^{r \delta } \mathcal{A}_\delta(\mathfrak{p}^\wedge)$.  It is clear that $\mathcal{A}_{\delta}(\mathfrak{p}^\wedge)$ is preserved under multiplication by the fixed polynomial $\zeta \mapsto \zeta^\gamma$, so we conclude that $\h^j a \star^j b \in \h^{r \delta} \mathcal{A}_{\delta}(\mathfrak{p}^\wedge)$, as required.
\end{proof}

\begin{lemma}\label{lem:star-product-asymptotics-applied-to-cal-A-delta-2}
  For $a,b \in \mathcal{A}_\delta(\mathfrak{p}^\wedge)$, we have
  \begin{equation*}
    a \star_{\h} b \in \mathcal{A}_\delta(\mathfrak{p}^\wedge),
  \end{equation*}
  where the star product is defined with respect to any pair of $\delta$-admissible cutoffs.
\end{lemma}
\begin{proof}
  The proof is similar to, but much simpler than, that of Lemma \ref{lem:star-product-asymptotics-applied-to-cal-A-delta}, hence left to the reader.  (Alternatively, we could appeal to Lemma \ref{lem:let-c-frakC-vs-frakA} and Proposition \ref{lem:standard:subcl-mathfr-theta}.)
\end{proof}

\subsubsection{Lifts of the class of localized Whittaker functions}
We next define (with $\mathfrak{C}(\dotsb)$ as in \S\ref{sec:spec-group-sett}) \index{classes!$\tilde{\mathcal{W}}_\delta(\mathfrak{p}), \tilde{\mathcal{W}}_\delta(\mathfrak{h})$}
\begin{equation*}
  \tilde{\mathcal{W}}_\delta(\mathfrak{p}) := \mathfrak{C}(P, - \theta, T, \delta),
\end{equation*}
\begin{equation*}
  \tilde{\mathcal{W}}_\delta(\mathfrak{h}) := \mathfrak{C}(H, - \theta, T, \delta).
\end{equation*}

\begin{lemma}\label{lem:W-delta-via-Fourier}
  For each $a \in \mathcal{A}_\delta(\mathfrak{p}^\wedge)$, each fixed $\delta_+ \in (\delta,1/2)$ and each $\chi \in \mathcal{X}_{\delta_+}(\mathfrak{p})$, the truncated rescaled inverse Fourier transform $\chi a_{\h}^\vee$, pushed forward via the exponential map, belongs to $\tilde{\mathcal{W}}_{\delta}(\mathfrak{p})$.  Conversely, every element of $\tilde{\mathcal{W}}_{\delta}(\mathfrak{p})$ arises in this way.
\end{lemma}
\begin{proof}
  This follows from part \eqref{itm:standard:fourier-description-of-bumps-arch} of Lemma \ref{lem:let-c-frakC-vs-frakA}.  As noted there, the proof consists of elementary Fourier analysis, i.e., repeated integration by parts in the definition of the Fourier integral.  The informal content here is that a symbol $a \in \mathcal{A}_\delta(\mathfrak{p}^\wedge)$ is a smooth bump on $\theta + \O(\h^\delta)$ if and only if the rescaled inverse Fourier transform $a_{\h}^\vee$ is concentrated on $\O(\h^{1-\delta})$ and oscillates like $e^{-\theta/\h}$.
\end{proof}

\begin{lemma}\label{lem:map-c_cinftyp-ni-tilde-W-delta-W}
  The maps $\phi \mapsto \pi(\phi) \Theta_{\pi,\psi_T}$, defined on the spaces of smooth, compactly-supported measures on $P$ and on $H$, induce surjective class maps
  \begin{equation*}
    \tilde{\mathcal{W}}_\delta(\mathfrak{p}) \twoheadrightarrow \mathcal{W}_\delta,
  \end{equation*}
  \begin{equation*}
    \tilde{\mathcal{W}}_\delta(\mathfrak{h}) \twoheadrightarrow \mathcal{W}_\delta.
  \end{equation*}
\end{lemma}
\begin{proof}
  This is a translation of Lemma \ref{lem:fix-delta-in}.
\end{proof}

\subsubsection{Localized Whittaker functions via localized symbols}

\begin{lemma}\label{lem:every-elem-mathc}
  For fixed $\delta_+ \in (\delta,1/2)$ and $(a,\chi) \in \mathcal{A}_\delta(\mathfrak{p}^\wedge) \times \mathcal{X}_{\delta_+}(\mathfrak{p})$, we have
  \begin{equation*}
    \Opp_{\h}(a:\chi) \Theta_{\pi,\psi} \in \mathcal{W}_\delta.
  \end{equation*}
  Conversely, every element of $\mathcal{W}_\delta$ may be expressed in this way.
\end{lemma}
\begin{proof}
  By Lemmas \ref{lem:W-delta-via-Fourier} and \ref{lem:map-c_cinftyp-ni-tilde-W-delta-W}, since $\Opp_{\h}(a : \chi) = \pi(\chi a_{\h}^\vee)$.
\end{proof}

\subsection{Localization with respect to $\mathfrak{p}$}\label{sec:local-with-resp-P}
In this section, we describe the asymptotic behavior of the class $\mathcal{W}_\delta$ under the mirabolic subgroup $P$, using the explicit description of the action of that subgroup on the Kirillov model.  The results of this section are essentially formal, but provide a convenient base for the more complicated arguments to follow.

More precisely, our first objective is to define on $\mathcal{W}_\delta$ the structure of an $S^\infty_\delta(\mathfrak{p}^\wedge)$-module (Definition \ref{defn:s-infinity-delta-module}).  We then verify that the latter is $(\mathfrak{p},\delta)$-localized at $\theta$ (Definition \ref{defn:localized-vectors}).

As the underlying $\mathfrak{U}(\mathfrak{g}_\mathbb{C})$-module for $\mathcal{W}_\delta$, we take the space $\pi^\infty$ of all smooth vectors.  For the operator assignment, we fix $\delta_{+,0} \in (\delta,1/2)$, choose a cutoff $\chi_0 \in \mathcal{X}_{\delta_{+,0}}(\mathfrak{p})$, and take
\begin{equation*}
  \Opp(a) := \Opp_{\h}(a:\chi_0).
\end{equation*}
We must check that
\begin{equation*}
  \Opp( S^\infty_\delta(\mathfrak{p}^\wedge)) \mathcal{W}_\delta \subseteq \mathcal{W}_\delta
\end{equation*}
For $v \in \mathcal{W}_\delta$ and $a \in S^\infty_\delta(\mathfrak{p}^\wedge)$, we know from Theorem \ref{thm:qualitative-extended-operator-assignment} (applied on the group $P$) that $\Opp(a) v$ lies in $\pi|_P^\infty$, but it is not immediately obvious that it lies in $\mathcal{W}_\delta$.  Let us verify that this is the case:
\begin{lemma}\label{lem:each-a-v}
  For each fixed $\delta_+ \in (\delta, 1/2)$ and all $(\chi,a,v) \in \mathcal{X}_{\delta_+}(\mathfrak{p}) \times S_\delta^\infty(\mathfrak{p}^\wedge) \times \mathcal{W}_\delta$, we have $\Opp_{\h}(a:\chi) v \in \mathcal{W}_\delta$.
\end{lemma}
\begin{proof}
  By Lemma \ref{lem:every-elem-mathc}, we may write
  \begin{equation}\label{eq:v-=-opp_hb:chi}
    v = \Opp_{\h}(b:\chi') \Theta_{\pi,\psi}
  \end{equation}
  for some fixed $\delta_+' \in (\delta_{+},1/2)$ and some $(b,\chi') \in \mathcal{A}_\delta(\mathfrak{p}^\wedge) \times \mathcal{X}_{\delta_+'}(\mathfrak{p})$.  Applying the definition of $\Opp_{\h}(a:\chi)$ followed by the composition formula (Lemma \ref{lem:action-opp-on-cal-W-delta}), we see that
  \begin{equation}\label{eq:oppa-v-Theta-stuff}
    \Opp_{\h}(a:\chi) v
    = \Opp_{\h}(a:\chi) \Opp_{\h}(b:\chi') \Theta_{\pi,\psi}
    = \Opp_{\h}(a \star_{\h} b:\chi'') \Theta_{\pi,\psi},
  \end{equation}
  for some $\chi'' \in \mathcal{X}_{\delta_+'}(\mathfrak{p})$.  By Lemma \ref{lem:star-product-asymptotics-applied-to-cal-A-delta-2}, we have $a \star_{\h} b \in \mathcal{A}_\delta(\mathfrak{p}^\wedge)$.
  By Lemma \ref{lem:every-elem-mathc}, it follows that $\Opp_{\h}(a \star_{\h} b:\chi'') \Theta_{\pi,\psi} \in \mathcal{W}_\delta$, hence $\Opp_{\h}(a:\chi) v \in \mathcal{W}_\delta$, as required.
\end{proof}

\begin{lemma}\label{lem:class-map-defined-cal-W-delta-module}
  The class map defined above equips $\mathcal{W}_\delta$ with the structure of an $S_\delta^\infty(\mathfrak{p}^\wedge)$-module.
\end{lemma}
\begin{proof}
  We check the desiderata of Definition \ref{defn:s-infinity-delta-module} one by one.  It is clear that $\mathcal{W}_\delta$ is an $\O(1)$-submodule of $\pi^\infty$.  The identity $\Opp(p) v = \sym(p_{\h}) v$ for $p \in \Sym(\mathfrak{g})$ is given by \eqref{eq:oppp-=-pisymp_h}.  In view of Lemma \ref{lem:each-a-v}, it remains to verify that for all $a,b \in S^\infty_\delta(\mathfrak{p}^\wedge)$ and $v \in \mathcal{W}_\delta$, we have
  \begin{equation*}
    \Opp(a) \Opp(b) v \equiv \Opp(a \star_{\h} b) v \pmod{\h^\infty \mathcal{W}_\delta}.
  \end{equation*}
  Choose cutoffs $\chi_1, \chi_2 \in \mathcal{X}_{\delta_+}(\mathfrak{p})$, with the support of $\chi_1$ (resp.\ $\chi_2$) taken large in terms of that of $\chi_0$ (resp.\ $\chi_1$).  Then by
  \eqref{eq:qualitative-composition-formula} and Lemma \ref{lem:fix-0-leq}, we have
  \begin{equation*}
    \Opp(a) \Opp(b) = \Opp_{\h}(a \star_{\h} b: \chi_1) = \Opp(a \star_{\h} b) + \Opp_{\h}(\mathcal{E}:\chi_2),
  \end{equation*}
  where $\mathcal{E} \in \h^\infty S^{-\infty}$.  By Lemma \ref{lem:each-a-v}, we have $\Opp_{\h}(\mathcal{E}:\chi_2) v \in \h^\infty \mathcal{W}_\delta$, whence the required conclusion.
\end{proof}

\begin{proposition}\label{prop:local-with-resp-p}
  The $S^\infty_\delta(\mathfrak{p}^\wedge)$-module $\mathcal{W}_\delta$ is $(\mathfrak{p},\delta)$-localized at $\theta$ (Definition \ref{defn:localized-vectors}).
\end{proposition}
\begin{proof}
  Fix $r \in \mathbb{Z}_{\geq 0}$.  Let $a \in S^\infty_0(\mathfrak{p}^\wedge)$ be a polynomial symbol that vanishes to order $r$ at $\theta$.  Let $v \in \mathcal{W}_\delta$.  We must verify that
  \begin{equation}\label{eq:oppa-v-in}
    \Opp(a) v \in \h^{r \delta} \mathcal{W}_\delta.
  \end{equation}

  To that end, choose $b \in \mathcal{A}_\delta(\mathfrak{p}^\wedge)$ as in the proof of Lemma \ref{lem:each-a-v}, so that \eqref{eq:v-=-opp_hb:chi} and \eqref{eq:oppa-v-Theta-stuff} hold.  By Lemma \ref{lem:star-product-asymptotics-applied-to-cal-A-delta}, we see that
  \begin{equation*}
    a \star_{\h} b \in \h^{r \delta} \mathcal{A}_\delta(\mathfrak{p}^\wedge).
  \end{equation*}
  The required membership \eqref{eq:oppa-v-in} then follows from Lemma \ref{lem:every-elem-mathc}.
\end{proof}

\section{Localization with respect to $G$}\label{sec:local-with-resp}

Here we refine the results of \S\ref{sec:local-with-resp-P} to describe the asymptotic behavior of $\mathcal{W}_\delta$ under the action of the full group $G$.  We continue to define $\theta := \theta_P(\psi) \in \mathfrak{p}^\wedge$ as \S\ref{sec:parameter-theta} and \S\ref{sec:class-bump-functions}, and now define $\tau := \tau(\pi, \psi, T) \in \mathfrak{g}^\wedge$ as Lemma \ref{lem:construction-tau-from-theta}, so that $\tau$ is characterized by
\begin{equation*}
  [\tau] = \h \lambda_\pi, \qquad
  \tau|_{\mathfrak{p}} = \theta.
\end{equation*}

\subsection{Statement of main result}
\begin{theorem}\label{thm:cal-W-delta-g-acts-localized}
  $\mathcal{W}_{\delta}$ is
  $(\mathfrak{g},\delta)$-localized at $\tau$.
\end{theorem}
In particular, for fixed $x \in \mathfrak{g}$, we have $\h \pi(x) \mathcal{W}_\delta \subseteq \mathcal{W}_\delta$.  This fact does not seem \emph{a priori} obvious (compare with \S\ref{sec:empha-priori-bounds}).

In \S\ref{sec:control-over-norms}--\ref{sec:proof-theor-refthm:b}, we apply Theorem \ref{thm:cal-W-delta-g-acts-localized} to prove Theorem \ref{thm:big-group-whittaker-behavior-under-G}.

\subsection{Reduction to factoring $\mathfrak{g}$ through $\mathfrak{p}$}
We explain here how the proof of Theorem \ref{thm:cal-W-delta-g-acts-localized} reduces to that of the following result, which describes the action on $\mathcal{W}_\delta$ of $\mathfrak{g}$ in terms of that of $S^\infty_0(\mathfrak{p}^\wedge)$.

In what follows, we define $\Opp(a)$, for a polynomial $a$ on $\mathfrak{g}^\wedge$, to be the $\h$-rescaling of the symmetrization map, as in \eqref{eq:oppp-:=-symp_h}; equivalently, $\Opp(a) := \Opp_{\h}(a:\chi)$ for any cutoff $\chi$ on $\mathfrak{g}$.

\begin{theorem}\label{thm:g-via-p}
  For each fixed $x \in \mathfrak{g}$, there is a symbol
  \[x^{\mathfrak{p}} \in S_0^\infty(\mathfrak{p}^\wedge)\] with the following properties.
  \begin{enumerate}[(i)]
  \item $\Opp(x) v \equiv \Opp(x^\mathfrak{p}) v \mod{\h^\infty \mathcal{W}_\delta}$ for each $v \in \mathcal{W}_\delta$.
  \item $x^\mathfrak{p}(\theta) = x(\tau) + \O(\h)$.
  \end{enumerate}
\end{theorem}

\begin{proof}[Proof of Theorem
  \ref{thm:cal-W-delta-g-acts-localized},
  assuming Theorem \ref{thm:g-via-p}]
  We must check that $\mathcal{W}_\delta$ is $(\mathfrak{g},\delta)$-localized at $\tau$.  By Lemma \ref{lem:localization-criterion-degree-1}, it is enough to show that for each $v \in \mathcal{W}_\delta$ and each fixed $x \in \mathfrak{g}$, we have $\Opp(x) v \equiv x(\tau) v \mod{\h^\delta \mathcal{W}_\delta}$.  Taking $x^\mathfrak{p}$ as in Theorem \ref{thm:g-via-p}, it is enough to show that
  \begin{equation}\label{eqn:required-estimate-reduce-to-going-x-to-x-p}
    \Opp(x^\mathfrak{p}) v
    \equiv x^\mathfrak{p}(\theta) v
    \pmod{\h^\delta \mathcal{W}_\delta}.
  \end{equation}
  We have seen in Proposition \ref{prop:local-with-resp-p} that $\mathcal{W}_\delta$ is an $S^\infty_\delta(\mathfrak{p}^\wedge)$-module that is $(\mathfrak{p},\delta)$-localized at $\tau$.  The symbol $x^\mathfrak{p} - x^\mathfrak{p}(\theta) \in S^\infty_0(\mathfrak{p}^\wedge)$ vanishes to order $\geq 1$ at $\theta$, so the required conclusion \eqref{eqn:required-estimate-reduce-to-going-x-to-x-p} follows from part \eqref{item:2-thm:S-infinity-modules-localized-stuff} of Theorem \ref{thm:S-infinity-modules-localized-stuff}.
\end{proof}

\subsection{Preliminary reductions}\label{sec:cuhn9jrbjg}
Here we reduce the proof of Theorem \ref{thm:g-via-p} to the crucial cases.

By linearity, we may assume that $x$ lies in some fixed spanning set for $\mathfrak{g}$.  Writing $e_{i j} \in \mathfrak{g}$ for the standard basis elements, we take our spanning set to be the union of the subsets
\begin{equation*}
  \mathfrak{p},
  \qquad \mathfrak{z},
  \qquad \{ e_{n j} : 1 \leq j \leq n-1 \}.
\end{equation*}
(To see that this is a spanning set, observe that $e_{i j} \in \mathfrak{p}$ for $i \neq n$, while $e_{nn} \in \mathfrak{p} + \mathfrak{z}$.)

If $x \in \mathfrak{p}$, then we may take $x^\mathfrak{p} = x \in S^1_0(\mathfrak{p}^\wedge)$.

Suppose that $x \in \mathfrak{z}$.  Then $x$ is (up to standard identifications) its own image under the Harish--Chandra isomorphism, and defines a regular function on $[\mathfrak{g}^\wedge]$.  The operator $\Opp(x)$ acts by the scalar $x(\h \lambda_\pi) \ll 1$, so we can take for $x^\mathfrak{p}$ the constant element of $S^0_0(\mathfrak{p}^\wedge)$ given by $x^\mathfrak{p}(\xi) := x(\h \lambda_\pi)$.  Note that the condition $[\tau] = \h \lambda_\pi$ implies that the restriction of $\tau$ to $\mathfrak{z}$ identifies with the linear form on $\mathfrak{z} \hookrightarrow \mathfrak{Z}(\mathfrak{g})$ given by evaluation at $\h \lambda_\pi$, so that $x(\h \lambda_\pi) = x(\tau)$.

The standard basis element $e_{n n}$ of $\mathfrak{g}$ lies in $\mathfrak{p} + \mathfrak{z}$: it is the difference $z - x$, where $z = e_{1 1} + \dotsb + e_{n n} \in \mathfrak{z}$ and $x = e_{1 1} + \dotsb + e_{n-1,n-1} \in \mathfrak{h} \subseteq \mathfrak{p}$.  We may thus take
\[
  e_{n n}^{\mathfrak{p}} := (e_{1 1} + \dotsb + e_{n n})(\h \lambda_\pi) - (e_{1 1} + \dotsb + e_{n-1,n-1}) \in S^1_0(\mathfrak{p}^\wedge).
\]
Then
\begin{equation}\label{eq:cui1fz583k}
  \Opp(e_{n n}^{\mathfrak{p}}) = \Opp(e_{n n}),
\end{equation}
since $(e_{1 1} + \dotsb + e_{n n})$ acts by the scalar $(e_{1 1} + \dotsb + e_{n n})(\h \lambda_\pi)$.

To complete the proof of Theorem \ref{thm:g-via-p}, it remains to consider $x$ of the form $e_{n j}$ for some $j \in \{1,\dotsc, n-1\}$.

\subsection{A priori bounds}\label{sec:empha-priori-bounds}

\begin{proposition}\label{prop:standard2:each-fixed-x}
  For each fixed $x \in \mathfrak{U}(\mathfrak{g})$, we have $\pi(x) \h^\infty \mathcal{W}_\delta \subseteq \h^\infty \mathcal{W}_\delta$.
\end{proposition}
\begin{remark}
  Since the notation is a bit concise, we pause to explicate the conclusion.  \emph{Let $W \in \mathcal{W}(\pi,\psi_T)$.  Suppose that for each fixed $m \geq 0$, we have $W \in \h^m \mathcal{W}_\delta$.  Then for each fixed $x \in \mathfrak{U}(\mathfrak{g})$ and fixed $m \geq 0$, we have $\pi(x) W \in \h^m \mathcal{W}_\delta$.}  Informally, the proposition says that the ``negligible'' subclass $\h^\infty \mathcal{W}_\delta$ of $\mathcal{W}_\delta$ ``remains negligible'' after taking $G$-derivatives.  In contrast to the corresponding assertion for $P$-derivatives, this is not obvious, and would be false without our running assumption that $\h \lambda_\pi = \O(1)$.  Of course, once we have proved Theorem \ref{thm:g-via-p}, we will know the much sharper assertion that $\pi(x) \mathcal{W}_\delta \subseteq \h^{-1} \mathcal{W}_\delta$.  That sharpening will be derived via a ``bootstrapping'' argument towards which Proposition \ref{prop:standard2:each-fixed-x} provides the first step.

  We note that the statement of Proposition \ref{prop:standard2:each-fixed-x} is weaker than the following statement: \emph{for each $W \in \mathcal{W}_\delta$ and fixed $x \in \mathfrak{U}(\mathfrak{g})$, we have $\pi(x) W \in \h^{-m} \mathcal{W}_\delta$ for some fixed $m$}.  We do not prove that stronger statement here, although it will eventually fall out as a byproduct of Theorem \ref{thm:g-via-p}.
\end{remark}
As we will explain, the proof of Proposition \ref{prop:standard2:each-fixed-x} is a direct quantification of Jacquet's proof of Theorem \ref{thm:jacquet-on-kirillov-model}.  We begin by recalling the crucial statement proved by Jacquet.
\begin{lemma}\label{lemma:cuhn9kcb25}
  Let $\psi$ be a nondegenerate unitary character of $N$.  There is a linear map
  \begin{equation*}
    \ell_{\psi} : \mathfrak{U}(\mathfrak{g}) \rightarrow \mathfrak{U}(\mathfrak{h}) \otimes \mathfrak{Z}(\mathfrak{g})
  \end{equation*}
  with the following property.  Let $\pi$ be a generic irreducible unitary representation of $G$, with central character $\lambda_\pi$.  Let
  \begin{equation*}
    \ell_\pi : \mathfrak{U}(\mathfrak{h}) \otimes \mathfrak{Z}(\mathfrak{g}) \rightarrow \mathfrak{U}(\mathfrak{g})
  \end{equation*}
  denote the linear map given on pure tensors by
  \begin{equation*}
    y \otimes z \mapsto \lambda_\pi(z) y.
  \end{equation*}
  Then for each $x \in \mathfrak{U}(\mathfrak{g})$,
  \begin{equation}\label{eq:pix-theta_pi-psi}
    \pi(x) \Theta_{\pi,\psi} = \pi(\ell_{\pi}(\ell_\psi(x))) \Theta_{\pi,\psi}.
  \end{equation}
\end{lemma}
\begin{proof}
  Besides notation, this is \cite[Lemma 3]{MR2733072}.  Indeed, the third inclusion of \emph{loc.\ cit.}\ asserts the following.
  \begin{itemize}
  \item Let $V$ be a $\mathfrak{U}(\mathfrak{g})$-module, and let $\nu \in V$ be an element on which $\mathfrak{n} := \Lie(N)$ acts via the differential of $\psi$.  Then
    \begin{equation}\label{eq:cuhn9lgz7j}
      \mathfrak{U}(\mathfrak{g}) \nu \subseteq \mathfrak{U}(\mathfrak{h}) \mathfrak{Z}(\mathfrak{g}) \nu.
    \end{equation}
  \end{itemize}
  To pass from \eqref{eq:cuhn9lgz7j} to the required conclusion, choose a basis $\{x_i\}$ of $\mathfrak{U}(\mathfrak{g})$.  Take for $V$ the universal Whittaker module $\mathfrak{U}(\mathfrak{g}) \otimes_{\mathfrak{U}(\mathfrak{n})} \mathbb{C}_\psi$, generated by a vector $\nu$ on which $\mathfrak{n}$ acts via the differential of $\psi$.  By \eqref{eq:cuhn9lgz7j}, we may write
  \begin{equation*}
    x_i \nu = \sum_j y_{i j} z_{i j} \nu,
  \end{equation*}
  where $y_{i j} \in \mathfrak{U}(\mathfrak{h})$, $z_{i j} \in \mathfrak{Z}(\mathfrak{g})$ and $j$ runs over some finite index set.  Define $\ell_\psi$ to be the linear map sending $x_i$ to $\sum_j y_{i j} \otimes z_{i j}$.  Recall (from \S\ref{sec:kirillov-model}) that $N$ acts on $\Theta_{\pi,\psi}$ via the character $\psi$.  In particular, $\mathfrak{n}$ acts on $\Theta_{\pi,\psi}$ via the differential of $\psi$, so there is a unique $\mathfrak{U}(\mathfrak{g})$-module map $V \rightarrow \pi^{-\infty}$ sending $\nu$ to $\Theta_{\pi,\psi}$.  Applying this map to the preceding formula and using that $\pi(z_{i j})$ is multiplication by the scalar $\lambda_\pi(z_{i j})$, we obtain
  \begin{equation*}
    \pi(x_i) \Theta_{\pi,\psi} = \sum_j \lambda_\pi(z_{i j}) \pi(y_{i j}) \Theta_{\pi,\psi}.
  \end{equation*}
  The required formula \eqref{eq:pix-theta_pi-psi} follows by linearity.
\end{proof}

We record a characterization of $\h^\infty \tilde{\mathcal{W}}_\delta(\mathfrak{h}) = T^{-\infty} \mathfrak{C}(H, -\theta_H(\psi), T, \delta)$, following from Lemma \ref{lem:standard2:let-v-p}:
\begin{lemma}\label{lem:scratch-research:hinfty-tild-cons}
  The class $\h^\infty \tilde{\mathcal{W}}_{\delta}(\mathfrak{h})$ consists of all $f \in C_c^\infty(H)$ that are supported on $1 + \O(\h^{1 - \delta_+})$ for some fixed $\delta_+ \in (\delta,1/2)$ and that satisfy, for all fixed $x,y \in \mathfrak{U}(\mathfrak{g})$, the derivative bound
  \begin{equation}\label{eq:x-ast-f}
    \|x \ast f \ast y\|_{\infty} \ll \h^\infty
  \end{equation}
\end{lemma}

We now apply arguments analogous to those of \cite[Proof of Lemma 5]{MR2733072}.  Recall that $T = 1/\h$ and that $x \in \mathfrak{U}(\mathfrak{g})$ is fixed.  The character $\psi$ of $N$ is fixed, so we may fix a functional $\ell_\psi$ for which the conclusion of Lemma \ref{lemma:cuhn9kcb25} holds.  Set $a_T := \diag(T^{n-1}, T^{n-2}, \dotsc, T, 1) \in G$.  Then $\pi(a_T) \Theta_{\pi,\psi}$ is a multiple of $\Theta_{\pi,\psi_T}$, so by applying \eqref{eq:pix-theta_pi-psi} to $\Ad(a_T)^{-1} x$, we see that
\begin{equation*}
  \pi(x) \Theta_{\pi,\psi_T} = \pi (\ell_\pi(\Ad(a_T) \ell_\psi(\Ad(a_T)^{-1} x))) \Theta_{\pi,\psi_T}.
\end{equation*}
Next,
\begin{align*}
  \pi(x) \pi(f) \Theta_{\pi,\psi_T} &= \int_{g \in H} f(g) \pi(x) \pi(g) \Theta_{\pi,\psi_T} \, d g \\
                      &= \int_{g \in H} f(g) \pi(g) \pi(\Ad(g)^{-1}  x) \Theta_{\pi,\psi_T} \, d g.
\end{align*}
Fix a monomial basis $\{x_i\}$ for $\mathfrak{U}(\mathfrak{g})$ and, for $g \in H$ as above, write
\begin{equation*}
  \Ad(g)^{-1} x = \sum_i c_i(g) x_i,
\end{equation*}
where $c_i$ lies in $C^\infty(H)$ and vanishes identically for $i$ outside some fixed finite set of indices.  Then, abbreviating
\begin{equation}\label{eq:cui0yocikh}
  y_{i,\pi,T} :=  \ell_\pi(\Ad(a_T) \ell_\psi(\Ad(a_T)^{-1} x_i)) \in \mathfrak{U}(\mathfrak{g}),
\end{equation}
we have
\begin{align*}
  \pi(\Ad(g)^{-1}  x) \Theta_{\pi,\psi_T}
  &=
    \sum_i c_i(g) \pi(x_i) \Theta_{\pi,\psi_T} \\
  &=
    \sum_i
    c_i(g) \pi (y_{i,\pi,T}) \Theta_{\pi,\psi_T}.
\end{align*}
Setting
\begin{equation*}
  f_i(g) := c_i(g) f(g), \quad
  f_i' := f_i \ast y_{i,\pi,T},
  \quad
  f' := \sum_i f_i',
\end{equation*}
we obtain
\begin{equation}\label{eq:pix-pif-theta_pi}
  \pi(x) \pi(f) \Theta_{\pi,\psi_T} = \pi(f') \Theta_{\pi,\psi_T}.
\end{equation}

Jacquet observes, in the ``$T=1$'' setting, that the map $f \mapsto f'$ on $C_c^\infty(H)$ is continuous.  In the present setting, the continuity is polynomially effective, in the following sense, with respect to $T = 1/\h$:
\begin{lemma}\label{lem:scratch-research:fix-x-in}
  Fix $x \in \mathfrak{U}(\mathfrak{g})$.  Let $f \in \h^\infty \tilde{\mathcal{W}}_\delta(\mathfrak{h})$.   Define $f' \in C_c^\infty(H)$ as above.  Then
  \begin{equation*}
    f' \in \h^\infty \tilde{\mathcal{W}}_\delta(\mathfrak{h}).
  \end{equation*}
\end{lemma}
\begin{proof}
  We appeal to the characterization of $\h^\infty \tilde{\mathcal{W}}_\delta(\mathfrak{h})$ given by Lemma \ref{lem:scratch-research:hinfty-tild-cons}.  Since $\supp(f') \subseteq \supp(f)$, the required support condition for $f'$ is satisfied.  Since $x$ is fixed, the functions $c_i \in C^\infty(H)$ vanish for $i$ outside some fixed finite set, and each function $c_i$ is fixed.  Therefore each $f_i$ likewise lies in $\h^\infty \tilde{\mathcal{W}}_\delta(\mathfrak{h})$, since multiplying by a fixed smooth function preserves the bounds \eqref{eq:x-ast-f} in view of the product rule for derivatives.

  Next, fix $i$.  Recall the definition \eqref{eq:cui0yocikh} of $y_{i, \pi, T}$.  Since
  \begin{itemize}
  \item $x_i$ has fixed degree,
  \item the maps $\Ad(a_T)^{\pm}$ and $\lambda_\pi$ do not increase degree, and
  \item $\ell_\psi$ is fixed and linear, hence sends fixed-degree elements to elements of degree $\O(1)$,
  \end{itemize}
  we see that $y_{i,\pi,T}$ has degree $\O(1)$.  Since $\h = 1/T$ and $\h \lambda_\pi = \O(1)$, we see that the operator norms of $\Ad(a_T)^{\pm}$ on fixed degree subspaces of $\mathfrak{U}(\mathfrak{g})$, and of $\lambda_\pi$ on fixed-degree subspaces of $\mathfrak{Z}(\mathfrak{g})$, are $\ll T^{\O(1)}$; consequently, the coefficients of $y_{i,\pi,T}$ are $\ll T^{\O(1)}$.  We may thus absorb the $y_{i,\pi,T}$ into the quantities $x$ and $y$ appearing in the estimate \eqref{eq:x-ast-f} without affecting that estimate.  We conclude that each $f_i'$ and hence $f'$ itself lies in $\h^\infty \tilde{\mathcal{W}}_\delta$.
\end{proof}

\begin{proof}[Proof of Proposition \ref{prop:standard2:each-fixed-x}]
  Let $v \in \h^\infty \mathcal{W}_\delta$.  By Lemma \ref{lem:standard2:overspill-frak-C-M}, we may find $\ell \geq 0$ with $\ell \ggg 1$ so that $\h^{-\ell} v \in \h^\infty \mathcal{W}_\delta \subseteq \mathcal{W}_\delta$.  We may thus write
  \begin{equation*}
    \h^{-\ell} v = \pi(f_0) \Theta_{\pi,\psi_T}
  \end{equation*}
  for some $f_0 \in \tilde{\mathcal{W}}_\delta$ (by Lemma \ref{lem:map-c_cinftyp-ni-tilde-W-delta-W}).  Then $f := \h^{\ell} f_0 \in \h^\infty \tilde{\mathcal{W}}_{\delta}$.  Fix $x \in \mathfrak{U}(\mathfrak{g})$.  Define $f'$ as above, so that \eqref{eq:pix-pif-theta_pi} holds, whence
  \begin{equation*}
    \pi(x) v = \pi(f') \Theta_{\pi,\psi_T}.
  \end{equation*}
  By Lemma \ref{lem:scratch-research:fix-x-in}, we have $f' \in \h^\infty \tilde{\mathcal{W}}_\delta$.  By Lemma \ref{lem:standard2:overspill-frak-C-M}, we may find $\ell ' \geq 0$ with $\ell ' \ggg 1$ so that $f ' \in \h^{\ell} \tilde{\mathcal{W}}_\delta$.  By Lemma \ref{lem:map-c_cinftyp-ni-tilde-W-delta-W}, it follows that
  \begin{equation*}
    \pi(f') \Theta_{\pi,\psi_T} \in \h^{\ell} \mathcal{W}_\delta \subseteq \h^\infty \mathcal{W}_\delta.
  \end{equation*}
  We conclude, as required, that $\pi(x) v \in \h^\infty \mathcal{W}_\delta$.
\end{proof}

\subsection{Capelli identities}\label{sec:cuhn9jqo9y}
Recall from \S\ref{sec:gener-line-groups} that
\begin{equation*}
  e = (e_{ij }) : \glLie_n(\mathbb{R}) \rightarrow \mathfrak{g} \hookrightarrow \Sym(\mathfrak{g})
\end{equation*}
denotes the ``identity map'' tailored to the symmetric algebra.  We denote by $E = (E_{i j})$ the corresponding ``identity map''
\[
  E = (E_{i j}) : \glLie_n(\mathbb{R}) \rightarrow \mathfrak{g} \hookrightarrow \mathfrak{U}(\mathfrak{g})
\]
but tailored now to the universal enveloping algebra.  Thus $E_{i j} \in \mathfrak{g} \hookrightarrow \mathfrak{U}(\mathfrak{g})$ is the standard basis element.  We note that $E_{i j}$ is the image of $e_{i j}$ under the symmetrization map $\sym : \Sym(\mathfrak{g}) \rightarrow \mathfrak{U}(\mathfrak{g})$.

We define the ``characteristic polynomial'' $\det(X+E) \in \mathfrak{U}(\mathfrak{g})[X]$, as follows.  We identify $X$ with the scalar matrix $\diag(X,\dotsc,X)$ and set
\[
  \rho := \diag \left( \tfrac{n-1}{2}, \tfrac{n-3}{2}, \dotsc , \tfrac{1-n}{2} \right),
\]
so that the sum
\[
  X + \rho + E
\]
defines an $n \times n$ matrix with entries in the ring $\mathfrak{U}(\mathfrak{g})[X]$.

For an $n \times n$ matrix $A = (A_{i j})$ with entries in a not-necessarily-commutative ring, we define the determinant via the following convention:
\begin{equation}\label{eq:cui1foghl9}
  \det(A) := \sum_{\sigma \in S(n)} (-1)^\sigma A_{n,\sigma(n)} A_{n-1,\sigma(n-1)} \dotsb A_{2,\sigma(2)} A_{1,\sigma(1)}.
\end{equation}
\begin{lemma}[Capelli, Howe--Umeda]\label{lemma:cui1foa3up}
  We have $\det(X + \rho + E) \in \mathfrak{Z}(\mathfrak{g})[X]$.  Moreover, writing
  \begin{equation*}
    \gamma : \mathfrak{Z}(\mathfrak{g}_\mathbb{C})[X] \xrightarrow{\cong } \Sym(\mathfrak{g}_\mathbb{C})^G[X]
  \end{equation*}
  for the coefficient-wise application of the Harish--Chandra isomorphism (see \S\ref{sec:cui1fnha7b}), we have
  \begin{equation}\label{eq:cui1fohnhu}
    \gamma(\det(X + \rho + E)) = \det(X + e),
  \end{equation}
  where $\det(X + e)$ is as defined in \S\ref{sec:char-polyn}.
\end{lemma}
\begin{proof}
  Up to cosmetic differences, this is contained in \cite{MR1116239}.  See for instance \cite[\S4.1]{MR2035111} for a formulation much like ours.  We note that the main point is that the coefficients of $\det(X + \rho + E)$ are $G$-invariant.  Taking this for granted, it is easy to compute the image under $\gamma$, as follows:

  Write $\mathfrak{t}_\mathbb{C}$ for the diagonal subalgebra, and $\mathfrak{n}_{\mathbb{C}}^{-}$ for the lower-triangular subspace.  Set $A := X + \rho + E$.  We observe that the contribution to $\det(A)$ of each non-identity permutation $\sigma$ lies in $\mathfrak{n}_\mathbb{C}^- \mathfrak{U}(\mathfrak{g}_\mathbb{C})$.  Indeed, if $j \in \{n, n - 1, \dotsc, 1\}$ is the largest index for which $\sigma(j) \neq j$, then $\sigma(j) < j$, hence
  \begin{equation*}
    A_{n, \sigma(n)} \dotsb A_{j, \sigma(j)}
    \in \Sym (\mathfrak{t}_{\mathbb{C}}) \mathfrak{n}_{\mathbb{C}}^-  = \mathfrak{n}_{\mathbb{C}}^- \Sym (\mathfrak{t}_{\mathbb{C}}) \subseteq \mathfrak{n}_{\mathbb{C}}^- \mathfrak{U}(\mathfrak{g}).
  \end{equation*}
  On the other hand, the identity permutation maps to the product over $j \in \{1, \dotsc, n\}$ of
  \[
    \gamma (X + E_{j j} + \tfrac{n +1 - 2 j}{2}) = X + e_{j j},
  \]
  where we have used that $\rho(E_{j j}) = \frac{n + 1 - 2 j}{2}$.  The restriction of that product to $\mathfrak{t}_\mathbb{C}^*$ coincides with the restriction of $\det(X+e)$.  Using that $\Sym(\mathfrak{t}_\mathbb{C})^W \cong \Sym(\mathfrak{g}_{\mathbb{C}})^G$, we deduce the claimed identity.
\end{proof}

\subsection{Commutative cofactor expansion}\label{sec:comm-cofact-expans}
Recall Lemma \ref{lem:construction-tau-from-theta} concerning the existence and uniquness of $\tau = \tau(\psi,\lambda)$.  Both as an excuse to introduce notation and for motivational purposes, we reformulate here the uniqueness part of the proof of that lemma a bit more verbosely.  The problem is to show that the entries $\tau_{j n}$ are determined by the conditions $[\tau] = \lambda$ and $\tau_P = \theta$.  Starting with the definition \eqref{eq:detx-+-e} of $\det(X + e)$, we group the permutations $\sigma \in S(n)$ according to $j = \sigma(n)$ to obtain the cofactor expansion
\begin{equation}\label{eq:detx-+-e-1}
  \det(X + e) = \sum_{j = 1}^n ( 1_{n=j} X + e_{n j} ) \mathfrak{d}_j(X),
\end{equation}
where $\mathfrak{d}_j(X)$ denotes the $(n,j)$ minor, obtained by restricting the sum and product in the definition \eqref{eq:detx-+-e} of $\det(X+e)$ to $\sigma(n) = j$ and $i \leq n-1$.  The coefficients of $\mathfrak{d}_j(X)$, viewing the latter as a one-variable polynomial, lie in $\Sym(\mathfrak{p})$, hence $\mathfrak{d}_j(X)$ itself lies in $\Sym(\mathfrak{p})[X]$.


We evaluate the expansion \eqref{eq:detx-+-e-1} on $\tau$.  Since $[\tau] = \lambda$, we have (see \eqref{eq:cui1fo922g})
\begin{equation*}
  \det(X + e)(\tau) = \mathcal{P}_\lambda(X)
\end{equation*}
Since $\tau_P = \theta$, we have $\mathfrak{d}_j(X)(\tau) = \mathfrak{d}_j(X)(\theta)$.  Since we have defined the duality between $\mathfrak{g}$ and $\mathfrak{g}^\wedge$ using the trace pairing, we have $e_{n j}(\tau) = \tau_{j n}$.  Recall also, from Lemma \ref{lem:construction-tau-from-theta}, that $\tau_{n n}$ coincides with $\mathfrak{c}_1(\lambda)$, the coefficient of $X^{n-1}$ in $\mathcal{P}_\lambda(X)$.  We may thus write our cofactor expansion in the form
\[
  \sum_{j=1}^{n-1} \tau_{j n} \mathfrak{d}_j(X)(\theta) = \mathcal{P}_\lambda(X) - (X + \mathfrak{c}_1(\lambda)) X^{n-1}.
\]

To see that this expansion determines $\tau$, we specialize the variable $X$ to the elements of any finite collection $\{x_1,\dotsc,x_{n-1}\}$ of distinct real numbers.  This gives a system of $n-1$ linear equations in the same number of unknowns.  With identifications as in \S\ref{sec:parameter-theta}, we have, e.g., for $n=4$,
\begin{equation*}
  X + \theta
  \leftrightarrow
  \begin{pmatrix}
    X & 0 & 0 & ? \\
    \eta_1  & X  & 0 & ? \\
    0 & \eta_2  & X & ? \\
    0 & 0 & \eta_3  & ?
  \end{pmatrix},
\end{equation*}
from which it follows that
\begin{equation}\label{eqn:frak-d-j-X-tau}
  \mathfrak{d}_j(X)(\theta)
  =
  \pm
  X^{j-1}
  \eta_j \eta_{j+1} \dotsb \eta_{n-1}.
\end{equation}
By the nonvanishing of the $\eta_j$ and the invertibility of the Vandermonde determinant, the system characterizes $\tau$, as required.

The proof of Theorem \ref{thm:g-via-p}, given below, amounts to a ``quantum analogue'' of the above calculations.

\subsection{Non-commutative cofactor expansion}
We may similarly write
\begin{equation}\label{eqn:non-comm-cofactor}
  \det(X + \rho + E)
  =
  \sum_{j=1}^n
  \left(1_{n=j}
    \left(X
      - \tfrac{n-1}{2} \right)
    +
    E_{n j} \right)
  \mathfrak{D}_j(X)
\end{equation}
for some ``non-commutative minors'' $\mathfrak{D}_j(X) \in \mathfrak{U}(\mathfrak{p})[X]$, defined in the evident way.

In what follows, we write $\mathfrak{U}^{\leq k}(\dotsb)$ for the $k$th graded part of a universal enveloping algebra, spanned by monomials of degree at most $k$.

By evaluating the commutative and non-commutator minors at a scalar $x$, we obtain elements
\begin{equation*}
  \mathfrak{d}_i(x) \in \Sym^{n-1}(\mathfrak{p}), \quad \mathfrak{D}_i(x) \in \mathfrak{U}^{\leq n-1}(\mathfrak{p}).
\end{equation*}
In fact, if $x \ll 1$, then the polynomial $\mathfrak{d}_i(x)$ defines an element of the symbol class $S^{n-1}_0(\mathfrak{p}^\wedge)$.  We may view $\mathfrak{d}_i(x)$ as the ``classical limit'' of $\mathfrak{D}_j(x)$ in many ways, the following of which suffices for our purposes.
\begin{lemma}\label{lem:non-comm-cofact}
  For each $i \in \{1, \dotsc, n\}$ and scalar $x \ll 1$ there is a symbol
  \[
    d_i(x) \in S_0^{n-1}(\mathfrak{p}^\wedge)
  \]
  with the following properties.
  \begin{enumerate}[(i)]
  \item $\Opp(d_i(x)) = \h^{n-1} \mathfrak{D}_i(\h^{-1} x)$.
  \item $d_i(x) - \mathfrak{d}_i(x) \in \h S_0^{n-2}(\mathfrak{p}^\wedge)$.
  \end{enumerate}
\end{lemma}
\begin{proof}
  Write $\sym : \Sym(\mathfrak{g})[X] \rightarrow \mathfrak{U}(\mathfrak{g})[X]$ for the map of polynomial rings induced by the symmetrization map on coefficients.
  Temporarily set $\mathfrak{S}^k := \Sym^k(\mathfrak{p})$ and $\mathfrak{U}^{\leq k} := \mathfrak{U}^{\leq k}(\mathfrak{p})$.  Observe that
  \[
    \mathfrak{d}_i(X) \in \oplus_{k=0}^{n-1} \mathfrak{S}^k X^{n-1-k},
  \]
  \[
    \mathfrak{D}_i(X) \in \oplus_{k=0}^{n-1} \mathfrak{U}^{\leq k} X^{n-1-k}.
  \]
  By comparing the definitions of $\mathfrak{d}_i$ and $\mathfrak{D}_i$, we see that the coefficient of $X^{n-1-k}$ in $\mathfrak{D}_i(X)$ has image in $\mathfrak{U}^{\leq k} / \mathfrak{U}^{\leq k-1} \cong \mathfrak{S}^k$ given by the corresponding coefficient of $\mathfrak{d}_i(X)$.  Thus
  \[
    \sym(\mathfrak{d}_i(X)) - \mathfrak{D}_i(X) \in \oplus_{k = 1}^{n-1} \mathfrak{U}^{\leq k-1} X^{n-1-k}.
  \]

  We iteratively construct polynomials $\mathfrak{d}_i^{(\ell)}(X)$ for $\ell \in \{0,\dotsc,n-1\}$ by setting $\mathfrak{d}_i^{(0)}(X) := \mathfrak{d}_i(X)$ and requiring that
  \[
    \mathfrak{d}_i^{(\ell)}(X) \in \oplus_{k=\ell}^{n-1} \mathfrak{S}^{k-\ell} X^{n-1-k}
  \]
  and
  \[
    \mathfrak{D}_i(X) = \sum_{\ell=0}^{n-1} \sym(\mathfrak{d}_i^{(\ell)}(X)).
  \]
  For $p \in \Sym(\mathfrak{p})$, write $\sym_{\h}(p) := \sym(p_{\h})$, where $p \mapsto p_{\h}$ multiplies the summand $\mathfrak{S}^k$ by $\h^k$.  Since the rescaled polynomial $\h^{n-1}\mathfrak{D}_i(\h^{-1} X)$ is obtained from $\mathfrak{D}_i(X)$ by multiplying the coefficient of $X^{n-1-k}$ by $\h^k$, we have
  \[
    \h^{n-1} \mathfrak{D}_i(\h^{-1} X) = \sym_{\h}(d_i(X)),
  \]
  where
  \[
    d_i(X) := \sum_{\ell =0 }^{n-1} \h^{\ell} \mathfrak{d}_i^{(\ell)}(X).
  \]
  By construction,
  \[
    d_i(X) \in \oplus_{k=0}^{n-1} S_0^k(\mathfrak{p}^\wedge) X^{n-1-k},
  \]
  \[
    d_i(X) - \mathfrak{d}_i(X) \in \oplus_{k=1}^{n-1} \h S_0^{k-1}(\mathfrak{p}^\wedge) X^{n-1-k}.
  \]
  The required conclusions follow by specializing to a scalar argument $x \ll 1$.
\end{proof}

\subsection{Completion of the proof}\label{sec:completion-proof-2}
We now prove Theorem \ref{thm:g-via-p}.  As noted following its statement, the crucial case is that of the Lie algebra elements $e_{n j}$ for $j \in \{1,\dotsc,n-1\}$.

Our point of departure is the following identity of polynomials in $X$, with coefficients in $\End(\pi)$:
\[
  \pi(\det(X + \rho + E)) = \det(X+e)(\lambda_\pi) = \mathcal{P}_{\lambda_\pi}(X).
\]
By rescaling, we obtain
\begin{equation}\label{eq:cui1fzb1w2}
  \h^{n} \pi(\det(\h^{-1} X + \rho + E)) = \mathcal{P}_{\h \lambda_\pi}(X).
\end{equation}
We choose real numbers $x_1,\dotsc,x_{n-1}$ with $x_i \ll 1$ and $x_i - x_j \gg 1$ for $i \neq j$.  For instance, the choice $x_j := j$ would do.  For each $j$, we use the cofactor expansion \eqref{eqn:non-comm-cofactor} to see that the left-hand side of \eqref{eq:cui1fzb1w2}, specialized to $X = x_j$, is given by
\begin{equation*}
  \sum_{i=1}^{n-1}
  \h \pi(E_{n i}) \h^{n-1}\pi(\mathfrak{D}_i(\h^{-1} x_j))
  + \h^n \pi \left(
    (\h^{-1}x_j - \tfrac{n-1}{2} + E_{n n})\mathfrak{D}_n(\h^{-1} x_j)
  \right).
\end{equation*}
By the formula \eqref{eq:oppp-=-pisymp_h} for quantizations of polynomial symbols, we have $\Opp(e_{n i}) = \h \pi(E_{n i})$.  With $d_i$ as in Lemma \ref{lem:non-comm-cofact}, we set
\begin{equation}\label{eq:aijdixj}
  a_{i j} := d_i(x_j) \in S^{n-1}_0(\mathfrak{p}^\wedge),
\end{equation}
so that $\Opp(a_{i j}) = \h^{n-1}\pi( \mathfrak{D}_i(\h^{-1} x_j))$.  We recall also, from \eqref{eq:cui1fz583k}, that $\Opp(e_{n n}^{\mathfrak{p}}) = \Opp(e_{n n})$.  Using these identities and \eqref{eq:cui1fzb1w2}, we obtain
\begin{equation*}
  \mathcal{P}_{\h \lambda_\pi}(x_j)
  =
  \sum_{i=1}^{n-1}
  \Opp(e_{n i}) \Opp(a_{i j})
  + \Opp (x_j - \h \tfrac{n-1}{2} + e_{n n}^{\mathfrak{p}}) \Opp(a_{n j}).
\end{equation*}
Since all involved symbols are polynomials, the composition formula for the star product holds without any remainder term (see part \eqref{enumerate:cui1f0676u} of Theorem \ref{thm:qualitative-extended-operator-assignment}), so the final term is equal to
\begin{equation*}
  \Opp \left( (x_j - \h \tfrac{n-1}{2} + e_{n n}^{\mathfrak{p}}) \star_{\h}(a_{n j}) \right).
\end{equation*}
Rearranging, we obtain
\begin{equation}\label{eqn:specialized-non-commutative-cofactor}
  \sum_{i=1}^{n-1}
  \Opp(e_{n i})
  \Opp(a_{i j})
  =
  \Opp(f_j),
\end{equation}
where
\[
  f_j := (\h \tfrac{n-1}{2} - x_j - e_{n n}^{\mathfrak{p}}) \star_{\h} a_{n j} + \mathcal{P}_{\h \lambda_\pi}(x_j) \in S^n_0(\mathfrak{p}^\wedge).
\]

We aim now to study $\Opp(e_{n i})$ by approximately inverting the relations \eqref{eqn:specialized-non-commutative-cofactor}.

By Lemma \ref{lem:non-comm-cofact}, we have $a_{i j}(\theta) = \mathfrak{d}_i(x_j)(\theta) + \O(\h)$.  Since the first derivatives of $a_{i j}$ are $\O(1)$ on $\theta + \O(1)$, we see in particular that
\[
  \xi = \theta + o(1) \implies a_{i j}(\xi) = \mathfrak{d}_i(x_j)(\theta) + o(1).
\]
For $r > 0$, set $B_\theta(r) := \{\xi : |\xi - \theta| < r\}$.
Let $A : \mathfrak{p}^\wedge \rightarrow M_{n - 1}(\mathbb{C})$ denote the matrix-valued symbol with entries $a_{i j}$.  By our evaluation \eqref{eqn:frak-d-j-X-tau} and the invertibility of the Vandermonde matrix, we may find a fixed $r > 0$ so that
\begin{equation}\label{eq:det-a_i-jxi_i}
  \det( A(\xi)) \gg 1
\end{equation}
for $\xi \in B_\theta(3 r)$.  We choose a symbol $c \in S^{-\infty}_0(\mathfrak{p}^\wedge)$ supported on $B_\theta(2 r)$ and taking the value $1$ on $B_\theta(r)$, and denote by $C$ the corresponding matrix-valued symbol with entries $1_{i=j} c$.

Let $v \in \mathcal{W}_\delta$.  Since $\mathcal{W}_\delta$ is an $S^\infty_\delta(\mathfrak{p}^\wedge)$-module that is $(\mathfrak{p},\delta)$-localized at $\theta$ and $c-1$ vanishes on $\theta + o(1)$, we see from Theorem \ref{thm:S-infinity-modules-localized-stuff} (part \eqref{item:1-thm:S-infinity-modules-localized-stuff}) that
\begin{equation*}
  \Opp(c) v \equiv v \mod{\h^\infty
    \mathcal{W}_\delta}.
\end{equation*}
Let $B^{(0)} : \mathfrak{p}^\wedge \rightarrow M_{n - 1}(\mathbb{C})$ denote the matrix-valued function supported on $B_\theta(3 r)$ and given there by
\begin{equation*}
  B^{(0)}(\xi) := A(\xi)^{-1} C(\xi).
\end{equation*}
This definition makes sense in view of the lower bound \eqref{eq:det-a_i-jxi_i}, which shows also that the entries $b_{j k}^{(0)}$ ($1 \leq j,k \leq n-1$) lie in $S^{-\infty}_0(\mathfrak{p}^\wedge)$.  These entries are characterized by the following properties:
\begin{itemize}
\item $\supp(b_{j k}^{(0)}) \subseteq \supp(c)$
\item $\sum_j a_{i j}(\xi) b_{j k}^{(0)}(\xi) = 1_{i=k} c(\xi)$ for all $\xi$.
\end{itemize}

For each fixed $\ell \geq 1$, we may inductively construct $b_{jk}^{(\ell)}$ so that we have the formal star product expansion
\begin{equation}\label{eqn:formal-star-prod-matrix}
  \sum_j a_{i j} \star_{\h} \sum_{\ell \geq 0} \h^{\ell} b_{j k}^{(\ell)} \sim  1_{i=k} c,
\end{equation}
where the notation $\sim$ signifies that the indicated expansion holds as a formal series in the variable $\h$ upon expanding $\star_{\h}$ formally into its homogeneous components.  In more detail, starting with the definition $B^{(0)} = A^{-1} C$ (juxtaposition denotes pointwise matrix multiplication), we set (on $B_\theta(3 r)$, extended by zero elsewhere)
\[
  B^{(1)} := - A^{-1} (A \star^1 B^{(0)}),
\]
\[
  B^{(2)} := -A^{-1} ( A \star^2 B^{(0)} + A \star^1 B^{(1)}),
\]
where, e.g., $A \star^1 B^{(0)}$ is the matrix-valued symbol with $(i,k)$th entry $\sum_j a_{i j} \star^1 b_{j k}^{(0)}$.  By \eqref{eq:aijdixj} and \eqref{eq:det-a_i-jxi_i}, we have $b_{j k}^{(\ell)} \in S^{-\infty}_0(\mathfrak{p}^\wedge)$ for each fixed $\ell$. 

We choose $b_{jk} := \sum_{\ell \leq L} \h^{\ell} b_{j k}^{(\ell)}$ for some $L \ggg 1$ chosen small enough that, by Lemma \ref{lem:borel-summation}, we have $b_{j k} \in S^{\infty}_0(\mathfrak{p}^\wedge)$.  Since $b_{j k}$ is supported on the fixed compact set $B_\theta(3 r)$, we in fact have $b_{j k} \in S_0^{- \infty}(\mathfrak{p}^\wedge)$.  By inserting the quantitative asymptotic expansions of $b_{j k}$ (see~\eqref{eq:cui1g9lsvk}) and of the star product (see Theorem~\ref{thm:basic-star-prod}) into the formal expansion \eqref{eqn:formal-star-prod-matrix}, we see that
\begin{equation*}
  \sum_j a_{i j} \star_{\h} b_{j k}
  \equiv 1_{i = k} c \mod{\h^\infty S^{-\infty}(\mathfrak{p}^\wedge)}.
\end{equation*}
It follows that for each $v \in \mathcal{W}_\delta$,
\begin{equation}\label{eqn:approx-invert-operators-aij-bjk}
  \sum_j \Opp(a_{i j})
  \Opp(b_{j k})
  v
  \equiv
  1_{i = k} \Opp(c) v
  \equiv 1_{i = k} v
  \mod{\h^\infty \mathcal{W}_\delta}.
\end{equation}

Let $v \in \mathcal{W}_\delta$ and $k \in \{1, \dotsc, n-1\}$.  To understand $\Opp(e_{n k}) v$, we evaluate in two ways the double sum
\begin{equation*}
  \Sigma :=
  \sum_{i}
  \sum_{j}
  \Opp(e_{n i})
  \Opp(a_{i j})
  \Opp(b_{j k})
  v,
\end{equation*}
where here and henceforth indices run over $\{1, \dotsc, n-1\}$.  By summing first over $i$ and appealing to \eqref{eqn:specialized-non-commutative-cofactor}, we see that
\[
  \Sigma = \sum_j \Opp(f_j) \Opp(b_{j k}) v,
\]
hence $\Sigma \in \mathcal{W}_\delta$.  On the other hand, summing first over $j$ and applying \eqref{eqn:approx-invert-operators-aij-bjk} yields
\begin{equation*}
  \Sigma -
  \Opp(e_{n k})
  v
  \in
  \sum_{i} \Opp(e_{n i})
  \h^\infty \mathcal{W}_\delta
  \subseteq \h^\infty \mathcal{W}_\delta,
\end{equation*}
where in the final step we applied the \emph{a priori} estimates given by Proposition \ref{prop:standard2:each-fixed-x}.
Comparing these two evaluations of $\Sigma$, we deduce that $\Opp(e_{n k}) v \in \mathcal{W}_\delta$, and more precisely that
\[
  \Opp(e_{n k}) v \equiv \sum_{j} \Opp(f_j) \Opp(b_{j k}) v \mod{\h^\infty \mathcal{W}_\delta}.
\]
Setting
\[
  e_{n k}^{\mathfrak{p}} := \sum_{0 \leq s < S} \h^s \sum_j f_j \star^s b_{j k}
\]
for some small enough natural number $S \ggg 1$, we deduce via Borel summation (Lemma \ref{lem:borel-summation}) that $e_{n k} ^{\mathfrak{p}} \in S^{-\infty}_0(\mathfrak{p}^\wedge)$ and, as in the verification of \eqref{eqn:approx-invert-operators-aij-bjk}, that
\begin{equation*}
  \Opp(e_{n k}) v \equiv \Opp(e_{n k}^{\mathfrak{p}}) v \pmod{\h^\infty \mathcal{W}_\delta},
\end{equation*}
as desired.

By comparing with the corresponding ``classical'' calculation of \S\ref{sec:comm-cofact-expans}, we see that $e_{n i}^{\mathfrak{p}}(\theta) = e_{n i}(\tau) + \O(\h)$.  Indeed, tracing through the construction and working modulo $\O(h)$, we see that for $x \in \{x_1,\dotsc,x_{n-1}\}$,
\[
  \sum_{i=1}^{n-1} e_{n i }^{\mathfrak{p}}(\theta) \mathfrak{d}_i(x)(\theta) = \mathcal{P}_{\h \lambda_\pi}(x) - (x + e_{n n}(\tau)) \mathfrak{d}_n(x)(\theta) + \O(\h).
\]
The $e_{n i}(\tau)$ satisfy the same linear relation, so we conclude as in \S\ref{sec:comm-cofact-expans} by the invertibility of the Vandermonde matrix.

The proof of Theorem \ref{thm:g-via-p}, hence that of Theorem \ref{thm:cal-W-delta-g-acts-localized}, is now complete.

\subsection{Control over norms}\label{sec:control-over-norms}

Let $\mathcal{V} \subseteq \pi^\infty$, as in Example \ref{example:class-V-wavelength-h}, denote the class consisting of all vectors $v$ for which $\|\pi(x_1 \dotsb x_m) v\| \ll \h^{-m}$ for all fixed $x_1,\dotsc,x_m \in \mathfrak{g}$.
\begin{theorem}\label{thm:control-over-norms}
  We have
  \begin{equation*}
    \mathcal{W}_\delta \subseteq \h^{-(1-\delta)\dim(B_H)/2} \mathcal{V}.
  \end{equation*}
\end{theorem}
\begin{proof}
  Let $W \in \mathcal{W}_\delta$.  By Lemma \ref{lem:let-w-in-square-integral}, we have
  \begin{equation}\label{eq:int-_n_h-backslash-1}
    \int_{N_H \backslash H} |W|^2 \ll \h^{-(1-\delta) \dim(B_H)}.
  \end{equation}
  By Theorem \ref{thm:cal-W-delta-g-acts-localized}, the required bounds persist for derivatives.  Indeed, our task is to check that for fixed $x_1,\dotsc,x_m \in \mathfrak{g}$, we have
  \begin{equation}\label{eq:int-_n_h-backslash-2}
    \int  _{N_H \backslash H} \lvert
    \pi(x_1 \dotsb x_m) W
    \rvert^2
    \ll \h^{-2 m - (1 - \delta) \dim(B_H)}.
  \end{equation}
  Theorem \ref{thm:cal-W-delta-g-acts-localized} implies in particular (cf.\ Remark \ref{rmk:r-=-0}) that the element $\h^m \pi(x_1 \dotsb x_m) W$ lies in $\mathcal{W}_\delta$.  Applying \eqref{eq:int-_n_h-backslash-1} to this element gives the required estimate \eqref{eq:int-_n_h-backslash-2}.
\end{proof}

\begin{corollary}\label{cor:g-localization-inside-cal-V}
  Let $v \in \mathcal{W}_\delta$.  Then $\h^{(1 - \delta) \dim(B_H)/2} v \in \mathcal{V}$ is $(\mathfrak{g},\delta)$-localized at $\tau$ inside $\mathcal{V}$.
\end{corollary}
\begin{proof}
  Let $p \in S^\infty_0(\mathfrak{g})$ be a polynomial that vanishes to some fixed order $\geq r$ at $\tau$.  We know by Theorem \ref{thm:cal-W-delta-g-acts-localized} that $v$ is $(\mathfrak{g},\delta)$-localized at $\tau$ inside $\mathcal{W}_\delta$, hence $\Opp(p) v \in \h^{r \delta} \mathcal{W}_\delta$.  By Theorem \ref{thm:control-over-norms}, we deduce that
  \begin{equation*}
    \Opp(p) \h^{(1-\delta) \dim(B_H)/2} v \in \h^{r \delta} \mathcal{V},
  \end{equation*}
  as required.
\end{proof}

\subsection{Approximate idempotents}\label{sec:appr-idemp}
For convenience, we recap below some of the running assumptions from before.
\begin{theorem}\label{thm:appr-idemp}
  Let $\pi$ be a generic irreducible unitary representation of the fixed general linear group $G = \GL_n(\mathbb{R})$.  Let $\h \lll 1$ be a wavelength parameter with $\h \lambda_\pi \ll 1$.  Let $\tau := \tau(\pi,\psi,T) \in \mathfrak{g}^\wedge_{\reg}$ be as given by Definition \ref{defn:parameter-tau-pi-T}.  Fix $\delta, \delta ', \delta ''$ with $0 < \delta ' < \delta < 1/2$ and $\delta ' < \delta '' < 2 \delta '$.  Recall the subclasses $\mathcal{W}_\delta$ and $\mathcal{V}$ of $\pi^\infty$ defined in \S\ref{sec:quant-prel} and \S\ref{sec:control-over-norms}.

  There exists $a \in S^\tau_{\delta ', \delta ''}$ with the following properties.
  \begin{enumerate}[(i)]
  \item Let $\chi$ be any $(\tau,\delta ',\delta '')$-admissible cutoff.  Set $\Opp(a) := \Opp_{\h}(a:\chi)$.  Then, for each $v \in \mathcal{W}_\delta$, we have
    \begin{equation}\label{eq:oppa-v-equiv}
      \Opp(a) v \equiv v \mod{\h^\infty \mathcal{V}}.
    \end{equation}
  \item With respect to $\tau$-coordinates as in \S\ref{sec:refined-classes}, we have
    \begin{equation}\label{eqn:a-desired-support-condition}
      \supp(a)
      \subseteq
      \left\{ \xi :
        |\xi ' - \tau '| \ll \h^{\delta '},
        |\xi '' - \tau ''| \ll \h^{\delta ''}
      \right\}.
    \end{equation}
  \end{enumerate}
\end{theorem}
\begin{proof}
  Choose $c \in S^{-\infty}_{\delta '}(\mathfrak{g}^\wedge)$ with $c(\xi) = 1$ for $\xi = \tau + o(\h^{\delta'})$ and $c(\xi) \neq 0$ only for $\xi = \tau + \O(\h^{\delta'})$.  By Corollary \ref{cor:g-localization-inside-cal-V}, the vector $\h^{(1 - \delta) \dim(B_H)/2} v$ is $(\mathfrak{g},\delta)$-localized at $\tau$ inside $\mathcal{V}$.  By Theorem \ref{thm:S-infinity-modules-localized-stuff} (part \eqref{item:1-thm:S-infinity-modules-localized-stuff}) applied to $c-1$, we deduce that
  \begin{equation*}
    \Opp(c) v \equiv v \mod{\h^\infty \mathcal{V}}.
  \end{equation*}
  This is nearly the desired conclusion, except that the support of $c$ is not as required.

  We observe that $c$ belongs to $S^{\tau}_{\delta ', \delta ''}$, since it is supported on elements $\xi$ with $\xi - \tau \ll \h^{\delta '} \lll \h^{\delta '' - \delta '}$ and satisfies stronger-than-required derivative bounds.  We use $\tau$-coordinates $\xi = (\xi ', \xi '') \in \mathfrak{g}_\tau^\perp \times \mathfrak{g}_\tau^{\perp \flat}$.  By choosing an orthonormal basis, we identify $\mathfrak{g}_{\tau}^{\perp\flat}$ with Euclidean space $\mathbb{R}^r$.  We fix $\eta_0 \in C_c^\infty(\mathbb{R}^r) \cong C_c^\infty(\mathfrak{g}_{\tau}^{\perp\flat})$ so that $\eta_0(u)=1$ for $|u| \leq 1$ and $\eta_0(u)=0$ for $|u| \geq 2$.  We set
  \[
    \eta(\xi) := \eta_0\!\left(\frac{\xi '' - \tau ''}{\h^{\delta ''}}\right), \quad
    a := \eta c, \quad b := (1-\eta)c.
  \]
  Then $c = a + b$, where $a$ and $b$ satisfy the support conditions
  \begin{equation*}
    \supp(a) \subseteq \{|\xi' - \tau'| \ll \h^{\delta'},\, |\xi'' - \tau''| \ll \h^{\delta''}\},
  \end{equation*}
  \begin{equation*}
    \supp(b) \subseteq \{|\xi' - \tau'| \ll \h^{\delta'},\, |\xi'' - \tau''| \gg \h^{\delta''}\}.
  \end{equation*}

  We claim that $a, b \in S^\tau_{\delta', \delta ''}$.  Indeed, since $a,b$ are supported on $\supp(c) \subseteq \tau + \O(\h^{\delta'}) \subseteq \tau + \O(\h^{\delta''-\delta'})$, we need only check the derivative bounds, which follow from the hypothesis $c \in S^{-\infty}_{\delta'}$, the construction of $\eta$, and the product rule.

  By Theorem \ref{thm:local-refin-symb}, we have
  \begin{equation*}
    \Opp(b) \in \h^\infty \Psi^{-\infty}.
  \end{equation*}
  By the easy inclusion
  \begin{equation}\label{eq:hinfty-psi-infty}
    \h^\infty \Psi^{-\infty} \mathcal{V} \subseteq \h^\infty \mathcal{V},
  \end{equation}
  verified below for completeness, we conclude as required that \eqref{eq:oppa-v-equiv} holds.
\end{proof}
\begin{proof}[Verification of \eqref{eq:hinfty-psi-infty}]
  Let $T \in \h^\infty \Psi^{-\infty}$ and $v \in \mathcal{V}$.  We must check that $T v \in \h^\infty \mathcal{V}$, i.e., that
  \begin{itemize}
  \item $T v$ lies in the space $\pi^\infty$ of smooth vectors, and
  \item for each fixed $N$ and $x_1,\dotsc,x_m \in \mathfrak{g}$, we have
    \begin{equation}\label{eq:pix_1-dotsb-x_m-1}
      \|\pi(x_1 \dotsb x_m) T v \| \ll \h^N.
    \end{equation}
  \end{itemize}
  By the definition of $\h^\infty \Psi^{-\infty}$ (\S\ref{sec:negligible-operators}), we have $T \in \underline{\Psi }^{\infty}$ and $T \in \h^N \Psi^{-\infty}$.  Recalling from Lemma \ref{lem:foll-incl-hold} that $\underline{\Psi }^{\infty} \pi^\infty \subseteq \pi^\infty$, we deduce that $T v \in \pi^\infty$.  By the definition of $\Psi^{-\infty}$, the operator norm of $\pi(x_1 \dotsb x_m) T$ is then $\O(\h^N)$.  Since $v \in \mathcal{V}$, we have in particular $\|v\| \ll 1$.  The required estimate \eqref{eq:pix_1-dotsb-x_m-1} follows.
\end{proof}

\subsection{Proof of Theorem \ref{thm:big-group-whittaker-behavior-under-G}}\label{sec:proof-theor-refthm:b}
Suppose given $T$, $\pi$, $\eps$ and $\delta$ as in the hypotheses of that theorem.  Then $\h := 1/T \lll 1$ and $\h \lambda_\pi = \O(1)$.  Let $\tau := \tau(\pi,\psi,T) \in \mathfrak{g}^\wedge_{\reg}$ as in Definition \ref{defn:parameter-tau-pi-T}, as usual.

Fix $\delta ', \delta ''$ as in the hypotheses of Theorem \ref{thm:appr-idemp}, with $(\delta ', \delta '')$ sufficiently close to $(\delta , 2 \delta)$; recalling that $\delta$ is chosen sufficiently close to $1/2$ relative to $\eps$, it follows that $(\delta ', \delta '')$ is close to $(1/2,1)$ relative to $\eps$.  Let $a \in S^{\tau}_{\delta ', \delta ''}$ be as furnished by the conclusion of that theorem.  Fix $\delta_+' \in (\delta ', 1/2)$ and $\delta_+'' \in (\delta '', 1)$.  Let $\chi \in \mathcal{X}_{\tau, \delta_+', \delta_+''}(\mathfrak{g})$.  Define
\begin{equation*}
  \omega_0 := \widetilde{\Opp}_{\h}(a:\chi),
\end{equation*}
a smooth compactly-supported measure on $G$; we identify it with a smooth compactly-supported function on $G$ by dividing by our fixed Haar measure.

Recalling the statement of Theorem \ref{thm:big-group-whittaker-behavior-under-G}, we must verify the following assertions concerning $\omega = \omega_0 \ast \omega_0^*$ and $\omega^\sharp = \int_{z \in Z} \pi|_Z (z) \omega(z g) \, d z$:
\begin{enumerate}[(i)]
\item \label{item:kirillov-final-verify-1} $\omega(g) \neq 0$ only if $g = 1 + o(1)$ and $\Ad^*(g) \tau = \tau + o(T^{-1/2})$.
\item \label{item:kirillov-final-verify-2} $\|\omega^\sharp\|_{\infty} \ll T^{n(n+1)/2+\eps}$.
\item \label{item:kirillov-final-verify-3} If $\tau$ lies in some fixed compact subset of $\mathfrak{g}^\wedge_{\stab}$, then $\int_H |\omega^\sharp| \ll T^{n/2 + \eps}$.
\item \label{item:kirillov-final-verify-4} For each $W \in \mathcal{W}_\delta = \mathfrak{W}(\pi,\psi,T,\delta)$ and fixed $x \in \mathfrak{U}(G)$, we have $\|\pi(x) (\pi(\omega_0)W - W)\| \ll T^{-\infty}$.
\end{enumerate}
We have $\Opp(a) = \pi(\omega_0)$, so assertion \eqref{item:kirillov-final-verify-4} is immediate from the conclusion $\Opp(a) W \equiv W \pmod{\h^\infty \mathcal{V}}$ of Theorem \ref{thm:appr-idemp} and the definition of the class $\mathcal{V}$.  The remaining assertions may be verified exactly as in \cite[\S14.9]{2020arXiv201202187N} (noting that our $\chi \in \mathcal{X}_{\tau, \delta_+', \delta_+''}(\mathfrak{g})$ plays the role of the modified cutoff ``$\tilde{\chi}$'' in the cited reference):
\begin{itemize}
\item Assertion \eqref{item:kirillov-final-verify-1} is given by \cite[(14.13)]{2020arXiv201202187N}.  In short, the support conditions on $\chi$ imply that for $g=\exp(x)$ in the support of $\omega_0$, we have $x' \ll \h^{1-\delta_+'}$ and $x'' \ll \h^{1-\delta_+''}$.  This support property is stable under $(x, y) \mapsto \log(\exp(x) \exp(y)^{-1})$, so it also holds on the support of $\omega = \omega_0 \ast \omega_0^*$.  Using that $\delta_+ ' < 1/2$ and $\delta_+'' < 1$, we deduce that $g = 1 + O(\h^\nu) = 1 + o(1)$ and $x ' \ll \h^{1/2 + \nu}$ for some fixed $\nu > 0$.  From this last estimate and a simple calculation with the exponential series (see \cite[Lemma 13.9]{2020arXiv201202187N}), we obtain further that $\Ad^*(g) \tau = \tau + \O(\h^{1/2 + \nu}) = \tau + o(T^{-1/2})$.
\item Assertion \eqref{item:kirillov-final-verify-2} is given by \cite[(14.14)]{2020arXiv201202187N}; for convenience, we record the proof.  The key estimate (compare \cite[(14.6)]{2020arXiv201202187N}) is obtained by Fourier inversion and repeated integration by parts, using the support properties \eqref{eqn:a-desired-support-condition} of $a$ and the defining derivative bounds \eqref{eq:partial_xi--alpha} for $S_{\delta ', \delta ''}^{\tau}$: in $\tau$-coordinates $x = (x',x'')$,
  \[
    |a_{\h}^{\vee}(x)|
    \ll
    \frac{
      \h^{(-1+\delta')(\dim(\mathfrak{g})-\rank(\mathfrak{g}))}
      \h^{(-1+\delta'')\rank(\mathfrak{g})}
    }{
      \langle \h^{-1+\delta'} x' \rangle^N
      \langle \h^{-1+\delta''} x'' \rangle^N
    }
  \]
  for each fixed $N \geq 0$.  Since $\chi \in \mathcal{X}_{\tau,\delta_+',\delta_+''}(\mathfrak{g})$ is supported on
  \[
    |x'| \ll \h^{1-\delta_+'},
    \quad
    |x''| \ll \h^{1-\delta_+''},
  \]
  with $\delta_+' > \delta'$ and $\delta_+'' > \delta''$, it follows that
  \[
    \|\omega_0\|_{L^1}
    \ll
    \h^{(\delta_+'-\delta')(\dim(\mathfrak{g})-\rank(\mathfrak{g}))}
    \h^{(\delta_+''-\delta'')\rank(\mathfrak{g})}
    \ll 1,
  \]
  \[
    \|\omega_0\|_{L^\infty}
    \ll
    \h^{(-1+\delta')(\dim(\mathfrak{g})-\rank(\mathfrak{g}))}
    \h^{(-1+\delta'')\rank(\mathfrak{g})}
    \ll T^{n(n+1)/2+\eps},
  \]
  using that $(\delta',\delta'')$ is sufficiently close to $(1/2,1)$ and $\dim(\mathfrak{g})-\rank(\mathfrak{g}) = n(n+1)$.  Thus
  \[
    \|\omega\|_{L^\infty}
    \leq
    \|\omega_0\|_{L^1}\|\omega_0\|_{L^\infty}
    \ll T^{n(n+1)/2+\eps}.
  \]
  Since $\supp(\omega)$ lies in a fixed compact subset of $G$, we have $\vol\{z \in Z : zg \in \supp(\omega)\} \ll 1$ for all $g \in G$, and so
  \[
    \|\omega^\sharp\|_{L^\infty}
    \leq
    \sup_{g \in G}
    \int_Z |\omega(zg)| \, d z
    \ll
    \|\omega\|_{L^\infty}.
  \]
  The claimed bound follows.
\item Assertion \eqref{item:kirillov-final-verify-3} is established in the final paragraph of \cite[\S14.9]{2020arXiv201202187N}.  The proof uses that the natural map $\mathfrak{g}^\wedge_{\stab} \rightarrow [\mathfrak{g}^\wedge] \times [\mathfrak{h}^\wedge]$ defines a principal $H$-bundle over its image, which gives
  \begin{equation*}
    V_\tau := \vol\{h \in H : \dist(\Ad^*(h)\tau,\tau) \leq T^{-1/2}\} \ll T^{-n^2/2}.
  \end{equation*}
  Combining this with (i) and (ii) gives
  \[
    \int_H |\omega^\sharp|
    \ll
    \|\omega^\sharp\|_{L^\infty}
    V_\tau \ll T^{n(n+1)/2+\eps} T^{-n^2/2}
    = T^{n/2+\eps},
  \]
  as required.
\end{itemize}

\part{Growth bounds for Eisenstein series}\label{part:local-l2-growth}
The aim of Part \ref{part:local-l2-growth} is to prove Theorem \ref{thm:growth-eisenstein-nonstandard}.  Throughout this section, $G$ denotes a general linear group over either $\mathbb{Z}$, $\mathbb{R}$, or a non-archimedean local field.  Recall from \S\ref{sec:gener-prel} the standard accompanying notation: $B = A N$ denotes the standard upper-triangular Borel subgroup of $G$, while $Q = A U$ denotes another Borel subgroup that contains the diagonal subgroup $A$, but is otherwise free for us to choose according to computational convenience.  For the purposes of this section, it is most convenient to take
\begin{equation*}
  U := N, \text{ so that } Q = B.
\end{equation*}
We accordingly write simply $\Delta$, $\rho$ and $\rho^\vee$ for $\Delta_N, \rho_N$ and $\rho_N^\vee$.  The meaning of ``dominant'' is likewise unambiguous now.  Recall also that $W$ (or $W_G$, to emphasize $G$) denotes the Weyl group for $(G,A)$.  Also, in the local case, $K$ denotes the standard maximal compact subgroup of $G$, while in the global case, $K = \prod_\mathfrak{p} K_\mathfrak{p}$ is the standard maximal compact subgroup of $G(\mathbb{A})$.

We employ the abbreviation
\begin{equation*}
  \left\langle s \right\rangle := 1 + |s|
\end{equation*}
for $s \in \mathfrak{a}_{\mathbb{C}}^*$.

The organization is as follows.
\begin{itemize}
\item Recall that Theorem \ref{thm:growth-eisenstein-nonstandard} provides nontrivial estimates for the local $L^2$-norms of Eisenstein series on $\GL_n$.  There are two sources of complication here: that we are working in higher rank, and that we are working with Eisenstein series.  In \S\ref{sec:simpler-variants}, we illustrate the proof scheme by applying it to simpler examples -- cusp forms on $\GL_n$, Eisenstein series on $\GL_2$ -- obtained by  ``turning off'' either of these difficulties.  The general proof is not quite a combination of the proofs of these special cases, but the latter do illustrate several key features of the former.
\item
  In \S\ref{sec:rank-selb-unfold}, we unfold the integrals $\int_{[\GL_n]} |\varphi|^2 \theta$, where $\varphi$ (resp.\ $\theta$) is a pseudo Eisenstein series attached to a minimal parabolic (resp.\ mirabolic) subgroup.
\item
  In \S\ref{sec:reduction-proof}, we develop the crucial global arguments that reduce the proof of Theorem \ref{thm:growth-eisenstein-nonstandard} to some technical problems.  These arguments include the construction of a mirabolic pseudo Eisenstein series $\theta$ taking large values in a given region in $[\GL_n]$, a crucial ``duality trick,'' and several alternative series representations of our minimal pseudo Eisenstein series carefully tailored to its estimation.  The resulting technical problems are then treated in \S\ref{sec:local-estimates} (local) and \S\ref{sec:completion-proof-1} (global).
\end{itemize}

\section{Simpler variants}\label{sec:simpler-variants}
As noted above, the purpose of this section is to illustrate the basic proof scheme of Part \ref{part:local-l2-growth} by application to some simpler variants of the main problem.  This section is not necessary for the logical purposes of the paper, but may serve as a useful guide to the arguments that follow.

\subsection{Cusp forms on $\GL_n$}\label{sec:simpler-cusp-forms-gl_n}
We refer here to \cite{MR3468028} for general background.  Let $\Psi \in L^2(\SL_n(\mathbb{Z}) \backslash \SL_n(\mathbb{R}))$ be a Hecke--Maass cusp form, i.e., a spherical vector in an everywhere-unramified cuspidal automorphic representation of $\PGL_n$ over $\mathbb{Q}$, viewed classically.  Its Fourier expansion has the shape
\begin{equation*}
  \Psi(g) = \sum_{m}
  \frac{\lambda (m)}{\delta_{N_n}^{1/2}(m)}
  \sum_{\gamma \in N_n(\mathbb{Z}) \backslash P_n(\mathbb{Z})}
  W(m \gamma g),
\end{equation*}
where $P_n \leq \GL_n$ denotes the mirabolic subgroup consisting of matrices with representatives having bottom row of the form $(0,\dotsc,0,1)$,  $N_n \leq P_n$ denotes the subgroup of upper-triangular unipotent matrices, $m$ runs over diagonal matrices of the form
\begin{equation*}
  m =
  \begin{pmatrix}
    m_1 \dotsb m_{n-1} & 0 & 0 & 0 \\
    0 & m_2 \dotsb m_{n-1} & 0 & 0 \\
    0 & 0 & \ddots & 0 \\
    0 & 0 & 0 & 1
  \end{pmatrix},
  \quad
  m_1, \dotsc, m_{n-1} \in \mathbb{Z}_{\geq 1},
\end{equation*}
$\delta_{N_n}$ is the modulus character for the upper-triangular Borel, and $\lambda$ denotes the normalized Hecke eigenvalue.  Set
\begin{equation*}
  \|\Psi \|^2 := \int_{\SL_n(\mathbb{Z}) \backslash \SL_n(\mathbb{R})} |\Psi|^2,
  \quad
  \|W\|^2 := \int_{N_{n-1}(\mathbb{R}) \backslash \GL_{n-1}(\mathbb{R})} |W|^2.
\end{equation*}
Let us assume that each archimedean parameter of $\pi$ has size $\ll T$.  We say that $\Psi$ is \emph{Whittaker $L^2$-normalized} if $\|W\|^2 = 1$.  In that case, we know (Lemma \ref{lem:nonstandard-norm-whittaker}) that
\begin{equation}\label{eq:psi-2-ll}
  \|\Psi \|^2 \ll T^{o(1)}.
\end{equation}

Let $\Omega$ be a fixed compact subset of $\SL_n(\mathbb{R})$.  Let $A_n \leq \SL_n(\mathbb{R})$ denote the subgroup of diagonal matrices, and let $t = \diag(t_1,\dotsc,t_n)$ be a dominant element of $A_n \cap \SL_n(\mathbb{R})$ with positive entries, i.e.,
\begin{equation*}
  t_1 \geq \dotsb \geq t_n > 0.
\end{equation*}
We set
\begin{equation*}
  \|\Psi \|_t^2 := \int_{g \in \Omega } |\Psi|^2(t g) \, d g.
\end{equation*}
By a simple argument with Siegel domains (Proposition \ref{prop:siegel-lower}), we obtain  the ``trivial'' bound
\begin{equation}\label{eq:psi-_r2-ll-2}
  \|\Psi \|_t^2 \ll \delta_{N_n}(t) \|\Psi \|^2.
\end{equation}
We improve upon this in the following cuspidal analogue of Theorem \ref{thm:growth-eisenstein-nonstandard}.
\begin{theorem}\label{thm:sub-gln:let-psi-be}
  Let $\Psi$ be as above and Whittaker $L^2$-normalized.  Let $t$ be as above.  Then
  \begin{equation}\label{eq:psi-_r2-ll-1}
    \|\Psi \|_t^2 \ll \min(t_1^{-1}, t_n)^n \delta_{N_n}(t) T^{o(1)}.
  \end{equation}
\end{theorem}
\begin{remark}
  We note that, since $t_1 \geq \dotsb \geq t_n$, we have $\min(t_1^{-1}, t_n) \leq 1$, and so \eqref{eq:psi-_r2-ll-1} actually does improve upon the ``trivial'' bound of $\ll \delta_{N_n}(t) T^{o(1)}$ provided that $t$ is not too close to the identity.
\end{remark}
\begin{remark}
  When $n = 2$, the bound \eqref{eq:psi-_r2-ll-1} reads $\lVert \Psi \rVert_t^2 \ll T^{o(1)}$ and is essentially optimal.  For $n \geq 3$, we do not know whether \eqref{eq:psi-_r2-ll-1} is optimal.  In any event, we recall from \S\ref{sec:high-level-overview} that \eqref{eq:psi-_r2-ll-1} serves our purposes well (see especially \eqref{eq:t-in-mathcalt} and \eqref{eq:cujyyxote2}).
\end{remark}
\begin{remark}
  We have referred to the above as a cuspidal analogue of Theorem \ref{thm:growth-eisenstein-nonstandard}, but the results are not quite analogous: Theorem \ref{thm:growth-eisenstein-nonstandard} concerns high frequency wave packets concentrated on representations of essentially bounded parameter, while in Theorem \ref{thm:sub-gln:let-psi-be}, the parameters of the representations vary considerably.  However, our arguments are insensitive to such differences.  We work here with a spherical vector to keep the presentation familiar and simple.
\end{remark}
Let us warm-up by considering the simplest case:
\begin{example}\label{example:growth-bound-sl2-cusp}
  Suppose for the moment that $n = 2$, so that $\Psi$ corresponds to a classical Maass form on the modular surface $\SL_2(\mathbb{Z}) \backslash \mathbb{H}$.  Writing
  \[
    t = \diag(t_1, t_2) = \diag(r^{1/2}, r^{-1/2}),
    \quad r \geq 1,
  \]
  we have $\delta_{N_2}(t) = r$.  By considering Siegel domains, we see that if $\Omega$ is large enough, then we may fix $c_1 > c_0 > 1$ so that
  \begin{equation*}
    I_{c_0}(t) \leq \lVert \Psi \rVert_t^2 \leq I_{c_1}(t),
  \end{equation*}
  where for $c > 1$,
  \begin{equation*}
    I_c(t) := \int_{x \in \mathbb{R} / \mathbb{Z}} \int_{y = r/c}^{c r} \lvert \Psi(x + i y) \rvert^2 \, \frac{d x \, d y}{y}.
  \end{equation*}
  Since $\min(t_1^{-1}, t_2)^2 \delta_{N_2}(t) = 1$, the content of Theorem \ref{thm:sub-gln:let-psi-be} is that for fixed $c$, we have
  \begin{equation}\label{eq:int-_x-in}
    I_c(t)
    \ll T^{o(1)}.
  \end{equation}
  This improves upon the ``trivial bound'' $I_c(t) \ll r T^{o(1)}$ that follows from the $L^2$-normalization, namely
  \begin{equation*}
    I_c(t) \ll r \int_{x \in \mathbb{R} / \mathbb{Z}}
    \int_{y = r/c}^{c r} \lvert \Psi(x + i y) \rvert^2 \, \frac{d x \, d y}{y^2}
    \ll r \lVert \Psi \rVert^2 \ll r T^{o(1)}.
  \end{equation*}
  The proof of \eqref{eq:int-_x-in} is an easy exercise, as we now explain.  Inserting the Fourier expansion and applying the Parseval relation gives
  \begin{equation*}
    I_c(t)
    =
    \sum_{0 \neq m \in \mathbb{Z} }
    \frac{|\lambda(m)|^2}{|m|}
    \int_{y = r/c}^{c r} |W(m y)|^2 \, d^\times y.
  \end{equation*}
  Let $S \geq 2$ be a parameter.  By applying the trivial estimate
  \begin{equation*}
    \int_{y=r/c}^{c r} |W(m y)|^2 \, d^\times y \leq \int_{\mathbb{R}^\times} |W|^2 = 1
  \end{equation*}
  when $|m| \leq S$, we obtain the majorization $I_c(t) \leq \Sigma_1 + \Sigma_2$, where
  \begin{equation*}
    \Sigma_1 :=  \sum_{0 \neq |m| \leq S}
    \frac{|\lambda(m)|^2}{|m|},
    \quad
    \Sigma_2 :=
    \sum_{|m| > S}
    \frac{|\lambda(m)|^2}{|m|}
    \int_{y = r/c}^{c r} |W(m y)|^2 \, d^\times y.
  \end{equation*}
  Rankin--Selberg theory and the estimates of Iwaniec \cite{MR1067982} give $\Sigma_1 \ll T^{o(1)} \log(S)$.  This estimate remains acceptable even if we take, say, $S = T^{1000}$.  In that case, the rapid decay of the Whittaker function $W$ at arguments larger than $T^{1+\eps}$ implies that $\Sigma_2 \ll T^{-100}$, say.  The required estimate follows.
\end{example}
We now turn to the general case:
\begin{proof}[Proof of Theorem \ref{thm:sub-gln:let-psi-be}]
  We begin with a key observation: it suffices to prove the apparently weaker bound
  \begin{equation}\label{eq:psi-_r2-ll-3}
    \|\Psi \|_t^2 \ll t_n^n \delta_{N_n}(t) T^{o(1)}.
  \end{equation}
  To obtain this reduction, we observe that the dual form
  \begin{equation*}
    \tilde{\Psi}(g) := \Psi(g^{-\transpose}),
  \end{equation*}
  obtained from $\Psi$ by applying the outer automorphism of $\GL_n$ defined by taking the inverse transpose, satisfies the same hypotheses as $\Psi$.  On the other hand, writing $w \in \SL_n(\mathbb{Z})$ for a long Weyl element, we have
  \begin{equation*}
    w t^{-\transpose} w = \diag(t_n^{-1}, \dotsc, t_1^{-1}).
  \end{equation*}
  By enlarging $\Omega$ suitably, we may and shall assume that it is invariant under left multiplication by $w$ and also under $g \mapsto g^{-\transpose}$.  Then
  \begin{equation*}
    \|{\tilde{\Psi}}\|_{w t^{-\transpose} w}^2
    =
    \int_{g \in \Omega }
    |{\tilde{\Psi}}|^2(w t^{-\transpose} w g)
    \, d g
    =
    \int_{g \in \Omega }
    |{\Psi}|^2(t (w g)^{-\transpose})
    \, d g
    =
    \|{\Psi}\|_{t}^2,
  \end{equation*}
  so if \eqref{eq:psi-_r2-ll-3} holds for $\tilde{\Psi}$, then
  \begin{equation*}
    \|\Psi \|_t^2 \ll
    ( w t^{-\transpose} w)_n^n \delta_{N_n}(w t^{-\transpose} w) T^{o(1)}
    =
    t_1^{-n} \delta_{N_n}(t) T^{o(1)}.
  \end{equation*}
  Thus if we can prove that \eqref{eq:psi-_r2-ll-3} holds for both $\Psi$ and $\tilde{\Psi}$, then \eqref{eq:psi-_r2-ll-1} will follow.

  In verifying \eqref{eq:psi-_r2-ll-3}, we may and shall assume that $t_n= T^{\O(1)}$, as otherwise the required estimates follow (in stronger form) from standard Sobolev-type bounds as in \cite[\S2]{michel-2009}.

  Our strategy is to construct a mirabolic pseudo Eisenstein series which
  \begin{itemize}
  \item takes a reasonably large value in a neighborhood of $t$, and yet
  \item whose integral against $|\Psi|^2$ is not too large.
  \end{itemize}
  To that end, we fix $C_0 \geq 2$ large enough in terms of $\Omega$ and construct an even Schwartz function $\phi$ on $\mathbb{R}^n$ with the following properties:
  \begin{enumerate}[(i)]
  \item $\phi$ is spherical, i.e., $\phi(v)$ depends only upon the Euclidean norm $\|v\|$.
  \item $0 \leq \phi \leq 1$.
  \item $\phi(v) = 1$ if $\|v\| / t_n \in (1/C_0, C_0)$.
  \item $\phi(v) = 0$ if $\|v\| / t_n \notin (1/2 C_0, 2 C_0)$.
  \end{enumerate}
  In particular, writing $e_n = (0,\dotsc,0,1)$ for the standard basis vector,
  \begin{equation}
    \phi(t_n e_n) = 1.
  \end{equation}
  We define
  \begin{equation*}
    \theta_\phi(g) = \sum_{\gamma \in P_n(\mathbb{Z}) \backslash \GL_n(\mathbb{Z})} \phi(e_n \gamma g)
    =
    \sum_{\gamma \in P_n(\mathbb{Z}) \cap \SL_n(\mathbb{Z}) \backslash \SL_n(\mathbb{Z})} \phi(e_n \gamma g).
  \end{equation*}
  Equivalently,
  \begin{equation*}
    \theta_\phi(g) = \sideset{}{^*}\sum_{v \in \mathbb{Z}^n}\phi(v g),
  \end{equation*}
  where $\sum^*$ indicates that the sum is restricted to primitive vectors.  Clearly
  \begin{equation*}
    \theta_\phi(t) \geq \phi(e_n t) = \phi(t_n e_n) = 1.
  \end{equation*}
  Since $C_0$ was chosen large enough in terms of $\Omega$, we have in fact $\theta_\phi(t g) \geq 1$ for all $g \in \Omega$.  Taking into account the formula for Haar measure on a Siegel domain (Proposition \ref{prop:siegel-lower}), we deduce that
  \begin{equation*}
    \frac{\|\Psi \|_t^2 }{\delta_{N_n}(t)} \ll
    \int_{\SL_n(\mathbb{Z}) \backslash \SL_n(\mathbb{R})} \theta_\phi |\Psi|^2.
  \end{equation*}
  Our task thereby reduces to verifying that
  \begin{equation}\label{eq:int_sl_nm-backsl-sl}
    \int_{\SL_n(\mathbb{Z}) \backslash \SL_n(\mathbb{R})} \theta_\phi |\Psi|^2 \ll t_n^n T^{o(1)}.
  \end{equation}
  A standard unfolding calculation shows that the LHS of \eqref{eq:int_sl_nm-backsl-sl} is a constant multiple of
  \begin{equation}\label{eq:sum-_m-fracl}
    \sum_{m}
    \frac{|\lambda(m)|^2}{\delta_{N_n}(m)}
    \int_{g \in N_n(\mathbb{R}) \backslash \SL_n(\mathbb{R})} \phi(e_n g)
    |W|^2(m  g) \, d g.
  \end{equation}
  We may explicate this last integral using the Iwasawa decomposition.  Write $A_{n-1}$ for the diagonal subgroup of $\GL_{n-1}(\mathbb{R})$.  For $a \in A_{n - 1}$, write $\underline{a} := \diag(a, \det(a)^{-1}) \in \mathrm{SL}_n(\mathbb{R})$.  Then with suitable measure normalization, the above integral is
  \begin{equation*}
    \int_{a \in A_{n-1}}
    \phi(e_n \det(a)^{-1})
    |W|^2(m \underline{a}) \, \frac{d a }{ \delta_{N_n}(\underline{a})}.
  \end{equation*}
  Using the identity
  \begin{equation*}
    \delta_{N_n}(\underline{a}) = |\det a|^{n} \delta_{N_{n-1}}(a),
  \end{equation*}
  we may rewrite this last integral as
  \begin{equation*}
    \int_{a \in A_{n-1}}
    |\det a|^{-n}
    \phi(e_n \det(a)^{-1})
    |W|^2(m \underline{a}) \, \frac{d a }{ \delta_{N_{n-1}}(a)}.
  \end{equation*}
  On the support of $\phi (e_n \det(a)^{-1})$, we have $\det(a) \asymp t_n^{-1}$, so the above integral is majorized by
  \begin{equation*}
    t_n^{n} \int_{a \in A_{n-1}} |W|^2(m \underline{a}) \, \frac{d a }{ \delta_{N_{n-1}}(a) }.
  \end{equation*}
  We now substitute $a \mapsto m^{-1} a$, which introduces the factor $\delta_{N_{n-1}}(m) = (\det m)^{-1} \delta_{N_{n}}(m)$.  The remaining integral
  \begin{equation*}
    \int_{a \in A_{n-1}} |W|^2(\underline{a}) \, \frac{d a}{\delta_{N_{n-1}}(a)}
  \end{equation*}
  is $\asymp 1$ thanks to the Whittaker $L^2$-normalization of $\Psi$.  We thereby obtain the bound
  \begin{equation*}
    \int  _{\SL_n(\mathbb{Z}) \backslash \SL_n(\mathbb{R})}
    \theta_\phi |\Psi|^2
    \ll
    t_n^{n}
    \sum_{m}
    \frac{|\lambda(m)|^2}{ \det m }.
  \end{equation*}
  Since $\det m = m_1 m_2^2 m_3^3 \dotsb$ and $\sum_{m_1 \in \mathbb{Z}_{\geq 1}} 1/m_1 = \infty$, the RHS of this estimate diverges, but only logarithmically; by truncating the $m$-sum as in Example \ref{example:growth-bound-sl2-cusp} and using \cite[Thm 2]{MR4093914} to estimate the Whittaker functions, we obtain an acceptable estimate.
\end{proof}

\subsection{Eisenstein series on $\GL_2$}\label{sec:cnjgglqbxq}
One can apply the arguments of \S\ref{sec:simpler-cusp-forms-gl_n} to a pseudo Eisenstein series $\Psi$.  The reduction to \eqref{eq:int_sl_nm-backsl-sl} proceeds exactly as in the cuspidal case.  The subsequent unfolding procedure produces a ``Whittaker'' term, like \eqref{eq:sum-_m-fracl}, together with $n-1$ additional ``degenerate'' terms, indexed by the standard maximal parabolic subgroups of $\GL_n$.  The estimation of the degenerate terms is a bit more involved.  In the simplest case $n=2$, there is a single degenerate term, coming from the constant term of the Eisenstein series.  The estimation is painless in that case, and serves as a poor guide to the general case, but nevertheless provides an opportunity to illustrate a few key features.

For clarity of presentation, we eschew the general notation and setup of Theorem \ref{thm:growth-eisenstein-nonstandard} (concerning seminorms $\nu_{\mathcal{D},\ell,T}$, completed pseudo Eisenstein series $\Phi[f_\infty]$, etc.), and restate the basic content of the ``$\SL_2$ variant'' of that estimate explicitly.

\begin{theorem}\label{thm:sub-gln:let-f-be}
  Let $T$ be a positive parameter with $T \ggg 1$.  Let $f$ be an even Schwartz function on $\mathbb{R}^2$ that, together with all derivatives, decays rapidly near the origin.  Regard $\mathbb{R}^2$ as a space of row vectors, equipped with the usual right action of $\SL_2(\mathbb{R})$.  Define $\Psi : \SL_2(\mathbb{Z}) \backslash \SL_2(\mathbb{R}) \rightarrow \mathbb{C}$ by the formula
  \begin{equation*}
    \Psi(g) :=
    \sum_{v \in \mathbb{Z}^2} f(v g)
    =
    \sum_{0 \neq v \in \mathbb{Z}^2} f(v g).
  \end{equation*}
  Assume the following conditions on $f$:
  \begin{enumerate}[(i)]
  \item \label{itm:sketch-sl2-assump-1} $f$ is equal to its symplectic Fourier transform, i.e., $\mathcal{F} f = f$, where
    \begin{equation*}
      \mathcal{F} f(x,y) := \int_{u, v \in \mathbb{R} } f(u,v ) e(x v - y u) \, d u \, d v, \quad  e(x) := e^{2 \pi i x}.
    \end{equation*}
  \item \label{itm:sketch-sl2-assump-2} $f$ has vanishing integral along each line through the origin, i.e., for all $0 \neq v \in \mathbb{R}^2$,
    \begin{equation}\label{eq:int-_t-in-1}
      \int_{t \in \mathbb{R} } f( t v) \, d t = 0.
    \end{equation}
  \item \label{itm:sketch-sl2-assump-3} $f$  has
    \begin{center}
      ``$L^2$-norm $\O(1)$, radial support $\asymp T^{1/2}$, radial frequency $\O(1)$ and angular frequency $\O(T)$''
    \end{center}
    in the following sense.  For elements $u$ of the universal enveloping algebra $\mathfrak{U}(\SL_2(\mathbb{R}))$, let $R(u)$ denote the associated differential operator on $C^\infty(\mathbb{R}^2)$.  Then for all fixed $m, \alpha \in \mathbb{Z}_{\geq 0}$, fixed $u \in \mathfrak{U}(\SL_2(\mathbb{R}))$ of degree at most $m$, and all $r > 0$, we have
    \begin{equation}\label{eq:int-_v-in-7}
      \int_{\theta \in \mathbb{R} / 2 \pi \mathbb{Z} }
      \left\lvert
        \frac{d^\alpha }{d z^\alpha } R(u) f(e^{z} r \cos \theta, e^z r \sin \theta )|_{z=0}
      \right\rvert^2 \, d \theta  \ll r^{-2} T^{2 m + o(1)} \max \left( \frac{r}{T^{1/2}}, \frac{T^{1/2}}{r} \right)^{-\infty}.
    \end{equation}
  \end{enumerate}
  Let $\Omega$ be a fixed compact subset of $\SL_2(\mathbb{R})$.  Let $t \geq 1$ with $t = T^{\O(1)}$.  Then
  \begin{equation}\label{eq:int-_g-in-1}
    \int_{g \in \Omega } |\Psi|^2 (\diag(t, 1/t) g) \, d g
    \ll
    T^{o(1)}.
  \end{equation}
\end{theorem}

\begin{remark}
  Using only the general $L^2$-theory of Eisenstein series for $\SL_2(\mathbb{Z})$ and our assumptions \eqref{itm:sketch-sl2-assump-1}, \eqref{itm:sketch-sl2-assump-2} and \eqref{itm:sketch-sl2-assump-3}, we may verify that $\int_{\SL_2(\mathbb{Z}) \backslash \SL_2(\mathbb{R}) } |\Psi|^2 \ll T^{o(1)}$.  This bound is essentially optimal (i.e., for some $f$ as above, it is sharp).  Arguing as in Example \ref{example:growth-bound-sl2-cusp}, we obtain the ``trivial bound''
  \begin{equation}\label{eq:int-_g-in-4}
    \int_{g \in \Omega } |\Psi|^2 (\diag(t, 1/t) g) \, d g
    \ll
    t^2 T^{o(1)},
  \end{equation}
  which the target bound \eqref{eq:int-_g-in-1} improves upon by the factor $t^2$.
\end{remark}
\begin{remark}
  The invariance condition $\mathcal{F} f = f$ under the symplectic Fourier transform, together with the assumed rapid decay of $f$ at $0$, implies that the Eisenstein series $\Psi$ is orthogonal to residual spectra, i.e., the constant function.  This feature is necessary: without such an assumption, the estimate \eqref{eq:int-_g-in-1} can fail and the ``trivial bound'' \eqref{eq:int-_g-in-4} can be sharp.  The failure becomes even worse in higher rank.
\end{remark}
\begin{remark}
  As mentioned earlier in Remark \ref{rmk:we-discuss-informal-cal-F-U-G-T}, the function $f$ relevant to our application (furnished by Theorem \ref{thm:main-local-results}, constructed in \S\ref{sec:constr-test-vect-proofs}) has the rough shape
  \begin{equation*}
    f(x,y) \approx T^{-1/4} 1_{x \asymp T^{1/2}} 1_{y \ll 1} e(T y / x),
  \end{equation*}
  or equivalently,
  \begin{equation*}
    f(r \cos \theta, r \sin \theta) \approx
    T^{-1/4} 1_{r \asymp T^{1/2}}
    1 _{\theta \ll T^{-1/2}} e(T \theta),
  \end{equation*}
  and is easily seen to satisfy all hypotheses of Theorem \ref{thm:sub-gln:let-f-be} except condition \eqref{itm:sketch-sl2-assump-2}, concerning vanishing of line integrals.  The latter assumption is satisfied by the modified function obtained by applying $\mathcal{P}_G(s)$ as a Mellin multiplier, which is the relevant function for defining $\Psi$ as a series (see Example \ref{exa:example-sl2-completed-eisenstein-series} for details).
\end{remark}
\begin{remark}
  Many of the assumptions in Theorem \ref{thm:sub-gln:let-f-be} could be relaxed.  For instance, instead of assuming that $\mathcal{F} f = f$, it would suffice to assume that $\mathcal{F} f$ vanishes at the origin, and less radial uniformity than \eqref{eq:int-_v-in-7} is required.  Since the stronger assumptions are valid in the setting of Theorem \ref{thm:growth-eisenstein-nonstandard} and streamline our discussion, we keep them here.
\end{remark}

The proof of Theorem \ref{thm:sub-gln:let-f-be} involves the Fourier expansion of $\Psi$.  Set
\begin{equation*}
  n(x) := \begin{pmatrix}
    1 & x \\
    0 & 1
  \end{pmatrix}.
\end{equation*}
The constant term in that Fourier expansion is given by
\begin{equation*}
  \Psi_N(g) = \int_{x \in \mathbb{R} / \mathbb{Z} } \Psi(n(x) g) \, d x.
\end{equation*}
By some elementary calculations and the Poisson summation formula, we verify that
\begin{equation*}
  \Psi_N(g) = \sum_{t \in \mathbb{Z} } (f + \mathcal{F} f ) ((0,t) g).
\end{equation*}
Since $\mathcal{F}f = f$, we may simplify $f + \mathcal{F} f$ to $2 f$.

The non-constant terms are indexed by nonzero integers $\ell$:
\begin{equation*}
  W(\ell,g) := \int_{x \in \mathbb{R} / \mathbb{Z} } \Psi(n(x) g) e(-\ell x) \, d x.
\end{equation*}
Their evaluation is a simple calculation that we leave to the reader:
\begin{equation}\label{eq:well-g-=}
  W(\ell, g)
  = \sum_{c | \ell} \int_{x \in \mathbb{R} } f((c, x) g) e( \tfrac{- \ell}{c} x) \, d x.
\end{equation}
Set $R := \diag(t,1/t)$.  After enlarging $\Omega$ if necessary, it is enough to bound
\begin{equation*}
  \int_{x \in  \mathbb{R} / \mathbb{Z} }
  \int_{g \in \Omega }
  |\Psi|^2(n(x) R g) \, d g \, d x.
\end{equation*}
By Parseval and the Fourier calculations given above, we may expand this double integral as a sum of two terms:
\begin{equation*}
  I_0 := 4 \int_{g \in \Omega } \left\lvert \sum_{t \in \mathbb{Z} } f((0,t) R g) \right\rvert^2 \, d g.
\end{equation*}
and
\begin{equation*}
  I_1 := \sum_{0 \neq \ell \in \mathbb{Z}}
  \int_{g \in \Omega }
  |W(\ell, R g)|^2 \, d g.
\end{equation*}

The term $I_1$ may be treated much like in the cuspidal case, following a spectral expansion of $f$.  We focus here on the treatment of $I_0$.  Define $f^\sharp : \mathbb{R}^2 \rightarrow \mathbb{C}$ by the formula
\begin{equation*}
  f^\sharp(v)  := \sum_{t \in \mathbb{Z} } f(t v).
\end{equation*}
Recall that $f(0) = 0$, so that $f^\sharp (0)$ makes sense and vanishes.  Then
\begin{equation*}
  I_0 = 4 \int_{g \in \Omega } | f^\sharp(e_2 R g)|^2 \, d g.
\end{equation*}

From the bounds \eqref{eq:int-_v-in-7}, the vanishing condition \eqref{eq:int-_t-in-1} and some Mellin analysis, one can check that, e.g., for all fixed $u \in \mathfrak{U}(\SL_2(\mathbb{R}))$ of degree at most $m$,
\begin{equation}\label{eq:int-_theta-in}
  \int_{\theta \in \mathbb{R} / 2 \pi \mathbb{Z} }
  \left\lvert R(u) f^\sharp  (r \cos \theta, r \sin \theta ) \right\rvert^2
  \, d \theta
  \ll r^{-2} T^{2 m + o(1)} \max \left( \frac{r}{T^{1/2}}, \frac{T^{1/2}}{r} \right)^{-\infty}.
\end{equation}
Indeed, focusing on the case $u=1$ for notational simplicity, to obtain \eqref{eq:int-_theta-in}, it suffices by basic Paley--Wiener theory to verify that the Mellin transform
\begin{equation*}
  \int_{r =0}^{\infty}
  f^\sharp (r \cos \theta, r \sin \theta )
  r^s \, \frac{d r}{|r|}
\end{equation*}
is majorized, in an $L^2$-sense with respect to $\theta$, by $T^{s/2}(1 + |s|)^{-\infty}$ whenever $\Re(s) \ll 1$.  In view of the definition of $f^\sharp$, that Mellin transform evaluates to
\begin{equation*}
  \zeta(s)
  \int_{r =0}^{\infty}
  f  (r \cos \theta, r \sin \theta )
  r^s \, \frac{d r}{|r|}.
\end{equation*}
The vanishing line integral condition on $f$ implies that the above remains holomorphic at the unique pole of $\zeta$.  Since $\zeta$ grows at most polynomially, we obtain the required bounds for $f^\sharp $ from the assumed bounds for $f$.

Using the Sobolev lemma for $K$, we may deduce from \eqref{eq:int-_theta-in} that
\begin{equation*}
  \sup_{\theta \in \mathbb{R} / 2 \pi \mathbb{Z} }
  \left\lvert f^\sharp  (r \cos \theta, r \sin \theta ) \right\rvert
  \ll T^{m + C_0} \max \left( \frac{r}{T^{1/2}}, \frac{T^{1/2}}{r} \right)^{-\infty}
\end{equation*}
for some fixed $C_0$.  From this, we see that the integral $I_0$ is always negligibly small.  Indeed, we have $\left\lvert e_2 R \right\rvert = t^{-1} \leq 1$, whereas $f^\sharp $ is concentrated on elements of size $\asymp T^{1/2}$.

The above analysis of $I_0$ was perhaps a bit underwhelming, but it illustrates a couple points that will recur in the general case.
\begin{enumerate}[(i)]
\item The first is that it matters where $f$ is concentrated.  In this example, the concentration condition implies that $I_0$ is negligible.  What typically happens (for the various ``degenerate'' terms intermediary between the constant term contribution and the ``Whittaker term'') is that the concentration properties of $f$ yield a satisfactory estimate in certain cases (depending upon $t$, $T$ and other variables), while in other cases, other methods succeed.
\item The second is the passage from $f$ to $f^\sharp$.  While it turns out to be simpler to estimate the Whittaker contribution using the series representation of $\Psi$ involving $f$ (the formula \eqref{eq:well-g-=} being much nicer than the expression for $W$ in terms of $f^\sharp$), we could have simplified our treatment of the constant term contribution by working with the series representation
  \begin{equation*}
    \Psi(g) =
    \sideset{}{^*}\sum_{v \in \mathbb{Z}^2} f^\sharp (v g),
  \end{equation*}
  where $\sum^\ast$ denotes a sum over primitive elements.  In general, for each of the degenerate terms that we encounter, there is a specific series representation of $\Psi$ best suited for its estimation (\S\ref{sec:seri-repr-via}).
\end{enumerate}

\section{Rankin--Selberg unfolding}\label{sec:rank-selb-unfold}

\subsection{Notation}\label{sec:notation-rs-unfolding}
For each natural number $n$, let $V_n$ denote the punctured affine $n$-plane, omitting the origin.  For a ring $R$, $V_n(R)$ consists of all elements of the free module $R^n$ that generate a direct summand.  For a field or local ring $R$, we have $V_n(R) = \{(x_1,\dotsc,x_n) \in R^n : x_j \in R^\times \text{ for some } j \in \{1,\dotsc,n\}\}$.  We regard points of $V_n$ as row vectors, with $\GL_n$ acting on the right by matrix multiplication.  Let $e_1,\dotsc,e_n \in V_n(\mathbb{Z})$ denote the standard basis for $\mathbb{Z}^n$, e.g., $e_n = (0,\dotsc,0,1)$.  We denote by $P_n \leq \GL_n$ the stabilizer of $e_n$, i.e., the ``mirabolic'' subgroup consisting of all matrices with bottom row $(0,\dotsc,0,1)$.  The map $g \mapsto e_n g$ identifies $V_n$ with $P_n \backslash \GL_n$.  We let $U_n \leq P_n$ denote the unipotent radical and $\GL_{n-1} \leq P_n$ the standard upper-left $(n-1) \times (n-1)$ Levi subgroup, so that $P_n = \GL_{n-1} \rtimes U_n$.  We denote by $N_n \leq \GL_n$ the standard upper-triangular maximal unipotent subgroup.  Then $N_n = N_{n-1} \rtimes U_n \leq P_n$.

We henceforth fix one such $n$.  For each $m < n$, we regard $\GL_m$, hence also each of its subgroups defined above, as included in $\GL_n$ in the usual way, as the upper-left $m \times m$ block.

Let $F$ be a number field with adele ring $\mathbb{A}$.  Let $\psi : \mathbb{A}/F \rightarrow \U(1)$ be a nontrivial unitary character.  We denote also simply by $\psi$ the corresponding standard nondegenerate character of $N_n(\mathbb{A})$, given by $\psi(u) = \psi(\sum_{j=1}^{n-1} u_{j,j+1})$.  It induces a nontrivial character of the abelian group $[U_m]$ for each $m \leq n$.

\subsection{Mirabolic Eisenstein series}
For each $\phi \in C_c^\infty(V_n(\mathbb{A}))$, we define $\theta_\phi : [\GL_n] \rightarrow \mathbb{C}$ by the formula
\begin{equation*}
  \theta_\phi(g) := \sum_{0 \neq v \in V_n(F)} \phi(v g).
\end{equation*}
Since $\GL_n(F)$ acts transitively on $V_n(F)-\{0\}$ with stabilizer $P_n(F)$, we may equivalently write
\begin{equation}\label{eq:theta-phi-via-p_n}
  \theta_\phi(g) = \sum_{\gamma \in P_n(F)\backslash \GL_n(F)} \phi(e_n \gamma g).
\end{equation}
\begin{example}
  Suppose that $F = \mathbb{Q}$ and that $\phi = \otimes \phi_\mathfrak{p}$, where $\phi_p$ is the characteristic function of $V_n(\mathbb{Z}_p)$ and $\phi_\infty \in C_c^\infty(V_n(\mathbb{R}))$.  Then for $g \in \GL_n(\mathbb{Z}) \backslash \GL_n(\mathbb{R}) \hookrightarrow [\GL_n]$, we have
  \begin{equation*}
    \theta_\phi(g) =
    \sum_{
      \substack{
        v = (v_1,\dotsc,v_n) \in \mathbb{Z}^n :
        \\
        \gcd(v_1,\dotsc,v_n) = 1
      }
    } \phi_\infty(v g).
  \end{equation*}
  This special case is in fact the only case relevant for the applications of this paper, but it seemed clearer to us to develop the preliminaries in their natural generality.
\end{example}
It is easy to see, by elementary lattice point counting arguments, that any such function $\theta_\phi$ lies in the space $\mathcal{T}([G])$ of automorphic functions of uniform moderate growth.

\subsection{Standard Rankin--Selberg unfolding}\label{sec:stand-rank-selb}
Let $\varphi \in \mathcal{S}([\GL_n])$.  We aim to study, for $\phi \in C_c^\infty(V_n(\mathbb{A}))$, the integrals
\begin{equation*}
  \int_{[\GL_n]} |\varphi|^2 \theta_\phi.
\end{equation*}
Note that since $\varphi$ is of rapid decay and $\theta_\phi$ of moderate growth, such integrals converge absolutely.

Suppose for the moment that $\varphi$ is cuspidal.  Then $\varphi$ admits a rapidly-convergent Whittaker expansion \cite[Thm 5.9]{MR348047}
\begin{equation}\label{eq:varphig-=-sum}
  \varphi(g) = \sum_{\gamma \in N_n(F) \backslash P_n(F) } W_{\varphi}(\gamma g),
  \quad
  W_{\varphi}(g) := \int_{[N_n]}
  \varphi(u g) \psi^{-1}(u) \, d u.
\end{equation}
Using \eqref{eq:theta-phi-via-p_n}, we obtain the following unfolding:
\begin{align*}
  \int_{[\GL_n]} |\varphi|^2 \theta_\phi  &=
                            \int_{P_n(F) \backslash \GL_n(\mathbb{A})}|\varphi|^2(g) \phi(e_n g) \, d g \\
                          &=
                            \int_{N_n(F) \backslash \GL_n(\mathbb{A})} W_{\varphi}(g) \overline{\varphi(g)} \phi(e_n g) \, d g \\
                          &=
                            \int_{g \in N_n(\mathbb{A}) \backslash \GL_n(\mathbb{A})}
                            \int_{u \in N_n(F) \backslash N_n(\mathbb{A})}
                            W_{\varphi}(u g) \overline{\varphi(u g)} \phi(e_n u g) \, d g \\
                          &=
                            \int_{N_n(\mathbb{A}) \backslash \GL_n(\mathbb{A})} |W_{\varphi}|^2(g) \phi(e_n g) \, d g.
\end{align*}
Calculations like this are the basis of the Rankin--Selberg method on $\GL_n \times \GL_n$ \cite{MR701565, MR618323}.

Our first aim in this section (Proposition \ref{prop:sub-gln:let-varphi-in-preilm-unfold}) is to carry out similar calculations, but without the cuspidality assumption on $\varphi$.  The resulting formula is similar to the above, but with some ``correction'' terms coming from other parts of the Fourier--Whittaker expansion of $\varphi$.  For the specific identity needed here, we could not quickly locate a direct reference.  Similar calculations, including a more general form of \eqref{eq:varphig-=-sum}, may be found in \cite[Prop 4.2]{MR3334892}.  Rankin--Selberg calculations in non-cuspidal settings go back to Zagier's work \cite{MR656029}.

We then apply the resulting formula to the case that $\varphi$ is an incomplete Eisenstein series defined with respect to a minimal parabolic.  Doing so involves several applications of standard formulas for constant terms and inner products of Eisenstein series.

\subsection{Partial Whittaker transforms}\label{sec:part-whitt-transf}
We require some further notation.

Let $j \in \{0, \dotsc, n-1\}$.

We define the subgroups $P_n^{(j)} \leq \GL_n$ by $P_n^{(j)} := P_{n-j} N_n$.  We denote by $U_n^{(j)}$ the unipotent radical of $P_n^{(j)}$.  Then
\begin{equation}\label{eq:cukcadro62}
  P_n^{(j)} = \GL_{n-j-1} \ltimes U_n^{(j)}
\end{equation}
and $U_n^{(j)} = U_{n-j} \dotsb U_{n-1} U_n$.  For example, $P_n^{(0)} = P_n$ and $U_n^{(0)} = U_n$, while $P_n^{(n-1)} = U_n^{(n-1)} = U_n^{(n-2)} = N_n$.

Generalizing the conventions of \S\ref{sec:haar-measures}, we equip the adelic points of any unipotent group $U$ with the Haar measure for which $[U] = U(F) \backslash U(\mathbb{A})$ has volume one.

We write $W^{(j)}_\varphi(g)$ for the ``$j$th partial Whittaker function,'' defined inductively as follows: $W_\varphi^{(0)} := \varphi$ and
\begin{equation*}
  W_{\varphi}^{(j+1)}(g) := \int_{u \in [U_{n-j}]} W_\varphi^{(j)}(u g) \psi^{-1}(u) \, d u.
\end{equation*}
Thus, for instance, $W_\varphi^{(n-1)} = W_\varphi$.   We note that $W_\varphi^{(j)}$ is left invariant under $P_n^{(j)}(F)$ and, for $j \geq 1$, transforms on the left under $U_n^{(j-1)}(\mathbb{A})$ via $\psi$.

We denote by $W^{(j)}_{\varphi,U_{n-j}}$ the constant term of $W^{(j)}_{\varphi}$ with respect to $U_{n-j}$:
\begin{equation*}
  W^{(j)}_{\varphi,U_{n-j}}(g) :=
  \int_{u \in [U_{n-j}]} W_\varphi^{(j)}(u g) \, d u.
\end{equation*}
For example, when $j = n-1$, the group $U_{n-j}$ is trivial and we have $W_{\varphi,U_{n-j}}^{(j)} = W_{\varphi}$, whereas for $j=0$, we have $W_{\varphi,U_{n}}^{(0)}(g) = \int_{u \in [U_n]} \varphi(u g) \, d u$ (which we  later denote by $\varphi_{U_n}(g)$).

We note that the restriction of $W_{\varphi,U_{n-j}}^{(j)}$ to $\GL_{n-j-1}(\mathbb{A})$ is left $\GL_{n-j-1}(F)$-invariant, hence defines a function on $[\GL_{n-j-1}]$.

Recall that we have chosen a Haar measure on each $\GL_m(\mathbb{A})$.  Using \eqref{eq:cukcadro62}, we obtain a left Haar measure $d p$ on each $P_n^{(j)}(\mathbb{A})$.  The corresponding right Haar measure is given by $\lvert \det p \rvert^{j + 1} \, d p$.  Explicitly, writing $m := n - j - 1$, for an integrable function $\phi$ on $[P_n^{(j)}]$, we have
\begin{align}
  \int_{p \in [P_n^{(j)}]} \phi(p)\, d p
  &=
  \int_{h \in [\GL_m]}
  \int_{u \in [U_n^{(j)}]}
    \phi(h u)\, d u\,
    d h \nonumber
  \\
  &= \label{eq:cujy0e5m3g}
    \int_{h \in [\GL_m]}
    \int_{u \in [U_n^{(j)}]}
    \phi(u h)\, d u\,
    \frac{d h}{|\det h|^{j + 1}}.
\end{align}

For each $j \in \{0,\dotsc,n-1\}$, let $\mathcal{C}_j$ denote the space of smooth functions $f : \GL_n(\mathbb{A}) \rightarrow \mathbb{C}$ satisfying
\begin{equation*}
  f(p g)=|\det p|^{j+1} f(g)
  \quad
  \text{ for } p \in P_n^{(j)}(\mathbb{A})
\end{equation*}
and of compact support modulo $P_n^{(j)}(\mathbb{A})$.  The construction of quotient integrals gives a $\GL_n(\mathbb{A})$-invariant functional
\begin{equation*}
  \mathcal{C}_j \rightarrow \mathbb{C},
  \qquad
  f \mapsto \int_{P_n^{(j)}(\mathbb{A}) \backslash \GL_n(\mathbb{A})}
  f(g)\, d g,
\end{equation*}
characterized by the identity
\begin{equation}\label{eq:culjy5hdwt}
  \int_{P_n^{(j)}(\mathbb{A}) \backslash \GL_n(\mathbb{A})}
  \left(
    \int_{P_n^{(j)}(\mathbb{A})}
    \phi(p g) \, d p
  \right)\, d g
  =
  \int_{\GL_n(\mathbb{A})} \phi(g)\, d g.
\end{equation}

For each $j \in \{0,\dotsc,n-2\}$, applying the above construction with ambient group $P_n^{(j)}(\mathbb{A})$ and subgroup $P_n^{(j+1)}(\mathbb{A})$ gives a quotient functional for $P_n^{(j+1)}(\mathbb{A}) \backslash P_n^{(j)}(\mathbb{A})$, defined on functions satisfying $f(p g) = \lvert \det p \rvert^{j + 2} f(g)$ for $p \in P_n^{(j + 1)}(\mathbb{A})$.

\begin{lemma}\label{lem:sub-gln:quotient-composition}
  Let $0 \leq j \leq n - 2$, and let $\Phi : \GL_n(\mathbb{A}) \rightarrow \mathbb{C}$ be a nonnegative measurable function satisfying
  \begin{equation*}
    \Phi(p g) = \lvert \det p \rvert^{j + 2} \Phi(g)
    \qquad
    \text{ for } p \in P_n^{(j + 1)}(\mathbb{A}).
  \end{equation*}
  Then
  \begin{equation}\label{eq:cumiu83k5d}
    \int_{P_n^{(j)}(\mathbb{A}) \backslash \GL_n(\mathbb{A})}
    \int_{P_{n - j - 1}(\mathbb{A}) \backslash \GL_{n - j - 1}(\mathbb{A})}
    \Phi(x g)
    \frac{d x \, d g}{\lvert \det x \rvert^{j + 1}}
    =
    \int_{P_n^{(j + 1)}(\mathbb{A}) \backslash \GL_n(\mathbb{A})}
    \Phi(g)\, d g.
  \end{equation}
\end{lemma}
\begin{proof}
  By the construction of the quotient integrals, it suffices to prove the identity for functions of the form
  \begin{equation*}
    \Phi(g)
    =
    \int_{P_n^{(j + 1)}(\mathbb{A})}
    f(p g)\, d p,
  \end{equation*}
  with $f$ nonnegative and measurable.  For such $\Phi$, the right-hand side of \eqref{eq:cumiu83k5d} is
  \begin{equation*}
    \int_{\GL_n(\mathbb{A})} f(g)\, d g.
  \end{equation*}
  It remains only to check that the iterated quotient integral on the left gives the same value.  Equivalently, by applying \eqref{eq:culjy5hdwt} first with $P_n^{(j)}$ and then with $P_n^{(j + 1)}$, it is enough to verify the corresponding identity inside $P_n^{(j)}(\mathbb{A})$:
  \begin{equation}\label{eq:cu4my82zab}
    \int_{P_{n - j - 1}(\mathbb{A}) \backslash \GL_{n - j - 1}(\mathbb{A})}
    \left(
      \int_{P_n^{(j + 1)}(\mathbb{A})}
      f(p x)\, d p
    \right)
    \frac{d x}{\lvert \det x \rvert^{j + 1}}
    =
    \int_{P_n^{(j)}(\mathbb{A})} f(p)\, d p.
  \end{equation}
  The measure on $P_n^{(j)}(\mathbb{A})$ is the product measure coming from
  \[
    P_n^{(j)} = \GL_{n - j - 1} \rtimes U_n^{(j)},
  \]
  and the measure on $P_n^{(j + 1)}(\mathbb{A})$ is the corresponding product measure coming from
  \[
    P_n^{(j + 1)} = P_{n - j - 1} \rtimes U_n^{(j)}.
  \]
  Thus \eqref{eq:cu4my82zab} follows from \eqref{eq:cujy0e5m3g}, applied with $\GL_{n - j - 1}$ in place of the Levi factor and then unfolded along $P_{n - j - 1}(\mathbb{A}) \backslash \GL_{n - j - 1}(\mathbb{A})$:
  \begin{align*}
    \int_{P_n^{(j)}(\mathbb{A})} f(p)\, d p
    &=
      \int_{\GL_{n - j - 1}(\mathbb{A})}
      \int_{U_n^{(j)}(\mathbb{A})}
      f(u x)\, d u\,
      \frac{d x}{\lvert \det x \rvert^{j + 1}}
    \\
    &=
      \int_{P_{n - j - 1}(\mathbb{A}) \backslash \GL_{n - j - 1}(\mathbb{A})}
      \int_{P_{n - j - 1}(\mathbb{A})}
      \int_{U_n^{(j)}(\mathbb{A})}
      f(u p x)\, d u\,
      \frac{d p}{\lvert \det p \rvert^{j + 1}}
      \frac{d x}{\lvert \det x \rvert^{j + 1}}
    \\
    &=
      \int_{P_{n - j - 1}(\mathbb{A}) \backslash \GL_{n - j - 1}(\mathbb{A})}
      \int_{P_{n - j - 1}(\mathbb{A})}
      \int_{U_n^{(j)}(\mathbb{A})}
      f(p u x)\, d u\, d p\,
      \frac{d x}{\lvert \det x \rvert^{j + 1}}
    \\
    &=
      \int_{P_{n - j - 1}(\mathbb{A}) \backslash \GL_{n - j - 1}(\mathbb{A})}
      \int_{P_n^{(j + 1)}(\mathbb{A})}
      f(p x)\, d p\,
      \frac{d x}{\lvert \det x \rvert^{j + 1}},
  \end{align*}
  where in the third step, we use that conjugation by $p$ on $U_n^{(j)}(\mathbb{A})$ has Jacobian $\lvert \det p \rvert^{j + 1}$.
\end{proof}

\subsection{Preliminary unfolding}

\begin{proposition}\label{prop:sub-gln:let-varphi-in-preilm-unfold}
  Let $\varphi \in \mathcal{S}([\GL_n])$ and $\phi \in C_c^\infty(V_n(\mathbb{A}))$.  We have
  \begin{equation}\label{eq:int-_gl_n-varphi2}
    \int_{[\GL_n]} |\varphi|^2 \theta_\phi  =
    \sum_{j=0}^{n-1} I_j,
  \end{equation}
  where
  \begin{equation}\label{eq:cullblm1fm}
    I_j := \int_{P_n^{(j)}(\mathbb{A}) \backslash \GL_n(\mathbb{A})} \phi(e_n g)
    \int_{h \in [\GL_{n-j-1}]}
    \left\lvert
      W_{\varphi,U_{n-j}}^{(j)} (h g)
    \right\rvert^2
    \, \frac{d h}{|\det h|^{j+1}}
    \, d g.
  \end{equation}
  Each of the above integrals converges absolutely.
\end{proposition}
\begin{remark}
  It is not hard to see that if $\varphi$ is cuspidal, then for each $j \neq n-1$, we have $W_{\varphi,U_{n-j}}^{(j)} = 0$, hence $I_j = 0$.  On the other hand,
  \begin{equation*}
    I_{n-1} = \int_{N_n(\mathbb{A}) \backslash \GL_n(\mathbb{A})} |W_{\varphi}|^2(g) \phi(e_n g) \, d g,
  \end{equation*}
  as in \S\ref{sec:stand-rank-selb}.
\end{remark}

The proof requires a preliminary lemma.  For each continuous function $\varphi$ on $[\GL_n]$, we denote by $\varphi_{U_n}$ its constant term with respect to $U_n$, i.e., $\varphi_{U_n}(g) := \int_{u \in [U_n]} \varphi(u g) \, d u$.
\begin{lemma}\label{lem:sub-gln:let-varphi-be-parseval-U-n}
  Let $\varphi$ be a continuous function on $[\GL_n]$.  Then
  \begin{equation*}
    \int_{u \in [U_n]} |\varphi|^2( u g) \, d u
    =
    |\varphi_{U_n}|^2(g)
    +
    \sum_{\gamma \in P_{n-1}(F) \backslash \GL_{n-1}(F)}
    \left\lvert
      W_{\varphi}^{(1)}(\gamma g)
    \right\rvert^2.
  \end{equation*}
\end{lemma}
\begin{proof}
  This computation is essentially contained in the proof of \eqref{eq:varphig-=-sum}.  We sketch it for the sake of convenience.  The quotient $[U_n]$ is a compact abelian group.  By considering the first $n-1$ rows of its final column, we may identify it with the space of column vectors $u = (u_1,\dotsc,u_{n-1})^\transpose$ with each $u_j \in \mathbb{A}/F$.  The Pontryagin dual $[U_n]^\wedge$ identifies with the space of row vectors $\xi = (\xi_1,\dotsc,\xi_{n-1})$, with the corresponding character $\psi_\xi$ of $[U_n]$ given by $\psi_\xi(u) = \psi(\xi u) := \psi(\sum_{j=1}^{n-1} \xi_j u_j)$.  In any event, Parseval's relation gives
  \begin{equation*}
    \int_{u \in [U_n]} |\varphi|^2(u g) \, d u
    = \sum_{\xi \in [U_n]^\wedge }
    |\varphi_\xi(g)|^2,
    \quad
    \varphi_\xi(g)
    :=
    \int_{u \in [U_n]}
    \varphi(u g) \psi_\xi^{-1}(u) \, d u.
  \end{equation*}
  The group $\GL_{n-1}(F)$ acts via conjugation on $[U_n]$, hence by duality on $[U_n]^\wedge$.  For $\gamma \in \GL_{n-1}(F)$, we have
  \begin{equation*}
    \varphi_{\xi}(\gamma g) = \varphi_{\xi \gamma}(g).
  \end{equation*}
  There are two orbits for the latter action: $\{0\}$ and $[U_n]^\wedge - \{0\}$.  The contribution from $\xi = 0$ is $\varphi_0 = \varphi_{U_n}$.  The representative $\xi = e_{n-1}^\transpose$ for $[U_n]^\wedge - \{0\}$ has stabilizer $P_{n-1}(F)$, and we have $\varphi_{e_{n-1}^\transpose} = W_{\varphi}^{(1)}$.  The required formula folllows.
\end{proof}

\begin{proof}[Proof of Proposition \ref{prop:sub-gln:let-varphi-in-preilm-unfold}]
  We first address the issue of convergence.  The proof given below shows that the identity \eqref{eq:int-_gl_n-varphi2} holds when $\phi$ is a nonnegative-valued continuous function of rapid decay.  For such $\phi$, the incomplete Eisenstein series $\theta_\phi$ remains of moderate growth, so the integral on the LHS of \eqref{eq:int-_gl_n-varphi2} converges.  Since each integrand on the RHS is nonnegative, it follows that each integral on the RHS converges.  We can apply the argument just indicated to the absolute value $|\phi|$ of $\phi$, giving the required absolute convergence in the stated generality of the proposition and also in each step of the proof that follows.

  We use the convention that $P_1$ and $\GL_0$ are trivial.  For $j \in \{0,\dotsc,n-1\}$, set
  \begin{equation*}
    J_j(g)
    :=
    \int_{p \in [P_{n-j}]}
    \left\lvert W_\varphi^{(j)}(p g) \right\rvert^2
    \frac{d p}{\lvert \det p \rvert^j}
  \end{equation*}
  and
  \begin{equation*}
    R_j
    :=
    \int_{P_n^{(j)}(\mathbb{A}) \backslash \GL_n(\mathbb{A})}
    \phi(e_n g) J_j(g)\, d g.
  \end{equation*}

  Using \eqref{eq:theta-phi-via-p_n} for $\theta_\phi$ and \eqref{eq:culjy5hdwt} with $j=0$, we unfold the left-hand side of \eqref{eq:int-_gl_n-varphi2} as $R_0$.  We see also that $R_{n-1}=I_{n-1}$.  To complete the proof, it suffices to show that for each $0 \leq j \leq n - 2$, we have
  \begin{equation*}
    R_j = I_j + R_{j+1}.
  \end{equation*}
  We begin by rewriting $J_j(g)$ using \eqref{eq:cujy0e5m3g} applied to $P_{n-j} = P_{n - j}^{(0)}$:
  \begin{equation*}
    J_j(g)
    =
    \int_{h \in [\GL_{n-j-1}]}
    \int_{u \in [U_{n-j}]}
    \left\lvert W_\varphi^{(j)}(u h g) \right\rvert^2\, d u\,
    \frac{d h}{\lvert \det h \rvert^{j+1}}.
  \end{equation*}
  Next, the same Parseval calculation as in Lemma \ref{lem:sub-gln:let-varphi-be-parseval-U-n}, applied to $W_\varphi^{(j)}$ and $U_{n-j}$, gives
  \begin{equation*}
    \int_{u \in [U_{n-j}]}
    \left\lvert W_\varphi^{(j)}(u g) \right\rvert^2\, d u
    =
    \left\lvert W_{\varphi,U_{n-j}}^{(j)}(g) \right\rvert^2
    +
    \sum_{\gamma \in P_{n-j-1}(F) \backslash \GL_{n-j-1}(F)}
    \left\lvert W_\varphi^{(j+1)}(\gamma g) \right\rvert^2.
  \end{equation*}
  It follows that $J_j(g) = A_j(g) + B_j(g)$, where
  \begin{equation*}
    A_j(g)
    :=
    \int_{h \in [\GL_{n-j-1}]}
    \left\lvert W_{\varphi,U_{n-j}}^{(j)}(h g) \right\rvert^2
    \frac{d h}{\lvert \det h \rvert^{j+1}},
  \end{equation*}
  \begin{equation*}
    B_j(g)
    :=
    \int_{h \in [\GL_{n-j-1}]}
    \sum_{\gamma \in P_{n-j-1}(F) \backslash \GL_{n-j-1}(F)}
    \left\lvert W_\varphi^{(j+1)}(\gamma h g) \right\rvert^2
    \frac{d h}{\lvert \det h \rvert^{j+1}}.
  \end{equation*}
  We now separately evaluate the contributions of $A_j$ and $B_j$ to $R_j$.  The contribution of $A_j$ is $I_j$, as follows immediately from the definition \eqref{eq:cullblm1fm} of the latter.  It remains to show that the contribution of $B_j$ is $R_{j + 1}$.  To that end, we first unfold the sum in the definition of $B_j$:
  \begin{align*}
    B_j(g)
    &=
      \int_{P_{n-j-1}(F) \backslash \GL_{n-j-1}(\mathbb{A})}
      \left\lvert W_{\varphi}^{(j+1)}(x g) \right\rvert^2
      \frac{d x}{\lvert \det x \rvert^{j+1}} \\
    &=
      \int_{x \in P_{n-j-1}(\mathbb{A}) \backslash \GL_{n-j-1}(\mathbb{A})}
      \int_{p \in[P_{n - j - 1}]}
      \left\lvert W_{\varphi}^{(j+1)}(p x g) \right\rvert^2
      \frac{d p}{\lvert \det p \rvert^{j+1}}
      \frac{d x}{\lvert \det x \rvert^{j+1}}
    \\
    &=
      \int_{x \in P_{n-j-1}(\mathbb{A}) \backslash \GL_{n-j-1}(\mathbb{A})}
      J_{j + 1}(x g)
      \frac{d x}{\lvert \det x \rvert^{j+1}}.
  \end{align*}
  We pause to observe that
  \begin{equation*}
    J_{j+1}(p g)  =  \lvert \det p \rvert^{j + 2} J_{j+1}(g) \text{ for } p \in P_n^{(j + 1)}(\mathbb{A}).
  \end{equation*}
  In verifying this, we may use that $P_n^{(j + 1)} = P_{n - j - 1} \rtimes U_n^{(j)}$ to reduce to the following cases:
  \begin{itemize}
  \item $p \in P_{n - j - 1}(\mathbb{A})$, in which case the claim follows from the definition of $J_{j + 1}$ and the relation between left and right Haar measures on $P_{n - j - 1}(\mathbb{A})$;
  \item $p \in U_n^{(j)}(\mathbb{A})$, in which case the claim follows from the fact that $W_\varphi^{(j + 1)}$ transforms on the left under $U_n^{(j)}(\mathbb{A})$ via $\psi$, so $J_{j + 1}$ is left invariant under $U_n^{(j)}(\mathbb{A})$.
  \end{itemize}
  Since $e_n p = e_n$ for $p \in P_n^{(j + 1)}(\mathbb{A})$, we also have $\phi(e_n p g) = \phi(e_n g)$ for $p \in P_n^{(j + 1)}(\mathbb{A})$.  Therefore, the function
  \[
    \Phi(g) := \phi(e_n g) J_{j+1}(g)
  \]
  satisfies the covariance required in Lemma \ref{lem:sub-gln:quotient-composition}.  By the conclusion of that lemma, we see that the contribution of $B_j$ to $R_j$ is precisely
  \[
    \int_{P_n^{(j + 1)}(\mathbb{A}) \backslash \GL_n(\mathbb{A})}
    \phi(e_n g) J_{j+1}(g)\, d g
    =
    R_{j+1}.
  \]
  Thus $R_j = I_j + R_{j+1}$, as required.
\end{proof}

\subsection{Reformulation in terms of constant terms}\label{sec:reform-terms-const-1}
The definitions of \S\ref{sec:part-whitt-transf} were convenient for writing the proof of Proposition \ref{prop:sub-gln:let-varphi-in-preilm-unfold}, but in our applications, a different description of the quantities $W_{\varphi,U_{n-j}}^{(j)}$ is convenient.

Set $G := \GL_n$.  We return to our standard practice of denoting by $B = N A$ the standard upper-triangular Borel subgroup, with Levi $A$ given by the diagonal subgroup.  Thus $N = N_n$ in the notation above.

For each $m \in \{1, \dotsc, n-1\}$, let $P_{(m,n-m)}$ denote the standard upper-triangular parabolic subgroup of $G$ with Levi component $M_{(m,n-m)}$ given by the block-diagonal product $\GL_m \times \GL_{n-m}$.  Let $U_{(m,n-m)}$ denote its unipotent radical.  By a \emph{standard maximal parabolic subgroup} $P = M U_P$ of $G := \GL_n$, we mean that $(P,M,U_P) = (P_{(m,n-m)}, M_{(m,n-m)}, U_{(m,n-m)})$ for some $m$.  In that case, we may factor $M = M' \times M''$, where $M' = \GL_m$ and $M'' \cong \GL_{n-m}$, with the latter included now as the lower-right block.  Set $N'' := N_n \cap M''$, and define $j$ by writing $m = n -j - 1$.  Then
\begin{equation*}
  W_{\varphi,U_{n-j}}^{(j)}(g)
  = \int_{u \in [N'']}
  \varphi_{P}(u g) \psi^{-1}(u) \, d u,
\end{equation*}
where $\varphi_P$ denotes the constant term $\varphi_{P}(g) := \int_{u \in [U_P]} \varphi(u g) \, d u$ as defined in \S\ref{sec:constant-terms}.  Proposition \ref{prop:sub-gln:let-varphi-in-preilm-unfold} thus admits the following reformulation: the difference $\int_{[\GL_n]} |\varphi|^2 \theta_\phi - \int_{N(\mathbb{A}) \backslash G(\mathbb{A})} \phi |W_\varphi|^2$ coincides with the following sum, taken over all standard maximal parabolic subgroups $P$ of $G$:
\begin{equation}\label{eq:sum_-p-=}
  \sum_{
    P = M U_P
  }
  \int_{g \in (M' N'' U_P \backslash G)(\mathbb{A})}
  \phi (g)
  \int_{h \in [M']}
  \left\lvert
    \int_{u \in [N'']}
    \varphi_P(h u g) \psi^{-1}(u) \, d u
  \right\rvert^2 \, \frac{d h}{\delta_{U_P}(h)} \, d g.
\end{equation}

\subsection{Partial Whittaker functions of pseudo Eisenstein series}
We specialize henceforth to $G = \GL_n$, with accompanying notation as in \S\ref{sec:reform-terms-const-1}.  Let $W$ denote the Weyl group for $(G,A)$.

\begin{lemma}\label{lem:sub-gln:let-phi-in-eval-W-eis}
  Let $\Phi \in \mathcal{S}([G]_B)$, $\varphi := \Eis_B^G(\Phi)$ (as defined in \S\ref{sec:cuhjp7akc2}).  Then $W_{\varphi} = W_{\Phi}$, where, with $w_G \in W \hookrightarrow G(F)$ the long Weyl element,
  \begin{equation*}
    W_\Phi(g) := \int_{u \in N(\mathbb{A})} \Phi(w_G^{-1} u g) \psi^{-1}(u) \, d u.
  \end{equation*}
\end{lemma}
\begin{proof}
  This is a case of the standard unfolding calculation recalled earlier in \S\ref{sec:whittaker-functions}.
\end{proof}

\begin{lemma}\label{lem:sub-gln:let-p-=}
  Let $P = M {U_P}$ be a standard maximal parabolic subgroup of $G$, with accompanying notation as above.  Let $\Phi \in \mathcal{S}([M]_{B_M})$, $\varphi := \Eis_{B_M}^{M}(\Phi)$.  Assume that for each $w' \in W_{M'}$, we have $M_{w'} \Phi \in \mathcal{S}([M]_{B_M})$.  Then, writing $w_{M''} \in W_{M''} \hookrightarrow M''(F)$ for the long Weyl element,
  \begin{align*}
    &\int_{h \in [M']}
      \left\lvert \int_{u \in [N'']}
      \varphi(h u) \psi^{-1}(u) \, d u
      \right\rvert^2
      \, d h \\
    &\quad =
      \frac{1}{|W_{M'}|}
      \int_{h \in [M']_{B_{M'}}}
      \left\lvert
      \sum_{w \in W_{M'}}
      \int_{u \in N''(\mathbb{A})}
      M_w \Phi (h w_{M''}^{-1} u)
      \psi^{-1}(u) \, d u
      \right\rvert^2
      \, d h.
  \end{align*}
\end{lemma}
\begin{proof}
  For $w \in W_{M'}$, we may define (\S\ref{sec:basic-definition-intertwining-ops}) an intertwining operator $M_w$ on both $\mathcal{S}([M']_{B_{M'}})$ and -- via the inclusion $W_{M'} \hookrightarrow W_M$ -- on $\mathcal{S}([M]_{B_{M}})$.  These operators are compatible in the sense that the diagram
  \begin{equation*}
    \begin{CD}
      \mathcal{S}([M]_{B_{M}}) @>M_w>>\mathcal{T}([M]_{B_M})\\
      @V\text{restriction}VV  @VV\text{restriction}V \\
      \mathcal{S}([M']_{B_{M'}}) @>>M_w> \mathcal{T}([M']_{B_{M'}})\\
    \end{CD}
  \end{equation*}
  commutes.

  Turning to the proof of the proposition, consider first the case of a pure tensor $\Phi = \Phi ' \otimes \Phi ''$, i.e., $\Phi(m' m'') = \Phi ' (m') \Phi ''(m'')$ for $(m', m'') \in M'(\mathbb{A}) \times M''(\mathbb{A})$, where
  \begin{itemize}
  \item  $\Phi ' \in \mathcal{S}([M']_{B_{M'}})$ satisfies $M_{w'} \Phi ' \in \mathcal{S}([M']_{B_{M'}})$ for all $w' \in W_{M'}$, and
  \item  $\Phi '' \in \mathcal{S}([M'']_{B_{M''}})$.
  \end{itemize}
  Then $\varphi = \varphi ' \otimes \varphi ''$ with $\varphi ' = \Eis_{B_{M'}}^{M'}(\Phi ')$, $\varphi '' = \Eis_{B_{M''}}^{M''}(\Phi '')$, so the claimed identity reads
  \begin{align*}
    &\int_{[M'] } |\varphi '|^2
      \cdot \left\lvert \int_{u \in [N'']} \varphi ''(u) \psi^{-1}(u) \, d u  \right\rvert^2 \\
    &\quad
      =
      \frac{1}{|W_{M'}|}
      \int_{[M']_{B_{M'}}} \left\lvert \sum_{w \in W_{M'}} M_w \Phi ' \right\rvert^2
      \cdot
      \left\lvert \int_{u \in N''(\mathbb{A} )}
      \Phi ''(w_{M''}^{-1} u) \psi^{-1}(u) \, d u
      \right\rvert^2
  \end{align*}
  and follows from Proposition \ref{lem:let-f_1-f_2} (resp.\ Lemma \ref{lem:sub-gln:let-phi-in-eval-W-eis}) applied to $M'$ (resp.\ $M''$).

  The general case may be deduced in a couple ways.
  \begin{itemize}
  \item We may reduce, as follows, to the case of pure tensors: both sides depend continuously upon $\Phi$, and as explained in \cite[\S2.5.9]{MR4426741}, for each pair of compact open subgroups $J' < M'(\mathbb{A}^\infty)$ and $J'' < M''(\mathbb{A}^\infty)$ with product $J := J' \times J'' < M(\mathbb{A}^\infty)$, we may identify $\mathcal{S}([M]_{B_M})^J$ with the completed tensor product of nuclear Fr{\'e}chet spaces $\mathcal{S}([M']_{B_{M'}})^{J'} \hat{\otimes} \mathcal{S}([M'']_{B_{M''}})^{J''}$.
  \item We can apply the proof given in the pure tensor case to the general case.  Indeed, by Lemma \ref{lem:sub-gln:let-phi-in-eval-W-eis}, we have
    \begin{equation}\label{eq:int-_u-in-1}
      \int_{u \in [N'']} \varphi(h u) \psi^{-1}(u) \, d u
      =
      \sum_{\gamma \in B_{M'}(F) \backslash M'(F)}
      \int_{u \in N''(\mathbb{A})}
      \Phi(\gamma h u)
      \psi^{-1}(u) \,  d u,
    \end{equation}
    where the RHS converges absolutely.  For $u_1, u_2 \in [N'']$, define $\Phi_1, \Phi_2 \in \mathcal{S}([M']_{B_{M'}})$ by $\Phi_j(m') := \Phi(u_j m') = \Phi(m' u_j)$.  Then for $j =1,2$ and $w \in W_{M'}$, our assumption on $\Phi_j$ (and the noted compatibility of the definitions of $M_w$) implies that $M_w \Phi_j \in \mathcal{S}([M']_{B_{M'}})$.  By Proposition \ref{lem:let-f_1-f_2}, it follows that
    \begin{align*}
      &\int_{h \in [M']}
        \sum_{\gamma_1, \gamma_2 \in B_{M'}(F) \backslash M'(F)}
        \Phi_1(\gamma_1 h u_1)
        \overline{\Phi_2}(\gamma_2 h u_2) \, d h \\
      &\quad =
        \frac{1}{|W_{M'}|}
        \sum_{w_1, w_2 \in W_{M'}}
        \int_{h \in [M']_{B_{M'}}}
        M_{w_1} \Phi_1(h)
        \overline{M_{w_2} \Phi_2(h)} \, d h.
    \end{align*}
    We conclude by inserting the definitions of $\Phi_1, \Phi_2$ and integrating over $u_1, u_2 \in N''(\mathbb{A})$ as in \eqref{eq:int-_u-in-1}.  We emphasize that at each stage of the proof, every iterated sum/integral converges absolutely, thanks in each case to the locally uniform convergence of pseudo Eisenstein series on $M$.
  \end{itemize}
\end{proof}

We note that Lemma \ref{lem:sub-gln:let-p-=} applies also to the apparently more general integrals obtained by multiplying the integrand by a positive-valued character evaluated at $h$, e.g.,
\begin{equation*}
  \int_{h \in [M']}
  \left\lvert \int_{u \in [N'']} \varphi(h u) \psi^{-1}(u) \,d u  \right\rvert^2
  \, \frac{d h}{\delta_{U_P}(h)}.
\end{equation*}
To obtain a formula for such integrals, we apply  Lemma \ref{lem:sub-gln:let-p-=} with $\Phi$ replaced by $\Phi \delta_{U_P}^{-1/2}$, noting that multiplication by $\delta_{U_P^{-1/2}}$ commutes with the intertwining operators $M_w$.

\subsection{Summary}\label{sec:summary-unification}
To unify the presentation, we extend the notation introduced in \S\ref{sec:reform-terms-const-1} from the case that $P$ is a standard maximal parabolic subgroup to the case $P = G$.  In that case, we set
\begin{equation*}
  M' := \{1\}, \quad M'' := G, \quad
  {U_P} = \{1\},
\end{equation*}
and take $N'$ and $N''$ as before.

For $P$ maximal parabolic or $P = G$ and $\Phi \in \mathcal{S}([G]_B)$, we define $W^P[\Phi] \in C^\infty(G(\mathbb{A}))$ by the formula \index{Jacquet integrals and distributions!$W^P[\Phi]$}
\begin{equation*}
  W^P[\Phi](g) := \int_{u \in N''(\mathbb{A}) } \Phi(w_{M''}^{-1} u g) \psi^{-1}(u) \, d u.
\end{equation*}
For example, $W^G[\Phi] := W_{\Phi}$ as defined in Lemma \ref{lem:sub-gln:let-phi-in-eval-W-eis}.

We recall the set $W(A,M)$, defined following \eqref{eq:eis_bgf_p-=-sum} in our preliminary discussion of constant terms of Eisenstein series.

\begin{proposition}\label{prop:sub-gln:let-phi-in-summary}
  Let $\phi \in C_c^\infty(V_n(\mathbb{A}))$ and $\Phi \in \mathcal{S}([G]_B)$.  Assume that $M_w \Phi \in \mathcal{S}([G]_B)$ for each $w \in W$.  Then
  \begin{align}
    \int_{[G]} \theta_\phi |\Eis(\Phi)|^2  \label{prop:=-sum_p-=}
    &=
      \sum_{P = M {U_P}}
      \int_{g \in (M' N'' {U_P} \backslash G)(\mathbb{A}) } \phi(e_n g) \\ \nonumber
    &\quad
      \frac{1}{|W_{M'}|}
      \int_{h \in [M']_{B_{M'}}}
      \left\lvert
      \sum_{w' \in W_{M'}}
      \sum_{w \in W(A,M)}
      W^P[M_{w' w} \Phi](h g)
      \right\rvert^2
      \, \frac{d h}{\delta_{U_P}(h)}
      \, d g,
  \end{align}
  where $P$ runs over
  \begin{equation*}
    \{\text{standard maximal parabolic subgroups $P$ of $G$}\} \cup \{G\}.
  \end{equation*}
\end{proposition}
\begin{proof}
  We first apply Proposition \ref{prop:sub-gln:let-varphi-in-preilm-unfold}, as reformulated in \S\ref{sec:reform-terms-const-1}, to evaluate the LHS of \eqref{prop:=-sum_p-=}.  The integral $\int_{N(\mathbb{A}) \backslash G(A)} \phi |W_\varphi|^2$ matches up with the $P = G$ contribution to the RHS of \eqref{prop:=-sum_p-=}, since in that case $M'$, $U_P$ and the set $W(A,M)$ are all trivial.  Following \eqref{eq:sum_-p-=}, it remains to check that for $P \neq G$ and all $g \in G(\mathbb{A})$,
  \begin{align*}
    &\int_{h \in [M']}
      \left\lvert
      \int_{u \in [N'']}
      \Eis(\Phi)_P(h u g) \psi^{-1}(u) \, d u
      \right\rvert^2 \, \frac{d h}{\delta_{U_P}(h)} \\
    &\quad =
      \frac{1}{|W_{M'}|}
      \int_{h \in [M']_{B_{M'}}}
      \left\lvert
      \sum_{
      \substack{
      w' \in W_{M'}  \\
    w \in W(A,M)
    }
    }
    W^P[M_{w' w} \Phi](h g)
    \right\rvert^2
    \, \frac{d h}{\delta_{U_P}(h)}.
  \end{align*}
  To that end, we first apply the constant term formula \eqref{eq:eis_bgf_p-=-sum} to evaluate
  \begin{equation*}
    \Eis(\Phi)_P(h u g)
    = \sum_{w \in W(A,M)} \Eis_{B_M}^M(M_w \Phi)(h u g).
  \end{equation*}
  We next apply the bilinearization of Lemma \ref{lem:sub-gln:let-p-=} to the elements $\Phi_1, \Phi_2 \in \mathcal{S}([M]_{B_M})$ defined, for $w_1, w_2 \in W(A,M)$, by
  \begin{equation*}
    \Phi_1(m) :=  M_{w_1} \Phi(m g), \quad
    \Phi_2(m) :=  M_{w_2} \Phi(m g).
  \end{equation*}
  Summing up the result of that application, we obtain
  \begin{align*}
    &\int_{h \in [M']}
      \left\lvert
      \int_{u \in [N'']}
      \Eis(\Phi)_P(h u g) \psi^{-1}(u) \, d u
      \right\rvert^2 \, \frac{d h}{\delta_{U_P}(h)} \\
    &\quad =
      \frac{1}{|W_{M'}|}
      \int_{h \in [M']_{B_{M'}}}
      \left\lvert
      \sum_{
      \substack{
      w' \in W_{M'}  \\
    w \in W(A,M)
    }
    }
    \int_{u \in N''(\mathbb{A})}
    M_{w' w} \Phi (h w_{M''}^{-1} u g)
    \psi^{-1}(u) \, d u
    \right\rvert^2
    \, \frac{d h}{\delta_{U_P}(h)}.
  \end{align*}
  Since $M'(\mathbb{A})$ and $M''(\mathbb{A})$ commute, we have
  \begin{align*}
    \int_{u \in N''(\mathbb{A})}
    M_{w' w} \Phi (h w_{M''}^{-1} u g)
    \psi^{-1}(u) \, d u
    &=
      \int_{u \in N''(\mathbb{A})}
      M_{w' w} \Phi (w_{M''}^{-1} u h g)
      \psi^{-1}(u) \, d u \\
    &= W^P[M_{w' w} \Phi](h g).
  \end{align*}
  The required identity follows.
\end{proof}

\section{Reduction of the proof}\label{sec:reduction-proof}
Throughout this section we consider fixed general linear groups over $\mathbb{R}$ or $\mathbb{Z}$.  We denote by $\psi$ the standard nondegenerate unitary character of $N$ (when working over $\mathbb{R}$) or of $N(\mathbb{A})$ (when  working over $\mathbb{Z}$), as described in \S\ref{sec:stand-nond-char}.  In particular, when applying the results of \S\ref{sec:rank-selb-unfold}, we employ this choice of $\psi$.

\subsection{A class of local vectors}
Let $T$ be a positive parameter with $T \ggg 1$.  Let $G$ be a fixed general linear group over $\mathbb{R}$.

\subsubsection{Basic definition}

\begin{definition}
  Let $\mathfrak{E}(N \backslash G, T)$ \index{classes!$\mathfrak{E}(N \backslash G, T)$, $\mathfrak{E}(N \backslash G, T)^W$} denote the class of all $f \in \mathcal{S}^e(N \backslash G)$ such that for each fixed compact subset $\mathcal{D}$ of $\mathfrak{a}^*$ and fixed $\ell \in \mathbb{Z}_{\geq 0}$, we have, with notation as in \eqref{eq:defn-frak-E-W},
  \begin{equation*}
    \nu_{\mathcal{D},\ell,T}(f) \ll T^{o(1)}.
  \end{equation*}
  We set
  \begin{equation}\label{eq:defn-frak-E-W}
    \mathfrak{E}(N \backslash G, T)^W := \mathfrak{E}(N \backslash G, T) \cap \mathcal{S}^e(N \backslash G)^W,
  \end{equation}
  with notation as in \S\ref{sec:fourier-transforms} (recall the general convention $U = N$ of Part \ref{part:local-l2-growth}, and Remark \ref{rmk:we-note-that-dependence-psi-w} concerning $\mathcal{S}^e(N \backslash G)^W$).
\end{definition}

The informal content of this definition is that an element $f \in \mathcal{S}^e(N \backslash G)$ belongs to $\mathfrak{E}(N \backslash G, T)$ if it is ``concentrated near $T^{-\rho^\vee} K$, with $A$-frequency $\O(1)$ and $G$-frequency $\O(T)$, and essentially $L^2$-normalized.''  This informal description is made precise by the following lemmas.

\begin{lemma}\label{lem:standard:frakE-equivalences}
  For $f \in \mathcal{S}^e(N \backslash G)$, the following are equivalent.
  \begin{enumerate}[(i)]
  \item\label{itm:lemma-class-frak-E-basic:1} $f \in \mathfrak{E}(N \backslash G, T)$.
  \item\label{itm:lemma-class-frak-E-basic:2}  For all $s \in \mathfrak{a}_{\mathbb{C}}^*$ with $\Re(s) \ll 1$, fixed $m \in \mathbb{Z}_{\geq 0}$ and fixed $y_1,\dotsc,y_m \in \Lie(G)$, we have
    \begin{equation}\label{eq:x_1-dotsb-x_m}
      \|R(y_1 \dotsb y_m) f[s]\| \ll T^{\rho^\vee(s) + m + o(1)} \langle s \rangle^{-\infty}.
    \end{equation}
  \item\label{itm:lemma-class-frak-E-basic:3} For all fixed $x \in \mathfrak{U}(A)$, fixed $m \in \mathbb{Z}_{\geq 0}$ and fixed $y_1,\dotsc,y_m \in \Lie(G)$, we have for all $(a,k) \in A \times K$ the estimate
    \begin{equation}\label{eq:int-_k-in}
      \int_{k \in K}
      \left\lvert
        L(x) R(y_1 \dotsb y_m) f(a k)
      \right\rvert^2 \, d k
      \ll
      \delta_N(a)
      \|T^{\rho^\vee} a\|^{-\infty}
      T^{2 m + o(1)}.
    \end{equation}
  \end{enumerate}
\end{lemma}
\begin{proof}
  The equivalence between \eqref{itm:lemma-class-frak-E-basic:1} and \eqref{itm:lemma-class-frak-E-basic:2} is immediate from the definitions (specifically, that $\Re(s)$ lies in a fixed compact set if and only if $\Re(s) \ll 1$).

  We now verify the equivalence between \eqref{itm:lemma-class-frak-E-basic:2} and \eqref{itm:lemma-class-frak-E-basic:3} (which may be understood as a form of Paley--Wiener theory for certain vector-valued functions on $A$).  For each $\sigma \in \mathfrak{a}^*$, we have the Mellin expansion
  \begin{align*}
    f(a k)
    &= \int_{\Re(s) = \sigma } |a|^{s + \rho} f[s](k) \, d \mu_A(s) \\
    &= \int_{\Re(s) = \sigma } |a|^{\rho^\vee} \left\lvert a T^{\rho^\vee } \right\rvert^{s} \frac{f[s](k)}{T^{\rho^\vee(s)}} \, d \mu_A(s),
  \end{align*}
  and also its variant for derivatives: for $(x,y) \in \mathfrak{U}(A) \times \mathfrak{U}(G)$,
  \begin{equation*}
    L(x) R(y) f(a k) =
    \int_{\Re(s) = \sigma }
    \langle x, s \rangle
    |a|^{\rho^\vee} \left\lvert a T^{\rho^\vee } \right\rvert^{s} \frac{R(y) f[s](k)}{T^{\rho^\vee(s)}} \, d \mu_A(s).
  \end{equation*}
  Here $x \mapsto \langle x, s \rangle$ denotes the algebra functional attached to $s \in \mathfrak{a}_{\mathbb{C}}^*$.  Abbreviate $\left\langle s \right\rangle := 1 + |s|$.  Fix $\ell \in \mathbb{Z}_{\geq 0}$ so that $\int_{\Re(s) = \sigma} \left\langle s \right\rangle^{-\ell} < \infty$.  Then by Cauchy--Schwarz,
  \begin{equation}\label{eq:lx-ry-fa}
    |L(x) R(y) f(a k)|^2 \leq
    \delta_N(a)
    \left\lvert a T^{\rho^\vee } \right\rvert^{2 \sigma }
    \int_{\Re(s) = \sigma }
    \left\langle s \right\rangle^{-\ell}
    \left\lvert \langle x, s \rangle \right\rvert^2
    \left\lvert \frac{R(y) f[s](k)}{T^{\rho^\vee(s)}} \right\rvert^2
    \, d \mu_A(s).
  \end{equation}
  Integrating this inequality over $k \in K$ and taking $\sigma \in \mathfrak{a}^*$ fixed but arbitrary, we see that condition \eqref{itm:lemma-class-frak-E-basic:2} implies condition \eqref{itm:lemma-class-frak-E-basic:3}.  A similar argument, using the Mellin identity
  \begin{equation*}
    \langle x,  s \rangle
    \frac{R(y) f[s](k)}{T^{\rho^\vee(s)}}
    =
    \int_{a \in A}
    \left\lvert a T^{\rho^\vee} \right\rvert^{-s}
    |a|^{-\rho^\vee} L(x) R(y) f(a k) \, d a
  \end{equation*}
  and Cauchy--Schwarz, gives the converse implication.
\end{proof}

\begin{lemma}\label{lem:f-in-mathfraken-square-integral-frakE}
  For $f \in \mathfrak{E}(N \backslash G, T)$, we have for each fixed $\ell,m \in \mathbb{Z}_{\geq 0}$, $x \in \mathfrak{U}(A)$ and $y_1,\dotsc,y_m \in \Lie(G)$ the estimate
  \begin{equation*}
    \int_{g \in N \backslash G } \|T^{\rho^\vee} g\|_{N \backslash G}^{\ell}
    \left\lvert L(x) R(y_1 \dotsb y_m) f \right\rvert^2(g) \, d g
    \ll
    T^{2 m + o(1)}.
  \end{equation*}
\end{lemma}
\begin{proof}
  The LHS of the above is
  \begin{equation*}
    \asymp
    \int_{a \in A}
    \|T^{\rho^\vee} a\|^{\ell}
    \int_{k \in K}
    \left\lvert L(x) R(y_1 \dotsb y_m) f \right\rvert^2(a k)
    \, d k
    \, \frac{d a}{ \delta_N(a)}.
  \end{equation*}
  For large enough fixed $d \geq 0$, we have
  \begin{equation*}
    \int_{a \in A} \|T^{\rho^\vee} a\|^{-d} \, d a
    =
    \int_{a \in A} \|a\|^{-d} \, d a
    \ll 1.
  \end{equation*}
  Thus the required estimate follows from \eqref{eq:int-_k-in}.
\end{proof}

\subsubsection{Duality}\label{sec:duality}

\begin{definition}
  Let $w_G \in W$ denote the long Weyl element.  For a function $f : N \backslash G \rightarrow \mathbb{C}$, define the dual function $\tilde{f} : N \backslash G \rightarrow \mathbb{C}$ by the formula
  \begin{equation*}
    \tilde{f}(g) := f(w_G g^{-\transpose}).
  \end{equation*}
\end{definition}
\begin{lemma}\label{lem:standard:let-f-in}
  We have
  \begin{align*}
    f \in \mathfrak{E}(N \backslash G, T) &\implies \tilde{f} \in \mathfrak{E}(N \backslash G, T), \\
    f \in \mathfrak{E}(N \backslash G, T)^W &\implies \tilde{f} \in \mathfrak{E}(N \backslash G, T)^W.
  \end{align*}
\end{lemma}
\begin{remark}
  We have noted that, informally, $\mathfrak{E}(N \backslash G,T)$ consists of functions concentrated near $T^{-\rho^\vee} K$, with $A$-frequency $\ll 1$ and $G$-frequency $\ll T$.  We have
  \begin{equation*}
    w_G (T^{-\rho^\vee})^{-\transpose} = T^{-\rho^\vee} w_G,
    \quad
    w_G (T^{-\rho^\vee} K)^{-\transpose} = T^{-\rho^\vee} K,
  \end{equation*}
  so the duality assignment $f \mapsto \tilde{f}$ preserves the concentration condition.  It is also twisted-equivariant, hence preserves the frequency conditions.  Finally, since it comes from an automorphism of the group, it should come as no surprise  that $f \mapsto \tilde{f}$ preserves the $W$-invariance condition.
\end{remark}
\begin{proof}
  By direct computation, we see that $\tilde{f}$ defines an element of $\mathcal{S}^e(N \backslash G)$, with Mellin components given by
  \begin{equation*}
    \tilde{f}[s](g) = f[- w_G s](w_G g^{-\transpose}).
  \end{equation*}
  Recalling the definitions of $\mathcal{F}_w$ (\S\ref{sec:fourier-transforms}) and $\nu_{\mathcal{D},\ell,T}$ (\S\ref{sec:local-sobolev-norms-prelims}), our task is check the following:
  \begin{enumerate}[(a)]
  \item \label{itm:local-standard:1b} If $f[w s] = \mathcal{F}_w f[s]$ for all $w \in W$, then $\tilde{f}[w s] = \mathcal{F}_{w} \tilde{f}[s]$ for all $w \in W$.
  \item \label{itm:local-standard:2b} For $\Re(s) \ll 1$, fixed $m \in \mathbb{Z}_{\geq 0}$, and fixed $x_1,\dotsc,x_m \in \Lie(G)$, we have
    \begin{equation*}
      \|R(x_1 \dotsb x_m) \tilde{f}[s]\| \ll T^{\rho^\vee(s) + m + o(1)} \langle s \rangle^{-\infty}.
    \end{equation*}
  \end{enumerate}
  Assertion \eqref{itm:local-standard:1b} reduces to the following identity: for $v \in \mathcal{I}(s)$, with dual $\tilde{v} \in \mathcal{I}(- w_G s)$ given by $\tilde{v}(g) := v(w_G g^{-\transpose})$, we have
  \begin{equation*}
    W[\tilde{v},\psi](g) =W[v, \psi^{-1}](w_G g^{-\transpose}).
  \end{equation*}
  To see this, we use that the involution $g \mapsto w_G g^{-\transpose} w_G$ of $G$ preserves $N$ and inverts $\psi$.  Assertion \eqref{itm:local-standard:2b} follows from \eqref{eq:x_1-dotsb-x_m} and the fact that $\langle -w_G s \rangle \asymp \langle s \rangle$ and $\rho^\vee(-w_G s) = \rho^\vee(s)$.
\end{proof}

\subsection{Construction of theta functions}
We fix a natural number $n$ and let $G := \GL_n, B, N, A, V_n$ be as in \S\ref{sec:notation-rs-unfolding} and the beginning of Part \ref{part:local-l2-growth}.

Let $t \in A(\mathbb{R})^0_{\geq 1}$, thus $t = \diag(t_1,\dotsc,t_n)$ with $t_1 \geq \dotsb \geq t_n > 0$.  For each prime number $p$, let $\phi_p \in C_c^\infty(V_n(\mathbb{Q}_p))$ denote the characteristic function of $V_n(\mathbb{Z}_p)$.  We fix a compact subset $\Omega$ of $G(\mathbb{R})$, and define the seminorm $\|.\|_{t,\Omega}$ on $\mathcal{S}([G])$ by \eqref{eq:varphi-2_r-:=}.

\begin{lemma}\label{lem:sub-gln:let-r-be}
  There exists $\phi_\infty \in C_c^\infty(V_n(\mathbb{R}))$ with the following properties:
  \begin{enumerate}[(i)]
  \item $0 \leq \phi_\infty \leq 1$.  In particular, setting $\phi := \otimes_{\mathfrak{p}} \phi_\mathfrak{p} \in C_c^\infty(V_n(\mathbb{A}))$, we have $\theta_{\phi} \geq 0$.
  \item $\phi_\infty$ is $K_\infty$-invariant, hence $\phi$ is invariant by $K = \prod_\mathfrak{p} K_\mathfrak{p}$.
  \item For each $v \in V_n(\mathbb{R})$, we have $\phi_\infty(v) \neq 0$ only if $\|v\| \asymp t_n$.
  \item For each $u \in N(\mathbb{R})$ and $g \in \Omega$, we have
    \begin{equation}\label{eq:theta_phiu-r-g}
      \theta_\phi(u t g) \geq 1.
    \end{equation}
  \end{enumerate}
\end{lemma}
\begin{proof}
  We fix a $[0,1]$-valued function $f \in C_c^\infty(\mathbb{R}^\times_+)$ taking the value $1$ on the interval $[1/c,c]$, where $c > 2$ is fixed sufficiently large in terms of $\Omega$, and take $\phi_\infty(v) := f(t_n^{-1} \|v\|)$.  The first three properties are then clear by construction.  By our choice of the finite components $\phi_p$, we have
  \begin{equation*}
    \theta_\phi(u t g) = \sum_{0 \neq v \in V_n(\mathbb{Z})}
    \phi_\infty(v u t g).
  \end{equation*}
  We consider the contribution of the basis vector $e_n \in V_n(\mathbb{Z})$.  We have $e_n u t g = t_n e_n g$, thus
  \begin{equation*}
    \theta_\phi(u t g) \geq \phi_\infty(e_n u t g)
    = f(\|e_n g\|).
  \end{equation*}
  By our construction of $c$ and the fact that $g \in \Omega$, we have $\|e_n g\| \in [1/c,c]$, hence $f( \|e_n g\|) = 1$.  The required inequality follows.
\end{proof}

\begin{corollary}\label{cor:sub-gln:let-psi-in}
  Let $\Psi \in \mathcal{S}([G])$ be right-invariant under $G(\mathbb{Z}_p)$ for each prime $p$, so that $\Psi$ identifies with an element of $\mathcal{S}(G(\mathbb{Z}) \backslash G(\mathbb{R}))$.  Let $\phi_\infty$ and $\phi$ be as in Lemma \ref{lem:sub-gln:let-r-be}.   Then
  \begin{equation*}
    \|\Psi \|_{t,\Omega}^2
    \ll
    \delta_N(t)
    \int_{[\GL_n]} \theta_\phi |\Psi|^2.
  \end{equation*}
\end{corollary}
\begin{proof}
  We recall from Proposition \ref{prop:siegel-lower} that for each nonnegative integrable function $\phi$ on $[G]$,
  \begin{equation*}
    \int_{\Omega \times \prod_p K_p} \varphi(t g) \, d g \ll \delta_N(t) \int_{[G]} \varphi.
  \end{equation*}
  Taking $\varphi := \theta_\phi |\Psi|^2$, we have by  \eqref{eq:theta_phiu-r-g} that $|\Psi|^2(t g) \leq \varphi(t g)$ for all $g \in \Omega$.  The required conclusion is then immediate from the definition of $\|.\|_{t,\Omega}$.
\end{proof}

\subsection{Recap}
We recall the setting of Theorem \ref{thm:growth-eisenstein-nonstandard}:
\begin{assumptions}\label{assumptions-part-5}
  Suppose given the following.
  \begin{itemize}
  \item $n$: a fixed natural number, $G := \GL_n$ over $\mathbb{Z}$.
  \item $\Omega$: a fixed compact subset of $G(\mathbb{R})$.
  \item $T \geq 1$.
  \item $Y \in A(\mathbb{R})^0$ with $T^{-\kappa} \leq |Y_j| \leq T^{\kappa}$ for all $j$ and some fixed $\kappa \in (0,1/2)$.
  \item $t \in A(\mathbb{R})^0_{\geq 1}$.
  \end{itemize}
\end{assumptions}
For $f \in \mathcal{S}(N(\mathbb{R}) \backslash G(\mathbb{R}))$, define  the pseudo Eisenstein series
\begin{equation*}
  \Psi[f,Y] := \Eis[\Phi[L(Y) f]].
\end{equation*}
In view of the homogeneity of the desired estimate with respect to $f$, it is enough to show the following:
\begin{theorem}[Theorem \ref{thm:growth-eisenstein-nonstandard}, reformulated]
  Retain Assumptions \ref{assumptions-part-5}.  For $f \in \mathfrak{E}(N(\mathbb{R}) \backslash G(\mathbb{R}), T)^W$, we have
  \begin{equation}\label{eq:psiy_r2-ll-dety}
    \|\Psi[f,Y]\|_{t,\Omega}^2 \ll \min ( \det(Y)^{-1} t_1^{-n}, \det(Y) t_n^n ) \delta_N(t) T^{o(1)}.
  \end{equation}
\end{theorem}

We now make some preliminary reductions.  It suffices to consider the case $T \ggg 1$, since if $T \ll 1$, then the required estimate follows readily from the continuity of the maps $\Phi : \mathcal{S}(N(\mathbb{R}) \backslash G(\mathbb{R})) \rightarrow \mathcal{S}([G]_B)$ and $\Eis : \mathcal{S}([G]_B) \rightarrow \mathcal{S}([G])$.  For the same reason, we may reduce to the case that each component $t_j$ is of the form $T^{\O(1)}$ (compare with the proof of Lemma \ref{lem:standard2:each-fixed-c}).  Finally, exactly as in \S\ref{sec:appl-crude-bounds}, we may reduce to the case that $\det(t Y) = T^{o(1)}$, noting that Lemmas \ref{lem:each-fixed-c_0-Psi-crude} and \ref{lem:standard2:each-fixed-c_0} remain valid in our setting, with the same proofs.  It will thus suffice to show the following:
\begin{theorem}\label{thm:sub-gln:let-r-as}
  Retain Assumptions \ref{assumptions-part-5}.  Assume moreover that $\det(t Y) = T^{o(1)}$ and $t_j = T^{\O(1)}$ for all $j \in \{1, \dotsc, n\}$.  Then
  \eqref{eq:psiy_r2-ll-dety} holds.
\end{theorem}

\begin{remark}
  For the purpose of proving our main result, it was unnecessary to carry out the above reductions; it suffices to note that such assumptions are in force in the one place where we apply Theorem \ref{thm:growth-eisenstein-nonstandard}, in \S\ref{sec:appl-growth-bounds}.  We have repeated them here to simplify the statement of Theorem \ref{thm:growth-eisenstein-nonstandard} by keeping it free of unnecessary hypotheses.
\end{remark}

\subsection{Application of duality}\label{sec:cq2z5f8a4z}
Retaining the setting of Theorem \ref{thm:sub-gln:let-r-as}, we record the preliminary reduction that it suffices to verify the apparently weaker estimate: for all $f \in  \mathfrak{E}(N(\mathbb{R}) \backslash G(\mathbb{R}), T)^W$,
\begin{equation}\label{eq:psiy_r2-ll-dety-1}
  \|\Psi[f,Y]\|_{t,\Omega}^2 \ll
  \det(Y)
  t_n^n \delta_N(t) T^{o(1)}.
\end{equation}

To obtain this reduction, we first define, for any function $\varphi : [G] \rightarrow \mathbb{C}$, the dual function $\tilde{\varphi} : [G] \rightarrow \mathbb{C}$ by
\begin{equation*}
  \tilde{\varphi}(g) := \varphi(g^{-\transpose}),
\end{equation*}
where, as in \S\ref{sec:duality}, $g \mapsto g^{-\transpose}$ denotes the outer automorphism given by taking the inverse transpose.  Temporarily denote by $w_0 := w_G$ the standard representative for the long Weyl element, which satisfies $w_0 = w_0^{-1} = w_0^\transpose$.  Set $Y^{-w_0} := w_0^{-1} Y^{-1} w_0 \in A(\mathbb{R})^0$ and $\tilde{\Phi}(g) := \Phi(w_0 g^{-\transpose})$.  Then a simple calculation gives
\begin{equation*}
  \widetilde{\Psi[f,Y]} = \Psi[\tilde{f}, Y^{-w_0}],
\end{equation*}
where $\tilde{f}$ is as defined in \S\ref{sec:duality}.  In deriving this, we use that
\begin{equation*}
  \mathcal{P}_G(-w_0 s) = \mathcal{P}_G(s),
  \quad
  \zeta(N,-w_0 s) = \zeta(N,s)
\end{equation*}
and that, for finite primes $p$,
\begin{equation*}
  f_p[-w_0 s](w_0 g^{-\transpose}) = f_p[s](g),
\end{equation*}
which follows from the fact that both sides defines normalized spherical elements of the same induced representation.

We now explain the deduction of Theorem \ref{thm:sub-gln:let-r-as} from \eqref{eq:psiy_r2-ll-dety-1}.  By enlarging $\Omega$ if necessary, we may assume that it is invariant under $g \mapsto g^{-\transpose}$ and also under left and right translation by $w_0$.  Then
\begin{equation*}
  \tilde{t} := w_0 t^{-\transpose} w_0 = \diag(t_n^{-1}, \dotsc, t_1^{-1}),
  \quad
  \delta_N(\tilde{t}) = \delta_N(t)
\end{equation*}
and
\begin{equation*}
  \|\Psi[f,Y]\|_{t,\Omega} = \|\Psi[\tilde{f}, Y^{-w_0}]\|_{\tilde{t},\Omega}.
\end{equation*}
We have $\det(Y^{-w_0}) = \det(Y)^{-1}$.  By Lemma \ref{lem:standard:let-f-in}, we have $\tilde{f} \in \mathfrak{E}(N(\mathbb{R}) \backslash G(\mathbb{R}), T)^W$.  By \eqref{eq:psiy_r2-ll-dety-1} applied to $\tilde{f}$ and $Y^{-w_0}$, it follows that
\begin{equation}\label{eq:psiy_r2-ll-dety-2}
  \|\Psi[f,Y]\|_{t,\Omega}^2 \ll
  \det(Y)^{-1} t_1^{-n} \delta_N(t) T^{o(1)}.
\end{equation}
Taking the minimum of the two estimates \eqref{eq:psiy_r2-ll-dety-1} and \eqref{eq:psiy_r2-ll-dety-2} yields the required estimate \eqref{eq:psiy_r2-ll-dety}.

\subsection{Reduction to bounds for Rankin--Selberg integrals}\label{sec:cq2z5fy3gx}
We recall the assignment
\begin{equation*}
  \Phi : \mathcal{S}^e(N(\mathbb{R}) \backslash G(\mathbb{R})) \rightarrow \mathcal{S}^e([G]_B)
\end{equation*}
\begin{equation*}
  f \mapsto \Phi[f]
\end{equation*}
characterized, as in \S\ref{sec:compl-eisenst-seri}, by the Mellin component identity
\begin{equation}\label{eq:phif_-=-mathc}
  \Phi[f][s] = \mathcal{P}_G(s) \zeta(N,s) \otimes_{\mathfrak{p}} f_\mathfrak{p}[s]
\end{equation}
where $f_\infty[s] := f[s]$, while $f_p[s]$ is normalized spherical for all finite primes $p$.

Our aim is to prove Theorem \ref{thm:sub-gln:let-r-as} via the criterion \eqref{eq:psiy_r2-ll-dety-1}.  In view of Corollary \ref{cor:sub-gln:let-psi-in}, it is enough to show that, with $\phi$ as in the statement of the latter and $f \in \mathfrak{E}(N \backslash G, T)^W$,
\begin{equation}\label{eq:int_g-thet-psiy2}
  \int_{[G]}
  \theta_\phi \left\lvert \Psi[f,Y] \right\rvert^2 \ll
  \det(Y)
  t_n^n  T^{o(1)},
\end{equation}
where $\Psi[f,Y] = \Eis[L(Y) \Phi[f]]$.

By Proposition \ref{prop:sub-gln:let-phi-in-summary}, we reduce further to verifying the same bound for each of the following integrals, defined for $P = MU$ either a standard maximal parabolic subgroup or $P = G$ itself:
\begin{align*}
  &\int_{g \in (M' N'' {U_P} \backslash G)(\mathbb{A}) } \phi(e_n g)
    \int_{h \in [M']_{B_{M'}}} \\
  &\quad  \quad
    \left\lvert
    \sum_{w' \in W_{M'}}
    \sum_{w \in W(A,M)}
    W^P[M_{w' w} L(Y) \Phi[f]](h g)
    \right\rvert^2
    \, \frac{d h}{\delta_{U_P}(h)}
    \, d g.
\end{align*}
Since $\phi$ is nonnegative, we can estimate the above expression by bounding the $h$-integral in magnitude and integrating the resulting bound.  In particular, using the triangle inequality, the sums over Weyl group elements may be taken outside of the absolute value.  By Lemma \ref{lem:standard2:let-f_infty-in} and our assumption that $f$ is $W$-invariant, we have
\begin{equation*}
  M_{w} L(Y) \Phi[f] = L({}^w Y) \Phi[f]
\end{equation*}
for all $w \in W_G$.  Since $\det(Y) = \det({}^w Y)$, and since ${}^w Y$ satisfies the hypotheses of Theorem \ref{thm:sub-gln:let-r-as} if and only if $Y$ does, we reduce to estimating the simpler integrals $\mathcal{Q}_P(\phi, L(Y) f)$, where for $f \in \mathcal{S}^e(N(\mathbb{R}) \backslash G(\mathbb{R}))$, we set
\begin{equation}\label{eq:mathcalq_pphi-f-:=}
  \mathcal{Q}_P(\phi, f)
  :=
  \int_{g \in (M' N'' {U_P} \backslash G)(\mathbb{A}) } \phi(e_n g)
  \int_{h \in [M']_{B_{M'}}} \\
  \left\lvert
    W^P[\Phi[f]](h g)
  \right\rvert^2
  \, \frac{d h}{\delta_{U_P}(h)}
  \, d g.
\end{equation}

We note that when $P = G$, the integrals $\mathcal{Q}_G(\phi,f)$ are what one might call ``pseudo local Rankin--Selberg integrals for $\GL_n \times \GL_n$:'' a spectral expansion of $\phi$ and $f$ would lead to an expression for $\mathcal{Q}_G(\phi,f)$ as an integral of ordinary local Rankin--Selberg integrals.  For $P \neq G$, they might be understood as degenerate variants of such integrals.

\subsection{Summary thus far}\label{sec:summary-thus-far}
The proof of Theorem \ref{thm:sub-gln:let-r-as} thereby reduces to that of the following proposition.  We summarize the running assumptions:
\begin{assumptions}\label{assumptions:furth-reduct-proof}
  Suppose given
  \begin{itemize}
  \item a positive parameter $T \ggg 1$,
  \item an element $f \in \mathfrak{E}(N(\mathbb{R}) \backslash G(\mathbb{R}), T)$,
  \item an element $r > 0$ with $r =T^{\O(1)}$ (playing the role of $t_n$), and
  \item an element $Y \in A(\mathbb{R})^0$ with $T^{-\kappa} \leq |Y_j| \leq T^{\kappa}$ for some fixed $\kappa < 1/2$ and all $j$.
  \end{itemize}
  Let $V := V_n$ denote the punctured affine $n$-plane, as in \S\ref{sec:notation-rs-unfolding}.  Let $\phi_p \in C_c^\infty(V(\mathbb{Q}_p))$ be the characteristic function of $V(\mathbb{Z}_p)$.  Let $\phi_\infty \in C_c^\infty(V(\mathbb{R}))$ be the element attached, as in Lemma \ref{lem:sub-gln:let-r-be}, to the quantity $r > 0$, and let $\phi = \otimes_{\mathfrak{p}} \phi_\mathfrak{p}$.
\end{assumptions}

We have reduced the proof of Theorem \ref{thm:sub-gln:let-r-as} to that of the following.
\begin{proposition}\label{prop:bound-after-theta-degenerate-case-2}
  Retain Assumptions \ref{assumptions:furth-reduct-proof}.  Let $P$ be either $G$ itself or a standard maximal parabolic subgroup of $G$.  Then
  \begin{equation*}
    \mathcal{Q}_P(\phi,L(Y) f)
    \ll
    \det(Y) r^{n} T^{o(1)}.
  \end{equation*}
\end{proposition}

\subsection{Interlude on basic vectors}\label{sec:interl-basic-vect}
We suspend the general setting of \S\ref{sec:reduction-proof} and pause to record some preliminaries, to be applied below, concerning ``basic vectors'' following Braverman--Kazhdan \cite[\S3.12]{MR1694894}.

\subsubsection{Local}\label{sec:basic-vectors-local}
Let $F$ be a non-archimedean local field, with ring of integers $\mathfrak{o}$, maximal ideal $\mathfrak{p}$ and $q = \# \mathfrak{o} / \mathfrak{p}$, as usual.  Let $G$ be a split reductive group over $\mathfrak{o}$, with maximal torus $A$ and Borel subgroup $B = N A$.

As in \S\ref{sec:local-zeta-functions}, we define the meromorphic function $\zeta_F(N,\cdot)$ on $\mathfrak{a}_{\mathbb{C}}^*$ by
\begin{equation*}
  \zeta_F(N,s) := \prod_{\alpha > 0} \zeta_F(1 + \alpha^\vee(s)),
\end{equation*}
where the product is over all positive (with respect to $N$) roots.  (We have switched notation from ``$U$'' there to ``$N$'' here, following the general convention ``$U = N$'' of Part \ref{part:local-l2-growth}.)

The meromorphic function $\zeta_F(N,s)$ depends only upon the unramified character $|.|^s$ of $A(F)$.  It is holomorphic in the domain where, for each $\alpha > 0$, we have $\Re(\alpha^\vee(s)) > -1$; in particular, it is holomorphic where $\Re(s)$ is dominant.  Thus for any dominant $\sigma \in \mathfrak{a}^*$, we may define a smooth function $\Theta : A(F) / A(\mathfrak{o}) \rightarrow \mathbb{C}$ by the contour integral \index{basic vectors!$\Theta = \Theta_G$}
\begin{equation*}
  \Theta(a) = \int_{(\sigma)} |a|^{\rho + s} \zeta_F(N,s) \, d \mu_{A}(s),
\end{equation*}
defined as in \eqref{eq:int-_sigma-fs}.  Here $d \mu_{A}(s)$ denotes the \emph{probability} Haar measure on $\{|.|^s \in \mathfrak{X}^e(A) : \Re(s)  = \sigma \}$, dual to the Haar measure on $A(F)$ assigning volume one to $A(\mathfrak{o})$.

If $G$ is a torus, then $\Theta|_{A(F)}$ is the characteristic function of $A(\mathfrak{o})$.  Otherwise, $\zeta_F(N,s)$ is nontrivial, hence has poles, and so the integral depends upon the choice of $\sigma$.  From another perspective, $\Theta$ is not compactly supported, but instead has some asymptotics at $\infty$ along the ``negative cone.''

By the Iwasawa decomposition, the function $\Theta$ extends uniquely to a function on $N(F) \backslash G(F) / G(\mathfrak{o})$, right-invariant under $G(\mathfrak{o})$.  We write $\Theta = \Theta^G$ when we wish to indicate which group $G$ is under consideration.

As in \cite[\S3.12]{MR1694894}, let $\Gamma^+ \subseteq \mathfrak{a}^*$ denote the set of all nonnegative integral combinations of coroots $\alpha^\vee$, and let $\mathcal{K} : \Gamma^+ \rightarrow \mathbb{R}_{\geq 0}$ denote the function given by
\begin{equation*}
  \mathcal{K}(\gamma) = \sum_{P \in \mathcal{P}_\gamma } q^{-|P|},
\end{equation*}
where $\mathcal{P}_\gamma$ denotes the set of representations of $\gamma$ as a sum of positive coroots and $|P|$ denotes the number of summands in a given representation $P \in \mathcal{P}_\gamma$.  The series representation $\zeta_F(s) = (1 - q^{-s})^{-1} = \sum_{n \geq 0} q^{-n s}$  yields
\begin{equation*}
  \zeta_F(N,s)
  =
  \sum_{\gamma \in \Gamma^+ }
  \mathcal{K}(\gamma)
  q^{- \gamma(s) },
\end{equation*}
valid, e.g., for $s$ of dominant real part.  We may identify $\Gamma^+$ with a subset of $A(F) / A(\mathfrak{o})$ by associating to $\gamma \in \Gamma^+$ the unique element $a_\gamma \in A(F) / A(\mathfrak{o})$ for which $|a_\gamma|^s = q^{\gamma(s)}$ for all $s \in \mathfrak{a}_{\mathbb{C}}^*$.  Then by Mellin inversion, we see that $\Theta$ is supported on the image of $\Gamma^+$ and given there by
\begin{equation*}
  \Theta(a_\gamma) = q^{ \langle \rho, \gamma  \rangle} \mathcal{K}(\gamma).
\end{equation*}
In particular, for every rational character $\chi : G(F) \rightarrow F^\times$, we have $\chi(g) \in \mathfrak{o}^\times$ for all $g$ in the support of $\Theta$.  For instance, when $G = \GL_n$, every element of the support of $\Theta$ has unit determinant.

\begin{example}
  If $G = \SL_2$ and $a_\gamma = \diag(\varpi^{-r}, \varpi^r)$, then $\mathcal{K}(\gamma) = 1_{r \geq 0} q^{-r}$ and $\langle \rho, \gamma  \rangle = r$, so $\Theta(a_\gamma) = 1_{r \geq 0}$.  Under the ``bottom row'' map $N(F) \backslash G(F) \cong F^2 - \{0\}$, $\Theta$ identifies with the characteristic function of $1_{\mathfrak{o}^2 - \{0\}}$ (cf.\ Example \ref{exa:example-sl2-completed-eisenstein-series}).  In most other examples, $\Theta$ is more complicated, e.g., non-constant on its support.
\end{example}

\subsubsection{Global}

Now let $G$ be a split reductive group over $\mathbb{Z}$, with maximal torus $A$ and Borel $B = N A$.  (The present discussion is the specialization to $(\mathbb{Q}, \{\infty \}$) of a more general one concerning any pair consisting of number field and a large enough finite set of places.)  For each prime number $p$, the discussion of \S\ref{sec:basic-vectors-local} provides us with an element
\begin{equation*}
  \Theta^G_p  : N(\mathbb{Q}_p) \backslash G(\mathbb{Q}_p) / G(\mathbb{Z}_p) \rightarrow \mathbb{C}.
\end{equation*}

\begin{lemma}\label{lem:sub-gln:let-m-in}
  Let $\eps > 0$.  Define $f_\infty : N(\mathbb{R}) \backslash G(\mathbb{R}) \rightarrow \mathbb{C}$ in Iwasawa coordinates $(a,k) \in A(\mathbb{R}) \times K_\infty$ by the formula
  \begin{equation*}
    f_\infty(a k) := |a|^{2 \rho} \|a\|^{-\eps}.
  \end{equation*}
  Set $f_p := \Theta_p^G$.  Define $f : N(\mathbb{A}) \backslash G(\mathbb{A}) \rightarrow \mathbb{R}_{\geq 0}$ by $f(g) := \prod_\mathfrak{p} f_\mathfrak{p}(g_\mathfrak{p})$.  Then the series $\sum_{\gamma \in N(\mathbb{Q}) \backslash G(\mathbb{Q})} f(\gamma)$ converges.
\end{lemma}
\begin{proof}
  This seems well-known (see, e.g., \cite[\S7.6]{MR1694894}).  Since we could not quickly locate a written proof, we record the necessary arguments for the sake of completeness.


  We may and shall assume that all quantities in our hypotheses are fixed, so that we may employ asymptotic notation effectively.

  We define the intermediate function $\Phi : [G]_B \rightarrow \mathbb{R}_{\geq 0}$ by the series
  \begin{equation*}
    \Phi(g) := \sum_{c \in A(\mathbb{Q})} f(c g).
  \end{equation*}
  Then
  \begin{equation*}
    \sum_{N(\mathbb{Q}) \backslash G(\mathbb{Q})} f(\gamma) = \sum_{B(\mathbb{Q}) \backslash G(\mathbb{Q})} \Phi(\gamma).
  \end{equation*}
  By the convergence of ``ordinary'' Eisenstein series (\S\ref{sec:eisenstein-series}), it will suffice to verify for some fixed strictly dominant element $\sigma \in \mathfrak{a}^*$ that for all $(a,k) \in [A] \times K$,
  \begin{equation}\label{eq:phia-k-ll}
    \Phi(a k) \ll |a|^{2 \rho + \sigma}.
  \end{equation}
  Since $f$ is right-invariant under $K$, so is $\Phi$, and we may thus reduce to the case $k = 1$.  We fix a small enough neighborhood $\Omega_\infty$ of the identity element in $A(\mathbb{R})^0$, and set $\Omega := \Omega_\infty \prod_p A(\mathbb{Z}_p) \subseteq A(\mathbb{A})$.  Then the cosets $c \Omega$, with $c \in A(\mathbb{Q})$, are disjoint from one another.  Moreover, for $u \in A(\mathbb{A})$ and $t \in \Omega$, we have $f(u t) \asymp f(u)$, hence $f(u) \asymp \int_{t \in \Omega } f(u t) \, d t$, where $d t$ denotes our Haar measure on $A(\mathbb{A})$.  We obtain

  \begin{equation*}
    \Phi(a) = \sum_{c \in A(\mathbb{Q})} f(c a) \asymp  \int_{t \in A(\mathbb{A})}  1_{A(\mathbb{Q}) \Omega}(t) f(t a) \, d t.
  \end{equation*}
  For $t \in A(\mathbb{Q}) \Omega$, the product rule implies that $|t|^{2 \rho + \sigma} \asymp 1$.  Thus
  \begin{align*}
    |a|^{-2 \rho-\sigma} \Phi(a)
    &\ll \int_{t \in A(\mathbb{A})} 1_{A(\mathbb{Q}) \Omega}(t) |t a|^{-2 \rho - \sigma} f(t a) \, d t \\
    &\leq \int_{t \in A(\mathbb{A})} |t|^{-2 \rho - \sigma} f(t) \, d t \\
    &= \prod_\mathfrak{p} \int_{t \in A(\mathbb{Q}_\mathfrak{p})} |t|^{-2 \rho - \sigma} f_\mathfrak{p}(t) \, d t.
  \end{align*}
  The $p$-adic integrals evaluate to $\zeta_p(N,\rho + \sigma)$, whose product converges (comfortably) in view of the strict dominance of $\rho + \sigma$.  The archimedean integral
  \begin{equation*}
    \int_{t \in A(\mathbb{R})} |t|^{- 2 \rho - \sigma} f_\infty(t) \, d t
    =
    \int_{t \in A(\mathbb{R})} |t|^{- \sigma} \|t\|^{-\eps} \, d t
  \end{equation*}
  converges provided that $\sigma$ is fixed sufficiently small in terms of $\eps$.  Thus $|a|^{- 2 \rho - \sigma} \Phi(a) \ll 1$, as required.
\end{proof}

\subsection{Series representations via partial Schwartz spaces}\label{sec:seri-repr-via}
We now resume the setting of \S\ref{sec:summary-thus-far} and derive alternative series representations for $\Phi[f]$ that will be very convenient for certain calculations (see the final paragraphs of \S\ref{sec:strong-approximation}).

Let $P = M {U_P}$, $M = M' \times M''$ be either $G$ or a standard maximal parabolic, with accompanying notation as in \S\ref{sec:reform-terms-const-1} and \S\ref{sec:summary-unification}.  We introduce the additional notation
\begin{equation*}
  A' \leq M', \quad A'' \leq M''
\end{equation*}
for the diagonal tori in the general linear groups $M'$ and $M''$.

Recall from \S\ref{eq:zeta-U-s-global-over-Z} that we have defined a meromorphic function
\begin{equation*}
  \zeta(N,s) = \prod_{\alpha > 0} \zeta(1 + \alpha^\vee(s))
\end{equation*}
of $s \in \mathfrak{a}_{\mathbb{C}}^*$.  We now introduce the analogous definition
\begin{equation*}
  \zeta(N'',s) := \prod_{
    \substack{
      \alpha > 0 :  \\
      \Lie(G)_\alpha \subseteq \Lie(N'')
    }
  }
  \zeta(1 + \alpha^\vee(s)),
\end{equation*}
obtained by restricting the product to those $\alpha$ coming from $N''$.  Equivalently, $\zeta(N'',s)$ is obtained via pullback under the map $\mathfrak{a}_{\mathbb{C}}^* \twoheadrightarrow (\mathfrak{a}_{\mathbb{C} }'')^*$, dual to the inclusion $A'' \hookrightarrow A$, of the analogous meromorphic function attached to $N'' \subseteq M''$.  The ratio $\zeta(N,s) / \zeta(N'', s)$ is then likewise a product of zeta functions:
\begin{equation*}
  \frac{\zeta (N, s)}{ \zeta (N'', s)}
  = \prod_{
    \substack{
      \beta > 0:   \\
      \Lie(G)_\beta \not\subseteq \Lie(N'')
    }
  }
  \zeta(1 + \beta^\vee(s)).
\end{equation*}

For each finite prime $p$, we denote by $f_p[s] \in \mathcal{I}_\mathfrak{p}(s)$ the spherical elements satisfying $f_p[s](1) = 1$, and define $\Theta^P_p \in C^\infty(N(\mathbb{Q}_p) \backslash G(\mathbb{Q}_p))$ by the Mellin integral \index{basic vectors!$\Theta^P_p$}
\begin{equation}\label{eq:thetap_p-:=-int}
  \Theta^P_p :=
  \int_{(\sigma)}
  \zeta_p(N'', s) f_p[s]
  \, d \mu_{A(\mathbb{Q}_p)}(s)
\end{equation}
for any dominant $\sigma \in \mathfrak{a}^*$.  For $s$ with dominant real part, the Mellin component $\Theta_p^P[s]$ is given by
\begin{equation*}
  \Theta_p^P[s] = \zeta_P(N'',s) f_p[s].
\end{equation*}
The function $\Theta_p^P$ descends to the quotient $N(\mathbb{Q}_p) \backslash G(\mathbb{Q}_p) / K_p$, hence is determined by its restriction to $A(\mathbb{Q}_p)/A(\mathbb{Z}_p)$.  In the special case $P = G$ (i.e., $M'' = G$), the function $\Theta_p^P = \Theta_p^G$ is as described in \S\ref{sec:basic-vectors-local}.  At the other extreme, if $\dim(M'') = 1$, then $\zeta(N'',s) = 1$, so $\Theta^P_p$ is the characteristic function of $K_p$.  The general case is a mixture of these two:
\begin{lemma}\label{lem:Theta-P-p}
  Let $a \in A(\mathbb{Q}_p)$.  Write $a = a' a''$, with $(a',a'') \in A'(\mathbb{Q}_p) \times A''(\mathbb{Q}_p)$.  Then

  \begin{equation}\label{eq:thet-=-1_am}
    \Theta_p^P(a) = 1_{A'(\mathbb{Z}_p)}(a')
    \Theta^{M''}_p(a''),
  \end{equation}
  where $\Theta^{M''}_p$ denotes the function defined like $\Theta^G_p$, but relative to the (typically smaller) general linear group $M''$.
\end{lemma}
\begin{proof}
  By definition, $f_p[s](a) = |a|^{\rho + s} = |a'|^{\rho + s} |a''|^{\rho + s}$.  Writing $\sigma = \sigma' + \sigma''$, the integral \eqref{eq:thetap_p-:=-int} factors as the product of the two integrals
  \begin{equation*}
    \int_{(\sigma')}
    |a'|^{\rho + s}
    \, d \mu_{A'(\mathbb{Q}_p)}(s),
    \quad
    \int_{(\sigma'')}
    |a''|^{\rho + s}
    \zeta_p(N'',s) \, d \mu_{A''(\mathbb{Q}_p)}(s).
  \end{equation*}
  The first integral evaluates to $1_{A'(\mathbb{Z}_p)}(a')$.  By definition,
  \begin{equation*}
    \Theta_p^{M''}(a'') = \int_{(\sigma'')} |a''|^{\rho'' + s} \zeta_p(N'',s) \, d \mu_{A''(\mathbb{Q}_p)}(s),
  \end{equation*}
  where $\rho''$ denotes the half-sum of positive roots for $N''$.  We have $|a''|^{\rho} = |a''|^{\rho''} |\det a''|^{-n'}$.  The required identity follows by noting, as in \S\ref{sec:basic-vectors-local}, that $|\det a''| = 1$ for $a''$ in the support of $\Theta_p^{M''}$.
\end{proof}

We now use the functions $\Theta^P_p$ to define, as follows, a linear map
\begin{equation*}
  \Theta^P : \mathcal{S}^e(N(\mathbb{R}) \backslash G(\mathbb{R})) \rightarrow C^\infty(N(\mathbb{A}) \backslash G(\mathbb{A}))
\end{equation*}
\begin{equation*}
  f \mapsto \Theta^P[f].
\end{equation*}
We first define $\Theta^P[f]_\infty \in \mathcal{S}^e(N(\mathbb{R}) \backslash G(\mathbb{R}))$ by requiring that the Mellin component $\Theta^P[f]_\infty [s] \in \mathcal{I}_\infty(s)$ be given by
\begin{equation*}
  \Theta^P[f]_\infty [s] = \mathcal{P}_G(s) \frac{\zeta(N,s)}{\zeta(N'', s)} f[s].
\end{equation*}
We emphasize that this product of $\zeta$ values, formed by Euler products over the finite primes, is being employed as a Mellin multiplier at the archimedean place.  Doing so turns out to have the effect of detecting cancellation from an otherwise large global sum (cf.\ the final paragraphs of \S\ref{sec:strong-approximation}).  We then define $\Theta^P[f] : N(\mathbb{A}) \backslash G(\mathbb{A}) \rightarrow \mathbb{C}$ by \index{Eisenstein series!modified inducing data $\Theta^P[f]$}
\begin{equation*}
  \Theta^P[f](g) := \Theta^P[f]_{\infty}(g_\infty) \prod_p \Theta^P_p (g_{p}),
\end{equation*}
noting that for given $g$, all but finitely many factors in this product are equal to $1$.  For notational convenience in what follows, we write
\begin{equation*}
  \Theta^P[f]_p := \Theta_p^P.
\end{equation*}

\begin{lemma}\label{lem:sub-gln:each-g-in}
  Let $f \in \mathcal{S}(N(\mathbb{R}) \backslash G(\mathbb{R}))$.  For each $g \in G(\mathbb{A})$ and $P$ as above, the series $\sum_{c \in A(\mathbb{Q}) } \Theta^P[f](c^{-1} g)$ converges absolutely, locally uniformly in $g$, and we have
  \begin{equation}\label{eq:sum-_c-in-1}
    \sum_{c \in A(\mathbb{Q}) } \Theta^P[f](c^{-1} g) = \Phi[f](g).
  \end{equation}
  We have
  \begin{equation}\label{eq:wpphif_inftyg-=-sum}
    W^P[\Phi[f]](g) = \sum_{c \in A(\mathbb{Q})}
    \int_{u \in N''(\mathbb{A})}
    \Theta^P[f]
    (
    c^{-1} w_{M''}^{-1} u g
    )
    \psi^{-1}(u) \, d u,
  \end{equation}
  where the RHS converges absolutely.
\end{lemma}
\begin{proof}
  To verify  the first convergence assertion, it suffices to show more precisely that the double sum
  \begin{equation}\label{eq:sum-_gamma-in-1}
    \sum_{\gamma \in {B}_{M''}(\mathbb{Q}) \backslash M''(\mathbb{Q})}
    \sum_{c \in A(\mathbb{Q})} \Theta^P[f](c^{-1} \gamma g)
  \end{equation}
  converges absolutely.  (In the special case $P = G$, we will see eventually that this double sum coincides with $\Eis[\Phi[f]]$; otherwise, it is an Eisenstein series on the smaller group $M''$ attached to the restriction of some translate of $\Phi[f]$.)  To that end, write $g = a k$ with $a \in A(\mathbb{A})$, $k \in K$.  The required convergence depends only upon the coset of $a$ modulo $A(\mathbb{Q})$.  Since $A$ is a split torus and $\mathbb{Q}$ has class number one, we have $A(\mathbb{A}) = A(\mathbb{Q}) A(\mathbb{R}) \prod_p A(\mathbb{Z}_p)$.  Since $\Theta^P[f]$ is left-invariant under $\prod_p A(\mathbb{Z}_p)$ and right-invariant under $K_p$, we may thus reduce to verifying the convergence and the identity in the special case that $a_p = 1$ and $k_p = 1$ for all $p$.  By translating $f $ on the right, we may reduce further to verifying convergence when $g_\infty = 1$, provided that we bound the series by a $G(\mathbb{R})$-continuous seminorm in $f$ so as to obtain the required local uniformity.  Our task is then to verify the convergence of
  \begin{equation*}
    \sum_{\gamma \in {B}_{M''}(\mathbb{Q}) \backslash M''(\mathbb{Q})}
    \sum_{c \in A(\mathbb{Q})}
    \Theta^P[f](c^{-1} \gamma).
  \end{equation*}
  Writing $c = c' c''$ with $(c',c'') \in A'(\mathbb{Q}) \times A''(\mathbb{Q})$, we can replace the sum over $c$ with a double sum over $c'$ and $c''$.  For $\gamma \in M''(\mathbb{Q})$, we see from Lemma \ref{lem:Theta-P-p} that $\Theta^P_p(c^{-1} \gamma) \neq 0$ only if each entry of $c'$ is a $p$-adic unit.  For this to happen for all $p$, the element $c'$ must lie in $A(\mathbb{Z})$, which is a finite group (of order $2^n$).  For the purpose of verifying convergence, we may thus drop the $c'$ sum.  We are left to verify the convergence of the remaining sum over $\gamma \in {B}_{M''}(\mathbb{Q}) \backslash M''(\mathbb{Q})$ and $c'' \in A''(\mathbb{Q})$, which may be unfolded as
  \begin{equation*}
    \sum_{\gamma \in N''(\mathbb{Q}) \backslash M''(\mathbb{Q})}
    \Theta^P[f](\gamma)
    =
    \sum_{\gamma \in N''(\mathbb{Q}) \backslash M''(\mathbb{Q})}
    \Theta^P[f]_\infty (\gamma)
    \prod_p
    \Theta^P_p(\gamma).
  \end{equation*}
  By  Lemma \ref{lem:Theta-P-p}, we have $\Theta_p^P(\gamma) = \Theta_p^{M''}(\gamma)$.  The restriction of  $\Theta^P[f]_\infty$ to $M''(\mathbb{R})$ is an element of $\mathcal{S}(N''(\mathbb{R}) \backslash M''(\mathbb{R}))$, hence decays rapidly with respect to the Iwasawa $A''(\mathbb{R})$-component.  The required convergence now follows from Lemma \ref{lem:sub-gln:let-m-in}, applied with $G$ replaced by $M''$.

  Let us now abbreviate $\Phi := \Phi[f]$.  We have seen that the series
  \begin{equation*}
    \Phi^P(g) := \sum_{c \in A(\mathbb{Q})} \Theta^P[f](c^{-1} g)
  \end{equation*}
  converges absolutely.  We verify now that $\Phi^P = \Phi$.  Let $\sigma \in \mathfrak{a}^*$ be strictly dominant.  Both $\Phi$ and $\Phi^P$ are left-invariant under $\prod_p A(\mathbb{Z}_p)$ and $\{a \in A(\mathbb{R}) : |a|^s = 1 \text{ for all } s \}$.  For all $s$, the integral
  \begin{equation*}
    \int_{a \in [A]} |a|^{\rho + s} \Phi(a^{-1} g) \, d a
  \end{equation*}
  converges absolutely and evaluates to $\Phi[s]$.  Thus by Mellin inversion for integrable functions, it will suffice to check that the integrals
  \begin{equation}\label{eq:int-_a-in}
    \int_{a \in [A]} |a|^{\rho + s} \Phi^P (a^{-1} g) \, d a
  \end{equation}
  converge absolutely and evaluate to $\Phi[s]$.

  To that end, we note the stronger assertion that
  \begin{equation}\label{eq:int-_a-in-1}
    \int_{a \in [A]} |a|^{\rho + \sigma} \sum_{c \in A(\mathbb{Q})}
    \left\lvert
      \Theta^P[f](c^{-1} a^{-1} g)
    \right\rvert
    \, d a < \infty.
  \end{equation}
  In establishing \eqref{eq:int-_a-in-1}, it suffices by strong approximation to consider the case $g \in G(\mathbb{R})$.  The above integral unfolds and then factors as a product of the archimedean integral
  \begin{equation*}
    \int_{a \in A(\mathbb{R})} |a|^{\rho + \sigma}
    \left\lvert
      \Theta^P[f](a^{-1} g)
    \right\rvert
    \, d a,
  \end{equation*}
  which converges because $\Theta^P[f] \in \mathcal{S}(U(\mathbb{R}) \backslash G(\mathbb{R}))$, and the $p$-adic integrals
  \begin{equation*}
    \int_{a \in A(\mathbb{Q}_p)} |a|^{\rho + \sigma}
    \Theta_p^P(a^{-1})
    \, d a,
  \end{equation*}
  which evaluate to the local zeta factor $\zeta_p(N,\sigma)$.  Since $\sigma$ is strictly dominant, the product of these factors converges, whence \eqref{eq:int-_a-in-1}.

  In particular, the integrals \eqref{eq:int-_a-in} converge absolutely.  By interchanging summation with integration, we obtain
  \begin{equation*}
    \int_{a \in [A]} |a|^{\rho + s}
    \Phi^P (a^{-1} g) \, d a
    =
    \int_{a \in A(\mathbb{A})}
    |a|^{\rho + s}
    \Theta^P[f] (a^{-1} g) \, d a.
  \end{equation*}
  This last integral factors as a convergent product of local integrals and evaluates to the convergent Euler product $\prod_{\mathfrak{p}} \Theta^P[f]_\mathfrak{p}[s](g_\mathfrak{p}) = \Phi[s](g)$, as required.

  It remains only to verify the identity of (partial) Whittaker functions \eqref{eq:wpphif_inftyg-=-sum}.  The absolute convergence of the double sum \eqref{eq:sum-_gamma-in-1} implies, in particular, the absolute convergence of the contribution of the open Bruhat cell in ${B}_{M''}(\mathbb{Q}) \backslash M''(\mathbb{Q})$, i.e., that of
  \begin{equation*}
    \sum_{b \in N''(\mathbb{Q})}
    \sum_{c \in A(\mathbb{Q})}
    \Theta^P[f]
    (
    c^{-1} w_{M''}^{-1} b g
    ),
  \end{equation*}
  locally uniformly in $g$, hence also the absolute convergence of the triple sum/integral
  \begin{equation*}
    I(g)
    :=
    \int_{u \in [N'']}
    \sum_{b \in N''(\mathbb{Q})}
    \sum_{c \in A(\mathbb{Q})}
    \Theta^P[f]
    (
    c^{-1} w_{M''}^{-1} b u g
    ) \psi^{-1}(u) \, d u
  \end{equation*}
  upon recalling that $[N'']$ is compact.  By evaluating the sum over $c$ using the identity \eqref{eq:sum-_c-in-1} established above and unfolding the integral/sum over $u$ and $b$, we see that $I(g) = W^P[\Phi[f]](g)$.  If we instead rearrange $I(g)$ as $\sum_{A(\mathbb{Q}) } \int_{[N'']} \sum_{N''(\mathbb{Q})}$ and unfold, we obtain the RHS of \eqref{eq:wpphif_inftyg-=-sum}.  The required identity \eqref{eq:wpphif_inftyg-=-sum} follows.
\end{proof}

\subsection{Strong approximation}\label{sec:strong-approximation}
We now use the series expressions derived in \S\ref{sec:seri-repr-via} to partially evaluate the integrals appearing in Proposition \ref{prop:bound-after-theta-degenerate-case-2}.  We denote in what follows by
\begin{equation*}
  K'  = \prod_\mathfrak{p} K_\mathfrak{p} ' \leq M'(\mathbb{A}), \quad K'' = \prod_\mathfrak{p} K_\mathfrak{p} '' \leq M''
\end{equation*}
the standard maximal compact subgroups, defined like $K =\prod_\mathfrak{p} K_\mathfrak{p} \leq G(\mathbb{A})$.
\begin{lemma}\label{lem:sub-gln:let-phi-in}
  Let $\phi \in C_c^\infty(V(\mathbb{A}))$ be nonnegative-valued and $K$-invariant.  Let $f \in \mathcal{S}(N(\mathbb{R}) \backslash G(\mathbb{R}))$.  For $(c,g) \in A(\mathbb{Q}) \times G(\mathbb{A})$, we introduce the temporary notation
  \begin{equation*}
    W(c,g) := \int_{u \in N''(\mathbb{A})} \Theta^P[f](c^{-1} w_{M''}^{-1} u g) \psi^{-1} (u) \, d u.
  \end{equation*}
  We have
  \begin{equation}\label{eq:mathc-f_infty-leq}
    \mathcal{Q}_P(\phi,f)
    \leq
    2^{n'}
    \int_{
      \substack{
        g \in N(\mathbb{R}) \backslash G(\mathbb{R}) \\
        a'' \in A''(\mathbb{A}^\infty)
      }
    }
    \phi(e_n g a'')
    \left\lvert
      \sum_{c \in A''(\mathbb{Q})}
      W(c,g a'')
    \right\rvert^2
    \, d g \frac{d a''}{\delta_N(a'')}.
  \end{equation}
\end{lemma}
\begin{proof}
  Abbreviate $\mathcal{W} := W^P[\Phi[f]]$, thus $\mathcal{W} : G(\mathbb{A}) \rightarrow \mathbb{C}$.  Note that $\phi(e_n g)$ depends only upon the bottom row of $g$.  Using the Iwasawa decompositions for $G$ and $M'$ and a careful computation of Haar measures, we may write the integral $I := \mathcal{Q}_P(\phi,f)$ to be estimated as
  \begin{equation*}
    I = \int_{a'' \in A''(\mathbb{A})}
    \phi(e_n a'')
    \int_{k \in K}
    \int_{a' \in [A']}
    \int_{k' \in K'}
    |\mathcal{W}(a' k' a'' k)|^2
    \, d k' \, \frac{d a'}{\delta_N(a')}
    \, d k \, \frac{d a''}{\delta_N(a'')}.
  \end{equation*}
  The integrand is invariant in the $a'$-variable under $\prod_p A'(\mathbb{Z}_p)$.  By strong approximation and the fact that $\mathbb{Q}$ has class number one, we can replace the integral over $[A']$ with one over $A'(\mathbb{Z}) \backslash A'(\mathbb{R})$.  Next, we observe that the integrand is invariant in the $k'$-variable under $\prod_p K_p'$ and in the $k$-variable under $\prod_p K_p$.  We may thus restrict the $k'$ and $k$ integrals to $K'_\infty$ and $K_\infty$, respectively.  Having restricted the integration in this way, we now open the series expression for $\mathcal{W}$ given by Lemma \ref{lem:sub-gln:each-g-in}, which reads
  \begin{equation*}
    \mathcal{W}(g) = \sum_{c \in A(\mathbb{Q})} W(c,g).
  \end{equation*}
  Specializing to $g = a' k' a'' k$ with $(a',k',a'',k) \in A'(\mathbb{Z}) \backslash A'(\mathbb{R}) \times K'_\infty \times A''(\mathbb{A}) \times K_\infty$ gives
  \begin{equation*}
    \Theta^P[f](c^{-1} w_{M''}^{-1} u g)
    =
    \Theta^P[f]_\infty(c^{-1} w_{M''}^{-1} u_\infty a' k' a''_\infty k)
    \prod_{p}
    \Theta^P_p(c^{-1} w_{M''}^{-1} u_p a''_p).
  \end{equation*}
  We further factor $c = c' c''$, with $(c', c'') \in A'(\mathbb{Q}) \times A''(\mathbb{Q})$.  For orientation, we note that the arguments of $\Theta^P[f]_\infty$ and $\Theta_p^P$ may be visualized as the block-diagonal matrices
  \begin{equation*}
    c^{-1} w_{M''}^{-1} u_\infty a' k' a''_\infty k
    =
    \begin{pmatrix}
      (c')^{-1} a' k' & 0 \\
      0 & (c'')^{-1} w_{M''}^{-1} u_\infty a''_\infty
    \end{pmatrix} k
  \end{equation*}
  and
  \begin{equation*}
    c^{-1} w_{M''}^{-1} u_p a_p''
    =
    \begin{pmatrix}
      (c')^{-1} & 0 \\
      0 & (c'')^{-1} w_{M''}^{-1} u_p a_p''
    \end{pmatrix}.
  \end{equation*}
  From this and \eqref{eq:thet-=-1_am}, we see that the product $\prod_p \Theta_p^P(c^{-1} w_{M''}^{-1} u_p a_p'')$ is nonzero only if $c' \in A'(\mathbb{Z})$, in which case that product is independent of the precise value of $c'$.  The set $A'(\mathbb{Z})$ has $2^{n'} = \O(1)$ elements.  By our assumption that the function $\phi$ is nonnegative and Cauchy--Schwarz, we obtain
  \begin{equation*}
    I \leq
    2^{n'}
    \int_{
      \substack{
        a'' \in A''(\mathbb{A})  \\
        k \in K_\infty  \\
        a' \in A'(\mathbb{Z}) \backslash A'(\mathbb{R}) \\
        k' \in K'_\infty
      }
    }
    \phi(e_n a'')
    \sum_{c' \in A'(\mathbb{Z})}
    \left\lvert
      \sum_{c'' \in A''(\mathbb{Q})}
      W(c' c'', a' k' a'' k)
    \right\rvert^2
    \, d k' \, d k \, \frac{d a' \, d a''}{\delta_N(a' a'')}.
  \end{equation*}
  We may further unfold the integral over $a'$ with the sum over $c'$, giving
  \begin{equation*}
    I
    \leq
    2^{n'}
    \int_{
      \substack{
        a'' \in A''(\mathbb{A})  \\
        k \in K_\infty  \\
        a' \in A'(\mathbb{R}) \\
        k' \in K'_\infty
      }
    }
    \phi(e_n a'')
    \left\lvert
      \sum_{c \in A''(\mathbb{Q})}
      W(c,a' k' a'' k)
    \right\rvert^2
    \, d k' \, d k \, \frac{d a' \, d a''}{\delta_N(a' a'')}.
  \end{equation*}
  Finally, we factor $A''(\mathbb{A}) = A''(\mathbb{R}) \times A''(\mathbb{A}^\infty)$ and recombine the archimedean integrals, giving the required estimate.
\end{proof}

With Lemma \ref{lem:sub-gln:let-phi-in}, the stage is now set for the final steps of the proof.  The integral $W(c,g)$ appearing in that lemma is Eulerian, so opening the square on the RHS of \eqref{eq:mathc-f_infty-leq} yields a sum, indexed by (say) $b,c \in A''(\mathbb{Q})$, of Eulerian integrals:
\begin{equation*}
  \int_{
    \substack{
      g \in N(\mathbb{R}) \backslash G(\mathbb{R}) \\
      a'' \in A''(\mathbb{A}^\infty)
    }
  }
  \phi(e_n a'')
  W(b, g a'')
  \overline{W(c, g a'')}
  \, d g \frac{d a''}{\delta_N(a'')}.
\end{equation*}
The estimation of such integrals is a local problem, which we address below in \S\ref{sec:local-estimates}.

The non-archimedean local integrals will turn out to detect the diagonal condition $b = c$, up to essentially bounded multiplicity, and we will not need to detect cancellation coming from the global sum over $b$ and $c$.  This pleasant feature is the motivation for why, for the analysis of $\mathcal{Q}_P(\phi, f)$, we represented $\Phi[f]$ by the series \eqref{eq:sum-_c-in-1} involving $\Theta^P[f]$, whose definition was carefully tailored to $P$.  For instance, had we instead used the series involving $\Theta^G[f]$ in the study of $\mathcal{Q}_P(\phi,f)$ for some $P \neq G$, we would be facing additional sums over $A'(\mathbb{Q})$ from which it would be necessary to extract cancellation.  It would still be possible to do so by further careful analysis, but in effect, that cancellation is ``baked in'' to our argument through the careful choice of series representation for $\Phi[f]$.

By Cauchy--Schwarz and the positivity of $\phi$, we may thus reduce to estimating the archimedean integrals in the case $b=c$.  That estimation is the most delicate part of the argument, requiring three different approaches depending upon the relative sizes of $T$, $r$ and $c$.

After completing the local analysis, we resume the global thread in \S\ref{sec:completion-proof-1} and complete the proof of Proposition \ref{prop:bound-after-theta-degenerate-case-2}, hence of our main results.

\section{Local estimates}\label{sec:local-estimates}

\subsection{Whittaker functions attached to Schwartz functions}\label{sec:whitt-funct-attach}

Let $F$ be a local field.  Set $G := \GL_n(F)$, with $n$ fixed.  As usual, let $A$ denote the diagonal subgroup, $N$ the subgroup of upper-triangular unipotent matrices, $K$ the standard maximal compact subgroup, and $Z$ the center.  We take for $\psi$ either the standard nondegenerate unitary character of $N$ or its inverse (\S\ref{eq:psin-=-psi_f}).

For $f \in C^\infty(N \backslash G)$ and $(c,g) \in A \times G$, we define
\begin{equation*}
  W(f,c,g) := \delta_N^{1/2}(c) \int_{u \in N} f(c^{-1} w^{-1}_G u g) \psi^{-1}(u) \, d u,
\end{equation*}
provided that it converges absolutely.  (We recall that $w_G^{-1} = w_G$, so the reader can ignore the inverses if desired.)

\begin{lemma}\label{lem:standard2:let-f-in}
  Let $f \in C^\infty(N \backslash G)$.  Suppose that for some strictly dominant element $\omega \in \mathfrak{a}^*$, there exists $C \geq 0$ so that for all $(a,k) \in A \times K$,
  \begin{equation}\label{eq:fa-k-leq}
    |f(a k)| \leq C |a|^{\rho + \omega}.
  \end{equation}
  Then the integral defining $W(f,c,g)$ converges absolutely  for all $(c,g)$.
\end{lemma}
\begin{proof}
  This reduces to fact that, if $f(a k) = |a|^{\rho + \omega}$ for some such $\omega$, then
  \begin{equation*}
    \int_{N} f(w^{-1}_G n) \, d n< \infty
  \end{equation*}
  (see \cite[\S10.1.2]{MR1170566}).
\end{proof}
We will primarily apply the above definition to elements $f$ of the Schwartz space $\mathcal{S}(N \backslash G)$, for which the hypotheses of Lemma \ref{lem:standard2:let-f-in} are valid (for all $\omega$, not necessarily strictly dominant).

\subsection{Bounds at finite places}\label{sec:cnjeonep68}
We retain the notation of \S\ref{sec:whitt-funct-attach}, but assume here that $F$ is a non-archimedean local field.

We pause to introduce some notation.  Let $\mu, \lambda \in \mathbb{Z}^n$.  We may identify these with weights for the diagonal torus in the complex algebraic group $\GL_n(\mathbb{C})$.  Assume that $\lambda$ is \emph{dominant} in the sense that $\lambda_1 \geq \dotsb \geq \lambda_n$.  It corresponds then to an irreducible algebraic representation $V_\lambda$, of $\GL_n(\mathbb{C})$, of highest weight $\lambda$.  We denote by $\mathfrak{M}_{\lambda}(\mu)$ the multiplicity of the weight $\mu$ in the representation $V_\lambda$.  For example, $\mathfrak{M}_{\lambda}(\lambda) = 1$, and more generally $\mathfrak{M}_{\lambda}(w \lambda) = 1$ for all permutations $w \in S(n)$.  We extend this definition to all pairs $(\mu,\lambda)$ by adopting the convention that $\mathfrak{M}_{\lambda}(\mu) := 0$ if $\lambda$ is not dominant.  We require the following standard properties.
\begin{enumerate}[(i)]
\item Suppose $\mathfrak{M}_{\lambda}(\mu) \neq 0$.  Then $\mu$ lies in the convex hull of the permutations of $\lambda$.  In particular,
  \begin{equation}\label{eq:mu_1-+-dotsb}
    \mu_1 + \dotsb + \mu_n = \lambda_1 + \dotsb + \lambda_n,
  \end{equation}
  and if $\lambda_n \geq 0$, then
  \begin{equation}\label{eq:0-leq-mu_1}
    0 \leq \mu_1,\dotsc,\mu_n \leq \lambda_1 + \dotsb + \lambda_n.
  \end{equation}
\item We have
  \begin{equation}\label{eq:mathfr-ll-1}
    \mathfrak{M}_{\lambda}(\mu) \ll (1 + |\lambda|)^{\O(1)},
  \end{equation}
  where $|\lambda|$ denotes the Euclidean norm and the implied constant depends at most upon $n$.  Such an estimate is classical, and can be deduced (in sharper forms) from, e.g., the Kostant partition formula \cite[p258, Theorem IV.3.2]{MR1410059}.
\end{enumerate}

Let $\ord : A \rightarrow \mathbb{Z}^n$ denote the map assigning to a diagonal matrix
\begin{equation*}
  a = \diag(a_1,\dotsc,a_n)
\end{equation*}
the corresponding $n$-tuple
\begin{equation*}
  \ord(a) = (\ord(a_1),\dotsc,\ord(a_n)),
\end{equation*}
where we write $\ord(y) := k$ if $y = \varpi^{k} u$ with $(k, u) \in \mathbb{Z} \times \mathfrak{o}^\times$.

We recall the formula for the normalized spherical Whittaker functions, due to Shintani \cite{MR407208} and Casselman--Shalika \cite[Thm 5.4]{MR581582}.
\begin{lemma}\label{lem:standard:spherical-whittaker-schur}
  Let $G=\GL_n(F)$, with $N$ the upper-triangular unipotent subgroup and $A$ the diagonal torus.  Let $\psi$ be unramified and let $\chi=\chi_s : A \rightarrow \mathbb{C}^\times$ be the unramified character $\chi_s(\diag(a_1,\dotsc,a_n))=\prod_i |a_i|^{s_i}$.  Let $W_s^0$ denote the corresponding normalized spherical Whittaker function, normalized by $W_s^0(1)=1$.  Then for $m \in A$ we have
  \begin{equation*}
    W_s^0(m)=\delta_N^{1/2}(m)\,S_{\ord(m)}(q^{-s_1},\dotsc,q^{-s_n}).
  \end{equation*}
\end{lemma}
\begin{proof}
  This is the Casselman--Shalika formula (after normalizing by $W_s^0(1)=1$), together with the identification of the resulting character with the Schur polynomial in the case of $\GL_n$ (cf.\ Shintani \cite{MR407208}).
\end{proof}

Let $\Theta \in C^\infty(N \backslash G)$ denote the basic vector, as in \S\ref{sec:basic-vectors-local}.  It is easy to see from the definition of $\Theta$ via a Mellin integral that it satisfies the estimate \eqref{eq:fa-k-leq}, so that $W(\Theta,c,g)$ is defined for $(c,g) \in A \times G$.

\begin{lemma}\label{lem:sub-gln:we-have-begin}
  For $c, m \in A$, we have
  \begin{equation*}
    W(\Theta,c,m) = \delta_N^{1/2}(m) \mathfrak{M}_{\ord(m)}(\ord(c)).
  \end{equation*}
\end{lemma}
\begin{proof}
  Let $\sigma \in \mathfrak{a}^*$ be strictly dominant.  We then have the convergent Mellin expansion $\Theta = \int_{\Re(s) = \sigma} \Theta[s] \, d \mu_A(s)$, where the Mellin components $\Theta[s] \in \mathcal{I}(s)$ are the spherical vectors taking the value $\zeta_F(N,s)$ at the identity.  In particular, those components satisfy the transformation property
  \begin{equation}\label{eq:f0sc-1-w_g}
    \Theta[s](c^{-1} w^{-1}_G u m)
    =
    \delta_N^{-1/2}(c) |c|^{- s} \Theta[s](w^{-1}_G u m).
  \end{equation}
  By the standard convergence lemma for Jacquet integrals \cite[\S10.1.2]{MR1170566} and the compactness of $\mathfrak{X}^e(A)$, the double integral
  \begin{equation*}
    W(\Theta, c, m) = \delta_N^{1/2}(c) \int_{u \in N} \psi^{-1}(u) \int_{\Re(s) = \sigma } \Theta[s](c^{-1} w^{-1}_G u m) \, d \mu_A(s) \, d u
  \end{equation*}
  converges absolutely, hence rearranges to the identity
  \begin{equation*}
    W(\Theta,c,m) =
    \int_{\Re(s) = \sigma }
    |c|^{-s}
    W_s^0(m)
    \, d s,
  \end{equation*}
  where we set
  \begin{equation*}
    W_s^0(m) := \int_{u \in N} \Theta[s](w^{-1}_G u m) \psi^{-1}(u) \, d u.
  \end{equation*}
  Using the Casselman--Shalika formula \eqref{eq:wf0chi-psi1-casselman-shalika} and the definition of the basic vector, we see that $W_s^0$ satisfies the normalization $W_s^0(1) = 1$.  Thus $W_s^0$ is the normalized spherical Whittaker function for the parameter $s$.  Moreover (see Example~\ref{example:cuhi8y36c4}), we have
  \begin{equation*}
    W_s^0(m) = \delta_N^{1/2}(m) S_{\ord(m)}(q^{-s_1}, \dotsc, q^{-s_n}),
  \end{equation*}
  where for $\lambda \in \mathbb{Z}^n$, we define $S_\lambda := 0$ unless $\lambda$ is dominant, in which case $S_\lambda$ is the corresponding Schur polynomial, given for $z = (z_1,\dotsc,z_n) \in (\mathbb{C}^\times)^n$ by
  \begin{equation*}
    S_\lambda(z) = \sum_{\mu \in \mathbb{Z}^n} \mathfrak{M}_{\lambda}(\mu) \prod_j z_j^{\mu_j}.
  \end{equation*}
  Thus
  \begin{equation*}
    W(\Theta,c,m) =
    \delta_N^{1/2}(m)
    \sum_{\mu \in \mathbb{Z}^n}
    \mathfrak{M}_{\ord(m)}(\mu)
    \int_{\Re(s) = \sigma } |c|^{-s} \prod_j q^{-s_j \mu_j} \, d s.
  \end{equation*}
  Writing $|c|^{-s} = \prod_j q^{s_j \ord(c_j)}$ and appealing to Mellin inversion, we see that the inner integral detects the condition $\mu = \ord(c)$.  The claimed formula follows.
\end{proof}

\begin{lemma}\label{lem:sub-gln:b-c-in}
  For $b, c \in A$, define
  \begin{equation*}
    I(b,c) :=
    \int_{g \in N(F) \backslash G(F)}
    1_{\|e_n g\| = 1}
    W  (\Theta, b, g) \overline{W }(\Theta, c, g)
    \, d g.
  \end{equation*}
  We have $I(b,c) = 0$ unless each entry $b_i,c_j$ is integral and $\sum_i \ord(b_i) = \sum_j \ord(c_j)$, in which case
  \begin{equation*}
    I(b,c) \ll (1 + \sum_j \ord(c_j))^{\O(1)}.
  \end{equation*}
\end{lemma}
\begin{proof}
  With suitable measure normalizations, we have
  \begin{equation*}
    I(b,c) =
    \int_{a \in A(F)}
    1_{|a_n| = 1}
    W  (\Theta, b, a) \overline{W }(\Theta, c, a)
    \, \frac{d^\times a}{\delta_N(a)}.
  \end{equation*}
  By Lemma \ref{lem:sub-gln:we-have-begin}, we have
  \begin{align*}
    I(b,c) &=
             \int_{a \in A(F)}
             1 _{|a_n| = 1}
             \mathfrak{M}_{\ord(a)} (\ord (b)) \mathfrak{M} _{\ord (a)} (\ord (c)) \, d a \\
           &=
             \sum_{
             \substack{
             \lambda  \in \mathbb{Z}^n:  \\
    \text{dominant}, \lambda_n = 0
    }
    }
    \mathfrak{M}_{\lambda }(\ord(b))
    \mathfrak{M}_{\lambda }(\ord(c)).
  \end{align*}
  The vanishing condition on $I(b,c)$ then follows from \eqref{eq:mu_1-+-dotsb} and \eqref{eq:0-leq-mu_1}, while the bound follows from \eqref{eq:mathfr-ll-1}.
\end{proof}

\subsection{Bounds at the real place: Whittaker functions}\label{sec:bounds-at-real}
We retain the notation of \S\ref{sec:whitt-funct-attach}, but assume here that $F = \mathbb{R}$.

\begin{lemma}\label{lem:integrate-by-parts-archimedean-whittaker-average-bound}
  There is a fixed $d_1 \in \mathbb{Z}$ so that for each $f \in \mathcal{S}^e(N \backslash G)$, $c \in A$ and $g \in G$,
  \begin{equation*}
    \int_{h \in N_H \backslash H}
    \int_{z \in Z}|W(f,c, h z g)|^2  \, d h \, d z \ll  \mathcal{S}_{d_1,0}(f)^2.
  \end{equation*}
\end{lemma}
\begin{proof}
  By replacing $f$ with $R(g) f$ and noting that $\mathcal{S}_{d_1,0}(R(g) f) = \mathcal{S}_{d_1,0}(f)$, we may reduce to the case $g = 1$.  We consider the Mellin expansion $f = \int_s f[s] \, d \mu_A(s)$, where the Mellin components satisfy the same transformation property \eqref{eq:f0sc-1-w_g} as in the non-archimedean case.  Then
  \begin{equation*}
    W(f,c,g) =
    \int_{\Re(s) = \sigma}
    |c|^{-s} W_s(f,g) \, d \mu_A(s),
  \end{equation*}
  where $W_s(f,g) = \int_{u \in N} f[s](w^{-1}_G u g) \psi^{-1}(u) \, d u$.  For $s \in \mathfrak{a}_{\mathbb{C}}^*$, write $\chi_s$ for the corresponding character of the symmetric algebra $\mathfrak{U}(A)$.  For $x \in \mathfrak{U}(A)$, we have
  \begin{equation*}
    (L(x) f)_s = \chi_s(x) f_s.
  \end{equation*}
  Thus, provided that $x$ is chosen so that $\chi_s(x) \neq 0$ whenever $\Re(s) = \sigma$, we have
  \begin{equation}\label{eq:wf-c-g}
    W(f,c,g) =
    \int_{\Re(s) = \sigma}
    |c|^{-s} \chi_s(x)^{-1} W_s(L(x) f,g) \, d \mu_A(s).
  \end{equation}
  Take $\sigma = 0$.  By the Parseval relation for the group $Z$, we have
  \begin{equation*}
    \int_{z \in Z} |W(f, c, g z)|^2 \, d z
    =
    \int_{t \in i \mathbb{R} }
    \left\lvert
      \int_{
        \substack{
          \Re(s) = 0,  \\
          \trace(s) = t
        }
      }
      |c|^{-s} \chi_s(x)^{-1} W_s(L(x) f, g ) \, d \mu_A^0(s)
    \right\rvert^2
    \, d t,
  \end{equation*}
  where $d \mu_A^0$ denotes a suitable Haar measure, regarding the integration domain as the unramified subset of a coset of $(A/Z)^\wedge$.

  Fix $m \in \mathbb{Z}_{\geq 0}$ large enough that
  \begin{equation*}
    \int_{
      \substack{
        \Re(s) = 0,  \\
        \trace(s) = t
      }
    }
    (1 + |s|^2)^{-2 m} \, d \mu_A^0(s) \ll 1
  \end{equation*}
  for all $t \in i \mathbb{R}$.  We may fix $x$ of degree $2 m$ so that $\chi_s(x) = (1 + |s|^2)^m$ whenever $\Re(s) = 0$ (indeed, take $x = (1 - \sum_{j} x_j^2)^{m}$, where $x_j$ is chosen so that $\chi_s(x_j) = s_j$ with the normalization $\lvert s \rvert^2 = \sum_j \lvert s_j \rvert^2 = -\sum_j s_j^2$ when each $\Re(s_j) = 0$).  Then, by Cauchy--Schwarz,
  \begin{equation*}
    \int_{z \in Z}
    |W(f,c,g z)|^2 \, d z \ll
    \int_{\Re(s) = 0}
    \left\lvert W_s(L(x) f, g)  \right\rvert^2 \, d \mu_A(s).
  \end{equation*}
  We now integrate over $h \in N_H \backslash H$.  It is known \cite[Appendix A]{MR2930996} that for each $s$ with $\Re(s) = 0$, we have
  \begin{equation}\label{eq:int-_g-in-norm}
    \int_{h \in N_H \backslash H} |W_s(L(x) f, h)|^2 \, d h
    \asymp
    \int_{K} |L(x) f_s|^2,
  \end{equation}
  where the implied constants depend only upon measure normalization.  Finally, the Parseval relation on $A$ gives
  \begin{equation*}
    \int_{\Re(s) = 0} \int_K |L(x) f_s|^2 \, d \mu_A(s) \asymp \int_{N \backslash G} |L(x) f|^2
    \ll \mathcal{S}_{2 m,0}(f)^2.
  \end{equation*}
  We conclude by taking $d_1 := 2 m$.
\end{proof}

\begin{lemma}\label{lem:sub-gln:let-f-in-1}
  There is a fixed $d_1 \in \mathbb{Z}_{\geq 0}$ with the following property.  Let $D > 0$.  Let $f \in \mathcal{S}(N \backslash G)$ be such that $f(g) \neq 0$ only if $|\det(g)| \asymp D$.  Let $c \in A$.  Then the integral
  \begin{equation*}
    I(c) :=
    \int_{g \in N \backslash G}
    \|e_n g\|^{-n}
    |W(f, c, g)|^2  \, d g
  \end{equation*}
  satisfies the bound
  \begin{equation}\label{eq:ic-ll-detc}
    I(c) \ll \frac{1}{D |\det(c)|} \mathcal{S}_{d_1,0}(f)^2.
  \end{equation}
\end{lemma}
\begin{proof}
  We rewrite
  \begin{equation*}
    I(c) =
    \int_{a \in A}
    \int_{k \in K}
    |a_n|^{-n}
    |W(f, c, a k)|^2 \, d k \, \frac{d a}{\delta_N(a)}.
  \end{equation*}
  We appeal to the identity
  \begin{equation*}
    |a_n|^n \delta_N(a) = |\det a| \delta_{N_H}(a),
  \end{equation*}
  corresponding to the following identity for $n$-tuples of exponents:
  \begin{align*}
    &\left( 0, \dotsc, 0, n \right) + \left( n-1, \dotsc, 3-n, 1 - n \right) \\
    &\quad =
      \left( 1, \dotsc, 1, 1 \right)
      +
      \left( n-2, \dotsc, 2-n, 0 \right).
  \end{align*}
  We factor $a = a_H z$, where $a_H \in A_H$ and $z \in Z$, and then rename $a_H$ to $a$.  We also substitute $k \mapsto k_H k$, with $k_H \in K_H$, and integrate over $K_H$.  Up to measure normalizations, we obtain
  \begin{equation*}
    I(c) =
    \int_{a \in A_H}
    \int_{k_H \in K_H}
    \int_{z \in Z}
    \int_{k \in K}
    |W(f, c, a k_H z k)|^2
    |\det(a z)|^{-1}
    \, d k \, d z \, d k_H \, \frac{d a}{\delta_{N_H}(a)}.
  \end{equation*}
  We can rewrite this further, again up to measure normalizations, as
  \begin{equation*}
    I(c) =
    \int_{h \in N_H \backslash H}
    \int_{z \in Z}
    \int_{k \in K}
    |W(f, c, h z k)|^2
    |\det(h z)|^{-1}
    \, d k \, d z \, d h.
  \end{equation*}
  On the support of the integrand, we have $|\det(h z)| \asymp |\det(c)| D$.  We may write $W(f,c,h z k) = W(R(k) f,c,h z)$.  We conclude via Lemma \ref{lem:integrate-by-parts-archimedean-whittaker-average-bound}, noting that $\mathcal{S}_{d_1,0}(f) \asymp \mathcal{S}_{d_1,0}(R(k) f)$.
\end{proof}

For $a \in A$, $\eps > 0$ and $\ell \in \mathbb{Z}_{\geq 0}$, we introduce the notation
\begin{equation*}
  \mathcal{G}_{\eps,\ell}(a) := \delta_N^{1/2-\eps}(a)
  \min (\{1\} \cup \{|a|^{-\alpha} : \alpha \in \Delta \})^{\ell}.
\end{equation*}

\begin{lemma}\label{lem:sub-gln:let-ell-geq}
  Let $\ell \geq 0$.  There exists $d_0 \geq 0$, depending only upon $\ell$ and $n$, so that for each $\eps > 0$, there exists $C_0 \geq 0$ with the following property.  Let $\pi$ be a generic tempered irreducible unitary representation of $G$, realized in its Whittaker model $\mathcal{W}(\pi,\psi)$ and with unitary structure defined via the Kirillov model.  Let $W \in \mathcal{W}(\pi,\psi)$ be a smooth vector.  Let $(a,k) \in A \times K$.  Then
  \begin{equation*}
    |W(a k)| \leq C_0 \mathcal{G}_{\eps,\ell}(a) \mathcal{S}_{d_0}(W).
  \end{equation*}
\end{lemma}
\begin{proof}
  This is \cite[Lemma 5.2]{JN19a}.\footnote{The cited result imposes hypotheses that we do not record here: that $\eta$ be small enough, and $N$ large enough.  These hypotheses are clearly unnecessary.  It also includes an additional factor $\prod_{i=1}^{n} (a_i + a_i^{-1})^\eps$, which can be absorbed into $\mathcal{G}_{\eps,\ell}$ by applying the cited result to some suitable $\eps ' < \eps$ and $\ell ' > \ell$.}
\end{proof}
\begin{remark}
  Lemma \ref{lem:sub-gln:let-ell-geq} is uniform in $\pi$, while standard references, such as \cite[Prop 3.5]{MR2533003}, \cite[\S15]{MR1170566} or \cite[Prop 9]{MR2058615}, give uniformity at most for $\pi$ in a fixed compact set.  However, it seems likely that a careful quantification of the proofs of those estimates would recover the conclusion of Lemma \ref{lem:sub-gln:let-ell-geq}.
\end{remark}

\begin{lemma}\label{lem:sub-gln:let-underlinen-bzp}
  Let $\underline{\ell} = (\ell_\alpha)_{\alpha \in \Delta}$ be a $\Delta$-tuple of nonnegative integers.  Set $\ell := \sum_{\alpha \in \Delta} \ell_\alpha$ and $u_{\underline{\ell}} := \prod_{\alpha \in \Delta} X_\alpha^{\ell_\alpha} \in \mathfrak{U}(N)$, where $X_\alpha$ is a basis element for the corresponding root subspace.  There exists $C \geq 0$ so that for all $W \in C^\infty(N \backslash G, \psi)$ and $(a,k) \in A \times K$, we have
  \begin{equation*}
    |W(a k)| \leq
    C
    |a|^{-\sum_{\alpha \in \Delta} \ell_\alpha \alpha}
    \sup_{x_1,\dotsc,x_\ell \in \mathcal{B}(G) \cup \{1\}}
    |R(x_1 \dotsb x_\ell) W(a k)|.
  \end{equation*}
\end{lemma}
\begin{proof}
  This is a step in the proof of \cite[Lemma 2.5.1]{2018arXiv181200053B}.
\end{proof}
\begin{corollary}\label{cor:sub-gln:each-ell-in}
  For each $\ell \in \mathbb{Z}_{\geq 0}$, there exists $C \geq 0$ so that for all $W \in C^\infty(N \backslash G, \psi)$ and $(a,k) \in A \times K$, we have
  \begin{equation*}
    |W(a k)| \leq C
    \min (\{ 1 \} \cup \{ |a|^{-\alpha} : \alpha \in \Delta \})^{\ell}
    \sup_{x_1,\dotsc,x_\ell \in \mathcal{B}(G) \cup \{1\}}
    |R(x_1 \dotsb x_\ell) W(a k)|.
  \end{equation*}
\end{corollary}

\begin{lemma}\label{lem:scratch-research:each-fixed-m}
  There is a fixed $d_0 \geq 0$ with the following property.  Let $\ell \in \mathbb{Z}_{\geq 0}$ be fixed.  Let $f \in \mathcal{S}(N \backslash G)$.  Let $c \in A$.  Let $k$ belong to a fixed compact subset of $G$.  Let $\eps > 0$ be fixed.  Then
  \begin{equation*}
    W(f,c,a k) \ll \mathcal{G}_{\eps,\ell}(a) \mathcal{S}_{d_0, d_0 + \ell}(f).
  \end{equation*}
\end{lemma}
\begin{proof}
  We let $x \in \mathfrak{U}(A)$ be as in the proof of Lemma \ref{lem:integrate-by-parts-archimedean-whittaker-average-bound}, and appeal once again to the Mellin expansion \eqref{eq:wf-c-g}:
  \begin{equation*}
    W(f,c,a k) =
    \int_{\Re(s) = 0}
    |c|^{-s} \chi_s(x)^{-1} W_s(L(x) f,a k) \, d \mu_A^0(s).
  \end{equation*}
  We let $d_0$ be as in Lemma \ref{lem:sub-gln:let-ell-geq} applied with $\ell = 0$, so that for each fixed $\eps > 0$,
  \begin{equation*}
    W_s(L(x) f, a k) \ll \mathcal{G}_{\eps,0}(a) \mathcal{S}_{d_0}(W_s(L(x) f, \cdot)).
  \end{equation*}
  Using the norm relation \eqref{eq:int-_g-in-norm}, we deduce that
  \begin{equation*}
    W_s(L(x) f, a k) \ll \mathcal{G}_{\eps,0}(a) \mathcal{S}_{d_0}(L(x) f_s).
  \end{equation*}
  Applying Cauchy--Schwarz and appealing to the Parseval relation as in the proof of Lemma \ref{lem:integrate-by-parts-archimedean-whittaker-average-bound}, we obtain
  \begin{equation}\label{eq:mathcalg_eps-0a-1}
    |\mathcal{G}_{\eps,0}(a)^{-1} W(f,c,a k)|^2 \ll
    \int_{\Re(s) = 0}
    \mathcal{S}_{d_0}(L(x) f_s)^2 \, d \mu_A^0(s)
    \asymp
    \mathcal{S}_{d_0}(L(x) f)^2.
  \end{equation}
  By increasing $d_0$ as necessary, we may assume that it exceeds the degree of $x$, so that
  \begin{equation*}
    \mathcal{S}_{d_0}(L(x) f)^2 \ll \mathcal{S}_{d_0,d_0}(f)^2.
  \end{equation*}

  To improve the estimate obtained thus far, we invoke Corollary \ref{cor:sub-gln:each-ell-in} and apply \eqref{eq:mathcalg_eps-0a-1} with $f$ replaced by $R(x_1 \dotsb x_{\ell}) f$ (noting that $f \mapsto W(f,c,\cdot)$ is $G$-equivariant).  We obtain
  \begin{equation*}
    |\mathcal{G}_{\eps,\ell}(a)^{-1} W(f,c,a k)|^2 \ll
    \sup_{x_1,\dotsc,x_\ell \in \mathcal{B}(G) \cup \{1\}}
    \mathcal{S}_{d_0,d_0}(R(x_1 \dotsb x_{\ell}) f)^2
    \ll
    \mathcal{S}_{d_0, d_0 + \ell}(f)^2,
  \end{equation*}
  as required.
\end{proof}

Let $T \geq 1$.  For $a \in A$, we set
\begin{equation*}
  \mathcal{G}_T(a) :=
  \min (\{ 1 \} \cup \{ T |a|^{-\alpha} : \alpha \in \Delta \}).
\end{equation*}
\begin{lemma}\label{lem:sub-gln:there-fixed-d_0}
  There is a fixed $d_0 \geq 0$ so that for all fixed $\ell \geq 0$, we have for all $f \in \mathcal{S}(N \backslash G)$ and $(c,a,k) \in A \times  A \times K$ the estimate
  \begin{equation}\label{eq:wf-c-k}
    W(f,c,a k) \ll
    T^{d_0} \mathcal{G}_T(a)^{\ell} \mathcal{S}_{d_0, d_0 + \ell, T}(f).
  \end{equation}

\end{lemma}
\begin{proof}
  We apply Lemma \ref{lem:scratch-research:each-fixed-m} with $\eps = 1/2$.  Let $d_0$ be as in the conclusion of that lemma.  Let $\ell \in \mathbb{Z}_{\geq 0}$ be fixed.  It is immediate from the definition \eqref{eq:mathcals_d_1-d_2-tf} that
  \begin{equation*}
    \mathcal{S}_{d_0, d_0 + \ell}(f) \leq T^{d_0 + \ell} \mathcal{S}_{d_0, d_0 + \ell, T}(f),
  \end{equation*}
  so the conclusion of that lemma gives
  \begin{equation}\label{eq:wf-c-k-1}
    W(f,c,a k) \ll
    T^{d_0 + \ell}
    \min (\{1\} \cup \{|a|^{-\alpha} : \alpha \in \Delta \})^{\ell} \mathcal{S}_{d_0, d_0 + \ell, T}(f).
  \end{equation}
  This estimate holds in particular with $\ell = 0$, in which case it specializes to
  \begin{equation}\label{eq:wf-c-k-2}
    W(f,c,a k) \ll
    T^{d_0}
    \mathcal{S}_{d_0,d_0,T}(f).
  \end{equation}
  We now distinguish two cases:
  \begin{itemize}
  \item If $\mathcal{G}_T(a) = 1$, then the desired bound \eqref{eq:wf-c-k} follows from \eqref{eq:wf-c-k-2}.
  \item If $\mathcal{G}_T(a) < 1$, then $\mathcal{G}_T(a) = T |a|^{-\beta}$ for some $\beta \in \Delta$.  Thus
    \begin{equation*}
      T
      \min
      (\{1\} \cup \{|a|^{-\alpha} : \alpha \in \Delta \})
      \leq T |a|^{-\beta} = \mathcal{G}_T(a),
    \end{equation*}
    and so \eqref{eq:wf-c-k} follows from \eqref{eq:wf-c-k-1}.
  \end{itemize}
\end{proof}

The function $\mathcal{G}_T$ is simplest when $n=1$: in that case, $\mathcal{G}_T(a) = 1$ for all $a$.  For this reason, we restrict to $n \geq 2$ in the following results.
\begin{lemma}\label{lem:sub-gln:let-a-in}
  Assume $n \geq 2$.  Let $a \in A$.  Then
  \begin{equation*}
    \mathcal{G}_T(a) \leq
    \min \left( 1,
      \frac{|a_n^n| T^{n(n-1)/2}}{ |\det a|}
    \right)^{ \frac{2}{n(n-1)} }.
  \end{equation*}
\end{lemma}
\begin{proof}
  We have $\mathcal{G}_T(a) = \min(1,T/X)$, where
  \begin{equation*}
    X := \max \{|a|^\alpha : \alpha \in \Delta \}
    = \max \left( \left\lvert \frac{a_1}{a_2} \right\rvert, \dotsc, \left\lvert \frac{a_{n-1}}{a_{n}} \right\rvert \right).
  \end{equation*}
  Clearly $|a_i / a_{i+1}| \leq X$ for each $i$.  Iterating, we obtain
  \begin{equation*}
    |a_{n-1}| \leq X |a_n|,
    \quad
    |a_{n-2}| \leq X^2 |a_{n}|,
    \quad
    \dotsc,
    \quad
    |a_1| \leq X^{n-1} |a_n|.
  \end{equation*}
  Multiplying these inequalities together gives
  \begin{equation*}
    |\det a| \leq X^{n(n-1)/2} |a_n^n|.
  \end{equation*}
  Therefore
  \begin{equation*}
    \mathcal{G}_T(a)^{n(n-1)/2}
    \leq \left( \frac{T}{X} \right)^{n(n-1)/2}
    \leq
    \frac{|a_n^n| T^{n(n-1)/2} }{|\det a|}.
  \end{equation*}
  We conclude by taking the $2/n(n-1)$st power of both sides.
\end{proof}

\begin{lemma}\label{lem:sub-gln:there-fixed-d_0-1}
  Assume $n \geq 2$.  There is a fixed $d_0 \geq 0$ with the following property.  Let $\ell_1, \ell_2 \geq 0$ be fixed.  Let $D > 0$.  Suppose that $f \in \mathcal{S}(N \backslash G)$ satisfies the support condition that $f(g) \neq 0$ only if $|\det(g)| \asymp D$.  Then for all $(c,a,k) \in A \times A \times K$, we have
  \begin{equation}\label{eq:wf-c-k-3}
    W(f,c,a k) \ll
    T^{d_0}
    \mathcal{G}_T(a)^{\ell_1}
    \min \left( 1,
      \frac{|a_n^n| T^{n(n-1)/2}}{D |\det c|}
    \right)^{\frac{2}{n(n-1)} \ell_2 }
    \mathcal{S}_{d_0, d_0+ \ell_1 + \ell_2, T}(f).
  \end{equation}
\end{lemma}
\begin{proof}
  We apply Lemma \ref{lem:sub-gln:there-fixed-d_0}.  To estimate some of the $\mathcal{G}_T(a)$ factors, we apply Lemma \ref{lem:sub-gln:let-a-in}, noting that our hypotheses on $f$ imply that $W(f,c,a k) \neq 0$ only if $|\det(c^{-1} a k)| \asymp D$, i.e., $|\det a| \asymp D |\det c|$.
\end{proof}

\subsection{Rankin--Selberg integrals}\label{sec:rank-selb-integr}
The definitions that follow apply to any local field $F$, but will be applied in this paper only when $F = \mathbb{R}$.

We denote by $V = F^{n} - \{0\}$ the punctured $n$-plane over $F$.  For $\phi \in C_c^\infty(V)$, $f \in \mathcal{S}(N \backslash G)$ and $c \in A$, we define
\begin{equation*}
  Q(\phi,f,c) := \int_{g \in N \backslash G}
  \phi(e_n g) |W(f,c,g)|^2 \, d g.
\end{equation*}
As noted earlier, these may be understood as ``pseudo local Rankin--Selberg integrals:'' a spectral expansion of $\phi$ and $f$ would yield an integral of ``ordinary'' local Rankin--Selberg integrals for $\GL_n \times \GL_n$.  We note that if $\phi$ is $K$-invariant, then
\begin{equation*}
  Q(\phi,f,c) =
  \int_{a \in A}
  \phi(a_n e_n)
  \int_{k \in K}
  |W(f,c,a k)|^2 \, d k \, \frac{d a}{\delta_N(a)}.
\end{equation*}

\subsection{Bounds for Rankin--Selberg integrals}\label{sec:bounds-at-real-1}
We now resume working over $F = \mathbb{R}$.  We let $T \geq 1$ and $r > 0$ be parameters.  For the following results, we let $\phi \in C_c^\infty(V)$ be an element satisfying the following conditions, motivated by the conclusion of Lemma \ref{lem:sub-gln:let-r-be}:
\begin{assumptions}\label{sec:bounds-at-real-assumptions-phi}
  We have:
  \begin{enumerate}[(i)]
  \item $0 \leq \phi \leq 1$.
  \item $\phi$ is $K$-invariant.
  \item $\phi(v) \neq 0$ only if $\|v\| \asymp r$.
  \end{enumerate}
\end{assumptions}
Thus for $g = a k$ with $(a,k) \in A \times K$, we have $\phi(e_n g) = \phi(a_n e_n)$, which vanishes unless $|a_n| \asymp r$.

In what follows, $f$ is an element of $\mathcal{S}(N \backslash G)$, $D$ is a positive real, and $c$ is an element of $A$.

\begin{lemma}\label{lem:standard:there-fixed-d_1-1}
  There is a fixed $d_1 \in \mathbb{Z}_{\geq 0}$ so that, if $f(g) \neq 0 \implies |\det g| \asymp D$, then
  \begin{equation*}
    Q(\phi,f,c) \ll \frac{r^n}{D |\det c|} \mathcal{S}_{d_1,0}(f)^2.
  \end{equation*}
\end{lemma}
\begin{proof}
  By the size and support properties of $\phi$, we have
  \begin{equation*}
    \phi(e_n g) \ll r^n \|e_n g\|^{-n}.
  \end{equation*}
  Thus $Q(\phi,f,c) \ll r^n I(c)$, where $I(c)$ is as in Lemma \ref{lem:sub-gln:let-f-in-1}.  The conclusion of that lemma gives the required estimate.
\end{proof}

\begin{lemma}\label{lem:sub-gln:there-fixed-x}
  There is a fixed $d_1 \in \mathbb{Z}_{\geq 0}$ with the following property.  Let $b \in \mathbb{Z}_{\geq 0}$ be fixed and large enough in terms of $n$.  There is a fixed $d_2 \in \mathbb{Z}_{\geq 0}$ so that
  \begin{equation*}
    Q(\phi,f,c)
    \ll
    \frac{r^{n b}}{|\det c|^b}
    \mathcal{S}_{d_1,d_2} ( |\det|^{-b/2} f)^2.
  \end{equation*}
  Here $|\det|^{-b/2} f$ denotes the function on $N \backslash G$ given by $g \mapsto |\det g|^{-b/2} f(g)$.
\end{lemma}
\begin{proof}
  We apply Lemma \ref{lem:scratch-research:each-fixed-m} with $\eps = 1/2$.  This gives, for some fixed $d_1 \in \mathbb{Z}_{\geq 0}$ and all fixed $\ell \in \mathbb{Z}_{\geq 0}$, the estimate
  \begin{equation*}
    W(f,c,a k) \ll \mathcal{G}(a)^{\ell} \mathcal{S}_{d_1, d_1+\ell}(f)
  \end{equation*}
  for all $f \in \mathcal{S}(N \backslash G)$ and $(a,k) \in A \times K$, where
  \begin{equation*}
    \mathcal{G}(a) := \min(\{1\} \cup \{|a|^{-\alpha} : \alpha \in \Delta \})
    =
    \min(1, |a_2/a_1|, \dotsc, |a_n / a_{n-1}|).
  \end{equation*}
  Fix $b \in \mathbb{Z}_{\geq 0}$.  We can apply the above estimate with $f$ replaced by $|\det|^{-b/2} f$.   Since
  \begin{equation*}
    W(f,c, a k)
    =
    |\det c^{-1} a|^{b/2}
    W(|\det|^{-b/2} f, c, a k),
  \end{equation*}
  we obtain
  \begin{equation}\label{eq:wf-c-k2}
    |W(f, c, a k)|^2 \ll
    \frac{|\det a|^{b}}{|\det c|^b}
    \mathcal{G}(a)^{2 \ell} \mathcal{S}_{d_1, d_1 + \ell} ( |\det|^{-b/2} f)^2.
  \end{equation}
  Thus
  \begin{equation*}
    Q(\phi,f,c) \ll
    \frac{J}{|\det c|^b}
    \mathcal{S}_{d_1, d_1 + \ell} ( |\det|^{-b/2} f)^2,
  \end{equation*}
  where
  \begin{equation*}
    J := \int_{a \in A}
    \phi(e_n a)
    |\det a|^b
    \mathcal{G}(a)^{2 \ell}
    \frac{d a}{\delta_N(a)}.
  \end{equation*}
  We factor $a = a_H z$ with $(a_H,z)  \in A_H \times Z$.  Thus $J$ factors as the product of integrals
  \begin{equation*}
    \left(
      \int_{z \in Z} \phi (e_n z) |\det z|^b  \, d z
    \right)
    \left(
      \int_{a \in A_H} |\det a|^b \mathcal{G}(a)^{2 \ell} \, \frac{d a}{ \delta_N(a)}
    \right).
  \end{equation*}
  Since $\phi$ is supported on $e_n z \asymp r$ and of size $\O(1)$, the first parenthetical integral is $\asymp r^{n b}$.  We verify readily (Lemma \ref{lem:sub-gln:let-r-in} below) that the second integral is finite, hence $\O(1)$, provided that $b$ (resp.\ $\ell$) is taken sufficiently large in terms of $n$ (resp.\ $b$).  Thus $J \ll r^{n b}$ and the required conclusion holds with $d_2 := d_1 + \ell$.
\end{proof}

The following calculation establishes the convergence required in the proof of the previous lemma:
\begin{lemma}\label{lem:sub-gln:let-r-in}
  Let $b \in \mathbb{Z}_{\geq 0}$ be large enough in terms of $n$, and let $\ell \in \mathbb{Z}_{\geq 0}$ be large enough in terms of $b$.  Then, with notation as in the proof above,
  \begin{equation*}
    \int_{a \in A_H} |\det a|^b \mathcal{G}(a)^{\ell} \, \frac{d a}{ \delta_N(a)} < \infty.
  \end{equation*}
\end{lemma}
\begin{proof}
  There is nothing to check if $n = 1$ (since then $A_H$ is trivial), so we may assume that $n \geq 2$.  It suffices to consider the same integral taken over the connected component $A_H^0$ of $A_H$, on which all entries are positive.  For $a \in A_H^0$, we define $y = (y_1,\dotsc,y_{n-1}) \in (\mathbb{R}^\times_+)^n$ by the formulas
  \begin{equation*}
    y_j := a_j / a_{j+1}.
  \end{equation*}
  Since $a_{n + 1} = 1$, we then have
  \begin{equation*}
    a_1 = y_1 \dotsb y_{n-1},
    \quad
    a_2 = y_2 \dotsb y_{n-1},
    \quad
    \dotsc
    \quad
    a_{n-1} = y_{n-1},
  \end{equation*}
  \begin{equation*}
    \det a = y_1 y_2^2 y_3^3 \dotsb y_{n-1}^{n-1},
  \end{equation*}
  \begin{equation*}
    \mathcal{G}(a) =
    \min(1,y_1^{-1},\dotsc,y_{n-1}^{-1})
    \leq
    \prod_{j=1}^{n-1} \max(1, y_j)^{-1/(n-1)}
  \end{equation*}
  and
  \begin{equation*}
    \delta_N(a) =  \prod_{j=1}^{n-1} y_j^{b_j}, \quad
    b_j := \sum_{k=1}^{j} (n + 1 - 2 k).
  \end{equation*}
  The integral in question then factors into a product of one-dimensional integrals:
  \begin{equation*}
    \prod_{j=1}^{n-1}
    \int_{y \in \mathbb{R}^\times_+}
    y^{b j} \max(1, y)^{-\ell/(n-1)}
    \, \frac{d^\times y}{y^{b_j}}.
  \end{equation*}
  We choose $b$ large enough that $b j > b_j$ for all $j$, so that each integral converges near zero, and then choose $\ell$ large enough that $\ell/(n-1) > b j - b_j$ for all $j$, so that each integral converges near infinity.  The proof is then complete.
\end{proof}

\begin{lemma}\label{lem:standard:there-fixed-d_1}
  Assume that $n \geq 2$.  Let $b \in \mathbb{Z}_{\geq 0}$ be fixed and large enough in terms of $n$.  There are fixed $d_1,d_2 \in \mathbb{Z}_{\geq 0}$ such that, if $f(g) \neq 0 \implies |\det g| \asymp D$, then for each fixed $\ell \in \mathbb{Z}_{\geq 0}$, we have
  \begin{equation*}
    Q(\phi,f,c) \ll
    \frac{T^{d_1} r^{n b}}{|\det c|^b}
    \min \left( 1,
      \frac{r^n T^{n(n-1)/2}}{D |\det c|}
    \right)^{\frac{4}{n(n-1)} \ell }
    \mathcal{S}_{d_1, d_2 + \ell, T}(|\det|^{-b/2} f)^2.
  \end{equation*}
\end{lemma}
\begin{proof}
  Take $d_0$ as in the conclusion of Lemma \ref{lem:sub-gln:there-fixed-d_0-1}.  By Lemma \ref{lem:sub-gln:let-r-in}, we may fix $d \in \mathbb{Z}_{\geq 0}$ large enough that
  \begin{equation*}
    \int_{a \in A_H} |\det a|^b \mathcal{G}(a)^{2 d} \, \frac{d a}{ \delta_N(a)} < \infty.
  \end{equation*}
  We apply Lemma \ref{lem:sub-gln:there-fixed-d_0-1} with $(\ell_1, \ell_2) = (d,\ell)$ and with $|\det|^{-b/2} f$ in place of $f$, giving
  \begin{align*}
    &|W(f,c,a k)|^2 \\
    &\ll
      \frac{ T^{2 d_0} |\det a|^b}{|\det c|^b}
      \mathcal{G}_T(a)^{2 d}
      \min \left( 1,
      \frac{|a_n^n| T^{n(n-1)/2}}{D |\det c|}
      \right)^{\frac{4}{n(n-1)} \ell }
      \mathcal{S}_{d_0, d_0 + d + \ell, T}(|\det |^{-b/2} f)^2.
  \end{align*}
  We now use that $\mathcal{G}_T(a) \leq T \mathcal{G}(a)$ and argue exactly as in the proof of Lemma \ref{lem:sub-gln:there-fixed-x}, from \eqref{eq:wf-c-k2} onwards.  We obtain the required estimate with $d_1 := 2 d_0 + 2 d$ (which satisfies $d_1 \geq d_0$, so we may increase the first index of the Sobolev norm without invalidating the estimate) and $d_2 := d_0 + d$.
\end{proof}

\subsection{Degenerate integrals}\label{sec:degenerate-integrals}
Let $P$ be $G$ or a maximal standard parabolic.

For convenience, we recall (from \S\ref{sec:reform-terms-const-1}, \S\ref{sec:summary-unification}, \S\ref{sec:seri-repr-via}, \S\ref{sec:strong-approximation}) the notation we have attached to such a group $P$.  We may factor $P = M {U_P}$, with $M$ a standard Levi (containing the diagonal subgroup) and ${U_P}$ the unipotent radical (consisting of upper-triangular unipotent matrices).  If $P \neq G$, then we may factor $M = M' \times M''$, where $M'$ (resp.\ $M''$) denotes the ``upper-left'' (resp.\ ``lower-right'') block of $M$.  In the case $P = G$, we adopt the convention that $M' := \{1\}$ and $M'' := G$.  We then define $A', N', K'$ and $A'', N'', K''$ to be the intersections with $M'$ and $M''$, respectively, of the corresponding subgroups $A,N,K$ of $G$.  Finally, we introduce the notation
\begin{equation*}
  n' := \rank(M'), \quad n'' := \rank(M''),
\end{equation*}
so that $n = n' + n''$.  We have $n' \in \{0,\dotsc,n-1\}$ and $n'' \in \{1, \dotsc, n\}$.

For $(f, c,g) \in \mathcal{S}(N \backslash G) \times A'' \times G$, we define the (partial) Whittaker transform
\begin{equation*}
  W_P(f,c,g) = \delta_{N}^{1/2}(c) \int_{u \in N''} f(c^{-1} w^{-1}_{M''} u g) \psi^{-1} (u) \,d u.
\end{equation*}
Such integrals converge absolutely for the same reasons as in \S\ref{sec:whitt-funct-attach}.  For $(\phi,f,c,g) \in C_c^\infty(V) \times \mathcal{S}(N \backslash G) \times A'' \times G$, we define
\begin{equation*}
  Q_P(\phi,f,c) := \int_{g \in N \backslash G}
  \phi(e_n g) |W_P(f,c,g)|^2 \, d g.
\end{equation*}
For example, $W_G(\dotsb)$ and $Q_G(\dotsb)$ are the integrals $W(\dotsb)$ and $Q(\dotsb)$ as defined in \S\ref{sec:whitt-funct-attach} and \S\ref{sec:rank-selb-integr}.

The definition of $Q$ makes sense for any general linear group $G$.  We write $W^{M''}$ and $Q^{M''}$ for the integrals defined like $W$ and $Q$, but with $M''$ playing the role of $G$.  That is, write $V'$ (resp.\ $V''$) for the punctured $n'$-plane (resp.\ $n''$-plane), included in $V$ as a ``punctured direct summand.''  For $(\phi, f, c, g) \in C_c^\infty(V'') \times \mathcal{S}(N'' \backslash M'') \times A'' \times M''$, we set
\begin{equation*}
  W^{M''}(f,c,g) := \delta_{N''}^{1/2}(c) \int_{u \in N''} f(c^{-1} w_{M''}^{-1} u g) \psi^{-1}(u) \, d u,
\end{equation*}
\begin{equation*}
  Q^{M''}(\phi,f,c) := \int_{g \in N'' \backslash M''} \phi (e_{n} g) |W^{M''}(f,c,g)|^2 \, d g.
\end{equation*}
(Here $e_n$ identifies with the $n''$th basis vector for $V''$.)

Given $\phi \in C_c^\infty(V)$ and $f \in \mathcal{S}(N \backslash G)$, the restrictions $\phi|_{V''}$ and $f|_{M''}$ lie in $C_c^\infty(V'')$ and $\mathcal{S}(N'' \backslash M'')$, respectively.  For such $\phi$ and $f$, we define $Q^{M''}(\phi,f,c)$ for $c \in A''$ by applying the above definitions to the restrictions.

We now record a general integration lemma, which we will shortly apply to evaluate $Q_P$ in terms of $Q^{M''}$.

\begin{lemma}\label{lem:standard:rewrite-N-G-via-M}
  Let $\Phi$ be an integrable function on $N \backslash G$.  Then
  \begin{equation*}
    \int_{N \backslash G} \Phi
    =
    \int_{a ' \in A'}
    \int_{k \in K}
    \int_{g \in N'' \backslash M''}
    |\det g|^{n'}
    \Phi (g a' k)
    \, d g
    \, d k \, \frac{d a '}{\delta_N(a')}.
  \end{equation*}
\end{lemma}
\begin{proof}
  We write $\int_{N \backslash G}$ out in Iwasawa coordinates as
  \begin{equation*}
    \int_{a' \in A'}
    \int_{k \in K}
    \int_{a'' \in A''}
    \int_{k'' \in K'' }
    \Phi (a'' k'' a' k)
    \, d k''
    \, \frac{d a''}{\delta_N(a'')}
    \, d k
    \, \frac{d a'}{\delta_N (a')}.
  \end{equation*}
  (We have extracted $k''$ from the $k$-integral, using that $a'$ and $k''$ commute.)  The identity $\delta_N(a'') = \delta_{N''}(a'') |\det a''|^{-n'}$ gives
  \begin{equation*}
    \int_{a' \in A'}
    \int_{k \in K}
    \int_{a'' \in A''}
    |\det a''|^{n'}
    \int_{k'' \in K'' }
    \Phi( a'' k'' a' k)
    \, d k''
    \, \frac{d a''}{\delta_{N''}(a'')}
    \, d k
    \, \frac{d a'}{\delta_N (a')}.
  \end{equation*}
  Recombining the $A''$ and $K''$ integrals, we arrive at
  \begin{equation*}
    \int_{a' \in A'}
    \int_{k \in K}
    \int_{g \in N'' \backslash M''}
    |\det g|^{n'}
    \Phi(g a' k)
    \, d g
    \, d k
    \, \frac{d a'}{\delta_N (a')},
  \end{equation*}
  as required.
\end{proof}

\begin{lemma}\label{lem:reduce-Q-P-to-Q}
  For $(\phi, f, c) \in C_c^\infty(V) \times \mathcal{S}(N \backslash G) \times A''$, we have
  \begin{equation*}
    Q_P(\phi, f, c)
    = \int_{a ' \in A'}
    \int_{k \in K}
    Q^{M''}\left(\phi, |\det|^{n'/2 }R(a' k) f, c\right) \, \frac{d a '}{\delta_{N}(a')} \, d k.
  \end{equation*}
\end{lemma}
\begin{proof}
  Using the identity $\delta_N(c) = \lvert \det c \rvert^{- n'} \delta_{N''}(c)$, we see that
  \begin{equation*}
    |\det g|^{n'/2}
    W_P(f, c, g)
    =
    W^{M''}(|\det|^{n'/2} f,c,g).
  \end{equation*}
  Integrating, we obtain
  \begin{equation*}
    \int_{g \in N'' \backslash M''}
    \phi(e_n g)
    |\det g|^{n'}
    \left| W_P(f,c,g) \right|^2 \, d g
    =
    Q^{M''}\left(\phi, |\det|^{n'/2} f, c\right).
  \end{equation*}
  We replace $f$ by its translate $R(a' k) f$, as in the desired identity, and integrate.  We conclude by applying Lemma \ref{lem:standard:rewrite-N-G-via-M} with $\Phi(g) := \phi(e_n g) \left| W_P(f,c,g) \right|^2$.
\end{proof}

We write $\mathcal{S}_{d_1,d_2}^{N'' \backslash M''}$ and $\mathcal{S}_{d_1,d_2,T}^{N'' \backslash M''}$ for the Sobolev norms on $\mathcal{S}(N'' \backslash M'')$ defined like in \S\ref{sec:local-sobolev-norms-prelims}, but with $N'' \backslash M''$ playing the role of $N \backslash G$.  Thus
\begin{align*}
  \mathcal{S}_{d_1,d_2}^{N'' \backslash M''}(f)^2
  &=
    \sum_{
    \substack{
    x_1, \dotsb, x_{d_1} \in \mathcal{B}(A'') \cup \{1\}  \\
  y_1, \dotsc, y_{d_2 } \in \mathcal{B}(M'') \cup \{1\}
  }
  }
  \| L(x_1 \dotsb x_{d_1}) R(y_1 \dotsb y_{d_2}) f \|^2, \\
  \mathcal{S}_{d_1,d_2,T}^{N'' \backslash M''}(f)
  &=
    \max_{0 \leq k \leq d_2}
    T^{-k}
    \mathcal{S}_{d_1,k}^{N'' \backslash M''}(f).
\end{align*}
We can apply these just as well to restrictions of elements of $\mathcal{S}(N \backslash G)$.

\begin{lemma}\label{lem:induct-sobolev-norms}
  Let $d_1, d_2 \in \mathbb{Z}_{\geq 0}$ be fixed.  Then for all $f \in \mathcal{S}(N \backslash G)$,
  \begin{equation*}
    \int_{a' \in A'}
    \int_{k \in K}
    \mathcal{S} _{d_1, d_2} ^{N'' \backslash M''}
    \left(
      |\det|^{n'/2}
      R(a' k) f
    \right)^2
    \, d k
    \, \frac{d a '}{\delta_N(a')}
    \ll
    \mathcal{S}_{d_1,d_2}(f)^2.
  \end{equation*}
\end{lemma}
\begin{proof}
  The LHS is, by definition,
  \begin{align*}
    &\sum_{
      \substack{
      x_1, \dotsb, x_{d_1} \in \mathcal{B}(A'') \cup \{1\}  \\
    y_1, \dotsc, y_{d_2 } \in \mathcal{B}(M'') \cup \{1\}
    }
    }
    \int_{a' \in A'}
    \int_{k \in K}
    \int_{g \in N'' \backslash M''} \\
    &\quad \left\lvert    L(x_1 \dotsb x_{d_1}) R(y_1 \dotsb y_{d_2}) |\det|^{n'/2} R(a' k) f
      \right\rvert^2(g)
      \, d g
      \, d k
      \, \frac{d a '}{\delta_N(a')}.
  \end{align*}
  Using the product rule for derivatives and the fact that $A'$ commutes with $M''$, we see that the above is
  \begin{align*}
    \ll &\sum_{
        \substack{
        x_1, \dotsb, x_{d_1} \in \mathcal{B}(A'') \cup \{1\}  \\
    y_1, \dotsc, y_{d_2 } \in \mathcal{B}(M'') \cup \{1\}
    }
    }
    \int_{a' \in A'}
    \int_{k \in K}
    \int_{g \in N'' \backslash M''} \\
      &\quad \left\lvert  |\det|^{n'/2} R(a')   L(x_1 \dotsb x_{d_1}) R(y_1 \dotsb y_{d_2}) R(k) f
        \right\rvert^2(g)
        \, d g
        \, d k \, \frac{d a '}{\delta_N(a')}.
  \end{align*}
  We have
  \begin{equation*}
    L(x_1 \dotsb x_{d_1}) R(y_1 \dotsb y_{d_2}) R(k)
    =
    R(k) L(x_1 \dotsb x_{d_1}) R(\Ad(k)^{-1} y_1 \dotsb \Ad(k)^{-1} y_{d_2}).
  \end{equation*}
  Since $K$ is compact, we obtain the further majorization
  \begin{align*}
    \ll &\sum_{
        \substack{
        x_1, \dotsb, x_{d_1} \in \mathcal{B}(A) \cup \{1\}  \\
    y_1, \dotsc, y_{d_2 } \in \mathcal{B}(M) \cup \{1\}
    }
    }
    \int_{a' \in A'}
    \int_{k \in K}
    \int_{g \in N'' \backslash M''} \\
      &\quad \left\lvert  |\det|^{n'/2} R(a' k)   L(x_1 \dotsb x_{d_1}) R(y_1 \dotsb y_{d_2}) f
        \right\rvert^2(g)
        \, d g
        \, d k
        \, \frac{d a '}{\delta_N(a')}.
  \end{align*}
  By Lemma \ref{lem:standard:rewrite-N-G-via-M}, this last expression is just $\mathcal{S}_{d_1,d_2}(f)^2$.
\end{proof}

\begin{corollary}\label{cor:standard:induct-sobolev-norms}
  Let $T \geq 1$.  Let $d_1, d_2 \in \mathbb{Z}_{\geq 0}$ be fixed.  Then for all $f \in \mathcal{S}(N \backslash G)$,
  \begin{equation*}
    \int_{a' \in A'}
    \int_{k \in K}
    \mathcal{S} _{d_1, d_2,T} ^{N'' \backslash M''}
    \left(
      |\det|^{n'/2}
      R(a' k) f
    \right)^2
    \, \frac{d a '}{\delta_N(a')} \, d k
    \ll
    \mathcal{S}_{d_1,d_2,T}(f)^2.
  \end{equation*}
\end{corollary}

\subsection{Bounds at the real place: degenerate integrals}\label{sec:bounds-at-real-2}
For the following results, we take $\phi \in C_c^\infty(V)$ satisfying Assumptions \ref{sec:bounds-at-real-assumptions-phi}, $f \in \mathcal{S}(N \backslash G)$, $D > 0$, $T \geq 1$ and $c \in A''$.  We establish estimates for the degenerate integrals $Q_P$ analogous to those established in \S\ref{sec:bounds-at-real-1} for the Rankin--Selberg integrals $Q$.  For the proofs, we express $Q_P$ in terms of $Q^{M''}$ using Lemma \ref{lem:reduce-Q-P-to-Q} to reduce to the bounds established in \S\ref{sec:bounds-at-real-1}.

\begin{lemma}\label{lem:sub-gln:there-fixed-x-3}
  There is a fixed $d_1 \in \mathbb{Z}_{\geq 0}$ so that, if $f(a' a'' k) \neq 0 \implies |\det a''| \asymp D$, then
  \begin{equation*}
    Q_P(\phi,f,c) \ll \frac{r^{n''}}{D |\det c|} \mathcal{S}_{d_1,0}(f)^2.
  \end{equation*}
\end{lemma}
\begin{proof}
  Let $(a', k) \in A' \times K$.  Then the restriction to $M''$ of $R(a' k) f$ is supported on elements of determinant $\asymp D$.  Thus by Lemma \ref{lem:standard:there-fixed-d_1-1} (applied to $M''$ in place of $G$), we have
  \begin{equation*}
    Q^{M''}(\phi, |\det|^{n'/2 }R(a' k) f, c)
    \ll
    \frac{r^{n''}}{D |\det c|}
    \mathcal{S}_{d_1,0}^{N'' \backslash M''}\left(|\det |^{n'/2} R(a' k) f\right)^2.
  \end{equation*}
  By Lemma \ref{lem:reduce-Q-P-to-Q}, it follows that
  \begin{equation*}
    Q_P(\phi, f, c) \ll \frac{r^{n''}}{D |\det c|} J,
  \end{equation*}
  where
  \begin{equation*}
    J :=
    \int_{a ' \in A'}
    \int_{k \in K}
    \mathcal{S}_{d_1,0}^{N'' \backslash M''}\left(|\det|^{n'/2} R(a' k) f\right)^2
    \, d k
    \, \frac{d a '}{\delta_N(a')}.
  \end{equation*}
  By Lemma \ref{lem:induct-sobolev-norms}, we have $J \ll \mathcal{S}_{d_1,0}(f)^2$.  The required bound follows.
\end{proof}

For the following results, we introduce the notation
\begin{equation*}
  |\det|_{M''} : N \backslash G \rightarrow \mathbb{R}_{>0}
\end{equation*}
for the unique function that is invariant under $A' \times K$ and whose restriction to $M''$ is $|\det|$.  Explicitly, each $g \in N \backslash G$ may be represented in the form $g = a' a'' k$ with $(a',a'',k) \in A' \times A'' \times K$, in which case
\begin{equation*}
  |\det|_{M''}(g) := |\det a''|.
\end{equation*}

\begin{lemma}\label{lem:sub-gln:there-fixed-x-1}
  There is a fixed $d_1 \in \mathbb{Z}_{\geq 0}$ with the following property.  Let $b \in \mathbb{Z}_{\geq 0}$ be fixed and large enough in terms of $n$.  There is a fixed $d_2 \in \mathbb{Z}_{\geq 0}$ so that
  \begin{equation*}
    Q_P(\phi,f,c)
    \ll
    \frac{r^{n'' b}}{|\det c|^b}
    \mathcal{S}_{d_1, d_2}\left(|\det|_{M''}^{-b/2} f\right)^2.
  \end{equation*}
\end{lemma}
\begin{proof}
  By Lemma \ref{lem:sub-gln:there-fixed-x}, we have
  \begin{equation*}
    Q^{M''}\left(\phi, |\det|^{n'/2} R(a' k) f, c\right)
    \ll
    \frac{r^{n'' b}}{|\det c|^b}
    \mathcal{S}_{d_1, d_2}^{N'' \backslash M''}\left(|\det|_{M''}^{n'/2 - b/2} R(a' k) f\right)^2.
  \end{equation*}
  We have
  \begin{equation*}
    |\det|_{M''}^{n'/2 - b/2} R(a' k) f
    =
    |\det|_{M''}^{n'/2} R(a' k)
    |\det|_{M''}^{- b/2}
    f.
  \end{equation*}
  Applying Lemma \ref{lem:reduce-Q-P-to-Q} and Lemma \ref{lem:induct-sobolev-norms} as in the proof of Lemma \ref{lem:sub-gln:there-fixed-x-3} yields the required estimate.
\end{proof}

\begin{lemma}\label{lem:sub-gln:there-fixed-x-2}
  Assume that $n'' \geq 2$.  Let $b \in \mathbb{Z}_{\geq 0}$ be fixed and large enough in terms of $n$.  There are fixed $d_1,d_2 \in \mathbb{Z}_{\geq 0}$ so that, if $f(a ' a'' k) \neq 0 \implies |\det a''| \asymp D$, then for each fixed $\ell \in \mathbb{Z}_{\geq 0}$, we have
  \begin{align*}
    &Q_P(\phi,f,c)
    \\
    &\ll \frac{T^{d_1} r^{n'' b}}{|\det c|^b}
      \min \left( 1,
      \frac{r^{n''} T^{n''(n''-1)/2}}{D |\det c|}
      \right)^{\frac{4}{n''(n''-1)} \ell }
      \mathcal{S}_{d_1, d_2 +\ell, T}\left(|\det|_{M''}^{-b/2} f\right)^2.
  \end{align*}
\end{lemma}
\begin{proof}
  We reduce to Lemma \ref{lem:standard:there-fixed-d_1} exactly as in the proofs of the previous two lemmas, using Corollary \ref{cor:standard:induct-sobolev-norms} to estimate the resulting Sobolev norms.
\end{proof}

\subsection{Basic estimates for $\mathfrak{E}(N \backslash G, T)$}\label{sec:cnq7s0k2or}
We now prepare to apply the above estimates to elements of the classes $\mathfrak{E}(N \backslash G, T)$ and their translates.  Here we undertake the necessary preparatory work by recording basic estimates for elements of such classes.  We assume here that $T \ggg 1$.

\begin{lemma}\label{lem:standard:let-f-in-2}
  Let $f \in \mathfrak{E}(N \backslash G, T)$.  We may find a subset $\Sigma$ of $A$ and a decomposition
  \begin{equation}\label{eq:f-=-fclubsuit}
    f = f^{\clubsuit} + \sum_{\mu \in \Sigma } f^{\mu}
  \end{equation}
  with the following properties.
  \begin{enumerate}[(i)]
  \item $f^{\clubsuit} \in T^{-\infty}\mathfrak{E}(N \backslash G, T)$ and each $f^{\mu} \in \mathfrak{E}(N \backslash G, T)$.
  \item \label{itm:standard:clubsuit-2} For all fixed $x \in \mathfrak{U}(A)$ and $y \in \mathfrak{U}(G)$, and all $g \in N \backslash G$, we have
    \begin{equation*}
      L(x) R(y) f^{\clubsuit}(g) \ll \|g\|_{N \backslash G}^{-\infty} T^{-\infty}.
    \end{equation*}
  \item \label{itm:standard:decomposition-support-3} $\# \Sigma \ll (\log T)^{\O(1)}$.
  \item \label{itm:standard:decomposition-support-4} For each $\mu \in \Sigma$, we have
    \begin{equation*}
      \|\mu T^{\rho^\vee}\| = T^{o(1)},
    \end{equation*}
    or equivalently,
    \begin{equation}\label{eq:mu_1-=-tfrac1}
      \mu_1 = T^{\frac{1-n}{2} + o(1)}, \quad
      \mu_2 = T^{\frac{3-n}{2} + o(1)},
      \quad
      \dotsc,
      \quad
      \mu_n = T^{\frac{n-1}{2} + o(1)}.
    \end{equation}
    In particular, $|\det \mu| = T^{o(1)}$.
  \item \label{itm:standard:decomposition-support-5} For each $\mu \in \Sigma$, $a \in A$ and $k \in K$ with $f^{\mu}(a k) \neq 0$, we have
    \begin{equation*}
      a_i \asymp \mu_i
    \end{equation*}
    for each $i \in \{1, \dotsc, n\}$.  In particular,
    \begin{equation*}
      f^{\mu}(g) \neq 0 \implies \det g \asymp \det \mu.
    \end{equation*}
  \end{enumerate}
\end{lemma}
\begin{proof}
  For $a \in A$, define $r(a)  = (r_1(a),\dotsc,r_n(a)) \in \mathbb{R}^n$ by writing
  \begin{equation*}
    a T^{\rho^\vee} = \diag(\pm T^{r_1(a)}, \dotsc, \pm T^{r_n(a)}).
  \end{equation*}
  Explicitly,
  \begin{equation*}
    r(a) =
    \left(
      \frac{\log |a_1|}{\log T} + \frac{n - 1}{2},
      \frac{\log |a_2|}{\log T} + \frac{n - 3}{2},
      \dotsc,
      \frac{\log |a_n|}{\log T} + \frac{1 - n}{2},
    \right).
  \end{equation*}
  Since $r$ is invariant under $A \cap K$, the function on $N \backslash G$ given in Iwasawa coordinates $(a,k) \in A \times K$ by the formula $a k \mapsto r(a)$ is well-defined.  Fix $\phi \in C_c^\infty(\mathbb{R}^n)$ taking the value $1$ in some neighborhood of the origin.  Then $\phi ' := 1 - \phi$ vanishes in some neighborhood of the origin.  For each $\eps > 0$, we define $\Phi_\eps , \Phi_\eps ' \in C^\infty(N \backslash G)$ in Iwasawa coordinates
  \begin{equation*}
    \Phi_\eps (a k ) := \phi(r(a)/\eps), \quad
    \Phi_\eps '(a k) = \phi '(r(a) / \eps),
  \end{equation*}
  so that $\Phi_\eps  + \Phi_\eps  ' = 1$.

  For each $\eps > 0$, define $f^{\clubsuit, \eps} : N \backslash G \rightarrow \mathbb{C}$ by
  \begin{equation*}
    f^{\clubsuit, \eps }(g) := \Phi'_\eps (g) f(g).
  \end{equation*}

  Suppose for the moment that $\eps > 0$ is fixed.

  It is not hard to see that for fixed $(x,y) \in \mathfrak{U}(A) \times \mathfrak{U}(G)$, we have for all $g \in N \backslash G$
  \begin{equation}\label{eq:lx-ry-phi_eps}
    L(x) R(y) \Phi_\eps (g) \ll T^{\O(1)}, \quad
    L(x) R(y) \Phi_\eps '(g) \ll T^{\O(1)}.
  \end{equation}
  Indeed, $\Phi_\eps$ and $\Phi_\eps'$ are defined in Iwasawa coordinates by composition with the fixed smooth functions $\phi$ and $\phi'$ and the map $a \mapsto r(a)$, and $r$ is smooth and $A\cap K$-invariant.  Since $\phi$ and $\phi'$ have fixed compact support, the functions $\Phi_\eps$ and $\Phi_\eps'$ are supported on $\|T^{\rho^\vee}a\| = T^{\O(1)}$, hence on $\|a\|_{N\backslash G} = T^{\O(1)}$.  The bounds \eqref{eq:lx-ry-phi_eps} then follow from repeated application of the chain rule.

  We now bound $f$ on the support of $\Phi_\eps '$.  Let $(a,k) \in A \times K$.  Suppose $\Phi_{\eps}'(a k) \neq 0$.  Then $r(a)$ lies outside some fixed neighborhood of the origin, hence $\|T^{\rho^\vee} a\| \geq T^{\eps '}$ for some fixed $\eps ' > 0$.  By the characterization of $\mathfrak{E}(N \backslash G, T)$ given in part \eqref{itm:lemma-class-frak-E-basic:3} of Lemma \ref{lem:standard:frakE-equivalences}, it follows that
  \begin{equation*}
    \int_{k \in K}
    \left\lvert
      L(x) R(y) f(a k)
    \right\rvert^2 \, d k
    \ll
    \delta_N(a)
    \|T^{\rho^\vee} a\|^{-\infty}
    T^{\O(1)}
    \ll \|a k\|_{N \backslash G}^{-\infty} T^{-\infty}
  \end{equation*}
  for all fixed $(x,y) \in \mathfrak{U}(A) \times \mathfrak{U}(G)$.  By the Sobolev lemma for $K$, we deduce that on the support of $\Phi_\eps'$,
  \begin{equation}\label{eq:lx-ry-fg}
    L(x) R(y) f(g) \ll \|g\|_{N \backslash G}^{-\infty} T^{-\infty}.
  \end{equation}

  Combining \eqref{eq:lx-ry-phi_eps} and \eqref{eq:lx-ry-fg} and appealing to the product rule for derivatives, we deduce that for all fixed $\eps > 0$ and fixed $(x,y) \in \mathfrak{U}(A) \times \mathfrak{U}(G)$, we have
  \begin{equation}\label{eq:rx-ly-f}
    L(x) R(y)
    f^{\clubsuit, \eps }(g)
    \ll \|g\|_{N \backslash G}^{-\infty} T^{-\infty}.
  \end{equation}
  By overspill, we may find some $\eps > 0$ with $\eps \lll 1$ for which the conclusion \eqref{eq:rx-ly-f} remains valid and for which the estimate  \eqref{eq:lx-ry-phi_eps} holds in the slightly stronger form
  \begin{equation*}
    L(x) R(y) \Phi_\eps (g) \ll T^{o(1)},
  \end{equation*}
  say.
  (Indeed, we know from \eqref{eq:rx-ly-f} that the set of $\eps > 0$ for which
  \begin{equation*}
    \sup_{
      \substack{
        d_1,d_2 \in \mathbb{Z} _{\geq 0} :   d_1, d_2 \leq 1/\eps \\
        x_1,\dotsc, x_{d_1} \in \mathcal{B}(A) \\
        y_1,\dotsc, y_{d_2} \in \mathcal{B}(G) \\
        g \in N \backslash G}
    }
    \|g\|_{N \backslash G}^{1/\eps}
    \left\lvert L(x_1 \dotsb x_{d_1}) R(y_1 \dotsb y_{d_2}) f^{\clubsuit,\eps}(g) \right\rvert \leq T^{-1/\eps}
  \end{equation*}
  and
  \begin{equation*}
    \sup_{
      \substack{
        d_1,d_2 \in \mathbb{Z} _{\geq 0} :   d_1, d_2 \leq 1/\eps \\
        x_1,\dotsc, x_{d_1} \in \mathcal{B}(A) \\
        y_1,\dotsc, y_{d_2} \in \mathcal{B}(G) \\
        g \in N \backslash G}
    }
    \left\lvert L(x_1 \dotsb x_{d_1}) R(y_1 \dotsb y_{d_2}) \Phi_\eps (g)  \right\rvert \leq T^{\eps}
  \end{equation*}
  contains all fixed $\eps > 0$, hence must contain some $\eps \lll 1$.)  We henceforth work with such an $\eps$, and set $f^{\clubsuit} := f^{\clubsuit, \eps }$.  It is clear from \eqref{eq:rx-ly-f} and Lemma \ref{lem:standard:frakE-equivalences} that $f^{\clubsuit}$ belongs to $T^{-\infty} \mathfrak{E}(N \backslash G, T)$ and satisfies the required estimate \eqref{itm:standard:clubsuit-2}.  Moreover, since \eqref{eq:rx-ly-f} persists after multiplying $f^{\clubsuit}$ by any fixed power of $T$, we in fact have $f^{\clubsuit} \in T^{-\infty}\mathfrak{E}(N \backslash G, T)$.

  It remains to show that the function $f - f^{\clubsuit}$, given by
  \begin{equation*}
    (f - f^{\clubsuit})(g) =  \Phi_\eps(g) f(g),
  \end{equation*}
  admits a decomposition $\sum_{\mu \in \Sigma} f^\mu$ of the required form.  For this, we appeal to dyadic decomposition, as follows.  Let $\Lambda$ denote the discrete subgroup of $A$ consisting of those elements each of whose coordinates lie in the discrete subgroup $\exp(\mathbb{Z})$ of $\mathbb{R}^\times$.  We may find a fixed test function $\eta \in C_c^\infty(A)$, invariant under $A \cap K$, so that $\sum_{\mu \in \Lambda} \eta(\mu^{-1} a) = 1$ for all $a \in A$.  We define $f^{\mu} \in C^\infty(N \backslash G)$ by, for $(a,k) \in A \times K$,
  \begin{equation*}
    f^{\mu}(a k) := \eta(\mu^{-1} a) \Phi_\eps(a) f(a k).
  \end{equation*}
  By construction, the identity \eqref{eq:f-=-fclubsuit} holds.  Set $\Sigma := \left\{ \mu \in \Lambda : f^\mu \neq 0 \right\}$.  From the support condition on $\Phi_\eps$, we see that $f^\mu \neq 0$ only if $\|\mu T^{\rho^\vee}\| = T^{o(1)}$, whence assertion \eqref{itm:standard:decomposition-support-4}.  In particular, the cardinality bound \eqref{itm:standard:decomposition-support-3} holds.  The required support condition \eqref{itm:standard:decomposition-support-5} for $f^\mu$ follows from fact that $\eta$ has fixed compact support.

  It remains only to verify that $f^{\mu} \in \mathfrak{E}(N \backslash G, T)$.
  The function $N \backslash G \ni a k \mapsto \eta(\mu^{-1} a)$ and its fixed $(A \times G)$-derivatives have magnitude $\O(1)$, uniformly in $\mu$.  By \eqref{eq:lx-ry-phi_eps}, the corresponding fixed derivatives of $\Phi_\eps$ are $\ll T^{\O(1)}$ (and for the particular $\eps \lll 1$ chosen above, even $\ll T^{o(1)}$).  Since $f \in \mathfrak{E}(N \backslash G, T)$, the characterization of $\mathfrak{E}$ given in part \eqref{itm:lemma-class-frak-E-basic:3} of Lemma \ref{lem:standard:frakE-equivalences}, together with the product rule for derivatives, implies that $f^{\mu} = \eta(\mu^{-1}\cdot)\,\Phi_\eps\,f$ also lies in $\mathfrak{E}(N \backslash G, T)$.
\end{proof}

Let $\sigma \in \mathfrak{a}^*$.  It defines a positive-valued character $A \rightarrow \mathbb{R}^\times_+$, $a \mapsto |a|^{\sigma}$.  We denote temporarily by
\begin{equation*}
  \phi_\sigma : N \backslash G \rightarrow \mathbb{R}^\times_+
\end{equation*}
the function given in Iwasawa coordinates by $\phi_\sigma(a k) := |a|^{\sigma}$ for all $(a,k) \in A \times K$.  For example, the functions $|\det|_{M''}^{b}$ with $b \in \mathbb{R}$ are of this form.

Recall from  \S\ref{sec:local-sobolev-norms-prelims} the definitions and basic properties of the Sobolev norms $\mathcal{S}_{d_1,d_2}$ and their rescaled variants $\mathcal{S}_{d_1,d_2,T}$.

\begin{lemma}\label{lem:lemma-Sobolev-estimates-frakE}
  Let $T \ggg 1$.  Let $f \in \mathfrak{E}(N \backslash G, T)$.  Let $Y \in A^0$ with each component satisfying $Y_j \in [T^{-1/2}, T^{1/2}]$.  Let $\sigma$ belong to some fixed compact subset of $\mathfrak{a}^*$.  Then:
  \begin{enumerate}[(i)]
  \item \label{item:sobolev-frakE-0} For all fixed $d_1,d_2 \in \mathbb{Z}_{\geq 0}$,
    \begin{equation}\label{eq:mathc-d_2-phi_s}
      \mathcal{S}_{d_1,d_2} \left(\phi_{\sigma} L(Y) f\right) \ll T^{- \langle \rho^\vee , \sigma  \rangle + d_2 + o(1)} Y^{-\sigma}.
    \end{equation}
  \item \label{item:sobolev-frakE-1} For all fixed $d_1, d_2 \in \mathbb{Z}_{\geq 0}$,
    \begin{equation}\label{eq:mathc-d_2-det_m}
      \mathcal{S}_{d_1,d_2} (\phi_\sigma  L(Y) f) \ll T^{\O(1)}.
    \end{equation}
  \item \label{item:sobolev-frakE-2} There is a fixed $C_0 \geq 0$ so that for all fixed $d_1, d_2 \in \mathbb{Z}_{\geq 0}$,
    \begin{equation*}
      \mathcal{S}_{d_1,d_2,T} (\phi_{\sigma} L(Y) f) \ll T^{C_0}.
    \end{equation*}
  \item \label{item:sobolev-frakE-3} For all fixed $d_1 \in \mathbb{Z}_{\geq 0}$,
    \begin{equation*}
      \mathcal{S}_{d_1,0}(L(Y) f) \ll T^{o(1)}.
    \end{equation*}
  \end{enumerate}
\end{lemma}
\begin{proof}
  We recall (from \S\ref{sec:local-sobolev-norms-prelims}) that the LHS of \eqref{eq:mathc-d_2-det_m} is the mean square, over $x = x_1 \dotsb x_{d_1}$ and $y = y_1 \dotsb y_{d_2}$ with $x_i \in \mathcal{B}(A) \cup \{1\}$ and $y_j \in \mathcal{B}(G) \cup \{1\}$, of the $L^2(N \backslash G)$-norms
  \begin{equation*}
    \|L(x) R(y) \phi_{\sigma} L(Y) f\|.
  \end{equation*}
  For any such $x$ and $y$, we have
  \begin{equation*}
    \left(
      L(x) R(y) \phi_{\sigma}
    \right)(g)
    \ll \phi_{\sigma}(g).
  \end{equation*}
  (Indeed, it is enough to check that the class of functions $\phi_{\sigma,f} : N \backslash G \rightarrow \mathbb{C}$ of the form $\phi_{\sigma,f}(a k) := |a|^{\sigma} f(k)$, where $f : (A \cap K) \backslash K \rightarrow \mathbb{C}$ has some fixed $K$-type and $L^\infty$-norm $\ll 1$, is preserved under $L(x)$ and $R(y)$ for $x \in \Lie(A)$ and $y \in \Lie(G)$.  The required conclusion concerning $L(x)$ follows from the homomorphism property of $|.|^{\sigma}$, that for $R(y)$ by considering the Iwasawa decomposition of $\Ad(k) y$.)  By the product rule for derivatives, we deduce that
  \begin{align*}
    &\mathcal{S}_{d_1,d_2} (\phi_\sigma  L(Y) f)^2 \\
    &\ll
      \max_{
      \substack{
      x,y:  \\
    \text{as above}
    }
    }
    \int_{g \in N \backslash G}
    \delta_N^{-1}(Y)
    \phi_\sigma(g)^2
    \left\lvert
    L(x) R(y) f
    \right\rvert^2(Y g) \, d g.
  \end{align*}
  Motivated by the concentration properties of $f$, we now write
  \begin{equation*}
    \phi_\sigma(g) = T^{-\langle \rho^\vee, \sigma  \rangle} Y^{-\sigma} \phi_\sigma(T^{\rho^\vee} Y g)
  \end{equation*}
  and substitute $g \mapsto Y^{-1} g$ in the integral, giving
  \begin{equation*}
    (T^{-\langle \rho^\vee, \sigma  \rangle} Y^{-\sigma})^2
    \int_{g \in N \backslash G}
    \phi_{\sigma} (T^{\rho^\vee} g)^2
    \left\lvert
      L(x) R(y) f
    \right\rvert^2(g) \, d g.
  \end{equation*}
  Our hypotheses on $\sigma$ give $\phi_{\sigma} (T^{\rho^\vee} g) \ll \|T^{\rho^\vee} g\|_{N \backslash G}^{\O(1)}$, so by Lemma \ref{lem:f-in-mathfraken-square-integral-frakE}, we deduce assertion \eqref{item:sobolev-frakE-0}.  Assertion \eqref{item:sobolev-frakE-1} follows by noting that the RHS of \eqref{eq:mathc-d_2-phi_s} is $T^{\O(1)}$ under our hypotheses on $\sigma$ and $Y$.  Assertion \eqref{item:sobolev-frakE-2} similarly follows upon recalling, from \eqref{eq:mathcals_d_1-d_2-tf}, the definition $\mathcal{S}_{d_1,d_2,T}(f') = \max_{0 \leq k \leq d_2} T^{-k} \mathcal{S}_{d_1,k}(f')$, upon taking $C_0$ large enough in terms of the fixed compact to which $\sigma$ belongs.  Assertion \eqref{item:sobolev-frakE-3} follows by specializing \eqref{eq:mathc-d_2-phi_s} to $(d_2,\sigma) = (0,0)$.
\end{proof}

\subsection{The main local estimate}\label{sec:20230514074956}

\begin{proposition}\label{prop:sub-gln:let-f-in}
  Let $f \in \mathfrak{E}(N \backslash G, T)$.  Let $Y \in A^0$ with each component $Y_j$ satisfying
  \begin{equation}\label{eq:t-kappa-leq}
    T^{-\kappa} \leq |Y_j| \leq T^{\kappa}
  \end{equation}
  for some fixed $\kappa < 1/2$.  Let $r > 0$ with $r = T^{\O(1)}$.  Retain Assumptions \ref{sec:bounds-at-real-assumptions-phi} concerning $\phi$.  Let $c \in A''$ with $|\det c| \geq 1$.  Then
  \begin{equation}\label{eq:qpphi-ly-f}
    Q_P(\phi, L(Y) f, c)
    \ll
    \min
    \left(
      \frac{r^n T^{o(1)} \det(Y)}{|\det (c)|} ,
      \frac{T^{\O(1)}}{|\det c|^{10}}
    \right)
    +
    \frac{T^{-\infty }}{|\det c|^{10}}
  \end{equation}
\end{proposition}
We have abused notation slightly in our use of ``$\min$.''  The estimate \eqref{eq:qpphi-ly-f} means, more precisely, that there is a fixed $C_0 \geq 0$ so that for each fixed $C_1 \geq 0$ and $\eps > 0$, we have
\begin{equation*}
  Q_P(\phi, L(Y) f, c)
  \ll
  \min
  \left(
    \frac{r^n T^{\eps} \det(Y)}{|\det (c)|},
    \frac{T^{C_0}}{|\det c|^{10}}
  \right)
  +
  \frac{T^{-C_1 }}{|\det c|^{10}}.
\end{equation*}

On a first reading, the reader is encouraged to pretend that $Y = 1$ (since, to deduce a qualitative subconvex bound, it is enough to treat the case that each $Y_j = T^{o(1)}$).  The variable $r$, on the other hand, should be regarded as being potentially quite small (e.g., of size $T^{(1-n)/2}$).

The proof of Proposition \ref{prop:sub-gln:let-f-in} requires all of the main lemmas of \S\ref{sec:bounds-at-real-2}, and proceeds in several steps.  We first decompose
\begin{equation}\label{eq:f-=-f}
  f = f^{\clubsuit} + \sum f^\mu
\end{equation}
according to Lemma \ref{lem:standard:let-f-in-2}.  Since the number of summands is $(\log T)^{\O(1)} \ll T^{o(1)}$, we may reduce via the Minkowski-type inequality
\begin{equation*}
  Q_P(\phi, f_1 + \dotsb + f_k, c)
  \leq
  k
  \sum_{j=1}^k
  Q_P(\phi, f_j, c)
\end{equation*}
to estimating the contribution to $Q_P(\phi,L(Y) f,c)$ of each term on the RHS of \eqref{eq:f-=-f}.  We treat these separately in the lemmas below.

We first treat the contribution of $f^{\clubsuit}$.  The following estimate is acceptable for the purposes of proving Proposition \ref{prop:sub-gln:let-f-in}.
\begin{lemma}\label{lemma:20230514162239}
  Let hypotheses be as in Proposition \ref{prop:sub-gln:let-f-in}; in particular, $\lvert \det c \rvert \geq 1$.  Then
  \begin{equation*}
    Q_P(\phi, L(Y) f^{\clubsuit},c) \ll
    \frac{T^{-\infty}}{|\det c|^{10}}.
  \end{equation*}
\end{lemma}
\begin{proof}
  In view of our assumption $r = T^{\O(1)}$, it is enough to show that for some fixed $b \geq 10$ taken sufficiently large in terms of $n$,
  \begin{equation*}
    Q_P(\phi, L(Y) f^{\clubsuit},c) \ll
    \frac{r^{n'' b}}{|\det c|^{b}} T^{-\infty},
  \end{equation*}
  say.  By Lemma \ref{lem:sub-gln:there-fixed-x-1}, we reduce to showing that any fixed $d_1,d_2 \in \mathbb{Z}_{\geq 0}$ and $b \in \mathbb{R}$,
  \begin{equation}\label{eq:estimate-clubsuit-contribution-trivially}
    \mathcal{S}_{d_1, d_2}(|\det|_{M''}^{-b/2} L(Y) f^{\clubsuit})
    \ll T^{-\infty},
  \end{equation}
  which follows in turn from Lemma \ref{lem:lemma-Sobolev-estimates-frakE}, part \eqref{item:sobolev-frakE-1} in view of the membership $f^{\clubsuit} \in T^{-\infty} \mathfrak{E}(N \backslash G, T)$.
\end{proof}

We next focus on an individual term $f^{\mu}$ in \eqref{eq:f-=-f}.  We will show that the contribution of $f^\mu$ is majorized by the first term on the RHS of \eqref{eq:qpphi-ly-f}.  To see this, we must establish upper bounds involving each of the two alternatives in the minimum defining that term.  The required bounds are given by Lemma \ref{lem:standard:we-have-begin-first-alternative} and Proposition \ref{prop:standard:we-have-begin-second-alternative}, below.

We start with the (simpler) second alternative in the minimum.
\begin{lemma}\label{lem:standard:we-have-begin-first-alternative}
  We have
  \begin{equation*}
    Q_P(\phi, L(Y) f^\mu , c)
    \ll
    \frac{T^{\O(1)}}{|\det c|^{10}}.
  \end{equation*}
\end{lemma}
\begin{proof}
  Using our assumptions $r = T^{\O(1)}$ and $\lvert \det c \rvert \geq 1$, it is enough to show that for some fixed $b \geq 10$ taken sufficiently large in terms of $n$, there is a fixed $C_0 \geq 0$ so that
  \begin{equation*}
    Q_P(\phi, L(Y) f^\mu, c) \ll
    \frac{r^{n'' b}}{|\det c|^{b}} T^{C_0}.
  \end{equation*}
  By Lemma \ref{lem:sub-gln:there-fixed-x-1}, it suffices to verify that for each fixed $d_{1}, d_2 \in \mathbb{Z}_{\geq 0}$ and $b \in \mathbb{R}$, we have
  \begin{equation}\label{eq:mathc-d_2-left}
    \mathcal{S}_{d_1,d_2} \left(
      |\det|_{M''}^{-b/2} L(Y) f^\mu
    \right)^2 \ll T^{\O(1)}.
  \end{equation}
  This estimate follows in turn from Lemma \ref{lem:lemma-Sobolev-estimates-frakE}, part \eqref{item:sobolev-frakE-1}.
\end{proof}

It remains to address the (crucial) first alternative in the minimum.
\begin{proposition}\label{prop:standard:we-have-begin-second-alternative}
  We have
  \begin{equation}\label{eq:qpphi-ly-fmu}
    Q_P(\phi, L(Y) f^\mu , c)
    \ll
    \frac{r^n T^{o(1)} \det(Y)}{|\det (c)|}.
  \end{equation}
\end{proposition}
The proof divides into a couple steps.  For both steps, it is relevant to recall that by the construction of $f^\mu$, we have
\begin{equation*}
  f^\mu (a' a'' k) \neq 0 \implies |\det a''| \asymp \prod_{j=n' + 1}^{n} \mu_j,
\end{equation*}
hence
\begin{equation}\label{eq:ly-fmua-a}
  L(Y) f^\mu(a' a'' k) \neq 0 \implies  \det(Y'') |\det a''| \asymp D , \quad D := \prod_{j=n'+1}^{n} \mu_j,
\end{equation}
where we write $Y = Y' Y''$ with $(Y', Y'') \in A' \times A''$.

\begin{lemma}\label{lemma:20230514162116}
  We have
  \begin{equation}\label{eq:q_pphi-f-c}
    Q_P(\phi, L(Y) f^\mu,c) \ll \frac{r^{n''} T^{o(1)}  \det (Y'') }{D |\det (c)| }.
  \end{equation}
\end{lemma}
\begin{proof}
  For fixed $d_1 \in \mathbb{Z}_{\geq 0}$ we have by \eqref{eq:mathcals_d_1-d_2la-f}, followed by Lemma \ref{lem:lemma-Sobolev-estimates-frakE}, part \eqref{item:sobolev-frakE-3} that

  \begin{equation*}
    \mathcal{S}_{d_1,0}(L(Y) f^{\mu})
    \ll
    T^{o(1)}.
  \end{equation*}
  The required estimate then follows from Lemma \ref{lem:sub-gln:there-fixed-x-3} and \eqref{eq:ly-fmua-a}.
\end{proof}

\begin{lemma}\label{lemma:20230514162141}
  There is a fixed $C_1 \geq 0$ so that for each fixed $\ell \in \mathbb{Z}_{\geq 0}$, we have
  \begin{equation}\label{eq:q_pphi-f-c-1}
    Q_P(\phi,L(Y) f^\mu,c) \ll T^{C_1}
    \min \left( 1,
      \frac{r^{n''} T^{n''(n''-1)/2} \det(Y'')}{D |\det (c)|}
    \right)^{\ell }.
  \end{equation}
\end{lemma}

\begin{proof}
  We apply Lemma \ref{lem:lemma-Sobolev-estimates-frakE}, part \eqref{item:sobolev-frakE-2} to see that for each fixed $b \in \mathbb{R}$, there is a fixed $C_0 \geq 0$ so that for all fixed $d_1,d_2 \in \mathbb{Z}_{\geq 0}$, we have
  \begin{equation}\label{eq:mathc-d_2-tdet_m}
    \mathcal{S}_{d_1, d_2, T}(|\det|_{M''}^{-b/2} L(Y) f^{\mu}) \ll T^{C_0}.
  \end{equation}
  In the case $n'' \geq 2$, we fix $b$ large enough (in terms of $n$) and Lemma \ref{lem:sub-gln:there-fixed-x-2} with this $b$.  The required estimate then follows (possibly after replacing $\ell$ by a multiple) from \eqref{eq:mathc-d_2-tdet_m} and our hypotheses $r = T^{\O(1)}$ and $|\det c| \geq 1$.

  In the case $n'' = 1$, we argue more directly (cf.\ \S\ref{sec:cnjgglqbxq}).  For $c \in A''$, we have $W_P(L(Y) f^\mu, c, g) = \delta_N^{1/2}(c Y^{-1}) f^\mu(Y c^{-1} g)$, hence
  \begin{equation}\label{eq:cnjfy7zeut}
    Q_P(\phi, L(Y) f^\mu, c) =
    \delta_N(c Y^{-1}) \int_{N \backslash G} \phi(e_n g) \lvert f^\mu(Y c^{-1} g) \rvert^2 \, d g.
  \end{equation}
  We write the integral using Iwasawa coordinates $g = a k$.  By the support conditions on $\phi$ and $f^\mu$, we see that the resulting integrand is supported on $a_n c_n \asymp r Y_n$ and $a_n \asymp \mu_n = D$, so the integral vanishes unless $r Y_n / D c_n \asymp 1$, as we may assume.  We have $c_n = \det (c)$ and $Y_n = \det(Y'')$, so $r Y_n / D c_n$ is precisely the ratio appearing in \eqref{eq:q_pphi-f-c-1}, which is thus $\asymp 1$.  Since each of $Y$, $a_n$ and $r$ is bounded polynomially in $T$, the same holds for $c_n$, hence also for $\delta_N(c Y^{-1})$.  We reduce to verifying that $Q_P(\phi, L(Y) f^\mu , c) \ll T^{\O(1)}$, which follows from (e.g.) \eqref{eq:q_pphi-f-c}.
\end{proof}

For the purposes of proving \eqref{eq:qpphi-ly-fmu}, the estimate \eqref{eq:q_pphi-f-c} is adequate provided that
\begin{equation*}
  \frac{r^{n''} T^{o(1)}  \det (Y'') }{D |\det (c)| } \ll \frac{r^n T^{o(1)} \det(Y)}{|\det (c)|},
\end{equation*}
which simplifies to
\begin{equation}\label{eq:rn-ll-d}
  r^{n''}\det(Y'') \ll D r^n T^{o(1)} \det(Y).
\end{equation}
In particular, in the ``nondegenerate'' case $P = G$, we have $\det(Y'') = \det(Y)$, $n'' = n$ and $D = T^{o(1)}$, so \eqref{eq:rn-ll-d} always holds and the proof is now complete.  We may thus suppose henceforth that
\begin{equation}\label{eq:p-neq-g}
  P \neq G,
  \quad
  n' \geq 1,
  \quad
  r^{n''} \det(Y'') \ggg D r^n \det(Y).
\end{equation}
In view of the estimate \eqref{eq:q_pphi-f-c-1}, our task reduces to the following calculation:

\begin{lemma}\label{lem:under-stat-cond}
  Under the stated conditions \eqref{eq:p-neq-g}, we have
  \begin{equation}\label{eq:fracrn-tnn-12d}
    \frac{r^{n''} T^{n''(n''-1)/2} \det (Y'')}{D}
    \ll T^{-\eps}
  \end{equation}
  for some fixed $\eps > 0$.
\end{lemma}
To explicate the sufficiency of this estimate, observe by \eqref{eq:q_pphi-f-c-1}, \eqref{eq:fracrn-tnn-12d} and the hypothesis $|\det c| \geq 1$ that
\begin{equation*}
  Q_P(\phi, L(Y) f,c) \ll |\det c|^{-\ell} T^{C_1 - \ell \eps + o(1)}
\end{equation*}
for each fixed $\ell \in \mathbb{Z}_{\geq 0}$.  Taking $\ell \geq 10$ sufficiently large then gives
\begin{equation*}
  Q_P(\phi,L(Y) f,c) \ll |\det c|^{-10} T^{-\infty},
\end{equation*}
which is adequate for the purposes of proving \eqref{eq:qpphi-ly-fmu} in view of our assumptions $r = T^{\O(1)}$ and $\det(Y'') = T^{\O(1)}$.
\begin{proof}[Proof of Lemma \ref{lem:under-stat-cond}]
  We first employ \eqref{eq:p-neq-g} to see that
  \begin{equation*}
    r^{n'} \lll \frac{1}{D \det(Y')}.
  \end{equation*}
  Our task \eqref{eq:fracrn-tnn-12d} thereby reduces to verifying that
  \begin{equation}\label{eq:fracdetynn-tnn-12dnn}
    \frac{T^{n''(n''-1)/2} \det (Y'') }{D^{n/n'} \det(Y')^{n''/n'}  }
    \ll T^{-\eps}.
  \end{equation}
  From our hypothesis \eqref{eq:t-kappa-leq} on $Y$ and the condition $\kappa < 1/2$, we know that
  \begin{equation*}
    \frac{\det(Y'')}{\det(Y')^{n''/n'}} \leq
    T^{n'' - \eps}
  \end{equation*}
  for some fixed $\eps > 0$ (depending only upon $\kappa$ and $n$).  On the other hand, from the definition \eqref{eq:ly-fmua-a} of $D$ and the estimates \eqref{eq:mu_1-=-tfrac1} for the $\mu_j$, we have
  \begin{equation*}
    D = T^{\sum_{j=n'+1}^{n} (2 j - n - 1)/2 + o(1)}.
  \end{equation*}
  By the combinatorial identity
  \begin{equation*}
    \sum_{j=n'+1}^{n}
    \frac{2j -n - 1}{2} =
    \frac{n ' n''}{2},
  \end{equation*}
  it follows that
  \begin{equation*}
    D^{n / n'} =
    T^{n n''/2 + o(1)}.
  \end{equation*}
  Since $n' \geq 1$, we have $n'' + 1 \leq n$, and so
  \begin{equation*}
    n''
    +
    \frac{n'' (n''-1)}{2}
    =
    \frac{n'' (n''+1)}{2}
    \leq
    \frac{n n''}{2},
  \end{equation*}
  which yields \eqref{eq:fracdetynn-tnn-12dnn}.
\end{proof}

\section{Completion of the proof}\label{sec:completion-proof-1}
Our aim is to complete the proof of Proposition \ref{prop:bound-after-theta-degenerate-case-2}.  We retain Assumptions \ref{assumptions:furth-reduct-proof}.  We resume the global discussion of \S\ref{sec:reduction-proof}, retaining the notation of that section.  In particular,
\begin{itemize}
\item $G$ is $\GL_n$ over $\mathbb{Z}$ (or for the present purposes, over $\mathbb{Q}$)
\item $P = M U_P$ is either $G$ itself or a standard maximal parabolic,
\item $M = M' \times M''$, $A', N' \leq M'$ and $A'', N'' \leq M''$ are as before, and
\item $\mathcal{Q}_P(\phi, f)$ was defined in \eqref{eq:mathcalq_pphi-f-:=}, motivated by Proposition \ref{prop:bound-after-theta-degenerate-case-2} (our remaining goal), and given a preliminary estimate in Lemma \ref{lem:sub-gln:let-phi-in}.
\end{itemize}
We also set $(n', n'') := (\rank(M'), \rank(M''))$, as in \S\ref{sec:degenerate-integrals}.

\subsection{Reduction of global bounds to archimedean bounds}
For $m \in \mathbb{N} := \mathbb{Z}_{\geq 1}$, we denote by $\tau(m)$ the number of divisors of $m$ and by $\Sigma^P(m)$ the set of all $c \in A''(\mathbb{Q})$, identified with $n''$-tuples of nonzero rational numbers, with the following properties:
\begin{enumerate}[(i)]
\item Each entry $c_j$ is an integer.
\item $|\det(c)| = \prod |c_j| = m$.
\end{enumerate}
Clearly $|\Sigma^P(m)| \ll \tau(m)^{\O(1)}$.

We have assumed that $f \in \mathfrak{E}(N(\mathbb{R}) \backslash G(\mathbb{R}), T)$, but the following lemma is valid for any $f \in \mathcal{S}(N(\mathbb{R}) \backslash G(\mathbb{R}))$.
\begin{lemma}\label{lem:scratch-research:we-have-begin}
  We have
  \begin{equation}\label{eq:mathc-f_infty-ll}
    \mathcal{Q}_P(\phi,f) \ll \sum_{m \in \mathbb{N}}
    \tau(m)^{\O(1)}
    \max_{c \in \Sigma^P(m)}
    Q_P(\phi_\infty, f, c).
  \end{equation}
\end{lemma}
\begin{proof}
  With $W_P(\dotsb)$ as in \S\ref{sec:degenerate-integrals}, let us introduce the notation
  \begin{equation*}
    W_\mathfrak{p}(c,g) :=
    W_P(\Theta^P[f]_\mathfrak{p}, c, g),
  \end{equation*}
  \begin{equation*}
    I_\infty(b,c) := \int_{N(\mathbb{R}) \backslash G(\mathbb{R}) } \phi_\infty(g) W_\infty(b,g) \overline{W_\infty(c,g)} \, d g,
  \end{equation*}
  \begin{equation*}
    I_p(b,c) := \int_{a \in A''(\mathbb{Q}_p)}
    1_{|a_n| = 1}
    W_p(b,a)
    \overline{W_p(c,a)}
    \, \frac{d a}{\delta_{N(\mathbb{Q}_p)}(a)}.
  \end{equation*}
  With this notation, Lemma \ref{lem:sub-gln:let-phi-in} reads
  \begin{equation*}
    \mathcal{Q}_P(\phi, f) \leq 2^{n '}
    \sum_{b, c \in A ''(\mathbb{Q})}
    \prod_{\mathfrak{p}} I_{\mathfrak{p}}(b, c),
  \end{equation*}
  so by the triangle inequality, we have
  \begin{equation*}
    \mathcal{Q}_P(\phi,f)
    \ll
    \sum_{b,c \in A''(\mathbb{Q})}
    \prod_\mathfrak{p}
    \left\lvert
      I_\mathfrak{p}(b,c)
    \right\rvert,
  \end{equation*}
  Using the support and invariance properties of $\Theta^P[f]_p$, we may rewrite
  \begin{equation*}
    I_p(b, c) = \int_{N \backslash G}
    1 _{ \lVert e_n g \rVert = 1}
    W_p(b, g) \overline{W_p(c, g)} \, d g.
  \end{equation*}
  By the bilinearization of Lemma 21.19, we obtain
  \begin{equation*}
    I_p(b, c) = \int_{N''(\mathbb{Q}_p) \backslash M''(\mathbb{Q}_p)}
    1 _{ \lVert e_n g \rVert = 1} W^{M''}(\Theta^P[f]_p, b, g)
    \overline{W^{M''}(\Theta^P[f]_p, c, g)}
    \, d g.
  \end{equation*}
  The remaining integral is the subject of Lemma \ref{lem:sub-gln:b-c-in}, but applied on the smaller general linear group $M''$ (cf.\ \eqref{eq:thet-=-1_am}).  The results of that lemma give that $\prod_p I_p(b,c) = 0$ unless $b,c \in \Sigma^P(m)$ for some $m \in \mathbb{N}$, in which case $\prod_p I_p(b,c) \ll \tau(m)^{\O(1)}$.  Thus
  \begin{equation*}
    \mathcal{Q}_P(\phi,f)
    \ll
    \sum_{m \in \mathbb{N} }
    \tau(m)^{\O(1)}
    \sum_{b,c \in \Sigma^P(m)} |I_\infty(b,c)|.
  \end{equation*}
  By the triangle inequality, Cauchy--Schwarz and the assumed nonnegativity of $\phi_\infty$, we have
  \begin{align*}
    \sum_{b,c \in \Sigma^P(m)}
    |I_\infty(b,c)|
    &\leq
      \sum_{b, c \in \Sigma^P(m)}
      \int_{g \in N(\mathbb{R}) \backslash G(\mathbb{R})}
      \phi_\infty(e_n g)
      |W_\infty(b,g) \overline{W_\infty}(c,g)|
      \, d g
    \\
    &\leq
      |\Sigma^P(m)|
      \sum_{ c \in \Sigma^P(m)}
      \int_{g \in N(\mathbb{R}) \backslash G(\mathbb{R})}
      \phi_\infty(e_n g)
      |W_\infty(c,g)|^2
      \, d g
  \end{align*}
  Since $\# \Sigma^P(m) \ll \tau(m)^{\O(1)}$, the desired estimate follows.
\end{proof}

\subsection{Proof of Proposition \ref{prop:bound-after-theta-degenerate-case-2}}
Assume now once again that $f \in \mathfrak{E}(N(\mathbb{R}) \backslash G(\mathbb{R}), T)$.  We must show that
\begin{equation*}
  \mathcal{Q}_P(\phi,L(Y) f)
  \ll
  r^{n} T^{o(1)} \det(Y).
\end{equation*}
The bound \eqref{eq:qpphi-ly-f} for $Q_P(\phi_\infty, L(Y) f, c)$ afforded by Proposition \ref{prop:sub-gln:let-f-in} gives, for some fixed $C_0 \geq 0$ and all $c \in \Sigma^P(m)$,
\begin{equation}\label{eq:qpphi-ly-f-1}
  Q_P(\phi_\infty, L(Y) f, c)
  \ll
  \min
  \left(
    \frac{r^n T^{o(1)} \det(Y) }{m},
    \frac{T^{C_0}}{m^{10}}
  \right)
  +
  \frac{T^{-\infty }}{m^{10}}.
\end{equation}
We insert this bound into Lemma \ref{lem:scratch-research:we-have-begin}.  The contribution to \eqref{eq:mathc-f_infty-ll} of the second term appearing on the RHS of \eqref{eq:qpphi-ly-f-1} is
\begin{equation*}
  \ll
  T^{-\infty}
  \sum_{m \in \mathbb{N} }
  \frac{\tau(m)^{\O(1)}}{m^{10}}
  \ll T^{-\infty}.
\end{equation*}
It remains to estimate the contribution of the first term.  To that end, we let $S \geq 2$ be an arbitrary parameter, and apply the first or second alternative in the minimum of \eqref{eq:qpphi-ly-f-1} according as $m \leq S$ or $m > S$.  We obtain
\begin{equation*}
  \mathcal{Q}_P(\phi, L(Y) f)
  \ll
  r^n T^{o(1)} \det(Y)
  \sum_{m \leq S}
  \frac{\tau(m)^{\O(1)}}{ m}
  +
  T^{C_0}
  \sum_{m > S}
  \frac{\tau(m)^{\O(1)}}{m^{10}}
  + T^{-\infty}.
\end{equation*}
We take $S = T^{C_1}$, where $C_1$ is fixed but sufficiently large in terms of $C_0$.  Then
\begin{equation*}
  \sum_{m \leq S} \frac{\tau(m)^{\O(1)}}{ m} \ll \log T \ll T^{o(1)},
\end{equation*}
while
\begin{equation*}
  \sum_{m > S} \frac{\tau(m)^{\O(1)}}{m^{10}} \ll S^{-8} = T^{-8 C_1}.
\end{equation*}
Since $C_1$ can be taken arbitrarily large, it follows that
\begin{equation*}
  \mathcal{Q}_P(\phi, L(Y) f)
  \ll
  r^n T^{o(1)} \det(Y) + T^{-\infty}.
\end{equation*}
This estimate is adequate in view of our assumptions that $r$ and $\det(Y)$ are of the form $T^{\O(1)}$.

\appendix
\section{Elementary formulations}\label{sec:elem-form}

Here we apply the algorithm described in \cite[\S2]{MR469763}, and illustrated earlier in \S\ref{sec:main-result:-reform} and in the proof of Lemma \ref{lem:four-transf-begin}, to reformulate some of the key supporting results of the paper in elementary language, free from nonstandard formalism.  (We also simplify the statements a bit using, e.g., Lemma \ref{lem:overspill-A-vs-T-eps}.)  We do this for ``information purposes'' only, with the hope of making more accessible the structure of our argument.  This discussion is not logically necessary for our aims, so we do not prove that we have carried out the reformulations accurately.

\begin{theorem}[Theorem \ref{thm:main-local-results}, reformulated]\label{thm:construction-vectors-original}
  Suppose given
  \begin{itemize}
  \item a general linear pair $(G,H)$ over $\mathbb{R}$, with accompanying notation $A$, $B$, $N$, $Q$, $U$, $Z$, $A_H$, $B_H$, $N_H$, $Q_H$, $U_H$, $Z_H$,
  \item a nondegenerate unitary character $\psi$ of $N$, equal to either the standard character or its inverse,
  \item $C_0 > c_0 > 0$,
  \item a compact subset $\Omega_H$ of $H$,
  \item an open neighborhood $E_G \subseteq G$ of the identity element, and
  \item $\eps > 0$.
  \end{itemize}
  There exists
  \begin{itemize}
  \item a compact neighborhood $\Omega_H '$ of $H$,
  \item nonnegative reals
    \begin{itemize}
    \item $C_{\mathcal{D},\ell}$, indexed by compact subsets $\mathcal{D} \subseteq \mathfrak{a}^*$ and $\ell \in \mathbb{Z}_{\geq 0}$,
    \item $C_{u, \ell}$, indexed by $u \in \mathfrak{U}(G)$ and $\ell \in \mathbb{Z}_{\geq 0}$,
    \end{itemize}
  \item $c_1 > 0$,
  \item $T_0 \geq 1$, and
  \item a finite subset $\mathcal{F}_{Z_H} \subseteq C_c^\infty(Z_H)$
  \end{itemize}
  such that
  \begin{itemize}
  \item for each $T \geq T_0$,
  \item for each generic irreducible unitary representation $\pi$ of $G$, each of whose archimedean $L$-function parameters has magnitude between $c_0 T$ and $C_0 T$,
  \end{itemize}
  there exists
  \begin{itemize}
  \item $f, f_{\ast} \in \mathcal{S}^e(U_H \backslash H)^{W_H}$, $\phi_{Z_H} \in \mathcal{F}_{Z_H}$,
  \item a smooth vector $W \in \mathcal{W}(\pi,\psi)$ of norm $\leq 1$ (with respect to the norm defined in \S\ref{sec:local-prelim-kirillov-model}), and
  \item $\omega_0 \in C_c^\infty(E_G)$
  \end{itemize}
  with the following properties.  Define
  \begin{equation*}
    \omega(g) := (\omega_0 \ast \omega_0^*)(g) := \int_{h \in G} \omega_0(g h^{-1}) \overline{\omega_0(h^{-1})} \, d h,
  \end{equation*}
  \begin{equation*}
    \omega^\sharp (g) := \int_{z \in Z}
    \pi|_Z (z) \omega(z g) \, d z,
  \end{equation*}
  where $\pi|_Z$ denotes the central character of $\pi$.  Then:
  \begin{enumerate}[(i)]
  \item $f_{\ast}(g) = \int_{z \in Z_H} f(g z^{-1}) \phi_{Z_H}(z) \, d z$.
  \item Each of the following assertions remains valid upon replacing $f$ with $f_{\ast}$.
  \item For each compact subset $\mathcal{D}$ of $\mathfrak{a}^*$ and $\ell \in \mathbb{Z}_{\geq 0}$, we have, with notation as in \S\ref{sec:local-sobolev-norms-prelims},
    \begin{equation*}
      \nu_{\mathcal{D},\ell,T}(f) \leq C_{\mathcal{D},\ell} T^\eps.
    \end{equation*}
  \item The local zeta integral $Z(W,f[s])$ is entire.  We have the lower bound
    \begin{equation*}
      |Z(W,f[0])| \geq c_1 T^{- \dim(H)/4}.
    \end{equation*}
    For each compact subset $\mathcal{D} \subseteq \mathfrak{a}^*$, $\ell \in \mathbb{Z}_{\geq 0}$, and $s \in \mathfrak{a}_{\mathbb{C}}^*$ with $\Re(s) \in \mathcal{D}$, we have the upper bound
    \begin{equation*}
      |Z(W,f[s])| \leq C_{\mathcal{D},\ell} T^{\rank(G) \trace(s)/2 - \frac{\dim(H)}{4} + \eps }
      (1 + |s|)^{-\ell}.
    \end{equation*}
  \item For each $u \in \mathfrak{U}(G)$ and $\ell \in \mathbb{Z}_{\geq 0}$, we have
    \begin{equation*}
      \|\pi(u) ( \pi(\omega_0) W - W)\| \leq C_{u,\ell} T^{-\ell}.
    \end{equation*}
  \item We have $\int_{H} |\omega^\sharp| \leq T^{\rank(H)/2 + \eps}$.
  \item Let $\Psi_1, \Psi_2 : H \rightarrow \mathbb{C}$ be continuous functions.  For $j=1,2$, define the convolutions
    \begin{equation*}
      (\Psi_j \ast \phi_{Z_H})(x) := \int_{z \in Z_H} \Psi_j(x z^{-1}) \phi_{Z_H}(z) \, d z.
    \end{equation*}
    Let $\gamma \in \bar{G} - H$.  Then the integral
    \begin{equation*}
      I := \int_{x, y \in \Omega_H }
      \left\lvert
        (\Psi_1 \ast \phi_{Z_H})(x)
        (\Psi_2 \ast \phi_{Z_H})(y)
        \omega^\sharp(x^{-1} \gamma y)
      \right\rvert
      \, d x \, d y
    \end{equation*}
    satisfies the estimate
    \begin{equation*}
      |I| \leq
      T^{\rank(H)/2 + \eps}
      \min\left( 1,
        \frac{T^{-1/2}}{ d_H(\gamma)}
      \right)
      \|\Psi_1\|_{L^2(\Omega_H')}
      \|\Psi_2\|_{L^2(\Omega_H')}.
    \end{equation*}
  \end{enumerate}
\end{theorem}

\begin{theorem}[Theorem \ref{thm:growth-eisenstein-nonstandard}, reformulated]
  Let $n$ be a natural number, $G := \GL_n$ over $\mathbb{Z}$.  Let $\Omega \subseteq G(\mathbb{R})$ be a compact set.  Let $\kappa \in (0, 1/2), \eps > 0$.  There exists $C \geq 0$, $\mathcal{D} \subseteq \mathfrak{a}^*$ compact, and $\ell \in \mathbb{Z}_{\geq 0}$ with the following property.  Let $T \geq 1$.  Let $Y \in A(\mathbb{R})^0$ with $T^{-\kappa} \leq |Y_j| \leq T^{\kappa}$ for all $j$.  Let $f \in \mathcal{S}^e(U(\mathbb{R}) \backslash G(\mathbb{R}))^W$.  Let $t \in A(\mathbb{R})^0_{\geq 1}$.  Then
  \begin{equation*}
    \frac{\|\Eis [ \Phi[L(Y) f] ] \|_{t,\Omega}^2}{\delta_N(t)}
    \leq C T^\eps \min \left( \det(Y)^{-1} t_1^{-n}, \det(Y) t_n^n \right) \nu_{\mathcal{D},\ell,T}(f).
  \end{equation*}
\end{theorem}

\begin{theorem}[Theorem \ref{thm:big-group-whittaker-behavior-under-G}, reformulated]
  Let $n$ be a natural number.  Set $G := \GL_n(\mathbb{R})$.  For each $\eps > 0$, there exists $\delta_0 < 1/2$ with the following property.  Let $\delta \in (\delta_0, 1/2)$, $\delta_+ \in (\delta, 1/2)$ and $\eta > 0$.  Let $\mathcal{X}$ be a finite subset of $\mathfrak{U}(G)$.  Let $m \geq 0$.  Let $\Lambda$ be a compact subset of $[\mathfrak{g}^\wedge]$.  There is a finite subset $\mathcal{Y}$ of $\mathbb{Z}_{\geq 0}^{\dim(Q_H)}$, $k \geq 0$ and $T_0 \geq 1$ with the following property.  Let $T \geq T_0$.  Let $\pi$ be a generic irreducible unitary representation of $G$ whose rescaled infinitesimal character satisfies $T^{-1} \lambda_\pi \in \Lambda$.  There exists $\omega_0 \in C_c^\infty(G)$ with the following properties.  Set $\bar{G} = G/Z$.  Define $\omega := \omega_0 \ast \omega_0^\ast \in C_c^\infty(G)$ and $\omega^\sharp \in C_c^\infty(\bar{G})$ by
  \begin{equation*}
    \omega(g) := \int_{h \in G} \omega_0(g h^{-1}) \overline{\omega_0(h^{-1})} \, d h,
  \end{equation*}
  \begin{equation*}
    \omega^\sharp (g) := \int_{z \in Z}
    \pi|_Z (z) \omega(z g) \, d z,
  \end{equation*}
  where $\pi|_Z$ denotes the central character of $\pi$.  Then:
  \begin{enumerate}[(i)]
  \item $\omega(g) \neq 0$ only if $|g -1| \leq \eta$ and $|\Ad^*(g) \tau - \tau| \leq \eta T^{-1/2}$, where $\tau = \tau(\pi,\psi,T)$ is given by \eqref{eq:tau-characterize-via-pi}.
  \item $\|\omega^\sharp\|_{\infty} \leq T^{n(n+1)/2+\eps}$.
  \item $\int_H |\omega^\sharp| \leq T^{n/2 + \eps}$.
  \item Let $W$ be an element of the $\psi_T$-Whittaker model $\mathcal{W}(\pi,\psi_T)$.  Suppose that $W|_{Q_H}$ is supported near the identity element, so that it may be described via pullback to the Lie algebra.  Let $\phi : \Lie(Q_H) \rightarrow \mathbb{C}$ be given by $\phi(x) := W(\exp(x))$.  Assume that $\phi$ is supported on $\{x : |x| \leq T^{\delta_+ - 1}\}$, where $|.|$ denotes (say) the Euclidean norm on $\Lie(Q_H)$ defined with respect to the matrix entries.  Assume that for each $\alpha \in \mathcal{Y}$ and all $x \in \Lie(Q_H)$, we have
    \begin{equation*}
      \left\lvert \partial^\alpha \phi(x) \right\rvert \leq T^{(1 - \delta) (\dim(Q_H) + |\alpha|)} \left( 1 + \frac{|x|}{T^{\delta-1}} \right)^{-k}.
    \end{equation*}
    Then for all $t\in \mathcal{X}$, we have
    \begin{equation*}
      \|\pi(t) (\pi(\omega_0)W - W)\| \leq T^{-m},
    \end{equation*}
    where $\|.\|$ denotes the norm on $\mathcal{W}(\pi,\psi_T)$ defined by integration over $N_H \backslash H$.
  \end{enumerate}
\end{theorem}

\printindex

\bibliography{refs}{} \bibliographystyle{alpha}


\end{document}